FROM STRINGS TO SETS

A technical report

Zlatan Damnjanovic

University of Southern California

We introduce an elementary theory of concatenation, $QT^+$, and construct a formal interpretation of Adjunctive Set Theory with Extensionality in $QT^+$. The proofs are given in complete detail.

# Table of contents

## 5. Coding sets by strings

## 6. Lexical precedence

## 10. Minimal, special and canonical set codes

# 1. A Theory of Concatenation

We consider a first-order theory with identity and a single binary function symbol *. Informally, we let the variables range over nonempty strings of a's and b's – or 0's and 1's – and let x*y be the string that consists of the digits of the string x followed by the digits of string y, subject to the following conditions:

(QT1)   (x*y)*z = x*(y*z)

(QT2)   ¬(x*y=a) & ¬(x*y=b)

(QT3)   (x*a=y*a → x=y) & (x*b=y*b → x=y) &

         & (a*x=a*y → x=y) & (b*x=b*y → x=y)

(QT4)   ¬(a*x=b*y) & ¬(x*a=y*b)

(QT5)   x=a v x=b v (∃y(a*y=x v b*y=x) & ∃z(z*a=x v z*b=x))

It is convenient to have a function symbol for a successor operation on strings:

(QT6)     Sx=y ↔ ((x=a & y=b) v (¬x=a & x*b=y)).

Of course, we may also think of a single letter a appended to x, x*a, as a successor of the string x.

Since the last axiom is basically a definition, adding it to the other five results in an inessential (i.e. conservative) extension.

We call the resulting theory QT⁺. We let M be an arbitrary model of QT⁺.

We introduce some obvious abbreviations:

$$xBy \equiv \exists z\ x*z=y \quad \text{and} \quad xEy \equiv \exists z\ z*x=y.$$

Often, we write xy for x*y and omit parentheses in x(yz) and (xy)z on account of (QT1).

Let $\quad xRy \equiv (x=a\ \&\ \neg y=a)\ \vee\ xBy.$

We make the following observations:

(1.1) QT ⊢ ¬aRa & ¬bRb

(1.2) QT ⊢ aRb

(1.3) QT ⊢ ∀x (xRb ↔ x=a)

(1.4) QT ⊢ ∀x ¬xRa & ∀x (¬a=x → aRx)

(1.5) QT ⊢ ∀x (xRa ∨ x=a ↔ x=a)

(1.6) QT ⊢ ∀x ∀y ∀z (¬x=a → ( xRy → (z*x)R(z*y)))

(1.1) and (1.2) are immediate; for (1.3), we have ← by (1.2), and → follows from (QT2) because QT ⊢ ∀x ¬xBb; on the other hand, (1.4) is immediate from the definition and (QT2); then (1.5) is immediate from (1.4).

For (1.6), assume M ⊨ xRy for x≠a. Then M ⊨ xBy, that is, M ⊨ ∃u x*u=y. So M ⊨ z*(x*u)=z*y, whence M ⊨ (z*x)*u=z*y by (QT1). But this means M ⊨ (z*x)B(z*y), so M ⊨ (z*x)R(z*y).

We proceed to prove:

(1.7)  QT ⊢ xRy & yRz → xRz.

Suppose M ⊨ x≠a and M ⊨ y≠a . Then M ⊨ xRy implies M ⊨ ∃u x*u=y, and M ⊨ yRz implies M ⊨ ∃v y*v=z; so M ⊨ (x*u)*v=z, whence M ⊨ x*(u*v)=z by (QT1). Thus M ⊨ ∃w x*w=z, whence M ⊨ xRz. If M ⊨ y=a, then M ⊨ ¬xRy by (1.4) whether M ⊨ x=a or not, and (1.7) holds trivially.

Note that (1.7) along with (1.1) implies that QT ⊢ ¬∃y(aRy & yRa), and QT ⊢ ¬∃y ∃z(aRy & yRz & zRa), and similarly for R-cycles of any given length. Likewise for b.

We observe:

(1.8)  $QT^+ \vdash \forall x\, xR(Sx)$.

(1.9)  $QT^+ \vdash \forall x\, \forall y\, (Sx=Sy \rightarrow x=y)$.

For (1.8), if $M \vDash x=a$, then $M \vDash Sx=b$, and we appeal to (2); and if $M \vDash x \neq a$, then $M \vDash Sx=x*b$, so $M \vDash xB(Sx)$, whence $M \vDash xR(Sx)$ by definition.

For (1.9), the result is immediate if $M \vDash x=a=y$; if $M \vDash x=a$ and $M \vDash y \neq a$, then $M \vDash Sx=b$ and $M \vDash Sy=y*b$ and (1.9) holds trivially because $M \vDash \neg b=y*b$ by (QT2), and likewise if $M \vDash y=a$ and $M \vDash x \neq a$. So we may assume $M \vDash x \neq a$ and $M \vDash y \neq a$. But then $M \vDash Sx=x*b$ and $M \vDash Sy=y*b$, and we have (1.9) by (QT3).

We also have

(1.10)  $QT^+ \vdash \forall x\, \neg(Sx=a)$.

By (QT4), $M \vDash \neg(x=a\ \&\ a=b)$. By (QT2), $M \vDash \neg(\neg x=a\ \&\ x*b=a)$.

So  $QT^+ \vdash \neg(Sx=a)$.

Note that

(1.11) $QT^+ \vdash \forall x\, (xR(Sa) \leftrightarrow xRa \lor x=a)$

By (1.3), $M \vDash xRb \leftrightarrow x=a$, whereas $M \vDash x=a \leftrightarrow xRa \lor x=a$ by (1.4). It follows that $M \vDash xR(Sa) \leftrightarrow xRa \lor x=a$ because $M \vDash Sa=b$.

Next, we show that

(1.12) $QT \vdash \forall x\, \forall y\, \forall z\, (x*z=y*b \rightarrow xRy \lor x=y)$.

Suppose $M \vDash x*z=y*b$. Then $M \vDash z \neq a$ by (QT4), and for the same reason $M \vDash \neg aEz$. We now use (QT5). If $M \vDash z=b$, then $M \vDash x*b=y*b$, so $M \vDash x=y$ by (QT3). Otherwise, $M \vDash bEz$, that is, $M \vDash \exists u\, u*b=z$, so $M \vDash x*(u*b)=y*b$, whence $M \vDash (x*u)*b=y*b$ by (QT1). But then $M \vDash x*u=y$ by (QT3), that is, $M \vDash xBy$. So $M \vDash xRy$.

We can now generalize (1.11) to

(1.13)  QT⁺ ⊢ ∀x ∀y (xR(Sy) ↔ xRy v x=y)

First note that M ⊨ yR(Sy) by (1.8), and that M ⊨ xRy implies M ⊨ xR(Sy) by (1.8) and (1.7). This proves ←.

Suppose M ⊨ xR(Sy). If M ⊨ y=a, the claim holds by (1.11). If M ⊨ y≠a, then M ⊨ Sy=y*b. If M ⊨ x=a, we have M ⊨ xRy by (1.4) and (QT2). Otherwise, M ⊨ y≠a and M ⊨ x≠a. Then M ⊨ xB(Sy) from M ⊨ xR(Sy). So M ⊨ ∃z x*z=y*b. But then M ⊨ xRy v x=y by (1.12).

From (11) we immediately derive

(1.14)  QT⁺ ⊢ ∀x ∀y ((xR(Sy) v x=Sy) ↔ (xRy v x=y v x=Sy)).

Note that we have

   aRb, bR(Sb), (Sb)R(SSb), (SSb)R(SSSb) … as well as

   aRb, bR(b*a), (b*a)R((b*a)*b), ((b*a)*b)R((b*a)*b)*a …

but not bR(a*b) or (a*b)R(b*a).

These elementary facts tell us that, provably in QT⁺,

   xRy v x=y  is a discrete preordering of strings.

## 2. Tractable strings

Consider the following problem. We know that neither a nor b are their own initial segments:

$$QT^+ \vdash \neg aRa \quad \text{and} \quad QT^+ \vdash \neg bRb.$$

But we <u>don't know whether our theory proves that no string</u> is an initial segment of itself: $\quad QT^+ \vdash^? \forall x \, \neg xRx$.

Let $I_0(x) \equiv \forall y \, (yRx \vee y=x \to \neg yRy)$.

Unsurprisingly, $QT \vdash I_0(a)$ by (1.1) and (1.4).

Likewise, $QT \vdash I_0(b)$ by (1.1) and (1.3).

But it also turns out that $I_0(x)$ is closed under S, provably in $QT^+$.

Assume $M \vDash I_0(x)$. Suppose $M \vDash yR(Sx) \vee y=Sx$.

Then $M \vDash (yRx \vee y=x) \vee y=Sx$ by (1.14).

If $M \vDash yRx \vee y=x$, then $M \vDash \neg yRy$ follows from the hypothesis $M \vDash I_0(x)$.

What if $M \vDash y=Sx$ ?

If $M \vDash x=a$, then we have $M \vDash Sx=b$, and $M \vDash \neg yRy$ by (1.1).

So we may assume $M \vDash \neg x=a$.

By the definition of S, then $M \vDash Sx=x*b$. We now derive a reductio.

Suppose $M \vDash yRy$, that is, $M \vDash (Sx)R(Sx)$. Then $M \vDash (x*b)R(x*b)$.

Because $M \vDash \neg a=x*b$ by (QT2), we have in fact that $M \vDash (x*b)B(x*b)$, which means that $M \vDash \exists u\ (x*b)*u=(x*b)$. Now $M \vDash u=b \vee bEu$ by (QT5) and (QT4). If $M \vDash u=b$, then $M \vDash (x*b)*b=x*b$, so $M \vDash x*b=x$ by (QT3), whence $M \vDash xBx$ and thus $M \vDash xRx$, which plainly contradicts the hypothesis $M \vDash I_0(x)$.
If $M \vDash bEu$, then $M \vDash \exists v\ v*b=u$. Thus we have $M \vDash (x*b)*(v*b)=x*b$, so

$M \vDash ((x*b)*v)*b=x*b$ by (QT1). Then $M \vDash (x*b)*v=x$ by (QT3), so

$M \vDash x*(b*v)=x$ again by (QT1). But then $M \vDash xBx$, and $M \vDash xRx$, again contradicting $M \vDash I_0(x)$.

So $M \vDash \neg yRy$ under the overall assumption $M \vDash yR(Sx) \vee y=Sx$, and thus, in addition to the atoms a and b, the formula $I_0(x)$ provably applies to strings that are closed under S.

We remark that $I_0$ is also closed under the operation $S_a(x)=x*a$. First note that the analogue (1.8$^a$) of (1.8) is immediate, and (1.10$^a$) holds by (QT3). (The superscript $^a$ indicates that S is replaced by $S_a$.) For (1.11$^a$), we have to show that $QT \vdash \forall x\ (xR(a*a) \leftrightarrow xRa \lor x=a)$. This obviously holds for x=a, by (1.4) and (QT2). Furthermore, the ← part is immediate given (QT2). For →, suppose that $M \vDash xR(a*a)$.

Suppose $M \vDash x \neq a$. Then $M \vDash xB(a*a)$, whence $M \vDash x \neq b$ by (QT4).

Likewise, $M \vDash \neg bBx$. By (QT5), $M \vDash aBx$.

So $M \vDash \exists u\ a*u=x$, and we have $M \vDash (a*u)B(a*a)$, that is, $M \vDash \exists v\ (a*u)*v=a*a$. By (QT1), this gives $M \vDash a*(u*v)=a*a$, whence $M \vDash u*v=a$ by (QT3). But this contradicts (QT2).

So we have shown that $M \vDash xR(a*a) \to x=a$, which then yields the → part of (1.11$^a$). Now, (1.12$^a$) and (1.13$^a$) are proved exactly as (1.12) and (1.13) if we just replace b with a, and then (1.14$^a$) is immediate.

In the proof of the closure of $I_0$ under $S_a$ where we appeal to (1.14$^a$), the only difference is in the case $M \vDash x=a$ under the hypothesis $M \vDash y=x*a$. Assume that $M \vDash yRy$, that is $M \vDash (a*a)R(a*a)$. Then $M \vDash (a*a)B(a*a)$, because M ⊨ a≠a*a by (QT2). So M ⊨ ∃u (a*a)*u=(a*a). Then M ⊨ a*(a*u)=a*a by (QT1), and we have M ⊨ a*u=a by (QT3). But this, of course, contradicts (QT2). So M ⊨ ¬yRy if M ⊨ x=a and M ⊨ y=x*a. The rest of the argument is the same if in the reductio we replace b by a.

So we also have that QT ⊢ ∀x ($I_0(x)$ → $I_0(x*a)$).

Obviously,

(2.1)   QT ⊢ ∀x ($I_0$(x) → ¬xRx).

But note that  QT ⊢ xRy  & yRx → xRx  by (1.7),  so (2.1) implies

$$QT ⊢ ∀x ∀y(I_0(x) \& yRx → ¬xRy),$$

whence we have

(2.2)   QT ⊢ ∀x ∀y($I_0$(x) → ¬(xRy & yRx)).

But (2.1) and (1.8) also imply

(2.3)   $QT^+$ ⊢ ∀x ($I_0$(x) → ¬x=Sx).

So if we were to restrict our attention to the strings in $I_0$, then (1.4), (1.7), (1.8), (1.9), (1.13) and (2.1) give us a discrete partial ordering under R and a successor operation, with a as the least element.   Thus we may write

$$x<y ≡ I_0(y) \& xRy \quad \text{and} \quad x≤y ≡ (I_0(y) \& xRy) \lor x=y.$$

We observe that QT ⊢ $I_0$(aa) & $I_0$(ab) while  QT ⊢ ¬(aa≤ab) & ¬(ab≤aa).

For, suppose  QT ⊢ aa<ab.  Then QT ⊢ aaBab, by the definition of R, so

QT ⊢ ∃u aau=ab.  Then  QT ⊢ au=b  by (QT3), which contradicts (QT2).

Of course,  QT ⊢ ¬(aa=ab) by (QT4).  Hence  QT ⊢ ¬(aa≤ab).  Likewise for

QT ⊢ ¬(ab≤aa).

But we do have  QT ⊢ a<aa & a<ab  because     QT ⊢ aBaa & aBab.

(2.4)   $QT^+ \vdash \forall x\, (I_0(x) \to x<Sx)$.

Suppose $M \vDash I_0(x)$. Then $M \vDash I_0(Sx)$. Since from (1.8), $M \vDash xR(Sx)$, this suffices for $M \vDash x<Sx$ by definition of $<$.

## 3. String concepts

Call a formula J(x) in the language with a, b, * and S as primitives a <u>string concept</u> if

$QT^+ \vdash J(a)$,

$QT^+ \vdash \forall x\, (J(x) \to J(x*a))$,

$QT^+ \vdash \forall x\, (J(x) \to J(Sx))$.

Note that a conjunction of string concepts is also a string concept.

We established that $I_0(x)$ is a string concept. Obviously, so is $x=x$. But we did not know if $x=x$ had the crucial property expressed in (2.1). We do know that of $I_0(x)$, provably in $QT^+$.

The formula $I_0(x)$ earned the title "string concept" by being closed under the string successor operations, provably in $QT^+$. But do we know that the strings in $I_0$ are closed under *? We don't. We need a string concept with that property.

(3.1) There is a string concept $I_1$ such that

$$QT^+ \vdash \forall x \forall y\, (I_1(x) \& I_1(y) \rightarrow I_1(x*y))$$

where $QT \vdash \forall x\, (I_1(x) \rightarrow I_0(x))$.

Let $\qquad I_1(x) \equiv I_0(x)\ \&\ \forall y(\, I_0(y) \rightarrow I_0(y*x))$.

We need to verify that $I_1(x)$ is indeed a string concept. First, that $QT^+ \vdash I_1(a)$:

we have $QT^+ \vdash I_0(a)$. Suppose $M \vDash I_0(y)$. Then $M \vDash I_0(y*a)$ because $I_0$ is closed under $S_a$, provably in $QT^+$. So indeed $QT^+ \vdash I_1(a)$.

As for $QT^+ \vdash I_1(b)$, that follows from $QT^+ \vdash I_0(a)$ and the closure of $I_0$ under $S$.

Next we show that $QT^+ \vdash \forall z\, (I_1(z) \rightarrow I_1(z*b))$.

So suppose $M \vDash I_1(z)$. We want $M \vDash I_1(z*b)$. We have $M \vDash I_0(z)$ from the hypothesis $M \vDash I_1(z)$, and so $M \vDash I_0(z*b)$ since $I_0$ is a string concept. Assume that $M \vDash I_0(y)$. From the hypothesis $M \vDash I_1(z)$ it then follows that $M \vDash I_0(y*z)$, and further that $M \vDash I_0((y*z)*b)$. By (QT1), this means that $M \vDash I_0(y*(z*b))$. So we have established that

$$M \vDash \forall y(\, I_0(y) \rightarrow I_0(y*(z*b))),$$

which, along with the previously obtained $M \vDash I_0(z*b)$, gives us $M \vDash I_1(z*b)$ under the hypothesis $M \vDash I_1(z)$, as required.

Similarly, $QT^+ \vdash \forall z\, (I_1(z) \rightarrow I_1(z*a))$.

This completes the argument that $I_1(x)$ is a string concept.

We now show that $I_1$ is closed under the concatenation operation $*$, that is,

$$QT^+ \vdash \forall x\, \forall y\, (I_1(x)\, \&\, I_1(y) \to I_1(x*y)).$$

Assume $M \vDash I_1(x)$ and $M \vDash I_1(y)$, namely

(a) $\quad M \vDash I_0(x)\, \&\, \forall z(I_0(z) \to I_0(z*x))$,

and (b) $\quad M \vDash I_0(y)\, \&\, \forall z(I_0(z) \to I_0(z*y))$.

From $M \vDash I_0(x)$ and (b) we obtain $M \vDash I_0(x*y)$. Assume now $M \vDash I_0(z)$. Then $M \vDash I_0(z*x)$ by (a), and further $M \vDash I_0((z*x)*y)$ by (b). But then $M \vDash I_0(z*(x*y))$ by (QT1). So we have that

$$M \vDash I_0(x*y)\, \&\, \forall z(I_0(z) \to I_0(z*(x*y))),$$

that is, $M \vDash I_1(x*y)$. This completes the proof of (3.1).

Note that we have not used any property specific to $I_0$ as a string concept in the above argument. Say that a string concept I is <u>stronger than</u> $I_0$ if $QT^+ \vdash \forall x\, (I(x) \rightarrow I_0(x))$ and write $I \subseteq I_0$. We have in fact proved something more general:

(3.2)  For any string concept $I \subseteq I_0$ there is a string concept $J \subseteq I$ such that

$$QT^+ \vdash \forall x\, \forall y\, (J(x)\, \&\, J(y) \rightarrow J(x*y)).$$

(3.3) Suppose J⊆I is a string concept where I⊆I$_0$. Then there is a string concept J$_\leq$⊆J such that

$$QT^+ \vdash \forall x\, (J_\leq(x)\ \&\ y \leq x \rightarrow J_\leq(y)).$$

Let J$_\leq$(x) ≡ ∀y≤x J(y).

We write ∀y≤x ... for ∀y(y≤x → ...), and analogously for ∃.

That the formula J$_\leq$(x) has the required property is immediate from the definition and (1.7).

We have QT$^+$ ⊢ J(a) by hypothesis, and QT ⊢ y≤a ↔ y=a by (1.5). So QT$^+$ ⊢ J$_\leq$(a).

Suppose M ⊨ J$_\leq$(x). Then M ⊨ ∀y≤x J(y). Suppose M ⊨ y≤Sx. By (1.14) and (1.8), M ⊨ y≤Sx ↔ y≤x ∨ y=Sx.

If M ⊨ y≤x, we have M ⊨ J(y) from the hypothesis M ⊨ J$_\leq$(x).

If M ⊨ y=Sx, then M ⊨ J(x) from the hypothesis M ⊨ J$_\leq$(x), whence M ⊨ J(Sx) from the principal hypothesis. Therefore, M ⊨ ∀y≤Sx J(y), that is, M ⊨ J$_\leq$(Sx).

That J$_\leq$(x) is closed under S$_a$ is proved in the same fashion except that we appeal to (1.14$^a$) instead.

This completes the argument that J$_\leq$(x) is a string concept and the proof of (3.3).

On account of (3.2) and (3.3), in establishing that a given string concept I may be strengthened to a string concept J with another property, we will not bother to argue that the formula J(x) is also closed with respect to * or downward closed with respect to ≤; by (3.2) and (3.3), we can always strengthen J(x) to one that is.

(3.4)   For any string concept I⊆I₀ there is a string concept J⊆I such that

$$QT \vdash \forall x \in J\ \neg xEx.$$

Let $J(x) \equiv I(x)\ \&\ \forall y\ \neg yx = x$.

Since $I(x)$ is a string concept, $QT \vdash I(a)$ and $QT \vdash I(b)$. By (QT2),

$QT \vdash \forall y\ \neg ya = a$ and $QT \vdash \forall y\ \neg yb = b$. Hence $QT \vdash J(a)$ and $QT \vdash J(b)$.

Assume $M \vDash J(x)$.

Then $M \vDash \neg xEx$. Suppose, for a reductio, that $M \vDash (x*a)E(x*a)$.

Then $M \vDash \exists y\ yxa = xa$. But then $M \vDash \exists y\ yx = x$ by (QT3), contradicting the hypothesis $M \vDash \neg xEx$.

Therefore $M \vDash \neg(x*a)E(x*a)$.

A completely analogous argument shows that $M \vDash \neg xEx \rightarrow \neg(x*b)E(x*b)$.

This completes the proof that $J(x)$ is a string concept and the proof of (3.4).

(3.5)  For any string concept $I \subseteq I_0$ there is a string concept $J \subseteq I$ such that

$$QT \vdash \forall x \in J \ \neg \exists x_1, x_2 \ (x_1 x x_2 = x).$$

Let $J(x) \equiv I(x) \ \& \ \neg \exists x_1, x_2 \ (x_1 x x_2 = x)$.

We have that $QT \vdash J(a)$ and $QT \vdash J(b)$ by (QT2) and the hypothesis that $I(x)$ is a string concept.

Assume $M \vDash J(x)$.

Suppose, for a reductio, that $M \vDash \exists x_1, x_2 \ (x_1(x*a)x_2 = x*a)$.

By (QT5) and (QT4), $M \vDash x_2 = a \ \lor \ aEx_2$.

If $M \vDash x_2 = a$, then $M \vDash x_1(x*a)a = x*a$, whence $M \vDash x_1 x*a = x$ by (QT3), contradicting the hypothesis $M \vDash J(x)$.

If $M \vDash aEx_2$, then $M \vDash \exists x_3 \ x_3 a = x_2$. Hence $M \vDash x_1(x*a)x_3 a = xa$. But then $M \vDash x_1(x*a)x_3 = x$ by (QT3), that is, $M \vDash x_1 x(ax_3) = x$, contradicting the hypothesis $M \vDash J(x)$.

Therefore $M \vDash J(x) \rightarrow J(x*a)$.

The argument for $M \vDash J(x) \rightarrow J(x*b)$ is completely analogous.

This completes the proof that $J(x)$ is a string concept and the proof of (3.5).

The axioms of QT, specifically (QT3), guarantee only cancellation of atoms. We shall need to be able to cancel entire strings.

(3.6) For any string concept $I \subseteq I_0$ there is a string concept $J \subseteq I$ such that

$$QT^+ \vdash \forall z \in J \; \forall x,y \; (x*z = y*z \to x=y).$$

Let $\quad J(z) \equiv I(z) \;\&\; \forall x,y \; (x*z = y*z \to x=y).$

To show that $QT^+ \vdash J(a)$, note we have that $QT^+ \vdash I(a)$, whereas $QT \vdash \forall x,y \; (x*a = y*a \to x=y)$ by (QT3). Likewise for $QT^+ \vdash J(b)$.

Assume $M \vDash J(z)$. Let $M \vDash x*(z*a) = y*(z*a)$. Then $M \vDash (x*z)*a = (y*z)*a$ by (QT1). But then $M \vDash x*z = y*z$ by (QT3), whence $M \vDash x=y$ by hypothesis $M \vDash J(z)$. Therefore,

$$M \vDash \forall z \; (J(z) \to J(z*a)).$$

Completely analogously, $M \vDash \forall z \; (J(z) \to J(z*b))$.

This completes the proof of (3.6).

This gives us right cancellation for arbitrary strings. Note that right cancellation rules out xEx for strings in $I_0$: that is, we have, for J(x) as in (3.6), that $\qquad QT^+ \vdash \forall x\, (J(x) \rightarrow \neg xEx)$.

For, suppose $M \vDash J(x)\ \&\ xEx$.

Then $M \vDash \exists x_1\, x = x_1 x$, so $M \vDash x_1 x = x_1 x_1 x$, whence $M \vDash x_1 = x_1 x_1$ by right cancellation. But then $M \vDash x_1 B x_1$.

If $M \vDash x_1 = a$, this contradicts (QT2); if $M \vDash \neg\, x_1 = a$, then $M \vDash x_1 R x_1$, which contradicts (2.1) because from $M \vDash x_1 B x\ \&\ J(x)$ we have $M \vDash J(x_1)$. Therefore, $M \vDash \neg xEx$.

For left cancellation, we construct a string concept similarly, but the argument for closure under S proceeds a little differently. We show

(3.7)  For any string concept $I \subseteq I_0$ there is a string concept $J \subseteq I$ such that

$$QT^+ \vdash \forall z \in J\ \forall x,y\ (z*x = z*y \rightarrow x = y).$$

Let    $J(z) \equiv I(z)\ \&\ \forall x,y\ (z*x = z*y \rightarrow x = y)$.

That $QT^+ \vdash J(a)$ is proved in the same way as in (3.6). We indicate the difference in the argument for closure under S.

Suppose  $M \vDash J(z)$. Then  $M \vDash I(z)\ \&\ \forall x,y\ (z*x = z*y \rightarrow x = y)$.

Assume  $M \vDash Sz*x = Sz*y$.

If  $M \vDash z = a$, then $M \vDash Sz = b$, and we have  $M \vDash x = y$  immediately by (QT3). Otherwise,  $M \vDash Sz = z*b$, and so the hypothesis says  $M \vDash (z*b)*x = (z*b)*y$. Then  $M \vDash z*(b*x) = z*(b*y)$. But then from the hypothesis  $M \vDash J(z)$ we get $M \vDash b*x = b*y$, whence  $M \vDash x = y$  by (QT3).  Therefore,

$$M \vDash I(Sz)\ \&\ \forall x,y\ (Sz*x = Sz*y \rightarrow x = y),$$

and so $QT^+ \vdash \forall z\ (J(z) \rightarrow J(Sz))$.  Again, the same argument with b replaced by a shows that J is closed under $S_a$. This completes the proof of (3.7).

We also need to prove in QT that initial segments of arbitrary strings may be totally ordered by the initial-segment-of relation B, making the partial ordering < in which a is the least element tree-like:

(3.8)  For any string concept $I \subseteq I_0$ there is a string concept $J \subseteq I$ such that

$$QT^+ \vdash \forall x \in J \; \forall u,v \; (uBx \; \& \; vBx \to u=v \lor uBv \lor vBu).$$

Let  $J(x) \equiv I(x) \; \& \; \forall y \leq x \; \forall u,v \; (uBy \; \& \; vBy \to u=v \lor uBv \lor vBu).$

For x=a, we have  $QT^+ \vdash J(a)$  because  $QT \vdash y \leq a \leftrightarrow y=a,$  and  $QT \vdash \forall u \; \neg uBa$  by (QT2).

For x=b, we have  $QT \vdash y \leq b \leftrightarrow y=a \lor y=b,$  and so  $QT^+ \vdash J(b)$  because  $QT \vdash \forall u \; \neg uBb$  again by (QT2).

Assume  $M \vDash J(x)$; we want to show that  $M \vDash J(Sx).$

If  $M \vDash x=a$, then  $M \vDash Sx=b$, and we already have the result.

For  $M \vDash \neg x=a,$  we need to show that  $M \vDash J(x*b)$ if $M \vDash J(x).$

So assume  $M \vDash uBy \; \& \; vBy$  where  $M \vDash y \leq (x*b).$  Then  $M \vDash \neg y=a$  by (QT2), so we have  $M \vDash yB(x*b)$  from  $M \vDash y \leq (x*b)$ by the definition of ≤.  Then

M ⊨ uB(x*b) & vB(x*b) from the hypothesis, that is, M ⊨ u*$u_1$=x*b and M ⊨ ∃$u_1$,$v_1$ v*$v_1$=x*b. From (QT4) and (QT5) we have from the former that M ⊨ $u_1$=b ∨ bE$u_1$, so M ⊨ u*b=x*b ∨ ∃$u_2$ u*($u_2$*b)=x*b. Applying (QT3) and (QT1) it follows that

    M ⊨ u=x ∨ uBx, and, similarly, M ⊨ v=x ∨ vBx.

Now, if M ⊨ u=x & v=x, then M ⊨ u=v; if M ⊨ u=x & vBx, then M ⊨ vBu, and if M ⊨ v=x & uBx, then M ⊨ uBv.

Thus the desired claim follows in any of these cases.

The remaining possibility is that M ⊨ uBx & vBx. Then

$$M ⊨ u=v ∨ uBv ∨ vBu$$

follows from the hypothesis M ⊨ J(x).

Therefore, $QT^+$ ⊢ ∀x (J(x) → J(Sx)).

An analogous argument shows that $QT^+$ ⊢ ∀x (J(x) → J($S_a$x)).

This completes the proof of (3.8).

(3.9)  For any string concept $I \subseteq I_0$ there is a string concept $J \subseteq I$ such that

$$QT \vdash \forall y\, \forall x \in J\, (\neg xaBy \vee \neg xbBy).$$

Let $J \equiv I_{3.7}$.

Assume $M \vDash xaBy\ \&\ xbBy$  where $M \vDash J(x)$.

$\Rightarrow M \vDash \exists z_1\, xaz_1 = y\ \&\ \exists z_2\, xbz_2 = y$,

$\Rightarrow M \vDash xaz_1 = xbz_2$,

$\Rightarrow$ by (3.7), $M \vDash az_1 = bz_2$, contradicting (QT4).

This suffices to prove (3.9).

(3.10) For any string concept $I \subseteq I_0$ there is a string concept $J \subseteq I$ such that

$$QT \vdash \forall x \in J \; \forall u,v \; (uEx \; \& \; vEx \rightarrow u=v \; \lor \; uEv \; \lor \; vEu).$$

Let $J(x) \equiv I(x) \; \& \; \forall u,v \; (uEx \; \& \; vEx \rightarrow u=v \; \lor \; uEv \; \lor \; vEu)$.

Let $x=a$.

$\Rightarrow$ by (QT2), $QT \vdash \neg uEx \; \& \; \neg vEx$,

$\Rightarrow QT \vdash J(a)$, trivially.

Likewise for $x=b$.

Assume $M \vDash uE(xa) \; \& \; vE(xa)$ where $M \vDash J(x)$ and $M \vDash QT$.

$\Rightarrow M \vDash \exists u_1, v_1 \; (u_1 u = xa \; \& \; v_1 v = xa)$,

$\Rightarrow$ by (QT4) and (QT5), $M \vDash (u=a \lor aEu) \; \& \; (v=a \lor aEv)$.

We distinguish four cases:

  (1)  $M \vDash u=a \; \& \; v=a$.

$\Rightarrow M \vDash u=v$, so the claim holds.

  (2)  $M \vDash u=a \; \& \; aEv$.

$\Rightarrow M \vDash uEv$.

  (3)  $M \vDash aEu \; \& \; v=a$.

$\Rightarrow M \vDash vEu$.

  (4)  $M \vDash aEu \; \& \; aEv$.

$\Rightarrow M \vDash \exists u_2, v_2 \; (u_2 a = u \; \& \; v_2 a = v)$,

$\Rightarrow M \vDash u_1 u_2 a = xa \; \& \; v_1 v_2 a = xa$,

$\Rightarrow$ by (QT3), $M \vDash u_1 u_2 = x \; \& \; v_1 v_2 = x$,

⇒ M ⊨ u₂Ex & v₂Ex,

⇒ from hypothesis M ⊨ J(x), M ⊨ u₂=v₂ v u₂Ev₂ v v₂Eu₂.

There are three subcases:

   (4a) M ⊨ u₂=v₂.

⇒ M ⊨ u₂a=v₂a,

⇒ by (QT3), M ⊨ u=v.

   (4b) M ⊨ u₂Ev₂.

⇒ M ⊨ ∃v₃ v₃u₂=v₂,

⇒ M ⊨ v₃(u₂a)=v₂a,

⇒ M ⊨ v₃u=v,

⇒ M ⊨ uEv.

   (4c) M ⊨ v₂Eu₂.

⇒ M ⊨ ∃u₃ u₃v₂=u₂,

⇒ M ⊨ u₃(v₂a)=u₂a,

⇒ M ⊨ u₃v=u,

⇒ M ⊨ vEu.

This completes the proof that M ⊨ J(x*a).

Exactly analogously for x*b.

It follows that J(x) is a string concept.

This completes the proof of (3.10).

Let

$$x \subseteq_p y \equiv x=y \lor xBy \lor xEy \lor \exists u,v\ (u*x)*v=y$$

meaning "x is a substring of y". We then write $\exists y \subseteq_p x$ (... for $\exists y(y \subseteq_p x\ \&\ ...$ and $\forall y \subseteq_p x$ (... for $\forall y(y \subseteq_p x \to ....$

We observe that

$$QT \vdash x \subseteq_p y\ \&\ y \subseteq_p z \to x \subseteq_p z.$$

We argue by cases, taking into account (QT1).

If $M \vDash x=y$ or $M \vDash y=z$, there is nothing to prove.

Suppose (1) $M \vDash xBy$.

(1a) If $M \vDash yBz$, then $M \vDash \exists x_1\ xx_1=y$ and $M \vDash \exists y_1\ yy_1=z$, so $M \vDash xx_1y_1=z$, that is

$M \vDash xBz$.

(1b) If $M \vDash yEz$, then $M \vDash \exists y_1\ y_1y=z$ and $M \vDash \exists x_1\ xx_1=y$, so $M \vDash y_1xx_1=z$.

Thus $M \vDash x \subseteq_p z$.

(1c) If $M \vDash \exists y_1,y_2\ y_1yy_2=z$, then from $M \vDash xx_1=y$ we have $M \vDash y_1xx_1y_2=z$.

So $M \vDash x \subseteq_p z$.

Suppose (2) $M \vDash xEy$.

(2a) If $M \vDash yBz$, then $M \vDash \exists x_1\, x_1x=y$ and $M \vDash \exists y_1\, yy_1=z$, so $M \vDash x_1xy_1=z$, thus $M \vDash x\subseteq_p z$.

(2b) If $M \vDash yEz$, then $M \vDash \exists x_1\, x_1x=y$ and $M \vDash \exists y_1\, y_1y=z$, so $M \vDash y_1x_1x=z$, thus $M \vDash xEz$.

(2c) If $M \vDash \exists y_1, y_2\, y_1yy_2=z$, then from $M \vDash x_1x=y$ we have $M \vDash y_1x_1xy_2=z$, so $M \vDash x\subseteq_p z$.

Finally, suppose (3) $M \vDash \exists x_1, x_2\, x_1xx_2=y$.

(3a) If $M \vDash yEz$, then $M \vDash \exists y_1\, y_1y=z$, so $M \vDash y_1x_1xx_2=z$, thus $M \vDash x\subseteq_p z$.

(3b) If $M \vDash yBz$, then $M \vDash \exists y_1\, yy_1=z$, so $M \vDash x_1xx_2y_1=z$, thus $M \vDash x\subseteq_p z$.

And (3c) if $M \vDash \exists y_1, y_2\, y_1yy_2=z$, then $M \vDash y_1x_1xx_2y_2=z$, so again $M \vDash x\subseteq_p z$.

(3.11) For any string concept $I \subseteq I_0$ there is a string concept $J \subseteq I$ such that

$$QT \vdash \forall x \in J \, \forall y \, (x \subseteq_p y \,\&\, y \subseteq_p x \rightarrow x = y).$$

Let $J(x) \equiv I_{3.4}(x) \,\&\, I_{3.5}(x)$.

Assume $M \vDash J(x) \,\&\, x \subseteq_p y \,\&\, y \subseteq_p x$.

Suppose $M \vDash x \neq y$.

From $M \vDash x \subseteq_p y$ we distinguish three cases:

(1) $M \vDash xBy$.

From $M \vDash y \subseteq_p x$, we distinguish three subcases:

(1a) $M \vDash yBx$.

$\Rightarrow M \vDash xBy \,\&\, yBx$,

$\Rightarrow$ by (QT2), $M \vDash x \neq a \,\&\, y \neq a$,

$\Rightarrow M \vDash xRy \,\&\, yRx$,

$\Rightarrow$ by transitivity of R, $M \vDash xRx$.

But this contradicts $M \vDash x \in J \subseteq I_0$.

(1b) $M \vDash yEx$.

$\Rightarrow M \vDash \exists x_1 \, x_1 y = x \,\&\, \exists y_1 \, xy_1 = y$,

$\Rightarrow M \vDash x_1 xy_1 = x$, which contradicts $M \vDash x \in J \subseteq I_{3.5}$.

(1c) $M \vDash \exists x_1, x_2 \, x_1 y x_2 = x$.

Also, $M \vDash \exists y_1 \, xy_1 = y$.

$\Rightarrow M \vDash x_1(xy_1)x_2 = x$,

$\Rightarrow M \vDash x_1 x(y_1 x_2) = x$, again contradicting $M \vDash x \in J \subseteq I_{3.5}$.

(2) $M \vDash xEy$.

Again, from $M \vDash y\subseteq_p x$, we distinguish three subcases:

(2a) $M \vDash yBx$.

Completely analogous to (1a) with x, y exchanging places.

(2b) $M \vDash yEx$.

$\Rightarrow M \vDash \exists x_1\, x_1y=x\ \&\ \exists y_1\, y_1x=y$,

$\Rightarrow M \vDash x_1y_1x=x$,

$\Rightarrow M \vDash xEx$, which contradicts $M \vDash x \in J\subseteq I_{3.4}$.

(2c) $M \vDash \exists x_1,x_2\, x_1yx_2=x$.

Also, $M \vDash \exists y_1\, y_1x=y$.

$\Rightarrow M \vDash x_1(y_1x)x_2=x$,

$\Rightarrow M \vDash (x_1y_1)xx_2=x$, contradicting $M \vDash x\in J\subseteq I_{3.5}$.

(3) $M \vDash \exists y_1,y_2\, y_1xy_2=y$.

The three subcases:

(3a) $M \vDash yBx$.

Completely analogous to (1c).

(3b) $M \vDash yEx$.

Completely analogous to (2c).

(3c) $M \vDash \exists x_1,x_2\, x_1yx_2=x$.

$\Rightarrow M \vDash x_1(y_1xy_2)x_2=x$,

$\Rightarrow M \vDash (x_1y_1)x(y_2x_2)=x$, contradicting $M \vDash x\in J\subseteq I_{3.5}$.

Since Cases (1)-(3) were all shown to lead to a contradiction, the only remaining possibility from $M \vDash x \subseteq_p y$ is $M \vDash x = y$.

This completes the proof of (3.11).

Note the following generalization of (3.4)-(3.5):

(3.12) For any string concept $I \subseteq I_0$ there is a string concept $J \subseteq I$ such that

$$QT \vdash \forall x \in J \, \forall y \, (\neg xy \subseteq_p x \,\&\, \neg yx \subseteq_p x).$$

Let $J(x) \equiv I_{3.4}(x) \,\&\, I_{3.5}(x)$, as in (3.11).

Suppose $M \vDash J(x)$.

Assume $M \vDash xy \subseteq_p x$.

Have $M \vDash xB(xy)$.

$\Rightarrow M \vDash x \subseteq_p xy$,

$\Rightarrow$ by (3.11), $M \vDash x = xy$,

$\Rightarrow M \vDash xBx$,

$\Rightarrow M \vDash xRx$, which contradicts $M \vDash x \in J \subseteq I_0$.

Therefore, $M \vDash \neg xy \subseteq_p x$.

Assume $M \vDash yx \subseteq_p x$.

Have $M \vDash x \subseteq_p yx$, $\Rightarrow$ by (3.11), $M \vDash x = yx$,

$\Rightarrow M \vDash xEx$, contradicting $\quad M \vDash x \in J \subseteq I_{3.4}$.

Therefore also $M \vDash \neg yx \subseteq_p x$.

This completes the proof of (3.12).

To guarantee downward closure under substrings, we need the following:

(3.13)   For any string concept I stronger than $I_0$ there is a string concept $J\subseteq I$ such that $QT^+ \vdash \forall x \in J \,\forall y\, (y\subseteq_p x \rightarrow J(y))$.

Let $I^{\subseteq p}(x) \equiv I(x) \,\&\, \forall z\leq x \,\forall y\, (y\subseteq_p z \rightarrow I(y))$, and let $J \equiv I^{\subseteq p}$.

For $QT^+ \vdash J(a)$, note that $QT^+ \vdash I(a)$ since I is a string concept, and $QT \vdash y\subseteq_p a \leftrightarrow y=a$ from (QT2). Hence $QT^+ \vdash \forall y\, (y\subseteq_p a \rightarrow I(y))$. But this suffices for $QT^+ \vdash J(a)$ because $QT \vdash z\leq a \leftrightarrow z=a$ by (1.5).

Likewise $QT^+ \vdash J(b)$, where we need only note that

$$QT \vdash z\leq b \leftrightarrow z=a \lor z=b$$

and appeal to $QT^+ \vdash I(b)$.

Suppose $M \vDash J(x)$.

If $M \vDash x=a$, we have $M \vDash Sx=b$, and so $M \vDash J(Sx)$ by what we just proved.

Otherwise $M \vDash Sx=x*b$. Suppose $M \vDash z\leq Sx$, and let $M \vDash y\subseteq_p z$. By (1.14), it is sufficient to consider the two cases $M \vDash z\leq x$ and $M \vDash z=Sx$.

If $M \vDash z\leq x$, then $M \vDash I(y)$ follows from the hypothesis $M \vDash J(x)$.

So let $M \models z=Sx=x*b$. Then, by definition,

$$M \models y \subseteq_p x*b \leftrightarrow y=x*b \vee yB(x*b) \vee yE(x*b) \vee \exists x_1, x_2\ x_1yx_2=x*b.$$

There are four cases to consider.

Case (i): $M \models y=x*b$.

We have $M \models I(x)$ from $M \models J(x)$, and $M \models I(x*b)$ because I is a string concept by hypothesis. But then $M \models I(y)$ as required.

Case (ii): $M \models yB(x*b)$.

Then $M \models \exists x_1\ yx_1=x*b$, and $M \models b=x_1 \vee bEx_1$ by (QT4) and (QT5).

If $M \models b=x_1$, then $M \models y*b=x*b$, and we have $M \models y=x$ by (QT3).

Then $M \models I(y)$ follows from the hypothesis $M \models J(x)$.

If $M \models \exists x_1'\ x*b=yx_1=y(x_1'b)$, then $M \models x=yx_1'$ by (QT3) and so $M \models y \subseteq_p x$.

But then $M \models I(y)$ from the hypothesis $M \models J(x)$.

Case (iii): $M \models yE(x*b)$.

Then $M \models \exists x_1\ x_1y=x*b$. By (QT4) and (QT5), we have $M \models y=b \vee bEy$.

If $M \models y=b$, we have $M \models I(y)$ from the hypothesis that I is a string concept.

If $M \vDash y_1b=y$, then $M \vDash x_1(y_1b)=x_1y=x*b$, whence $M \vDash x_1y_1=x$ by (QT3).

So $M \vDash y_1 \subseteq_p x$, and $M \vDash I(y_1)$ from the hypothesis $M \vDash J(x)$. But then $M \vDash I(y_1b)$ from the hypothesis that I is a string concept, and thus $M \vDash I(y)$.

Case (iv): $M \vDash x_1yx_2=x*b$.

Then $M \vDash b=x_2 \lor bEx_2$ by (QT4) and (QT5). Then $M \vDash x_1yx_2=x_1yb$ or $M \vDash \exists x_2'\ x_1yx_2=x_1y(x_2'b)$, whence $M \vDash xb=x_1yb$ or $M \vDash xb= x_1y(x_2'b)$.

But then $M \vDash x=x_1y$ or $M \vDash x= x_1yx_2'$ by (QT3).

In either case $M \vDash y \subseteq_p x$ and we have $M \vDash I(y)$ from the hypothesis $M \vDash J(x)$.

From Cases (i)-(iv) we thus have $M \vDash \forall y\ (y \subseteq_p Sx \to I(y))$, which is what was needed to show that $M \vDash \forall z \leq Sx\ \forall y\ (y \subseteq_p z \to I(y))$.

So we proved that $M \vDash J(Sx)$ if $M \vDash J(x)$ as required.

That $J(x)$ is closed under $S_a$ is established in a similar fashion, with (1.14) replaced by (1.14a). This completes the proof that $J(x)$ is a string concept with the required properties and the proof of (3.13).

## §4. Tallies

Let  Digit(y) ≡ x=a v y=a.

We now set

   $Tally_a(x) \equiv \forall y \subseteq_p x\ (Digit(y) \to y=a)$

   $Tally_b(x) \equiv \forall y \subseteq_p x\ (Digit(y) \to y=b)$

for the strings that consist exclusively of a's and b's, respectively.

We want to show that tallies are closed under concatenation. The first step towards that is:

(4.1)  $QT^+ \vdash Tally_b(y) \to Tally_b(Sy)$.

Assume that $M \vDash Tally_b(y)$. Then  $M \vDash \neg y=a$  and  so $M \vDash Sy=y*b$. Assume that $M \vDash \neg Tally_b(y*b)$. Then $M \vDash \exists z \subseteq_p (y*b)(Digit(z)\ \&\ \neg z=b)$. By (QT5),

$M \vDash a \subseteq_p (y*b)$. But  $M \vDash \neg a=(y*b)$  by (QT2), and  $M \vDash \neg aB(y*b)$  by (QT4) because  $M \vDash Tally_b(y)$ , and  $M \vDash \neg aE(y*b)$,  again by (QT4). So we must have $M \vDash \exists y_1, y_2\ y_1 a y_2 = y*b$. Then  $M \vDash y_2=b\ v\ bEy_2$ by (QT4) and (QT5). Either way we derive  $M \vDash a \subseteq_p y$, contradicting the hypothesis $M \vDash Tally_b(y)$.

Therefore   $M \vDash Tally_b(y*b)$.

This completes the proof of (4.1).

(4.2)   $QT^+ \vdash Tally_b(y) \leftrightarrow y=b \lor \exists y_1 (Tally_b(y_1) \& y=Sy_1)$.

Assume $M \vDash Tally_b(y)$. Then by, (QT4) and (QT5),

$$M \vDash y=b \lor bEy.$$

If $M \vDash bEy$, then $M \vDash \exists y_1\, y=y_1*b$. But $M \vDash Tally_b(y)$ implies $M \vDash Tally_b(y_1)$. Then $M \vDash y_1 \neq a$, so $M \vDash y_1*b = Sy_1$. Therefore

$$M \vDash Tally_b(y) \rightarrow y=b \lor \exists y_1 (Tally_b(y_1) \& y=Sy_1).$$

Conversely, suppose $M \vDash y=b$. Let $M \vDash x \subseteq_p y \& Digit(x)$. Then, by (QT2), $M \vDash x \subseteq_p y$ implies $M \vDash x=y$, so $M \vDash \forall x \subseteq_p y\, (Digit(x) \rightarrow x=b)$, that is, $M \vDash y=b \rightarrow Tally_b(y)$.

On the other hand,

$$M \vDash Tally_b(y_1) \& y=Sy_1 \rightarrow Tally_b(y)$$

follows by (4.1).

This completes the proof of (4.2).

(4.3)   QT ⊢ Tally$_a$(w)  →  w=a v w=aa v w=aaa v (aaaBw & aaaEw).

Suppose  M ⊨ Tally$_a$(w) & ¬w=a & ¬w=aa.

By (QT4) and (QT5) it follows that  M ⊨ aBw & aEw,  i.e. that

$$M \vDash \exists z_1\ az_1=w\ \&\ \exists z_2\ z_2a=w.$$

Again by (QT4) and (QT5),  M ⊨ z$_1$=a v aEz$_1$  because M ⊨ Tally$_a$(w).  But M ⊨ z$_1$=a  contradicts the hypothesis that  M ⊨ ¬w=aa.  Therefore, M ⊨ aEz$_1$, so  M ⊨ ∃z$_1$' a(z$_1$'a)=w.  Again, by (QT4) and (QT5),

$$M \vDash z_1'=a \ v\ (aBz_1'\ \&\ aEz_1')$$

because  M ⊨ Tally$_a$(w).  In either case M ⊨ aaBw.

But then, again by (QT4) and (QT5),

$$M \vDash aaa=w\ v\ \exists z_1''\ aa(z_1''a)a=w,$$

so     M ⊨ aaa=w v (aaaBw & aaaEw)  because  M ⊨ z$_1$''=a v (aBz$_1$'' & aEz$_1$'').

This completes the proof of (4.3).

(4.4)   QT ⊢ ∀x (Tally$_b$(x) → ∀v⊆$_p$x ¬Tally$_a$(v)).

Assume  M ⊨ Tally$_b$(x) & v⊆$_p$x.

Suppose for a reductio that M ⊨ Tally$_a$(v).

⇒  M ⊨ v=a  v ∃v$_1$ (Tally$_a$(v$_1$) & av$_1$=v),

⇒  M ⊨ a⊆$_p$v⊆$_p$x,  which contradicts  M ⊨ Tally$_b$(x).

This proves (4.4).

(4.5) For any string concept I⊆I₀ there is a string concept J⊆I such that

$$QT^+ \vdash \forall z \in J\, \forall y\, (Tally_b(y)\, \&\, Tally_b(z) \rightarrow Tally_b(y*z)).$$

Let $J(z) \equiv I(z)\, \&\, \forall y\, (Tally_b(y)\, \&\, Tally_b(z) \rightarrow Tally_b(y*z))$.

We have that $QT^+ \vdash J(a)$ because $QT^+ \vdash I(a)$ and $QT \vdash \neg Tally_b(y*a)$.

For z=b, we have $QT^+ \vdash I(b)$. From $M \vDash Tally_b(y)$ we have $M \vDash Tally_b(y*b)$ by (4.1). Hence $QT^+ \vdash J(b)$.

Suppose $M \vDash J(z)$. Then $M \vDash I(z)$, whence $M \vDash I(Sz)$.

Suppose $M \vDash Tally_b(y)\, \&\, Tally_b(Sz)$.

We want $M \vDash Tally_b(y*Sz)$.

If z=a, then $M \vDash Sz=b$, and $M \vDash Tally_b(y*Sz)$ was proved in $QT^+ \vdash J(b)$.

If z≠a, we have $M \vDash y*Sz = y*(z*b) = (y*z)*b = S(y*z)$ since $M \vDash y*z \neq a$ by (QT2).

Now, $QT \vdash Tally_b(Sz) \rightarrow Tally_b(z)$.

So from the hypothesis $M \vDash J(z)$, we have $M \vDash Tally_b(y*z)$.

Assume, for a reductio, that $M \vDash \neg Tally_b(S(y*z))$, i.e. that $M \vDash a \subseteq_p S(y*z)$.

Then $M \vDash a \subseteq_p (y*z)*b$. Then, just as in the proof of (4.1), we derive

$M \vDash a \subseteq_p (y*z)$, so $M \vDash \neg \text{Tally}_b(y*z)$, a contradiction.

Therefore, $M \vDash J(Sz)$, and we have $QT^+ \vdash J(z) \to J(Sz)$. That J is closed under $S_a$ is immediate.

This completes the proof of (4.5).

(4.6)  For any string concept I stronger than $I_0$ there is a string concept $J \subseteq I$

such that  $QT^+ \vdash \forall z \in J \, \forall x \, (Tally_b(x) \, \& \, Tally_b(z) \to x \leq z \, \lor \, z \leq x)$.

Let $J(z) \equiv I(z) \, \& \, \forall x \, \forall y \leq z \, (Tally_b(x) \, \& \, Tally_b(y) \to x \leq y \, \lor \, y \leq x)$.

We have $QT^+ \vdash I(a)$. That the right hand conjunct of $J(a)$ is provable in QT follows from (1.5) and the definition of $Tally_b$. So $QT^+ \vdash J(a)$.

Suppose $M \vDash J(z)$. We want $M \vDash J(Sz)$. We do have $M \vDash I(z)$ from $M \vDash J(z)$, so $M \vDash I(Sz)$ since I is a string concept.

Assume $M \vDash y \leq Sz$ along with $M \vDash Tally_b(x) \, \& \, Tally_b(y)$.

If $M \vDash y \leq z$, the claim follows from the hypothesis $M \vDash J(z)$.

By (1.14), it is sufficient to establish the same under the condition $M \vDash y = Sz$. We may assume that $M \vDash \neg z = a$, for otherwise $M \vDash Sz = b$, and $M \vDash y = a$ by (1.3), in which case the desired conditional holds by the definition of $Tally_b$ since $M \vDash \neg Tally_b(y)$. So $M \vDash y = z*b \, \& \, Tally_b(z*b)$.

We have that $QT \vdash Tally_b(z*b) \to Tally_b(z)$, so the hypothesis $M \vDash J(z)$ entails $M \vDash x \leq z \, \lor \, z \leq x$. But we need $M \vDash x \leq y \, \lor \, y \leq x$.

So assume $M \vDash x \leq z$, whence $M \vDash x \leq z*b$ by (1.7) and (1.8), and thus $M \vDash x \leq y$. Assume, on the other hand, $M \vDash z \leq x$.

If $M \vDash z=x$, then $M \vDash y=x*b$, so $M \vDash x \leq y$ by (1.8).

If $M \vDash z<x$, then $M \vDash zBx$, and $M \vDash \exists u\ z*u=x$.

But we have that $M \vDash \text{Tally}_b(u)$ because $QT \vdash \text{Tally}_b(x)\ \&\ u \subseteq_p x \rightarrow \text{Tally}_b(u)$.

Then $M \vDash \text{Tally}_b(u)$ implies $M \vDash u=b$, or $M \vDash b*w=u$, by (QT5).

In the former case, $M \vDash z*b=z*u=x$, and we then have $M \vDash y=x$.

In the latter case, $M \vDash z*(b*w)=x$, so $M \vDash (z*b)*w=x$ by (QT1), so

$M \vDash (z*b)Bx$, and thus $M \vDash (z*b)<x$, that is $M \vDash y<x$.

We have proved that $M \vDash y=z*b\ \&\ \text{Tally}_b(z*b) \rightarrow x \leq y \vee y \leq x$ under the hypothesis $M \vDash \text{Tally}_b(x)\ \&\ \text{Tally}_b(y)$. But this gives us

   $M \vDash \text{Tally}_b(x)\ \&\ \text{Tally}_b(Sz) \rightarrow x \leq Sz \vee Sz \leq x$,

which is what we needed to complete our argument that

   $M \vDash \forall y \leq Sz\ (\text{Tally}_b(x)\ \&\ \text{Tally}_b(y) \rightarrow x \leq y \vee y \leq x)$.

Along with $M \vDash I(Sz)$, this establishes that $M \vDash J(Sz)$ if $M \vDash J(z)$. The proof that $J$ is closed under $S_a$ is much easier, appealing to (1.14[a]) and the fact that $T \vdash \neg \text{Tally}_b(z*a)$.

This completes the proof of (4.6).

According to (4.6), the b-tallies in J are linearly ordered by ≤, provably in $QT^+$. Claims analogous to (4.3)-(4.6) with $Tally_a$ in place of $Tally_b$ are similarly proved. We also have:

(4.7)   $QT^+ \vdash \forall v,u\ (Tally_b(v)\ \&\ u<v \rightarrow Su \leq v)$.

Assume $M \vDash Tally_b(v)\ \&\ u<v$. Then $M \vDash I_0(v)\ \&\ uRv$.

If $M \vDash u=a$, then $M \vDash Su=b$, and $M \vDash (Su)Bv$ by (QT4) and (QT5), so

$M \vDash Su \leq v$. So we may assume $M \vDash \neg u=a$, in which case the hypothesis implies $M \vDash Tally_b(u)$. From $M \vDash u<v$, we have $M \vDash \exists w\ uw=v$. By (QT5),

$$M \vDash w=a \vee w=b \vee aBw \vee bBw.$$

From $M \vDash Tally_b(v)$ we have $M \vDash \neg w=a\ \&\ \neg aBw$. Hence

$$M \vDash ub=v \vee \exists w'\ u(bw')=v,$$

and so $M \vDash Su=ub \leq v$, as required.

This completes the proof of 4.7.

(4.8) For any string concept I⊆I₀ there is a string concept J⊆I such that

$$QT^+ \vdash \forall u \in J \ (Tally_b(u) \to u*b = b*u).$$

Let $J(x) \equiv I(x) \ \& \ \forall u \leq x \ (Tally_b(u) \to u*b = b*u)$.

Then $QT^+ \vdash J(a)$ because $QT \vdash u \leq a \to u = a$ by (1.3) and $QT \vdash \neg Tally_b(a)$. Suppose that $M \vDash J(x)$, and assume $M \vDash u \leq Sx \ \& \ Tally_b(u)$. If $M \vDash u \leq x$, the claim $M \vDash u*b = b*u$ holds by the hypothesis $M \vDash J(x)$. By (1.14), we need only consider $M \vDash u = Sx$.

If $M \vDash x = a$, then $M \vDash Sx = b$, and we have $M \vDash b*b = b*b$ for free.

Otherwise, $M \vDash \neg x = a$, and $M \vDash Sx = x*b$. Then $M \vDash Sx*b = (x*b)*b$. But from the hypothesis $M \vDash J(x)$ we have in particular that $M \vDash x*b = b*x$. Therefore

$$M \vDash (x*b)*b = (b*x)*b = b*(x*b) = b*Sx,$$

appealing to (QT1). So $M \vDash Sx*b = b*Sx$, as required. That $J(x)$ is closed under $S_a$ is seen as in (4.6).

This completes the proof of (4.8).

(4.9)    $QT^+ \vdash Tally_b(y) \rightarrow (x<y \leftrightarrow Sx<Sy)$.

Assume $M \vDash Tally_b(y)$ along with $M \vDash x<y$.

$\Rightarrow M \vDash I_0(y)$,

$\Rightarrow M \vDash I_0(Sy)$.

From (4.7), $M \vDash Sx \leq y$, and from (2.4), $M \vDash y<Sy$.

$\Rightarrow$ from (1.7), $M \vDash Sx<Sy$.

Conversely, assume $M \vDash Sx<Sy$.

$\Rightarrow M \vDash I_0(Sy) \& (Sx)R(Sy)$,

$\Rightarrow$ by (1.13), $M \vDash (Sx)Ry \lor Sx=Sy$.

Suppose, for a reductio, that $M \vDash Sx=Sy$.

$\Rightarrow M \vDash (Sx)R(Sx)$, contradicting $M \vDash I_0(Sy) \& (Sx)R(Sy)$.

Suppose that $M \vDash (Sx)Ry$.

$\Rightarrow$ by (1.8), $M \vDash xR(Sx)$,

$\Rightarrow$ by (1.7), $M \vDash xRy$,

$\Rightarrow$ by (1.8), $M \vDash yR(Sy)$,

$\Rightarrow$ from $M \vDash I_0(Sy)$, since $I_0$ is a downward closed with respect to R, $M \vDash I_0(y)$,

$\Rightarrow$ from $M \vDash xRy$, by definition of $<$, $M \vDash x<y$, as required.

This completes the proof of (4.9).

What we proved a moment ago, in (4.8), is the stepping stone to:

(4.10)   For any string concept $I \subseteq I_{4.8}$ there is a string concept $J \subseteq I$ such that

$$QT^+ \vdash \forall u,v \in J \ (Tally_b(u) \ \& \ Tally_b(v) \ \rightarrow \ u*v = v*u).$$

Let  $J(x) \equiv I(x) \ \& \ \forall u \in I \ \forall v \leq x \ (Tally_b(u) \ \& \ Tally_b(v) \ \rightarrow \ u*v = v*u)$.

Note we have  $QT^+ \vdash J(a)$  because $QT \vdash v \leq a \rightarrow v = a$  by (1.3)  and

$QT \vdash \neg Tally_b(a)$. We have that  $QT^+ \vdash I(b)$.

Suppose that $M \vDash v \leq b$ and $M \vDash I(u)$. Then  $M \vDash v = a$  or  $M \vDash v = b$, and we have

$QT^+ \vdash J(b)$ because by hypothesis $I \subseteq I_{4.8}$.

To show closure under S, suppose  $M \vDash J(x)$  and assume $M \vDash v \leq Sx$ along with

$M \vDash Tally_b(u) \ \& \ Tally_b(v)$  where $M \vDash I(u)$.

If  $M \vDash v \leq x$, the result follows from the hypothesis $M \vDash J(x)$. By (1.14), it is sufficient to show the same under the condition $M \vDash v = Sx$. We may assume that $M \vDash x \neq a$  for otherwise  $M \vDash Sx = b$  and we already have that  $M \vDash J(b)$.

Then  $M \vDash v = x*b$ and we have  $M \vDash u*v = u*(x*b) = (u*x)*b$  by (QT1). But we

have $M \vDash u*x=x*u$ from $M \vDash J(x)$. So $M \vDash (u*x)*b=(x*u)*b=x*(u*b)$ by (QT1). But $M \vDash u*b=b*u$ since $QT^+ \vdash J(b)$. Therefore

$$M \vDash x*(u*b)=x*(b*u)=(x*b)*u=v*u,$$

again appealing to (QT1). But then $M \vDash u*v=v*u$, as required. Again, that $J(x)$ is closed under $S_a$ is established as in (4.6).

This completes the proof of (4.10).

(4.11)  $QT^+ \vdash Tally_a(x)$ & $Tally_a(y)$ & $x \subseteq_p y \rightarrow xa \subseteq_p ya$.

Assume  $M \vDash Tally_a(x)$ & $Tally_a(y)$ & $x \subseteq_p y$.

We distinguish the cases based on the definition of $\subseteq_p$.

(i)  $M \vDash x=y$.

$\Rightarrow M \vDash xa=ya$,

$\Rightarrow M \vDash xa \subseteq_p ya$.

(ii)  $M \vDash xBy$.

$\Rightarrow M \vDash \exists x_1\ xx_1=y$,

$\Rightarrow$ from  $M \vDash Tally_a(y)$,  $M \vDash Tally_a(x_1)$,

$\Rightarrow$ by (QT5) and (QT4),  $M \vDash x_1=a \lor aBx_1$,

$\Rightarrow M \vDash xa=y \lor \exists x_2\ x(ax_2)=y$,

$\Rightarrow M \vDash xaBya$,

$\Rightarrow M \vDash xa \subseteq_p ya$.

(iii)  $M \vDash xEy$.

$\Rightarrow M \vDash \exists x_1\ x_1x=y$,

$\Rightarrow M \vDash x_1xa=ya$,

$\Rightarrow M \vDash xaEya$,

$\Rightarrow M \vDash xa \subseteq_p ya$.

(iv)  $M \vDash \exists x_1,x_2\ x_1xx_2=y$.

$\Rightarrow$ from  $M \vDash Tally_a(y)$,  $M \vDash Tally_a(x_2)$,

$\Rightarrow$ by (QT5) and (QT4),  $M \vDash x_2=a \lor aBx_2$,

⇒ M ⊨ $x_1xa=y$ ∨ $\exists x_3\ x_1x(ax_3)=y$,

⇒ M ⊨ $x_1(xa)a=ya$ ∨ $x_1(xa)x_3a=ya$,

⇒ M ⊨ $xa\subseteq_p ya$.

This completes the proof of (4.11).

Let
$$MaxT_a(z,x) \equiv Tally_a(z) \mathbin{\&} \forall v\subseteq_p x\, (Tally_a(v) \to v\subseteq_p z)),$$
meaning "z is an a-tally as long as any a-tally in x" if x has any a digits. Similarly for $MaxT_b(z,x)$.

(4.12) For any string concept $I\subseteq I_0$ there is a string concept $I_{MaxT}\subseteq I$ such that
$$QT^+ \vdash \forall x\in I_{MaxT}\, \exists z(MaxT_a(z,x) \mathbin{\&} (\neg Tally_b(x) \to z\subseteq_p x)).$$

Let $\quad I_{MaxT}(x) \equiv I(x) \mathbin{\&} \exists z(MaxT_a(z,x) \mathbin{\&} (\neg Tally_b(x) \to z\subseteq_p x))$.

For x=a, we have $QT \vdash I(a)$. Since $QT \vdash \forall v\, (v\subseteq_p a \leftrightarrow v=a)$, and $QT \vdash Tally_a(a)$, we have $QT \vdash I_{MaxT}(a)$.

On the other hand, note that, from (4.4),
$$QT^+ \vdash Tally_b(x) \to MaxT_a(a,x).$$

Then, for x=b, since $QT \vdash I(b)$ and $QT \vdash Tally_b(b)$, we have that $QT \vdash I_{MaxT}(b)$ follows immediately.

Let $M \vDash I_{MaxT}(x)$.

Consider x*b. We may assume $M \vDash \neg Tally_b(x*b)$.

Then $M \vDash \neg Tally_b(x)$ by (4.1).

Hence from hypothesis $M \vDash I_{MaxT}(x)$ we have $M \vDash \exists z_x\subseteq_p x\, MaxT_a(z_x,x)$.

Suppose $M \vDash v\subseteq_p x*b \mathbin{\&} Tally_a(v)$.

Then, by the definition of $\subseteq_p$, we have that
$$M \vDash v=x*b \;\vee\; vB(x*b) \;\vee\; vE(x*b) \;\vee\; \exists u,w\, uvw=x*b.$$

We distinguish the four cases:

(a) $M \vDash v=x*b$.

Have $QT \vdash Tally_a(v) \rightarrow \forall y \subseteq_p v\ (Digit(y) \rightarrow y=a)$ and

$QT \vdash v=x*b \rightarrow b \subseteq_p v\ \&\ Digit(b)$. Hence, from the hypothesis (a), $M \vDash b=a$.

$\Rightarrow M \vDash \neg v=x*b$ since $QT \vdash \neg b=a$.

(b) $M \vDash vE(x*b)$.

$\Rightarrow M \vDash \exists x_1\ x_1 v = x*b$,

$\Rightarrow$ by axioms of QT, $M \vDash v=b \lor bEv$,

$\Rightarrow M \vDash b \subseteq_p v\ \&\ Digit(b)$,

$\Rightarrow$ from hypothesis, $M \vDash Tally_a(v) \rightarrow b \subseteq_p v\ \&\ Digit(b)$,

$\Rightarrow$ as in (a), $M \vDash b=a$,

$\Rightarrow M \vDash \neg vE(x*b)$.

So both (a) and (b) are ruled out.

(c) $M \vDash vB(x*b)$.

$\Rightarrow M \vDash \exists x_1\ v x_1 = x*b$,

$\Rightarrow$ by (QT5) and (QT4), $M \vDash x_1=b \lor bEx_1$,

$\Rightarrow M \vDash vb=x*b \lor \exists x_2\ (vx_2)b = x*b$,

$\Rightarrow$ by (QT3), $M \vDash v=x \lor vx_2=x$,

$\Rightarrow M \vDash v \subseteq_p x$.

Letting $z=z_x$, we have from $M \vDash MaxT_a(z_x,x)$ that

$$M \vDash Tally_a(z)\ \&\ v \subseteq_p x*b\ \&\ (Tally_a(v) \rightarrow v \subseteq_p z),$$

as required.

(d) $M \vDash \exists u,w \; uvw=x*b$.

$\Rightarrow$ by (QT5) and (QT4), $M \vDash w=b \lor bEw$,

$\Rightarrow M \vDash uvb=x*b \lor \exists w_1 \; (uv)w_1 b=x*b$,

$\Rightarrow M \vDash uv=x \lor uvw_1=x$,

$\Rightarrow M \vDash v \subseteq_p x$.

Now the argument proceeds exactly as in (c).

So we have shown that $M \vDash \exists z \subseteq_p x*b \; (Tally_a(z) \;\&\; \forall v \subseteq_p x*b \; (Tally_a(v) \to v \subseteq_p z))$ and so $M \vDash \exists z \subseteq_p x*b \; MaxT_a(z, x*b)$.

Therefore $M \vDash I_{MaxT}(x) \to I_{MaxT}(x*b)$.

Consider $x*a$. Then $QT \vdash \neg Tally_b(x*a)$.

We distinguish two cases.

(i) $M \vDash Tally_b(x)$.

Assume $M \vDash v \subseteq_p x*a \;\&\; Tally_a(v)$.

Again, by the definition of $\subseteq_p$, we have that

$\qquad M \vDash v=x*a \;\lor\; vB(x*a) \;\lor\; vE(x*a) \;\lor\; \exists u,w \; uvw=x*a$.

(ia) $M \vDash v=x*a$.

Have $QT^+ \vdash v=x*a \;\&\; Tally_a(v) \to \exists v_1 \; (v=v_1 a \;\&\; Tally_a(v_1))$.

$\Rightarrow M \vDash v=x*a \to \exists v_1 \; (v_1=x \;\&\; Tally_a(v_1))$,

$\Rightarrow$ by (i), $M \vDash Tally_b(x) \;\&\; Tally_a(x)$,

$\Rightarrow$ since $QT \vdash \neg b=a$, by (QT5), $M \vDash (x=b \lor \exists y \; x=b*y) \;\&\; (x=a \lor \exists y \; x=a*y)$,

contradicting (QT2) and (QT4).

(ib)  $M \vDash vB(x*a)$.

$\Rightarrow$ $M \vDash \exists x_1\ vx_1 = x*a$,

$\Rightarrow$ by (QT5) and (QT4),  $M \vDash x_1 = a\ \vee\ aEx_1$,

$\Rightarrow$ $M \vDash va = xa\ \vee\ \exists x_2\ vx_2a = x*a$,

$\Rightarrow$ by (QT3),  $M \vDash v = x\ \vee\ vx_2 = x$,

$\Rightarrow$ $M \vDash v \subseteq_p x$,

$\Rightarrow$ $M \vDash \mathrm{Tally}_a(v)\ \&\ \mathrm{Tally}_b(x)$, contradicting $QT \vdash \neg b = a$.

So both (ia) and (ib) are ruled out.

(ic)  $M \vDash vE(x*a)$.

$\Rightarrow$ $M \vDash \exists x_1\ x_1v = x*a$,

$\Rightarrow$ by (QT5) and (QT4),  $M \vDash v = a\ \vee\ aEv$.

If $M \vDash v = a$, then of course $M \vDash v \subseteq_p a$. Then $M \vDash \mathrm{MaxT}_a(a, x*b)$.

If $M \vDash aEv$, then $M \vDash \exists v_1\ (v_1 a = v\ \&\ x_1 v_1 a = xa)$, so $M \vDash x_1 v_1 = x$  where $M \vDash v_1 \subseteq_p v$.

So $M \vDash \mathrm{Tally}_a(v_1)$. But then also $M \vDash v_1 \subseteq_p x$, so $M \vDash \mathrm{Tally}_b(v_1)$, a contradiction as in (ia).

(id)  $M \vDash \exists u,w\ uvw = x*a$.

$\Rightarrow$ by (QT5) and (QT4),  $M \vDash w = a\ \vee\ aEw$,

$\Rightarrow$ $M \vDash uva = xa\ \vee\ \exists x_1\ uvx_1 a = xa$,

$\Rightarrow$ by (QT3),  $M \vDash uv = x\ \vee\ \exists x_1\ uvx_1 = x$,

$\Rightarrow$ $M \vDash \mathrm{Tally}_a(v)\ \&\ v \subseteq_p x\ \&\ \mathrm{Tally}_b(x)$,  a contradiction as in (ib).

So in Case (i) $M \vDash \mathrm{Tally}_b(x)$ we have

$$M \vDash (\mathrm{Tally}_a(z)\ \&\ \forall v \subseteq_p x*a\ (\mathrm{Tally}_a(v) \to v \subseteq_p z))$$

where z=a.

Suppose (ii) $M \vDash \neg Tally_b(x)$.

Then from the hypothesis $M \vDash I_{MaxT}(x)$ we have $M \vDash \exists z_x \subseteq_p x\; MaxT_a(z_x,x)$.

Assume $M \vDash v \subseteq_p x*a\; \&\; Tally_a(v)$.

Again,   $M \vDash v=x*a\; \lor\; vB(x*a)\; \lor\; vE(x*a)\; \lor\; \exists u,w\; uvw=x*a$.

  (iia) $M \vDash v=x*a$.

$\Rightarrow$ by (4.2), $QT^+ \vdash Tally_a(v) \to \exists v_1\; (v=v_1 a\; \&\; Tally_a(v_1))$,

$\Rightarrow$ $M \vDash v_1 a = x*a$,

$\Rightarrow$ by (QT3), $M \vDash v_1 = x$,

$\Rightarrow$ $M \vDash v_1 \subseteq_p x\; \&\; Tally_a(v_1)$,

$\Rightarrow$ since $M \vDash MaxT_a(z_x,x)$, $M \vDash v_1 \subseteq_p z_x$,

$\Rightarrow$ by (4.11), $M \vDash Tally_a(v_1)\; \&\; Tally_a(z_x)\; \&\; v_1 \subseteq_p z_x \to v_1 a \subseteq_p z_x a$,

$\Rightarrow$ $M \vDash v \subseteq_p z_x a$.

  (iib) $M \vDash vB(x*a)$.

$\Rightarrow$ $M \vDash \exists x_1\; vx_1 = x*a$,

$\Rightarrow$ by (QT5) and (QT4), $M \vDash x_1 = a\; \lor\; aEx_1$,

$\Rightarrow$ $M \vDash va = xa\; \lor\; \exists x_2\; vx_2 a = x*a$,

$\Rightarrow$ by (QT3), $M \vDash v=x\; \lor\; vx_2 = x$,

$\Rightarrow$ $M \vDash v \subseteq_p x$,

$\Rightarrow$ from $M \vDash MaxT_a(z_x,x)$, $M \vDash v \subseteq_p z_x$,

$\Rightarrow$ $M \vDash v \subseteq_p z_x a$.

(iic) $M \vDash vE(x*a)$.

$\Rightarrow$ $M \vDash \exists x_1\ x_1v = x*a$,

$\Rightarrow$ by (QT5) and (QT4), $M \vDash v=a \lor aEv$.

If $M \vDash v=a$, then of course $M \vDash v \subseteq_p z_x a$.

If $M \vDash aEv$, then $M \vDash \exists v_1 (v_1a=v\ \&\ x_1v_1a=xa)$, so by (QT3), $M \vDash x_1v_1=x$ whence $M \vDash v_1 \subseteq_p z_x$ by the choice of $z_x$. Then $M \vDash v \subseteq_p z_x a$ as in (iia).

(iid) $M \vDash \exists u,w\ uvw = x*a$.

$\Rightarrow$ from (QT5) and (QT4), $M \vDash w=a \lor aEw$,

$\Rightarrow$ $M \vDash uva=xa \lor \exists x_1\ uvx_1a=xa$,

$\Rightarrow$ by (QT3), $M \vDash uv=x \lor uvx_1=x$,

$\Rightarrow$ $M \vDash Tally_a(v)\ \&\ v \subseteq_p x$,

$\Rightarrow$ from $M \vDash MaxT_a(z_x,x)$, $M \vDash v \subseteq_p z_x$,

$\Rightarrow$ $M \vDash v \subseteq_p z_x a$.

So we have in Case (ii) $M \vDash \neg Tally_b(x)$ that

$$M \vDash z_x a \subseteq_p xa\ \&\ Tally_a(z_x a)\ \&\ \forall v \subseteq_p x*a\ (Tally_a(v) \to v \subseteq_p z_x a)),$$

which finally shows that $M \vDash I_{MaxT}(x) \to I_{MaxT}(x*a)$.

This completes the proof that $I_{MaxT}(x)$ is a string concept, and the proof of (4.12).

Let $\text{Max}^+\text{T}_a(t,y) \equiv \text{MaxT}_a(t,y) \mathbin{\&} \neg t \subseteq_p y$.

Then we have:

(4.13) For any string concept $I \subseteq I_0$ there is a string concept $I_{\text{Max}+\text{T}} \subseteq I$ such that

$$QT \vdash \forall x \in I_{\text{Max}+\text{T}} \, \exists z \, \text{Max}^+\text{T}_a(z,x).$$

Let $\quad I_{\text{Max}+\text{T}}(x) \equiv I_{\text{MaxT}}(x) \mathbin{\&} I_{3.12}(x)$.

Assume $M \vDash I_{\text{Max}+\text{T}}(x)$.

$\Rightarrow M \vDash I_{\text{MaxT}}(x)$,

$\Rightarrow M \vDash \exists z \, (\text{MaxT}_a(z,x) \mathbin{\&} (\neg \text{Tally}_b(x) \rightarrow z \subseteq_p x))$.

Fix such a z.

If $M \vDash \text{Tally}_b(x)$, then $M \vDash \text{Max}^+\text{T}_a(z,x)$ because

$$QT \vdash \text{Tally}_b(x) \mathbin{\&} \text{Tally}_a(z) \rightarrow \neg z \subseteq_p x.$$

If $M \vDash \neg \text{Tally}_b(x)$, then $M \vDash z \subseteq_p x$. Let $z' = za$. Assume now that $M \vDash v \subseteq_p x$ where $M \vDash \text{Tally}_a(v)$. From the hypothesis $M \vDash \text{MaxT}_a(z,x)$ we have $M \vDash v \subseteq_p z$, hence $M \vDash v \subseteq_p za$. So $M \vDash \text{MaxT}_a(za,x)$. Suppose $M \vDash z' \subseteq_p x$. Then $M \vDash z' = za \subseteq_p z$ by hypothesis $M \vDash \text{MaxT}_a(z,x)$. But this contradicts (3.12). (By 3.13 we may assume that $I_{3.12}$ is downward closed under $\subseteq_p$.)

Therefore $M \vDash \neg z' \subseteq_p x$, and we have $M \vDash Max^+T_a(z',x)$.

This completes the proof of (4.13).

It is easily seen that

$QT \vdash Max^+T_a(z,x)$ & $Tally_a(t)$ & $t \subseteq_p x \rightarrow t<z$.

(4.14) For any string concept $I \subseteq I_0$ there is a string concept $J \subseteq I$ such that

$QT^+ \vdash \forall x \in J \ \forall v,w,y,z \ (Tally_a(z) \ \& \ zBx \ \& \ x=vyw \ \& \ Tally_b(y) \rightarrow z=v \lor zBv)$.

Let $J(x) \equiv I_{3.7}(x) \ \& \ I_{3.8}(x)$.

We may assume that $J$ is downwards closed under $\subseteq_p$.

Assume $M \vDash Tally_a(z) \ \& \ zBx \ \& \ x=vyw \ \& \ Tally_b(y)$ where $M \vDash J(x)$.

$\Rightarrow M \vDash zBx \ \& \ vBx$,

$\Rightarrow$ by (3.8), $M \vDash zBv \lor z=v \lor vBz$.

Suppose, for a reductio, that $M \vDash vBz$.

$\Rightarrow M \vDash \exists v_1 \ (vv_1=z \ \& \ Tally_a(v_1))$,

$\Rightarrow$ from $M \vDash zBx$, $M \vDash \exists x_1 \ zx_1=x$,

$\Rightarrow M \vDash (vv_1)x_1=zx_1=x=vyw$,

$\Rightarrow$ by (3.7), $M \vDash v_1x_1=yw$, which contradicts $M \vDash Tally_a(v_1) \ \& \ Tally_b(y)$ by (QT4).

Therefore, $M \vDash \neg vBz$, and so $M \vDash zBv \lor z=v$, as required.

This completes the proof of (4.14).

(4.15) For any string concept $I \subseteq I_0$ there is a string concept $J \subseteq I$ such that

$QT^+ \vdash \forall x \in J \; \forall v,w,y,z \; (Tally_a(z) \& zEx \& x=vyw \& Tally_b(y) \rightarrow z=w \lor zEw)$.

Let $J(x) \equiv I_{3.6}(x) \& I_{3.10}(x)$.

We may assume that $J$ is downwards closed under $\subseteq_p$.

Assume $M \vDash Tally_a(z) \& zEx \& x=vyw \& Tally_b(y)$ where $M \vDash J(x)$.

$\Rightarrow M \vDash zEx \& wEx$,

$\Rightarrow$ by (3.10), $M \vDash zEw \lor z=w \lor wEz$.

Suppose, for a reductio, that $M \vDash wEz$.

$\Rightarrow M \vDash \exists w_1 (w_1 w = z \& Tally_a(w_1))$,

$\Rightarrow$ from $M \vDash zEx$, $M \vDash \exists x_1 \; x = x_1 z$,

$\Rightarrow M \vDash x_1(w_1 w) = x_1 z = x = vyw$,

$\Rightarrow$ by (3.6), $M \vDash x_1 w_1 = vy$, which contradicts $M \vDash Tally_a(w_1) \& Tally_b(y)$ by (QT4).

Therefore, $M \vDash \neg wEz$, and so $M \vDash zEw \lor z=w$, as required.

This completes the proof of (4.15).

(4.16) For any string concept $I \subseteq I_0$ there is a string concept $J \subseteq I$ such that

$QT \vdash \forall x \in J \, \forall t_1,v,t_2,w_1,w,w_2 \, (Tally_b(t_1) \,\&\, Tally_b(t_2) \,\&\, t_1vt_2=x=w_1awaw_2 \to$

$$\to w \subseteq_p v).$$

Let $J \equiv I_{4.14b} \,\&\, I_{4.15b}$.

Assume $M \vDash t_1vt_2=x=w_1awaw_2$ where $M \vDash Tally_b(t_1) \,\&\, Tally_b(t_2)$ and $M \vDash J(x)$.

By $(4.14^b)$ and $(4.15^b)$ we have that

$$M \vDash (w_1=t_1 \vee t_1Bw_1) \,\&\, (w_2=t_2 \vee t_2Ew_2),$$

$\Rightarrow M \vDash (w_1=t_1 \,\&\, w_2=t_2) \vee (w_1=t_1 \,\&\, t_2Ew_2) \vee (t_1Bw_1 \,\&\, w_2=t_2) \vee$

$$\vee (t_1Bw_1 \,\&\, t_2Ew_2),$$

$\Rightarrow M \vDash t_1vt_2=t_1awat_2 \vee \exists w_4 \, t_1vt_2=t_1awa(w_4t_2) \vee \exists w_3 \, t_1vt_2=(t_1w_3)awat_2 \vee$

$$\vee \exists w_3,w_4 \, t_1vt_2=(t_1w_3)awa(w_4t_2),$$

$\Rightarrow$ by (3.6) and (3.7), $M \vDash v=awa \vee \exists w_4 \, v=awaw_4 \vee \exists w_3 \, v=w_3awa \vee$

$$\vee \exists w_3,w_4 \, v=w_3awaw_4,$$

$\Rightarrow M \vDash w \subseteq_p v$, as required.

This completes the proof of (4.16).

The following proposition generalizes somewhat (4.14) and (4.15):

(4.17) For any string concept $I \subseteq I_0$ there is a string concept $J \subseteq I$ such that

$QT^+ \vdash \forall x \in J \; \forall u,w,y,z \; (Tally_a(u) \;\&\; u \subseteq_p x \;\&\; x = wyz \;\&\; Tally_b(y) \rightarrow u \subseteq_p w \;\vee\; u \subseteq_p z)$.

Let $I\hat{\;}(x) \equiv I_{4.14}(x) \;\&\; I_{4.15}(x)$.

We let $J(x)$ abbreviate

$I\hat{\;}(x) \;\&\; \forall u,w,y,z \; (Tally_a(u) \;\&\; u \subseteq_p x \;\&\; x = wyz \;\&\; Tally_b(y) \rightarrow u \subseteq_p w \;\vee\; u \subseteq_p z)$.

That $QT^+ \vdash J(a)$ and $QT^+ \vdash J(b)$ follows trivially from (QT2).

Assume $M \vDash J(x)$.

Then

$M \vDash I\hat{\;}(x) \;\&\; \forall u,w,y,z \; (Tally_a(u) \;\&\; u \subseteq_p x \;\&\; x = wyz \;\&\; Tally_b(y) \rightarrow u \subseteq_p w \;\vee\; u \subseteq_p z)$.

Assume $M \vDash Tally_a(u) \;\&\; u \subseteq_p (wyz)*a$.

There are four cases:

  (a1) $M \vDash u = (wyz)*a$.

Contradicts (QT4) because by hypothesis $M \vDash Tally_a(u) \;\&\; Tally_b(y)$.

  (a2) $M \vDash uB((wyz)*a)$.

$\Rightarrow$ since $M \vDash Tally_a(u) \;\&\; Tally_b(y)$, by (4.14), $M \vDash u = w \;\vee\; uBw$,

$\Rightarrow M \vDash u \subseteq_p w$.

  (a3) $M \vDash uE((wyz)*a)$.

$\Rightarrow M \vDash \exists u_1 \; u_1 u = wyz*a$,

$\Rightarrow$ by (QT5) and (QT4), $M \vDash u = a \;\vee\; aEu$.

If $M \vDash u = a$, then $M \vDash u \subseteq_p (z*a)$, as required.

If $M \vDash aEu$, then $M \vDash \exists u_2 (u_2 a = u$ & $Tally_a(u_2))$, so $M \vDash u_1(u_2 a) = wyz^*a$. But then, by (QT3), $M \vDash u_1 u_2 = wyz$. Since $M \vDash Tally_b(y)$ & $Tally_a(u_2)$, it follows by (4.15) that $M \vDash u_2 = z \lor u_2 E z$.

If $M \vDash u_2 = z$, then $M \vDash u = u_2 a = z^*a$, so $M \vDash u \subseteq_p (z^*a)$.

If $M \vDash u_2 E z$, then $M \vDash \exists z_1\ z = z_1 u_2$. Then $M \vDash z_1 u_2 a = za$.

Hence $M \vDash u_2 a \subseteq_p (z^*a)$, and so $M \vDash u \subseteq_p (z^*a)$.

  (a4)  $M \vDash \exists u_1, u_2\ (wyz)^*a = u_1 u u_2$.

$\Rightarrow$ by (QT5) and (QT4), $M \vDash u_2 = a \lor aEu_2$.

If $M \vDash u_2 = a$, then $M \vDash u_1 u a = wyz^*a$, so by (QT3), $M \vDash u_1 u = wyz$.

Then $M \vDash u \subseteq_p wyz = x$.

From the hypothesis $M \vDash J(x)$ we then have that

$$M \vDash u \subseteq_p w \lor u \subseteq_p z,$$

whence $M \vDash u \subseteq_p w \lor u \subseteq_p (z^*a)$ since, of course, $M \vDash z \subseteq_p (z^*a)$.

If $M \vDash aEu_2$, then $M \vDash \exists u_3\ u_3 a = u_2$.

So $M \vDash u_1 u(u_3 a) = wyz^*a$, whence $M \vDash u_1 u u_3 = wyz = x$ by (QT3).

We then have $M \vDash u \subseteq_p w \lor u \subseteq_p z$ from the hypothesis $M \vDash J(x)$, whence $M \vDash u \subseteq_p w \lor u \subseteq_p (z^*a)$ as required.

So from (1)-(4) we have that $QT^+ \vdash \forall x (J(x) \to J(x^*a))$.

Assume now $M \vDash \Gamma(u)$ & $Tally_a(u)$ & $u \subseteq_p wb(x^*b)$.

  (b1)  $M \vDash u = (wyz)^*b$.

Contradicts (QT4) because by hypothesis $M \vDash Tally_a(u)$.

  (b2)  $M \vDash uB(wb(x^*b))$.

Exactly as (a2).

   (b3)  $M \vDash uE(wyz*b)$.

$\Rightarrow$  $M \vDash \exists u_1\ u_1u=wyz*b$,  which contradicts the hypothesis $M \vDash \text{Tally}_a(u)$.

   (b4)  $M \vDash \exists u_1,u_2\ wyz*b=u_1uu_2$.

$\Rightarrow$ by (QT5) and (QT4), $M \vDash u_2=b\ \vee\ bEu_2$.

If $M \vDash u_2=b$, then $M \vDash u_1ub=wyz*b$, so $M \vDash u_1u=wyz$ by (QT3).

Then  $M \vDash u \subseteq_p wyz = x$.

From the hypothesis $M \vDash J(x)$ we then have that

$$M \vDash u \subseteq_p w\ \vee\ u \subseteq_p z$$

and so $M \vDash u \subseteq_p w\ \vee\ u \subseteq_p z*b$.

If $M \vDash bEu_2$, then $M \vDash \exists u_3\ u_3b=u_2$.

So $M \vDash u_1u(u_3b)=wyz*b$, whence $M \vDash u_1uu_3=wyz=x$ by (QT3).

Then again     $M \vDash u \subseteq_p w\ \vee\ u \subseteq_p z$

from the hypothesis $M \vDash J(x)$, whence

$$M \vDash u \subseteq_p w\ \vee\ u \subseteq_p (z*b)$$

as required.

So from (b1)-(b4) we have that $QT^+ \vdash \forall x\ (J(x) \rightarrow J(x*b))$.

This completes the proof that $J(x)$ is a string concept.

This completes the proof of (4.17).

(4.18) For any string concept $I \subseteq I_0$ there is a string concept $J \subseteq I$ such that

$QT^+ \vdash \forall x \in J \ \forall t_1, t_2, w, w_1, w_2, y, z \ \neg(Tally_b(t_1) \ \& \ Tally_b(t_2) \ \&$

$\& \ w_1 a t_1 w a t_2 y = x = w_2 a t_2 z \ \& \ Max^+T_b(t_1, w_2) \ \& \ Max^+T_b(t_2, z))$.

Let $J \equiv I_{3.10} \ \& \ I_{4.15b}$.

Assume, for a reductio, that $M \vDash w_1 a t_1 w a t_2 y = x = w_2 a t_2 z$ where

$M \vDash Tally_b(t_1) \ \& \ Tally_b(t_2) \ \& \ Max^+T_b(t_1, w_2) \ \& \ Max^+T_b(t_2, z)$ and $M \vDash J(x)$.

$\implies \ M \vDash (t_2 y)Ex \ \& \ zEx$,

$\implies$ by (3.10), $M \vDash (t_2 y)Ez \ \lor \ t_2 y = z \ \lor \ zE(t_2 y)$.

If $M \vDash (t_2 y)Ez \ \lor \ t_2 y = z$, then $M \vDash t_2 \subseteq_p t_2 y \subseteq_p z$, contradicting $M \vDash Max^+T_b(t_2, z)$.

Therefore, $M \vDash zE(t_2 y)$.

$\implies \ M \vDash \exists z_1 \ t_2 y = z_1 z$,

$\implies \ M \vDash w_1 a t_1 w a (z_1 z) = w_2 a t_2 z$,

$\implies$ by (3.6), $M \vDash w_1 a t_1 w a z_1 = w_2 a t_2$,

$\implies$ by (4.15$^b$), $M \vDash z_1 = t_2 \ \lor \ t_2 E z_1$,

$\implies \ M \vDash w_1 a t_1 w a t_2 = w_2 a t_2 \ \lor \ \exists z_2 \ w_1 a t_1 w a z_2 t_2 = w_2 a t_2$,

$\implies$ by (3.6), $M \vDash w_1 a t_1 w = w_2 \ \lor \ w_1 a t_1 w a z_2 = w_2 a$,

$\implies \ M \vDash z_2 = a \ \lor \ a E z_2$,

$\implies \ M \vDash w_1 a t_1 w = w_2 \ \lor \ w_1 a t_1 w a = w_2 \ \lor \ \exists z_3 \ w_1 a t_1 w a z_3 = w_2$,

$\implies \ M \vDash t_1 \subseteq_p w_2$, contradicting $M \vDash Max^+T_b(t_1, w_2)$.

This completes the proof of (4.18).

(4.19) For any string concept $I \subseteq I_0$ there is a string concept $J \subseteq I$ such that

$QT^+ \vdash \forall x \in J \; \forall t_1, t_2, w_1, w_2, y_1, y_2 \; (w_1 a t_1 y_1 = x = w_2 a t_2 y_2 \; \&$

$\& \; Max^+T_b(t_1, w_2) \; \& \; Max^+T_b(t_2, w_1) \rightarrow (w_2 a)B(w_1 a) \; v \; w_1 = w_2).$

Let $J \equiv I_{3.8} \; \& \; I_{4.14b}$.

Assume $M \vDash w_1 a t_1 y_1 = x = w_2 a t_2 y_2$ where $M \vDash Max^+T_b(t_1, w_2) \; \& \; Max^+T_b(t_2, w_1)$ and $M \vDash J(x)$.

$\Rightarrow \; M \vDash (w_1 a)Bx \; \& \; (w_2 a)Bx,$

$\Rightarrow$ by (3.8), $M \vDash (w_1 a)B(w_2 a) \; v \; w_1 = w_2 \; v \; (w_2 a)B(w_1 a).$

Suppose, for a reductio, that $M \vDash (w_1 a)B(w_2 a).$

$\Rightarrow \; M \vDash \exists w \; w_1 a w = w_2 a,$

$\Rightarrow \; M \vDash w_1 a t_1 y_1 = x = (w_1 a w) t_2 y_2,$

$\Rightarrow$ by (3.7), $M \vDash t_1 y_1 = x = w t_2 y_2,$

$\Rightarrow$ from $M \vDash w_1 a w = w_2 a$, $M \vDash w = a \; v \; a E w_2.$

Now, $M \vDash w \neq a$ because $M \vDash Tally_b(t_1).$

$\Rightarrow M \vDash aEw,$

$\Rightarrow M \vDash \exists w' \; (w'a = w \; \& \; t_1 y_1 = x = (w'a)t_2 y_2),$

$\Rightarrow$ by (4.14$^b$), $M \vDash w' = t_1 \; v \; t_1 B w',$

$\Rightarrow M \vDash t_1 \subseteq_p w'.$

But we have $M \vDash w_1 a w' a = w_2 a.$

$\Rightarrow M \vDash w_1 a w' = w_2,$

$\Rightarrow M \vDash t_1 \subseteq_p w' \subseteq_p w_2$, contradicting hypothesis $M \vDash Max^+T_b(t_1, w_2).$

Therefore   $M \vDash w_1a = w_2a \lor (w_2a)B(w_1a)$,  whence

$$M \vDash w_1 = w_2 \lor (w_2a)B(w_1a),$$

as required.

This completes the proof of (4.19).

(4.20) For any string concept $I \subseteq I_0$ there is a string concept $J \subseteq I$ such that

$QT \vdash \forall x \in J\ \forall t_1,t_2,w_1,w_2,y_1,y_2\ (w_1at_1y_1=x=w_2at_2y_2\ \&$

$\&\ Max^+T_b(t_1,w_2)\ \&\ Max^+T_b(t_2,w_1) \rightarrow w_1=w_2).$

Let $J \equiv I_{4.19}$.

Assume $M \vDash w_1at_1y_1=x=w_2at_2y_2$ where $M \vDash Max^+T_b(t_1,w_2)\ \&\ Max^+T_b(t_2,w_1)$ and $M \vDash J(x)$.

$\Rightarrow\ M \vDash Tally_b(t_2)$,

$\Rightarrow$ by (4.19), $M \vDash (w_2a)B(w_1a)\ \vee\ w_1=w_2$,

$\Rightarrow$ from $M \vDash Max^+T_b(t_1,w_2)$, $M \vDash Tally_b(t_1)$,

$\Rightarrow$ by (4.19), $M \vDash (w_1a)B(w_2a)\ \vee\ w_1=w_2$.

Now, $M \vDash \neg((w_1a)B(w_2a)\ \&\ (w_2a)B(w_1a))$. Otherwise, $M \vDash (w_1a)B(w_1a)$, contradicting $M \vDash (w_1a) \in J \subseteq I_0$.

Likewise, $M \vDash \neg((w_1a)B(w_2a)\ \&\ w_1=w_2)$ and $M \vDash \neg((w_2a)B(w_1a)\ \&\ w_1=w_2)$.

Therefore, $M \vDash w_1=w_2$, as required.

This completes the proof of (4.20).

(4.21) For any string concept $I \subseteq I_0$ there is a string concept $J \subseteq I$ such that

$\quad$ QT ⊢ $\forall x \in J \, \forall t_1,t_2,t_3,u,v,w$ (Max$^+$T$_b$(t$_1$,aua) & Tally$_b$(t$_2$) &

$$\& \; x=t_1auat_2 \; \& \; x=wat_3avat_3 \; \rightarrow \; \neg MaxT_b(t_3,x)).$$

Let $J \equiv I_{4.16}$.

Assume $M \vDash$ Max$^+$T$_b$(t$_1$,aua) & Tally$_b$(t$_2$) where $M \vDash x \in J$.

Let $M \vDash t_1auat_2 = x = wat_3avat_3$.

Suppose for a reductio that $M \vDash MaxT_b(t_3,x)$.

$\Rightarrow$ by (4.16), $M \vDash t_3 \subseteq_p u$,

$\Rightarrow M \vDash t_1 \leq t_3$,

$\Rightarrow M \vDash t_1 \subseteq_p u$, contradicting $M \vDash$ Max$^+$T$_b$(t$_1$,aua).

This completes the proof of (4.21).

(4.22) For any string concept $I \subseteq I_0$ there is a string concept $J \subseteq I$ such that

$$QT \vdash \forall x \in J \, \forall t,u(x=tauat \, \& \, Max^+T_b(t,aua) \rightarrow MaxT_b(t,x)).$$

Let $J \equiv I_{4.17b}$.

Assume $M \vDash x=tauat \, \& \, Max^+T_b(t,aua)$ where $M \vDash J(x)$.

$\Longrightarrow$ $M \vDash Tally_b(t)$.

Let $M \vDash t' \subseteq_p x \, \& \, Tally_b(t')$.

$\Longrightarrow$ $M \vDash t' \subseteq_p x = tauat$,

$\Longrightarrow$ by (4.17$^b$), $M \vDash t' \subseteq_p t \, \vee \, t' \subseteq_p uat$,

$\Longrightarrow$ by (4.17$^b$), $M \vDash t' \subseteq_p t \, \vee \, t' \subseteq_p u$,

$\Longrightarrow$ $M \vDash t' \subseteq_p t \, \vee \, t' \subseteq_p aua$,

$\Longrightarrow$ from $M \vDash Max^+T_b(t,aua)$, $M \vDash t' \subseteq_p t$,

$\Longrightarrow$ $M \vDash MaxT_b(t,aua)$.

This completes the proof of (4.22).

(4.23) For any string concept $I \subseteq I_0$ there is a string concept $J \subseteq I$ such that

$QT \vdash \forall z_1, z_2 \in J\ \forall u_1, u_2, v_1, v_2\ (Tally_a(z_1)\ \&\ Tally_a(z_2)\ \&\ z_1 u_1 v_1 = z_2 u_2 v_2\ \&$

$\&\ bBu_1\ \&\ bBu_2 \rightarrow z_1 = z_2).$

Let $J(x) \equiv I_{3.7}(x)\ \&\ I_{4.6a}(x)$.

Assume $M \vDash z_1 u_1 v_1 = z_2 u_2 v_2$

where $M \vDash Tally_a(z_1)\ \&\ Tally_a(z_2)\ \&\ bBu_1\ \&\ bBu_2$.

$\Rightarrow$ by $(4.6^a)$, $M \vDash z_1 \leq z_2\ \vee\ z_2 \leq z_1$.

Suppose $M \vDash z_1 < z_2$.

$\Rightarrow$ because $M \vDash Tally_a(z_1)\ \&\ Tally_a(z_2)$, $M \vDash \exists z_3 (z_1 z_3 = z_2\ \&\ Tally_a(z_3))$,

$\Rightarrow M \vDash z_1 u_1 v_1 = (z_1 z_3) u_2 v_2$,

$\Rightarrow$ by (3.7), $M \vDash u_1 v_1 = z_3 u_2 v_2$,

which contradicts (QT4) because $M \vDash bBu_1\ \&Tally_a(z_3)$.

Suppose $M \vDash z_2 < z_1$.

$\Rightarrow$ because $M \vDash Tally_a(z_1)\ \&\ Tally_a(z_2)$, $M \vDash \exists z_3(z_2 z_3 = z_1\ \&\ Tally_a(z_3))$,

$\Rightarrow M \vDash (z_2 z_3) u_1 v_1 = z_2 u_2 v_2$,

$\Rightarrow$ by (3.7), $M \vDash z_3 u_1 v_1 = u_2 v_2$,

again a contradiction because $M \vDash bBu_2\ \&Tally_a(z_3)$.

Hence $M \vDash \neg(z_1 < z_2)\ \&\ \neg(z_2 < z_1)$.

It follows that $M \vDash z_1 = z_2$.

This completes the proof of (4.23).

(4.24) For any string concept $I \subseteq I_0$ there is a string concept $J \subseteq I$ such that

$QT \vdash \forall z_3, z_4 \in J \; \forall u_1, u_2, v_1, v_2 \; (Tally_a(z_3) \& Tally_a(z_4) \& v_1 u_1 z_3 = v_2 u_2 z_4 \; \& \; \& \; bEu_1 \& bEu_2 \rightarrow z_3 = z_4)$

Let $J(x) \equiv I_{3.6}(x) \& I_{4.6a}(x) \& I_{4.10a}(x)$.

Assume $M \vDash v_1 u_1 z_3 = v_2 u_2 z_4$

where $M \vDash Tally_a(z_3) \& Tally_a(z_4) \& bEu_1 \& bEu_2$.

$\Rightarrow$ by (4.6$^a$), $M \vDash z_3 \leq z_4 \vee z_4 \leq z_3$.

Suppose $M \vDash z_3 < z_4$.

$\Rightarrow$ because $M \vDash Tally_a(z_3) \& Tally_a(z_4)$, $M \vDash \exists z_5 (z_3 z_5 = z_4 \& Tally_a(z_5))$,

$\Rightarrow M \vDash v_1 u_1 z_3 = v_2 u_2 (z_3 z_5)$,

$\Rightarrow$ by (4.10$^a$), $M \vDash z_3 z_5 = z_5 z_3$,

$\Rightarrow M \vDash v_1 u_1 z_3 = v_2 u_2 (z_5 z_3)$,

$\Rightarrow$ by (3.6), $M \vDash v_1 u_1 = v_2 u_2 z_5$, a contradiction because $M \vDash bEu_1 \& Tally_a(z_5)$.

Suppose $M \vDash z_4 < z_3$.

$\Rightarrow$ because $M \vDash Tally_a(z_3) \& Tally_a(z_4)$, $M \vDash \exists z_5 (z_4 z_5 = z_3 \& Tally_a(z_5))$,

$\Rightarrow M \vDash v_1 u_1 (z_4 z_5) = v_2 u_2 z_4$,

$\Rightarrow$ by (4.10$^a$), $M \vDash z_4 z_5 = z_5 z_4$,

$\Rightarrow M \vDash v_1 u_1 (z_5 z_4) = v_2 u_2 z_4$,

$\Rightarrow$ by (3.6), $M \vDash v_1 u_1 z_5 = v_2 u_2$, a contradiction because $M \vDash bEu_2 \& Tally_a(z_5)$.

Therefore, $M \vDash z_3 = z_4$.

This completes the proof of (4.24).

## 5. Coding sets by strings

We wish to have strings $y_1, y_2, y_3, \ldots$ coded by

$$t_1 a y_1 a t_2 a y_2 a t_3 a y_3 a \ldots$$

so that the successive $y_i$'s are embedded in "frames" separated by sufficiently long b-tallies $t_i$, where $t_i$'s are such that $\neg t_i \subseteq_p a y_i a$ and $t_j < t_i$ for $j<i$. Each marker $t_i$ is supposed to be longer than any b-tally in $y_i$ and any preceding marker. In this sense the coding will work like a step-ladder: starting with the b-tally that precedes the first occurrence of the letter a in w, each next longer b-tally in the code is a successive step of the ladder marking off a frame that corresponds to another member of the coded set.

So far we have been talking about codes of sets of strings <u>informally</u>. What specific assumptions about concatenated strings suffice to formally validate the claim that the coding works? It will pay off to make these assumptions as weak as possible.

We will show that all the necessary reasoning can be carried out in $QT^+$.

We first define when a string u is a <u>preframe indexed by</u> t:

$$\text{Pref}(u,t) \equiv \text{Tally}_b(t) \ \& \ \exists y \subseteq_p u \ (aya=u \ \& \ \text{Max}^+T_b(t,u));$$

when $t_1ut_2$ is (the) <u>first frame</u> in the string x, $\text{Firstf}(x,t_1,u,t_2)$:

$\text{Pref}(u,t_1) \ \& \ \text{Tally}_b(t_2) \ \& \ ((t_1=t_2 \ \& \ t_1ut_2=x) \ \vee \ (t_1<t_2 \ \& \ (t_1ut_2a)Bx));$

when $t_1ut_2$ is (the) <u>last frame</u> in x, $\text{Lastf}(x,t_1,u,t_2)$:

$\text{Pref}(u,t_1) \ \& \ \text{Tally}_b(t_2) \ \& \ t_1=t_2 \ \& \ (t_1ut_2=x \ \vee \ \exists w \ (wat_1ut_2=x \ \& \ \text{Max}^+T_b(t_1,w)));$

and when $t_1ut_2$ is an <u>intermediate frame</u> in x, $\text{Intf}(x,w,t_1,u,t_2)$:

$\text{Pref}(u,t_1) \ \& \ \text{Tally}_b(t_2) \ \& \ t_1<t_2 \ \& \ \exists w_1(wat_1ut_2aw_1=x) \ \& \ \text{Max}^+T_b(t_1,w).$

We then define when a string $u$ is <u>$t_1,t_2$-framed</u> in x:

$\text{Fr}(x,t_1,u,t_2) \equiv \text{Firstf}(x,t_1,u,t_2) \ \vee \ \exists w \ \text{Intf}(x,w,t_1,u,t_2) \ \vee \ \text{Lastf}(x,t_1,u,t_2),$

We say that $t_1$ is the <u>initial</u>, and $t_2$ <u>terminal tally marker</u> in the frame.

We then define "t <u>envelops</u> x", $\text{Env}(t,x)$, to be the conjunction of the following five conditions:

(a) $\text{MaxT}_b(t,x)$                          "t is a longest b-tally in x",

(b) $\exists u \subseteq_p x \ \exists t_1,t_2 \ \text{Firstf}(x,t_1,u,t_2)$     "x has a first frame",

(c) $\exists u \subseteq_p x \ \text{Lastf}(x,t,u,t)$             "x has a last frame with t as its initial

                                                       and terminal marker"

(d) $\forall u \subseteq_p x \ \forall t_1,t_2,t_3,t_4 \ (\text{Fr}(x,t_1,u,t_2) \ \& \ \text{Fr}(x,t_3,u,t_4) \to t_1=t_3)$

      "different initial tally markers frame distinct strings",

(e) $\forall u_1, u_2 \subseteq_p x \; \forall t', t_1, t_2 \; (Fr(x,t',u_1,t_1) \; \& \; Fr(x,t',u_2,t_2) \; \rightarrow \; u_1 = u_2)$

"distinct strings are framed by different initial tally markers"

We then say  x <u>is a set code</u>  if x is aa or else x is enveloped by some b-tally:

$Set(x) \equiv x = aa \; \lor \; \exists t \subseteq_p x \; Env(t,x)$.

Finally, we say that a string y <u>is a member of the set coded by</u> string x if x is enveloped by some b-tally t and the juxtaposition of the string y with single tokens of digit a is framed in x:

$y \; \varepsilon \; x \equiv \exists t \subseteq_p x \; (Env(t,x) \; \& \; \exists u \subseteq_p x \; \exists t_1, t_2 (Fr(x,t_1,u,t_2) \; \& \; u = aya))$.

If a set of strings X  is extended by adjoining  a string  y to obtain another set  $Y = X \cup \{y\}$, then  a code for Y can be picked so that a given code of X  will be its <u>initial segment</u>.  This will have to be formally proved in QT$^+$.

We let   $x \sim y \equiv Set(x) \; \& \; Set(y) \; \& \; \forall w(w \; \varepsilon \; x \leftrightarrow w \; \varepsilon \; y)$.

(5.1) $QT^+ \vdash \forall z \in I_0\ \forall t',t'',y(\text{Firstf}(z,t',aya,t'')\ \&\ (t'ayat''a)Bz \rightarrow t<t'\ \&\ z \neq t'ayat'')$.

Assume $M \vDash \text{Firstf}(z,t',aya,t'')\ \&\ (t'ayat''a)Bz$ where $M \vDash I_0(z)$.

$\Rightarrow M \vDash \text{Pref}(ava,t')\ \&\ \text{Tally}_b(t')\ \&$

$\&\ ((t'=t''\ \&\ z=t'ayat'') \lor (t'<t''\ \&\ (t'ayat''a)Bz))$.

Suppose, for a reductio, that $M \vDash z=t'ayat''$.

$\Rightarrow$ from $M \vDash (t'ayat''a)Bz$, $M \vDash zBz$, contradicting $M \vDash I_0(z)$.

$\Rightarrow M \vDash z \neq t'ayat''$,

$\Rightarrow M \vDash \neg(t'=t''\ \&\ z=t'ayat'')$,

$\Rightarrow M \vDash t<t'\ \&\ z \neq t'ayat''$.

This completes the proof of (5.1).

(5.2) For any string I⊆I₀ there is a string concept  J⊆I such that

   $QT^+ \vdash \forall z \in J \, \forall t,x,u \subseteq_p z$ (Env(t,x) & z=xt'utt' & Pref(u,tt') → $MaxT_b$(tt',z)).

Let  $J(x) \equiv I_{4.17b}(x)$.

Assume  M ⊨ z=xt'utt'  where  M ⊨ Env(t,x) & Pref(u,tt')  and M ⊨ J(x).

First note that from M ⊨ Env(t,x) we have  M ⊨ $MaxT_b$(t,x),  and from M ⊨ Pref(u,tt')  we have

   M ⊨ $\exists u_0$ u=a$u_0$a &  $Max^+T_b$(tt',u).

Let M ⊨ $Tally_b$(t'') & t''$\subseteq_p$z=xt'utt'=xt'(a$u_0$a)tt'.

⇒ by (4.17$^b$),  M ⊨ t''$\subseteq_p$xt'  v  t''$\subseteq_p u_0$  v  t''$\subseteq_p$tt'.

From  M ⊨ $Max^+T_b$(tt',u)  and  M ⊨ t''$\subseteq_p u_0 \subseteq_p$u, we have  M ⊨ t''$\subseteq_p$tt'.

On the other hand, from M ⊨ Env(t,x) we have  M ⊨ (at)Ex.

⇒  M ⊨ $\exists x_1$ $x_1$at=x,

⇒ from  M ⊨ t''$\subseteq_p$xt'=($x_1$at)t',  by (4.17$^b$),

   M ⊨ t''$\subseteq_p x_1$  v  t''$\subseteq_p$tt',

⇒ from  M ⊨ $MaxT_b$(t,x) & t''$\subseteq_p x_1 \subseteq_p$x,  M ⊨ t''$\subseteq_p$tt'.

Hence in any case  M ⊨ t''$\subseteq_p$tt', and therefore  M ⊨ $MaxT_b$(tt',z), as required.

This completes the proof of (5.2).

(5.3) For any string concept I⊆I₀ there is a string concept J⊆I such that

$$QT^+ \vdash \forall x \in J \; \forall t,t',t'',v,v',w' \subseteq_p x \; (MaxT_b(t,x) \; \& \; Lastf(x,t,v',t) \; \& \; wat'vt''=x \;\&$$
$$\& \; Tally_b(t'') \; \& \; Pref(v,t') \; \& \; Max^+T_b(t',w) \; \to \; t'=t=t'' \; \& \; v=v').$$

Let J(x) abbreviate

$$I^{\subseteq p}_{3.6}(x) \; \& \; I^{\subseteq p}_{3.7}(x) \; \& \; I_{3.10}(x) \; \& \; I^{\subseteq p}_{4.5}(x) \; \& \; I^{\subseteq p}_{4.10}(x) \; \& \; I_{4.14b}(x) \; \& \; I^{\subseteq p}_{4.16}(x) \; \& \; I^{\subseteq p}_{4.23b}(x) \;\&$$
$$\& \; I_{4.24b}(x).$$

Assume  M ⊨ MaxT_b(t,x) & Lastf(x,t,v',t) & wat'vt''=x

where  M ⊨ Tally_b(t'') & Pref(v,t') & Max⁺T_b(t',w)  and  M ⊨ J(x).

From  M ⊨ Lastf(x,t,v',t)  we have

  M ⊨ Pref(v',t) & (x=tv't ∨ ∃w'(w'atv't=x & Max⁺T_b(t,w'))).

Suppose  M ⊨ x=tv't.

⇒  M ⊨ wat'vt''=x=tv't,

⇒ by (4.14ᵇ),  M ⊨ w=t ∨ tBw,

⇒  M ⊨ t⊆_p w,

⇒ from M ⊨ Max⁺T_b(t',w),  M ⊨ t<t'.

But from M ⊨ MaxT_b(t,x),  M ⊨ t'≤t,

⇒  M ⊨ t<t'≤t,  contradicting  M ⊨ t ∈ I ⊆ I₀.

Therefore  M ⊨ ¬(x=tv't).

So we may assume that  M ⊨ ∃w'(w'atv't=x=wat'vt'').

From  M ⊨ Pref(v',t) & Pref(v,t')  we have  M ⊨ ∃v'' v'=av''a & ∃v₀ v=av₀a,

⇒  M ⊨ w'at(av''a)t=w'atv't=x=wat'vt''=wat'(av₀a)t'',

⇒ by (4.24$^b$), M ⊨ t=t'',

⇒ M ⊨ w'atv't=wat'vt.

From M ⊨ MaxT$_b$(t,x), M ⊨ t'≤t, and we also have that

$$M ⊨ (tv't)Ex \ \& \ (t'vt)Ex,$$

⇒ by (3.10), M ⊨ (tv't)E(t'vt) v tv't=t'vt v (t'vt)E(tv't).

Case 1. M ⊨ (tv't)E(t'vt).

⇒ M ⊨ ∃x$_1$ x$_1$tv't=t'vt.

Now, M ⊨ Tally$_b$(x$_1$) v ¬Tally$_b$(x$_1$).

(1a) M ⊨ Tally$_b$(x$_1$).

⇒ by (4.5), M ⊨ Tally$_b$(x$_1$t), and M ⊨ x$_1$t(av''a)t=t'vt,

⇒ by (4.23$^b$), M ⊨ x$_1$t=t',

⇒ by (4.10), M ⊨ tx$_1$=t',

⇒ M ⊨ t<t'≤t, contradicting M ⊨ t ∈ I ⊆ I$_0$.

(1b) M ⊨ ¬Tally$_b$(x$_1$).

⇒ M ⊨ a⊆$_p$x$_1$.

We distinguish cases based on the definition of ⊆$_p$:

We have that M ⊨ ¬(x$_1$=a) & ¬aBx$_1$ because M ⊨ Tally$_b$(t'). Suppose M ⊨ aEx$_1$.

⇒ M ⊨ ∃x$_2$ x$_2$a=x$_1$,

⇒ M ⊨ (x$_2$a)tv't=t'vt,

⇒ M ⊨ x$_2$at(av''a)t=t'vt,

⇒ by (4.16), M ⊨ t⊆$_p$v,

⇒ M ⊨ t'≤t⊆$_p$v, contradicting M ⊨ Pref(v,t').

The same argument yields a contradiction if $M \vDash \exists x_2, x_3 \; x_2 a x_3 = x_1$.

This completes Case 1.

<u>Case 2.</u> $M \vDash tv't = t'vt$.

$\Rightarrow$ by (4.23$^b$), $M \vDash t = t'$,

$\Rightarrow M \vDash tv't = tvt$,

$\Rightarrow$ by (3.6) and (3.7), $M \vDash v' = v$.

Hence $M \vDash t' = t = t'' \; \& \; v' = v$, as required.

<u>Case 3.</u> $M \vDash (t'vt)E(tv't)$.

$\Rightarrow M \vDash \exists x_1 \; x_1 t'vt = tv't$.

Again, $M \vDash \text{Tally}_b(x_1) \vee \neg \text{Tally}_b(x_1)$.

(3a) $M \vDash \text{Tally}_b(x_1)$.

$\Rightarrow M \vDash w'a(x_1 t'vt) = w'atv't = x = wat'vt$,

$\Rightarrow$ by (3.6), $M \vDash w'ax_1 = wa$, a contradiction.

(3b) $M \vDash \neg \text{Tally}_b(x_1)$.

$\Rightarrow M \vDash a \subseteq_p x_1$.

Again, $M \vDash \neg(x_1 = a) \; \& \; \neg aBx_1$. Suppose that $M \vDash aEx_1$.

$\Rightarrow M \vDash \exists x_2 \; (x_2 a)t'vt = tv't$,

$\Rightarrow M \vDash x_2 at'(av_0 a)t = tv't$,

$\Rightarrow$ by (4.16), $M \vDash t' \subseteq_p v'$,

$\Rightarrow$ from $M \vDash \text{Pref}(v, t')$, $M \vDash t' < t$.

But also from $M \vDash x_2 at'vt = tv't$, we have by (4.14$^b$), that

$$M \vDash x_2 = t \vee tBx_2,$$

$\Rightarrow M \vDash t \subseteq_p x_2$.

Now, we have $M \vDash x_2 a t'vt = tv't$ and $M \vDash w'atv't = x = wat'vt$,

$\Rightarrow M \vDash w'a(x_2 a t'vt) = watv't$,

$\Rightarrow$ by (3.6), $M \vDash w'ax_2 = w$,

$\Rightarrow M \vDash t \subseteq_p x_2 \subseteq_p w$,

$\Rightarrow M \vDash Max^+ T_b(t',w)$, $M \vDash t < t'$,

$\Rightarrow M \vDash t < t' < t$, contradicting $M \vDash t \in I \subseteq I_0$.

An analogous argument derives a contradiction if $M \vDash \exists x_2, x_3\ x_2 a x_3 = x_1$.

This completes Case 3 and the proof of (5.3).

(5.4) For any string concept $I \subseteq I_0$ there is a string concept $J \subseteq I$ such that

$QT^+ \vdash \forall z \in J \ \forall x,v,t,t',t''(Env(t,x) \ \& \ xBz \ \& \ Firstf(z,t',ava,t'')) \rightarrow$

$\rightarrow \exists t_2 \ Firstf(x,t',ava,t_2))$.

Let $J \equiv I_{3.8} \ \& \ I_{4.5} \ \& \ I_{4.16} \ \& \ I_{4.23b} \ \& \ I_{4.24b}$.

Assume $M \vDash Firstf(z,t',ava,t'')$ and let $M \vDash xBz \ \& \ Env(t,x)$ where $M \vDash J(z)$.

$\Longrightarrow M \vDash Pref(ava,t') \ \& \ ((t'=t'' \ \& \ z=t'avat'') \lor (t'<t'' \ \& \ (t'avat''a)Bz))$,

$\Longrightarrow$ from $M \vDash Env(t,x)$, $M \vDash \exists v',t_1,t_2 \ Firstf(x,t_1,av'a,t_2)$,

$\Longrightarrow M \vDash Pref(av'a,t_1) \ \& \ ((t_1=t_2 \ \& \ x=t_1av'at_2) \lor (t_1<t_2 \ \& \ (t_1av'at_2a)Bx))$,

$\Longrightarrow M \vDash (t_1a)Bx \ \& \ (t'a)Bz \ \& \ Tally_b(t_1) \ \& \ Tally_b(t')$,

$\Longrightarrow M \vDash (t_1a)Bz \ \& \ (t'a)Bz$,

$\Longrightarrow M \vDash \exists z_2,z_3 \ t_1az_2=z=t'az_3$,

$\Longrightarrow$ by $(4.23^b)$, $M \vDash t_1=t'$.

There are two cases:

(1) $M \vDash t'=t'' \ \& \ z=t'avat''$.

(1a) $M \vDash t_1=t_2 \ \& \ x=t_1av'at_2$.

$\Longrightarrow$ from $M \vDash xBz$,

$\Longrightarrow M \vDash \exists z_1 \ xz_1=z$,

$\Longrightarrow M \vDash t_1av'at_2z_1=xz_1=z=t'avat''$.

Suppose

(1ai) $M \vDash Tally_b(z_1)$.

$\Longrightarrow$ by $(4.5)$, $M \vDash Tally_b(t_2z_1)$,

$\Rightarrow$ by (4.24), $M \vDash t_2 z_1 = t''$,

$\Rightarrow M \vDash t_2 < t''$,

$\Rightarrow M \vDash t' = t_1 = t_2 < t'' = t'$, contradicting $M \vDash t' \in J \subseteq I_0$.

Suppose

  (1aii)  $M \vDash \neg Tally_b(z_1)$.

$\Rightarrow M \vDash a \subseteq_p z_1$,

$\Rightarrow$ since $M \vDash Tally_b(t'')$, $M \vDash \neg(z_1 = a)$ & $\neg(aEz_1)$.

Assume, for a reductio, that $M \vDash aBz_1$. Then $M \vDash \exists z_2\ az_2 = z_1$.

$\Rightarrow M \vDash t_1(av'a)t_2 az_2 = t'avat''$,

$\Rightarrow$ by (4.16), $M \vDash t_2 \subseteq_p v$,

$\Rightarrow M \vDash t' = t_1 = t_2 \subseteq_p v$, contradicting $M \vDash Pref(ava, t')$.

Hence $M \vDash \neg(aBz_1)$.

Likewise $M \vDash \neg \exists z_2, z_3\ z_1 = z_2 a z_3$.

Therefore, $M \vDash \neg(a \subseteq_p z_1)$.

  (1b)  $M \vDash t_1 < t_2$ & $(t_1 av'at_2 a)Bx$.

$\Rightarrow M \vDash (t_1 av'at_2 a)Bz$,

$\Rightarrow M \vDash \exists z_2\ (t_1 av'at_2 a)z_2 = z = t'avat''$,

$\Rightarrow$ by (4.16), $M \vDash t_2 \subseteq_p v$,

$\Rightarrow M \vDash t' = t_1 < t_2 \subseteq_p v$, contradicting $M \vDash Pref(ava, t')$.

(2) $M \vDash t' < t''$ & $(t'avat''a)Bz$.

$\Rightarrow$ from $M \vDash x = t_1 av'at_2\ v\ (t_1 av'at_2 a)Bx$, $M \vDash (t_1 av'at_2)Bz$ & $(t'avat'')Bz$,

$\Rightarrow$ since $M \vDash t_1 = t'$, by (3.8),

$M \vDash (t_1av'at_2)B(t_1avat'') \lor t_1av'at_2=t_1avat'' \lor (t_1avat'')B(t_1av'at_2)$.

(2a) $M \vDash (t_1av'at_2)B(t_1avat'')$.

$\Rightarrow M \vDash \exists y_1 (t_1av'at_2)y_1 = t_1avat''$.

Suppose

(2ai) $M \vDash Tally_b(y_1)$.

$\Rightarrow$ by (4.5), $M \vDash Tally_b(t_2y_1)$,

$\Rightarrow$ from $M \vDash t_1(av'a)t_2y_1 = t_1avat''$, by (4.24$^b$), $M \vDash t_2y_1 = t''$,

$\Rightarrow M \vDash t_1av'at'' = t_1avat''$,

$\Rightarrow$ by (3.6) and (3.7), $M \vDash v = v'$,

$\Rightarrow M \vDash Firstf(x,t',ava,t_2)$, as required.

Suppose

(2aii) $M \vDash \neg Tally_b(y_1)$.

Exactly analogous to (1aii) with $y_1$ in place of $z_1$.

(2b) $M \vDash t_1av'at_2 = t_1avat''$.

$\Rightarrow$ by (4.24$^b$), $M \vDash t_2 = t''$, $\Rightarrow M \vDash t_1av'at_2 = t_1avat_2$,

$\Rightarrow$ by (3.6) and (3.7), $M \vDash v = v'$,

$\Rightarrow M \vDash Firstf(x,t',ava,t_2)$, as required.

(2c) $M \vDash (t_1avat'')B(t_1av'at_2)$.

Exactly analogous to (2a).

This completes the proof of (5.4).

(5.5) For any string concept $I \subseteq I_0$ there is a string concept $J \subseteq I$ such that

$QT^+ \vdash \forall x \in J \ \forall t,u,t_1,t_2 (x=tauat \ \& \ MaxT_b(t,x) \rightarrow$

$\rightarrow \forall v,t_1,t_2 \ (Fr(x,t_1,aua,t_2) \rightarrow t_1=t=t_2))$.

Let $J \equiv I_{4.14b} \ \& \ I_{4.23} \ \& \ I_{4.24b}$.

Assume $M \vDash x=tauat \ \& \ MaxT_b(t,x)$ where $M \vDash J(x)$.

$\Rightarrow \ M \vDash Tally_b(t)$.

Let $M \vDash Fr(x,t_1,ava,t_2)$.

Suppose, for a reductio, that $M \vDash \exists w_1 \ Intf(x,w_1,t_1,ava,t_2)$.

$\Rightarrow \ M \vDash \exists w_2 \ w_1at_1avat_2aw_2=x=tauat \ \& \ t_1<t_2 \ \& \ Max^+T_b(t_1,w_1)$.

$\Rightarrow$ from $M \vDash MaxT_b(t,x)$, $M \vDash t_1<t_2 \leq t$,

$\Rightarrow$ by (4.14), $M \vDash w_1=t \ v \ tBw_1$,

$\Rightarrow \ M \vDash t \subseteq_p w_1$,

$\Rightarrow$ from $M \vDash Max^+T_b(t_1,w_1)$, $M \vDash t<t_1$,

$\Rightarrow \ M \vDash t<t_1<t$, contradicting $M \vDash t \in J \subseteq I_0$.

Therefore, $M \vDash \neg \exists w_1 \ Intf(x,w_1,t_1,ava,t_2)$.

$\Rightarrow \ M \vDash Firstf(x,t_1,ava,t_2) \ v \ Lastf(x,t_1,ava,t_2)$.

Suppose, for a reductio, that

$M \vDash \exists w \ (wat_1avat_2=x=tauat \ \& \ t_1=t_2 \ \& \ Max^+T_b(t_1,w))$.

$\Rightarrow$ by (4.14$^b$), $M \vDash w=t \ v \ tBw$,

$\Rightarrow \ M \vDash t \subseteq_p w$,

$\Rightarrow$ from $M \vDash Max^+T_b(t,x)$, $M \vDash t<t_1<t$, a contradiction again.

Also, suppose that $M \vDash (t_1 a v a t_2 a)Bx$ & $t_1 < t_2$.

$\implies M \vDash \exists w\ t_1 a v a t_2 a w = x = t a u a t$,

$\implies$ by (4.23$^b$), since $M \vDash Tally_b(t_1)$ & $Tally_b(t)$, $M \vDash t_1 = t$,

$\implies$ from $M \vDash Max^+T_b(t,x)$ & $t_2 \subseteq_p x$, $M \vDash t_2 < t$,

$\implies M \vDash t < t_2 < t$, contradicting $M \vDash t \in J \subseteq I_0$.

Hence from $M \vDash Firstf(x,t_1,ava,t_2) \lor Lastf(x,t_1,ava,t_2)$ we have that

$\qquad M \vDash x = t_1 a v a t_2$.

$\implies M \vDash t_1 a v a t_2 = x = t a u a t$ & $Tally_b(t_1)$ & $Tally_b(t_2)$,

$\implies$ by (4.23$^b$) and (4.24$^b$), $M \vDash t_1 = t = t_2$, as required.

This completes the proof of (5.5).

(5.6) For any string concept $I \subseteq I_0$ there is a string concept $J \subseteq I$ such that

$QT^+ \vdash \forall z \in J \ \forall x,t,v,t_1,t_2,t',w,t''$ (Env(t,x) & Fr(x,t_1,ava,t_2) & z=xt'wt'' &

& aBw & aEw & Tally$_b$(t') & Tally$_b$(t'') $\rightarrow \ \exists t_3$ Fr(z,t_1,ava,t_3)).

Let $J \equiv I_{4.5}$.

Assume $M \vDash$ Env(t,x) & Fr(x,t_1,ava,t_2) & z=xt'wt''

where $M \vDash$ aBw & aEw & Tally$_b$(t') & Tally$_b$(t'') and $M \vDash$ J(x).

We distinguish three cases:

(1) $M \vDash$ Firstf(x,t_1,ava,t_2).

$\Rightarrow M \vDash$ Pref(ava,t_1) & Tally$_b$(t_2) &

& ((t_1=t_2 & x=t_1avat_2) v ((t_1<t_2 & (t_1avat_2a)Bx)).

If $M \vDash$ t_1=t_2 & x=t_1avat_2, then $M \vDash$ z=(t_1avat_2)t'wt''.

$\Rightarrow$ since $M \vDash$ aBw, $M \vDash$ (t_1ava(t_2t')a)Bz,

$\Rightarrow$ from $M \vDash$ Tally$_b$(t_2) & Tally$_b$(t'), by (4.5), $M \vDash$ Tally$_b$(t_2t'),

$\Rightarrow M \vDash$ Pref(ava,t_1) & Tally$_b$(t_2t') & (t_1ava(t_2t')a)Bz,

$\Rightarrow M \vDash$ Firstf(z,t_1,ava,t_2t').

If $M \vDash$ t_1<t_2 & (t_1avat_2a)Bx, then $M \vDash$ (t_1avat_2a)Bz, whence

$M \vDash$ Firstf(z,t_1,ava,t_2).

(2) $M \vDash \exists w_1$Intf(x,w_1,t_1,ava,t_2).

$\Rightarrow M \vDash \exists w_2$ x=w_1at_1avat_2aw_2 & Pref(ava,t_1) & Tally$_b$(t_2) & t_1<t_2 &

& Max$^+$T$_b$(t_1,w_1),

$\Rightarrow M \vDash$ z=(w_1at_1avat_2aw_2)t'wt'',

$\Rightarrow$ M ⊨ $\exists w_3$ z=$w_1$a$t_1$avat$_2$a$w_3$,

$\Rightarrow$ M ⊨ Intf(z,$w_1$,$t_1$,ava,$t_2$).

(3)  M ⊨ Lastf(x,$t_1$,ava,$t_2$).

$\Rightarrow$ from M ⊨ Env(t,x), M ⊨ $t_1$=$t_2$=t & Pref(ava,$t_1$) & Tally$_b$(t) &

$\qquad\qquad\qquad$ & $\exists w_1$(x=$w_1$atavat & Max$^+$T$_b$(t,$w_1$)),

$\Rightarrow$ M ⊨ z=($w_1$atavat)t'wt'',

$\Rightarrow$ M ⊨ $\exists w_2$ z=$w_1$atavatt'a$w_2$,

$\Rightarrow$ since M ⊨ Tally$_b$(t) & Tally$_b$(t'), by (4.5), M ⊨ Tally$_b$(tt') & t<tt',

$\Rightarrow$ M ⊨ Pref(ava,$t_1$) & Tally$_b$(tt') & t<tt' &

$\qquad\qquad\qquad$ & $\exists w_2$ z=$w_1$atavatt'a$w_2$ & Max$^+$T$_b$(t,$w_1$),

$\Rightarrow$ M ⊨ Intf(z,$w_1$,t,ava,tt').

This completes the proof of (5.6).

(5.7) For any string concept $I \subseteq I_0$ there is a string concept $J \subseteq I$ such that

$QT^+ \vdash \forall x \in J \; \forall t,t',u \; (Firstf(x,t,aua,t') \; \& \; Lastf(x,t,aua,t) \rightarrow x=tauat \; \& \; t=t')$.

Let $J \equiv I_{4.16}$.

Assume $M \vDash Firstf(x,t,aua,t') \; \& \; Lastf(x,t,aua,t)$ where $M \vDash J(x)$.

$\Longrightarrow$ $M \vDash Pref(aua,t) \; \& \; Tally_b(t) \; \& \; (x=tauat \; \lor \; \exists w(x=watauat \; \& \; Max^+T_b(t,w)))$.

Suppose $M \vDash \exists w(x=watauat \; \& \; Max^+T_b(t,w))$.

$\Longrightarrow$ $M \vDash Firstf(x,t,aua,t'), \; M \vDash (ta)Bx,$

$\Longrightarrow$ $M \vDash \exists x_1 \; tax_1 = x = watauat,$

$\Longrightarrow$ by $(4.14^b)$, $M \vDash w=t \; \lor \; tBw,$

$\Longrightarrow$ $M \vDash t \subseteq_p w,$

$\Longrightarrow$ from $M \vDash Max^+T_b(t,w)$, $M \vDash t<t$, contradicting $M \vDash t \in I \subseteq I_0$.

Therefore, $M \vDash x=tauat$.

$\Longrightarrow$ from $M \vDash Firstf(x,t,aua,t')$,

$\quad\quad\quad M \vDash Tally_b(t') \; \& \; ((t=t' \; \& \; x=tauat') \; \lor \; (t<t' \; \& \; (tauat'a)Bx))$.

Suppose $M \vDash t<t' \; \& \; (tauat'a)Bx$.

$\Longrightarrow$ $M \vDash \exists x_2 \; tauat'ax_2 = x = tauat,$

$\Longrightarrow$ by $(4.16)$, $M \vDash t' \subseteq_p u,$

$\Longrightarrow$ from $M \vDash Pref(aua,t)$, $M \vDash Max^+T_b(t,u),$

$\Longrightarrow$ $M \vDash t'<t,$

$\Longrightarrow$ $M \vDash t'<t<t'$, contradicting $M \vDash t' \in I \subseteq I_0$.

Therefore, $M \vDash t=t' \; \& \; tauat' = x = tauat$.

This completes the proof of (5.7).

(5.8) For any string concept $I \subseteq I_0$ there is a string concept $J \subseteq I$ such that

$QT^+ \vdash \forall x \in J \; \forall t,t_1,t_2,u,w \; (\text{Firstf}(x,t_1,aua,t_2) \; \& \; \text{Lastf}(x,t,awa,t') \; \& \; t_1 < t_2 \leq t \;\rightarrow$

$\rightarrow \exists w_1(x = w_1 ataw'at \; \& \; \text{Max}^+T_b(t,w_1) \; \& \; ((t_1au)Bw_1 \vee t_1au = w_1)))$.

Let $J \equiv I_{3.8} \; \& \; I_{4.17b} \; \& \; I_{4.23b}$.

Assume $M \vDash \text{Firstf}(x,t_1,aua,t_2) \; \& \; \text{Lastf}(x,t,awa,t')$ where $M \vDash t_1 < t_2 \leq t \; \& \; J(x)$.

$\Longrightarrow M \vDash \text{Pref}(aua,t_1) \; \& \; \text{Tally}_b(t_2) \; \&$

$\& \; ((t_1 = t_2 \; \& \; x = t_1auat_2) \vee (t_1 < t_2 \; \& \; (t_1auat_2a)Bx))$,

and $M \vDash \text{Pref}(awa,t') \; \& \; \text{Tally}_b(t) \; \&$

$\& \; (x = tawat \vee \exists w_1(x = w_1 atawat \; \& \; \text{Max}^+T_b(t,w_1))$.

We first claim that $M \vDash x \neq tawat$.

Suppose otherwise. Then $M \vDash (t_1a)Bx \; \& \; (ta)Bx \; \& \; \text{Tally}_b(t_1) \; \& \; \text{Tally}_b(t)$.

$\Longrightarrow M \vDash \exists x_1,x_2 \; t_1ax_1 = x = tax_2$,

$\Longrightarrow$ by (4.23$^b$), $M \vDash t_1 = t$,

$\Longrightarrow M \vDash t_1 < t = t_1$, contradicting $M \vDash t_1 \in I \subseteq I_0$.

Therefore $M \vDash x \neq tawat$.

$\Longrightarrow M \vDash \exists w_1(x = w_1 atawat \; \& \; \text{Max}^+T_b(t,w_1))$,

$\Longrightarrow M \vDash (t_1au)Bx \; \& \; w_1Bx$,

$\Longrightarrow$ by (3.8), $M \vDash (t_1au)Bw_1 \vee t_1au = w_1 \vee w_1B(t_1au)$.

Suppose, for a reductio, that $M \vDash w_1B(t_1au)$.

$\Longrightarrow M \vDash \exists w_2 \; t_1au = w_1w_2$,

$\Longrightarrow$ from $M \vDash t_1 \in I \subseteq I_0 \; \& \; t_1 < t_2$, $M \vDash \neg(t_1 = t_2)$,

$\Rightarrow M \vDash t_1 < t_2$ & $(t_1auat_2a)Bx$,

$\Rightarrow M \vDash \exists x_1 (w_1w_2)at_2ax_1 = x = w_1atawat$,

$\Rightarrow$ by (3.7), $M \vDash w_2at_2ax_1 = atawat$,

$\Rightarrow M \vDash w_2 = a \lor aBw_2$.

If $M \vDash w_2 = a$, then $M \vDash aat_2ax_1 = atawat$, whence $M \vDash at_2ax_1 = tawat$, which contradicts $M \vDash \text{Tally}_b(t)$.

$\Rightarrow M \vDash aBw_2$,

$\Rightarrow M \vDash \exists w_3(aw_3 = w_2$ & $aw_3at_2ax_1 = atawat)$,

$\Rightarrow M \vDash w_3at_2ax_1 = tawat$,

$\Rightarrow M \vDash w_3 = t \lor tBw_3$,

$\Rightarrow M \vDash t \subseteq_p w_3$,

$\Rightarrow M \vDash t \subseteq_p w_3 \subseteq_p w_2 \subseteq_p t_1au$,

$\Rightarrow$ by (4.17$^b$), $M \vDash t \subseteq_p t_1 \lor t \subseteq_p u$.

If $M \vDash t \subseteq_p u \subseteq_p aua$, then from $M \vDash \text{Pref}(aua,t_1)$ we have $M \vDash \text{Max}^+T_b(t_1,aua)$, so $M \vDash t < t_1$.

Therefore, $M \vDash t \leq t_1$.

$\Rightarrow M \vDash t \leq t_1 < t_2 \leq t$, contradicting $M \vDash t \in I \subseteq I_0$.

Hence $M \vDash (t_1au)Bw_1 \lor t_1au = w_1$, as required.

This completes the proof of (5.8).

(5.9) For any string concept $I \subseteq I_0$ there is a string concept $J \subseteq I$ such that

$QT^+ \vdash \forall x \in J \ \forall t_1, t_2, t', t'', u, w \ (Env(t_2, x) \ \& \ x = t_1 w t_2 \ \& \ aBw \ \& \ aEw \ \& \ Pref(aua, t') \ \&$

$\& \ t' < t'' \ \& \ Tally_b(t'') \ \& \ Max^+T_b(t', w') \rightarrow x \neq w't'auat''at_2 \ \& \ x \neq w'at'auat''at_2)$.

Let $J \equiv I_{3.6} \ \& I_{4.17b} \ \& \ I_{4.23b}$.

Assume $M \vDash Env(t_2, x) \ \& \ x = t_1 w t_2$

where $M \vDash aBw \ \& \ aEw \ \& \ Pref(aua, t') \ \& \ t' < t'' \ \& \ Tally_b(t'') \ \& \ Max^+T_b(t', w')$ and $M \vDash J(x)$.

Assume for a reductio that $M \vDash x = w'at'auat''at_2 \ \vee \ x = w'at'auat''at_2$.

$\Rightarrow$ from $M \vDash aBw$, $M \vDash \exists w_1 \ aw_1 = w$,

$\Rightarrow$ from $M \vDash aEw$, $M \vDash \exists w_2 \ w = w_2 a$,

$\Rightarrow$ from $M \vDash Env(t_2, x)$, $M \vDash \exists v \ Lastf(x, t_2, ava, t_2)$,

$\Rightarrow M \vDash Pref(ava, t_2) \ \& \ (x = t_2 avat_2 \ \vee \ \exists w(x = wat_2 avat_2 \ \& \ Max^+T_b(t_2, w)))$.

There are two cases:

(1) $M \vDash x = t_2 avat_2$.

$\Rightarrow M \vDash t_2 avat_2 = t_1 aw_1 t_2$,

$\Rightarrow$ since $M \vDash Tally_b(t_1) \ \& \ Tally_b(t_2)$, by (4.23$^b$), $M \vDash t_1 = t_2$.

Suppose, for a reductio, that

(1a) $M \vDash x = w'at'auat''at_2$.

$\Rightarrow M \vDash t_2 avat_2 = x = w'at'auat''at_2$,

$\Rightarrow$ by (4.14$^b$), $M \vDash w' = t_2 \ \vee \ t_2 Bw'$,

$\Rightarrow M \vDash t_2 \subseteq_p w'$,

$\Rightarrow$ from $M \vDash Max^+T_b(t',w')$, $M \vDash t_2<t'$,

$\Rightarrow$ from $M \vDash Env(t_2,x)$, $M \vDash MaxT_b(t_2,x)$,

$\Rightarrow M \vDash t' \leq t_2$,

$\Rightarrow M \vDash t_2 < t' \leq t_2$, contradicting $M \vDash t_2 \in I \subseteq I_0$.

 (1b)  $M \vDash x=w'at'auat''aat_2$.

The same argument applies as in (1a).

(2)  $M \vDash \exists w(x=wat_2avat_2 \,\&\, Max^+T_b(t_2,w))$.

Again, suppose that

 (2a)  $M \vDash x=w'at'auat''at_2$.

$\Rightarrow M \vDash wat_2avat_2=x=w'at'auat''at_2$,

$\Rightarrow$ by (3.6), $M \vDash wat_2av=w'at'auat''$,

$\Rightarrow M \vDash v=t'' \lor t''Ev$,

$\Rightarrow M \vDash t'' \subseteq_p v$,

$\Rightarrow$ from $M \vDash Pref(ava,t_2)$, $M \vDash t''<t_2$,

$\Rightarrow$ from $M \vDash wat_2av=x=w'at'auat''$, $M \vDash t_2 \subseteq_p w'at'auat''$,

$\Rightarrow$ repeatedly by (4.17[b]), $M \vDash t_2 \subseteq_p w' \lor t_2 \subseteq_p t' \lor t_2 \subseteq_p u \lor t_2 \subseteq_p t''$.

If $M \vDash t_2 \subseteq_p w'$, then from $M \vDash Max^+T_b(t',w')$, $M \vDash t_2<t'$; if $M \vDash t_2 \subseteq_p u$, then from $M \vDash Pref(aua,t')$, again $M \vDash t_2<t'$.

$\Rightarrow$ from $M \vDash t_2 \subseteq_p w'at'auat''$, $M \vDash t_2 \leq t'<t'' \leq t_2$, contradicting $M \vDash t_2 \in I \subseteq I_0$.

 (2b)  $M \vDash x=w'at'auat''aat_2$.

We proceed analogously to (2a) and derive $M \vDash wat_2av=w'at'auat''a$,

whence $M \vDash v=a \lor aEv$. Now, $M \vDash v \neq a$, because otherwise

M ⊨ wat₂aa=w'at'auat''a  and  M ⊨ wat₂a=w'at'auat'',  a contradiction.

$\Rightarrow$ M ⊨ ∃v₁(v=v₁a & wat₂a(v₁a)=w'at'auat''a),

$\Rightarrow$ M ⊨ wat₂av₁=w'at'auat''.

From this point on the same argument applies as in (2a) with v₁ in place of v.

This concludes the proof of (5.9).

(5.10) For any string concept $I \subseteq I_0$ there is a string concept $J \subseteq I$ such that

$QT^+ \vdash \forall z \in J \; \forall t (Env(t,z) \rightarrow \forall t',t'' (Tally_b(t') \; \& \; Tally_b(t'') \rightarrow z \neq t'at'' \; \& \; z \neq t'aat''))$.

Let $J \equiv I_{4.16} \; \& \; \& \; I_{4.23b} \; \& \; I_{4.24b}$.

Assume $M \vDash Env(t,z)$ where $M \vDash J(z)$.

Let $M \vDash Tally_b(t') \; \& \; Tally_b(t'')$, and assume for a reductio that

$$M \vDash z=t'at'' \; v \; z=t'aat''.$$

$\Rightarrow$ from $M \vDash Env(t,z)$, $M \vDash \exists u \; Lastf(z,t,aua,t)$,

$\Rightarrow M \vDash Pref(aua,t) \; \& \; (z=tauat \; v \; \exists w(z=watauat \; \& \; Max^+T_b(t,w)))$.

There are two cases:

(1) $M \vDash z=tauat$.

$\Rightarrow M \vDash tauat=z=t'at'' \; v \; tauat=z=t'aat''$,

$\Rightarrow$ by (4.23$^b$) and (4.24$^b$), $M \vDash t'=t=t''$,

$\Rightarrow M \vDash tauat=tat \; v \; tauat=taat$,

$\Rightarrow$ by (3.6) and (3.7), $M \vDash aua=a \; v \; au=a$, a contradiction.

(2) $M \vDash \exists w \; z=watauat$.

$\Rightarrow M \vDash watauat=z=t'at'' \; v \; watauat=z=t'aat''$,

$\Rightarrow$ by (4.16), $M \vDash t \subseteq_p a \; v \; t \subseteq_p aa$, a contradiction because $M \vDash Tally_b(t)$.

This completes the proof of (5.10).

(5.11) For any string concept $I \subseteq I_0$ there is a string concept $J \subseteq I$ such that

$QT^+ \vdash \forall z \in J\ \forall t(Env(t,z) \to \exists w, t'(Tally_b(t') \ \&\ z=t'wt'\ \&\ aBw\ \&\ aEw))$.

Let $J \equiv I_{4.14b}\ \&\ I_{4.15b}\ \&\ I_{5.10}$.

Assume $M \vDash Env(t,z)$ where $M \vDash J(z)$.

$\Rightarrow M \vDash \exists t'(Tally_b(t')\ \&\ (t'a)Bz))\ \&\ Tally_b(t)\ \&\ (at)Ez$,

$\Rightarrow M \vDash \exists z_1, z_2\ (t'a)z_2 = z = z_1(at)$,

$\Rightarrow$ by (4.14$^b$), $M \vDash z_1 = t'\ v\ t'Bz_1$,

$\Rightarrow$ by (4.15$^b$), $M \vDash z_2 = t\ v\ tEz_2$.

Suppose, for a reductio, that $M \vDash z_1 = t'\ v\ z_2 = t$. Then $M \vDash z = t'at$, which contradicts (5.10).

Therefore, $M \vDash t'Bz_1\ v\ tEz_2$.

$\Rightarrow M \vDash \exists z_3, z_4\ (t'z_3 = z_1\ v\ z_4 t = z_2)$,

$\Rightarrow M \vDash (t'z_3)at = z\ v\ z = (t'a)z_4 t$,

$\Rightarrow M \vDash t'z_3 at = z = t'az_2\ v\ t'az_4 t = z = z_1 at$,

$\Rightarrow$ by (3.6) and (3.7), $M \vDash z_3 at = az_2\ v\ t'az_4 = z_1 a$,

$\Rightarrow M \vDash (z_3 = a\ v\ aBz_3)\ v\ (z_4 = a\ v\ aEz_4)$,

$\Rightarrow M \vDash (z = t'(z_3 a)t\ \&\ aB(z_3 a)\ \&\ aE(z_3 a))\ v\ (z = t'(az_4)t\ \&\ aB(az_4)\ \&\ aE(az_4))$,

as required.

This completes the proof of (5.11).

(5.12) $QT^+ \vdash Env(t,x)$ & $t=b \rightarrow \exists u(x=bauab$ & $Pref(aua,b))$.

Assume $M \vDash Env(t,x)$ & $t=b$.

$\Rightarrow$ $M \vDash MaxT_b(b,x)$,

$\Rightarrow$ $M \vDash \exists u,t',t''\ Firstf(x,t',aua,t'')$,

$\Rightarrow$ $M \vDash Tally_b(t')$ & $Tally_b(t'')$,

$\Rightarrow$ from $M \vDash MaxT_b(b,x)$, $M \vDash t' \leq b$ & $t'' \leq b$,

$\Rightarrow$ from (3) and the definition of $\leq$, $M \vDash t'=b=t''$,

$\Rightarrow$ from $M \vDash Firstf(x,t',aua,t'')$, $M \vDash x=t'auat'' \lor t'<t''$,

$\Rightarrow$ from $M \vDash t'=b=t''$ & $\neg(b<b)$, $M \vDash \neg(t'<t'')$,

$\Rightarrow$ $M \vDash x=t'auat''=bauab$,

$\Rightarrow$ from $M \vDash Firstf(x,b,aua,b)$, $M \vDash Pref(aua,b)$, as required.

This completes the proof of (5.12).

(5.13) For any string concept $I\subseteq I_0$ there is a string concept $J\subseteq I$ such that

$QT^+ \vdash \forall z\in J\ \forall x,t_0,t,t',u(Env(t_0,x)\ \&\ z=xt'auat\ \&\ Tally_b(t')\ \&\ Lastf(z,t,aua,t) \rightarrow$

$\rightarrow t'<t\ \&\ \exists t^*,w^*(\ Tally_b(t^*)\ \&\ aBw^*\ \&\ aEw^*\ \&\ z=t^*w^*tauat))$.

Let $J \equiv I_{3.6}\ \&\ I_{4.5}\ \&\ I_{4.10}\ \&\ I_{5.11}$.

Assume $M \vDash Env(t_0,x)\ \&\ z=xt'auat\ \&\ Tally_b(t')$ where

$M \vDash Lastf(z,t,aua,t)\ \&\ J(z)J$.

$\Rightarrow M \vDash z=tauat\ \lor\ \exists w_1(z=w_1atauat\ \&\ Max^+T_b(t,w_1))$.

Suppose, for a reductio, that $M \vDash z=tauat$.

$\Rightarrow M \vDash xt'auat=z=tauat$,

$\Rightarrow$ by (3.6), $M \vDash xt'=t$.

But from $M \vDash Env(t_0,x)$, $M \vDash (at_0)Ex\ \&\ Tally_b(t_0)$.

$\Rightarrow M \vDash \exists x_1\ x=x_1at_0$,

$\Rightarrow M \vDash a\subseteq_p x\subseteq_p t$, which contradicts $M \vDash Tally_b(t)$.

Therefore, $M \vDash z\neq tauat$.

$\Rightarrow M \vDash \exists w_1\ z=w_1atauat$,

$\Rightarrow M \vDash (x_1at_0)t'auat=z=w_1atauat$,

$\Rightarrow$ by (3.6), $M \vDash x_1at_0t'=w_1at$,

$\Rightarrow$ from $M \vDash Lastf(z,t,aua,t)$, $M \vDash Tally_b(t)$,

$\Rightarrow$ from $M \vDash Tally_b(t_0)\ \&\ Tally_b(t')$, by (4.5), $M \vDash Tally_b(t_0t')$,

$\Rightarrow$ by (4.24$^b$), $M \vDash t_0t'=t$,

$\Rightarrow$ by (4.10), $M \vDash t't_0=t_0t'=t$, $M \vDash t'<t$.

We also have from  M ⊨ Env($t_0$,x),  by (5.11),

$$M \vDash \exists t^*,w^*(\text{Tally}_b(t^*) \ \& \ aBw^* \ \& \ aEw^* \ \& \ z=t^*w^*t_0).$$

Hence  M ⊨ z=xt'auat=t*w*$t_0$t'auat=t*w*tauat,  as required.

This completes the proof of (5.13).

(5.14) For any string concept I⊆I₀ there is a string concept J⊆I such that

$$QT^+ \vdash \forall x \in J\ \forall t,t',t'',v,v',w' \subseteq_p x\ (MaxT_b(t,x)\ \&\ Lastf(x,t,v',t)\ \&\ wat'vt''=x\ \&$$
$$\&\ Tally_b(t'')\ \&\ Pref(v,t')\ \&\ Max^+T_b(t',w)\ \to\ t'=t=t''\ \&\ v=v').$$

Let J(x) abbreviate

$I_{3.6}^{\subseteq p}(x)\ \&\ I_{3.7}^{\subseteq p}(x)\ \&\ I_{3.10}(x)\ \&\ I_{4.5}^{\subseteq p}(x)\ \&\ I_{4.10}^{\subseteq p}(x)\ \&\ I_{4.14b}(x)\ \&\ I_{4.16}^{\subseteq p}(x)\ \&\ I_{4.23b}^{\subseteq p}(x)\ \&$

$\&\ I_{4.24b}(x).$

Assume  M ⊨ MaxT_b(t,x) & Lastf(x,t,v',t) & wat'vt''=x

where  M ⊨ Tally_b(t'') & Pref(v,t') & Max⁺T_b(t',w)  and  M ⊨ J(x).

From  M ⊨ Lastf(x,t,v',t)  we have

  M ⊨ Pref(v',t) & (x=tv't ∨ ∃w'(w'atv't=x & Max⁺T_b(t,w'))).

Suppose  M ⊨ x=tv't.

⇒ M ⊨ wat'vt''=x=tv't,

⇒ by (4.14ᵇ),  M ⊨ w=t ∨ tBw,

⇒ M ⊨ t⊆_p w,

⇒ from M ⊨ Max⁺T_b(t',w),  M ⊨ t<t'.

But from M ⊨ MaxT_b(t,x),  M ⊨ t'≤t,

⇒ M ⊨ t<t'≤t,  contradicting  M ⊨ t ∈ I ⊆ I₀.

Therefore  M ⊨ ¬(x=tv't).

So we may assume that  M ⊨ ∃w'(w'atv't=x=wat'vt'').

From  M ⊨ Pref(v',t) & Pref(v,t') we have  M ⊨ ∃v'' v'=av''a & ∃v₀ v=av₀a,

⇒ M ⊨ w'at(av''a)t=w'atv't=x=wat'vt''=wat'(av₀a)t'',

⇒ by (4.24$^b$), M ⊨ t=t'',

⇒ M ⊨ w'atv't=wat'vt.

From M ⊨ MaxT$_b$(t,x), M ⊨ t'≤t, and we also have that

$$M ⊨ (tv't)Ex\ \&\ (t'vt)Ex,$$

⇒ by (3.10), M ⊨ (tv't)E(t'vt) v tv't=t'vt v (t'vt)E(tv't).

Case 1. M ⊨ (tv't)E(t'vt).

⇒ M ⊨ ∃x$_1$ x$_1$tv't=t'vt.

Now, M ⊨ Tally$_b$(x$_1$) v ¬Tally$_b$(x$_1$).

(1a) M ⊨ Tally$_b$(x$_1$).

⇒ by (4.5), M ⊨ Tally$_b$(x$_1$t), and M ⊨ x$_1$t(av''a)t=t'vt,

⇒ by (4.23$^b$), M ⊨ x$_1$t=t',

⇒ by (4.10), M ⊨ tx$_1$=t',

⇒ M ⊨ t<t'≤t, contradicting M ⊨ t ∈ I ⊆ I$_0$.

(1b) M ⊨ ¬Tally$_b$(x$_1$).

⇒ M ⊨ a⊆$_p$x$_1$.

We distinguish cases based on the definition of ⊆$_p$:

We have that M ⊨ ¬(x$_1$=a) & ¬aBx$_1$ because M ⊨ Tally$_b$(t'). Suppose M ⊨ aEx$_1$.

⇒ M ⊨ ∃x$_2$ x$_2$a=x$_1$,

⇒ M ⊨ (x$_2$a)tv't=t'vt,

⇒ M ⊨ x$_2$at(av''a)t=t'vt,

⇒ by (4.16), M ⊨ t⊆$_p$v,

⇒ M ⊨ t'≤t⊆$_p$v, contradicting M ⊨ Pref(v,t').

The same argument yields a contradiction if $M \vDash \exists x_2, x_3 \ x_2 a x_3 = x_1$.

This completes Case 1.

<u>Case 2.</u> $M \vDash tv't = t'vt$.

$\Rightarrow$ by (4.23$^b$), $M \vDash t = t'$,

$\Rightarrow M \vDash tv't = tvt$,

$\Rightarrow$ by (3.6) and (3.7), $M \vDash v' = v$.

Hence $M \vDash t' = t = t''$ & $v' = v$, as required.

<u>Case 3.</u> $M \vDash (t'vt)E(tv't)$.

$\Rightarrow M \vDash \exists x_1 \ x_1 t'vt = tv't$.

Again, $M \vDash \text{Tally}_b(x_1) \lor \neg \text{Tally}_b(x_1)$.

(3a) $M \vDash \text{Tally}_b(x_1)$.

$\Rightarrow M \vDash w'a(x_1 t'vt) = w'atv't = x = wat'vt$,

$\Rightarrow$ by (3.6), $M \vDash w'ax_1 = wa$, a contradiction.

(3b) $M \vDash \neg \text{Tally}_b(x_1)$.

$\Rightarrow M \vDash a \subseteq_p x_1$.

Again, $M \vDash \neg(x_1 = a)$ & $\neg a B x_1$. Suppose that $M \vDash a E x_1$.

$\Rightarrow M \vDash \exists x_2 \ (x_2 a) t'vt = tv't$, $\Rightarrow M \vDash x_2 at'(av_0 a)t = tv't$,

$\Rightarrow$ by (4.16), $M \vDash t' \subseteq_p v'$,

$\Rightarrow$ from $M \vDash \text{Pref}(v, t')$, $M \vDash t' < t$.

But also from $M \vDash x_2 at'vt = tv't$, we have by (4.14$^b$), that

$$M \vDash x_2 = t \lor t B x_2,$$

$\Rightarrow M \vDash t \subseteq_p x_2$.

Now, we have $M \vDash x_2at'vt=tv't$ and $M \vDash w'atv't=x=wat'vt$,

$\Rightarrow$ $M \vDash w'a(x_2at'vt)=watv't$,

$\Rightarrow$ by (3.6), $M \vDash w'ax_2=w$,

$\Rightarrow M \vDash t\subseteq_p x_2 \subseteq_p w$,

$\Rightarrow M \vDash \text{Max}^+T_b(t',w)$, $M \vDash t<t'$,

$\Rightarrow M \vDash t<t'<t$, contradicting $M \vDash t \in I \subseteq I_0$.

An analogous argument derives a contradiction if $M \vDash \exists x_2, x_3 \; x_2ax_3=x_1$.

This completes Case 3 and the proof of (5.14).

(5.15)  For any string concept $I \subseteq I_0$ there is a string concept $J \subseteq I$ such that

$QT^+ \vdash \forall z \in J \; \forall t, t_1, t_2, t_3, t_4, u_1, u_2 \subseteq_p z \; ((Firstf(z, t_1, u_1, t_2) \; \& \; Firstf(z, t_3, u_2, t_4)) \; v$

$v \; ((Lastf(z, t_1, u_1, t_2) \; \& \; Lastf(z, t_3, u_2, t_4))) \rightarrow u_1 = u_2).$

Let $J(x)$ abbreviate

$I_{3.6}(x) \; \& \; I_{3.7}^{\subseteq p}(x) \; \& \; I_{3.10}(x) \; \& \; I_{4.5}^{\subseteq p}(x) \; \& \; I_{4.14b}^{\subseteq p}(x) \; \& \; I_{4.16}^{\subseteq p}(x) \; \& \; I_{4.23b}^{\subseteq p}(x) \; \& \; I_{4.24b}(x).$

Assume first that  $M \vDash Lastf(z, t_1, u_1, t_2) \; \& \; Lastf(z, t_3, u_2, t_4)$

where  $M \vDash MaxT_b(t, z)$  and  $M \vDash J(z)$.

$\Rightarrow M \vDash Pref(u_1, t_1) \; \& \; t_1 = t_2 \; \& \; (z = t_1 u_1 t_2 \; v \; \exists w_1 \; (z = w_1 a t_1 u_1 t_2 \; \& \; Max^+T_b(t_1, w_1)))$

and

$\Rightarrow M \vDash Pref(u_2, t_3) \; \& \; t_3 = t_4 \; \& \; (z = t_3 u_2 t_4 \; v \; \exists w_2 \; (z = w_2 a t_3 u_2 t_4 \; \& \; Max^+T_b(t_3, w_2))).$

From $M \vDash Pref(u_1, t_1) \; \& \; Pref(u_2, t_3)$  we have

$M \vDash Tally_b(t_1) \; \& \; \exists u' \; u_1 = au'a \; \& \; Tally_b(t_3) \; \& \; \exists u'' \; u_2 = au''a.$

<u>Case 1</u>.  $M \vDash z = t_1 u_1 t_2 \; \& \; z = t_3 u_2 t_4.$

$\Rightarrow M \vDash t_1(au'a)t_1 = z = t_3(au''a)t_3,$

$\Rightarrow$ by (4.24$^b$),  $M \vDash t_1 = t_3,$

$\Rightarrow M \vDash t_1(au'a)t_1 = t_1(au''a)t_1,$

$\Rightarrow$ by (3.6) and (3.7),  $M \vDash u_1 = au'a = au''a = u_2.$

<u>Case 2</u>.  $M \vDash z = t_1 u_1 t_2 \; \& \; \exists w_2 \; (z = w_2 a t_3 u_2 t_4 \; \& \; Max^+T_b(t_3, w_2)).$

$\Rightarrow M \vDash t_1(au'a)t_1 = z = w_2 a t_3(au''a)t_3,$

$\Rightarrow$ by (4.24$^b$),  $M \vDash t_1 = t_3,$

$\Rightarrow$ by (4.14$^b$),  $M \vDash w_2 = t_1 \; v \; t_1 B w_2,$

$\Rightarrow$ M ⊨ $t_1 \subseteq_p w_2$,

$\Rightarrow$ M ⊨ $t_3 \subseteq_p w_2$, contradicting M ⊨ $Max^+T_b(t_3,w_2)$.

Case 3.  M ⊨ $\exists w_1 (z=w_1at_1u_1t_2 \& Max^+T_b(t_1,w_1)) \& z=t_3u_2t_4$.

Exactly analogous to Case 2.

Case 4.  M ⊨ $\exists w_1 (z=w_1at_1u_1t_2 \& Max^+T_b(t_1,w_1)) \&$

$$\& \exists w_2 (z=w_2at_3u_2t_4 \& Max^+T_b(t_3,w_2)).$$

$\Rightarrow$ M ⊨ $w_1at_1(au'a)t_1 = w_2at_3(au''a)t_3$,

$\Rightarrow$ by (4.24$^b$), M ⊨ $t_1 = t_3$,

$\Rightarrow$ M ⊨ $w_1at_1(au'a)t_1 = z = w_2at_1(au''a)t_1$,

$\Rightarrow$ M ⊨ $(t_1u_1t_1)Ez \& (t_1u_2t_1)Ez$,

$\Rightarrow$ by (3.10), M ⊨ $(t_1u_1t_1)E(t_1u_2t_1) \lor t_1u_1t_1=t_1u_2t_1 \lor (t_1u_2t_1)E(t_1u_1t_1)$.

We derive a contradiction from M ⊨ $(t_1u_1t_1)E(t_1u_2t_1)$ and M ⊨ $(t_1u_2t_1)E(t_1u_1t_1)$ exactly as in Case 1 of (5.14), and we derive M ⊨ $u_1=u_2$ from M ⊨ $t_1u_1t_1=t_1u_2t_1$ by reasoning exactly as in Case 2 there.

This completes the proof that M ⊨ $Lastf(z,t_1,u_1,t_2) \& Lastf(z,t_3,u_2,t_4) \rightarrow u_1=u_2$ and the proof of (5.15).

(5.16) For any string concept $I \subseteq I_0$ there is a string concept $J \subseteq I$ such that

$QT^+ \vdash \forall x \in J \ \forall z,t_0,t_1,t_2,t,v(Env(t,z) \ \& \ xBz \ \& \ Env(t_0,x) \ \& \ Firstf(x,t_1,ava,t_2) \rightarrow$

$\rightarrow \exists t_3 Firstf(z,t_1,ava,t_3))$.

Let $J \equiv I_{5.4} \ \& \ I_{5.15}$.

Assume $M \vDash Env(t,z) \ \& \ xBz \ \& \ Env(t_0,x)$

along with $M \vDash Firstf(x,t_1,ava,t_2) \ \& \ J(x)$.

$\implies$ from $M \vDash Env(t,z)$, $M \vDash \exists v',t',t_3 \ Firstf(z,t',av'a,t_3)$,

$\implies$ from $M \vDash Env(t_0,x) \ \& \ xBz$, by (5.4), $M \vDash \exists t_4 \ Firstf(x,t',av'a,t_4)$,

$\implies$ by (5.15), $M \vDash v'=v \ \& \ t'=t_1$,

$\implies M \vDash \exists t_3 \ Firstf(z,t_1,ava,t_3)$, as required.

This completes the proof of (5.16).

(5.17)   $QT^+ \vdash \forall z \, [Set(z) \to (z=aa \lor \exists y \, y \, \varepsilon \, z) \,\&\, \neg(z=aa \,\&\, \exists y \, y \, \varepsilon \, z)]$.

Suppose $M \vDash Set(z)$.

$\implies$  $M \vDash z=aa \lor \exists t \subseteq_p z \, Env(t,z)$.

If $M \vDash z=aa$, we immediately have $M \vDash z=aa \lor \exists y \, y \, \varepsilon \, z$.

So assume $M \vDash \exists t \subseteq_p z \, Env(t,z)$.

$\implies$  $M \vDash \exists u \subseteq_p z \, \exists t_1, t_2 \leq t \, Firstf(z, t_1, aua, t_2)$,

$\implies$  $M \vDash Pref(aua, t_1)$,

$\implies$  $M \vDash \exists t \subseteq_p z \, \exists u \subseteq_p z \, \exists t_1, t_2 \leq t \, Fr(z, t_1, aua, t_2)$,

$\implies$  $M \vDash \exists y \, y \, \varepsilon \, z$, as required.

Suppose, for a reductio, that $M \vDash z=aa \,\&\, \exists y \, y \, \varepsilon \, z$.

$\implies$  $M \vDash \exists t \subseteq_p z \, \exists u \subseteq_p z \, \exists t_1, t_2 \leq t \, Fr(z, t_1, aua, t_2)$,

$\implies$  $M \vDash Pref(aua, t_1)$,

$\implies$  $M \vDash Tally_b(t_1) \,\&\, t_1 \subseteq_p t \subseteq_p z$,

$\implies$  $M \vDash b \subseteq_p aa$, whence a contradiction follows.

Therefore, also $M \vDash \neg(z=aa \,\&\, \exists y \, y \, \varepsilon \, z)$.

This completes the proof of (5.17).

We thus have:

(5.18)   $QT^+ \vdash \forall z \, [Set(z) \to (\exists y \, y \, \varepsilon \, z \leftrightarrow z \neq aa)]$.

(5.19) For any string concept $I \subseteq I_0$ there is a string concept $J \subseteq I$ such that

$QT^+ \vdash \forall x \in J \ \forall u, t_1, t_2, t_3, t_4 (Set(x) \ \& \ (Firstf(x,t_1,aua,t_2) \lor Lastf(x,t_1,aua,t_2)) \rightarrow$

$\rightarrow \neg \exists w \ Intf(x,w,t_3,aua,t_4))$.

Let $J \equiv I_{3.12} \ \& \ I_{4.16}$.

Assume $M \vDash Set(x)$ and let $M \vDash \exists w \ Intf(x,w,t_3,aua,t_4)$ where $M \vDash J(x)$.

$\Rightarrow M \vDash \exists t \ Env(t,x)$.

Assume $M \vDash Firstf(x,t_1,aua,t_2)$.

$\Rightarrow M \vDash Fr(x,t_1,aua,t_2) \ \& \ Fr(x,t_3,aua,t_4)$,

$\Rightarrow$ by (d) of $M \vDash Env(t,x)$, $M \vDash t_1 = t_3$,

$\Rightarrow$ from $M \vDash Intf(x,w,t_3,aua,t_4)$, $M \vDash \exists w' \ wat_3auat_4aw' = x$,

$\Rightarrow$ from $M \vDash Firstf(x,t_1,aua,t_2)$,

$M \vDash (t_1 = t_2 \ \& \ x = t_1auat_2) \lor (t_1 < t_2 \ \& \ (t_1auat_2a)Bx)$.

<u>Case 1.</u> $M \vDash t_1 = t_2 \ \& \ x = t_1auat_2$.

$\Rightarrow M \vDash wat_3auat_4aw' = x = t_1auat_2$,

$\Rightarrow$ by (4.16), $M \vDash t_3aua \subseteq_p (t_3auat_4) \subseteq_p aua$, contradicting (3.12).

<u>Case 2.</u> $M \vDash t_1 < t_2 \ \& \ (t_1auat_2a)Bx$.

$\Rightarrow M \vDash \exists x_1 \ t_1auat_2ax_1 = x = wat_3auat_4aw'$,

$\Rightarrow$ by (4.14$^b$), $M \vDash w = t_1 \lor t_1Bw$,

$\Rightarrow M \vDash t_1 \subseteq_p w$,

$\Rightarrow$ from $M \vDash Max^+T_b(t_3,w)$, $M \vDash t_1 < t_3$,

$\Rightarrow M \vDash t_1 < t_3 = t_1$, which contradicts $M \vDash t_1 \in J \subseteq I_0$.

Therefore, $M \vDash \neg Firstf(x,t_1,aua,t_2)$.

Assume $M \vDash Lastf(x,t_1,aua,t_2)$.

$\Rightarrow M \vDash t_1=t_2=t$,

$\Rightarrow M \vDash Fr(x,t_1,aua,t_2) \& Fr(x,t_3,aua,t_4)$,

$\Rightarrow$ by (d) of $M \vDash Env(t,x)$, $M \vDash t_1=t_3$,

$\Rightarrow M \vDash t_3=t$,

$\Rightarrow$ from $M \vDash Lastf(x,t,aua,t)$, $M \vDash x=tauat \;\lor\; \exists w_1\; w_1atauat=x$.

If $M \vDash x=tauat$, we reason as in Case 1 above.

If $M \vDash w_1atauat=x=wat_3auat_4aw'$, then by (a) of $M \vDash Env(t,x)$,

$M \vDash MaxT_b(t,x)$.

$\Rightarrow M \vDash t_3<t_4\leq t$,

$\Rightarrow M \vDash t=t_3<t_4\leq t$, contradicting $M \vDash t\in J\subseteq I_0$.

Therefore, also $M \vDash \neg Lastf(x,t_1,aua,t_2)$.

This completes the proof of (5.19).

(5.20) For any string concept $I \subseteq I_0$ there is a string concept $J \subseteq I$ such that

$QT^+ \vdash \forall x \in J\ \forall t',v,t'',t_1,u,t_2\ (Fr(x,t'ava,t'') \& Firstf(x,t_1,aua,t_2) \rightarrow$

$\rightarrow (Firstf(x,t'ava,t'') \& t_1=t') \lor t_1<t')$.

Let $J \equiv I_{4.14b}\ \&\ I_{4.23b}$.

Assume $M \models Fr(x,t'ava,t'') \& Firstf(x, t_1,aua,t_2)$ where $M \models J(x)$.

$\Rightarrow M \models Pref(aua,t_1) \& ((t_1=t_2 \& x=t_1auat_2) \lor (t_1<t_2 \& (t_1auat_2a)Bx))$.

We distinguish three cases:

Case 1. $M \models Firstf(x,t'ava,t'')$.

$\Rightarrow M \models Pref(ava,t') \& ((t'=t'' \& x=t'avat'') \lor (t'<t'' \& (t'avat''a)Bx))$,

$\Rightarrow M \models (t_1auab)Bx \& (t'avab)Bx \& Tally_b(t_1) \& Tally_b(t')$,

$\Rightarrow$ by $(4.23^b)$, $M \models t_1=t'$, as required.

Case 2. $M \models \exists w\ Intf(x,w,t'ava,t'')$.

$\Rightarrow M \models Pref(ava,t') \& \exists w_1,w_2\ (x=w_1at'avat''aw_2 \& Max^+T_b(t',w_1))$,

$\Rightarrow M \models Tally_b(t_1) \& t_1B(w_1at'avat''aw_2)$,

$\Rightarrow$ by $(4.14^b)$, $M \models w_1=t_1 \lor t_1Bw_1$,

$\Rightarrow M \models t_1\subseteq_p w_1, \Rightarrow M \models t_1<t'$.

Case 3. $M \models Lastf(x,t'ava,t'')$.

$\Rightarrow M \models Pref(ava,t') \& \exists w\ (x=wat'avat'' \& Max^+T_b(t',w))$.

Then, just as in Case 2, $M \models t_1<t'$.

This completes the proof of (5.20).

(5.21) For any string concept $I\subseteq I_0$ there is a string concept $J\subseteq I$ such that

$$QT^+ \vdash \forall x \in J\ \forall u,t_1,t_2\ (Set(x)\ \&\ Firstf(x,t_1,aua,t_2)\ \&\ x=t_1auat_2\ \rightarrow$$

$$\rightarrow\ \forall w\ (w\ \varepsilon\ x \leftrightarrow w=u)).$$

Let $J \equiv I_{5.5}\ \&\ I_{5.15}\ \&\ I_{5.19}\ \&\ I_{5.20}$.

Assume $M \vDash Set(x)\ \&\ Firstf(x,t_1,aua,t_2)\ \&\ x=t_1auat_2$ where $M \vDash J(x)$.

$\Rightarrow\ M \vDash \exists t\ Env(t,x)$,

$\Rightarrow\ M \vDash MaxT_b(t,x)$,

$\Rightarrow$ from $M \vDash Fr(x,t_1,aua,t_2)$, by (5.5), $M \vDash t_1=t=t_2$,

$\Rightarrow$ from $M \vDash Firstf(x,t_1,aua,t_2)$, $M \vDash Pref(aua,t_1)$,

$\Rightarrow\ M \vDash Lastf(x,t_1,aua,t_2)$.

Suppose $M \vDash w=u$.

$\Rightarrow$ from $M \vDash Firstf(x,t_1,aua,t_2)$, $M \vDash w\ \varepsilon\ x$.

Suppose $M \vDash w\ \varepsilon\ x\ \&\ w \neq u$.

$\Rightarrow\ M \vDash \exists t_3,t_4\ Fr(x,t_3,awa,t_4)$.

If $M \vDash Firstf(x,t_3,awa,t_4)\ v\ Lastf(x,t_3,awa,t_4)$, then, by (5.15), $M \vDash w=u$, contradicting the hypothesis.

So we may assume that $M \vDash \exists w_1\ Intf(x,w_1,t_3,awa,t_4)$.

$\Rightarrow$ by (5.19), $M \vDash \neg\ Firstf(x,t_3,awa,t_4)$,

$\Rightarrow$ by (5.20), $M \vDash t=t_1<t_3$,

$\Rightarrow$ from $M \vDash MaxT_b(t,x)\ \&\ t_3\subseteq_p x$, $M \vDash t_3 \leq t$,

$\Rightarrow\ M \vDash t=t_1<t_3\leq t$, contradicting $M \vDash t\in J\subseteq I_0$.

Therefore, $M \vDash w = u$, as required.

This completes the proof of (5.21).

(5.22) For any string concept $I \subseteq I_0$ there is a string concept $J \subseteq I$ such that

$$QT^+ \vdash \forall x \in J\ \forall u,t\ (Env(t,x)\ \&\ \forall z\ (z\ \varepsilon\ x \leftrightarrow z=u) \leftrightarrow$$

$$\leftrightarrow x=tauat\ \&\ Firstf(x,t,aua,t)\ \&\ Lastf(x,t,aua,t)).$$

Let $J \equiv I_{4.22}\ \&\ I_{5.21}$.

Assume $M \vDash Env(t,x)\ \&\ \forall z\ (z\ \varepsilon\ x \leftrightarrow z=u)$ where $M \vDash J(x)$.

$\Rightarrow M \vDash u\ \varepsilon\ x$,

$\Rightarrow$ from $M \vDash Env(t,x)$, $M \vDash \exists w,t_1,t_2\ Firstf(x,t_1,awa,t_2)\ \&\ \exists v Lastf(x,t,ava,t)$,

$\Rightarrow M \vDash w=u=v$,

$\Rightarrow M \vDash Lastf(x,t,aua,t)$.

On the other hand, we also have $M \vDash Fr(x,t_1,aua,t_2)\ \&\ Fr(x,t,aua,t)$.

$\Rightarrow$ by (d) of $M \vDash Env(t,x)$, $M \vDash t=t_1$.

From $M \vDash Firstf(x,t_1,aua,t_2)$, $M \vDash (t_1=t_2\ \&\ x=t_1auat_2)\ v\ (t_1<t_2\ \&\ (t_1auat_2a)Bx)$,

$\Rightarrow$ from $M \vDash Env(t,x)$, $M \vDash MaxT_b(t,x)$,

$\Rightarrow M \vDash t_2 \leq t_1 = t$.

Since we may assume that $J \subseteq I_0$ is closed under $\subseteq_p$, it follows that

$$M \vDash \neg(t_1<t_2).$$

$\Rightarrow M \vDash t_1=t=t_2\ \&\ x=tauat$,

$\Rightarrow M \vDash Firstf(x,t,aua,t)$.

Conversely, assume $M \vDash x=tauat\ \&\ Firstf(x,t,aua,t)\ \&\ Lastf(x,t,aua,t)$.

First, we show that $M \vDash Env(t,x)$.

From $M \vDash Firstf(x,t,aua,t)$, we have $M \vDash Pref(aua,t)$.

$\implies$ M ⊨ Max$^+$T$_b$(t,aua),

$\implies$ by (4.22), M ⊨ MaxT$_b$(t,x).

This gives part (a) of M ⊨ Env(t,x). Parts (b) and (c) follow immediately from the hypothesis. Parts (d) and (e) follow immediately from M ⊨ x=tauat and (a) by (5.5). Hence M ⊨ Env(t,x).

$\implies$ by (5.21), M ⊨ $\forall$z (z ε x $\leftrightarrow$ z=u).

This completes the proof of (5.22).

(5.23)  For any string concept $I\subseteq I_0$ there is a string concept $J\subseteq I$ such that

$QT^+ \vdash \forall z\in J\ \forall t,u,z'(Set(z)\ \&\ z=tauat\ \&\ Pref(aua,t)\ \&\ Set(z')\ \&\ z\sim z'\ \&$

    $\&\ \forall w,t_1,t_2\ (Firstf(z,t_1,awa,t_2) \to \exists t_3\ Firstf(z',t_1,awa,t_3)) \to z=z')$.

Let $J \equiv I_{5.15}\ \&\ I_{5.22}$.

Assume  $M \vDash Set(z)\ \&\ z=tauat\ \&\ Pref(aua,t)\ \&\ Set(z')\ \&\ z\sim z'\ \&$

    $\&\ \forall w,t_1,t_2\ (Firstf(z,t_1,awa,t_2) \to \exists t_3\ Firstf(z',t_1,awa,t_3))$

where $M \vDash J(z)$.

$\Rightarrow$ from $M \vDash Pref(aua,t)\ \&\ z=tauat$, $M \vDash Firstf(z,t,aua,t)$,

$\Rightarrow$ from $M \vDash Set(z)$, by (5.21), $M \vDash \forall w(w\ \varepsilon\ z \leftrightarrow w=u)$,

$\Rightarrow$ from $M \vDash Set(z')\ \&\ z\sim z'$, $M \vDash \forall w(w\ \varepsilon\ z' \leftrightarrow w=u)$,

$\Rightarrow$ by (5.22), $M \vDash \exists t',u'(z'=t'au'at'\ \&\ Firstf(z',t',au'a,t')\ \&\ Lastf(z',t',au'a,t'))$,

$\Rightarrow$ by hypothesis, from $M \vDash Firstf(z,t,aua,t)$, $M \vDash \exists t_3\ Firstf(z',t,aua,t_3)$,

$\Rightarrow$ by (5.15), $M \vDash t=t'\ \&\ u=u'$,

$\Rightarrow$ $M \vDash z=tauat=t'au'at'=z'$, as claimed.

This completes the proof of (5.23).

(5.24) For any string concept I⊆I₀ there is a string concept J⊆I such that

QT⁺ ⊢ ∀y∈J ∀y',y'',t,t₁,t₂,u,v(Env(t,y) & Lastf(y,t,aua,t) & Fr(y,t₁,ava,t₂) & u≠v &

& y'By & Env(t₁,y') & Env(t'',y'') & Lastf(y'',t,aua,t) → ¬ y''By).

Let J ≡ I₄.₁₈ & I₅.₄ & I₅.₂₂.

Assume  M ⊨ Env(t,y) & y'By & Env(t₁,y') & Env(t'',y'')  where

M ⊨ Lastf(y,t,aua,t) & Fr(y,t₁,ava,t₂) & u≠v & Lastf(y'',t,aua,t)  and  M ⊨ J(y).

Assume, for a reductio, that  M ⊨ y''By.

⟹ M ⊨ ∃z y''z=y,

⟹ from  M ⊨ Lastf(y,t,aua,t),

   M ⊨ Pref(aua,t) & (y=tauat ∨ ∃w₁(y=w₁atauat & Max⁺T_b(t,w₁))).

If  M ⊨ y=tauat,  then from  M ⊨ Pref(aua,t)  we also have  M ⊨ Firstf(y,t,aua,t),

whence  M ⊨ ∀w (w ε y ↔ w=u)  by (5.22).  But then  M ⊨ v ε y  follows from hypothesis  M ⊨ Fr(y,t₁,ava,t₂),  whence  M ⊨ v=u,  contradicting hypothesis.

Therefore  M ⊨ ∃w₁(y=w₁atauat & Max⁺T_b(t,w₁)).

⟹ from  M ⊨ Lastf(y'',t,aua,t),

   M ⊨ Pref(aua,t) & (y''=tauat ∨ ∃w₂(y''=w₂atauat & Max⁺T_b(t,w₂))).

(1)  M ⊨ y''=tauat.

⟹ M ⊨ Firstf(y'',t,aua,t),

⟹ from  M ⊨ Env(t,y),  M ⊨ ∃v',t',t''Firstf(y,t',av'a,t''),

⟹ from  M ⊨ Env(t₁,y') & y'By, by (5.4),  M ⊨ ∃t₃Firstf(y',t',av'a,t₃),

⟹ from  M ⊨ Env(t₁,y'),  M ⊨ MaxT_b(t₁,y'),

$\Rightarrow$ M ⊨ t'≤t$_1$,

$\Rightarrow$ from  M ⊨ Env(t'',y'') & y''By, by (5.4),  M ⊨ ∃t$_4$Firstf(y'',t',av'a,t$_4$),

$\Rightarrow$ by (5.15),  M ⊨ t=t',

$\Rightarrow$ from  M ⊨ Env(t,y) & u≠v,  M ⊨ t≠t$_1$,

$\Rightarrow$ from  M ⊨ Env(t,y),  M ⊨ MaxT$_b$(t,y),

$\Rightarrow$  M ⊨ t$_1$<t,

$\Rightarrow$ M ⊨ t=t'≤t$_1$<t, contradicting  M ⊨ t∈I⊆I$_0$.

(2)  M ⊨ ∃w$_2$(y''=w$_2$atauat & Max$^+$T$_b$(t,w$_2$)).

$\Rightarrow$ M ⊨ (w$_2$atauat)z=y=w$_1$atauat.

Assume, for a reductio, that  M ⊨ Tally$_b$(z).

$\Rightarrow$ by (4.24$^b$),  M ⊨ tz=t,

$\Rightarrow$ M ⊨ tBt, contradicting  M ⊨ t∈I⊆I$_0$.

Therefore,  M ⊨ ¬Tally$_b$(z).

$\Rightarrow$ M ⊨ a⊆$_p$z,

$\Rightarrow$ from  M ⊨ Tally$_b$(t),  M ⊨ z≠a & ¬(aEz),

$\Rightarrow$ M ⊨ ∃z$_1$ az$_1$=z  v  ∃z$_1$,z$_2$ z=z$_2$az$_1$,

$\Rightarrow$ M ⊨ w$_2$at(au)at(az$_1$)=w$_1$at(auat)  v  w$_2$at(au)at(z$_2$az$_1$)=w$_1$at(auat)

where  M ⊨ Max$^+$T$_b$(t,w$_1$) & Max$^+$T$_b$(t,w$_2$).

But this contradicts (4.18).

This completes the proof of (5.24).

(5.25) For any string concept $I \subseteq I_0$ there is a string concept $J \subseteq I$ such that

$$QT^+ \vdash \forall x \in J\ \forall z,t',v,t''\ [Fr(z,t',ava,t'')\ \&\ \exists t_1,w\ (x=t_1wz\ \&\ aBw\ \&\ aEw\ \&\ Max^+T_b(t',t_1w)) \to Fr(x,t',ava,t'')].$$

Let $J \equiv I_{4.17b}$.

Assume $M \vDash Fr(z,t',ava,t'')$ and let $M \vDash x=t_1wz$

where $M \vDash aBw\ \&\ aEw\ \&\ Max^+T_b(t',t_1w)$ and $M \vDash J(x)$.

From the hypothesis we have $M \vDash Pref(ava,t')\ \&\ Tally_b(t'')$ along with one of the following four cases:

<u>Case 1.</u> $M \vDash t'=t''\ \&\ t'avat''=z$.

$\Rightarrow M \vDash t_1wz=x=t_1wt'avat''$,

$\Rightarrow$ since $M \vDash aEw$, $M \vDash (at'avat'')Ex$ and $M \vDash \exists w_1\ w_1a=w$,

$\Rightarrow M \vDash t_1w_1at'avat''=x\ \&\ Max^+T_b(t',t_1w)$,

$\Rightarrow M \vDash Lastf(x,t',ava,t'')$.

<u>Case 2.</u> $M \vDash t'=t''\ \&\ \exists w'(w'at'avat''=z\ \&\ Max^+T_b(t',w'))$.

$\Rightarrow M \vDash t_1wz=x=t_1w(w'at'avat'')$,

$\Rightarrow M \vDash (at'avat'')Ex$.

Let $M \vDash s \subseteq_p t_1ww'\ \&\ Tally_b(s)$.

$\Rightarrow$ since $M \vDash aEw$, $M \vDash s \subseteq_p t_1w_1aw'$ where $M \vDash w_1a=w$.

$\Rightarrow$ by (4.17$^b$), $M \vDash s \subseteq_p t_1\ \lor\ s \subseteq_p w'$,

$\Rightarrow$ from $M \vDash Max^+T_b(t',t_1w)\ \&\ Max^+T_b(t',w')$, $M \vDash s<t'$.

Therefore, $M \vDash Max^+T_b(t',t_1ww')$.

This suffices to show that $M \vDash \text{Lastf}(x,t',ava,t'')$.

Case 3. $M \vDash t'<t''$ & $(t'avat'')Bz$.

$\Rightarrow M \vDash \exists z_1 (t'avat''a)z_1=z$,

$\Rightarrow$ since $M \vDash w_1a=w$, $M \vDash t_1wz=x=t_1w_1at'avat''az_1$,

$\Rightarrow$ since $M \vDash \text{Max}^+T_b(t',t_1w)$,

  $M \vDash \text{Pref}(ava,t')$ & $\text{Tally}_b(t'')$ & $t'<t''$ & $\exists z_1(w'at'avat''az_1=x$ & $\text{Max}^+T_b(t',w'))$

where $M \vDash w'=t_1w_1$.

But then $M \vDash \text{Intf}(x,w't',ava,t'')$.

Case 4. $M \vDash t'<t''$ & $\exists w_2,w_3(w_2at'avat''aw_3=z$ & $\text{Max}^+T_b(t',w_2))$.

$\Rightarrow M \vDash t_1wz=x=(t_1w_1aw_2)at'avat''aw_3$.

Let $M \vDash s\subseteq_p t_1w_1aw_2$ & $\text{Tally}_b(s)$.

$\Rightarrow$ by $(4.17^b)$, $M \vDash s\subseteq_p t_1w_1 \vee s\subseteq_p w_2$.

If $M \vDash s\subseteq_p t_1w_1$, then $M \vDash s<t'$ from $M \vDash \text{Max}^+T_b(t',t_1w)$; if $M \vDash s\subseteq_p w_2$, then $M \vDash s<t'$ from $M \vDash \text{Max}^+T_b(t',w_2)$.

Therefore, $M \vDash \text{Max}^+T_b(t',t_1w_1aw_2)$.

We then have

 $M \vDash \text{Pref}(ava,t')$ & $\text{Tally}_b(t'')$ & $t'<t''$ & $\exists w_3(w'at'avat''aw_3=x$ & $\text{Max}^+T_b(t',w'))$

where $M \vDash w'=t_1w_1aw_2$.

It follows that $M \vDash \text{Intf}(x,w't',ava,t'')$.

This completes the proof of (5.25).

(5.26) For any string concept $I \subseteq I_0$ there is a string concept $J \subseteq I$ such that

$QT^+ \vdash \forall x \in J \, \forall t, t_1, t_2, t_3, t_4, v_0, w, z [Env(t,x) \,\&\, x = t_1 w z \,\&\, aBw \,\&\, aEw \,\&\, x' = t_1 w t_2 \,\&$

$\&\, Max^+T_b(t_2, x') \,\&\, Firstf(z, t_3, av_0 a, t_4) \,\&\, t_2 < t_3 \,\&\, Set(z) \rightarrow Env(t,z)].$

Let $J \equiv I_{4.17b} \,\&\, I_{4.24b}$.

Assume $M \vDash Env(t,x) \,\&\, x = t_1 w z \,\&\, aBw \,\&\, aEw$

where $M \vDash x' = t_1 w t_2 \,\&\, Max^+T_b(t_2, x') \,\&\, Firstf(z, t_3, av_0 a, t_4) \,\&\, t_2 < t_3 \,\&\, Set(z)$

and $M \vDash J(x)$.

$\Rightarrow$ from $M \vDash Set(z) \,\&\, Firstf(z, t_3, av_0 a, t_4)$, $M \vDash v_0 \, \varepsilon \, z$,

$\Rightarrow$ by (5.18), $M \vDash z \neq aa$,

$\Rightarrow M \vDash \exists t' \, Env(t', z)$,

$\Rightarrow M \vDash \exists v' \, Lastf(z, t', av'a, t')$.

We claim that $M \vDash Max^+T_b(t', t_1 w)$.

Assume that $M \vDash Tally_b(t_0) \,\&\, t_0 \subseteq_p t_1 w$.

$\Rightarrow M \vDash t_0 \subseteq_p x'$,

$\Rightarrow$ from $M \vDash Max^+T_b(t_2, x')$, $M \vDash t_0 \leq t_2$,

$\Rightarrow$ from $M \vDash Env(t', z)$, $M \vDash Max^+T_b(t', z)$,

$\Rightarrow$ from $M \vDash Firstf(z, t_3, av_0 a, t_4)$, $M \vDash t_3 \leq t'$,

$\Rightarrow M \vDash t_0 \leq t_2 < t_3 \leq t'$,

$\Rightarrow M \vDash Max^+T_b(t', t_1 w)$, as claimed.

So we have

$M \vDash Lastf(z, t', av'a, t') \,\&\, x = t_1 w z \,\&\, aBw \,\&\, aEw \,\&\, Max^+T_b(t', t_1 w)$,

$\Rightarrow$ by proof of (5.25), parts (1) and (2),  $M \vDash \text{Lastf}(x,t',av'a,t')$,

$\Rightarrow$ from hypothesis  $M \vDash \text{Env}(t,x)$,  $M \vDash \exists v \text{Lastf}(x,t,ava,t)$,

$\Rightarrow$ $M \vDash (at')Ex \ \& \ (at)Ex$,

$\Rightarrow$ $M \vDash \exists x_1 \ x = x_1 at' \ \& \ \exists x_2 \ x = x_2 at$,

$\Rightarrow$ $M \vDash x_1 at' = x = x_2 at$,

$\Rightarrow$ since $M \vDash \text{Tally}_b(t') \ \& \ \text{Tally}_b(t)$, by (4.24$^b$),  $M \vDash t = t'$,

$\Rightarrow$ from  $M \vDash \text{Env}(t',z)$,  $M \vDash \text{Env}(t,z)$,  as required.

This completes the proof of (5.26).

(5.27) For any string concept I⊆I₀ there is a string concept J⊆I such that

$QT^+ \vdash \forall x \in J \; \forall u,v,w_1,w_3,t_1,t_2,t_3,t_4(Set(x) \; \& \; Intf(x,w_1,t_1,aua,t_2) \;\&$

$\& \; Intf(x,w_3,t_3,ava,t_4) \; \& \; u \neq v \; \& \; t_2 \neq t_3 \; \& \; t_4 \neq t_1 \;\rightarrow\; (t_1auat_2) \subseteq_p w_3 \;\vee\; (t_3avat_4) \subseteq_p w_1).$

Let $J \equiv I_{3.8} \; \& \; I_{4.14b} \; \& \; I_{4.15b} \; \& \; I_{4.24b}$.

Assume $M \vDash Set(x) \; \& \; Intf(x,w_1,t_1,aua,t_2) \; \& \; Intf(x,w_3,t_3,ava,t_4)$

where $M \vDash u \neq v \; \& \; t_2 \neq t_3 \; \& \; t_4 \neq t_1$ and $M \vDash J(x)$.

$\Longrightarrow M \vDash \exists t \; Env(t,x)$ and $M \vDash \exists w_2,w_4 \; w_1at_1auat_2aw_2 = x = w_3at_3avat_4aw_4$.

$\Longrightarrow M \vDash (w_1at_1)Bx \; \& \; (w_3at_3)Bx$,

$\Longrightarrow$ by (3.8), $M \vDash (w_1at_1)B(w_3at_3) \;\vee\; w_1at_1 = w_3at_3 \;\vee\; (w_3at_3)B(w_1at_1)$.

<u>Case 1.</u> $M \vDash w_1at_1 = w_3at_3$.

$\Longrightarrow$ since $M \vDash Tally_b(t_1) \; \& \; Tally_b(t_3)$, by (4.24$^b$), $M \vDash t_1 = t_3$.

But $M \vDash Fr(x,t_1,aua,t_2) \; \& \; Fr(x,t_3,ava,t_4)$. Hence from $M \vDash Env(t,x)$, $M \vDash u=v$,

which contradicts the hypothesis.

<u>Case 2.</u> $M \vDash (w_1at_1)B(w_3at_3)$.

$\Longrightarrow M \vDash \exists w_5 \; w_1at_1w_5 = w_3at_3$,

$\Longrightarrow M \vDash w_1at_1w_5(avat_4aw_4) = w_3at_3avat_4aw_4 = x = w_1at_1auat_2aw_2$,

$\Longrightarrow$ by (3.7), $M \vDash w_5avat_4aw_4 = auat_2aw_2$,

$\Longrightarrow M \vDash w_5 = a \;\vee\; aBw_5$.

We cannot have $M \vDash w_5 = a$ because $M \vDash Tally_b(t_3)$.

$\Longrightarrow M \vDash aBw_5$,

$\Longrightarrow M \vDash \exists w_6 \; aw_6 = w_5$,

$\Rightarrow$ $M \vDash w_1at_1(aw_6)=w_3at_3$,

$\Rightarrow$ by (4.15$^b$), $M \vDash w_6=t_3 \lor t_3Ew_6$.

(2a) $M \vDash w_6=t_3$.

$\Rightarrow$ $M \vDash w_1at_1at_3=w_3at_3$,

$\Rightarrow$ by (3.6), $M \vDash w_1at_1=w_3$,

$\Rightarrow$ $M \vDash t_1\subseteq_p w_3$,

$\Rightarrow$ from $M \vDash \text{Max}^+T_b(t_3,w_3)$, $M \vDash t_1<t_3$,

$\Rightarrow$ $M \vDash w_3at_3avat_4aw_4=w_3auat_2aw_2$,

$\Rightarrow$ by (3.7), $M \vDash t_3avat_4aw_4=uat_2aw_2$,

$\Rightarrow$ by (4.14$^b$), $M \vDash u=t_3 \lor t_3Bu$,

$\Rightarrow$ $M \vDash t_3\subseteq_p u$,

$\Rightarrow$ $M \vDash t_1<t_3\subseteq_p u$, contradicting $M \vDash \text{Pref}(aua,t_1)$.

(2b) $M \vDash t_3Ew_6$.

$\Rightarrow$ $M \vDash \exists w_7\; w_7t_3=w_6$,

$\Rightarrow$ $M \vDash w_1at_1a(w_7t_3)=w_3at_3$,

$\Rightarrow$ by (3.6), $M \vDash w_1at_1aw_7=w_3a$,

$\Rightarrow$ $M \vDash w_7=a \lor aEw_7$,

$\Rightarrow$ $M \vDash w_1at_1aa=w_3a \lor \exists w_8\; w_1at_1a(w_8a)=w_3a$,

$\Rightarrow$ $M \vDash w_1at_1a=w_3 \lor \exists w_8\; w_1at_1aw_8=w_3$,

$\Rightarrow$ $M \vDash t_1\subseteq_p w_3$,

$\Rightarrow$ since $M \vDash \text{Max}^+T_b(t_3,w_3)$, $M \vDash t_1<t_3$,

$\Rightarrow$ $M \vDash (w_1at_1aw_7)t_3avat_4aw_4=w_3at_3avat_4aw_4=x=w_1at_1auat_2aw_2$,

$\Rightarrow$ by (3.7), $M \vDash w_7t_3avat_4aw_4=uat_2aw_2$,

$\Rightarrow$ by (3.8), $M \vDash (w_7t_3)Bu \vee w_7t_3=u \vee uB(w_7t_3)$.

If $M \vDash (w_7t_3)Bu \vee w_7t_3=u$, then $M \vDash t_3\subseteq_p u$, whence $M \vDash t_1<t_3\subseteq_p u$,

contradicting $M \vDash \text{Pref}(aua,t_1)$.

Hence $M \vDash uB(w_7t_3)$.

$\Rightarrow M \vDash \exists w_8 \, uw_8=w_7t_3$,

$\Rightarrow M \vDash (uw_8)avat_4aw_4=uat_2aw_2$,

$\Rightarrow$ by (3.7), $M \vDash w_8avat_4aw_4=at_2aw_2$,

$\Rightarrow M \vDash w_8a \vee aBw_8$.

If $M \vDash w_8a$, then $M \vDash aavat_4aw_4=at_2aw_2$, whence $M \vDash avat_4aw_4=t_2aw_2$,

which is a contradiction because $M \vDash \text{Tally}_b(t_2)$.

Hence $M \vDash aBw_8$.

$\Rightarrow M \vDash \exists w_9 \, aw_9=w_8$,

$\Rightarrow M \vDash (aw_9)avat_4aw_4=at_2aw_2$,

$\Rightarrow M \vDash w_9avat_4aw_4=t_2aw_2$,

$\Rightarrow$ by (4.14$^b$), $M \vDash w_9=t_2 \vee t_2Bw_9$.

  (2bi)  $M \vDash w_9=t_2$.

$\Rightarrow M \vDash t_2avat_4aw_4=t_2aw_2$,

$\Rightarrow M \vDash w_1at_1auat_2avat_4aw_4=w_1at_1auat_2aw_2=x=w_3at_3avat_4aw_4$,

$\Rightarrow$ by (3.6), $M \vDash w_1at_1auat_2=w_3at_3$,

$\Rightarrow$ by (4.24$^b$), since $M \vDash \text{Tally}_b(t_2) \, \& \, \text{Tally}_b(t_3)$, $M \vDash t_2=t_3$, contradicting the hypothesis.

(2bii)   $M \vDash t_2Bw_9$.

$\Rightarrow$ $M \vDash \exists w_{10}\ t_2w_{10}=w_9$,

$\Rightarrow$ $M \vDash (t_2w_{10})avat_4aw_4=t_2aw_2$,

$\Rightarrow$ by (3.7), $M \vDash w_{10}avat_4aw_4=aw_2$,

$\Rightarrow$ $M \vDash w_{10}=a\ \vee\ aBw_{10}$,

$\Rightarrow$ $M \vDash w_1at_1auat_2(w_{10}avat_4aw_4)=w_1at_1auat_2aw_2=x=w_3at_3avat_4aw_4$,

$\Rightarrow$ by (3.6), $M \vDash w_1at_1auat_2w_{10}=w_3at_3$,

But we cannot have $M \vDash w_{10}=a$ because $M \vDash \text{Tally}_b(t_3)$.

$\Rightarrow$ $M \vDash aBw_{10}$,

$\Rightarrow$ $M \vDash \exists w_{11}\ aw_{11}=w_{10}$,

$\Rightarrow$ $M \vDash w_1at_1auat_2(aw_{11})=w_3at_3$,

$\Rightarrow$ by (4.24$^b$), $M \vDash w_{11}=t_3\ \vee\ t_3Ew_{11}$,

$\Rightarrow$ $M \vDash w_1at_1auat_2at_3=w_3at_3\ \vee\ \exists w_{12}\ w_1at_1auat_2a(w_{12}t_3)=w_3at_3$,

$\Rightarrow$ by (3.7), $M \vDash \vDash w_1at_1auat_2=w_3\ \vee\ w_1at_1auat_2aw_{12}=w_3a$,

$\Rightarrow$ $M \vDash w_{12}=a\ \vee\ aEw_{12}$,

$\Rightarrow$ $M \vDash w_1at_1auat_2=w_3\ \vee\ w_1at_1auat_2a=w_3\ \vee\ \exists w_{13}\ w_1at_1auat_2aw_{13}=w_3$,

$\Rightarrow$ $M \vDash (t_1auat_2)\subseteq_p w_3$, as required.

<u>Case 3.</u> $M \vDash (w_3at_3)B(w_1at_1)$.

By an argument entirely analogous to that of Case 2 we prove that

$$M \vDash (t_3avat_4)\subseteq_p w_1.$$

This completes the proof of (5.27).

(5.28) For any string concept $I \subseteq I_0$ there is a string concept $J \subseteq I$ such that

$QT^+ \vdash \forall x \in J \, \forall t',t'',t_1,t_2,u,v(Intf(x,w',t',ava,t'') \& Intf(x,w_1,t_1,aua,t_2) \& t''=t_1 \rightarrow$

$\rightarrow w'at'avat''=w_1at_1)$.

Let $J \equiv I_{3.8} \& I_{4.17b}$.

Assume $M \vDash Intf(x,w',t',ava,t'') \& Intf(x,w_1,t_1,aua,t_2)$ and $M \vDash t''=t_1$ where $M \vDash J(x)$.

$\Rightarrow M \vDash Pref(ava,t') \& t'<t''=t_1 \& Max^+T_b(t',w') \& Pref(aua,t_1) \&$

$\& Max^+T_b(t_1,w_1) \& \exists w'',w_2 \, w'at'avat''aw''=x=w_1at_1auat_2aw_2$,

$\Rightarrow M \vDash (w'at'avat'')Bx \& (w_1at_1)Bx$,

$\Rightarrow$ by (3.8),

$M \vDash (w'at'avat'')B(w_1at_1) \lor w'at'avat''=w_1at_1 \lor (w_1at_1)B(w'at'avat'')$.

We first show that

(1) $M \vDash \neg(w_1at_1)B(w'at'ava)$.

Suppose for a reductio that $M \vDash (w_1at_1)B(w'at'ava)$.

$\Rightarrow M \vDash t_1 \subseteq_p w_1at_1 \subseteq_p w'at'ava \subseteq_p w'at'avab$,

$\Rightarrow$ by (4.17$^b$), $M \vDash t_1 \subseteq_p w'at' \lor t_1 \subseteq_p vab$,

$\Rightarrow$ by (4.17$^b$), $M \vDash t_1 \subseteq_p w' \lor t_1 \subseteq_p t' \lor t_1 \subseteq_p v \lor t_1 \subseteq_p b$.

From $M \vDash Max^+T_b(t',w')$ we have that $M \vDash t_1 \subseteq_p w'$ implies $M \vDash t_1<t'$, which along with $M \vDash t'<t_1$ yields $M \vDash t_1<t_1$, contradicting $M \vDash t_1 \in I_0 \subseteq I$.

Therefore $M \vDash \neg(t_1 \subseteq_p w') \& \neg(t_1 \subseteq_p t')$.

If $M \vDash t_1 \subseteq_p v$, then from $M \vDash Pref(ava,t')$ we have $M \vDash Max^+T_b(t',ava)$, so $M \vDash t_1 \subseteq_p v$ implies $M \vDash t_1<t'$, again a contradiction.

Therefore also $M \vDash \neg(t_1 \subseteq_p v)$.

Hence $M \vDash t_1 \subseteq_p b$.

Since we may assume that $M \vDash b < t_1$, we have that $M \vDash \neg(w_1 a t_1) B(w'at'ava)$.

(1A) Suppose now, for a reductio, that $M \vDash (w_1 a t_1) B(w'at'avat'')$.

$\Rightarrow M \vDash \exists x_3\ w_1 a t_1 x_3 = w'at'avat''$,

$\Rightarrow$ by (4.15$^b$), $M \vDash t_1 x_3 = t''\ v\ t''E(t_1 x_3)$.

If $M \vDash t_1 x_3 = t''$, then $M \vDash t_1 < t''$, contradicting the hypothesis $M \vDash t'' = t_1$.

$\Rightarrow M \vDash t''E(t_1 x_3)$,

$\Rightarrow M \vDash \exists x_4\ t_1 x_3 = x_4 t''$.

We claim that $M \vDash aEx_4$.

For, we cannot have $M \vDash x_4 = a$ because $M \vDash \text{Tally}_b(t_1)$.

Suppose that $M \vDash x_4 = b\ v\ bEx_4$.

$\Rightarrow M \vDash t_1 x_3 = bt''\ v\ \exists x_5\ t_1 x_3 = (x_5 b)t''$,

$\Rightarrow M \vDash w_1 abt'' = w'at'avat''\ v\ w_1 a x_5 bt'' = w'at'avat''$,

$\Rightarrow$ by (3.6), $M \vDash w_1 ab = w'at'ava\ v\ w_1 a x_5 b = w'at'ava$,

a contradiction either way.

Therefore $M \vDash \neg x_4 = a\ \&\ \neg x_4 = b\ \&\ \neg bEx_4$, so by (QT5) we have that $M \vDash aEx_4$.

$\Rightarrow M \vDash \exists x_5\ x_5 a = x_4$,

$\Rightarrow M \vDash t_1 x_3 = x_5 a t''$,

$\Rightarrow$ by (4.14$^b$), $M \vDash x_5 = t_1\ v\ t_1 B x_5$,

$\Rightarrow M \vDash t_1 \subseteq_p x_5 \subseteq_p x_4$,

$\Rightarrow$ from $M \vDash t_1 x_3 = x_4 t''$, $M \vDash w_1 a x_4 t'' = w'at'avat''$,

$\Rightarrow$ by (3.6), $M \vDash w_1ax_4=w'at'ava$,

$\Rightarrow$ $M \vDash t_1\subseteq_p x_4 \subseteq_p w'at'ava$, whence a contradiction follows as in the above argument for (1).

Therefore $M \vDash \neg(w_1at_1)B(w'at'avat'')$.

We also have that

(2)   $M \vDash \neg(w'at'avat'')Bw_1$.

Suppose, for a reductio, that $M \vDash (w'at'avat'')Bw_1$.

$\Rightarrow$ $M \vDash t''\subseteq_p w_1$,

$\Rightarrow$ from $M \vDash Max^+T_b(t_1,w_1)$, $M \vDash t''<t_1$, contradicting the hypothesis $M \vDash t''=t_1$.

Therefore $M \vDash \neg(w'at'avat'')Bw_1$.

(2A) Suppose now, for a reductio, that $M \vDash (w'at'avat'')B(w_1at_1)$.

$\Rightarrow$ $M \vDash \exists x_3\ w_1at_1=w'at'avat''x_3$,

$\Rightarrow$ by (4.15$^b$), $M \vDash t''x_3=t_1 \vee t_1E(t''x_3)$.

We now argue, analogously to (1A) and show that $M \vDash \exists x_4\ t''x_3=x_4t_1$, that $M \vDash aEx_4$, and that $M \vDash t''\subseteq_p x_4$. Then from $M \vDash t''x_3=x_4t_1$ we have

$\qquad M \vDash w'at'avax_4t_1=w_1at_1$,

$\Rightarrow$ by (3.6), $M \vDash w'at'avax_4=w_1a$,

$\Rightarrow$ $M \vDash t_1\subseteq_p x_4\subseteq_p w_1$, and we derive a contradiction as in (2).

Therefore, $M \vDash \neg(w'at'avat'')B(w_1at_1)$.

It follows that $M \vDash w'at'avat''=w_1at_1$, as required.

This completes the proof of (5.28).

(5.29) For any string concept $I \subseteq I_0$ there is a string concept $J \subseteq I$ such that

$QT^+ \vdash \forall x \in J \; \forall t',u,v,w',w_1,t',t'' \subseteq_p x \; (MaxT_b(t,x) \; \& \; z=xt'utt' \; \& \; Intf(x,w',t',u,t'') \;$ &

& $Lastf(x,t,v,t) \; \& \; w_1atvt=x \; \& \; Max^+T_b(t_1,w_1) \to$

$\to \; w'at'ut''=w_1 \; v \; (w'at'ut'')Bw_1 \; v \; w'at'u=w_1a).$

Let $J(x) \equiv I_{3.6}^{\subseteq p}(x) \; \& \; I_{3.8}(x) \; \& \; I_{4.14b}^{\subseteq p}(x) \; \& \; I_{4.15b}^{\subseteq p}(x) \; \& \; I_{4.17b}^{\subseteq p}(x) \; \& \; I_{5.14b}^{\subseteq p}(x).$

Assume $M \vDash MaxT_b(t,x) \; \& \; Intf(x,w',t',u,t'') \; \& \; Lastf(x,t,v,t)$

where $M \vDash J(x)$ and $M \vDash w_1atvt=x \; \& \; Max^+T_b(t_1,w_1)$.

$\Rightarrow M \vDash Pref(u,t') \; \& \; \exists w'' \; w'at'ut''aw''=x=w_1atvt$ and $M \vDash Max^+T_b(t',w')$,

$\Rightarrow M \vDash \exists u_0 \; u=au_0a \; \& \; (w_1at)Bx \; \& \; (w'at'ut'')Bx$,

$\Rightarrow$ by (3.8), $M \vDash (w_1at)B(w'at'ut'') \; v \; w_1at=w'at'ut'' \; v \; (w'at'ut'')B(w_1at).$

<u>Case 1.</u> $M \vDash (w_1at)B(w'at'ut'')$.

$\Rightarrow M \vDash \exists x_1 \; w_1atx_1=w'at'ut''$,

$\Rightarrow$ by (4.15$^b$), $M \vDash tx_1=t'' \; v \; t''E(tx_1)$.

(1a) $M \vDash tx_1=t''$.

$\Rightarrow M \vDash t<t''$,

$\Rightarrow$ from $M \vDash MaxT_b(t,x)$, $M \vDash t<t'' \leq t$, contradicting $M \vDash t \in I \subseteq I_0$.

(1b) $M \vDash t''E(tx_1)$.

$\Rightarrow M \vDash \exists x_2 \; tx_1=x_2t''$.

We then show, by an argument exactly analogous to that in (1A) in (5.28), that

$M \vDash aEx_2$.

$\Rightarrow M \vDash \exists x_3 \; x_3a=x_2$,

$\Rightarrow$ M ⊨ $tx_1=x_3at''$,

$\Rightarrow$ by (4.14$^b$), M ⊨ $x_3=t$ v $tBx_3$,

$\Rightarrow$ M ⊨ $t\subseteq_p x_3 \subseteq_p x_2$.

From M ⊨ $tx_1=x_2t''$, we have

$$M \vDash w_1atx_1=w_1ax_2t''=w'at'ut'',$$

$\Rightarrow$ by (3.6), M ⊨ $w_1ax_2=w'at'u$,

$\Rightarrow$ M ⊨ $t\subseteq_p x_2 \subseteq_p w'at'u=w'at'(au_0a)$,

$\Rightarrow$ by (4.17$^b$), M ⊨ $t\subseteq_p w'$ v $t\subseteq_p t'au_0a$,

$\Rightarrow$ by (4.17$^b$), M ⊨ $t\subseteq_p w'$ v $t\subseteq_p t'$ v $t\subseteq_p u_0a\subseteq_p u$.

From M ⊨ $Max^+T_b(t',w')$ we have from M ⊨ $t\subseteq_p w'$ that M ⊨ $t<t'$. On the other hand, from M ⊨ $t\subseteq_p t'$ we have M ⊨ $t\leq t'$, and from M ⊨ Pref(u,t') we have M ⊨ $Max^+T_b(t',u)$. So from M ⊨ $t\subseteq_p u$ we also have M ⊨ $t<t'$. But from M ⊨ $MaxT_b(t,x)$ we then have M ⊨ $t\leq t'<t''\leq t$, which contradicts M ⊨ $t \in I \subseteq I_0$.

<u>Case 2.</u>  M ⊨ $w_1at=w'at'ut''$.

$\Rightarrow$ M ⊨ $w_1at=w'at'(au_0a)t''$,

$\Rightarrow$ by (4.24$^b$), M ⊨ $t=t''$,

$\Rightarrow$ M ⊨ $w_1at=w'at'ut$,

$\Rightarrow$ by (3.6), M ⊨ $w'at'u=w_1a$.

<u>Case 3.</u>  M ⊨ $(w'at'ut'')B(w_1at)$.

$\Rightarrow$ M ⊨ $\exists x_2\ w'at'ut''x_2=w_1at$,

$\Rightarrow$ by (4.15$^b$), M ⊨ $t''x_2=t$ v $tE(t''x_2)$.

If M ⊨ $t''x_2=t$, then M ⊨ $w'at'ut=w'at'ut''x_2=w_1at$, whence by (3.6),

$$M \vDash w'at'u = w_1 a.$$

If $M \vDash tE(t''x_2)$, then $M \vDash \exists x_4\ t''x_2 = x_4 t$.

We claim that $M \vDash aEx_4$.

Now, we have $M \vDash \neg(x_4 = a)$ because $M \vDash \text{Tally}_b(t'')$.

If $M \vDash x_4 = b\ v\ bEx_4$, then

$$M \vDash t''x_2 = bt\ v\ \exists x_5\ t''x_2 = (x_5)bt,$$

$\Rightarrow M \vDash w'at'ubt = w_1 at\ v\ w'at'ux_5 bt = w_1 at$,

$\Rightarrow$ by (3.6), $M \vDash w'at'ub = w_1 a\ v\ w'at'ux_5 b = w_1 a$, a contradiction either way.

Therefore $M \vDash aEx_4$.

$\Rightarrow M \vDash \exists x_6\ x_6 a = x_4$,

$\Rightarrow M \vDash t''x_2 = x_6 at$,

$\Rightarrow$ by (4.14$^b$), $M \vDash x_6 = t''\ v\ t''Bx_6$.

But $\quad M \vDash w_1 at = w'at'ut''x_2 = w'at'ux_4 t = w'at'u(x_6 a)t$,

$\Rightarrow M \vDash w_1 at = w'at'ut''at\ v\ \exists x_7\ w_1 at = w'at'u(t''x_7)at$,

$\Rightarrow$ by (3.6), $M \vDash w'at'ut'' = w_1\ v\ w'at'ut''x_7 = w_1$,

$\Rightarrow M \vDash w'at'ut'' = w_1\ v\ (w'at'ut'')Bw_1.$

This completes the proof of (5.29).

(5.30) For any string concept I⊆I₀ there is a string concept J⊆I such that

$QT^+ \vdash \forall z \in J \; \forall x,t,t',u,v,w_1,t_1,t_2 \subseteq_p z \; (Env(t,x) \; \& \; z=xt'utt' \; \& \; Lastf(z,tt',u,tt') \; \&$

$\& \; Intf(z,w_1,t_1,v,t_2) \; \rightarrow \; Intf(x,w_1,t_1,v,t_2) \; v \; Lastf(x,t_1,v,t_1)).$

Let $J(x) \equiv I_{3.7}^{\subseteq p}(x) \; \& \; I_{4.23b}^{\subseteq p}(x) \; \& \; I_{5.2}(x) \; \& \; I_{5.3}(x) \; \& \; I_{5.29}(x)$.

Assume $M \vDash Intf(z,w_1,t_1,v,t_2) \; \& \; Lastf(z,tt',u,tt')$

where $M \vDash Env(t,x) \; \& \; z=xt'utt'$ and $M \vDash z \in J$.

Then

$M \vDash \exists w_2 \; w_1 a t_1 v t_2 a w_2 = z = xt'utt' \; \& \; Pref(v,t_1) \; \& \; Tally_b(t_2) \; \& \; t_1 < t_2 \; \&$

$\& \; Max^+T_b(t_1,w_1).$

From $M \vDash Lastf(z,tt',u,tt')$ we have

$M \vDash Pref(u,tt') \; \& \; (z=tt'utt' \; v \; \exists w(watt'utt'=z)).$

Now, from $M \vDash tt'utt'=z=xt'utt'$ we have by (3.6) that $M \vDash t=x$. But this contradicts the hypothesis $M \vDash Env(t,x)$ since $M \vDash Tally_b(t)$.

So we may assume that

$M \vDash \exists w \; watt'utt' = z = w_1 a t_1 v t_2 a w_2.$

From $M \vDash Env(t,x)$ we have $M \vDash MaxT_b(t,x)$, hence by (5.2),

$M \vDash MaxT_b(tt',z).$

We may then apply (5.29) to obtain

$M \vDash w_1 a t_1 v t_2 = w \; v \; (w_1 a t_1 v t_2) B w \; v \; w_1 a t_1 v = wa.$

<u>Case 1.</u> $M \vDash w_1 a t_1 v t_2 = w$.

⇒ by (3.7), $M \vDash tt'utt' = w_2$,

$\Rightarrow M \vDash w_1at_1vt_2at=wat$.

But  $M \vDash xt'utt'=z=watt'utt'$,

$\Rightarrow$ by (3.6),  $M \vDash x=wat$,

$\Rightarrow M \vDash w_1at_1vt_2at=x$,

$\Rightarrow M \vDash Pref(v,t_1)$ & $Tally_b(t_2)$ & $t_1<t_2$ & $\exists w_3\, x=w_1at_1vt_2aw_3$ &

$\qquad\qquad\qquad\qquad\qquad\qquad\qquad$ & $Max^+T_b(t_1,w_1)$,

$\Rightarrow M \vDash Intf(x,w_1,t_1,v,t_2)$.

<u>Case 2.</u>  $M \vDash (w_1at_1vt_2)Bw$.

$\Rightarrow M \vDash \exists w'\, w_1at_1vt_2w'=w$,

$\Rightarrow M \vDash w_1at_1vt_2w'att'utt'=watt'utt'=z=xt'utt'$,

$\Rightarrow$ by (3.6),  $M \vDash w_1at_1vt_2w'at=x$.

Also, $M \vDash w_1at_1vt_2w'att'utt'= z=w_1at_1vt_2aw_2$,

$\Rightarrow$ by (3.7),  $M \vDash w'att'utt'=aw_2$,

$\Rightarrow M \vDash w'=a \lor aBw'$,

$\Rightarrow M \vDash \exists w''\, w_1at_1vt_2aw''=x$,

$\Rightarrow M \vDash Pref(v,t_1)$ & $Tally_b(t_2)$ & $t_1<t_2$ & $\exists w''\, x=w_1at_1vt_2aw''$ &

$\qquad\qquad\qquad\qquad\qquad\qquad\qquad$ & $Max^+T_b(t_1,w_1)$,

$\Rightarrow M \vDash Intf(x,w_1,t_1,v,t_2)$.

<u>Case 3.</u>  $M \vDash w_1at_1v=wa$.

We have that

$\Rightarrow M \vDash (wa)t_2aw_2=(w_1at_1v)t_2aw_2=z=xt'utt'=watt'utt'$,

$\Rightarrow$ by (3.7),  $M \vDash t_2aw_2=tt'utt'$.

From $M \vDash \text{Pref}(u,tt')$ we have $M \vDash \exists u_0\ u=au_0a$,

$\Rightarrow M \vDash t_2aw_2=tt'(au_0a)tt'$,

$\Rightarrow$ by (4.23$^b$), $M \vDash t_2=tt'$.

From $M \vDash watt'utt'=xt'utt'$ by (3.6), $M \vDash wat=x$,

$\Rightarrow M \vDash w_1at_1vt=wat=x$.

From $M \vDash \text{Env}(t,x)$ we have $M \vDash \text{MaxT}_b(t,x)\ \&\ \exists v'\ \text{Lastf}(x,t,v',t)$.

Since we have $M \vDash \text{Pref}(v,t_1)\ \&\ \text{Max}^+\text{T}_b(t_1,w_1)$, we may apply (5.3) to obtain

$\quad M \vDash t_1=t\ \&\ v'=v$.

So we have $M \vDash \text{Pref}(v,t_1)\ \&\ \exists w_1(w_1at_1vt_1=x\ \&\ \text{Max}^+\text{T}_b(t_1,w_1))$,

that is, $M \vDash \text{Lastf}(x,t_1,v,t_1)$, as required.

This completes the proof of (5.30).

(5.31) For any string concept $I\subseteq I_0$ there is a string concept $J\subseteq I$ such that

$QT^+ \vdash \forall x,z\in J \ \forall t_1,u,t_2,t',v,w_1(x=t_1auat_1ax_1$ & $Pref(aua,t_1)$ & $Tally_b(t_2)$ &

& $z=t'avat_2ax_1$ & $Pref(ava,t')$ & $t'<t_2 \rightarrow$

$\rightarrow \forall t_3,t_4,w(t_1<t_3$ & $(\exists w'Intf(z,w',t_3,awa,t_4)$ v $Lastf(z,t_3,awa,t_4)) \rightarrow$

$\rightarrow Fr(x,t_3,awa,t_4))]$.

Let $J \equiv I_{3.7}$ & $I_{4.17b}$ & $I_{4.23b}$.

Assume

$M \vDash x=t_1auat_1ax_1$ & $Pref(aua,t_1)$ & $Tally_b(t_2)$ & $z=t'avat_2ax_1$ & $Pref(ava,t')$

where $M \vDash J(x)$ & $J(z)$.

Let $M \vDash t_1<t_3$ & $t'<t_2$.

(1) Assume $M \vDash \exists w'Intf(z,w',t_3,awa,t_4)$.

$\Rightarrow M \vDash Pref(awa,t_3)$ & $Tally_b(t_4)$ & $t_3<t_4$ &

& $\exists w''(z=w'at_3awat_4aw''$ & $Max^+T_b(t_3,w'))$,

$\Rightarrow M \vDash t'avat_2ax_1=z=w'at_3awat_4aw''$,

$\Rightarrow M \vDash (t'av)Bw'$ v $t'av=w'$ v $w'B(t'av)$.

(1a) $M \vDash w'B(t'av)$.

$\Rightarrow M \vDash \exists w_2 \ w'w_2=t'av$,

$\Rightarrow M \vDash (w'w_2)at_2ax_1=(t'av)at_2ax_1=z=w'at_3awat_4aw''$,

$\Rightarrow$ by (3.7), $M \vDash w_2at_2ax_1=at_3awat_4aw''$,

$\Rightarrow M \vDash w_2a$ v $aBw_2$.

If $M \vDash w_2a$, then $M \vDash aat_2ax_1=at_3awat_4aw''$, whence $M \vDash at_2ax_1=t_3awat_4aw''$,

a contradiction because $M \vDash \text{Tally}_b(t_3)$.

Therefore $M \vDash aBw_2$.

$\Rightarrow M \vDash \exists w_3(aw_3=w_2 \:\&\: (aw_3)at_2ax_1=at_3awat_4aw")$,

$\Rightarrow M \vDash w_3at_2ax_1=t_3awat_4aw"$,

$\Rightarrow$ by (4.14$^b$), $M \vDash w_3=t_3 \lor t_3Bw_3$,

$\Rightarrow M \vDash t_3 \subseteq_p w_3 \subseteq_p w_2 \subseteq_p t'av$,

$\Rightarrow$ by (4.17$^b$), $M \vDash t_3 \subseteq_p t' \lor t_3 \subseteq_p v$.

If $M \vDash t_3 \subseteq_p ava$, then from hypothesis $M \vDash \text{Pref}(ava,t')$ we have

$M \vDash \text{Max}^+T_b(t',ava)$, hence $M \vDash t_3 < t'$.

Hence in any case $M \vDash t_3 \leq t'$.

But from $M \vDash \text{Pref}(ava,t')$ we have

Therefore, we have $M \vDash \text{Tally}_b(t')$.

$\Rightarrow$ from $M \vDash t'(ava)t_2ax_1=z=w'at_3awat_4aw"$, by (4.14$^b$), $M \vDash w'=t' \lor t'Bw'$,

$\Rightarrow M \vDash t' \subseteq_p w'$,

$\Rightarrow M \vDash t_3 \subseteq_p t' \subseteq_p w'$, contradicting $M \vDash \text{Max}^+T_b(t_3,w')$.

(1b) $M \vDash t'av=w' \lor (t'av)Bw'$.

$\Rightarrow M \vDash t'av=w' \lor \exists w_4 \: (t'av)w_4=w'$,

$\Rightarrow M \vDash t'avat_3awat_4aw"=w'at_3awat_4aw" \lor$

$\qquad\qquad\qquad \lor \exists w_4 \: (t'avw_4)a\:t_3awat_4aw"=w'at_3awat_4aw"=z$,

$\Rightarrow M \vDash (t'av)at_3awat_4aw"=z=(t'av)at_2ax_1 \lor$

$\qquad\qquad\qquad \lor (t'avw_4)a\:t_3awat_4aw"=z=(t'av)at_2ax_1$,

$\Rightarrow$ by (3.7), $M \vDash at_3awat_4aw"=at_2ax_1 \lor w_4a\:t_3awat_4aw"=at_2ax_1$,

$\implies$ M ⊨ $(t_1au)at_3awat_4aw'' = (t_1au)at_2ax_1$ v

v $(t_1au)w_4a\ t_3awat_4aw'' = (t_1au)at_2ax_1 = x$.

We now claim that  M ⊨ $Max^+T_b(t_3, t_1auw_4)$.

Assume  M ⊨ $Tally_b(t'')$ & $t'' \subseteq_p t_1auw_4$.

$\implies$ by (4.17$^b$),  M ⊨ $t'' \subseteq_p t_1$ v $t'' \subseteq_p uw_4$.

From  M ⊨ $w_4at_3awat_4aw'' = at_2ax_1$, we have  M ⊨ $w_4 = a$ v $aBw_4$.

If  M ⊨ $w_4 = a$,  then  M ⊨ $aat_3awat_4aw'' = at_2ax_1$,  whence

M ⊨ $at_3awat_4aw'' = t_2ax_1$,  a contradiction because  M ⊨ $Tally_b(t_2)$.

Hence  M ⊨ $aBw_4$.  So  M ⊨ $\exists w_5\ aw_5 = w_4$.

So if  M ⊨ $t'' \subseteq_p uw_4$,  then, by (4.17$^b$),   M ⊨ $t'' \subseteq_p u$ v $t'' \subseteq_p w_5$.

If  M ⊨ $t'' \subseteq_p u \subseteq_p aua$,  then from hypothesis  M ⊨ $Pref(aua, t_1)$  we have

M ⊨ $Max^+T_b(t_1, aua)$,  hence  M ⊨ $t'' < t_1$.

If  M ⊨ $t'' \subseteq_p w_5 \subseteq_p w_4 \subseteq_p w'$,  then  M ⊨ $t'' < t_3$  from M ⊨ $Max^+T_b(t_3, w')$.

Hence we have  M ⊨ $t'' < t_1$ v $t'' < t_3$.

$\implies$ from hypothesis  $t_1 < t_3$, M ⊨ $t'' \leq t_1 < t_3$.

Therefore  M ⊨ $Max^+T_b(t_3, t_1auw_4)$  as claimed.

A fortiori,  M ⊨ $Max^+T_b(t_3, t_1au)$.

So we have that

 M ⊨ $Pref(awa, t_3)$ & $Tally_b(t_4)$ & $t_3 < t_4$ &

& $\exists w_6\ (x = w_6at_3awat_4aw''$ & $Max^+T_b(t_3, w_6))$,

$\implies$ M ⊨ $\exists w_6\ Intf(x, w_6, t_3, awa, t_4)$,  as required.

(2)  Assume  M ⊨ $Lastf(z, t_3, awa, t_4)$.

$\Rightarrow$ $M \vDash \text{Pref}(awa,t_3)$ & $\text{Tally}_b(t_4)$ & $t_3=t_4$ &

& $\exists w'(z=w'at_3awat_4$ & $\text{Max}^+T_b(t_3,w'))$.

Suppose, for a reductio, that $M \vDash z=t_3awat_4$.

$\Rightarrow$ $M \vDash t_3awat_4=z=t'(ava)t_2ax_1$,

$\Rightarrow$ by (4.23$^b$), $M \vDash t_3=t'$,

$\Rightarrow$ from $M \vDash t_3awat_4=t'(ava)t_2ax_1$, by (4.16), $M \vDash t_2\subseteq_p vat_2\subseteq_p w$,

$\Rightarrow$ from hypothesis $M \vDash t'<t_2$, $M \vDash t_3=t'\subseteq_p t_2\subseteq_p w$,

contradicting $M \vDash \text{Pref}(awa,t_3)$.

Therefore $M \vDash z \neq t_3awat_4$.

Hence $M \vDash \exists w'(z=w'at_3awat_4$ & $\text{Max}^+T_b(t_3,w'))$.

Now exactly the same argument applies as in (1), omitting $aw''$, to obtain

$M \vDash \exists w_6(x=w_6at_3awat_4$ & $\text{Max}^+T_b(t_3,w_6))$,

whence from $M \vDash \text{Pref}(awa,t_3)$ & $\text{Tally}_b(t_4)$ & $t_3=t_4$ we have

$M \vDash \text{Lastf}(x,t_3,awa,t_4)$.

This completes the proof of (5.31).

(5.32) For any string concept I⊆I₀ there is a string concept J⊆I such that

$QT^+ \vdash \forall x \in J\ \forall t, t_1, t_2, x', v, w_1 (Env(t,x)\ \&\ Intf(x,w_1,t_1,ava,t_2)\ \&\ x'=w_1at_1avat_1 \rightarrow$

$\rightarrow \exists t', w' (Tally_b(t')\ \&\ x=x't'w't\ \&\ aBw'\ \&\ aEw'))$.

Let $J \equiv I_{3.4}\ \&\ I_{4.10}\ \&\ I_{4.19}\ \&\ I_{4.17b}\ \&\ I_{5.22}$.

Assume $M \vDash Env(t,x)\ \&\ Intf(x,w_1,t_1,ava,t_2)$ where $M \vDash J(x)$.

Let $x'=w_1at_1avat_1$.

$\Rightarrow$ from $M \vDash Intf(x,w_1,t_1,ava,t_2)$,

  $M \vDash Pref(ava,t_1)\ \&\ t_1<t_2\ \&\ \exists w_2\ x=w_1at_1avat_2aw_2\ \&\ Max^+T_b(t_1,w_1)$,

$\Rightarrow$ from $M \vDash Env(t,x)$, $M \vDash \exists v"Lastf(x,t,av"a,t)$,

$\Rightarrow M \vDash Pref(av"a,t)\ \&\ (x=tav"at\ \vee\ \exists w"(x=w"atav"at\ \&\ Max^+T_b(t,w")))$.

Suppose, for a reductio, that $M \vDash x=tav"at$.

$\Rightarrow$ from $M \vDash Pref(av"a,t)$, $M \vDash Firstf(x,t,av"at)\ \&\ Lastf(x,t,av"a,t)$,

$\Rightarrow$ by (5.22), $M \vDash \forall w(w\ \varepsilon\ x \leftrightarrow w=v")$,

$\Rightarrow$ from $M \vDash Intf(x,w_1,t_1,ava,t_2)$, $M \vDash v=v"$,

$\Rightarrow M \vDash Lastf(x,t,ava,t)$.

But this contradicts (5.19).

Therefore, $M \vDash x \neq tav"at$.

$\Rightarrow M \vDash \exists w"(x=w"atav"at\ \&\ Max^+T_b(t,w"))$,

$\Rightarrow$ from $M \vDash Pref(ava,t_1)$, $M \vDash Max^+T_b(t_1,ava)$,

$\Rightarrow$ from $M \vDash Max^+T_b(t_1,w_1)\ \&\ t_1<t_2$, by (4.17$^b$), $M \vDash Max^+T_b(t_2,w_1at_1ava)$,

$\Rightarrow$ from $M \vDash Env(t,x)$, $M \vDash MaxT_b(t,x)$,

$\Rightarrow$ $M \vDash t_2 \leq t$,

$\Rightarrow$ $M \vDash \text{Max}^+T_b(t_1, w_1at_1ava)$.

Now, we also have that $M \vDash w_1at_1avat_2aw_2 = x = w"atav"at$.

$\Rightarrow$ by (4.19), $M \vDash (w_1at_1ava)B(w"a) \vee w_1at_1ava = w"a$.

(1)  $M \vDash (w_1at_1ava)B(w"a)$.

$\Rightarrow$ $M \vDash \exists w_3 \ (w_1at_1ava)w_3 = w"a$,

$\Rightarrow$ $M \vDash (w_1at_1ava)t_2aw_2 = x = w"atav"at = (w_1at_1ava)w_3tav"at$,

$\Rightarrow$ by (3.7), $M \vDash t_2aw_2 = w_3tav"at$,

$\Rightarrow$ by (4.14$^b$), $M \vDash w_3t = t_2 \vee t_2B(w_3t)$.

 (1a)  $M \vDash w_3t = t_2$.

We have $M \vDash t_2 \leq t$.

If $M \vDash t_2 = t$, then $M \vDash w_3t = t$. But then $M \vDash tEt$, contradicting (3.4).

On the other hand, if $M \vDash t_2 < t$, then, since $M \vDash \text{Tally}_b(t_2)$, we have that $M \vDash \exists t_3 \ t_2t_3 = t$.

$\Rightarrow$ since $M \vDash \text{Tally}_b(t_2)$ & $\text{Tally}_b(t_3)$, by (4.10), $M \vDash t_3t_2 = t$,

$\Rightarrow$ $M \vDash w_3t_3t_2 = w_3t = t_2$,

$\Rightarrow$ $M \vDash t_2Et_2$, contradicting (3.4).

 (1b)  $M \vDash t_2B(w_3t)$.

$\Rightarrow$ $M \vDash \exists w_4 \ t_2w_4 = w_3t$,

$\Rightarrow$ $M \vDash t_2aw_2 = (t_2w_4)av"at$,

$\Rightarrow$ by (3.7), $M \vDash aw_2 = w_4av"at$,

$\Rightarrow$ $M \vDash w_4 = a \vee aBw_4$,

$\Rightarrow$ $M \vDash t_2aw_2=t_2aav''at \lor \exists w_5\ t_2aw_2=t_2(aw_5)av''at$,

$\Rightarrow$ from $M \vDash Tally_b(t_2)\ \&\ t_1<t_2$, $M \vDash \exists t'(Tally_b(t')\ \&\ t_1t'=t_2)$,

$\Rightarrow$ $M \vDash x=w_1at_1ava(t_1t')aw_2=(w_1at_1avat_1)t'aav''at\ \lor$

$\lor\ x=w_1at_1ava(t_1t')aw_2=(w_1at_1avat_1)t'aw_5av''at$,

$\Rightarrow$ either way, $M \vDash \exists t',w'(Tally_b(t')\ \&\ x=x't'w't\ \&\ aBw'\ \&\ aEw')$.

(2) $M \vDash w_1at_1ava=w''a$.

$\Rightarrow$ $M \vDash (w_1at_1ava)t_2aw_2=x=w''atav''at=(w_1at_1ava)tav''at$,

$\Rightarrow$ by (3.7), $M \vDash t_2aw_2=tav''at$,

$\Rightarrow$ from $M \vDash Tally_b(t_2)\ \&\ t_1<t_2$, $M \vDash \exists t'(Tally_b(t')\ \&\ t_1t'=t_2)$,

$\Rightarrow$ by (4.23$^b$), $M \vDash t_1t'=t_2=t$,

$\Rightarrow$ $M \vDash w_1at_1avat_1t'aw_2=x=w_1at_1avat_1t'av''at$,

$\Rightarrow$ $M \vDash x=x't'av''at$,

$\Rightarrow$ $M \vDash \exists t',w'(Tally_b(t')\ \&\ x=x't'w't\ \&\ aBw'\ \&\ aEw')$, as required.

This completes the proof of (5.32).

(5.33) For any string concept $I \subseteq I_0$ there is a string concept $J \subseteq I$ such that

$$QT^+ \vdash \forall y \in J\ \forall t, t_1, t_2, w_1, y'(\text{Env}(t,y)\ \&\ \text{Intf}(y, w_1, t_1, aua, t_2)\ \&\ y'By\ \&\ \text{Env}(t_1, y')\ \&$$
$$\&\ \text{Lastf}(y', t_1, aua, t_1) \rightarrow y \neq t_1 auat_1).$$

Let $J \equiv I_{5.4}\ \&\ I_{5.15}\ \&\ I_{5.19}$.

Assume $M \vDash \text{Env}(t,y)\ \&\ y'By\ \&\ \text{Env}(t_1, y')\ \&\ \text{Lastf}(y', t_1, aua, t_1)$ where $M \vDash \text{Intf}(y, w_1, t_1, aua, t_2)$ and $M \vDash J(y)$.

$\Rightarrow$ from $M \vDash \text{Env}(t,y)$, $M \vDash \exists v', t', t''\text{Firstf}(y, t', av'a, t'')$,

$\Rightarrow$ from $M \vDash \text{Env}(t_1, y')\ \&\ y'By$, by (5.4), $M \vDash \exists t_3 \text{Firstf}(y', t', av'a, t_3)$,

$\Rightarrow$ from $M \vDash \text{Intf}(y, w_1, t_1, aua, t_2)$, by (5.19), $M \vDash \neg \text{Firstf}(y, t_1, aua, t_2)$,

$\Rightarrow$ by (5.20), $M \vDash t' < t_1$,

$\Rightarrow$ from $M \vDash \text{Intf}(y, w_1, t_1, aua, t_2)$, $M \vDash \text{Pref}(aua, t_1)$.

Assume, for a reductio, that $M \vDash y' = t_1 auat_1$.

$\Rightarrow$ from $M \vDash \text{Pref}(aua, t_1)$, $M \vDash \text{Firstf}(y', t_1, aua, t_1)$,

$\Rightarrow$ from $M \vDash \text{Firstf}(y', t', av'a, t_3)\ \&\ \text{Firstf}(y', t_1, aua, t_1)$, by (5.15), $M \vDash v' = u$,

$\Rightarrow$ from $M \vDash \text{Env}(t_1, y')$, $M \vDash t' = t_1$,

$\Rightarrow M \vDash t' < t_1 = t'$, contradicting $M \vDash t' \in I \subseteq I_0$.

This completes the proof of (5.33).

(5.34) For any string concept $I \subseteq I_0$ there is a string concept $J \subseteq I$ such that

$QT^+ \vdash \forall x \in J \ \forall u,v,t,t_1,t_2,t',t'',w_1$ (Env(t,x) & Firstf(x,t',ava,t'') &

& Intf(x,$w_1$,$t_1$,aua,$t_2$) $\to$ $t'' \leq t_1$ ).

Let $J \equiv I_{3.8}$ & $I_{4.16}$ & $I_{4.24b}$.

Assume $M \vDash$ Env(t,x) where $M \vDash$ Firstf(x,t',ava,t'') & Intf(x,$w_1$,$t_1$,aua,$t_2$) and $M \vDash J(x)$.

$\Rightarrow M \vDash$ Pref(aua,$t_1$) & Tally$_b$($t_2$) & $t_1 < t_2$ &

& $\exists w_2$ x=$w_1$a$t_1$auau$t_2$a$w_2$ & Max$^+$T$_b$($t_1$,$w_1$)

where $M \vDash$ Pref(ava,t') & Tally$_b$(t'') &

& ((t'=t'' & x=t'avat'') v (t'<t'' & (t'avat''a)Bx)).

<u>Case 1.</u> $M \vDash$ t'=t'' & x=t'avat''.

$\Rightarrow M \vDash$ t'avat''=$w_1$a$t_1$aua$t_2$a$w_2$,

$\Rightarrow$ by (4.14$^b$), $M \vDash w_1$=t' v t'B$w_1$,

$\Rightarrow M \vDash$ t'$\subseteq_p w_1$,

$\Rightarrow$ from $M \vDash$ Max$^+$T$_b$($t_1$,$w_1$), $M \vDash$ t'<$t_1$.

On the other hand, by (4.16), we have $M \vDash t_1$aua$t_2 \subseteq_p$ v.

$\Rightarrow M \vDash$ t'<$t_1 \subseteq_p t_1$aua$t_2 \subseteq_p$ v, contradicting $M \vDash$ Pref(ava,t').

<u>Case 2.</u> $M \vDash$ t'<t'' & (t'avat''a)Bx.

$\Rightarrow M \vDash \exists x_1$ t'avat''a$x_1$=x=$w_1$a$t_1$aua$t_2$a$w_2$,

$\Rightarrow$ by (4.14$^b$), $M \vDash w_1$=t' v t'B$w_1$.

(2a) $M \vDash w_1 = t'$.

$\Rightarrow M \vDash t'avat''ax_1 = x = t'at_1auat_2aw_2$,

$\Rightarrow$ by (3.7), $M \vDash vat''ax_1 = t_1auat_2aw_2$,

$\Rightarrow$ by (4.14$^b$), $M \vDash x_0 = t_1 \lor t_1Bv$,

$\Rightarrow M \vDash t_1 \subseteq_p v$,

$\Rightarrow$ from $M \vDash w_1 = t'$ & $Max^+T_b(t_1, w_1)$, $M \vDash t' < t_1$,

$\Rightarrow M \vDash t' < t_1 \subseteq_p v$, contradicting $M \vDash Pref(ava, t')$.

(2b) $M \vDash t'Bw_1$.

$\Rightarrow M \vDash \exists w_3\ t'w_3 = w_1$,

$\Rightarrow M \vDash t'avat''ax_1 = t'w_3at_1auat_2aw_2$,

$\Rightarrow$ by (3.7), $M \vDash avat''ax_1 = w_3at_1auat_2aw_2$,

$\Rightarrow$ by (3.8), $M \vDash (avat'')B(w_3at_1) \lor avat'' = w_3at_1 \lor (w_3at_1)B(avat'')$.

(2bi) $M \vDash avat'' = w_3at_1$.

$\Rightarrow$ by (4.24$^b$), $M \vDash t'' = t_1$, as required.

(2bii) $M \vDash (w_3at_1)B(avat'')$.

$\Rightarrow M \vDash \exists x_4\ w_3at_1x_4 = avat''$,

$\Rightarrow M \vDash (t'w_3)at_1x_4 = t'avat''$,

$\Rightarrow M \vDash w_1at_1x_4 = t'avat''$,

$\Rightarrow M \vDash w_1at_1x_4ax_1 = t'avat''ax_1 = w_1at_1auat_2aw_2$,

$\Rightarrow$ by (3.7), $M \vDash x_4ax_1 = auat_2aw_2$,

$\Rightarrow M \vDash x_4 = a \lor \exists x_5\ ax_5 = x_4$.

We cannot have $M \vDash x_4 = a$ because $M \vDash Tally_b(t'')$.

$\Rightarrow$ $M \vDash \exists x_5\ ax_5=x_4$,

$\Rightarrow$ $M \vDash w_1at_1ax_5=t'avat''$, by (4.16), $M \vDash t_1 \subseteq_p v$,

$\Rightarrow$ from $M \vDash w_1at_1x_4=t'avat''$ by (4.14$^b$), $M \vDash w_1=t'\ v\ t'Bw_1$,

$\Rightarrow$ $M \vDash t' \subseteq_p w_1$,

$\Rightarrow$ from $M \vDash Max^+T_b(t_1,w_1)$, $M \vDash t'<t_1$,

$\Rightarrow$ $M \vDash t'<t_1 \subseteq_p v$, contradicting $M \vDash Pref(ava,t')$.

(2biii) $M \vDash (avat'')B(w_3at_1)$.

$\Rightarrow$ $M \vDash \exists x_2\ avat''x_2=w_3at_1$,

$\Rightarrow$ $M \vDash t'avat''x_2=t'w_3at_1=w_1at_1$,

$\Rightarrow$ $M \vDash (t'avat''x_2)auat_2aw_2=w_1at_1auat_2aw_2=x=t'avat''ax_1$,

$\Rightarrow$ by (3.7), $M \vDash x_2auat_2aw_2=ax_1$,

$\Rightarrow$ $M \vDash x_2=a\ v\ aBx_2$.

We cannot have $M \vDash x_2=a$ because $M \vDash Tally_b(t_1)$.

$\Rightarrow$ $M \vDash aBx_2$,

$\Rightarrow$ $M \vDash \exists x_4\ ax_4=x_2$,

$\Rightarrow$ $M \vDash t'avat''ax_4=w_1at_1$,

$\Rightarrow$ by (4.15$^b$), $M \vDash x_4=t_1\ v\ t_1Ex_4$,

$\Rightarrow$ $M \vDash t'avat''at_1=w_1at_1\ v\ \exists x_5\ t'avat''ax_5t_1=w_1at_1$,

$\Rightarrow$ by (3.7), $M \vDash t'avat''=w_1\ v\ t'avat''ax_5=w_1a$,

$\Rightarrow$ from $M \vDash t'avat''ax_5=w_1a$, $M \vDash x_5=a\ v\ aEx_5$,

$\Rightarrow$ $M \vDash t'avat''=w_1\ v\ t'avat''aa=w_1a\ v\ \exists x_6\ t'avat''ax_6a=w_1a$,

$\Rightarrow$ $M \vDash t'avat''=w_1\ v\ t'avat''a=w_1\ v\ t'avat''ax_6=w_1$,

$\Rightarrow M \vDash t'avat'' \subseteq_p w_1,$

$\Rightarrow M \vDash t'' \subseteq_p t'avat'' \subseteq_p w_1,$

$\Rightarrow$ from $M \vDash Max^+T_b(t_1, w_1),\ M \vDash t'' < t_1,$ as required.

This completes the proof of (5.34).

(5.35) For any string concept I⊆I₀ there is a string concept J⊆I such that

QT⁺ ⊢ ∀x,z∈J ∀t,u,v,t₁,t₂,t',x₁[Env(t,x) & x=t₁auat₂ax₁ & Firstf(x,t₁,aua,t₂) &

& t'<t₂ & ¬(v ε x) & Max⁺T_b(t',ava) & z=t'avat₂ax₁ → Env(t,z) &

& Firstf(z,t',ava,t₂) & ∀w(w ε z ↔ (w ε x & w≠u) v w=v)].

Let J ≡ I₅.₁₅ & I₅.₁₉ & I₅.₃₁ & I₅.₃₄.

Assume M ⊨ Env(t,x) & x=t₁auat₂ax₁ & z=t'avat₂ax₁

where M ⊨ Firstf(x,t₁,aua,t₂) & t'<t₂ & ¬(v ε x) & Max⁺T_b(t',ava)

and M ⊨ J(x) & J(z).

⟹ from M ⊨ Max⁺T_b(t',ava), M ⊨ Tally_b(t'),

⟹ M ⊨ Tally_b(t') & Max⁺T_b(t',ava),

⟹ M ⊨ Pref(ava,t'),

⟹ from M ⊨ Firstf(x,t₁,aua,t₂), M ⊨ Tally_b(t₂),

⟹ M ⊨ Pref(ava,t') & Tally_b(t₂) & (t'<t₂ & (t'avat₂a)Bz),

⟹ M ⊨ Firstf(z,t',ava,t₂).

This also establishes part (b) of M ⊨ Env(t,z).

To show (a) M ⊨ MaxT_b(t,z), assume that M ⊨ Tally_b(t'') & t''⊆_p z.

⟹ M ⊨ t''⊆_p z=t'avat₂ax₁,

⟹ by (4.17ᵇ), M ⊨ t''⊆_p t'av v t''⊆_p t₂ax₁.

If M ⊨ t''⊆_p t'av, then, by (4.17ᵇ), M ⊨ t''⊆_p t' v t''⊆_p v.

If M ⊨ t''⊆_p t', then from M ⊨ Max⁺T_b(t',ava), M ⊨ t''<t'. So from M ⊨ t''⊆_p t'av

we have M ⊨ t''<t'<t₂, whereas M ⊨ t₂≤t from M ⊨ MaxT_b(t,x). Hence

$$M \vDash t''\subseteq_p t'av \rightarrow t''<t.$$

If $M \vDash t''\subseteq_p t_2ax_1\subseteq_p x$, then from $M \vDash MaxT_b(t,x)$, $M \vDash t''\leq t'$.

Hence $M \vDash t''\leq t'$ in either case.

So $M \vDash MaxT_b(t,z)$, as required.

Next we show that (c) $M \vDash \exists w'\, Lastf(z,t,w',t)$.

From $M \vDash Env(t,x)$, we have $M \vDash \exists w'\, Lastf(x,t,w',t)$.

$\Rightarrow M \vDash Pref(w',t)\ \&\ Tally_b(t)$.

Suppose that $M \vDash t_1auat_2=x$.

$\Rightarrow M \vDash xBx$, contradicting $M \vDash x\in I\subseteq I_0$.

$\Rightarrow M \vDash \neg(t_1=t_2\ \&\ x=t_1auat_2x)$,

$\Rightarrow$ from $M \vDash Firstf(x,t_1,aua,t_2)$, $M \vDash t_1<t_2$,

$\Rightarrow$ by (5.8), $M \vDash \exists w_1(x=w_1atw't\ \&\ Max^+T_b(t,w_1)\ \&\ ((t_1au)Bw_1 \vee t_1au=w_1))$,

$\Rightarrow M \vDash \exists x_3\, t_1ax_3=w_1 \vee t_1au=w_1$,

$\Rightarrow M \vDash (t_1aux_3)atw't=w_1atw't \vee (t_1au)atw't=w_1atw't$,

$\Rightarrow M \vDash t_1aux_3atw't=x=t_1(auat_2)ax_1 \vee t_1auatw't=t_1auat_2ax_1$,

$\Rightarrow$ by (3.7), $M \vDash x_3atw't=at_2ax_1 \vee atw't=at_2ax_1$,

$\Rightarrow M \vDash z=(t'av)at_2ax_1=t'avx_3atw't \vee z=(t'av)at_2ax_1=t'avatw't$.

We claim that $M \vDash Max^+T_b(t,t'avx_3)$.

Suppose $M \vDash Tally_b(t'')\ \&\ t''\subseteq_p t'avx_3$.

$\Rightarrow$ by (4.17[b]), $M \vDash t''\subseteq_p t' \vee t''\subseteq_p vx_3$,

$\Rightarrow$ from $M \vDash x_3atw't=at_2ax_1$, $M \vDash x_3=a \vee aBx_3$.

So if $M \vDash t''\subseteq_p vx_3$, then $M \vDash t''\subseteq_p va \vee \exists x_4\, t''\subseteq_p vax_4$.

Now, if $M \vDash t''\subseteq_p va\subseteq_p ava$, then from $M \vDash \text{Pref}(ava,t')$ we have $M \vDash t''<t'$.

If $M \vDash t''\subseteq_p vax_4$, then, by (4.17$^b$), $M \vDash t''\subseteq_p v\subseteq_p ava \lor t''\subseteq_p x_4\subseteq_p x_3\subseteq_p w_1$.

Then if $M \vDash t''\subseteq_p v$ or $M \vDash t''\subseteq_p t'$, we have $M \vDash t''\leq t'<t_2\leq t$, whereas if

$M \vDash t''\subseteq_p x_4$, then $M \vDash t''<t$ from $M \vDash \text{Max}^+T_b(t,w_1)$.

Hence in any case $M \vDash \text{Max}^+T_b(t,t'avx_3)$, as claimed.

A fortiori, $M \vDash \text{Max}^+T_b(t,t'av)$.

It follows that $M \vDash \text{Pref}(w',t) \,\&\, \text{Tally}_b(t) \,\&\, \exists w_1 (z=w_1 atw' \,\&\, \text{Max}^+T_b(t,w_1))$.

Hence $M \vDash \text{Lastf}(z,t,w',t)$.

For part (d) of $M \vDash \text{Env}(t,z)$, assume $M \vDash \text{Fr}(z,t_3,av'a,t_4) \,\&\, \text{Fr}(z,t_5,av'a,t_6)$.

Assume (d1) $M \vDash \text{Firstf}(z,t_3,av'a,t_4)$.

$\Rightarrow$ by (5.19), $M \vDash \neg\exists w_5 \text{Intf}(z,w_5,t_5,av'a,t_6)$,

$\Rightarrow M \vDash \text{Firstf}(z,t_5,av'a,t_6) \lor \text{Lastf}(z,t_5,av'a,t_6)$.

  (d1a) $M \vDash \text{Firstf}(z,t_5,av'a,t_6)$.

$\Rightarrow M \vDash (t_3a)Bz \,\&\, (t_5a)Bz$,

$\Rightarrow M \vDash \exists z_1,z_2 (t_3az_1=z=t_5az_2 \,\&\, \text{Tally}_b(z_1) \,\&\, \text{Tally}_b(z_2))$,

$\Rightarrow$ by (4.23$^b$), $M \vDash t_3=t_5$, as required.

  (d1b) $M \vDash \text{Lastf}(z,t_5,av'a,t_6)$.

$\Rightarrow$ from (5.31), $M \vDash \text{Fr}(x,t_5,av'a,t_6)$,

$\Rightarrow$ from $M \vDash \text{Firstf}(z,t_3,av'a,t_4) \,\&\, \text{Firstf}(z,t',ava,t_2)$, by (5.15), $M \vDash v'=v$,

$\Rightarrow M \vDash v \,\varepsilon\, x$, contradicting the principal hypothesis.

Assume (d2) $M \vDash \exists w_5 \text{Intf}(z,w_5,t_3,av'a,t_4)$.

$\Rightarrow$ by (5.19), $M \vDash \neg\text{Firstf}(z,t_5,av'a,t_6) \,\&\, \neg\text{Lastf}(z,t_5,av'a,t_6)$,

$\Rightarrow$ M ⊨ $\exists w_6 \text{Intf}(z,w_6,t_5,\text{av'a},t_6)$,

$\Rightarrow$ by (5.31), M ⊨ $\text{Fr}(x,t_3,\text{av'a},t_4)$ & $\text{Fr}(x,t_5,\text{av'a},t_6)$,

$\Rightarrow$ by (d) of M ⊨ $\text{Env}(t,x)$, M ⊨ $t_3 = t_5$, as required.

Assume (d3) M ⊨ $\text{Lastf}(z,t_3,\text{av'a},t_4)$.

$\Rightarrow$ by (5.19), M ⊨ $\neg \exists w_6 \text{Intf}(z,w_5,t_5,\text{av'a},t_6)$,

$\Rightarrow$ M ⊨ $\text{Firstf}(z,t_5,\text{av'a},t_6)$ v $\text{Lastf}(z,t_5,\text{av'a},t_6)$.

(d3a) M ⊨ $\text{Firstf}(z,t_5,\text{av'a},t_6)$.

Exactly analogous to (d1b).

(d3b) M ⊨ $\text{Lastf}(z,t_5,\text{av'a},t_6)$.

$\Rightarrow$ by (5.31), M ⊨ $\text{Fr}(x,t_3,\text{av'a},t_4)$ & $\text{Fr}(x,t_5,\text{av'a},t_6)$,

$\Rightarrow$ by (d) of M ⊨ $\text{Env}(t,x)$, M ⊨ $t_3 = t_5$, as required.

This completes the proof of part (d) of M ⊨ $\text{Env}(t,z)$.

For part (e), assume that M ⊨ $\text{Fr}(z,t_3,\text{av'a},t_5)$ & $\text{Fr}(z,t_3,\text{av''a},t_6)$.

Assume (e1) M ⊨ $\text{Firstf}(z,t_3,\text{av'a},t_5)$.

(e1a) M ⊨ $\text{Firstf}(z,t_3,\text{av''a},t_6)$.

$\Rightarrow$ from M ⊨ $\text{Firstf}(z,t',\text{ava},t_2)$, by (5.15), M ⊨ $v' = \text{ava} = v''$, as required.

(e1b) M ⊨ $\exists w_5 \text{Intf}(z,w_5,t_3,\text{av''a},t_6)$.

$\Rightarrow$ M ⊨ $\exists w_6 (z = w_5 a t_3 \text{av''} a t_6 a w_6$ & $\text{Max}^+ T_b(t_3,w_5))$,

$\Rightarrow$ from M ⊨ $\text{Firstf}(z,t_3,\text{av'a},t_5)$, M ⊨ $(t_3 a)Bz$ & $\text{Tally}_b(t_3)$,

$\Rightarrow$ M ⊨ $\exists z_1\, t_3 a z_1 = z = w_5 a t_3 \text{av''} a t_6 a w_6$,

$\Rightarrow$ by (4.14$^b$), M ⊨ $w_5 = t_3$ v $t_3 B w_5$,

$\Rightarrow$ M ⊨ $t_3 \subseteq_p w_5$, contradicting M ⊨ $\text{Max}^+ T_b(t_3,w_5)$.

(e1c)  $M \vDash \text{Lastf}(z,t_3,av''a,t_6)$.

From (c) above,  $M \vDash \text{Lastf}(z,t,aw'a,t)$.

$\Rightarrow$ by (5.15),  $M \vDash v'=w'=v''$,  as required.

(e2)  $M \vDash \exists w_5 \text{Intf}(z,w_5,t_3,av'a,t_5)$.

(e2a)  $M \vDash \text{Firstf}(z,t_3,av''a,t_6)$.

Exactly analogous to (e1b).

(e2b)  $M \vDash \exists w_7 \text{Intf}(z,w_7,t_3,av''a,t_6)$.

$\Rightarrow$ by (5.31),  $M \vDash \text{Fr}(x,t_3,av'a,t_5) \& \text{Fr}(x,t_3,av''a,t_6)$,

$\Rightarrow$ by part (e) of  $M \vDash \text{Env}(t,x)$,  $M \vDash v'=v''$,  as required.

(e2c)  $M \vDash \text{Lastf}(z,t_3,av''a,t_6)$.

Same as (e2b), appealing to (5.31) and the principal hypothesis.

(e3)  $M \vDash \text{Lastf}(z,av'a,t_5)$.

(e3a)  $M \vDash \text{Firstf}(z,t_3,av''a,t_6)$.

Exactly analogous to (e1c).

(e3b)  $M \vDash \exists w_5 \text{Intf}(z,w_5,t_3,av''a,t_6)$.

Exactly analogous to (e2c).

(e3c)  $M \vDash \text{Lastf}(z,t_3,av''a,t_6)$.

As in (e3b), we appeal directly to (5.31) and the hypothesis  $M \vDash \text{Env}(t,x)$.

This completes the proof of part (e) of $M \vDash \text{Env}(t,z)$ and the proof of

$M \vDash \text{Env}(t,z)$.

We now show that

$$M \vDash \forall w(w \,\varepsilon\, z \leftrightarrow (w \,\varepsilon\, x \,\&\, w \neq u) \vee w=v).$$

Assume $M \vDash w \, \varepsilon \, z$.

$\Rightarrow M \vDash \exists t_3, t_4 \, Fr(z, t_3, awa, t_4)$.

(1) $M \vDash Firstf(z, t_3, awa, t_4)$.

$\Rightarrow$ from $M \vDash Firstf(z, t', ava, t_2)$, by (5.15), $M \vDash w = v$.

(2) $M \vDash \exists w_1 Intf(z, w_1, t_3, awa, t_4) \lor Lastf(z, t_3, awa, t_4)$.

$\Rightarrow$ by the proof of (5.31), $M \vDash \exists w_1 Intf(x, w_1, t_3, awa, t_4) \lor Lastf(x, t_3, awa, t_4)$,

$\Rightarrow M \vDash Fr(x, t_3, awa, t_4)$,

$\Rightarrow M \vDash w \, \varepsilon \, x$.

Suppose, for a reductio, that $M \vDash w = u$.

From the hypothesis we have that $M \vDash Firstf(x, t_1, aua, t_2)$.

$\Rightarrow$ by (d) of $M \vDash Env(t, x)$, $M \vDash t_1 = t_3$.

(2a) $M \vDash \exists w_1 Intf(x, w_1, t_3, awa, t_4)$.

$\Rightarrow$ from $M \vDash w = u$, $M \vDash Intf(x, w_1, t_3, aua, t_4)$.

But this contradicts $M \vDash Firstf(x, t_1, aua, t_2)$ by (5.19).

(2b) $M \vDash Lastf(x, t_3, awa, t_4)$.

$\Rightarrow M \vDash t_3 = t_4 \, \& \, (x = t_3 awat_4 \lor \exists w_1(x = w_1 at_3 awat_4 \, \& \, Max^+T_b(t_3, w_1)))$.

(2bi) $M \vDash x = t_3 awat_4$.

$\Rightarrow M \vDash t_1 auat_2 ax_1 = x = t_3 awat_4$,

$\Rightarrow$ from $M \vDash Env(t, x)$, $M \vDash t_3 = t_4 = t$,

$\Rightarrow$ from $M \vDash Tally_b(t_1) \, \& \, Tally_b(t_3)$, by (4.23[b]), $M \vDash t_1 = t_3$.

On the other hand, we have as in (c) above, that $M \vDash t_1 < t_2$.

$\Rightarrow M \vDash t_1 < t_2 \leq t_3 = t_1$, contradicting $M \vDash t_1 \in I \subseteq I_0$.

(2bii) $M \vDash \exists w_1(x=w_1at_3awat_4 \,\&\, Max^+T_b(t_3,w_1))$.

$\Rightarrow$ from $M \vDash Firstf(x,t_1,aua,t_2)$, $M \vDash (t_1a)Bx$,

$\Rightarrow M \vDash \exists x_1 \, t_1ax_1=x=w_1at_3awat_4$,

$\Rightarrow$ by (4.14$^b$), $M \vDash w_1=t_1 \lor t_1Bw_1$,

$\Rightarrow M \vDash t \subseteq_p w_1$,

$\Rightarrow M \vDash t_3=t_1 \subseteq_p w_1$, contradicting $M \vDash Max^+T_b(t_3,w_1)$.

Therefore, $M \vDash w \neq u$.

So we proved that $M \vDash \forall w(w \, \varepsilon \, z \to (w \, \varepsilon \, x \,\&\, w \neq u) \lor w=v)$.

Conversely, assume $M \vDash w=v$.

$\Rightarrow$ from $M \vDash Firstf(z,t',ava,t_2)$, $M \vDash w \, \varepsilon \, z$.

Assume $M \vDash w \, \varepsilon \, x \,\&\, w \neq u$.

$\Rightarrow M \vDash \exists t_3,t_4 \, Fr(x,t_3,awa,t_4)$.

(i) $M \vDash Firstf(x,t_3,awa,t_4)$.

$\Rightarrow$ from $M \vDash Firstf(x,t_1,aua,t_2)$, by (5.15), $M \vDash w=u$, contradicting the hypothesis.

(ii) $M \vDash \exists w_1 Intf(x,w_1,t_3,awa,t_4)$.

$\Rightarrow M \vDash Pref(awa,t_3) \,\&\, Tally_b(t_4) \,\&\, t_3<t_4 \,\&$

$\qquad\qquad\qquad\qquad \&\, \exists w_2 \, x=w_1at_3awat_4aw_2 \,\&\, Max^+T_b(t_3,w_1)$,

$\Rightarrow M \vDash t_1auat_2ax_1=x=w_1at_3awat_4aw_2$,

$\Rightarrow$ from $M \vDash Firstf(x,t_1,aua,t_2)$, by (5.34), $M \vDash t_2 \leq t_3$,

$\Rightarrow$ by principal hypothesis, $M \vDash t'<t_2 \leq t_3$,

$\Rightarrow M \vDash Pref(aua,t_1) \,\&\, ((t_1=t_2 \,\&\, x=t_1auat_2) \lor (t_1<t_2 \,\&\, (t_1auat_2a)Bx))$.

If $M \vDash t_1=t_2$ & $x=t_1auat_2$, then $M \vDash t_1auat_2=x=t_1auat_2ax_1$. But then $M \vDash xBx$, contradicting $M \vDash x \in I \subseteq I_0$.

Therefore $M \vDash \neg(t_1=t_2$ & $x=t_1auat_2)$.

$\Rightarrow M \vDash t_1<t_2$ & $(t_1auat_2a)Bx$.

With $M \vDash t'<t_3$ & $t_1<t_2$, we now apply (5.31), reversing the roles of x and z and derive from $M \vDash \exists w_1 Intf(x,w_1,t_3,awa,t_4)$ that $M \vDash Fr(z,t_3,awa,t_4)$.

$\Rightarrow M \vDash w\ \varepsilon\ z$, as required.

(iii) $M \vDash Lastf(x,t_3,awa,t_4)$.

$\Rightarrow M \vDash Pref(awa,t_3)$ & $Tally_b(t_4)$ & $t_3=t_4=t$ & $(x=t_3awat_4\ v$

& $\exists w_1\ x=w_1at_3awat_4$ & $Max^+T_b(t_3,w_1))$.

If $M \vDash x=t_3awat_4$, then $M \vDash x=t_3awat_4=x=t_1auat_2ax_1$, whence by (4.23$^b$), $M \vDash t_1=t_3$. From hypothesis $M \vDash Firstf(x,t_1,aua,t_2)$ we have $M \vDash t_1<t_2$, as in (ii) above. But then from $M \vDash MaxT_b(t,x)$ it follows that $M \vDash t_1<t_2\leq t=t_3=t_1$, contradicting $M \vDash t_1 \in I \subseteq I_0$.

Therefore, $M \vDash \exists w_1\ x=w_1at_3awat_4$ & $Max^+T_b(t_3,w_1))$.

$\Rightarrow M \vDash w_1at_3awat_4=x=t_1auat_2ax_1$,

$\Rightarrow$ as in (ii), $M \vDash t_1<t_2$,

$\Rightarrow$ by (5.34), $M \vDash t_2\leq t_3$,

$\Rightarrow M \vDash t_1<t_2\leq t_3$,

$\Rightarrow$ by principal hypothesis, $M \vDash t'<t_2$,

$\Rightarrow$ by (5.31), $M \vDash Fr(z,t_3,awa,t_4)$,

$\Rightarrow M \vDash w\ \varepsilon\ z$.

This completes the proof that  M ⊨ ∀w((w ε x & w≠u) v w=v  →  w ε z),

and the proof of (5.35).

(5.36) For any string concept $I \subseteq I_0$ there is a string concept $J \subseteq I$ such that

$QT^+ \vdash \forall x \in J\ \forall t, t_1, t_2, t', t'', w_1, u, w_2, z [\text{Env}(t,x)\ \&\ x = w_1 a(t_1 aua) t_2 a w_2\ \&$

$\&\ \text{Firstf}(x, t', ava, t'')\ \&\ \text{Intf}(x, w_1, t_1, aua, t_2)\ \&\ z = w_1 a t_2 a w_2 \rightarrow$

$\rightarrow \text{Firstf}(z, t', ava, t_2) \vee \text{Firstf}(z, t', ava, t'')]$

.

Let $J \equiv I_{5.34}$.

Assume $M \vDash \text{Env}(t,x)\ \&\ \text{Firstf}(x, t', ava, t'')\ \&\ \text{Intf}(x, w_1, t_1, aua, t_2)$ where

$M \vDash x = w_1 a(t_1 aua) t_2 a w_2$ and $M \vDash J(x)$.

Let $M \vDash z = w_1 a t_2 a w_2$.

$\Longrightarrow M \vDash \text{Pref}(ava, t')\ \&\ \text{Tally}_b(t'')\ \&\ ((t' = t''\ \&\ x = t' avat'')\ \vee$

$\vee\ (t' < t''\ \&\ (t' avat'' a) Bx))$,

and

$M \vDash \text{Pref}(aua, t_1)\ \&\ \text{Tally}_b(t_2)\ \&\ t_1 < t_2\ \&\ x = w_1 a(t_1 aua) t_2 a w_2\ \&\ \text{Max}^+ T_b(t_1, w_1)$.

We follow the pattern of the proof of (5.34).

<u>Case 1.</u> $M \vDash t' = t''\ \&\ x = t' avat''$.

Exactly the same as in (5.34).

<u>Case 2.</u> $M \vDash t' < t''\ \&\ (t' avat'' a) Bx$.

$\Longrightarrow M \vDash \exists x_1\ t' avat'' a x_1 = x = w_1 a t_1 auat_2 a w_2$,

$\Longrightarrow$ by (4.14$^b$), $M \vDash w_1 = t' \vee t' B w_1$.

(2a) $M \vDash w_1 = t'$.

We derive a contradiction exactly as in (5.34).

(2b) $M \vDash t' B w_1$.

$\Longrightarrow M \vDash \exists w_3\ t' w_3 = w_1$,

$\Rightarrow$ M ⊨ t'avat''ax$_1$=t'w$_3$at$_1$auat$_2$aw$_2$,

$\Rightarrow$ by (3.7), M ⊨ avat''ax$_1$=w$_3$at$_1$auat$_2$aw$_2$,

$\Rightarrow$ M ⊨ (avat'')B(w$_3$at$_1$auat$_2$aw$_2$) & (w$_3$at$_1$)B(w$_3$at$_1$auat$_2$aw$_2$),

$\Rightarrow$ by (3.8), M ⊨ (avat'')B(w$_3$at$_1$) v avat''=w$_3$at$_1$ v (w$_3$at$_1$)B(avat'').

(2bi)  M ⊨ avat''=w$_3$at$_1$.

$\Rightarrow$ by (4.24$^b$), M ⊨ t''=t$_1$,

$\Rightarrow$ by (3.6), M ⊨ ava=w$_3$a,

$\Rightarrow$ M ⊨ t'ava=t'w$_3$a=w$_1$a,

$\Rightarrow$ M ⊨ (t'ava)t$_2$aw$_2$=(w$_1$a)t$_2$aw$_2$=z.

We have that M ⊨ t'<t''=t$_1$<t$_2$. Thus we obtain

$$M ⊨ Pref(ava,t') \& Tally_b(t_2) \& (t'<t_2 \& (t'avat_2a)Bz),$$

whence M ⊨ Firstf(z,t',ava,t$_2$).

(2bii)  M ⊨ (w$_3$at$_1$)B(avat'').

We derive a contradiction exactly as in (5.34).

(2biii)  M ⊨ (avat'')B(w$_3$at$_1$).

We follow the argument given in (5.34) to obtain

$$M ⊨ t'avat''=w_1 \text{ v } t'avat''a=w_1 \text{ v } t'avat''ax_6=w_1,$$

$\Rightarrow$ M ⊨ (t'avat'')at$_2$aw$_2$=w$_1$at$_2$aw$_2$ v (t'avat''a)at$_2$aw$_2$=w$_1$at$_2$aw$_2$ v

v (t'avat''ax$_6$)at$_2$aw$_2$=w$_1$at$_2$aw$_2$,

$\Rightarrow$ M ⊨ z=t'avat''at$_2$aw$_2$ v z=t'avat''aat$_2$aw$_2$ v z=t'avat''ax$_6$at$_2$aw$_2$.

In any case, we have

$$M ⊨ Pref(ava,t') \& Tally_b(t'') \& (t'<t'' \& (t'avat''a)Bz),$$

whence $M \vDash \text{Firstf}(z,t',ava,t'')$.

This completes the proof of (5.36).

(5.37) For any string concept $I \subseteq I_0$ there is a string concept $J \subseteq I$ such that

$QT^+ \vdash \forall x \in J \ \forall t, t_1, t_2, u, v, w_1, w_2, z [Env(t,x) \ \& \ x = w_1 a(t_1 aua) t_2 a w_2 \ \&$

$\& \ Intf(x, w_1, t_1, aua, t_2) \ \& \ z = w_1 a t_2 a w_2 \ \& \ Lastf(x, t, ava, t) \rightarrow Lastf(z, t, ava, t)]$.

Let $J \equiv I_{3.6} \ \& \ I_{3.10} \ \& \ I_{4.17b} \ \& \ I_{4.23b}$.

Assume $M \vDash Env(t,x) \ \& \ Intf(x, w_1, t_1, aua, t_2) \ \& \ Lastf(x, t, ava, t)$ where

$M \vDash x = w_1 a(t_1 aua) t_2 a w_2$ and $M \vDash J(x)$.

Let $M \vDash z = w_1 a t_2 a w_2$.

$\Rightarrow M \vDash Pref(aua, t_1) \ \& \ Tally_b(t_2) \ \& \ t_1 < t_2 \ \& \ \exists w_3 \ x = w_1 a(t_1 aua) t_2 a w_3 \ \&$

$\& \ Max^+ T_b(t_1, w_1)$

and $M \vDash Pref(ava, t') \ \& \ \exists w'(x = w' atavat \ \& \ Max^+ T_b(t, w'))$.

$\Rightarrow M \vDash w_1 a(t_1 aua) t_2 a w_2 = x = w' atavat$,

$\Rightarrow M \vDash (t_2 a w_2) Ex \ \& \ (tavat) Ex$,

$\Rightarrow$ by (3.10), $M \vDash (t_2 a w_2) E(tavat) \ \vee \ t_2 a w_2 = tavat \ \vee \ (tavat) E(t_2 a w_2)$.

(1) $M \vDash (t_2 a w_2) E(tavat)$.

$\Rightarrow M \vDash \exists w_3 \ tavat = w_3(t_2 a w_2)$,

$\Rightarrow M \vDash (w'a) w_3(t_2 a w_2) = (w'a) tavat = x = w_1 a(t_1 aua) t_2 a w_2$,

$\Rightarrow$ by (3.6), $M \vDash w' a w_3 = w_1 a(t_1 aua)$,

$\Rightarrow M \vDash w_3 = a \ \vee \ aEw_3$.

Now, $M \vDash w_3 \neq a$ because $M \vDash Tally_b(t)$.

$\Rightarrow M \vDash aEw_3$,

$\Rightarrow M \vDash \exists w_4 \ w_3 = w_4 a$,

$\Rightarrow$ M ⊨ (w$_4$a)t$_2$aw$_2$=tavat,

$\Rightarrow$ by (4.14$^b$), M ⊨ w$_4$=t v tBw$_4$.

(1a)  M ⊨ w$_4$=t.

$\Rightarrow$ M ⊨ tat$_2$aw$_2$=tavat,

$\Rightarrow$ M ⊨ (w'a)tat$_2$aw$_2$=w'atavat=x=w$_1$at$_1$auat$_2$aw$_2$,

$\Rightarrow$ by (3.6), M ⊨ w'at=w$_1$at$_1$au,

$\Rightarrow$ by (4.15$^b$), M ⊨ u=t v tEu,

$\Rightarrow$ M ⊨ t⊆$_p$u,

$\Rightarrow$ from M ⊨ Env(t,x),  M ⊨ t$_1$<t$_2$≤t,

$\Rightarrow$ M ⊨ t$_1$⊆$_p$t⊆$_p$u, contradicting M ⊨ Pref(aua,t$_1$).

(1b)  M ⊨ tBw$_4$.

$\Rightarrow$ M ⊨ ∃w$_5$ tw$_5$=w$_4$,

$\Rightarrow$ M ⊨ (tw$_5$)at$_2$aw$_2$=tavat,

$\Rightarrow$ M ⊨ (w'a)tw$_5$at$_2$aw$_2$=w'atavat=x=w$_1$at$_1$auat$_2$aw$_2$,

$\Rightarrow$ by (3.6), M ⊨ w'atw$_5$=w$_1$at$_1$au,

$\Rightarrow$ M ⊨ (tw$_5$)E(w'atw$_5$) & (t$_1$au)E(w'atw$_5$),

$\Rightarrow$ by (3.10), M ⊨ (tw$_5$)E(t$_1$au) v tw$_5$=t$_1$au v (t$_1$au)E(tw$_5$).

(1bi)  M ⊨ (tw$_5$)E(t$_1$au) v tw$_5$=t$_1$au.

$\Rightarrow$ M ⊨ ∃w$_6$ t$_1$au=w$_6$(tw$_5$) v tw$_5$=t$_1$au,

$\Rightarrow$ M ⊨ t⊆$_p$t$_1$au,

$\Rightarrow$ by (4.17$^b$), M ⊨ t⊆$_p$t$_1$ v t⊆$_p$u,

$\Rightarrow$ from M ⊨ Env(t,x),  M ⊨ t$_1$<t$_2$≤t.

So $M \vDash t \subseteq_p t_1 < t$ contradicts $M \vDash t \in I \subseteq I_0$. If $M \vDash t \subseteq_p u$, then $M \vDash t_1 \subseteq_p t \subseteq_p u$, contradicting $M \vDash \text{Pref}(aua, t_1)$.

(1bii) $M \vDash (t_1 au)E(tw_5)$.

$\Rightarrow M \vDash \exists w_6 \; (tw_5) = w_6(t_1 au)$,

$\Rightarrow M \vDash w_6 t_1 aua = tw_5 a = w_4 a = w_3$,

$\Rightarrow M \vDash (w_6 t_1 aua) t_2 aw_2 = w_3 t_2 aw_2 = tavat$,

$\Rightarrow M \vDash (w'a) w_6 t_1 auat_2 aw_2 = (w'a) tavat = x = w_1 at_1 auat_2 aw_2$,

$\Rightarrow$ by (3.6), $M \vDash w'aw_6 = w_1 a$,

$\Rightarrow M \vDash w_6 = a \lor aEw_6$.

Now, $M \vDash w_6 \neq a$ because $M \vDash \text{Tally}_b(t)$.

$\Rightarrow M \vDash aEw_6$,

$\Rightarrow M \vDash \exists w_7 \; w_6 = w_7 a$,

$\Rightarrow M \vDash (w_7 a) t_1 au = tw_5$,

$\Rightarrow$ by (4.14[b]), $M \vDash w_7 = t \lor tBw_7$,

$\Rightarrow M \vDash t \subseteq_p w_7$.

But $M \vDash w'a(w_7 a) = w'aw_6 = w_1 a$,

$\Rightarrow M \vDash w'aw_7 = w_1$,

$\Rightarrow M \vDash t \subseteq_p w_7 \subseteq_p w_1$,

$\Rightarrow$ from $M \vDash \text{Env}(t,x)$, $M \vDash t_1 < t_2 \leq t$,

$\Rightarrow M \vDash t_1 \subseteq_p t \subseteq_p w_1$, contradicting $M \vDash \text{Max}^+T_b(t_1, w_1)$.

(2) $M \vDash t_2 aw_2 = tavat$.

$\Rightarrow M \vDash t_2 aw_2 = tavat$,

$\Rightarrow$ by (4.23$^b$), $M \vDash t_2=t$,

$\Rightarrow$ from $M \vDash \text{Max}^+T_b(t_1,w_1)$ & $t_1<t_2$, $M \vDash \text{Max}^+T_b(t_2,w_1)$,

$\Rightarrow$ from $M \vDash t_2aw_2=tavat$, $M \vDash z=w_1a(t_2aw_2)=w_1atavat$,

$\Rightarrow$ $M \vDash \text{Pref}(ava,t)$ & $\exists w_1(z=w_1atavat$ & $\text{Max}^+T_b(t,w_1))$,

$\Rightarrow$ $M \vDash \text{Lastf}(z,t,ava,t)$.

(3)  $M \vDash (tavat)E(t_2aw_2)$.

$\Rightarrow$ $M \vDash \exists w_3\ t_2aw_2=w_3(tavat)$,

$\Rightarrow$ $M \vDash w_3=a \lor aEw_3$.

$\Rightarrow$ $M \vDash (w_1at_1aua)w_3tavat=(w_1at_1aua)t_2aw_2=x=w'atavat$,

$\Rightarrow$ by (3.6), $M \vDash w_1at_1auaw_3=w'a$.

Now, $M \vDash w_3 \neq a$ because $M \vDash \text{Tally}_b(t_2)$.

$\Rightarrow$ $M \vDash aEw_3$,

$\Rightarrow$ $M \vDash \exists w_4\ w_3=w_4a$,

$\Rightarrow$ $M \vDash (w_4a)tavat=t_2aw_2$,

$\Rightarrow$ $M \vDash (w_1a)w_4atavat=(w_1a)t_2aw_2=z$.

We claim that $M \vDash \text{Max}^+T_b(t,wt_1aw_4)$.

Assume $M \vDash \text{Tally}_b(t_0)$ & $t_0 \subseteq_p w_1aw_4$. Then by (4.17$^b$), $M \vDash t_0 \subseteq_p w_1 \lor t_0 \subseteq_p w_4$.

If $M \vDash t_0 \subseteq_p w_1$, then from $M \vDash \text{Max}^+T_b(t_1,w_1)$, $M \vDash t_0<t_1<t_2\leq t$.

If $M \vDash t_0 \subseteq_p w_4$, note that from $M \vDash w_1at_1aua(w_4a)=w'a$, we have

$M \vDash w_1at_1auaw_4=w'$. Hence $M \vDash w_4 \subseteq_p w'$, so $M \vDash t_0 \subseteq_p w'$, and $M \vDash t_0<t$ from

$M \vDash \text{Max}^+T_b(t,w')$.

Therefore, in either case $M \vDash t_0<t$. So $M \vDash \text{Max}^+T_b(t,w_1aw_4)$ as claimed.

Thus we have $M \vDash \text{Pref}(ava,t)\ \&\ \exists w_0(z=w_0 atavat\ \&\ \text{Max}^+T_b(t,w_0))$ where $w_0=w_1 aw_4$. Hence $M \vDash \text{Lastf}(z,t,ava,t)$.

This completes the proof of (5.37).

(5.38) For any string concept $I\subseteq I_0$ there is a string concept $J\subseteq I$ such that

$QT^+ \vdash \forall x\in J\ \forall t,t',t_1,t_2,t_3,t_4,u,v,w,w_1,w',w'',z[Env(t,x)\ \&\ x=w'at'avat''aw''\ \&$

$\&\ x=w_1at_1auat_2\ \&\ t''aw''=t_1auat_2\ \&\ Env(t',z)\ \&\ Intf(x,w',t',ava,t'')\ \&$

$\&\ Lastf(x,t_1,aua,t_2)\ \&\ z=w'at'avat'' \to$

$\to (Firstf(x,t_3,awa,t_4)\ \vee\ \exists w_2 Intf(x,w_2,t_3,awa,t_4)\ \to\ \exists t_5\ Fr(z,t_3,awa,t_5))]$.

Let $J \equiv I_{4.16}\ \&\ I_{5.8}\ \&\ I_{5.30}$.

Assume $M \vDash Env(t,x)\ \&\ Intf(x,w',t',ava,t'')\ \&\ Lastf(x,t_1,aua,t_2)$ where

$M \vDash x=w'at'avat''aw''\ \&\ x=w_1at_1auat_2\ \&\ t''aw''=t_1auat_2$ and $M \vDash J(x)$.

Assume also that $M \vDash z=w'at'avat''\ \&\ Env(t',z)$ and

$M \vDash Firstf(x,t_3,awa,t_4)\ \vee\ \exists w_2 Intf(x,w_2,t_3,awa,t_4)$.

$\Rightarrow$ from $M \vDash Env(t,x)$, $M \vDash MaxT_b(t,x)\ \&\ t_1=t_2=t$,

$\Rightarrow$ from $M \vDash Intf(x,w',t',ava,t'')$,

$M \vDash Pref(ava,t')\ \&\ Tally_b(t'')\ \&\ t'<t''\ \&\ \exists w''\ x=w'at'avat''aw''\ \&\ Max^+T_b(t',w')$,

$\Rightarrow\ M \vDash t'<t''\leq t_1$,

$\Rightarrow M \vDash \exists t_0(Tally_b(t_0)\ \&\ t_1=t't_0)$.

(1) $M \vDash Firstf(x,t_3,awa,t_4)$.

$\Rightarrow M \vDash Pref(awa,t_3)\ \&\ Tally_b(t_4)\ \&$

$\&\ ((t_3=t_4\ \&\ x=t_3awat_4)\ \vee\ (t_3<t_4\ \&\ (t_3awat_4a)Bx))$.

(1a) $M \vDash t_3=t_4\ \&\ x=t_3awat_4$.

$\Rightarrow$ from $M \vDash Pref(awa,t_3)$, $M \vDash Max^+T_b(t_3,awa)$,

$\Rightarrow M \vDash t_3awat_4=x=w_1atauat$.

But this contradicts (4.21).

(1b)  $M \vDash t_3<t_4$ & $(t_3awat_4a)Bx$.

$\Rightarrow$ from  $M \vDash \exists x_2 \, (t_3awat_4a)x_2=x$,

$\Rightarrow M \vDash t_3<t_4\leq t$,

$\Rightarrow$ from  $M \vDash Firstf(x,t_3,awa,t_4)$ & $Lastf(x,t,aua,t)$ & $t_3<t_4\leq t$,

$\Rightarrow$ by (5.8),  $M \vDash \exists w_3 \, (x=w_3atauat$ & $((t_3aw)Bw_3 \lor t_3aw=w_3))$,

$\Rightarrow M \vDash w_1at_1auat_2=w_3atauat$ & $t_1=t_2=t$,

$\Rightarrow$ by (3.6),  $M \vDash w_1=w_3$,

$\Rightarrow M \vDash (t_3aw)Bw_1 \lor t_3aw=w_1$.

(1bi)  $M \vDash t_3aw=w_1$.

$\Rightarrow M \vDash (t_3aw)at_1auat_2=w_1at_1auat_2=x$,

$\Rightarrow M \vDash zt_0auat_2=(w'at'avat')t_0auat_2=w'at'avat_1auat_2=w'at'avat''aw''=x$,

$\Rightarrow M \vDash t_3awat't_0auat_2=zt_0auat_2$,

$\Rightarrow$ by (3.6),  $M \vDash t_3awat'=z=w'at'avat'$,

$\Rightarrow$ by (4.16),  $M \vDash t'\subseteq_p w$,

$\Rightarrow$ from  $M \vDash Env(t',z)$, $M \vDash MaxT_b(t',z)$,

$\Rightarrow M \vDash t_3\leq t'\subseteq_p w$, contradicting  $M \vDash Pref(awa,t_3)$.

(1bii)  $M \vDash (t_3aw)Bw_1$.

$\Rightarrow M \vDash \exists w_4 \, (t_3aw)w_4=w_1$,

$\Rightarrow$ from  $M \vDash w'at'avat''aw''=x=w_1at_1auat_2$ & $t''aw''=t_1auat_2$, by (3.6),

$\quad M \vDash w'at'av=w_1$,

$\Rightarrow M \vDash t_3aww_4=w'at'av$,

$\implies M \vDash (t_3aww_4)at_1auat_2=w_1at_1auat_2=x$,

$\implies M \vDash t_3aww_4at_1auat_2=t_3awat_4ax_2$,

$\implies$ by (3.7), $M \vDash w_4at_1auat_2=at_4ax_2$,

$\implies M \vDash w_4=a \lor aBw_4$.

If $M \vDash w_4=a$, then $M \vDash aat_1auat_2=at_4ax_2$, whence $M \vDash at_1auat_2=t_4ax_2$, contradicting $M \vDash Tally_b(t_4)$.

$\implies M \vDash aBw_4$,

$\implies M \vDash \exists w_5 \, aw_5=w_4$,

$\implies M \vDash t_3aw(aw_5)at_1auat_2=t_3awat_4ax_2$,

$\implies$ by (3.7), $M \vDash w_5at_1auat_2=t_4ax_2$,

$\implies$ by (4.14$^b$), $M \vDash w_5=t_4 \lor t_4Bw_5$.

(1bii1) $M \vDash w_5=t_4$.

$\implies M \vDash w_1=(t_3aw)w_4=t_3awaw_5=t_3awat_4$,

$\implies M \vDash (t_3awat_4a)Bz$,

$\implies M \vDash Pref(awa,t_3) \, \& \, Tally_b(t_4) \, \& \, (t_3<t_4 \, \& \, (t_3awat_4a)Bz)$,

$\implies M \vDash Firstf(z,t_3,awa,t_4)$,

$\implies M \vDash Fr(z,t_3,awa,t_4)$, as required.

(1bii2) $M \vDash t_4Bw_5$.

$\implies M \vDash \exists w_6 \, t_4w_6=w_5$,

$\implies M \vDash w_1=(t_3aw)w_4=t_3awaw_5=t_3awat_4w_6$,

$\implies M \vDash t_3awat_4ax_2=x=w_1at_1auat_2=(t_3awat_4w_6)at_1auat_2$,

$\implies$ by (3.7), $M \vDash ax_2=w_6at_1auat_2$,

$\implies$ M ⊨ $w_6$=a v a$Bw_6$,

$\implies$ M ⊨ $t_3$awa$t_4$a=$w_1$ v $\exists w_7$ (a$w_7$=$w_6$ & $t_3$awa$t_4$a$w_7$=$w_1$),

$\implies$ either way, M ⊨ Pref(awa,$t_3$) & Tally$_b$($t_4$) & ($t_3$<$t_4$ & ($t_3$awa$t_4$a)Bz),

$\implies$ M ⊨ Firstf(z,$t_3$,awa,$t_4$), $\implies$ M ⊨ Fr(z,$t_3$,awa,$t_4$), as required.

(2) M ⊨ $\exists w_2$ Intf(x,$w_2$,$t_3$,awa,$t_4$).

$\implies$ M ⊨ Env(t',z) & x=z$t_0$aua$t'_0$ & Lastf(x,$t'_0$,aua,$t'_0$) & Intf(x,$w_2$,$t_3$,awa,$t_4$),

$\implies$ by (5.30), M ⊨ Intf(z,$w_2$,$t_3$,awa,$t_4$) v Lastf(z,$t_3$,awa,$t_3$),

$\implies$ M ⊨ $\exists t_5$ Fr(z,$t_3$,awa,$t_5$), as required.

This completes the proof of (5.38).

(5.39) For any string concept $I \subseteq I_0$ there is a string concept $J \subseteq I$ such that

$QT^+ \vdash \forall x,x' \in J \; \forall t_1,t_2,t_3,t_4,w,v,z(Env(t_2,x') \; \& \; x'=t_1wt_2 \; \& \; aBw \; \& \; aEw \; \& \; x=t_1wz \; \&$

$\& \; Env(t,z) \; \& \; t_2<t_3 \; \& \; Firstf(z,t_3,ava,t_4) \; \rightarrow \; \forall t',t'',u(Fr(x,t',aua,t'') \rightarrow$

$\rightarrow \; \exists t_6 \; Fr(x',t',aua,t_6) \; v \; Fr(z,t_3,aua,t_4) \; v \; \exists t_6 \; Fr(z,t',aua,t_6)))$ .

Let $J \equiv I_{5.4} \; \& \; I_{5.9} \; \& \; I_{5.15} \; \& \; I_{5.28} \; \& \; I_{5.30}$.

Assume $M \vDash Env(t_2,x') \; \& \; x'=t_1wt_2 \; \& \; aBw \; \& \; aEw$ along with

$M \vDash x=t_1wz \; \& \; Env(t,z) \; \& \; t_2<t_3 \; \& \; Firstf(z,t_3,ava,t_4)$ and $M \vDash J(x) \; \& \; J(x')$.

Suppose that $M \vDash Fr(x,t',aua,t'')$.

We distinguish three cases:

(1) $M \vDash Firstf(x,t',aua,t'')$.

We first observe that $M \vDash x'Bx$.

For, from $M \vDash Firstf(z,t_3,ava,t_4)$ we have $M \vDash (t_3a)Bz \; \& \; Tally_b(t_3)$.

$\Rightarrow M \vDash \exists z_1 \; t_3az_1=z$,

$\Rightarrow M \vDash t_1wz=x=t_1w(t_3az_1)$,

$\Rightarrow$ from $M \vDash Env(t_2,x')$, $M \vDash Tally_b(t_2)$,

$\Rightarrow$ from hypothesis $M \vDash t_2<t_3$, $M \vDash \exists t_0 \; t_3=t_2t_0$,

$\Rightarrow M \vDash x=t_1wt_2t_0az_1=x't_0az_1$,

$\Rightarrow M \vDash x'Bx$, as claimed.

But then from $M \vDash Env(t_2,x') \; \& \; x'Bx \; \& \; Firstf(x,t',aua,t'')$ we have, by (5.4),

that $M \vDash \exists t_6 Firstf(x',t',aua,t_6)$. Hence $M \vDash \exists t_6 Fr(x',t',aua,t_6)$, as required.

(2) $M \vDash \exists w' \; Intf(x,w',t',aua,t'')$.

$\implies$ $M \vDash \text{Pref}(aua,t')$ & $\text{Tally}_b(t'')$ & $\exists w''\ x=w'at'auat''aw''$ & $\text{Max}^+T_b(t',w')$,

$\implies$ from $M \vDash \text{Firstf}(z,t_3,ava,t_4)$,

$M \vDash \text{Pref}(ava,t_3)$ & $\text{Tally}_b(t_4)$ & $((t_3=t_4\ \&\ z=t_3avat_4)\ \vee\ (t_3<t_4\ \&\ (t_3avat_4a)Bz))$,

$\implies$ from $M \vDash aEw$, $M \vDash \exists w_2\ w=w_2a$,

$\implies$ $M \vDash x=t_1w_2az=t_1w_2at_3avat_4\ \vee\ \exists z_2\ x=t_1w_2az=t_1w_2a(t_3avat_4az_2)$.

Assume $M \vDash t_0 \subseteq_p t_1w_2 \subseteq_p x'$ & $\text{Tally}_b(t_0)$. Then from $M \vDash \text{Env}(t_2,x')$ we have

$M \vDash t_0 \leq t_2 < t_3$. Hence $M \vDash \text{Max}^+T_b(t_3,t_1w_2)$.

$\implies$ $M \vDash (x=t_1w_2at_3avat_4\ \&\ \text{Lastf}(x,t_3,ava,t_4))\ \vee\ \text{Intf}(x,t_1w_2,t_3,ava,t_4)$.

  (2a)  $M \vDash x=t_1w_2at_3avat_4$ & $\text{Lastf}(x,t_3,ava,t_4)$.

$\implies$ $M \vDash t_1w_2az=x=t_1w_2at_3avat_4$,

$\implies$ by (3.7), $M \vDash z=t_3avat_4$,

$\implies$ since $M \vDash t_3=t_4$, $M \vDash z=t_3avat_3$,

$\implies$ $M \vDash t_1wz=x=t_1wt_3avat_3=t_1w(t_2t_0)ava(t_2t_0)=x't_0avat_2t_0$,

$\implies$ $M \vDash \text{Env}(t_2,x')$ & $x=x't_0avat_2t_0$ & $\text{Lastf}(x,t_2t_0,ava,t_2t_0)$ & $\text{Intf}(x,w',t',aua,t'')$,

$\implies$ by (5.30), $M \vDash \text{Intf}(x,w',t',aua,t'')\ \vee\ \text{Lastf}(x',t',aua,t')$,

$\implies$ $M \vDash \exists t_6 \text{Fr}(x',t',aua,t_6)$, as required.

  (2b)  $M \vDash \text{Intf}(x,t_1w_2,t_3,ava,t_4)$.

We distinguish three cases:

  (2bi)  $M \vDash t_4=t'\ \vee\ t''=t_3$.

    (2bi1)  $M \vDash t_4=t'$.

$\implies$ by (5.28), $M \vDash t_1w_2at_3avat_4=w'at'$,

$\implies$ by (3.6), $M \vDash t_1w_2at_3av=w'$.

Also, $M \vDash (t_1w_2at_3avat_4)auat''aw''=(w'at')auat''aw''=x=t_1w_2az$.

$\Rightarrow$ by (3.7), $M \vDash z=t_3avat_4auat''aw''$.

Assume $M \vDash t_0 \subseteq_p t_3av$ & $Tally_b(t_0)$. Then $M \vDash t_0 \subseteq_p t_3av \subseteq_p t_1w_2at_3av=w'$,

whence from $M \vDash Max^+T_b(t',w')$, $M \vDash t_0<t'$.

Letting $w_0=t_3av$, we then have

$M \vDash Pref(aua,t')$ & $Tally_b(t'')$ & $t'<t''$ & $\exists w_0(z=w_0at'auat''aw''$ & $Max^+T_b(t',w_0))$

whence $M \vDash \exists w_0\, Intf(z,w_0,t',aua,t'')$, so $M \vDash Fr(z,t',aua,t'')$, as required.

(2bi2) $M \vDash t''=t_3$.

$\Rightarrow$ by (5.28), $M \vDash w'at'auat''=t_1w_2at_3$,

$\Rightarrow$ by (3.6), $M \vDash w'at'aua=t_1w_2a=t_1w$,

$\Rightarrow M \vDash w'at'auat_2=t_1wt_2=x'$,

$\Rightarrow$ from $M \vDash Env(t_2,x')$, $M \vDash \exists v'Lastf(x',t_2,av'a,t_2)$,

$\Rightarrow M \vDash Pref(av'a,t_2)$ & $(x'=t_2av'at_2 \lor \exists w_4\,(x'=w_4at_2av'at_2$ & $Max^+T_b(t_2,w_4))$.

If $M \vDash x'=t_2av'at_2$, then $M \vDash t_2av'at_2=x'=w'at'auat_2$.

$\Rightarrow$ by (4.14$^b$), $M \vDash w'=t_2 \lor t_2Bw'$,

$\Rightarrow M \vDash t_2 \subseteq_p w'$,

$\Rightarrow$ from $M \vDash Max^+T_b(t',w')$, $M \vDash t_2<t'$.

But, by (4.16), $M \vDash t' \subseteq_p v'$,

$\Rightarrow$ from $M \vDash Pref(av'a,t_2)$, $M \vDash t'<t_2$,

$\Rightarrow M \vDash t_2<t'<t_2$, contradicting $M \vDash t' \in I \subseteq I_0$.

Therefore, $M \vDash \exists w_4\,(x'=w_4at_2av'at_2$ & $Max^+T_b(t_2,w_4))$.

$\Rightarrow M \vDash w'at'auat_2=x'=w_4at_2av'at_2$,

$\Rightarrow$ $M \vDash MaxT_b(t_2,x')$ & $Lastf(x',t_2,av'a,t_2)$ & $w'at'auat_2=x'$ & $Tally_b(t_2)$ &

& $Pref(av'a,t_2)$ & $Max^+T_b(t',w')$, by (5.3), $M \vDash t'=t_2$ & $u=v'$,

$\Rightarrow$ $M \vDash Lastf(x',t',aua,t')$, $\Rightarrow$ $M \vDash \exists t_6 Fr(x',t',aua,t_6)$, as required.

(2bii) $M \vDash t_4 \neq t'$ & $t'' \neq t_3$ & $u=v$.

$\Rightarrow$ from $M \vDash Firstf(z,t_3,ava,t_4)$, $M \vDash Firstf(z,t_3,aua,t_4)$,

$\Rightarrow$ $M \vDash Fr(z,t_3,aua,t_4)$, as required.

(2biii) $M \vDash t_4 \neq t'$ & $t'' \neq t_3$ & $u \neq v$.

$\Rightarrow$ $M \vDash (w'at')Bx$ & $(t_1w_2at_3)Bx$,

$\Rightarrow$ by (3.8), $M \vDash (w'at')B(t_1w_2at_3)$ v $w'at'=t_1w_2at_3$ v $(t_1w_2at_3)B(w'at')$.

(2biii1) $M \vDash (w'at')B(t_1w_2at_3)$.

We reason as in (5.27), Case 2, and derive

$M \vDash w'at'auat''=t_1w_2$ v $w'at'auat''a=t_1w_2$ v $\exists w_5$ $w'at'auat''aw_5=t_1w_2$.

Since $M \vDash Pref(aua,t')$ & $Tally_b(t'')$ & $t'<t''$ & $Max^+T_b(t',w')$ & $Env(t_2,x')$, we have, by (5.9), that

$M \vDash x' \neq w'at'auat''at_2$ & $x' \neq w'at'auat''aat_2$,

$\Rightarrow$ $M \vDash w'at'auat''at_2 \neq t_1w_2$ & $w'at'auat''aat_2 \neq t_1w_2$,

$\Rightarrow$ $M \vDash \exists w_5$ $w'at'auat''aw_5=t_1w_2$,

$\Rightarrow$ $M \vDash \exists w_6$ $w'at'auat''aw_6=t_1w_2at_2=x'$ & $Max^+T_b(t',w')$,

$\Rightarrow$ $M \vDash Intf(x',w',t',aua,t'')$,

$\Rightarrow$ $M \vDash Fr(x',t',aua,t'')$, as required.

(2biii2) $M \vDash (t_1w_2at_3)B(w'at')$.

We reason as in (5.27), Case 3, analogously to (2biii1) and obtain

$M \vDash t_1w_2at_3avat_4 = w' \lor t_1w_2at_3avat_4a = w' \lor \exists w_5\, t_1w_2at_3avat_4aw_5 = w'$.

Suppose that $M \vDash t_1w_2at_3avat_4 = w'$.

$\Rightarrow M \vDash (t_1w_2at_3avat_4)at'auat''aw'' = w'at'auat''aw'' = x = t_1w_2az$,

$\Rightarrow$ by (3.7), $M \vDash z = t_3avat_4at'auat''aw''$,

$\Rightarrow$ from $M \vDash \text{Intf}(x,w',t',aua,t'')$, $M \vDash \text{Max}^+T_b(t',w')$,

$\Rightarrow$ since $M \vDash t_3avat_4 \subseteq_p w'$, $M \vDash \text{Max}^+T_b(t',t_3avat_4)$.

So we have

$M \vDash \text{Pref}(aua,t') \,\&\, \text{Tally}_b(t'') \,\&\, t' < t'' \,\&\, z = t_3avat_4at'auat''aw'' \,\&$

$\&\, \text{Max}^+T_b(t',t_3avat_4)$,

$\Rightarrow M \vDash \text{Intf}(z,t_3avat_4,t',aua,t'')$, $\Rightarrow M \vDash \text{Fr}(z,t',aua,t'')$, as required.

Exactly analogous arguments apply if

$M \vDash t_1w_2at_3avat_4a = w' \lor \exists w_5\, t_1w_2at_3avat_4aw_5 = w'$.

(2biii3) $M \vDash w'at' = t_1w_2at_3$.

$\Rightarrow$ by (4.24$^b$), $M \vDash t'' = t_3$,

$\Rightarrow$ from $M \vDash \text{Pref}(aua,t')$, $M \vDash \text{Pref}(aua,t_3)$,

$\Rightarrow M \vDash w'at'auat''aw'' = x = (t_1w_2at_3)auat''aw''$,

$\Rightarrow$ from $M \vDash x = t_1w_2az$, by (3.7), $M \vDash z = t_3auat''aw''$,

$\Rightarrow M \vDash \text{Pref}(aua,t_3) \,\&\, \text{Tally}_b(t'') \,\&\, (t_3 < t'' \,\&\, (t_3auat''a)Bz)$,

$\Rightarrow M \vDash \text{Firstf}(z,t_3,aua,t'')$,

$\Rightarrow M \vDash \text{Fr}(z,t',aua,t'')$, as required.

(3) $M \vDash \text{Lastf}(x,t',aua,t'')$.

We first observe that $M \vDash \text{Max}^+T_b(t,t_1w)$. (This follows by the same argument as $M \vDash \text{Max}^+T_b(t,t_1w_2)$ in (2) above.)

Now, from $M \vDash \text{Env}(t,z)$, $M \vDash \exists v'\text{Lastf}(z,t,av'a,t)$.

$\implies$ from $M \vDash \text{Max}^+T_b(t,t_1w)$, by the proof given in (5.25), parts (1) and (3),

$$M \vDash \text{Lastf}(x,t,av'a,t).$$

$\implies$ from $M \vDash \text{Lastf}(x,t',aua,t'')$, by (5.15), $M \vDash t=t'=t''$ & $u=v'$,

$\implies M \vDash \text{Lastf}(z,t',aua,t''),$

$\implies M \vDash \text{Fr}(z,t',aua,t''),$ as required.

This completes the proof of (5.39).

(5.40)  For any string concept $I \subseteq I_0$ there is a string concept $J \subseteq I$ such that

$QT^+ \vdash \forall x \in J \; \forall t,t',t'',v(Env(t,x) \& Firstf(x,t',ava,t'') \rightarrow x \neq t'avat''at \& x \neq t'avat''at)$.

Let $J \equiv I_{4.16} \& I_{5.20}$.

Assume $M \vDash Env(t,x) \& Firstf(x,t',ava,t'')$ where $M \vDash J(x)$.

$\Rightarrow M \vDash \exists u \; Lastf(x,t,aua,t)$,

$\Rightarrow M \vDash Pref(aua,t) \& (x=tauat \lor \exists w(x=watauat \& Max^+T_b(t,w)))$.

(1) Suppose, for a reductio, that $M \vDash x=t'avat''at$.

  (1a)  $M \vDash x=tauat$.

$\Rightarrow M \vDash t'avat''at=x=tauat$,

$\Rightarrow$ by (4.16), $M \vDash t'' \subseteq_p aua$,

$\Rightarrow$ since $M \vDash Tally_b(t') \& Tally_b(t)$, by (4.23$^b$), $M \vDash t'=t$,

$\Rightarrow$ from $M \vDash Firstf(x,t',ava,t'')$ by (5.20),  $M \vDash t' \leq t''$,

$\Rightarrow M \vDash t=t' \leq t'' \subseteq_p aua$, contradicting $M \vDash Pref(aua,t)$.

  (1b)  $M \vDash \exists w(x=watauat \& Max^+T_b(t,w))$.

$\Rightarrow M \vDash watauat=t'avat''at$,

$\Rightarrow$ by (3.6),  $M \vDash watau=t'avat''$,

$\Rightarrow$ by (4.16),  $M \vDash t \subseteq_p v$,

$\Rightarrow$ from  $M \vDash Env(t,x)$,  $M \vDash t' \leq t$,

$\Rightarrow M \vDash t' \leq t \subseteq_p v$, contradicting $M \vDash Pref(ava,t')$ and $M \vDash Firstf(x,t',ava,t'')$.

Therefore,  $M \vDash x \neq t'avat''at$.

(2) Suppose, for a reductio, that  $M \vDash x=t'avat''aat$.

(2a)  $M \vDash x=tauat$.

The same argument as in (1a).

(2b)  $M \vDash \exists w(x=watauat \ \& \ Max^+T_b(t,w))$.

$\Rightarrow M \vDash watauat=t'avat''aat$,

$\Rightarrow$ by (3.6), $M \vDash watau=t'avat''a$,

$\Rightarrow M \vDash u=a \ \vee \ aEu$,

$\Rightarrow M \vDash wataa=t'avat''a \ \vee \ \exists u_1(u=u_1a \ \& \ watau_1a=t'avat''a)$,

$\Rightarrow M \vDash wata=t'avat'' \ \vee \ \exists u_1 \ watau_1=t'avat''$.

Now, we have that $M \vDash wata \neq t'avat''$.

$\Rightarrow M \vDash watau_1=t'avat''$, whence a contradiction follows as in (1b) with $u_1$ in place of u.

Therefore, $M \vDash x \neq t'avat''aat$.

This completes the proof of (5.40).

(5.41) For any string concept $I \subseteq I_0$ there is a string concept $J \subseteq I$ such that

$QT^+ \vdash \forall y \in J \ \forall t, t_0, t_1, t_2, t', w_1, w_2, u, y', y''$ (Env(t,y) & $y = t_0 w_1 t_1 t_2 w_2 t$ & $Tally_b(t_0)$ &

& $aBw_1$ & $aEw_1$ & $aBw_2$ & $aEw_2$ & $y' = t_0 w_1 t_1$ & $y'' = t_1 t_2 w_2 t$ & $Env(t_1, y')$ &

& $Env(t, y'')$ & $Firstf(y'', t_1 t_2, aua, t') \rightarrow \neg \exists w (w \ \varepsilon \ y' \ \& \ w \ \varepsilon \ y''))$.

Let $J \equiv I_{5.6} \ \& \ I_{5.20} \ \& \ I_{5.25}$.

Assume $M \vDash Env(t,y) \ \& \ y = t_0 w_1 t_1 t_2 w_2 t$ where

$M \vDash Tally_b(t_0) \ \& \ aBw_1 \ \& \ aEw_1 \ \& \ aBw_2 \ \& \ aEw_2 \ \& \ y' = t_0 w_1 t_1 \ \& \ y'' = t_1 t_2 w_2 t \ \&$

& $Env(t_1, y')$ & $Env(t, y'')$ & $Firstf(y'', t_1 t_2, aua, t')$

and $M \vDash J(y)$.

Suppose, for a reductio, that $M \vDash \exists w (w \ \varepsilon \ y' \ \& \ w \ \varepsilon \ y'')$.

$\Rightarrow M \vDash \exists t_3, t_4 \ Fr(y', t_3, awa, t_4) \ \& \ \exists t_5, t_6 \ Fr(y'', t_5, awa, t_6)$,

$\Rightarrow M \vDash Env(t_1, y'), \ M \vDash MaxT_b(t_1, y')$,

$\Rightarrow M \vDash t_3 \leq t_1$,

$\Rightarrow$ from $M \vDash Firstf(y'', t_1 t_2, aua, t') \ \& \ Fr(y'', t_5, awa, t_6)$, by (5.20), $M \vDash t_1 t_2 \leq t_5$,

$\Rightarrow M \vDash t_3 \leq t_1 < t_1 t_2 \leq t_5$.

We claim that $M \vDash Max^+ T_b(t_5, t_0 w_1)$.

From $M \vDash MaxT_b(t_1, y') \ \& \ t_1 < t_1 t_2 \leq t_5$ and $M \vDash t_1 \in I \subseteq I_0$, we have

$M \vDash MaxT_b(t_5, y')$. Hence from $M \vDash t_0 w_1 \subseteq_p y'$, $M \vDash Max^+ T_b(t_5, t_0 w_1)$ as claimed.

Now, we have that

$M \vDash Env(t_1, y') \ \& \ Fr(y', t_3, awa, t_4) \ \& \ y' t_2 w_2 t = y \ \& \ aBw_2 \ \& \ aEw_2 \ \& \ Tally_b(t_0) \ \&$

& $Tally_b(t)$,

$\Rightarrow$ by (5.6),  $M \vDash \exists t_7 Fr(y,t_3,awa,t_7)$.

On the other hand, we also have

$\quad\quad M \vDash Fr(y'',t_5,awa,t_6)$ & $y=t_0w_1y''$ & $aBw_1$ & $aEw_1$ & $Max^+T_b(t_5,t_0w_1)$.

$\Rightarrow$ by (5.25),  $M \vDash Fr(y,t_5,awa,t_6)$,

$\Rightarrow M \vDash Fr(y',t_3,awa,t_7)$ & $Env(t,y)$,   $M \vDash t_5=t_3$,

$\Rightarrow M \vDash t_5=t_3<t_5$,  contradicting  $M \vDash t_5 \in I \subseteq I_0$.

Therefore,  $M \vDash \neg \exists w(w \varepsilon y'$ & $w \varepsilon y'')$.

This completes the proof of (5.41).

(5.42) For any string concept $I \subseteq I_0$ there is a string concept $J \subseteq I$ such that

$QT^+ \vdash \forall x \in J\ \forall u,t_1,t_2,t_3,t_4\ (Env(t,x)\ \&\ Fr(x,t_1,aua,t_2)\ \&\ Fr(x,t_3,aua,t_4) \rightarrow$

$\rightarrow t_1=t_3\ \&\ t_2=t_4).$

Let $J \equiv I_{4.20}\ \&\ I_{5.15}\ \&\ I_{5.19}\ \&\ I_{5.20}$.

Assume that $M \vDash Env(t,x)\ \&\ Fr(x,t_1,aua,t_2)\ \&\ Fr(x,t_3,aua,t_4)$ where $M \vDash J(x)$.

$\Rightarrow M \vDash t_1=t_3$,

$\Rightarrow M \vDash Fr(x,t_1,aua,t_2)\ \&\ Fr(x,t_1,aua,t_4)$.

(1) $M \vDash Firstf(x,t_1,aua,t_2)$.

<u>Claim</u>: $M \vDash Firstf(x,t_1,aua,t_4)$.

For, if $M \vDash \neg Firstf(x,t_1,aua,t_4)$, then, by (5.20), $M \vDash t_1 < t_1$, contradicting

$M \vDash t_1 \in I \subseteq I_0$. Hence, $M \vDash Firstf(x,t_1,aua,t_4)$, as claimed.

$\Rightarrow$ by (5.15), $M \vDash t_2=t_4$.

(2) $M \vDash \exists w_1\ Intf(x,w_1,t_1,aua,t_2)$.

$\Rightarrow$ from $M \vDash Env(t,x)$, by (5.19),

$\quad M \vDash \neg Firstf(x,t_1,aua,t_4)\ \&\ \neg Lastf(x,t_1,aua,t_4)$,

$\Rightarrow M \vDash \exists w_3\ Intf(x,w_3,t_1,aua,t_4)$,

$\Rightarrow M \vDash \exists w_2\ x=w_1at_1auat_2aw_2\ \&\ Tally_b(t_2)\ \&\ Max^+T_b(t_1,w_1)\ \&$

$\quad\quad\quad\ \&\ \exists w_4\ x=w_3at_1auat_4aw_4\ \&\ Tally_b(t_4)\ \&\ Max^+T_b(t_1,w_3)$,

$\Rightarrow M \vDash w_1at_1auat_2aw_2 = w_3at_1auat_4aw_4$,

$\Rightarrow$ by (4.20), $M \vDash w_1=w_3$,

$\Rightarrow M \vDash (w_1at_1aua)t_2aw_2 = (w_1at_1aua)t_4aw_4$,

$\Rightarrow$ by (3.7), $M \vDash t_2aw_2=t_4aw_4$,

$\Rightarrow$ since $M \vDash Tally_b(t_2)$ & $Tally_b(t_4)$, by (4.23$^b$), $M \vDash t_2=t_4$, as required.

(3)  $M \vDash Lastf(x,t_1,aua,t_2)$.

$\Rightarrow$ $M \vDash t_1=t_2$,

$\Rightarrow$ by (5.19), $M \vDash \neg \exists w_1 \, Intf(x,w_1,t_1,aua,t_4)$,

$\Rightarrow$ $M \vDash Firstf(x,t_1,aua,t_4) \vee Lastf(x,t_1,aua,t_4)$.

 (3a)  $M \vDash Firstf(x,t_1,aua,t_4)$.

$\Rightarrow$ $M \vDash (t_1=t_4 \,\&\, x=t_1auat_4) \vee (t_1<t_4 \,\&\, (t_1auat_4a)Bx)$.

   (3ai)  $M \vDash t_1=t_4 \,\&\, x=t_1auat_4$.

$\Rightarrow$ $M \vDash t_4=t_1=t_2$, as required.

   (3aii)  $M \vDash t_1<t_4 \,\&\, (t_1auat_4a)Bx$.

$\Rightarrow$ $M \vDash \exists x_1 \, t_1auat_4ax_1=x$.

Suppose, for a reductio, that $M \vDash x=t_1auat_2$.

$\Rightarrow$ $M \vDash (t_1aua)t_2=(t_1aua)t_4ax_1$,

$\Rightarrow$ by (3.7), $M \vDash t_2=t_4ax_1$, a contradiction because $M \vDash Tally_b(t_2)$.

Therefore, $M \vDash x \neq t_1auat_2$.

$\Rightarrow$ from $M \vDash Lastf(x,t_1,aua,t_2)$, $M \vDash \exists w_1 \, (x=w_1at_1auat_2 \,\&\, Max^+T_b(t_1,w_1))$,

$\Rightarrow$ $M \vDash t_1auat_4ax_1=x=w_1at_1auat_2$,

$\Rightarrow$ since $M \vDash Tally_b(t_1)$, by (4.14$^b$), $M \vDash t_1 \subseteq_p w_1$,

which contradicts $M \vDash Max^+T_b(t_1,w_1)$.

 (3b)  $M \vDash Lastf(x,t_1,aua,t_4)$.

$\Rightarrow$ $M \vDash t_4=t_1=t_2$, as required.

This completes the proof of (5.42).

(5.43) For any string concept $I \subseteq I_0$ there is a string concept $J \subseteq I$ such that

$QT^+ \vdash \forall x \in J \, \forall t_1, v_0, t_2, t_3, v, t_4 \, (Firstf(x, t_1, av_0 a, t_2) \, \& \, Lastf(x, t_3, ava, t_4) \rightarrow$

$\rightarrow \forall w_1, t', u, t'', w_2 (Intf(x, w_1, t', aua, t'') \, \& \, x = w_1 a t' a u a t'' a w_2 \, \&$

$\& \, z = t' a u a t'' a w_2 \, \& \, t' \leq t_3 \rightarrow Lastf(z, t_3, ava, t_4)))$.

Let $J \equiv I_{3.8} \, \& \, I_{3.10} \, \& \, I_{5.15} \, \& \, I_{5.19} \, \& \, I_{5.20}$.

Assume $M \vDash Firstf(x, t_1, av_0 a, t_2) \, \& \, Lastf(x, t_3, ava, t_4)$,

and let $M \vDash Intf(x, w_1, t', aua, t'')$ where $M \vDash J(x)$ and $M \vDash t' \leq t_3$.

$\Rightarrow$ by (5.19), $M \vDash \neg Firstf(x, t', aua, t'')$,

$\Rightarrow$ by (5.20), $M \vDash t_1 < t'$,

$\Rightarrow$ from $M \vDash Firstf(x, t_1, av_0 a, t_2)$,

$M \vDash Pref(av_0 a, t_1) \, \& \, Tally_b(t_2) \, \&$

$\& \, ((t_1 = t_2 \, \& \, x = t_1 av_0 at_2) \vee (t_1 < t_2 \, \& \, (t_1 av_0 at_2 a) Bx))$.

Suppose, for a reductio, that $M \vDash x = t_1 av_0 at_2$.

$\Rightarrow M \vDash t_1 av_0 at_2 = x = w_1 at' aua t'' aw_2$,

$\Rightarrow$ by (4.16), $M \vDash t' aut'' \subseteq_p v_0$,

$\Rightarrow M \vDash t_1 \subseteq_p t' \subseteq_p t' aut'' \subseteq_p v_0$, contradicting $M \vDash Pref(av_0 a, t_1)$.

Therefore $M \vDash \neg(x = t_1 av_0 at_2)$.

So we have $M \vDash (t_1 av_0 at_2 a) Bx$.

$\Rightarrow$ from $M \vDash Lastf(x, t_3, ava, t_4)$,

$M \vDash Pref(ava, t_3) \, \& \, Tally_b(t_4) \, \& \, t_3 = t_4 = t \, \&$

$\& \, (t_3 avat_4 = x \vee \exists w'(x = w' at_3 avat_4 \, \& \, Max^+ T_b(t_3, w')))$.

We distinguish the two cases:

(1) $M \vDash t_3avat_4=x$.

$\Rightarrow$ from $M \vDash Pref(ava,t_3)$ & $Tally_b(t_4)$ & $t_3=t_4$, $M \vDash Firstf(x,t_3,ava,t_4)$,

$\Rightarrow$ from $M \vDash Firstf(x,t_1,av_0a,t_2)$ by (5.15), $M \vDash v_0=v$,

$\Rightarrow$ from $M \vDash (t_1av_0at_2a)Bx$, $M \vDash \exists x_1\ (t_1av_0at_2a)x_1=x$,

$\Rightarrow$ $M \vDash t_1av_0at_2ax_1=x=t_3avat_4$,

$\Rightarrow$ since $M \vDash Tally_b(t_1)$ & $Tally_b(t_3)$, by (4.23$^b$), $M \vDash t_1=t_3$,

$\Rightarrow$ $M \vDash t_1avat_2ax_1=t_1avat_4$,

$\Rightarrow$ by (3.7), $M \vDash t_2ax_1=t_4$,

$\Rightarrow$ $M \vDash a\subseteq_p t_4$, contradicting $M \vDash Tally_b(t_4)$.

Therefore $M \vDash \neg(t_3avat_4=x)$.

So we have

(2) $M \vDash \exists w'(x=w'at_3avat_4\ \&\ Max^+T_b(t_3,w'))$.

$\Rightarrow$ $M \vDash w'at_3avat_4=x=w_1at'auat''aw_2$,

$\Rightarrow$ $M \vDash (at_3avat_4)Ex$ & $(at'auat''aw_2)Ex$,

$\Rightarrow$ by (3.10), $M \vDash at_3avat_4=at'auat''aw_2$ v

$\qquad\qquad$ v $(at'auat''aw_2)E(at_3avat_4)$ v $(at_3avat_4)E(at'auat''aw_2)$.

There are three subcases:

(2a) $M \vDash at_3avat_4=at'auat''aw_2$.

$\Rightarrow$ $M \vDash t_3avat_4=t'auat''aw_2=z$,

$\Rightarrow$ $M \vDash Pref(ava,t_3)$ & $Tally_b(t_4)$ & $t_3=t_4$ & $t_3avat_4=z$,

$\Rightarrow$ $M \vDash Lastf(z,t_3,ava,t_4)$, as required.

(2b)  M ⊨ (at'auat''aw₂)E(at₃avat₄).

⟹ M ⊨ ∃w₃ at₃avat₄=w₃at'auat''aw₂,

⟹ M ⊨ (w₃a)B(at₃avat₄) & (at₃)B(at₃avat₄),

⟹ by (3.8),  M ⊨ w₃a=at₃ v (w₃a)B(at₃) v (at₃)B(w₃a).

   (2bi)  M ⊨ w₃a=at₃  is ruled out because  M ⊨ Tally$_b$(t₃).

   (2bii)  M ⊨ (w₃a)B(at₃).

⟹ M ⊨ ∃w₄ at₃=w₃aw₄,

⟹ by (QT4) and (QT5), M ⊨ w₃=a v aBw₃,

⟹ M ⊨ aaw₄=at₃ v ∃w₅ at₃=(aw₅)aw₄,

⟹ M ⊨ aw₄=t₃ v t₃=w₅aw₄,

⟹ M ⊨ a⊆$_p$t₃, a contradiction again.

   (2bii)  M ⊨ (at₃)B(w₃a).

⟹ M ⊨ ∃w₄ (at₃)w₄=w₃a,

⟹ M ⊨ at₃w₄t'auat''aw₂=at₃avat₄,

⟹ M ⊨ w₁at'auat''aw₂=x=w'at₃avat₄=w'at₃w₄t'auat''aw₂,

⟹ by (3.7),  M ⊨ w₁a=w'at₃w₄,

⟹ by (QT4) and (QT5), M ⊨ w₄=a v aEw₄,

⟹  M ⊨ w₁a=w'at₃a v ∃w₆ w₁a=w'at₃(w₆a),

⟹  M ⊨ w₁=w'at₃ v ∃w₆ w₁=w'at₃w₆,

⟹ M ⊨ t₃⊆$_p$w₁.

But from  M ⊨ Intf(x,w₁,t',aua,t'')  we have  M ⊨ Max⁺T$_b$(t',w₁).

Hence M ⊨ t₃<t'. But this contradicts the hypothesis  M ⊨ t'≤t₃  because

$M \vDash t' \in J \subseteq I_0$.

Therefore $M \vDash \neg(at_3)B(w_3a)$.

(2c)  $M \vDash (at_3avat_4)E(at'auat''aw_2)$.

$\Rightarrow$ $M \vDash \exists x_3\ at'auat''aw_2 = x_3at_3avat_4$,

$\Rightarrow$ by (QT4) and (QT5), $M \vDash x_3 = a \lor aBx_3$,

$\Rightarrow$ $M \vDash aat_3avat_4 = at'auat''aw_2 \lor \exists x_4\ at'auat''aw_2 = (ax_4)at_3avat_4$,

$\Rightarrow$ $M \vDash at_3avat_4 = t'auat''aw_2 \lor t'auat''aw_2 = x_4at_3avat_4$.

Now, $M \vDash at_3avat_4 = t'auat''aw_2$ is ruled out because $M \vDash Tally_b(t')$.

So $M \vDash z = t'auat''aw_2 = x_4at_3avat_4$, whence $M \vDash (at_3avat_4)Ez$.

Assume $M \vDash s \subseteq_p x_4\ \&\ Tally_b(s)$.

By hypothesis (2), $M \vDash w'at_3avat_4 = x = w_1az = w_1ax_4at_3avat_4$,

$\Rightarrow$ by (3.6), $M \vDash w' = w_1ax_4$,

$\Rightarrow$ $M \vDash x_4 \subseteq_p w'$,

$\Rightarrow$ $M \vDash s \subseteq_p w'$,

$\Rightarrow$ by hypothesis (2), $M \vDash s < t_3$.

So we have shown that $M \vDash Max^+T_b(t_3, x_4)$.

We thus have

$\qquad M \vDash Pref(ava, t_3)\ \&\ Tally_b(t_4)\ \&\ t_3 = t_4\ \&\ wat_3avat_4 = z\ \&\ Max^+T_b(t_3, w)$

where $M \vDash w = x_4$. But then $M \vDash Lastf(z, t_3, ava, t_4)$, which completes the proof of (5.43).

(5.44) For any string concept I⊆I₀ there is a string concept J⊆I such that

$QT^+ \vdash \forall x \in J \, \forall t,z,t',u,t'',w_1,w_2(Env(t,x) \& Intf(x,w_1,t',aua,t'')$ &

$\& \; x=w_1at'auat''aw_2 \; \& \; z=t'auat''aw_2 \; \rightarrow \; Env(t,z)).$

Let $J \equiv I_{5.25} \& I_{5.43}$.

Assume $M \vDash Env(t,x)$, and let

$M \vDash Intf(x,w_1,t',aua,t'') \; \& \; x=w_1at'auat''aw_2 \; \& \; z=t'auat''aw_2$

where $M \vDash J(x)$.

Assume $M \vDash s \subseteq_p z \; \& \; Tally_b(s)$.

$\Rightarrow M \vDash s \subseteq_p x$,

$\Rightarrow$ from $M \vDash Env(t,x)$, $M \vDash MaxT_b(t,x)$,

$\Rightarrow M \vDash s \subseteq_p t$.

Hence we have (a) $M \vDash MaxT_b(t,z)$.

We now show (b) $M \vDash Firstf(z,t',aua,t'')$.

From $M \vDash Intf(x,w_1,t',aua,t'')$ we have $M \vDash Pref(aua,t') \; \& \; Tally_b(t'') \; \& \; t'<t''$

whereas $M \vDash (t'auat''a)Bz$. But then $M \vDash Firstf(z,t',aua,t'')$, as required.

From $M \vDash Env(t,x)$ we have that $M \vDash \exists v_0,t_1,t_2 \, Firstf(x,t_1,av_0a,t_2)$ and

$M \vDash \exists v,t_3,t_4 \, Lastf(x,t_3,ava,t_4)$. Then $M \vDash t_3=t_4=t$, and from (a) above we have

that $M \vDash t' \leq t_3$. Along with the hypothesis $M \vDash Intf(x,w_1,t',aua,t'')$ and with

$M \vDash z=t'auat''aw_2$ we have from (5.43) that

(c) $M \vDash Lastf(z,t_3,ava,t_4)$.

We now move on to prove

(d)  $M \models \forall v, t_4, t_5$ (Fr(z,$t_4$,ava,$t_5$) & Fr(z,$t_6$,ava,$t_7$) $\rightarrow$ $t_4$=$t_6$).

Assume  $M \models$ Fr(z,$t_4$,ava,$t_5$).

By (b) above we have  $M \models$ Firstf(z,t',aua,t'')  where  $M \models w_1$at'auat''a$w_2$=$w_1$az.

$\Rightarrow$ by (5.20),  $M \models$ (Firstf(z,$t_4$,ava,$t_5$) & $t_4$=t') v  t'<$t_4$,

$\Rightarrow$ $M \models$ t'$\leq t_4$,

$\Rightarrow$ from  $M \models$ Firstf(x,$t_1$,a$v_0$a,$t_2$),  $M \models$ ($t_1v_0$)Bx & Pref(a$v_0$a,$t_1$),

$\Rightarrow$ since  $M \models$ Tally$_b$($t_1$), by (4.23$^b$),  $M \models w_1$=$t_1$ v $t_1$B$w_1$,

$\Rightarrow$ $M \models$ ($w_1$=$t_1$ & $t_1$az=x)  v  ($w_1 \neq t_1$ & $\exists w_3$($t_1 w_3$)az=x).

Suppose, for a reductio, that  $M \models w_1$=$t_1$ & $t_1$az=x.

$\Rightarrow$ $M \models t_1$az=x=$t_1$at'auat''a$w_2$,

$\Rightarrow$ from  $M \models (t_1 a v_0 a)$Bx,  $M \models \exists x_1 (t_1 a v_0 a)x_1$=x=$t_1$at'auat''a$w_2$,

$\Rightarrow$ by (3.7),  $M \models v_0 a x_1$=t'auat''a$w_2$,

$\Rightarrow$ since  $M \models$ Tally$_b$(t'), by (4.23$^b$),  $M \models v_0$=t' v t'B$v_0$,

$\Rightarrow$ $M \models$ t'$\subseteq_p v_0$,

$\Rightarrow$ from hypothesis  $M \models$ Intf(x,$w_1$,t',aua,t''),  $M \models$ Fr(x,t',aua,t''),

$\Rightarrow$ from  $M \models$ Firstf(x,$t_1$,a$v_0$a,$t_2$), by (5.20),  $M \models t_1 \leq$ t',

$\Rightarrow$ $M \models t_1 \subseteq_p$ t'$\subseteq_p v_0$, contradicting  $M \models$ Pref(a$v_0$a,$t_1$).

Therefore  $M \models \neg$($w_1$=$t_1$ & $t_1$az=x).

So the remaining possibility is $M \models w_1 \neq t_1$ & $\exists w_3(t_1 w_3)$az=x.

$\Rightarrow$ from  $M \models (t_1 a v_0 a)$Bx,  $M \models \exists x_1\ t_1 w_3$az=x=$t_1 a v_0 a x_1$,

$\Rightarrow$ by (3.7),  $M \models w_3$az=a$v_0$a$x_1$,

$\Rightarrow$ by (QT4) and (QT5),  $M \models w_3$=a v a$Bw_3$,

$\Rightarrow$ either way, $M \vDash \exists w(t_1wz=x \mathbin{\&} aBw \mathbin{\&} aEw)$ where $M \vDash w=w_3a$.

We now argue that $M \vDash \text{Max}^+T_b(t',t_1w)$.

Assume $M \vDash s\subseteq_p t_1w \mathbin{\&} \text{Tally}_b(s)$.

$\Rightarrow M \vDash s\subseteq_p t_1w = t_1w_3a$,

   where $M \vDash t_1w_3az = t_1w_3at'auat''aw_2 = x = w_1at'auat''aw_2$,

$\Rightarrow$ by (3.7), $M \vDash t_1w_3 = w_1$,

$\Rightarrow M \vDash s\subseteq_p w_1$,

$\Rightarrow$ from $M \vDash \text{Intf}(x,w_1,t',aua,t'')$, $M \vDash \text{Max}^+T_b(t',w_1)$,

$\Rightarrow M \vDash s\subseteq_p t'$.

Hence $M \vDash \text{Max}^+T_b(t',t_1w)$.

But from $M \vDash t'\leq t_4$, we then have $M \vDash \text{Max}^+T_b(t_4,t_1w)$.

$\Rightarrow$ from $M \vDash \text{Fr}(z,t_4,ava,t_5) \mathbin{\&} t_1wz=x \mathbin{\&} aBw \mathbin{\&} aEw$, by (5.25),

$$M \vDash \text{Fr}(x,t_4,ava,t_5).$$

An exactly analogous argument shows that $M \vDash \text{Fr}(x,t_6,ava,t_7)$.

Then from (d) of $M \vDash \text{Env}(t,x)$, it follows that $M \vDash t_4=t_6$, as required.

It remains to show that

(e) $M \vDash \forall u_1,u_2,t_3,t_4,t_5\ (\text{Fr}(z,t_3,au_1a,t_4) \mathbin{\&} \text{Fr}(z,t_3,au_2a,t_5) \rightarrow u_1=u_2)$.

Assume $M \vDash \text{Fr}(z,t_3,au_1a,t_4) \mathbin{\&} \text{Fr}(z,t_3,au_2a,t_5)$.

We reason as in the proof of (d) above and establish that

$$M \vDash \text{Fr}(x,t_3,au_1a,t_4) \mathbin{\&} \text{Fr}(x,t_3,au_2a,t_5),$$

whence $M \vDash u_1=u_2$ follows from (e) of $M \vDash \text{Env}(t,x)$.

Having proved (a)-(e), we have shown that $M \vDash \text{Env}(t,z)$.

This completes the proof of (5.44).

(5.45) For any string concept $I \subseteq I_0$ there is a string concept $J \subseteq I$ such that

$QT^+ \vdash \forall x \in J \; \forall t',t'',t,t_1,t_2,t_3,t_4,u,v,v',w',w_1,z,x_1(Env(t,x) \;\&\; x=t_1auat_2ax_1 \;\&$

$\&\; Firstf(x,t_1,aua,t_2) \;\&\; Intf(x,w',t',ava,t'') \;\&\; Intf(x,w_1,t_3,av'a,t_4) \;\&$

$\&\; w'at'=t_1auat_2 \;\&\; z=t_2ax_1 \;\rightarrow\; Firstf(z,t',ava,t'') \;\&\; \exists t_5 \; Fr(z,t_3,av'a,t_5))$.

Let $J \equiv I_{4.23b} \;\&\; I_{5.28} \;\&\; I_{5.34}$.

Assume

$M \vDash Firstf(x,t_1,aua,t_2) \;\&\; Intf(x,w',t',ava,t'') \;\&\; Intf(x,w_1,t_3,av'a,t_4) \;\&\; Env(t,x)$

where $M \vDash x=t_1auat_2ax_1 \;\&\; z=t_2ax_1 \;\&\; w'at'=t_1auat_2$ and $M \vDash J(x)$.

$\Rightarrow M \vDash Pref(v,t') \;\&\; Tally_b(t'') \;\&\; t'<t'' \;\&\; \exists w'' \; x=w'at'avat''aw'' \;\&\; Max^+T_b(t',w')$,

$\Rightarrow M \vDash Max^+T_b(t',v)$ and $M \vDash x=w'at'ax_1$,

$\Rightarrow$ from $M \vDash w'at'=t_1auat_2$ by (4.24$^b$), $M \vDash t'=t_2$,

$\Rightarrow M \vDash w'az=w'a(t_2ax_1)=w'at'ax_1=x=w'at'avat''aw''$,

$\Rightarrow$ by (3.7), $M \vDash z=t'avat''aw''$,

$\Rightarrow M \vDash Pref(v,t') \;\&\; Tally_b(t'') \;\&\; t'<t'' \;\&\; (t'avat''a)Bz$,

$\Rightarrow M \vDash Firstf(z,t',ava,t'')$.

Note that from the hypothesis $M \vDash Intf(x,w_1,t_3,av'a,t_4)$ we have that

$M \vDash Pref(v',t_3) \;\&\; Tally_b(t_4) \;\&\; t_3<t_4 \;\&\; \exists w_2 \; x=w_1at_3av'at_4aw_2 \;\&\; Max^+T_b(t_3,w_1)$.

We distinguish two cases:

(1) $M \vDash v=v'$.

$\Rightarrow$ from the hypothesis, $M \vDash Fr(x,t',ava,t'') \;\&\; Fr(x,t_3,ava,t_4)$,

$\Rightarrow$ from $M \vDash Env(t,x)$, $M \vDash t'=t_3$,

$\Rightarrow$ from $M \vDash \text{Firstf}(z,t',ava,t'')$, $M \vDash \text{Firstf}(z,t_3,av'a,t'')$,

$\Rightarrow$ $M \vDash \exists t_5\ \text{Fr}(z, t_3, av'a, t_5)$, as required.

(2) $M \vDash v \neq v'$.

Two subcases:

(2a) $M \vDash t''=t_3$.

$\Rightarrow$ from $M \vDash \text{Intf}(x,w',t',ava,t'')\ \&\ \text{Intf}(x,w_1,t_3,av'a,t_4)$ by (5.28),

$$M \vDash w'at'avat'' = w_1 at_3,$$

$\Rightarrow$ by (4.24$^b$), $M \vDash t''=t_3$,

$\Rightarrow$ by (3.6), $M \vDash w'at'av = w_1$,

$\Rightarrow$ from $M \vDash w'at'avat'' = w_1 at_3$

and $M \vDash w'at'avat''aw'' = x = w_1 at_3 av'at_4 aw_2$ by (3.7), $M \vDash w'' = v'at_4 aw_2$,

$\Rightarrow$ from $M \vDash t''=t_3$, $M \vDash t''aw'' = x = t_3 av'at_4 aw_2$,

$\Rightarrow$ $M \vDash z = t'avat''aw'' = (t'av)at_3 av'at_4 aw_2$.

We now argue that $M \vDash \text{Max}^+T_b(t'',t'av)$.

Assume $M \vDash \text{Tally}_b(t_0)\ \&\ t_0 \subseteq_p t'av$.

$\Rightarrow$ by (4.17$^b$), $M \vDash t_0 \subseteq_p t'\ \vee\ t_0 \subseteq_p v$,

$\Rightarrow$ from $M \vDash \text{Max}^+T_b(t',v)$, $M \vDash t_0 \leq t' < t''$.

Therefore, we have that

$$M \vDash \text{Pref}(av'a,t_3)\ \&\ \text{Tally}_b(t_4)\ \&\ t_3 < t_4\ \&$$
$$\&\ \exists w_2\ z = (t'av)at_3 av'at_4 aw_2\ \&\ \text{Max}^+T_b(t_3,t'av),$$

$\Rightarrow$ $M \vDash \text{Intf}(z,t'av,t_3,av'a,t_4)$,

$\Rightarrow$ $M \vDash \text{Fr}(z,t_3,av'a,t_4)$, as required.

(2b)  $M \vDash t''{\neq}t_3$.

  (2bi)  $M \vDash t_4{=}t'$.

$\Rightarrow$ from  $M \vDash \text{Intf}(x,w',t',ava,t'')$ & $\text{Intf}(x,w_1,t_3,av'a,t_4)$  by (5.28),

$$M \vDash w_1at_3av'at_4 = w'at',$$

$\Rightarrow$ from hypothesis  $M \vDash w'at' = t_1auat_2$, $M \vDash w_1at_3av'at_4 = t_1auat_2$,

$\Rightarrow$ by (4.16),  $M \vDash t_3 \subseteq_p u$.

But from  $M \vDash \text{Firstf}(x,t_1,aua,t_2)$ & $\text{Intf}(x,w_1,t_3,av'a,t_4)$  by (5.20),  $M \vDash t_1 {\leq} t_3 \subseteq_p u$.

$\Rightarrow$ $M \vDash t_1 \subseteq_p u$,  contradicting  $M \vDash \text{Firstf}(x,t_1,aua,t_2)$.

  (2bii)  $M \vDash t_4 {\neq} t'$.

Hence we have  $M \vDash v{\neq}v'$ & $t''{\neq}t_3$ & $t_4{\neq}t'$.

$\Rightarrow$ from  $M \vDash \text{Intf}(x,w',t',ava,t'')$ & $\text{Intf}(x,w_1,t_3,av'a,t_4)$  by (5.27),

$$M \vDash (t_3av'at_4) \subseteq_p w' \; \vee \; (t'avat'') \subseteq_p w_1.$$

  (2biia)  $M \vDash (t_3av'at_4) \subseteq_p w'$.

$\Rightarrow$ $M \vDash t_3 \subseteq_p (t_3av'at_4) \subseteq_p w'$,

$\Rightarrow$ from  $M \vDash \text{Max}^+T_b(t',w')$,  $M \vDash t_3 < t' = t_2$.

But from $M \vDash \text{Firstf}(x,t_1,aua,t_2)$ & $\text{Intf}(x,w_1,t_3,av'a,t_4)$, by (5.34),  $M \vDash t_2 {\leq} t_3$.

$\Rightarrow$ $M \vDash t_2 {\leq} t_3 < t_2$,  contradicting  $M \vDash t_2 \in I \subseteq I_0$.

  (2biib)  $M \vDash (t'avat'') \subseteq_p w_1$.

We have two principal subcases to distinguish:

   (2biib1)  $M \vDash (t'avat'') = w_1 \vee (t'avat'')Bw_1$.

$\Rightarrow$ $M \vDash (t'a)Bx$,

$\Rightarrow$ from  $M \vDash (t_1a)Bx$ & $(t'a)Bx$,

$\Rightarrow$ $M \vDash \exists x_3, x_4\ t_1 a x_3 = t'ax_4$,

$\Rightarrow$ by (4.23$^b$), $M \vDash t_1 = t'$,

$\Rightarrow$ $M \vDash t_1 < t_2 = t' = t_1$, contradicting $M \vDash t_1 \in I \subseteq I_0$.

(2biib2)  $M \vDash (t'avat'')Ew_1$.

$\Rightarrow$ $M \vDash \exists w_3\ w_3(t'avat'') = w_1$,

$\Rightarrow$ $M \vDash (w_3 t'avat'')at_3 av'at_4 aw_2 = x = w'at'avat''aw''$,

$\Rightarrow$ $M \vDash w_3 B x\ \&\ (w'a)Bx$,

$\Rightarrow$ by (3.8), $M \vDash w_3 B(w'a)\ \vee\ w_3 = w'a\ \vee\ (w'a)Bw_3$.

(2biib2a)  $M \vDash w_3 B(w'a)$.

$\Rightarrow$ $M \vDash \exists w_5\ w_3 w_5 = w'a$,

$\Rightarrow$ $M \vDash (w_3 t'avat'')at_3 av'at_4 aw_2 = (w_3 w_5)t'avat''aw''$,

$\Rightarrow$ by (3.7), $M \vDash t'avat''at_3 av'at_4 aw_2 = w_5 t'avat''aw''$,

$\Rightarrow$ from $M \vDash w_3 w_5 = w'a$, $M \vDash w_5 = a\ \vee\ aEw_5$.

Now, $M \vDash w_5 \neq a$ because $M \vDash \text{Tally}_b(t')$.

$\Rightarrow$ $M \vDash aEw_5$,

$\Rightarrow$ $M \vDash \exists w_6\ (w_6 a = w_5\ \&\ t'avat''at_3 av'at_4 aw_2 = (w_6 a)t'avat''aw'')$,

$\Rightarrow$ by (4.14$^b$), $M \vDash w_6 = t'\ \vee\ t'Bw_6$,

$\Rightarrow$ $M \vDash t' \subseteq_p w_6$.

But $M \vDash w_3(w_6 a) = w_3 w_5 = w'a$. So $M \vDash w_3 w_6 = w'$, whence $M \vDash t' \subseteq_p w_6 \subseteq_p w'$, contradicting $M \vDash \text{Max}^+ T_b(t',v)$.

(2biib2b)  $M \vDash (w'a)Bw_3$.

$\Rightarrow$ $M \vDash \exists w_6\ (w'a)w_6 = w_3$,

$\Rightarrow$ $M \vDash (w'aw_6)t'avat''at_3av'at_4aw_2 = w'at'avat''aw'' = x = w_1at_3av'at_4aw_2$,

$\Rightarrow$ by (3.6), $M \vDash w'aw_6t'avat'' = w_1$,

$\Rightarrow$ from $M \vDash Max^+T_b(t_3,w_1)$, $M \vDash Max^+T_b(t_3,w'aw_6t'avat'')$,

$\Rightarrow$ $M \vDash Max^+T_b(t_3,w_6t'avat'')$,

$\Rightarrow$ by (3.7), $M \vDash w_6t'avat''at_3av'at_4aw_2 = t'avat''aw'' = z$.

But then we have

$\quad M \vDash Pref(av'a,t_3)$ & $Tally_b(t_4)$ & $t_3<t_4$ &

$\qquad\qquad$ & $\exists w_2\ z = (w_6t'avat'')at_3av'at_4aw_2$ & $Max^+T_b(t_3,w_6t'avat'')$,

$\Rightarrow$ $M \vDash Intf(z,w_6t'avat'',t_3,av'a,t_4)$,

$\Rightarrow$ $M \vDash Fr(z,t_3,av'a,t_4)$, as required.

$\quad$ (2biib2c) $M \vDash w_3 = w'a$.

$\Rightarrow$ $M \vDash (w'a)t'avat''at_3av'at_4aw_2 = x = w'at'avat''aw''$,

$\Rightarrow$ by (3.7), $M \vDash t'avat''at_3av'at_4aw_2 = t'avat''aw'' = z$.

Also, from $M \vDash w_3t'avat'' = w_1$, $M \vDash w'at'avat'' = w_1$,

$\Rightarrow$ from $M \vDash Max^+T_b(t_3,w_1)$, $M \vDash Max^+T_b(t_3,w'at'avat'')$,

$\Rightarrow$ $M \vDash Max^+T_b(t_3,t'avat'')$.

But then we have

$\quad M \vDash Pref(av'a,t_3)$ & $Tally_b(t_4)$ & $t_3<t_4$ &

$\qquad\qquad$ & $\exists w_2\ z = (t'avat'')at_3av'at_4aw_2$ & $Max^+T_b(t_3,t'avat'')$,

$\Rightarrow$ $M \vDash Intf(z,t'avat'',t_3,av'a,t_4)$,

$\Rightarrow$ $M \vDash Fr(z,t_3,av'a,t_4)$.

$\quad$ (2biib3) $M \vDash \exists w_3,w_4\ w_3(t'avat'')w_4 = w_1$.

$\Rightarrow$ M ⊨ ($w_3$t'avat''$w_4$)at$_3$av'at$_4$a$w_2$=$w_1$at$_3$av'at$_4$a$w_2$=x=w'at'avat''aw'',

$\Rightarrow$ M ⊨ $w_3$Bx & (w'a)Bx,

$\Rightarrow$ by (3.8),  M ⊨ $w_3$B(w'a) v $w_3$=w'a v (w'a)B$w_3$.

 (2biib3a)  M ⊨ $w_3$B(w'a).

The same argument as in (2biib2a).

 (2biib3b)  M ⊨ (w'a)B$w_3$.

Analogous argument to (2biib2b) with $w_6$t'av'at''$w_4$ in place of $w_6$t'avat'' throughout.

 (2biib3c)  M ⊨ $w_3$=w'a.

Analogous to (2biib2c) with w'at'av'at''$w_4$ in place of w'at'avat'' throughout. This completes the proof of (5.45).

(5.46)  For any string concept I⊆I₀ there is a string  concept J⊆I such that

$\quad$ QT⁺ ⊢ ∀x,y∈J ∀t,t₂,t₃(Env(t₂,x) & Env(t,y) & (t₃a)By & Tally_b(t₃) & t₂<t₃ &

$\quad\quad$ & ¬∃u(u ε x  &  u ε y)  → ∃z∈J (Env(t,z) & ∀u(u ε z ↔ u ε x v u ε y)).

Let J ≡ I₅.₁₁ & I₅.₃₉ & I₅.₄₀.

Assume  M ⊨ Env(t₂,x) & Env(t,y) & ¬∃u(u ε x  &  u ε y)  along with

M ⊨ (t₃a)By & Tally_b(t₃) & t₂<t₃  and  M ⊨ J(x) & J(y).

⟹ by (5.11),  M ⊨ ∃w,t₁(Tally_b(t₁) & x=t₁wt₂ & aBw & aEw),

⟹ from M ⊨ Env(t,y),  M ⊨ ∃v,t₅ Firstf(y,t₅,ava,t₄),

⟹ M ⊨ Tally_b(t₅) & (t₅a)By,

⟹ by (4.23^b),  M ⊨ t₃=t₅,

⟹ M ⊨ Firstf(y,t₃,ava,t₄).

Let  z=t₁wy.

Since we may assume that J is closed under ⊆_p and *, we have that  M ⊨ J(z).

We first show that

$\quad\quad\quad$ M ⊨ ∃z',t₀(Tally_b(t₀) & z=xt₀z't & aBz' & aEz').

From  M ⊨ t₂<t₃ & Tally_b(t₃),  M ⊨ ∃t₀(Tally_b(t₀) & t₃=t₂t₀).

From  M ⊨ Firstf(y,t₃,ava,t₄),  we have

 M ⊨ Pref(ava,t₃) & Tally_b(t₄) & ((t₃=t₄ & y=t₃avat₄) v (t₃<t₄ & (t₃avat₄a)By)),

⟹  M ⊨ z=t₁wt₃avat₄ v ∃z₁ z=t₁w(t₃avat₄az₁),

⟹  M ⊨ z=t₁w(t₂t₀)avat₄ v ∃z₁ z=t₁w(t₂t₀)avat₄az₁.

If  M ⊨ t₃=t₄ & y=t₃avat₄, then  M ⊨ Lastf(y,t₃,ava,t₄), and from  M ⊨ Env(t,y)

we have $M \vDash t_3=t_4=t$.

From $M \vDash Env(t,y)$, $M \vDash (at)Ey$. Hence $M \vDash \exists z_2\ y=z_2at$.

So if $M \vDash z=t_1w(t_2t_0)avat_4az_1$, then $M \vDash (t_1wt_2)t_0avat_4az_1=z=t_1wz_2at$.

$\implies$ by (4.15$^b$), $M \vDash t=z_1\ v\ tEz_1$,

$\implies M \vDash z=xt_0avat_4at\ v\ \exists z_3(z_1=z_3t\ \&\ z=xt_0avat_4az_3t)$.

Then if $M \vDash t_1wz_2at=t_1wy=z=xt_0avat_4az_3t$, then by (3.6),

$$M \vDash t_1wz_2a=xt_0avat_4az_3.$$

$\implies M \vDash z_3=a\ v\ aEz_3$.

Therefore, we have that

$M \vDash z=xt_0(ava)t_4\ v\ z=xt_0(ava)t_4at\ v$

$$v\ z=xt_0(ava)t_4aat\ v\ \exists z_4\ z=xt_0(ava)t_4a(z_4a)t.$$

$\implies$ by (5.40), from $M \vDash Env(t,y)\ \&\ Firstf(y,t_3,ava,t_4)$,

$$M \vDash y \neq t_3avat_4at\ \&\ y \neq t_3avat_4aat,$$

$\implies M \vDash t_1wy=z \neq t_1wt_3avat_4at=t_1w(t_2t_0)avat_4at=xt_0(ava)t_4at$ and

$\qquad M \vDash t_1wy=z \neq t_1wt_3avat_4aat=t_1w(t_2t_0)avat_4aat=xt_0(ava)t_4aat$,

$\implies M \vDash z=xt_0(ava)t_4\ v\ \exists z_4\ z=xt_0(ava)t_4a(z_4a)t$,

and so either way we have $M \vDash \exists z',t_0(Tally_b(t_0)\ \&\ z=xt_0z't\ \&\ aBz'\ \&\ aEz')$, as claimed.

Now, from $M \vDash Env(t_2,x)$, $M \vDash \exists v',t_5,t_6\ Firstf(x,t_5,av'a,t_6)$.

$\implies$ from $M \vDash Tally_b(t_0)\ \&\ Tally_b(t)$, by the proof of (5.6)(1),

(1b) $\quad M \vDash \exists v',t_5,t_6\ Firstf(z,t_5,av'a,t_6)$.

We show that

(1a)   $M \vDash MaxT_b(t,z)$.

Assume  $M \vDash Tally_b(t_0)$ & $t_0 \subseteq_p z = t_1wy = t_1w_2ay$.

$\Rightarrow$ by (4.17$^b$),  $M \vDash t_0 \subseteq_p t_1w_2$  v  $t_0 \subseteq_p y$.

We have that  $M \vDash Max^+T_b(t,t_1w_2)$  as in the proof of (5.39)(2). On the other hand, from $M \vDash Env(t,y)$, $M \vDash MaxT_b(t,y)$.  Then $M \vDash t_0 \leq t$, as required.

For   (1c)  $M \vDash \exists v''\, Lastf(z,t,av''a,t)$,

we observe, as in the proof of (5.39),(3), that  $M \vDash Max^+T_b(t,t_1w)$.  Then $M \vDash \exists v''\, Lastf(z,t,av''a,t)$  follows by the argument given there.

For  (1d)  $M \vDash \forall u,t_5,t_6,t_7,t_8(Fr(z,t_5,aua,t_6)$ & $Fr(z,t_7,aua,t_8) \to t_5 = t_7)$,

assume  $M \vDash Fr(z,t_5,aua,t_6)$ & $Fr(z,t_7,aua,t_8)$.

We distinguish nine cases, by applying (5.39):

 (1di)  $M \vDash \exists t_9\, Fr(x,t_5,aua,t_9)$ & $\exists t_{10} Fr(x,t_7,aua,t_{10})$.

$\Rightarrow$ from  $M \vDash Env(t_2,x)$,  $M \vDash t_5 = t_7$.

 (1dii-iii)  $M \vDash \exists t_9\, Fr(x,t_5,aua,t_9)$ & $((t_7 = t_3$ & $Fr(y,t_3,aua,t_4))$ v

$\qquad\qquad\qquad\qquad\qquad\qquad$ v $\exists t_{10} Fr(y,t_7,aua,t_{10}))$.

$\Rightarrow$ $M \vDash u\,\varepsilon\, x$ & $u\,\varepsilon\, y$,  contradicting the hypothesis.

 (1div)  $M \vDash (t_5 = t_3$ & $Fr(y,t_3,aua,t_4))$ & $\exists t_{10} Fr(x,t_7,aua,t_{10})$.

Same as (1dii-iii).

 (1dv)  $M \vDash (t_5 = t_3$ & $Fr(y,t_3,aua,t_4))$ & $(t_7 = t_3$ & $Fr(y,t_3,aua,t_4))$.

Nothing to prove.

 (1dvi)  $M \vDash (t_5 = t_3$ & $Fr(y,t_3,aua,t_4))$ & $\exists t_{10} Fr(y,t_7,aua,t_{10})$.

$\Rightarrow$ from  $M \vDash Env(t,y)$,  $M \vDash t_3 = t_7$.

(1dvii) $M \vDash \exists t_9\, Fr(y,t_5,aua,t_9)\ \&\ \exists t_{10} Fr(x,t_7,aua,t_{10})$.

Same as (1dii-iii).

(1dviii) $M \vDash \exists t_9\, Fr(y,t_5,aua,t_9)\ \&\ (t_7=t_3\ \&\ Fr(y,t_3,aua,t_4))$.

$\Rightarrow$ from $M \vDash Env(t,y)$, $M \vDash t_5=t_3=t_7$.

(1dix) $M \vDash \exists t_9\, Fr(y,t_5,aua,t_9)\ \&\ \exists t_{10} Fr(y,t_7,aua,t_{10})$.

$\Rightarrow$ from $M \vDash Env(t,y)$, $M \vDash t_5=t_7$.

This proves (1d).

For (1e) $M \vDash \forall u,u',t_5,t_6,t_7(Fr(z,t_5,aua,t_6)\ \&\ Fr(z,t_5,au'a,t_7) \rightarrow u=u')$,

assume $M \vDash Fr(z,t_5,aua,t_6)\ \&\ Fr(z,t_5,au'a,t_7)$.

Again applying (5.39) we have:

(1ei) $M \vDash \exists t_8\, Fr(x,t_5,aua,t_8)\ \&\ \exists t_9 Fr(x,t_5,au'a,t_9)$.

$\Rightarrow$ from $M \vDash Env(t_2,x)$, $M \vDash u=u'$.

(1eii) $M \vDash \exists t_8\, Fr(x,t_5,aua,t_8)\ \&\ ((t_5=t_3\ \&\ Fr(y,t_3,au'a,t_4))$.

$\Rightarrow$ from $M \vDash Env(t_2,x)$, $M \vDash t_5 \leq t_2$,

$\Rightarrow M \vDash t_5 \leq t_2 < t_3 = t_5$, contradicting $M \vDash t_5 \in I \subseteq I_0$.

(1eiii) $M \vDash \exists t_8\, Fr(x,t_5,aua,t_8)\ \&\ \exists t_9\, Fr(y,t_5,au'a,t_9)$.

By (5.20), from $M \vDash Firstf(y,t_3,ava,t_4)$, $M \vDash t_3 \leq t_5$. From $M \vDash Env(t_2,x)$,

$M \vDash t_5 \leq t_2$. Hence $M \vDash t_5 \leq t_2 < t_3 \leq t_5$, again contradicting $M \vDash t_5 \in I \subseteq I_0$.

(1eiv) $M \vDash (t_5=t_3\ \&\ Fr(y,t_3,aua,t_4))\ \&\ \exists t_9 Fr(x,t_5,au'a,t_9)$.

Just like (1eii).

(1ev-vi) $M \vDash (t_5=t_3\ \&\ Fr(y,t_3,aua,t_4))\ \&\ (Fr(y,t_3,au'a,t_4)\ \vee\ \exists t_9 Fr(x,t_5,au'a,t_9))$.

$\Rightarrow$ from $M \vDash Env(t,z)$, $M \vDash u=u'$.

(1evii)  $M \vDash \exists t_8\, Fr(y,t_5,aua,t_8)\ \&\ \exists t_9\, Fr(x,t_5,au'a,t_9)$.

Same as (1eiii).

(1eviii)  $M \vDash \exists t_8 Fr(y,t_5,aua,t_8)\ \&\ (t_5=t_3\ \&\ Fr(y,t_3,au'a,t_4))$.

Same as (1ev-vi).

(1eix)  $M \vDash \exists t_8\, Fr(y,t_5,aua,t_8)\ \&\ \exists t_9 Fr(y,t_5,au'a,t_9)$.

$\Rightarrow$ from $M \vDash Env(t,y)$, $M \vDash u=u'$.

This proves (1e).

We have thus established that $M \vDash Env(t,z)$.

We now proceed to show that $M \vDash \forall u(u\ \varepsilon\ z \to u\ \varepsilon\ x\ v\ u\ \varepsilon\ y)$.

Assume $M \vDash u\ \varepsilon\ z$.

$\Rightarrow M \vDash \exists t',t''\, Fr(z,t',aua,t'')$,

$\Rightarrow$ by (5.39), $M \vDash \exists t_5,t_6(Fr(x,t_5,aua,t_6)\ v\ Fr(y,t_5,aua,t_6))$,

$\Rightarrow$ since $M \vDash Env(t_2,x)\ \&\ Env(t,y)$, $M \vDash u\ \varepsilon\ x\ v\ u\ \varepsilon\ y$, as required.

Finally, we show that $M \vDash \forall u(u\ \varepsilon\ x\ v\ u\ \varepsilon\ y \to u\ \varepsilon\ z)$.

Assume $M \vDash u\ \varepsilon\ y$.

$\Rightarrow M \vDash \exists t',t''\, Fr(y,t',aua,t'')$.

We claim that $M \vDash Max^+T_b(t',t_1w)$.

By hypothesis we have $M \vDash Firstf(y,t_3,ava,t_4)$.

$\Rightarrow$ by (5.20), $M \vDash t_3 \leq t'$,

$\Rightarrow$ from $M \vDash Env(t_2,x)$, $M \vDash MaxT_b(t_2,x)\ \&\ t_2<t_3\leq t'$,

$\Rightarrow M \vDash Max^+T_b(t',x)$,

$\Rightarrow$ since $M \vDash t_1w \subseteq_p x$, $M \vDash Max^+T_b(t',t_1w)$.

Now, we have   $M \vDash Fr(z,t',aua,t'')$ & $z=t_1wy$ & $aBw$ & $aEw$ & $Max^+T_b(t',t_1w)$.

$\implies$ by (5.25),  $M \vDash Fr(z,t',aua,t'')$,

$\implies M \vDash u\ \varepsilon\ z$,  as required.

On the other hand, assume  $M \vDash u\ \varepsilon\ x$.

$\implies M \vDash \exists t',t''\ Fr(x,t',aua,t'')$.

Hence we have

$M \vDash Env(t_2,x)$ & $Fr(x,t',aua,t'')$ & $xt_0z't=z$ & $aBz'$ & $aEz'$ & $Tally_b(t_0)$ & $Tally_b(t)$.

$\implies$ by (5.6),  $M \vDash \exists t_6\ Fr(z,t',aua,t_6)$,

$\implies M \vDash u\ \varepsilon\ z$,  as required.

Therefore, we have shown that  $M \vDash \forall u(u\ \varepsilon\ x\ \vee\ u\ \varepsilon\ y\ \rightarrow\ u\ \varepsilon\ z)$.

This completes the proof of (5.46).

(5.47) For any string concept $I\subseteq I_0$ there is a string concept $J\subseteq I$ such that

$QT^+ \vdash \forall x \in J \ \forall w_1,t_1,t_2,u(Intf(x,w_1,t_1,aua,t_2) \rightarrow$

$\rightarrow \forall z,w_2,t_3,v,t_4( x=w_1at_1auat_2aw_2 \ \& \ z= w_1at_2aw_2 \ \& \ Fr(z,t_3,ava,t_4) \rightarrow$

$\rightarrow Fr(x,t_3,ava,t_4) \ v \ Fr(z,t_3,ava,t_1))]$.

Let $J \equiv I_{3.10} \ \& \ I_{4.16} \ \& \ I_{4.17b} \ \& \ I_{4.23b} \ \& \ I_{4.24b}$.

Assume $M \vDash Intf(x,w_1,t_1,aua,t_2)$ where $M \vDash x=w_1at_1auat_2aw_2$ and $M \vDash J(x)$.

Let $M \vDash z=w_1at_2aw_2$. Note that z is obtained from x by excising $t_1aua$.

Assume $M \vDash Fr(z,t_3,ava,t_4)$.

There are three cases:

(1) $M \vDash Firstf(z,t_3,ava,t_4)$.

$\Rightarrow M \vDash Pref(ava,t_3) \ \& \ Tally_b(t_4) \ \&$

$\& \ ((t_3=t_4 \ \& \ z=t_3avat_4) \ v \ (t_3<t_4 \ \& \ (t_3avat_4a)Bz))$.

(1a) $M \vDash t_3=t_4 \ \& \ z=t_3avat_4$.

$\Rightarrow M \vDash w_1at_2aw_2=z=t_3avat_4$,

$\Rightarrow$ by (4.16), since $M \vDash Tally_b(t_3) \ \& \ Tally_b(t_4)$, $M \vDash t_2\subseteq_p v$,

$\Rightarrow$ from $M \vDash Pref(ava,t_3)$, $M \vDash t_2<t_3$.

On the other hand, by (4.14$^b$), $M \vDash w_1=t_3 \ v \ t_3Bw_1$,

$\Rightarrow M \vDash t_3\subseteq_p w_1$,

$\Rightarrow$ from hypothesis $M \vDash Intf(x,w_1,t_1,aua,t_2)$, $M \vDash t_3<t_1<t_2$, a contradiction.

Therefore, $M \vDash \neg(t_3=t_4 \ \& \ z=t_3avat_4)$.

(1b) $M \vDash t_3<t_4 \ \& \ (t_3avat_4a)Bz$.

$\Rightarrow$ $M \vDash (t_3avat_4a)Bz$ & $(w_1at_2a)Bz$,

$\Rightarrow$ $M \vDash t_3avat_4a=w_1at_2a$ v $(t_3avat_4a)B(w_1at_2a)$ v $(w_1at_2a)B(t_3avat_4a)$.

(1bi) $M \vDash t_3avat_4a=w_1at_2a$.

$\Rightarrow$ $M \vDash t_3avat_4=w_1at_2$,

$\Rightarrow$ by (4.24$^b$), $M \vDash t_4=t_2$,

$\Rightarrow$ by (3.6), $M \vDash t_3av=w_1$,

$\Rightarrow$ $M \vDash t_3 \subseteq_p w_1$,

$\Rightarrow$ from hypothesis $M \vDash Intf(x,w_1,t_1,aua,t_2)$, $M \vDash Max^+T_b(t_1,w_1)$,

$\Rightarrow$ $M \vDash t_3<t_1$.

Hence we have $M \vDash Pref(ava,t_3)$ & $Tally_b(t_1)$ & $t_3<t_1$ & $(t_3avat_1a)Bx$,

$\Rightarrow$ $M \vDash Firstf(x,t_3,ava,t_1)$,

$\Rightarrow$ $M \vDash Fr(x,t_3,ava,t_1)$, as required.

(1bii) $M \vDash (t_3avat_4a)B(w_1at_2a)$.

$\Rightarrow$ $M \vDash \exists x_1\ t_3avat_4ax_1=w_1at_2a$,

$\Rightarrow$ $M \vDash x_1=a$ v $aEx_1$.

If $M \vDash x_1=a$, then $M \vDash t_3avat_4aa=w_1at_2a$, whence $M \vDash t_3avat_4a=w_1at_2$, a contradiction since $M \vDash Tally_b(t_2)$ from hypothesis $M \vDash Intf(x,w_1,t_1,aua,t_2)$. Therefore, $M \vDash aEx_1$.

$\Rightarrow$ $M \vDash \exists x_2\ t_3avat_4a(x_2a)=w_1at_2a$,

$\Rightarrow$ $M \vDash t_3avat_4ax_2=w_1at_2$,

$\Rightarrow$ by (4.15$^b$), $M \vDash x_2=t_2$ v $t_2Ex_2$,

$\Rightarrow$ $M \vDash t_3avat_4at_2=w_1at_2$ v $\exists x_3\ t_3avat_4ax_3t_2=w_1at_2$,

$\Rightarrow$ by (3.6), $M \vDash t_3avat_4a=w_1a \lor t_3avat_4ax_3=w_1a$,

$\Rightarrow M \vDash (t_3avat_4a)B(w_1at_1auat_2aw_2)$,

$\Rightarrow M \vDash (t_3avat_4a)Bx$,

$\Rightarrow M \vDash Pref(ava,t_3) \& Tally_b(t_4) \& t_3<t_4 \& (t_3avat_4a)Bx$,

$\Rightarrow M \vDash Firstf(x,t_3,ava,t_4)$,

$\Rightarrow M \vDash Fr(x,t_3,ava,t_4)$, as required.

(1biii) $M \vDash (w_1at_2a)B(t_3avat_4a)$.

$\Rightarrow M \vDash \exists x_1\ t_3avat_4a=w_1at_2ax_1$,

$\Rightarrow M \vDash x_1=a \lor aEx_1$.

If $M \vDash x_1=a$, then $M \vDash t_3avat_4a=w_1at_2aa$, whence $M \vDash t_3avat_4=w_1at_2a$, a contradiction because $M \vDash Tally_b(t_4)$.

Therefore, $M \vDash aEx_1$.

$\Rightarrow M \vDash \exists x_2\ t_3avat_4a=w_1at_2a(x_2a)$,

$\Rightarrow M \vDash t_3avat_4=w_1at_2ax_2$,

$\Rightarrow$ by (4.24$^b$), $M \vDash x_2=t_4 \lor t_4Ex_2$,

$\Rightarrow M \vDash t_3avat_4=w_1at_2at_4 \lor \exists x_3\ t_3avat_4=w_1at_2ax_3t_4$,

$\Rightarrow$ by (4.16), either way $M \vDash t_2\subseteq_p v$,

$\Rightarrow$ from $M \vDash t_3avat_4=w_1at_2ax_2$, by (4.14$^b$), $M \vDash w_1=t_3 \lor t_3Bw_1$,

$\Rightarrow M \vDash t_3\subseteq_p w_1$.

From hypothesis $M \vDash Intf(x,w_1,t_1,aua,t_2)$, $M \vDash Max^+T_b(t_1,w_1) \& t_1<t_2$,

$\Rightarrow M \vDash t_3<t_1<t_2\subseteq_p v$, contradicting $M \vDash Pref(ava,t_3)$.

This completes the proof in Case 1.

(2) $M \vDash \exists w_3\ \text{Intf}(z,w_3,t_3,ava,t_4)$.

$\Longrightarrow$ $M \vDash \exists w_4\ w_1at_2aw_2=z=w_3at_3avat_4aw_4$,

$\Longrightarrow$ $M \vDash (w_1a)Bz\ \&\ (w_3a)Bz$,

$\Longrightarrow$ by (3.8), $M \vDash w_1a=w_3a\ \vee\ (w_1a)B(w_3a)\ \vee\ (w_3a)B(w_1a)$.

(2a) $M \vDash w_1a=w_3a$.

$\Longrightarrow$ $M \vDash w_1at_2aw_2=w_1at_3avat_4aw_4$,

$\Longrightarrow$ by (3.7), $M \vDash t_2aw_2=t_3avat_4aw_4$,

$\Longrightarrow$ $M \vDash w_1at_1auat_2aw_2=x=w_1at_1auat_3avat_4aw_4$,

$\Longrightarrow$ by ($4.23^b$), $M \vDash t_2=t_3$.

We claim that $M \vDash \text{Max}^+T_b(t_3,w_1at_1au)$.

Assume $M \vDash s\subseteq_p(w_1at_1au)\ \&\ \text{Tally}_b(s)$.

$\Longrightarrow$ by ($4.17^b$), $M \vDash s\subseteq_pw_1\ \vee\ s\subseteq_pt_1au$.

If $M \vDash s\subseteq_pw_1$, then $M \vDash s<t_3$ because from $M \vDash \text{Intf}(z,w_3,t_3,ava,t_4)$ we have

$M \vDash \text{Max}^+T_b(t_3,w_3)$ and $M \vDash w_1a=w_3a$.

If $M \vDash s\subseteq_pt_1au$, then, by ($4.17^b$), $M \vDash s\subseteq_pt_1\ \vee\ s\subseteq_pu$.

$\Longrightarrow$ from $M \vDash \text{Intf}(x,w_1,t_1,aua,t_2)$, $M \vDash \text{Pref}(aua,t_1)\ \&\ t_1<t_2=t_3$,

$\Longrightarrow$ $M \vDash s\subseteq_pt_1<t_3$.

Therefore, $M \vDash \text{Max}^+T_b(t_3,w_1at_1au)$, as claimed.

Thus we have

$M \vDash \text{Pref}(ava,t_3)\ \&\ \text{Tally}_b(t_4)\ \&\ t_3<t_4\ \&\ x=w'at_3avat_4aw_4\ \&\ \text{Max}^+T_b(t_3,w')$

where $w'= w_1at_1au$.

$\Longrightarrow M \vDash \exists w'\ \text{Intf}(x,w',t_3,ava,t_4)$,

$\Rightarrow$ $M \vDash Fr(x,t_3,ava,t_4)$.

(2b) $M \vDash (w_1a)B(w_3a)$.

$\Rightarrow$ $M \vDash \exists w_5\ w_1aw_5=w_3a$,

$\Rightarrow$ $M \vDash w_1at_2aw_2=z=(w_1aw_5)t_3avat_4aw_4$,

$\Rightarrow$ by (3.7), $M \vDash t_2aw_2=w_5t_3avat_4aw_4$,

$\Rightarrow$ $M \vDash w_1at_1auat_2aw_2=x=w_1at_1auaw_5t_3avat_4aw_4$,

$\Rightarrow$ from $M \vDash w_1aw_5=w_3a$, $M \vDash w_5=a\ v\ aEw_5$.

We cannot have $M \vDash w_5=a$ because $M \vDash Tally_b(t_2)$.

$\Rightarrow$ $M \vDash \exists w_6\ w_5=w_6a$,

$\Rightarrow$ $M \vDash x=(w_1at_1auaw_6)at_3avat_4aw_4$,

$\Rightarrow$ from $M \vDash w_1aw_6a=w_1aw_5=w_3a$, $M \vDash w_1aw_6=w_3$,

$\Rightarrow$ from $M \vDash t_2aw_2=w_5t_3avat_4aw_4$, $M \vDash t_2aw_2=w_6at_3avat_4aw_4$,

$\Rightarrow$ by $(4.14^b)$, $M \vDash w_6=t_2\ v\ t_2Bw_6$,

$\Rightarrow$ $M \vDash t_2\subseteq_p w_6 \subseteq_p w_3$.

We claim that $M \vDash Max^+T_b(t_3, w_1at_1auaw_6)$.

Assume $M \vDash s\subseteq_p(w_1at_1auaw_6)\ \&\ Tally_b(s)$.

$\Rightarrow$ by $(4.17^b)$, $M \vDash s\subseteq_p w_1\ v\ s\subseteq_p t_1auw_6$.

Note that $M \vDash s\subseteq_p w_1 \subseteq_p w_3$.

On the other hand, if $M \vDash s\subseteq_p t_1auaw_6$, then, by $(4.17^b)$,

$$M \vDash s\subseteq_p t_1au\ v\ s\subseteq_p w_6.$$

Again, note that $M \vDash s\subseteq_p w_6 \subseteq_p w_3$.

And if $M \vDash s\subseteq_p t_1au$, then $M \vDash s\subseteq_p t_1$ as in (2a) above.

But $M \vDash t_1 < t_2 \subseteq_p w_3$.

So we have in case $M \vDash s \subseteq_p w_3$.

But from $M \vDash \text{Intf}(z,w_3,t_3,ava,t_4)$ we have $M \vDash \text{Max}^+T_b(t_3,w_3)$.

Thus $M \vDash s < t_3$. Therefore $M \vDash \text{Max}^+T_b(t_3,w_1at_1auaw_6)$, as claimed.

So we have

$M \vDash \text{Pref}(ava,t_3)$ & $\text{Tally}_b(t_4)$ & $t_3 < t_4$ & $x = w'at_3avat_4aw_4$ & $\text{Max}^+T_b(t_3,w')$

where $w' = w_1at_1auaw_6$.

$\implies M \vDash \exists w' \, \text{Intf}(x,w',t_3,ava,t_4)$,

$\implies M \vDash \text{Fr}(x,t_3,ava,t_4)$.

(2c) $M \vDash (w_3a)B(w_1a)$.

$\implies M \vDash \exists w_5 \, w_3aw_5 = w_1a$,

$\implies M \vDash w_3aw_5t_2aw_2 = z = w_3at_3avat_4aw_4$,

$\implies$ by (3.7), $M \vDash w_5t_2aw_2 = t_3avat_4aw_4$,

$\implies$ from $M \vDash w_3aw_5 = w_1a$, $M \vDash w_5 = a \lor aEw_5$.

Again, we cannot have $M \vDash w_5 = a$ because $M \vDash \text{Tally}_b(t_3)$. So $M \vDash aEw_5$.

$\implies M \vDash \exists w_6 \, w_5 = w_6a$,

$\implies M \vDash w_6at_2aw_2 = t_3avat_4aw_4$,

$\implies$ by (4.14$^b$), $M \vDash w_6 = t_3 \lor t_3Bw_6$,

$\implies M \vDash t_3at_2aw_2 = t_3avat_4aw_4 \lor \exists w_7 \, (t_3w_7)at_2aw_2 = t_3avat_4aw_4$,

$\implies$ by (3.7), $M \vDash t_2aw_2 = vat_4aw_4 \lor \exists w_7 \, w_7at_2aw_2 = avat_4aw_4$.

Note that $M \vDash t_3 \subseteq_p w_6 \subseteq_p w_1$.

$\implies$ from $M \vDash \text{Max}^+T_b(t_1,w_1)$, $M \vDash t_3 < t_1 < t_2$.

(2ci)  $M \vDash t_2aw_2=vat_4aw_4$.

$\implies$ by (4.14$^b$),  $M \vDash v=t_2 \lor t_2Bv$,

$\implies M \vDash t_2\subseteq_p v$,

$\implies M \vDash t_3<t_2\subseteq_p v$,  contradicting  $M \vDash \text{Pref}(ava,t_3)$.

(2cii)  $M \vDash \exists w_7\ w_7at_2aw_2=avat_4aw_4$.

$\implies$ by (3.8), $M \vDash ava=w_7a \lor (ava)B(w_7a) \lor (w_7a)B(ava)$.

(2cii1)  $M \vDash (w_7a)B(ava)$.

$\implies M \vDash \exists w_8\ w_7aw_8=ava$,

$\implies M \vDash w_7at_2aw_2=(w_7aw_8)t_4aw_4$,

$\implies$ by (3.7), $M \vDash t_2aw_2=w_8t_4aw_4$,

$\implies$ from $M \vDash w_7aw_8=ava$,  $M \vDash w_8=a \lor aEw_8$.

We cannot have  $M \vDash w_8=a$  because  $M \vDash \text{Tally}_b(t_2)$.

$\implies M \vDash aEw_8$,

$\implies M \vDash \exists w_9\ (w_8=w_9a\ \&\ t_2aw_2=(w_9a)t_4aw_4)$,

$\implies$ by (4.14$^b$),  $M \vDash w_9=t_2 \lor t_2Bw_9$,

$\implies$ either way, $M \vDash t_2\subseteq_p w_9\subseteq_p w_8\subseteq_p v$,

$\implies M \vDash t_3<t_2\subseteq_p v$,  contradicting  $M \vDash \text{Pref}(ava,t_3)$.

(2cii2)  $M \vDash ava=w_7a$.

$\implies M \vDash avat_2aw_2=avat_4aw_4$,

$\implies$ by (3.7), $M \vDash t_2aw_2=t_4aw_4$,

$\implies M \vDash w_3at_3ava(t_4aw_4)=z=w_1at_2aw_2=w_1a(t_4aw_4)$,

$\implies$ by (3.6), $M \vDash w_3at_3ava=w_1a$,

$\Rightarrow$ $M \vDash (w_1a)t_1auat_2aw_2=x=(w_3at_3ava)t_1auat_2aw_2=w_3at_3avat_1aua(t_4aw_4)$,

$\Rightarrow$ $M \vDash \text{Pref}(ava,t_3)$ & $\text{Tally}_b(t_1)$ & $t_3<t_1$ &

$\qquad\qquad\qquad\qquad$ & $\exists w'\ x=w_3at_3avat_1aw'$ & $\text{Max}^+T_b(t_3,w_3)$,

$\Rightarrow$ $M \vDash \text{Intf}(x,w_3,t_3,ava,t_1)$,

$\Rightarrow$ $M \vDash \text{Fr}(x,t_3,ava,t_1)$.

(2cii3)  $M \vDash (ava)B(w_7a)$.

$\Rightarrow$ $M \vDash \exists w_8\ avaw_8=w_7a$,

$\Rightarrow$ $M \vDash avaw_8at_2aw_2=avat_4aw_4$,

$\Rightarrow$ by (3.7), $M \vDash w_8at_2aw_2=t_4aw_4$,

$\Rightarrow$ by (4.14$^b$), $M \vDash w_8=t_4 \lor t_4Bw_8$.

(2cii3a) $M \vDash w_8=t_4$.

$\Rightarrow$ $M \vDash t_4at_2aw_2=t_4aw_4$,

$\Rightarrow$ by (3.7), $M \vDash t_2aw_2=w_4$,

$\Rightarrow$ $M \vDash w_3at_3avat_4a(w_4)=z= w_1a(t_2aw_2)$,

$\Rightarrow$ by (3.6), $M \vDash w_3at_3avat_4a= w_1a$,

$\Rightarrow$ $M \vDash w_1at_1auat_2aw_2=x=w_3at_3avat_4at_1auat_2aw_2$,

$\Rightarrow$ $M \vDash \text{Pref}(ava,t_3)$ & $\text{Tally}_b(t_4)$ & $t_3<t_4$ &

$\qquad\qquad\qquad\qquad$ & $\exists w'\ x=w_3at_3avat_4aw'$ & $\text{Max}^+T_b(t_3,w_3)$,

$\Rightarrow$ $M \vDash \text{Intf}(x,w_3,t_3,ava,t_4)$,

$\Rightarrow$ $M \vDash \text{Fr}(x,t_3,ava,t_4)$, as required.

    (2cii3b)  $M \vDash t_4Bw_8$.

$\Rightarrow$  $M \vDash \exists w_9\ t_4w_9=w_8$,

$\Rightarrow$  $M \vDash (t_4w_9)at_2aw_2=t_4aw_4$,

$\Rightarrow$  by (3.7), $M \vDash w_9at_2aw_2=aw_4$,

$\Rightarrow$  $M \vDash w_1(at_2aw_2)=z=(w_3at_3avat_4)aw_4=(w_3at_3avat_4)w_9(at_2aw_2)$,

$\Rightarrow$  $M \vDash \text{Pref}(ava,t_3)\ \&\ \text{Tally}_b(t_4)\ \&\ t_3<t_4\ \&$

$$\&\ \exists w'\ x=w_3at_3avat_4aw'\ \&\ \text{Max}^+T_b(t_3,w_3),$$

$\Rightarrow$  $M \vDash \text{Intf}(x,w_3,t_3,ava,t_4)$,

$\Rightarrow$  $M \vDash \text{Fr}(x,t_3,ava,t_4)$, as required.

This completes the proof in Case 2.

(3)  $M \vDash \text{Lastf}(z,t_3,ava,t_4)$.

$\Rightarrow$  $M \vDash \text{Pref}(ava,t_3)\ \&\ \text{Tally}_b(t_4)\ \&\ t_3=t_4\ \&$

$$\&\ (z=t_3avat_4\ \vee\ \exists w(z=wat_3avat_4\ \&\ \text{Max}^+T_b(t_3,w))).$$

We derive a contradiction from $M \vDash t_3=t_4\ \&\ z=t_3avat_4$ exactly as in 1(a).

$\Rightarrow$  $M \vDash w_1at_2aw_2=z=wat_3avat_4\ \&\ \text{Max}^+T_b(t_3,w)$,

$\Rightarrow$  $M \vDash (w_1a)Bz\ \&\ (wa)Bz$,

$\Rightarrow$  by (3.8), $M \vDash w_1a=wa\ \vee\ (w_1a)B(wa)\ \vee\ (wa)B(w_1a)$.

  (3a)  $M \vDash w_1a=wa$.

$\Rightarrow$  $M \vDash w_1at_2aw_2=z=w_1at_3avat_4$,

$\Rightarrow$  by (3.7), $M \vDash t_2aw_2=t_3avat_4$,

$\Rightarrow$  $M \vDash w_1at_1aua(t_2aw_2)=x=w_1at_1aua(t_3avat_4)$.

We reason as in (2a), omitting $aw_4$ throughout and replacing $w_3$ with $w_1$, and derive $M \vDash Max^+T_b(t_3,w_1at_1au)$. We then conclude $M \vDash Lastf(x,w_3,t_3,ava,t_4)$, so $M \vDash Fr(x,t_3,ava,t_4)$, as required.

(3b)  $M \vDash (w_1a)B(wa)$.

Again we proceed as in (2b), omitting $aw_4$ throughout and replacing $w_3$ with w. We then obtain $M \vDash Lastf(x,w_3,t_3,ava,t_4)$, so $M \vDash Fr(x,t_3,ava,t_4)$, as required.

(3c)  $M \vDash (wa)B(w_1a)$.

$\Rightarrow$ $M \vDash \exists w_5\ waw_5=w_1a$,

$\Rightarrow$ $M \vDash waw_5t_2aw_2=z=wat_3avat_4$,

$\Rightarrow$ by (3.7), $M \vDash w_5t_2aw_2=t_3avat_4$,

$\Rightarrow$ $M \vDash w_5=a \vee aEw_5$.

But we cannot have $M \vDash w_5=a$ because $M \vDash Tally_b(t_3)$. Hence $M \vDash aEw_5$.

$\Rightarrow$ $M \vDash \exists w_6\ w_6a=w_5$,

$\Rightarrow$ $M \vDash w_6at_2aw_2=t_3avat_4$,

$\Rightarrow$ by (4.16), $M \vDash t_2\subseteq_p v$,

$\Rightarrow$ from $M \vDash w_6at_2aw_2=t_3avat_4$, $M \vDash w_6=t_3 \vee t_3Bw_6$.

(3ci)  $M \vDash w_6=t_3$.

$\Rightarrow$ $M \vDash w_1a=wat_3a$,

$\Rightarrow$ $M \vDash w_1=wat_3$,

$\Rightarrow$ $M \vDash t_3\subseteq_p w_1$,

From hypothesis $M \vDash Intf(x,w_1,t_1,aua,t_2)$, $M \vDash Max^+T_b(t_1,w_1)\ \&\ t_1<t_2$.

$\Rightarrow$ $M \vDash t_3<t_1<t_2\subseteq_p v$, contradicting $M \vDash Pref(ava,t_3)$.

(3cii) $M \vDash t_3 B w_6$.

$\Rightarrow$ $M \vDash \exists w_7\ t_3 w_7 = w_6$,

$\Rightarrow$ $M \vDash t_3 w_7 a = w_6 a = w_5$,

$\Rightarrow$ $M \vDash w_1 a = w a w_5 = w a t_3 w_7 a$,

$\Rightarrow$ $M \vDash w_1 = w a t_3 w_7$,

$\Rightarrow$ $M \vDash t_3 \subseteq_p w_1$,

$\Rightarrow$ from $M \vDash \text{Max}^+ T_b(t_1, w_1)$ & $t_1 < t_2$, $M \vDash t_3 < t_1 < t_2 \subseteq_p v$, again contradicting $M \vDash \text{Pref}(ava, t_3)$.

This completes the proof in Case 3 and the proof of (5.47).

(5.48) For any string concept $I \subseteq I_0$ there is a string concept $J \subseteq I$ such that

$QT^+ \vdash \forall x \in J \ \forall t, w_1, t_1, t_2, u(Env(t,x) \ \& \ Intf(x,w_1,t_1,aua,t_2) \rightarrow$

$\rightarrow \forall z, w_2 (x = w_1 a(t_1 aua)t_2 a w_2 \ \& \ z = w_1 a t_2 a w_2 \rightarrow Env(t,z))]$.

Let $J \equiv I_{5.37} \ \& \ I_{5.47}$.

Assume $M \vDash Env(t,x) \ \& \ Intf(x,w_1,t_1,aua,t_2)$ and let $M \vDash x = w_1 a(t_1 aua)t_2 a w_2$

where $M \vDash J(x)$. Let $M \vDash z = w_1 a t_2 a w_2$.

We first show that (a) $M \vDash MaxT_b(t,z)$.

Assume $M \vDash s \subseteq_p z = w_1 a t_2 a w_2 \ \& \ Tally_b(s)$.

$\Rightarrow$ by (4.14$^b$), $M \vDash s \subseteq_p w_1 \ \vee \ s \subseteq_p t_2 a w_2$,

$\Rightarrow$ since $M \vDash w_1 B x \ \& \ (t_2 a w_2)Ex$, $M \vDash w_1 \subseteq_p x \ \& \ t_2 a w_2 \subseteq_p x$,

$\Rightarrow$ either way, $M \vDash s \subseteq_p x$,

$\Rightarrow$ from $M \vDash Env(t,x)$, $M \vDash MaxT_b(t,x)$,

$\Rightarrow M \vDash s \subseteq_p t$,

$\Rightarrow M \vDash MaxT_b(t,z)$, as required.

Note that from the hypothesis $M \vDash Env(t,x)$ we have

$M \vDash \exists v, t, t' Firstf(x, t', ava, t'')$.

$\Rightarrow$ by (5.36), $M \vDash Firstf(z, t', ava, t_2) \vee Firstf(z, t', ava, t'')$,

Hence we have (b) $M \vDash \exists v, t_3, t_4 Firstf(z, t_3, ava, t_4)$.

Again, from the hypothesis $M \vDash Env(t,x)$ we have $M \vDash \exists v Lastf(x,t,ava,t)$,

$\Rightarrow$ by (5.37), $M \vDash Lastf(z,t,ava,t)$, as required.

For (d) of $M \vDash Env(t,z)$, assume that $M \vDash Fr(z,t_3,ava,t_4) \ \& \ Fr(z,t_5,ava,t_6)$.

$\Rightarrow$ by (5.47),

$M \vDash (Fr(x,t_3,ava,t_4) \lor Fr(x,t_3,ava,t_1)) \& (Fr(x,t_5,ava,t_6) \lor Fr(x,t_5,ava,t_1))$.

But then there are four cases:

(d1) $M \vDash Fr(x,t_3,ava,t_4) \& (Fr(x,t_5,ava,t_6))$,

(d2) $M \vDash Fr(x,t_3,ava,t_4) \& (Fr(x,t_5,ava,t_1))$,

(d3) $M \vDash Fr(x,t_3,ava,t_1) \& (Fr(x,t_5,ava,t_6))$,

(d4) $M \vDash Fr(x,t_3,ava,t_1) \& (Fr(x,t_5,ava,t_1))$.

In each of the cases, it follows from the hypothesis $M \vDash Env(t,x)$ that

$M \vDash t_3=t_5$, as required.

For (e) of $M \vDash Env(t,z)$, assume that $M \vDash Fr(z,t_3,au_1a,t_4) \& Fr(z,t_3,au_2a,t_5)$.

Then, again by (5.47), we have that

$M \vDash (Fr(x,t_3,au_1a,t_4) \lor Fr(x,t_3,au_1a,t_1)) \& (Fr(x,t_3,au_2a,t_5) \lor Fr(x,t_3,au_2 a,t_1))$

and the four cases:

(e1) $M \vDash Fr(x,t_3,au_1a,t_4) \& (Fr(x,t_3,au_2a,t_5))$,

(e2) $M \vDash Fr(x,t_3,au_1a,t_4) \& (Fr(x,t_3,au_2a,t_1))$,

(e3) $M \vDash Fr(x,t_3,au_1a,t_1) \& (Fr(x,t_3,au_2a,t_5))$,

(d4) $M \vDash Fr(x,t_3,au_1a,t_1) \& (Fr(x,t_3,au_2a,t_1))$.

In each case, the hypothesis $M \vDash Env(t,x)$ implies that $M \vDash u_1=u_2$, as required.

This completes the proof that $M \vDash Env(t,z)$ and the proof of (5.48).

(5.49) For any string concept $I \subseteq I_0$ there is a string concept $J \subseteq I$ such that

$QT^+ \vdash \forall x \in J \ \forall t,w_1,t_1,t_2,u,z,w_2 (x=w_1a(t_1aua)t_2aw_2 \ \& \ Intf(x,w_1,t_1,aua,t_2) \ \&$

$\& \ z=w_1at_2aw_2 \ \& \ Env(t,z) \rightarrow \forall t',v,t''(Fr(z,t',ava,t'') \rightarrow t'<t_1 \ v \ t_1<t')]$.

Let $J \equiv I_{3.8} \ \& \ I_{4.14b} \ \& \ I_{4.23b}$.

Assume $M \vDash Intf(x,w_1,t_1,aua,t_2) \ \& \ Env(t,z)$ where

$M \vDash x=w_1a(t_1aua)t_2aw_2 \ \& \ z=w_1at_2aw_2$ and $M \vDash J(x)$.

$\Rightarrow$ from $M \vDash Intf(x,w_1,t_1,aua,t_2)$, $M \vDash t_1<t_2 \ \& \ Max^+T_b(t_1,w_1)$.

Assume also that $M \vDash Fr(z,t',ava,t'')$.

There are three cases:

(1) $M \vDash Firstf(z,t',ava,t'')$.

$\Rightarrow M \vDash (t'ava)Bz$,

$\Rightarrow M \vDash \exists z_1 \ t'avaz_1=z=w_1at_2aw_2$,

$\Rightarrow$ by (4.14$^b$), $M \vDash w_1=t' \ v \ t'Bw_1$,

$\Rightarrow M \vDash t' \subseteq_p w_1$,

$\Rightarrow$ since $M \vDash MaxT_b(t_1,w_1)$, $M \vDash t'<t_1$, as required.

(2) $M \vDash Lastf(z,t',ava,t'')$.

$\Rightarrow M \vDash t'=t''$,

$\Rightarrow$ from $M \vDash Env(t,z)$, $M \vDash t'=t''=t \ \& \ MaxT_b(t,z)$,

$\Rightarrow M \vDash t_1<t_2 \leq t=t'$, as required.

(3) $M \vDash \exists w' Intf(z,w',t',ava,t'')$.

$\Rightarrow M \vDash \exists w'' \ w'at'avat''aw''=z=w_1at_2aw_2$,

$\Rightarrow$ M ⊨ (w'a)Bz & (w₁a)Bz,

$\Rightarrow$ by (3.8), M ⊨ w'a=w₁a v (w'a)B(w₁a) v (w₁a)B(w'a).

  (3a)  M ⊨ w'a=w₁a.

$\Rightarrow$ M ⊨ w₁at'avat''aw''=w₁at₂aw₂,

$\Rightarrow$ by (3.7), M ⊨ t'avat''aw''=t₂aw₂,

$\Rightarrow$ by (4.23$^b$), M ⊨ t'=t₂,

$\Rightarrow$ since M ⊨ t₁<t₂, M ⊨ t₁<t', as required.

  (3b)  M ⊨ (w'a)B(w₁a).

$\Rightarrow$ M ⊨ ∃w₃ w'aw₃=w₁a,

$\Rightarrow$ M ⊨ w'at'avat''aw''=(w'aw₃)t₂aw₂,

$\Rightarrow$ by (3.7), M ⊨ t'avat''aw''=w₃t₂aw₂,

$\Rightarrow$ from M ⊨ w'aw₃=w₁a, M ⊨ w₃=a v aEw₃.

We cannot have because M ⊨ Tally$_b$(t').

$\Rightarrow$ M ⊨ ∃w₄ (w₃=w₄a & t'avat''aw''=(w₄a)t₂aw₂),

$\Rightarrow$ M ⊨ w'aw₄a=w'aw₃=w₁a,

$\Rightarrow$ M ⊨ w'aw₄=w₁,

$\Rightarrow$ by (4.14$^b$), M ⊨ w₄=t' v t'Bw₄,

$\Rightarrow$ M ⊨ t'⊆$_p$w₄⊆$_p$w₁,

$\Rightarrow$ from M ⊨ Max$^+$T$_b$(t₁,w₁), M ⊨ t'<t₁, as required.

  (3c)  M ⊨ (w₁a)B(w'a).

$\Rightarrow$ M ⊨ ∃w₃ w₁aw₃=w'a,

$\Rightarrow$ M ⊨ (w₁aw₃)t'avat''aw''=w₁at₂aw₂,

$\Rightarrow$ by (3.7), $M \vDash w_3 t'avat''aw''=t_2aw_2$,

$\Rightarrow$ from $M \vDash w_1aw_3=w'a$, $M \vDash w_3=a \lor aEw_3$.

We cannot have because $M \vDash Tally_b(t_2)$.

$\Rightarrow$ $M \vDash \exists w_4 (w_3=w_4a \ \& \ (w_4a)t'avat''aw''=t_2aw_2)$,

$\Rightarrow$ by (4.14$^b$), $M \vDash w_4=t_2 \lor t_2Bw_4$,

$\Rightarrow$ $M \vDash t_2 \subseteq_p w_4$,

$\Rightarrow$ $M \vDash w_1aw_4a=w_1aw_3=w'a$,

$\Rightarrow$ $M \vDash w_1aw_4=w'$,

$\Rightarrow$ $M \vDash w_4 \subseteq_p w'$,

$\Rightarrow$ $M \vDash t_2 \subseteq_p w_4 \subseteq_p w'$,

$\Rightarrow$ from $M \vDash Intf(z,w',t',ava,t'')$, $M \vDash Max^+T_b(t',w')$,

$\Rightarrow$ $M \vDash t_1<t_2<t'$, as required.

This completes the proof of (5.49).

(5.50) For any string concept $I \subseteq I_0$ there is a string concept $J \subseteq I$ such that

$QT^+ \vdash \forall x \in J\ \forall w_1,t_1,t_2,u,z,w_2(x=w_1a(t_1aua)t_2aw_2\ \&\ Intf(x,w_1,t_1,aua,t_2)\ \&$

$\&\ z=w_1at_2aw_2 \rightarrow \forall t',v,t''(Intf(x,w',t',ava,t'')\ \&\ t' \neq t_1 \rightarrow$

$\rightarrow Fr(z,t',ava,t'')\ v\ Fr(z,t',ava,t_2)]$.

Let $J \equiv I_{3.8}\ \&\ I_{4.14b}\ \&\ I_{4.23b}$.

Assume $M \vDash Intf(x,w_1,t_1,aua,t_2)\ \&\ Intf(x,w',t',ava,t'')\ \&\ t' \neq t_1$ where

$M \vDash x=w_1a(t_1aua)t_2aw_2\ \&\ z=w_1at_2aw_2$ and $M \vDash J(x)$.

$\Rightarrow M \vDash Tally_b(t_1)\ \&\ Tally_b(t')$,

$\Rightarrow$ from $M \vDash t' \neq t_1$ by (4.6), $M \vDash t_1 < t'\ v\ t' < t_1$,

$\Rightarrow M \vDash \exists w''\ w_1a(t_1aua)t_2aw_2 = x = w'at'avat''aw''$ where

$M \vDash Pref(aua,t_1)\ \&\ Tally_b(t_2)\ \&\ t_1 < t_2\ \&\ Max^+T_b(t_1,w_1)$

and $M \vDash Pref(ava,t')\ \&\ Tally_b(t'')\ \&\ t' < t''\ \&\ Max^+T_b(t',w')$,

$\Rightarrow M \vDash (w_1a)Bx\ \&\ (w'a)Bx$,

$\Rightarrow$ by (3.8), $M \vDash w_1a=w'a\ v\ (w_1a)B(w'a)\ v\ (w'a)B(w_1a)$.

(1) $M \vDash w_1a=w'a$.

$\Rightarrow M \vDash w_1at_1auat_2aw_2=w_1at'avat''aw''$,

$\Rightarrow$ by (3.7), $M \vDash t_1auat_2aw_2=t'avat''aw''$,

$\Rightarrow$ since $M \vDash Tally_b(t_1)$, by $(4.23^b)$, $M \vDash t_1=t'$, contradicting the hypothesis.

(2) $M \vDash (w_1a)B(w'a)$.

$\Rightarrow M \vDash \exists w_3\ w_1aw_3=w'a$,

$\Rightarrow M \vDash w_1at_1auat_2aw_2=(w_1aw_3)t'avat''aw''$,

⟹ by (3.7), $M \vDash t_1auat_2aw_2=w_3t'avat''aw''$,

⟹ from $M \vDash w_1aw_3=w'a$,

⟹ $M \vDash w_3=a \lor aEw_3$.

We cannot have $M \vDash w_3=a$ because $M \vDash \text{Tally}_b(t_1)$.

⟹ $M \vDash aEw_3$,

⟹ $M \vDash \exists w_4 \; w_3=w_4a$,

⟹ $M \vDash w_1aw_4a=w_1aw_3=w'a$,

⟹ $M \vDash w_1aw_4=w'$,

⟹ $M \vDash t_1auat_2aw_2=(w_4a)t'avat''aw''$,

⟹ by (4.14$^b$), $M \vDash w_4=t_1 \lor t_1Bw_4$.

 (2i) $M \vDash w_4=t_1$.

⟹ $M \vDash t_1auat_2aw_2=t_1at'avat''aw''$,

⟹ by (3.7), $M \vDash uat_2aw_2=t'avat''aw''$,

⟹ by (4.14$^b$), $M \vDash u=t' \lor t'Bu$,

⟹ $M \vDash t'\subseteq_p u$,

⟹ $M \vDash t_1\subseteq_p w_4\subseteq_p w'$,

⟹ from $M \vDash \text{Max}^+T_b(t',w')$, $M \vDash t_1<t'$,

⟹ $M \vDash t_1\subseteq_p t'\subseteq_p u$, contradicting $M \vDash \text{Pref}(aua,t_1)$.

 (2ii) $M \vDash t_1Bw_4$.

⟹ $M \vDash \exists w_5 \; t_1w_5=w_4$,

⟹ $M \vDash t_1auat_2aw_2=(t_1w_5)at'avat''aw''$,

⟹ by (3.7), $M \vDash auat_2aw_2=w_5at'avat''aw''$,

$\Rightarrow$ by (3.8), $M \vDash aua=w_5a \lor (aua)B(w_5a) \lor (w_5a)B(aua)$.

(2iia) $M \vDash aua=w_5a$.

$\Rightarrow M \vDash (aua)t_2aw_2=(aua)t'avat''aw''$,

$\Rightarrow$ by (3.7), $M \vDash t_2aw_2=z=t'avat''aw''$,

$\Rightarrow M \vDash Pref(ava,t') \& Tally_b(t'') \& t'<t'' \& (t'avat''a)Bz$,

$\Rightarrow M \vDash Firstf(z,t',ava,t'')$,

$\Rightarrow M \vDash Fr(z,t',ava,t'')$.

(2iib) $M \vDash (aua)B(w_5a)$.

$\Rightarrow M \vDash \exists w_6\ auaw_6=w_5a$,

$\Rightarrow M \vDash auat_2aw_2=(auaw_6)t'avat''aw''$,

$\Rightarrow$ by (3.7), $M \vDash t_2aw_2=w_6t'avat''aw''$,

$\Rightarrow$ from $M \vDash auaw_6=w_5a$, $M \vDash w_6=a \lor aEw_6$.

We cannot have $M \vDash w_6=a$ because $M \vDash Tally_b(t_2)$.

$\Rightarrow M \vDash aEw_6$,

$\Rightarrow M \vDash \exists w_7\ w_6=w_7a$,

$\Rightarrow M \vDash t_2aw_2=z=(w_7a)t'avat''aw''$.

We claim that $M \vDash Max^+T_b(t',w_7)$.

Assume $M \vDash s\subseteq_p w_7 \& Tally_b(s)$.

$\Rightarrow$ from $M \vDash auaw_7a=auaw_6=w_5a$, $M \vDash auaw_7=w_5$,

$\Rightarrow M \vDash s\subseteq_p w_7 \subseteq_p w_5 \subseteq_p w_4 \subseteq_p w'$.

Therefore, $M \vDash Max^+T_b(t',w_7)$, as claimed.

$\Rightarrow M \vDash Pref(ava,t') \& Tally_b(t'') \& t'<t'' \& z=(w_7a)t'avat''aw'' \& Max^+T_b(t',w_7)$,

$\Rightarrow$ M ⊨ Intf(z,$w_7$, t',ava,t''),

$\Rightarrow$ M ⊨ Fr(z,t',ava,t'').

(2iic)  M ⊨ ($w_5$a)B(aua).

$\Rightarrow$ M ⊨ $\exists w_6$ aua=$w_5$a$w_6$,

$\Rightarrow$ M ⊨ ($w_5$a$w_6$)$t_2$a$w_2$=$w_5$at'avat''aw'',

$\Rightarrow$ by (3.7),  M ⊨ $w_6 t_2$a$w_2$=t'avat''aw'',

$\Rightarrow$ from  M ⊨ aua=$w_5$a$w_6$,  M ⊨ $w_6$=a  v  aE$w_6$.

We cannot have  M ⊨ $w_6$=a   because  M ⊨ Tally$_b$(t').

$\Rightarrow$ M ⊨ aE$w_6$,

$\Rightarrow$ M ⊨ $\exists w_7$ ($w_6$=$w_7$a & $w_7$a$t_2$a$w_2$=t'avat''aw''),

$\Rightarrow$ by (4.14$^b$),  M ⊨ $w_7$=t' v t'B$w_7$,

$\Rightarrow$ M ⊨ t'$\subseteq_p w_7 \subseteq_p w_6 \subseteq_p$ aua.

But  M ⊨ $t_1 \subseteq_p w_4 \subseteq_p$ w',  so from  M ⊨ Max$^+$T$_b$(t',w')  we have   M ⊨ $t_1$<t'.

$\Rightarrow$ M ⊨ $t_1 \subseteq_p$ t'$\subseteq_p$ aua,  contradicting  M ⊨ Pref(aua,$t_1$).

(3)  M ⊨ (w'a)B($w_1$a).

$\Rightarrow$ M ⊨ $\exists w_3$ $w_1$a=w'a$w_3$,

$\Rightarrow$ M ⊨ (w'a$w_3$)$t_1$aua$t_2$a$w_2$=x=w'at'avat''aw'',

$\Rightarrow$ by (3.7),  M ⊨ $w_3 t_1$aua$t_2$a$w_2$=t'avat''aw'',

$\Rightarrow$ from M ⊨ w'a$w_3$=$w_1$a,

$\Rightarrow$ M ⊨ $w_3$=a v aE$w_3$.

We cannot have  M ⊨ $w_3$=a  because  M ⊨ Tally$_b$(t').

$\Rightarrow$ M ⊨ aE$w_3$,

$\Rightarrow$ M ⊨ $\exists w_4\ w_3=w_4a$,

$\Rightarrow$ M ⊨ $w'aw_4a=w'aw_3=w_1a$,

$\Rightarrow$ M ⊨ $w'aw_4=w_1$,

$\Rightarrow$ M ⊨ $(w_4a)t_1auat_2aw_2=t'avat''aw''$,

$\Rightarrow$ by (4.14$^b$), M ⊨ $w_4=t'\ \vee\ t'Bw_4$.

(3i) M ⊨ $w_4=t'$.

$\Rightarrow$ M ⊨ $t'at_1auat_2aw_2=t'avat''aw''$,

$\Rightarrow$ by (3.7), M ⊨ $t_1auat_2aw_2=vat''aw''$,

$\Rightarrow$ by (4.14$^b$), M ⊨ $v=t_1\ \vee\ t_1Bv$, $\Rightarrow$ M ⊨ $t_1\subseteq_p v$,

$\Rightarrow$ from M ⊨ $t'\subseteq_p w_4 \subseteq_p w_1$ & Max$^+T_b(t_1,w_1)$, M ⊨ $t'<t_1$,

$\Rightarrow$ M ⊨ $t'\subseteq_p t_1\subseteq_p v$, contradicting M ⊨ Pref(ava,t').

(3ii) M ⊨ $t'Bw_4$.

$\Rightarrow$ M ⊨ $\exists w_5\ t'w_5=w_4$,

$\Rightarrow$ M ⊨ $(t'w_5)at_1auat_2aw_2=t'avat''aw''$,

$\Rightarrow$ by (3.7), M ⊨ $w_5at_1auat_2aw_2=avat''aw''$,

$\Rightarrow$ by (3.8), M ⊨ $ava=w_5a\ \vee\ (ava)B(w_5a)\ \vee\ (w_5a)B(ava)$.

(3iia) M ⊨ $ava=w_5a$.

$\Rightarrow$ M ⊨ $t'w_5a=t'ava=w_4a$,

$\Rightarrow$ M ⊨ $w_1a=w'a(w_4a)=w'a(t'ava)$,

$\Rightarrow$ M ⊨ $(w_1a)t_2aw_2=z=(w'at'ava)t_2aw_2$,

$\Rightarrow$ M ⊨ $w_1=w'at'av$,

$\Rightarrow$ M ⊨ $t'\subseteq_p w_1$,

$\Rightarrow$ from $M \vDash Max^+T_b(t_1,w_1)$, $M \vDash t'<t_1<t_2$,

$\Rightarrow M \vDash Pref(ava,t')$ & $Tally_b(t_2)$ & $t'<t_2$ & $z=w'at'avat_2aw_2$ & $Max^+T_b(t',w')$,

$\Rightarrow M \vDash Intf(z,w', t',ava,t_2)$,

$\Rightarrow M \vDash Fr(z,t',ava,t_2)$.

(3iib)  $M \vDash (ava)B(w_5a)$.

$\Rightarrow M \vDash \exists w_6\ avaw_6=w_5a$,

$\Rightarrow M \vDash (avaw_6)at_1auat_2aw_2=avat''aw''$,

$\Rightarrow$ by (3.7), $M \vDash w_6at_1auat_2aw_2=t''aw''$,

$\Rightarrow$ from $M \vDash avaw_6=w_5a$, $M \vDash w_6=a \lor aEw_6$.

We cannot have $M \vDash w_6=a$ because $M \vDash Tally_b(t'')$.

$\Rightarrow M \vDash aEw_6$,

$\Rightarrow M \vDash \exists w_7\ w_6=w_7a$,

$\Rightarrow$ from $M \vDash avaw_6=w_5a$, $M \vDash t'w_5a=t'avaw_6=t'avaw_7a$,

$\Rightarrow M \vDash (w_7a)at_1auat_2aw_2=t''aw''$,

$\Rightarrow$ by (4.14[b]), $M \vDash t''=w_7 \lor t''Bw_7$.

(3iib1)  $M \vDash t''=w_7$.

$\Rightarrow M \vDash t'w_5a=t'avat''a$,

$\Rightarrow M \vDash w_4a=t'avat''a$,

$\Rightarrow M \vDash w_3=t'avat''a$,

$\Rightarrow M \vDash w'aw_3=w'at'avat''a$,

$\Rightarrow$ from $M \vDash w_1a=w'aw_3$, $M \vDash w_1a=w'at'avat''a$,

$\Rightarrow M \vDash (w_1a)t_2aw_2=z=(w'at'avat''a)t_2aw_2$,

$\Rightarrow M \vDash \text{Pref}(ava,t') \ \& \ \text{Tally}_b(t'') \ \& \ t'<t'' \ \& \ z=w'at'avat''at_2aw_2 \ \& \ \text{Max}^+T_b(t',w')$,

$\Rightarrow M \vDash \text{Intf}(z,w', t',ava,t'')$,

$\Rightarrow M \vDash \text{Fr}(z,t',ava,t'')$.

(3iib2) $M \vDash t''Bw_7$.

$\Rightarrow M \vDash \exists w_8 \ t''w_8=w_7$,

$\Rightarrow$ since $M \vDash t'w_5a=t'avaw_7a$, $M \vDash t'w_5a=t'ava(t''w_8)a$,

$\Rightarrow$ from $M \vDash t'w_5=w_4$, $M \vDash w_4a=t'avat''w_8a$,

$\Rightarrow M \vDash w_3=t'avat''w_8a$,

$\Rightarrow M \vDash w'aw_3=w'at'avat''w_8a$,

$\Rightarrow M \vDash w_1a=w'at'avat''w_8a$,

$\Rightarrow$ from $M \vDash t''w_8=w_7 \ \& \ w_7a=w_6 \ \& \ avaw_6=w_5a$, $M \vDash w_5a=avaw_7a=avat''w_8a$,

$\Rightarrow$ from $M \vDash w_5at_1auat_2aw_2=avat''aw''$, $M \vDash (avat''w_8a)t_1auat_2aw_2=avat''aw''$,

$\Rightarrow$ by (3.7), $M \vDash w_8at_1auat_2aw_2=aw''$,

$\Rightarrow M \vDash w_8=a \ \lor \ aBw_8$,

$\Rightarrow M \vDash w_8=a \ \lor \ \exists w_9 \ aw_9=w_8$,

$\Rightarrow M \vDash w_1a=w'at'avat''aa \ \lor \ w_1a=w'at'avat''(aw_9)a$,

$\Rightarrow M \vDash (w_1a)t_2aw_2=z=(w'at'avat''aa)t_2aw_2 \ \lor$

$\lor \ (w_1a)t_2aw_2=z=w'at'avat''aw_9at_2aw_2$,

Either way we have

$M \vDash \text{Pref}(ava,t') \ \& \ \text{Tally}_b(t'') \ \& \ t'<t'' \ \& \ \exists w(z=w'at'avat''aw \ \& \ \text{Max}^+T_b(t',w'))$,

$\Rightarrow M \vDash \text{Intf}(z,w', t',ava,t'')$,

$\Rightarrow M \vDash \text{Fr}(z,t',ava,t'')$.

(3iic)  $M \vDash (w_5a)B(ava)$.

$\Longrightarrow M \vDash \exists w_6\ ava=w_5aw_6$,

$\Longrightarrow M \vDash w_5at_1auat_2aw_2=(w_5aw_6)t''aw''$,

$\Longrightarrow$ by (3.7), $M \vDash t_1auat_2aw_2=w_6t''aw''$,

$\Longrightarrow$ from $M \vDash ava=w_5aw_6$, $M \vDash w_6=a \vee aEw_6$.

We cannot have $M \vDash w_6=a$ because $M \vDash Tally_b(t'')$.

$\Longrightarrow M \vDash aEw_6$,

$\Longrightarrow M \vDash \exists w_7\ w_6=w_7a$,

$\Longrightarrow M \vDash t_1auat_2aw_2=w_7at''aw''$,

$\Longrightarrow$ by (4.14$^b$), $M \vDash w_7=t_1 \vee t_1Bw_7$,

$\Longrightarrow M \vDash t_1 \subseteq_p w_7 \subseteq_p w_6 \subseteq_p ava$.

But $M \vDash t' \subseteq_p w_4 \subseteq_p w_1$, so from $M \vDash Max^+T_b(t_1,w_1)$ we have $M \vDash t'<t_1$.

$\Longrightarrow M \vDash t' \subseteq_p t_1 \subseteq_p ava$, contradicting $M \vDash Pref(ava,t')$.

This completes the proof of (5.50).

(5.51) For any string concept $I\subseteq I_0$ there is a string concept $J\subseteq I$ such that

$QT^+ \vdash \forall x\in J\ \forall w_1,t_1,t_3,u,z,w,t',v't''((w_1at_1)Bx\ \&\ Intf(x,w,t_3,aua,t_1)\ \&$

$\&\ z=wat_3auat_3\ \&\ Fr(z,t',ava,t'') \rightarrow Fr(x,t',ava,t'')\ v\ Fr(x,t',ava,t_1))$.

Let $J \equiv I_{4.16}\ \&\ I_{4.24b}$.

Assume $M \vDash Intf(x,w,t_3,aua,t_1)\ \&\ Fr(z,t',ava,t'')$

where $M \vDash (w_1at_1)Bx\ \&\ z=wat_3auat_3$ and $M \vDash J(x)$.

(1) $M \vDash Firstf(z,t',ava,t'')$.

$\Rightarrow M \vDash Pref(ava,t')\ \&\ Tally_b(t'')\ \&$

$\&\ ((t'=t''\ \&\ z=t'avat'')\ v\ (t'<t''\ \&\ (t'avat''a)Bz))$.

 (1a) $M \vDash t'=t''\ \&\ z=t'avat''$.

$\Rightarrow M \vDash wat_3auat_3=z=t'avat''$,

$\Rightarrow$ since $M \vDash Tally_b(t')\ \&\ Tally_b(t'')$, by (4.16), $M \vDash t_3\subseteq_p v$,

$\Rightarrow$ from $M \vDash Pref(ava,t')$, $M \vDash t_3<t'$,

$\Rightarrow$ from $M \vDash wat_3auat_3=t'avat''$, by (4.14$^b$), $M \vDash w=t'\ v\ t'Bw$,

$\Rightarrow M \vDash t'\subseteq_p w$,

$\Rightarrow$ from $M \vDash Intf(x,w,t_3,aua,t_1)$, $M \vDash Max^+T_b(t_3,w)$,

$\Rightarrow M \vDash t'<t_3$, contradicting $M \vDash t'\in I_0$.

 (1b) $M \vDash t'<t''\ \&\ (t'avat''a)Bz$.

$\Rightarrow$ from $M \vDash Intf(x,w,t_3,aua,t_1)$, $M \vDash t_3<t_1\ \&\ Tally_b(t_1)\ \&\ \exists w_2\ x=wat_3auat_1aw_2$,

$\Rightarrow M \vDash \exists t_4\ t_3t_4=t_1$,

$\Rightarrow M \vDash x=wat_3aua(t_3t_4)aw_2$,

$\Rightarrow$ from $M \vDash z=wat_3auat_3$, $M \vDash zBx$,

$\Rightarrow$ from $M \vDash (t'avat''a)Bz$, $M \vDash Pref(ava,t')$ & $Tally_b(t'')$ & $t'<t''$ & $(t'avat''a)Bx$,

$\Rightarrow$ $M \vDash Firstf(x,t',ava,t'')$,

$\Rightarrow$ $M \vDash Fr(x,t',ava,t'')$, as required.

(2) $M \vDash \exists w'\ Intf(z,w',t',ava,t'')$.

$\Rightarrow$ $M \vDash Pref(ava,t')$ & $Tally_b(t'')$ & $t'<t''$ & $\exists w''(z=w'at'avat''aw''$ &

& $Max^+T_b(t',w'))$,

$\Rightarrow$ as in (1b), $M \vDash \exists t_4\ t_3t_4=t_1$,

$\Rightarrow$ $M \vDash w'at'avat''aw''=z=wat_3auat_3$,

$\Rightarrow$ $M \vDash wat_3auat_3(t_4aw_2)=x=w'at'avat''aw''(t_4aw_2)$,

$\Rightarrow$ $M \vDash Pref(ava,t')$ & $Tally_b(t'')$ & $t'<t''$ & $x=w'at'avat''aw_3$ & $Max^+T_b(t',w')$,

where $w_3=w''t_4aw_2$,

$\Rightarrow$ $M \vDash Intf(x,w',t',ava,t'')$,

$\Rightarrow$ $M \vDash Fr(x,t',ava,t'')$.

(3) $M \vDash Lastf(z,t',ava,t'')$.

$\Rightarrow$ $M \vDash Pref(ava,t')$ & $Tally_b(t'')$ & $t'=t''$ &

& $(z=t'avat''\ \vee\ \exists w'(z=w'at'avat''$ & $Max^+T_b(t',w'))$.

Now, $M \vDash z=t'avat''$ is ruled out as in (1a).

$\Rightarrow$ $M \vDash w'at'avat''=z=wat_3auat_3$,

$\Rightarrow$ by (4.24[b]), $M \vDash t'=t''=t_3$,

$\Rightarrow$ as in (1b), $M \vDash x=wat_3auat_1aw_2=wat_3aua(t_3t_4)aw_2=(w'at'avat'')t_4aw_2=$

$=w'at'ava(t_3t_4)aw_2=w'at'avat_1aw_2$,

$\implies$ M ⊨ Intf(x,w',t',ava,$t_1$),

$\implies$ M ⊨ Fr(x,t',ava,$t_1$).

This completes the proof of (5.51).

(5.52) For any string concept $I \subseteq I_0$ there is a string concept $J \subseteq I$ such that

$QT^+ \vdash \forall x \in J \; \forall t,t',t'',t_1,t_3,u,z,v,w(Env(t,x) \; \& \; Firstf(x,t',ava,t'') \;\&$

$\& \; Intf(x,w,t_3,aua,t_1) \; \& \; z=wat_3auat_3 \; \rightarrow \; Firstf(z,t',ava,t''))$.

Let $J \equiv I_{3.6} \; \& \; I_{4.14b} \; \& \; I_{5.34}$.

Assume $M \vDash Env(t,x) \; \& \; Firstf(x,t',ava,t'') \; \& \; Intf(x,w,t_3,aua,t_1)$ where $M \vDash J(x)$.

Let $M \vDash z=wat_3auat_3$.

$\Rightarrow$ by (5.34), $M \vDash t'' \leq t_3 < t_1$,

$\Rightarrow$ from $M \vDash Firstf(x,t',ava,t'')$,

$M \vDash Pref(ava,t') \; \& \; Tally_b(t'') \; \& \; ((t'=t'' \; \& \; x=t'avat'') \; v \; (t'<t'' \; \& \; (t'avat''a)Bx))$.

(1) $M \vDash t'=t'' \; \& \; x=t'avat''$.

This is ruled out by the same argument as in Case 1(a) of (5.51).

(2) $M \vDash t'<t'' \; \& \; (t'avat''a)Bx$.

$\Rightarrow$ from $M \vDash Intf(x,w,t_3,aua,t_1)$, $M \vDash Tally_b(t_3) \; \& \; \exists w_2 \; x=wat_3auat_1aw_2$,

$\Rightarrow M \vDash (wat_3)Bx \; \& \; (t'avat''a)Bx$,

$\Rightarrow$ by (3.8), $M \vDash (wat_3)B(t'avat''a) \; v \; wat_3=t'avat''a \; v \; (t'avat''a)B(wat_3)$.

Since $M \vDash Tally_b(t_3)$, $M \vDash \neg(wat_3=t'avat''a)$.

(2a) $M \vDash (wat_3)B(t'avat''a)$.

$\Rightarrow M \vDash \exists x_1 \; wat_3x_1=t'avat''a$,

$\Rightarrow M \vDash x_1=a \; v \; aEx_1$,

$\Rightarrow M \vDash wat_3a=t'avat''a \; v \; \exists x_2 \; wat_3(x_2a)=t'avat''a$.

(2ai) $M \vDash wat_3a=t'avat''a$.

$\Rightarrow$ from $M \vDash z = wat_3auat_3$, $M \vDash (t'avat''a)Bz$,

$\Rightarrow$ $M \vDash \text{Pref}(ava, t')$ & $\text{Tally}_b(t'')$ & $t' < t''$ & $(t'avat''a)Bz$,

$\Rightarrow$ $M \vDash \text{Firstf}(z, t', ava, t'')$, as required.

(2aii) $M \vDash \exists x_2 \, wat_3(x_2a) = t'avat''a$.

$\Rightarrow$ $M \vDash wat_3x_2 = t'avat''$,

$\Rightarrow$ by (3.10), $M \vDash x_2Et'' \lor x_2 = t'' \lor t''Ex_2$.

(2aii1) $M \vDash x_2 = t''$.

$\Rightarrow$ by (3.6), $M \vDash wat_3 = t'ava$, a contradiction because $M \vDash \text{Tally}_b(t_3)$.

(2aii2) $M \vDash t''Ex_2$.

$\Rightarrow$ $M \vDash \exists x_3 \, x_3t'' = x_2$,

$\Rightarrow$ $M \vDash wat_3(x_3t'') = t'avat''$,

$\Rightarrow$ by (3.6), $M \vDash wat_3x_3 = t'ava$,

$\Rightarrow$ $M \vDash x_3 = a \lor aEx_3$,

$\Rightarrow$ $M \vDash wat_3a = t'ava \lor \exists x_4 \, wat_3(x_4a) = t'ava$,

$\Rightarrow$ $M \vDash wat_3at'' = t'avat'' \lor wat_3x_4at'' = t'avat''$,

$\Rightarrow$ by (4.16), $M \vDash t_3 \subseteq_p v \lor t_3 \subseteq_p t_3x_4 \subseteq_p v$,

$\Rightarrow$ $M \vDash t_3 \subseteq_p v$,

$\Rightarrow$ from $M \vDash \text{Pref}(ava, t')$, $M \vDash t_3 < t'$,

$\Rightarrow$ $M \vDash t_3 < t' < t'' \leq t_3$, contradicting $M \vDash t_3 \in I \subseteq I_0$.

(2aii3) $M \vDash x_2Et''$.

$\Rightarrow$ $M \vDash \exists t_4 \, (\text{Tally}_b(t_4) \, \& \, t_4 < t'' \, \& \, t_4x_2 = t'')$,

$\Rightarrow$ $M \vDash wat_3x_2 = t'avat_4x_2$,

$\Rightarrow$ by (3.6), $M \vDash wat_3 = t'avat_4$,

$\Rightarrow$ since $M \vDash Tally_b(t_3)$ & $Tally_b(t_4)$, by (4.24$^b$), $M \vDash t_3 = t_4$,

$\Rightarrow$ $M \vDash t_3 = t_4 < t'' \leq t_3$, contradicting $M \vDash t_3 \in I \subseteq I_0$.

(2b) $M \vDash (t'avat''a)B(wat_3)$.

$\Rightarrow$ from $M \vDash (wat_3a)Bz$, $M \vDash (t'avat''a)Bz$,

$\Rightarrow$ $M \vDash Firstf(z,t',ava,t'')$, as required.

This completes the proof of (5.52).

(5.53) For any string concept $I \subseteq I_0$ there is a string concept $J \subseteq I$ such that

$QT^+ \vdash \forall x \in J \, \forall t, w_1, t_1 \, (Env(t,x) \,\&\, (w_1 a t_1)Bx \rightarrow$

$\rightarrow \forall z, w, u, t_3 \, (Intf(x, w, t_3, aua, t_1) \,\&\, z = wat_3auat_3 \rightarrow Env(t_3, z)))$.

Let $J \equiv I_{4.17b} \,\&\, I_{5.51} \,\&\, I_{5.52}$.

Assume $M \vDash Intf(x, w, t_3, aua, t_1) \,\&\, (w_1 a t_1)Bx$ where $M \vDash Env(t,x) \,\&\, J(x)$.

Assume that $M \vDash z = wat_3auat_3$.

Suppose that $M \vDash s \subseteq_p z \,\&\, Tally_b(s)$.

$\Rightarrow M \vDash s \subseteq_p wat_3auat_3$,

$\Rightarrow$ by (4.17$^b$), $M \vDash s \subseteq_p w \,\vee\, s \subseteq_p t_3 auat_3$.

If $M \vDash s \subseteq_p w$, then $M \vDash s < t_3$ because from $M \vDash Intf(x, w, t_3, aua, t_1)$ we have

$M \vDash Max^+T_b(t_3, w)$. If $M \vDash s \subseteq_p t_3 auat_3$, then, again by (4.17$^b$),

$M \vDash s \subseteq_p t_3 \,\vee\, s \subseteq_p uat_3$. Once again by (4.17$^b$), $M \vDash s \subseteq_p u \,\vee\, s \subseteq_p t_3$. But then from $M \vDash s \subseteq_p u$ we have $M \vDash s < t_3$ since $M \vDash Pref(aua, t_3)$.

Therefore, $M \vDash s \subseteq_p t_3$ in any case, and so we have

  (a) $M \vDash MaxT_b(t_3, z)$.

From $M \vDash Env(t,x)$ we have $M \vDash \exists v_0, t', t'' Firstf(x, t', av_0 a, t'')$. By (5.52), we have

  (b) $M \vDash Firstf(z, t', av_0 a, t'')$.

That $M \vDash Lastf(z, t_3, aua, t_3)$ follows immediately from the hypothesis

$M \vDash Intf(x, w, t_3, aua, t_1)$ and the choice of z.

For (d) of $M \vDash Env(t_3, z)$, assume that $M \vDash Fr(z, t', ava, t_4) \,\&\, Fr(z, t'', ava, t_5)$.

$\Rightarrow$ by (5.51),

$$M \vDash (Fr(x,t',ava,t_4) \vee Fr(x,t'',ava,t_1)) \& (Fr(x,t'',ava,t_5) \vee Fr(x,t'',ava,t_1)).$$

There are four cases:

(d1)  $M \vDash Fr(x,t',ava,t_4) \& Fr(x,t'',ava,t_5)$,

(d2)  $M \vDash Fr(x,t',ava,t_4) \& Fr(x,t'',ava,t_1)$,

(d3)  $M \vDash Fr(x,t',ava,t_1) \& Fr(x,t'',ava,t_5)$,

(d4)  $M \vDash Fr(x,t',ava,t_1) \& Fr(x,t'',ava,t_1)$.

In any of these cases it follows from (d) of $M \vDash Env(t,x)$ that $M \vDash t'=t''$. This establishes (d) of $M \vDash Env(t_3,z)$.

Finally, for (e) of $M \vDash Env(t_3,z)$, suppose that

$$M \vDash Fr(z,t',av_1a,t_4) \& Fr(z,t',av_2a,t_5).$$

We again apply (5.51) and obtain

$$M \vDash (Fr(x,t',av_1a,t_4) \vee Fr(x,t',av_1a,t_1)) \& (Fr(x,t',av_2a,t_5) \vee Fr(x,t',av_2a,t_1)).$$

Again we have four cases:

(e1)  $M \vDash Fr(x,t',av_1a,t_4) \& Fr(x,t',av_2a,t_5)$,

(e2)  $M \vDash Fr(x,t',av_1a,t_4) \& Fr(x,t',av_2a,t_1)$,

(e3)  $M \vDash Fr(x,t',av_1a,t_1) \& Fr(x,t',av_2a,t_5)$,

(e4)  $M \vDash Fr(x,t',av_1a,t_1) \& Fr(x,t',av_2a,t_1)$.

In each case we have from (e) of $M \vDash Env(t,x)$ that $M \vDash v_1=v_2$. This establishes (e) of $M \vDash Env(t_3,z)$ and concludes the proof of (5.51).

(5.54)  For any string concept I⊆I₀ there is a string  concept J⊆I such that

$QT^+ \vdash \forall x \in J \; \forall t,t_1,t_2,w_1,v \; (Env(t,x) \& Intf(x,w_1,t_1,ava,t_2) \rightarrow$

$\rightarrow \exists x' \in J \; \exists t_3,w \; (Env(t_1,x') \& x'Bx \& Lastf(x',t_1,ava,t_1) \&$

$\& \; \forall x''(Env(t_1,x') \& x'Bx \& Lastf(x',t_1,ava,t_1) \rightarrow x''=x')))$.

Let $J \equiv I_{4.19} \& I_{5.53}$.

Assume $M \vDash Env(t,x) \& Intf(x,w_1,t_1,ava,t_2)$  where $M \vDash J(x)$.

$\Rightarrow M \vDash Tally_b(t_2) \& t_1<t_2 \& \exists w_2 \; x=w_1at_1avat_2aw_2 \& t_1<t_2 \& Max^+T_b(t_1,w_1)$.

Let $x'=w_1at_1avat_1$.

$\Rightarrow$ from $M \vDash t_1<t_2$, $M \vDash \exists t_3 \; (Tally_b(t_3) \& t_1t_3=t_2)$,

$\Rightarrow M \vDash x=(w_1at_1avat_1)t_3aw_2=x't_3aw_2$,

$\Rightarrow M \vDash x'Bx$,

$\Rightarrow$ by (5.53),  $M \vDash Env(t_1,x') \& Lastf(x',t_1,ava,t_1)$.

Assume now that  $M \vDash Env(t_1,x'') \& x''Bx \& Lastf(x'',t_1,ava,t_1)$.

$\Rightarrow M \vDash x''=t_1avat_1 \lor \exists w_3 \; (x''=w_3at_1avat_1 \& Max^+T_b(t_1,w_3))$.

Suppose, for a reductio, that $M \vDash x''=t_1avat_1$.

$\Rightarrow$ from $M \vDash x''Bx$, $M \vDash \exists x_1 \; x''x_1=x$,

$\Rightarrow M \vDash (t_1avat_1)x_1=x''x_1=x=w_1at_1avat_2aw_2$,

$\Rightarrow$ since $M \vDash Tally_b(t_1)$, by $(4.14^b)$, $M \vDash w_1=t_1 \lor t_1Bw_1$,

$\Rightarrow M \vDash t_1 \subseteq_p w_1$, contradicting $M \vDash Max^+T_b(t_1,w_1)$.

Therefore, $M \vDash \exists w_3 \; (x''=w_3at_1avat_1 \& Max^+T_b(t_1,w_3))$.

$\Rightarrow M \vDash (w_3at_1avat_1)x_1=x''x_1=x=w_1at_1avat_2aw_2$ & $Max^+T_b(t_1,w_3)$ &

& $Max^+T_b(t_1,w_1)$,

$\Rightarrow$ by (4.19), $M \vDash w_1=w_3$,

$\Rightarrow M \vDash x'=w_1at_1avat_1=w_3at_1avat_1=x''$, as required.

This completes the proof of (5.54).

(5.55) For any string concept I⊆I₀ there is a string concept J⊆I such that

$QT^+ \vdash \forall x \in J \; \forall t, t_1, t_2, w_1, v \; (Env(t,x) \; \& \; Intf(x, w_1, t_1, ava, t_2) \rightarrow$

$\rightarrow \exists x' \in J \; \exists t_3, w \; (Env(t_1, x') \; \& \; Lastf(x', t_1, ava, t_1) \; \& \; x = x' t_3 wt \; \&$

$\& \; Tally_b(t_3) \; \& \; aBw \; \& \; aEw \; \&$

$\& \; \forall x''(Env(t_1, x'') \; \& \; Lastf(x'', t_1, ava, t_1) \; \& \; x = x'' t_3 wt \; \& \; Tally_b(t_3) \; \& \; aBw \; \& \; aEw \rightarrow$

$\rightarrow x' = x'')))$.

Let $J \equiv I_{5.22} \; \& \; I_{5.54}$.

Assume $M \vDash Env(t,x) \; \& \; Intf(x, w_1, t_1, ava, t_2)$ where $M \vDash J(x)$.

$\Rightarrow M \vDash \exists w_2 \; x = w_1 a t_1 a v a t_2 a w_2 \; \& \; t_1 < t_2 \; \& \; Max^+ T_b(t_1, w_1)$,

$\Rightarrow M \vDash \exists t_3 \; (t_2 = t_1 t_3 \; \& \; Tally_b(t_3))$,

$\Rightarrow M \vDash x = (w_1 a t_1 a v a t_1) t_3 a w_2$.

Let $x' = w_1 a t_1 a v a t_1$.

$\Rightarrow$ since $M \vDash J(x)$ and we may assume J is downward closed w.r. to $\subseteq_p$,

$M \vDash J(x')$,

$\Rightarrow$ by (5.54), $M \vDash Env(t_1, x')$,

$\Rightarrow M \vDash x = x' t_3 a w_2$.

We now argue that $M \vDash \exists w \; (t_3 a w_2 = t_3 wt \; \& \; aBw \; \& \; aEw)$.

$\Rightarrow$ from $M \vDash Env(t,x)$, $M \vDash \exists v' \; Lastf(x, t, av'a, t) \; \& \; MaxT_b(t,x)$,

$\Rightarrow M \vDash t_1 < t_2 \leq t$,

$\Rightarrow$ from $M \vDash Env(t,x)$, since $M \vDash t \in I_0$, $M \vDash v \neq v'$,

$\Rightarrow$ from $M \vDash Lastf(x, t, av'a, t)$,

$M \vDash \text{Pref}(av'a,t)$ & $(x=tav'at \lor \exists w'(x=w'atav'at$ & $\text{Max}^+T_b(t,w')))$.

If $M \vDash x=tav'at$, then from $M \vDash \text{Pref}(av'a,t)$ we have

$$M \vDash \text{Firstf}(x,t,av'a,t) \text{ \& } \text{Lastf}(x,t,av'a,t),$$

whence, by (5.22), $M \vDash \forall u (u \varepsilon x \leftrightarrow u=v')$. But from $M \vDash \text{Intf}(x,w_1,t_1,ava,t_2)$, $M \vDash v \varepsilon x$, and we obtain $M \vDash v=v'$, a contradiction.

Therefore $M \vDash \neg(x=tav'at)$.

$\Rightarrow M \vDash \exists w'(x=w'atav'at$ & $\text{Max}^+T_b(t,w'))$,

$\Rightarrow M \vDash (w'at)Bx$ & $x'Bx$,

$\Rightarrow$ by (3.8), $M \vDash (w'at)Bx' \lor w'at=x' \lor x'B(w'at)$.

(a) $M \vDash (w'at)Bx'$.

$\Rightarrow M \vDash t \subseteq_p x'$,

$\Rightarrow M \vDash \text{Env}(t_1,x')$, $M \vDash \text{MaxT}_b(t_1,x')$,

$\Rightarrow M \vDash t \leq t_1 < t_2 \leq t$, contradicting $M \vDash t \in I \subseteq I_0$.

(b) $M \vDash w'at=x'$.

$\Rightarrow M \vDash w'at=x'=w_1at_1avat_1$,

$\Rightarrow$ from $M \vDash \text{Tally}_b(t)$ & $\text{Tally}_b(t_1)$, by (4.24$^b$), $M \vDash t=t_1$, again contradicting $M \vDash t \in I \subseteq I_0$ as in (a).

Therefore, we have

(c) $M \vDash x'B(w'at)$.

$\Rightarrow M \vDash \exists x_1\, x'x_1=w'at$,

$\Rightarrow M \vDash (x'x_1)av'at=(w'at)av'at=x=x't_3aw_2$,

$\Rightarrow$ by (3.7), $M \vDash x_1av'at=t_3aw_2$,

$\Rightarrow$ by (4.14$^b$), $M \vDash x_1 = t_3 \vee t_3 B x_1$,

$\Rightarrow M \vDash t_3 a v' a t = t_3 a w_2 \vee \exists x_2\ t_3 x_2 a v' a t = t_3 a w_2$.

If $M \vDash t_3 x_2 a v' a t = t_3 a w_2$, then, by (3.7), $M \vDash x_2 a v' a t = a w_2$,

$\Rightarrow M \vDash x_2 = a \vee a B x_2$,

$\Rightarrow M \vDash t_3 a w_2 = t_3 a v' a t \vee t_3 a w_2 = t_3 a a v' a t \vee \exists x_3\ t_3 a w_2 = t_3 (a x_3) a v' a t$,

$\Rightarrow M \vDash \exists w\ (x = x' t_3 w t\ \&\ a B w\ \&\ a E w)$, as required.

The uniqueness follows immediately from (5.54).

This completes the proof of (5.55).

(5.56) For any string concept $I \subseteq I_0$ there is a string concept $J \subseteq I$ such that

$QT^+ \vdash \forall y \in J \, \forall t,t_1,t_2,t',t'',u,v,y' \, (Env(t,y) \, \& \, Fr(y,t_1,ava,t_2) \, \& \, y'By \, \& \, Env(t_1,y') \, \&$

$\& \, Lastf(y',t_1,ava,t_1) \, \& \, Fr(y,t',aua,t'') \, \& \, t' \leq t_1 \rightarrow u \, \varepsilon \, y').$

Let $J \equiv I_{5.24} \, \& \, I_{5.27} \, \& \, I_{5.33} \, \& \, I_{5.53}$.

Assume $M \vDash Env(t,y) \, \& \, Fr(y,t_1,ava,t_2) \, \& \, y'By$ along with

$M \vDash Env(t_1,y') \, \& \, Lastf(y',t_1,ava,t_1)$ and $M \vDash J(y)$.

Let $M \vDash Fr(y,t',aua,t'') \, \& \, t' \leq t_1$.

If $M \vDash t'=t_1$, then from $M \vDash Fr(y,t_1,ava,t_2) \, \& \, Fr(y,t',aua,t'') \, \& \, Env(t,y)$ we have

$M \vDash v=u$, whence from $M \vDash Lastf(y',t_1,ava,t_1)$, $M \vDash u \, \varepsilon \, y'$, as required.

So we may assume that $M \vDash t' < t_1$.

$\Rightarrow$ from $M \vDash Env(t,y) \, \& \, t_1 \in I \subseteq I_0$, $M \vDash t_1 \neq t'$,

$\Rightarrow M \vDash u \neq v$.

From hypothesis $M \vDash Fr(y,t',aua,t'')$ we distinguish three cases:

(1) $M \vDash Firstf(y,t',aua,t'')$.

$\Rightarrow$ from $M \vDash Env(t,y) \, \& \, y'By$, by (5.4), $M \vDash \exists t_3 \, Firstf(y',t',aua,t_3)$,

$\Rightarrow M \vDash u \, \varepsilon \, y'$, as required.

(2) $M \vDash \exists w_1 Intf(y,w_1,t',aua,t'')$.

$\Rightarrow M \vDash Pref(aua,t') \, \& \, t' < t'' \, \& \, Tally_b(t'') \, \& \, \exists w_2 \, y=w_1at'auat''aw_2 \, \&$

$\& \, Max^+T_b(t',w_1).$

From hypothesis $M \vDash Fr(y,t_1,ava,t_2)$ we distinguish three subcases:

(2a) $M \vDash Fr(y,t_1,ava,t_2)$.

$\Rightarrow$ from $M \vDash Env(t,y)$ & $Intf(y,w_1,t',aua,t'')$, by (5.19),

$\qquad M \vDash \neg Firstf(y,t',aua,t'')$,

$\Rightarrow$ by (5.20), $M \vDash t_1 < t'$,

$\Rightarrow M \vDash t' < t_1 < t'$, contradicting $M \vDash t' \in I \subseteq I_0$.

(2b) $M \vDash \exists w_3\, Intf(y,w_3,t_1,ava,t_2)$.

$\Rightarrow M \vDash Pref(ava,t_1)$ & $t_1 < t_2$ & $Tally_b(t_2)$ & $\exists w_4\, y=w_3 a t_1 a v a t_2 a w_4$ &

$\qquad\qquad\qquad\qquad\qquad$ & $Max^+T_b(t_1,w_3)$.

$\Rightarrow$ from $M \vDash Lastf(y',t_1,ava,t_1)$,

$\qquad M \vDash y' = t_1 a v a t_1 \lor \exists w'(y' = w' a t_1 a v a t_1$ & $Max^+T_b(t_1,w'))$.

Now we have

$\quad M \vDash Env(t,y)$ & $Intf(y,w_3,t_1,ava,t_2)$ & $y'By$ & $Env(t_1,y')$ & $Lastf(y',t_1,ava,t_1)$.

$\Rightarrow$ by (5.33), $M \vDash y' \neq t_1 a v a t_1$,

$\Rightarrow M \vDash \exists w'(y' = w' a t_1 a v a t_1$ & $Max^+T_b(t_1,w'))$.

(2bi) $M \vDash t'' = t_1$.

$\Rightarrow$ from hypothesis $M \vDash y'By$, $M \vDash \exists z\, y'z = y$,

$\Rightarrow M \vDash (w' a t_1 a v a t_1)z = y = w_1 a t' a u a t'' a w_2$,

$\Rightarrow$ from $M \vDash t' < t_1$ & $Max^+T_b(t'',w_1)$, $M \vDash Max^+T_b(t_1,w_1)$.

Assume $M \vDash Tally_b(t_0)$ & $t_0 \subseteq_p w_1 a t' a u$.

$\Rightarrow$ by (4.17$^b$), $M \vDash t_0 \subseteq_p w_1 \lor t_0 \subseteq_p t' a u$.

If $M \vDash t_0 \subseteq_p t' a u$, then, again by (4.17$^b$), $M \vDash t_0 \subseteq_p t' \lor t_0 \subseteq_p u$. Hence from

$M \vDash Pref(aua,t')$, $M \vDash t_0 \leq t' < t'' = t_1$. Together with $M \vDash Max^+T_b(t_1,w_1)$ we then

have that $M \vDash Max^+T_b(t_1, w_1 a t' a u)$.

Hence from

$M \vDash w'at_1(avat_1z)=w_1at'auat''aw_2$ & $Max^+T_b(t'',w')$ & $Max^+T_b(t_1,w_1at'au)$,

we have, by (4.20), that $M \vDash w'=w_1at'au$.

$\Rightarrow$ from $M \vDash y'=w'at_1avat_1$, $M \vDash (w_1at'au)at_1avat_1$,

$\Rightarrow M \vDash Pref(aua,t')$ & $Tally_b(t_1)$ & $t'<t_1$ & $\exists w''\ y'=w_1at'auat_1aw''$ &

$\phantom{x}$ & $Max^+T_b(t',w_1)$,

$\Rightarrow M \vDash Intf(y',w_1,t',aua,t_1)$,

$\Rightarrow M \vDash u\ \varepsilon\ y'$, as required.

$\quad$ (2bii) $M \vDash t''\neq t_1$.

$\quad\quad$ (2bii1) $M \vDash t_2=t'$.

$\Rightarrow M \vDash t_1<t_2$ & $t'<t''$ & $t'<t_1$,

$\Rightarrow M \vDash t'<t_1<t_2=t'$, contradicting $M \vDash t'\in I\subseteq I_0$.

$\quad\quad$ (2bii2) $M \vDash t_2\neq t'$.

$\Rightarrow M \vDash Env(t,y)$ & $Intf(y',w_1,t',aua,t_1)$ & $Intf(y,w_3,t_1,ava,t_2)$ & $t''\neq t_1$ & $t_2\neq t'$,

$\Rightarrow$ by (5.27), $M \vDash t'auat''\subseteq_p w_3\ v\ t_1avat_2\subseteq_p w_1$.

If $M \vDash t_1avat_2\subseteq_p w_1$, then $M \vDash t'<t_1\subseteq_p w_1$, contradicting $M \vDash Max^+T_b(t',w_1)$.

Therefore, $M \vDash t'auat''\subseteq_p w_3$, and from the proof of (5.27), part (2iib), we in fact have

$\quad\quad M \vDash w_1at'auat''=w_3\ v\ w_1at'auat''a=w_3\ v\ \exists w_5\ w_1at'auat''aw_5=w_3$.

But we have

$M \vDash w_3at_1avat_2aw_4=y=y'z=(w'at_1avat_1)z$ & $Max^+T_b(t_1,w_3)$ & $Max^+T_b(t',w_1)$,

$\Rightarrow$ by (4.20), $M \vDash w'=w_3$,

$\Rightarrow$ M ⊨ $(w_1at'aut'')at_1avat_1=y'$ v $(w_1at'aut''a)at_1avat_1=y'$ v

$$v\ \exists w_5\ (w_1at'aut''aw_5)at_1avat_1=y',$$

$\Rightarrow$ M ⊨ $Pref(aua,t')$ & $Tally_b(t'')$ & $t'<t''$ & $\exists w''\ y'=w_1at'auat''aw''$ &

$$\&\ Max^+T_b(t',w_1),$$

$\Rightarrow$ M ⊨ $Intf(y',w_1,t',aua,t'')$,

$\Rightarrow$ M ⊨ $u\ \varepsilon\ y'$, as required.

(2c) M ⊨ $Lastf(y,t_1,ava,t_2)$.

$\Rightarrow$ M ⊨ $Env(t,y)$ & $Lastf(y,t_1,ava,t_2)$ & $Fr(y,t_1,ava,t_2)$ & $u\neq v$ & $Env(t_1,y')$ &

$$\&\ Lastf(y',t_1,ava,t_1)\ \&\ \exists w_1 Intf(y,w_1,t',aua,t'').$$

Let $z'=w_1at'auat'$. Then M ⊨ $z'By$.

$\Rightarrow$ by (5.53), M ⊨ $Env(t,z')$.

Along with hypothesis M ⊨ $y'By$, we obtain a contradiction by (5.24).

(3) M ⊨ $Lastf(y,w_1,t',aua,t'')$.

$\Rightarrow$ from M ⊨ $Env(t,y)$, M ⊨ $t'=t''$ & $MaxT_b(t,y)$,

$\Rightarrow$ M ⊨ $t_1\leq t=t'$,

$\Rightarrow$ M ⊨ $t_1\leq t'<t_1$, contradicting M ⊨ $t_1\in I\subseteq I_0$.

This completes the proof of (5.56).

(5.57)  For any string concept $I \subseteq I_0$ there is a string concept $J \subseteq I$ such that

$QT^+ \vdash \forall y \in J \, \forall x,t,t',t'',u,v(Env(t,y) \, \& \, Fr(y,t',aua,t'') \, \& \, Env(t',x) \, \&$

$\& \, Lastf(y,t,ava,t) \, \& \, xBy \rightarrow \neg Lastf(x,t',ava,t'))$.

Let $J \equiv I_{5.7} \, \& \, I_{5.16} \, \& \, I_{5.19} \, \& \, I_{5.55}$.

Assume  $M \vDash Env(t,y) \, \& \, Fr(y,t',aua,t'') \, \& \, Env(t',x)$  where

$M \vDash Lastf(y,t,ava,t) \, \& \, xBy$  and  $M \vDash J(y)$.

Suppose, for a reductio, that  $M \vDash Lastf(x,t',ava,t')$.  Then from  $M \vDash xBy$,

$M \vDash \exists z \, xz=y$.  We distinguish three cases:

(1)  $M \vDash Firstf(y,t',aua,t'')$.

$\Rightarrow$  $M \vDash Pref(aua,t') \, \& \, Tally_b(t'') \, \& \, ((t'=t'' \, \& \, y=t'auat'') \, v$

$v \, (t'<t'' \, \& \, (t'auat''a)By))$,

$\Rightarrow$  from  $M \vDash Env(t',x) \, \& \, Firstf(y,t',aua,t'') \, \& \, xBy$, by (5.4),

$M \vDash \exists t_1 Firstf(x,t',aua,t_1)$,

$\Rightarrow$  from  $M \vDash Lastf(x,t',ava,t')$,   $M \vDash Fr(x,t',aua,t_1) \, \& \, Fr(x,t',ava,t')$,

$\Rightarrow$  from  $M \vDash Env(t',x)$,  $M \vDash u=v$,

$\Rightarrow$  from  $M \vDash Firstf(x,t',ava,t_1) \, \& \, Lastf(x,t',ava,t')$, by (5.7),

$M \vDash t'=t_1 \, \& \, x=t'avat'$.

(1a)  $M \vDash t'=t'' \, \& \, y=t'auat''$.

$\Rightarrow$  from  $M \vDash u=v \, \& \, xz=y$,   $M \vDash (t'avat')z=y=t'avat'$,

$\Rightarrow$  $M \vDash yBy$,  contradicting  $M \vDash y \in I \subseteq I_0$.

(1b)  $M \vDash t'<t'' \, \& \, (t'auat''a)By$.

$\Rightarrow$ from $M \vDash \text{Env}(t,y)$, $M \vDash \text{MaxT}_b(t,y)$,

$\Rightarrow$ $M \vDash t' < t'' \leq t$,

$\Rightarrow$ from $M \vDash \text{Lastf}(y,t,\text{ava},t)$ & $u=v$, $M \vDash \text{Fr}(y,t',\text{aua},t'')$ & $\text{Fr}(y,t,\text{aua},t)$,

$\Rightarrow$ from $M \vDash \text{Env}(t,y)$, $M \vDash t'=t$,

$\Rightarrow$ $M \vDash t=t'<t$, contradicting $M \vDash t \in I \subseteq I_0$.

(2) $M \vDash \exists w' \text{Intf}(y,w',t',\text{aua},t'')$.

$\Rightarrow$ from $M \vDash \text{Env}(t,y)$, by (5.55),

$M \vDash \exists x_1,t_2,w_1 (\text{Env}(t',x_1)$ & $y=x_1 t_2 w_1 t$ & $\text{Lastf}(x_1,t',\text{aua},t')$ & $\text{Tally}_b(t_1)$ &

& $aBw_1$ & $aEw_1$),

$\Rightarrow$ from $M \vDash \text{Intf}(y,w',t',\text{aua},t'')$, $M \vDash t'<t''$,

$\Rightarrow$ from $M \vDash \text{Env}(t,y)$, $M \vDash \text{MaxT}_b(t,y)$,

$\Rightarrow$ $M \vDash t'<t'' \leq t$,

$\Rightarrow$ from $M \vDash \text{Lastf}(x_1,t',\text{aua},t')$,

$M \vDash \text{Pref}(\text{aua},t')$ & $(x_1 = t' \text{auat}' \lor \exists w_2 (x_1 = w_2 \text{at'auat'}$ & $\text{Max}^+ T_b(t',w_2)))$.

Suppose, for a reductio, that $M \vDash x_1 = t' \text{auat}'$.

$\Rightarrow$ from $M \vDash \text{Pref}(\text{aua},t')$, $M \vDash \text{Firstf}(x,t',\text{aua},t')$,

$\Rightarrow$ from $M \vDash \text{Env}(t,y)$ & $x_1 By$ & $\text{Env}(t',x_1)$, by (5.16),

$M \vDash \exists t_3 \text{Firstf}(y,t',\text{aua},t_3)$.

But this contradicts $M \vDash \text{Intf}(y,w',t',\text{aua},t'')$, by (5.19).

$\Rightarrow$ $M \vDash \exists w_2 (x_1 = w_2 \text{at'auat'}$ & $\text{Max}^+ T_b(t',w_2))$,

$\Rightarrow$ from $M \vDash \text{Lastf}(x,t',\text{ava},t')$,

$M \vDash \text{Pref}(\text{ava},t')$ & $(x = t' \text{avat}' \lor \exists w''(x = w'' \text{at'avat'}$ & $\text{Max}^+ T_b(t',w'')))$.

(2a) $M \vDash x=t'avat'$.

$\Rightarrow$ from $M \vDash \text{Pref}(ava,t')$, $M \vDash \text{Firstf}(x,t',ava,t')$,

$\Rightarrow$ from $M \vDash \text{Env}(t,y)$ & $xBy$ & $\text{Env}(t',x)$, by (5.16), $M \vDash \exists t_4 \text{Firstf}(y,t',ava,t_4)$,

$\Rightarrow$ from $M \vDash \text{Lastf}(y,t,ava,t)$ & $\text{Env}(t,y)$, $M \vDash t=t'$,

$\Rightarrow$ $M \vDash t=t'<t$, contradicting $M \vDash t \in I \subseteq I_0$.

(2b) $M \vDash \exists w''(x=w''at'avat'$ & $\text{Max}^+T_b(t',w''))$.

$\Rightarrow$ from $M \vDash x_1By$ & $xz=y$, $M \vDash \exists z_1(w_2at'auat')z_1=x_1z_1=y=xz=(w''at'avat')z$,

$\Rightarrow$ from $M \vDash \text{Max}^+T_b(t',w_2)$ & $\text{Max}^+T_b(t',w'')$, by (4.20), $M \vDash w_2=w''$,

$\Rightarrow$ $M \vDash (w''at'a)uat'z_1=(w''at'a)vat'z$,

$\Rightarrow$ by (3.7), $M \vDash uat'z_1=vat'z$,

$\Rightarrow$ from $M \vDash \text{Pref}(aua,t')$ & $\text{Pref}(ava,t')$, $M \vDash \text{Max}^+T_b(t',u)$ & $\text{Max}^+T_b(t',v)$,

$\Rightarrow$ by (4.20), $M \vDash u=v$,

$\Rightarrow$ from $M \vDash \text{Lastf}(y,t,ava,t)$, $M \vDash \text{Fr}(y,t',aua,t'')$ & $\text{Fr}(y,t,aua,t)$,

$\Rightarrow$ from $M \vDash \text{Env}(t,y)$, $M \vDash t'=t$,

$\Rightarrow$ $M \vDash t=t'<t$, contradicting $M \vDash t \in I \subseteq I_0$.

(3) $M \vDash \text{Lastf}(y,t',aua,t'')$.

$\Rightarrow$ from $M \vDash \text{Lastf}(y,t,ava,t)$, by (5.15), $M \vDash u=v$,

$\Rightarrow$ from $M \vDash \text{Env}(t,y)$, $M \vDash t'=t''=t$,

$\Rightarrow$ from $M \vDash \text{Lastf}(y,t',aua,t'')$,

$M \vDash \text{Pref}(aua,t')$ & $(y=t'auat'' \lor \exists w'(y=w'at'auat''$ & $\text{Max}^+T_b(t',w')))$.

Now, $M \vDash y=t'auat''$ is ruled out as inconsistent with hypothesis

$M \vDash \text{Lastf}(x,t',ava,t')$ as in (1a).

Therefore, $M \vDash \exists w'(y=w'at'auat'' \ \& \ Max^+T_b(t',w'))$.

$\Rightarrow$ from $M \vDash u=v \ \& \ t'=t''=t$, $M \vDash w'atavat=y$,

$\Rightarrow$ from $M \vDash Lastf(x,t',ava,t')$,

$\quad M \vDash Pref(ava,t') \ \& \ (x=t'avat' \ \lor \ \exists w''(x=w''at'avat' \ \& \ Max^+T_b(t',w'')))$.

(3a) $M \vDash x=t'avat'$.

$\Rightarrow$ from $M \vDash xz=y$, $M \vDash (t'avat')z=y=w'atavat$,

$\Rightarrow$ since $M \vDash Tally_b(t')$, by $(4.14^b)$, $M \vDash w'=t' \ \lor \ t'Bw'$,

$\Rightarrow$ $M \vDash t' \subseteq_p w'$, contradicting $M \vDash Max^+T_b(t',w')$.

(3b) $M \vDash \exists w''(x=w''at'avat' \ \& \ Max^+T_b(t',w''))$.

$\Rightarrow$ from $M \vDash t=t' \ \& \ xz=y$, $M \vDash (w''atavat)z=y=w'atavat$,

$\Rightarrow$ from $M \vDash Max^+T_b(t,w'') \ \& \ Max^+T_b(t,w')$, by (4.20), $M \vDash w''=w'$,

$\Rightarrow$ $M \vDash w'atavatz=w'atavat=x$,

$\Rightarrow$ $M \vDash xBx$, contradicting $M \vDash x \in I \subseteq I_0$.

This completes the proof of (5.57).

(5.58) For any string concept $I \subseteq I_0$ there is a string concept $J \subseteq I$ such that

$QT^+ \vdash \forall x \in J \ \forall t_1, t_2, t_3, u, v (\text{Pref}(aua, t_1) \ \& \ \text{Pref}(ava, t_2) \ \& \ t_1 < t_2 \ \& \ t_2 = t_3 \ \& \ u \neq v \ \&$

$\& \ x = t_1 a u a t_2 a v a t_3 \rightarrow \text{Set}(x) \ \& \ \forall w (w \ \varepsilon \ x \leftrightarrow (w = u \lor w = v)).$

Let $J \equiv I_{5.46}$.

Assume $M \vDash \text{Pref}(aua, t_1) \ \& \ \text{Pref}(ava, t_2) \ \& \ t_1 < t_2 \ \& \ t_2 = t_3$ where $M \vDash u \neq v$ and $M \vDash J(x)$. Let $M \vDash x = t_1 a u a t_2 a v a t_3$.

Letting $x' = t_1 a u a t_1$ we have $M \vDash \text{Firstf}(x', t_1, aua, t_1) \ \& \ \text{Lastf}(x', t_1, aua, t_1)$.

$\Rightarrow$ by (5.22), $M \vDash \text{Env}(t_1, x') \ \& \ u \ \varepsilon \ x' \ \& \ \forall w (w \ \varepsilon \ x' \leftrightarrow w = u)$.

Likewise, for $z = t_2 a v a t_3$, $M \vDash \text{Firstf}(z, t_2, ava, t_3) \ \& \ \text{Lastf}(z, t_2, ava, t_3)$.

$\Rightarrow$ by (5.22), $M \vDash \text{Env}(t_2, z) \ \& \ v \ \varepsilon \ z \ \& \ \forall w (w \ \varepsilon \ z \leftrightarrow w = v)$,

$\Rightarrow$ by (5.46), $M \vDash \text{Env}(t_2, x) \ \& \ \forall w (w \ \varepsilon \ x \leftrightarrow (w \ \varepsilon \ z \lor w \ \varepsilon \ x') \leftrightarrow (w = u \lor w = v))$,

as required.

This completes the proof of (5.58).

Let   Free(x,t,aua,t')

abbreviate

  Fr(x,t,aua,t') & ¬Firstf(x,t,aua,t') &

    & ∀v,t$_1$,t$_2$ (Fr(x,t$_1$,ava,t$_2$) & v<$_x$u & ¬∃w(v<$_x$w & w<$_x$u)  →  t=t$_1$b),

Roughly, tauat' is a frame in x with initial tally marker t, and the frame immediately preceding tauat' has as its initial tally marker t$_1$, the immediate predecessor of the b-tally t.

Let   Bound(x,t,aua,t')

abbreviate

  Fr(x,t,aua,t') & ¬Firstf(x,t,aua,t') &

    & ∀v,t$_1$,t$_2$ (Fr(x,t$_1$,ava,t$_2$) & v<$_x$u & ¬∃w(v<$_x$w & w<$_x$u)  →  t$_1$b<t).

Roughly, tauat' is a frame in x with initial tally marker t strictly longer than the immediate successor  t$_1$b of the initial marker t$_1$ of the frame immediately preceding  tauat' in x.

We also define   $Free^+(x,t_1,ava,t_2)$

to abbreviate

$Fr(x,t_1,ava,t_2)$ & $\forall u,t',t''(Fr(x,t',aua,t'')$ & $u\leq_x v \rightarrow$

$\rightarrow Firstf(x,t',aua,t'')$ v $Free(x,t',aua,t''))$.

In other words, every frame that precedes $t_1avat_2$ in x that is not first frame is Free.

Finally, we let

$Free^-(x,t_1,ava,t_2) \equiv Free(x,t_1,ava,t_2)$ & $\neg Free^+(x,t_1,ava,t_2)$.

(5.59) For any string concept $I \subseteq I_0$ there is a string concept $J \subseteq I$ such that

$QT^+ \vdash \forall z \in J\ \forall y,v,t',t''(Env(t'',z)\ \&\ z=t'ayat''avat''\ \&\ Firstf(z,t',aya,t'')\ \&$

$\&\ Lastf(z,t'',ava,t'') \to \forall w,t_1,t_2(Free(z,t_1,awa,t_2) \lor Bound(z,t_1,awa,t_2) \to$

$\to w=v\ \&\ t_1=t''))$.

Let $J \equiv I_{5.46}$.

Assume $M \vDash Env(t'',z)\ \&\ z=t'ayat''avat''$ where

$M \vDash Firstf(z,t',aya,t'')\ \&\ Lastf(z,t'',ava,t'')$ and $M \vDash J(z)$.

Assume now that $M \vDash Free(z,t_1,awa,t_2)$.

$\Rightarrow M \vDash Fr(z,t_1,awa,t_2)\ \&\ \neg Firstf(z,t_1,awa,t_2)$,

$\Rightarrow$ from $M \vDash Firstf(z,t',aya,t'')\ \&\ Env(t'',z)$, $M \vDash w \neq y\ \&\ Tally_b(t'')$,

$\Rightarrow$ from $M \vDash Lastf(z,t'',ava,t'')$,

$\qquad M \vDash z=t''avat'' \lor \exists w_1(z=w_1at''avat''\ \&\ Max^+T_b(t'',w_1))$.

Suppose, for a reductio, that $M \vDash z=t''avat''$.

$\Rightarrow M \vDash t'aya(t''avat'')=z=t''avat''$,

$\Rightarrow$ by (3.6), $M \vDash t'ayat''=t''$, contradicting $M \vDash Tally_b(t'')$.

$\Rightarrow M \vDash \exists w_1(z=w_1at''avat''\ \&\ Max^+T_b(t'',w_1))$,

$\Rightarrow M \vDash t'aya(t''avat'')=z=w_1at''avat''$,

$\Rightarrow$ by (3.6), $M \vDash t'ay=w_1$,

$\Rightarrow$ from $M \vDash Firstf(z,t',aya,t'')$, $M \vDash Tally_b(t')$,

$\Rightarrow$ from $M \vDash Max^+T_b(t'',w_1)$, $M \vDash t'<t''$,

$\Rightarrow$ from $M \vDash t' \in I \subseteq I_0$, $M \vDash t' \neq t''$,

⇒ from $M \vDash Fr(z,t',aya,t'') \& Fr(z,t'',ava,t'') \& Env(t'',z)$, $M \vDash y \neq v$,

⇒ by (5.58), $M \vDash \forall u(u \, \varepsilon \, z \leftrightarrow (u=y \vee u=v))$,

⇒ from $M \vDash w \neq y$, $M \vDash w=v$,

⇒ from $M \vDash Fr(z,t_1,awa,t_2) \& Fr(z,t'',ava,t'') \& Env(t'',z)$, $M \vDash t_1=t''$, as required.

The same argument applies if $M \vDash Bound(z,t_1,awa,t_2)$.

This completes the proof of (5.59).

# 6. Lexical Precedence

Let $Rt_L(z,x,y)$ abbreviate

$(((zaBx \lor za=x) \& (zbBy \lor zb=y)) \lor ((zbBx \lor zb=x) \& (zaBy \lor za=y)))$

which we read "z is the (left) root of x and y".

(6.1) For any string concept $I \subseteq I_0$ there is a string concept $J \subseteq I$ such that

$$QT^+ \vdash \forall x,y \in J \, (\exists z \, Rt_L(z,x,y) \to \neg xBy \, \& \, \neg yBx).$$

Let $J \equiv I_{3.7} \& I_{3.8}$.

We may assume that J is closed under $\subseteq_p$.

Assume $M \vDash \exists z \, Rt_L(z,x,y)$ where $M \vDash x,y \in J$.

$\Rightarrow M \vDash ((zaBx \lor za=x) \& (zbBy \lor zb=y)) \lor$

$\lor ((zbBx \lor zb=x) \& (zaBy \lor za=y))$.

Suppose, for a reductio, that $M \vDash xBy$.

There are two cases:

(1) $M \vDash (zaBx \lor za=x) \& (zbBy \lor zb=y)$.

  (1a) $M \vDash zaBx \, \& \, zbBy$.

$\Rightarrow$ from $M \vDash xBy$, by (3.8), $M \vDash xBzb \lor zbBx \lor x=zb$,

$\Rightarrow$ from $M \vDash zaBx$, by (3.9), $M \vDash \neg zbBx$,

$\Rightarrow M \vDash xBzb \lor x=zb$,

$\Rightarrow$ from $M \vDash zaBx$, $M \vDash zaBzb$,

⇒ M ⊨ ∃w zaw=zb,

⇒ by (3.7),  M ⊨ aw=b,  contradicting (QT2).

   (1b)  M ⊨ zaBx  & zb=y.

⇒ from  M ⊨ xBy,  M ⊨ zaBzb.

But this yields a contradiction just as in (1a).

   (1c)  M ⊨ za=x  & zbBy.

⇒ from  M ⊨ xBy,  M ⊨ zaBy,

⇒ by (3.9),  M ⊨ ¬zbBy,  contradicting the hypothesis.

   (1d)  M ⊨ za=x  & zb=y.

⇒ from  M ⊨ xBy,  M ⊨ zaBy,

⇒ M ⊨ zaBzb,  whence a contradiction follows as in (1a).

(2)  M ⊨ (zbBx ∨ zb=x) & (zaBy ∨ za=y).

Analogous to (1).

Therefore  M ⊨ ¬xBy.

An analogous argument shows that  M ⊨ ¬yBx.

This completes the proof of (6.1).

(6.2) For any string concept $I \subseteq I_0$ there is a string concept $J \subseteq I$ such that

$$QT^+ \vdash \forall x,y \in J \; \forall z_1,z_2 \; (Rt_L(z_1,x,y) \; \& \; Rt_L(z_2,x,y) \rightarrow z_1=z_2).$$

Let $J \equiv I_{3.7} \; \& \; I_{3.8}$.

We may assume that J is closed under $\subseteq_p$.

Assume that $M \vDash Rt_L(z_1,x,y) \; \& \; Rt_L(z_2,x,y)$ where $M \vDash J(x) \; \& \; J(y)$.

$\Rightarrow M \vDash (((z_1aBx \; \vee \; z_1a=x) \; \& \; (z_1bBy \; \vee \; z_1b=y)) \; \vee$

$\vee \; ((z_1bBx \; \vee \; z_1b=x) \; \& \; (z_1aBy \; \vee \; z_1a=y))) \; \&$

$\& \; (((z_2aBx \; \vee \; z_2a=x) \; \& \; (z_2bBy \; \vee \; z_2b=y)) \; \vee$

$\vee \; ((z_2bBx \; \vee \; z_2b=x) \; \& \; (z_2aBy \; \vee \; z_2a=y)))$.

There are four cases to be considered, with altogether 64 primary subcases:

(1) $M \vDash (z_1aBx \; \vee \; z_1a=x) \; \& \; (z_1bBy \; \vee \; z_1b=y) \; \& \; (z_2aBx \; \vee \; z_2a=x) \; \&$

$\& \; (z_2bBy \; \vee \; z_2b=y)$.

(2) $M \vDash (z_1aBx \; \vee \; z_1a=x) \; \& \; (z_1bBy \; \vee \; z_1b=y) \; \& \; (z_2bBx \; \vee \; z_2b=x) \; \&$

$\& \; (z_2aBy \; \vee \; z_2a=y)$.

(3) $M \vDash (z_1bBx \; \vee \; z_1b=x) \; \& \; (z_1aBy \; \vee \; z_1a=y) \; \& \; (z_2aBx \; \vee \; z_2a=x) \; \&$

$\& \; (z_2bBy \; \vee \; z_2b=y)$.

(4) $M \vDash (z_1bBx \; \vee \; z_1b=x) \; \& \; (z_1aBy \; \vee \; z_1a=y) \; \& \; (z_2bBx \; \vee \; z_2b=x) \; \&$

$\& \; (z_2aBy \; \vee \; z_2a=y)$.

We illustrate the proof by considering case (1) with its 16 primary subcases.

(1a) $M \vDash z_1aBx \; \& \; z_1bBy \; \& \; z_2aBx \; \& \; z_2bBy$.

$\Rightarrow$ M ⊨ $z_1Bx$ & $z_2Bx$,

$\Rightarrow$ by (3.8), M ⊨ $z_1=z_2$ v $z_1Bz_2$ v $z_2Bz_1$.

Suppose, for a reductio, that M ⊨ $z_1Bz_2$.

$\Rightarrow$ M ⊨ $\exists w_1\ z_1w_1=z_2$,

$\Rightarrow$ M ⊨ $w_1=a$ v $w_1=b$ v $aBw_1$ v $bBw_1$.

   (1ai) M ⊨ $w_1=b$ v $bBw_1$.

$\Rightarrow$ M ⊨ $z_2=z_1b$ v $\exists w_2\ z_1bw_2=z_2$,

$\Rightarrow$ from M ⊨ $z_2aBx$, M ⊨ $(z_1b)aBx$ v $(z_1bw_2)aBx$,

$\Rightarrow$ M ⊨ $z_1bBx$,

$\Rightarrow$ from (3.9), M ⊨ $\neg z_1aBx$, contradicting hypothesis.

   (1aii) M ⊨ $w_1=a$ v $aBw_1$.

$\Rightarrow$ M ⊨ $z_2=z_1a$ v $\exists w_2\ z_1aw_2=z_2$,

$\Rightarrow$ from M ⊨ $z_2bBy$, M ⊨ $(z_1a)bBy$ v $(z_1aw_2)bBy$,

$\Rightarrow$ M ⊨ $z_1aBy$,

$\Rightarrow$ from (3.9), M ⊨ $\neg z_1bBx$, contradicting hypothesis.

Exactly analogously if M ⊨ $z_2Bz_1$.

Hence M ⊨ $z_1=z_2$, as required.

  (1b) M ⊨ $z_1aBx$ & $z_1bBy$ & $z_2aBx$ & $z_2b=y$.

$\Rightarrow$ M ⊨ $z_1Bx$ & $z_2Bx$,

$\Rightarrow$ by (3.8), M ⊨ $z_1=z_2$ v $z_1Bz_2$ v $z_2Bz_1$.

Suppose, for a reductio, that M ⊨ $z_1Bz_2$.

$\Rightarrow$ M ⊨ $\exists w_1\ z_1w_1=z_2$,

⇒ by (QT5),  $M \vDash w_1=a \lor w_1=b \lor aBw_1 \lor bBw_1$.

  (1bi)  $M \vDash w_1=b \lor bBw_1$.

Exactly as in (1ai).

  (1bii)  $M \vDash w_1=a \lor aBw_1$.

Exactly as (1aii).

Therefore,  $M \vDash \neg z_1Bz_2$.

Suppose that  $M \vDash z_2Bz_1$.

⇒ $M \vDash \exists w_2\ z_2w_2=z_1\ \&\ z_1By$,

⇒ $M \vDash \exists w_1\ z_1w_1=y$,

⇒ $M \vDash (z_2w_2)w_1=y$,

⇒ from  $M \vDash z_2b=y$,  $M \vDash z_2w_2w_1=z_2b$,

⇒ by (3.7),  $M \vDash w_2w_1=b$,  contradicting (QT2).

Therefore,  $M \vDash \neg z_2Bz_1$.

Hence, by (3.8),  $M \vDash z_1=z_2$,  as required.

  (1c)  $M \vDash z_1aBx\ \&\ z_1bBy\ \&\ z_2a=x\ \&\ z_2bBy$.

Same as (1a) with  $z_2a=x$  replacing  $z_2aBx$ throughout the argument.

  (1d)  $M \vDash z_1aBx\ \&\ z_1bBy\ \&\ z_2a=x\ \&\ z_2b=y$.

Same as (1b) with the modification as in (1c).

  (1e)  $M \vDash z_1aBx\ \&\ z_1b=y\ \&\ z_2aBx\ \&\ z_2bBy$.

Analogous to (1b) with appropriate modifications.

  (1f)  $M \vDash z_1aBx\ \&\ z_1b=y\ \&\ z_2aBx\ \&\ z_2b=y$.

⇒ $M \vDash z_1b=z_2b$,

$\Rightarrow$ by (QT3), $M \vDash z_1 = z_2$, as required.

(1g) $M \vDash z_1aBx$ & $z_1b=y$ & $z_2a=x$ & $z_2bBy$.

Analogous to (1d).

(1h) $M \vDash z_1aBx$ & $z_1b=y$ & $z_2a=x$ & $z_2b=y$.

Same as (1f).

(1i) $M \vDash z_1a=x$ & $z_1bBy$ & $z_2aBx$ & $z_2bBy$.

Analogous to (1b).

(1j) $M \vDash z_1a=x$ & $z_1bBy$ & $z_2aBx$ & $z_2b=y$.

$\Rightarrow$ $M \vDash yBx$, contradicting (6.1).

(1k) $M \vDash z_1a=x$ & $z_1bBy$ & $z_2a=x$ & $z_2bBy$.

$\Rightarrow M \vDash z_1a = z_2a$,

$\Rightarrow$ by (QT3), $M \vDash z_1 = z_2$, as required.

(1l) $M \vDash z_1a=x$ & $z_1bBy$ & $z_2a=x$ & $z_2b=y$.

Same as (1k).

(1m) $M \vDash z_1a=x$ & $z_1b=y$ & $z_2aBx$ & $z_2bBy$.

Analogous to (1d).

(1n) $M \vDash z_1a=x$ & $z_1b=y$ & $z_2aBx$ & $z_2b=y$.

Analogous to (1f).

(1o) $M \vDash z_1a=x$ & $z_1b=y$ & $z_2a=x$ & $z_2bBy$.

Same as (1k).

(1p) $M \vDash z_1a=x$ & $z_1b=y$ & $z_2a=x$ & $z_2b=y$.

Same as (1k).

This completes the argument in case (1). Cases (2)-(4) are treated analogously.

This completes the proof of (6.2).

The proposition just proved (6.2) establishes the uniqueness of the left root of x, y in an appropriately chosen string concept.   We now prove:

LEFT ROOT LEMMA (6.3)  For any string concept I⊆$I_0$ there is a string concept  $I_{RtL}$⊆I such that

  QT ⊢ ∀x ∈ $I_{RtL}$ (∃z zBx  → ∀y (x≠y  →  y=a v y=b v

                 v (aBx & bBy) v (bBx & aBy) v xBy v yBx v ∃z $Rt_L$(z,x,y))).

Let  $I_{RtL}$(x) abbreviate

I(x)  &  (∃z zBx  → ∀y (x≠y→  y=a v y=b v

                v (aBx & bBy) v (bBx & aBy) v xBy v yBx v ∃z $Rt_L$(z,x,y))).

Trivially,  by (QT2),  $QT^+$ ⊢ $I_{RtL}$(a)  and  $QT^+$⊢ $I_{RtL}$(b).

Assume  M ⊨ $I_{RtL}$(x).

We consider x*a.

Assume  M ⊨ ∃z zB(x*a), and fix y such that  M ⊨ (x*a)≠y. We may assume that  M ⊨ y≠a & y≠b & ¬(aB(x*a) & bBy) & ¬(bB(x*a) & aBy).

Case (1):  M ⊨ aEy.

⇒  M ⊨ ∃$y_1$ y=$y_1$*a,

⇒   M ⊨ x≠$y_1$.

We distinguish two subcases.

   (1a):  M ⊨ ¬∃z zBx.

⇒  by (QT5),  M ⊨ x=a v x=b.

If $M \vDash x=a$, then $M \vDash y_1 \neq a$. By (QT5), $M \vDash y_1=b \lor bBy_1 \lor aBy_1$.

But $M \vDash x=a$ & $(y_1=b \lor bBy_1)$ contradicts the hypothesis

$M \vDash \neg(aB(x*a)$ & $bBy)$. Therefore $M \vDash aBy_1$, so $M \vDash xBy_1$ and $M \vDash xBy$, as required.

If $M \vDash x=b$, then $M \vDash y_1 \neq b$. By (QT5), $M \vDash y_1=a \lor aBy_1 \lor bBy_1$. Then the argument proceeds analogously to the case $M \vDash x=a$.

(1b): $M \vDash \exists z\, zBx$.

$\Rightarrow$ by hypothesis $M \vDash I_{RtL}(x)$, from $M \vDash x \neq y_1$,

$M \vDash y_1=a \lor y_1=b \lor (aBx$ & $bBy_1) \lor (bBx$ & $aBy_1) \lor xBy_1 \lor y_1Bx \lor \exists z\, Rt_L(z,x,y_1)$.

(1b1): $M \vDash y_1=a$.

$\Rightarrow M \vDash \exists z\, zBx$ & $aBy$,

$\Rightarrow$ by (QT5) and (QT2), $M \vDash aBx \lor bBx$,

$\Rightarrow$ from hypothesis $M \vDash \neg(bB(x*a)$ & $aBy)$, $M \vDash aBx$ & $aBy$,

$\Rightarrow M \vDash \exists x_1\, ax_1=x$,

$\Rightarrow$ by (QT5), $M \vDash aBx_1 \lor a=x_1 \lor b=x_1 \lor bBx_1$.

(1b1a): $M \vDash aBx_1$.

$\Rightarrow M \vDash \exists x_2\, a(ax_2)=x$ & $aa=y$,

$\Rightarrow M \vDash yB(x*a)$, as required.

(1b1b): $M \vDash a=x_1$.

$\Rightarrow M \vDash aa=x$ & $aa=y$,

$\Rightarrow M \vDash yB(x*a)$, as required.

(1b1c):  $M \vDash b=x_1$.

$\Rightarrow$ $M \vDash ab=x$ & $aa=y$,

$\Rightarrow$ $M \vDash Rt_L(a,x*a,y)$,

$\Rightarrow$ $M \vDash \exists z\ Rt_L(z,x*a,y)$, as required.

(1b1d):  $M \vDash bBx_1$.

$\Rightarrow$ $M \vDash \exists x_2\ a(bx_2)=x$ & $aa=y$,

$\Rightarrow$ $M \vDash Rt_L(a,x*a,y)$,

$\Rightarrow$ $M \vDash \exists z\ Rt_L(z,x*a,y)$, as required.

(1b2):  $M \vDash y_1=b$.

$\Rightarrow$ $M \vDash \exists z\ zBx$ & $y=ba$,

$\Rightarrow$ by (QT5) and (QT2), $M \vDash aBx \lor bBx$,

$\Rightarrow$ from hypothesis $M \vDash \neg(aB(x*a)\ \&\ bBy)$, $M \vDash bBx\ \&\ bBy$,

$\Rightarrow$ $M \vDash \exists x_1\ bx_1=x$,

$\Rightarrow$ by (QT5), $M \vDash aBx_1 \lor a=x_1 \lor b=x_1 \lor bBx_1$.

(1b2a):  $M \vDash aBx_1$.

$\Rightarrow$ $M \vDash \exists x_2\ b(ax_2)=x$ & $ba=y$,

$\Rightarrow$ $M \vDash yB(x*a)$, as required.

(1b2b):  $M \vDash a=x_1$.

$\Rightarrow$ $M \vDash ba=x$ & $ba=y$,

$\Rightarrow$ $M \vDash yB(x*a)$, as required.

(1b2c):  $M \vDash b=x_1$.

$\Rightarrow$ $M \vDash bb=x$ & $ba=y$,

⇒ $M \vDash Rt_L(b,x*a,y)$,

⇒ $M \vDash \exists z\ Rt_L(z,x*a,y)$, as required.

 (1b1d):  $M \vDash bBx_1$.

⇒ $M \vDash \exists x_2\ b(bx_2)=x\ \&\ ba=y$,

⇒ $M \vDash Rt_L(b,x*a,y)$,

⇒ $M \vDash \exists z\ Rt_L(z,x*a,y)$, as required.

 (1b3):  $M \vDash aBx\ \&\ bBy_1$.

⇒ $M \vDash aB(x*a)\ \&\ bBy$, contradicting the hypothesis $M \vDash \neg(aB(x*a)\ \&\ bBy)$.

 (1b4):  $M \vDash bBx\ \&\ aBy_1$.

⇒ $M \vDash bB(x*a)\ \&\ aBy$, contradicting the hypothesis $M \vDash \neg(bB(x*a)\ \&\ aBy)$.

 (1b5):  $M \vDash xBy_1$.

⇒ $M \vDash \exists x_1\ xx_1=y_1$,

⇒ by (QT5), $M \vDash aBx_1\ v\ a=x_1\ v\ b=x_1\ v\ bBx_1$.

 (1b5a):  $M \vDash aBx_1$.

⇒ $M \vDash (x*a)By_1$,

⇒ $M \vDash (x*a)By$, as required.

 (1b5b):  $M \vDash a=x_1$.

⇒ $M \vDash (x*a)By_1$,

⇒ $M \vDash (x*a)By$, as required.

 (1b5c):  $M \vDash b=x_1$.

⇒ $M \vDash xb=y_1$,

⇒ $M \vDash xbBy$,

$\Rightarrow$ $M \vDash Rt_L(x, x*a, y)$,

$\Rightarrow$ $M \vDash \exists z\ Rt_L(z, x*a, y)$, as required.

    (1b5d): $M \vDash bBx_1$.

$\Rightarrow$ $M \vDash xbBy_1$,

$\Rightarrow$ $M \vDash xbBy$,

$\Rightarrow$ $M \vDash Rt_L(x, x*a, y)$,

$\Rightarrow$ $M \vDash \exists z\ Rt_L(z, x*a, y)$, as required.

   (1b6): $M \vDash y_1Bx$.

$\Rightarrow$ $M \vDash \exists y_2\ y_1y_2 = x$,

$\Rightarrow$ by (QT5), $M \vDash aBy_2\ \lor\ a=y_2\ \lor\ b=y_2\ \lor\ bBy_2$.

    (1b6a): $M \vDash aBy_2$.

Have $M \vDash y = y_1a$,

$\Rightarrow$ $M \vDash yBx$,

$\Rightarrow$ $M \vDash yB(x*a)$, as required.

    (1b6b): $M \vDash a = y_2$.

$\Rightarrow$ $M \vDash y_1y_2 = y_1a = y$,

$\Rightarrow$ $M \vDash yB(x*a)$, as required.

    (1b6c): $M \vDash b = y_2$.

Have $M \vDash y_1Bx\ \&\ y_1b = x\ \&\ y_1a = y$.

$\Rightarrow$ $M \vDash y_1bB(x*a)$,

$\Rightarrow$ $M \vDash Rt_L(y_1, x*a, y)$,

$\Rightarrow$ $M \vDash \exists z\ Rt_L(z, x*a, y)$, as required.

  (1b6d): $M \vDash bBy_2$.

$\Rightarrow$ $M \vDash y_1bBx$ & $y_1a=y$,

$\Rightarrow$ $M \vDash y_1bB(x*a)$,

$\Rightarrow$ $M \vDash Rt_L(y_1,x*a,y)$,

$\Rightarrow$ $M \vDash \exists z\, Rt_L(z,x*a,y)$, as required.

 (1b7): $M \vDash \exists z\, Rt_L(z,x,y_1)$.

$\Rightarrow$ $M \vDash zB(x*a)$ & $zBy$,

$\Rightarrow$ from $M \vDash Rt_L(z,x,y_1)$,

$M \vDash ((zaBx \;\vee\; za=x)$ & $(zbBy_1 \;\vee\; zb=y_1))\;\vee$

$$\vee\; ((zbBx \;\vee\; zb=x)\; \& \;(zaBy_1 \vee za=y_1)).$$

 (1b7a): $M \vDash (zaBx \;\vee\; za=x)$ & $(zbBy_1 \;\vee\; zb=y_1)$.

$\Rightarrow$ $M \vDash zaB(x*a)$ & $zbBy$,

$\Rightarrow$ $M \vDash Rt_L(z,x*a,y)$,

$\Rightarrow$ $M \vDash \exists z\, Rt_L(z,x*a,y)$, as required.

 (1b7b): $M \vDash (zbBx \;\vee\; zb=x)$ & $(zaBy_1 \vee za=y_1)$.

$\Rightarrow$ $M \vDash zbB(x*a)$ & $zaBy$,

$\Rightarrow$ $M \vDash Rt_L(z,x*a,y)$,

$\Rightarrow$ $M \vDash \exists z\, Rt_L(z,x*a,y)$, as required.

This completes the proof of $M \vDash I_{RtL}(x*a)$ in Case (1).

Case (2): $M \vDash bEy$.

$\Rightarrow$ $M \vDash \exists y_1\, y=y_1*b$.

Consider x.

(2a):  $M \vDash \neg \exists z\, zBx$.

$\Rightarrow$ by (QT5) and (QT2), $M \vDash x=a \lor x=b$.

We may assume that $M \vDash x=y_1$; otherwise we proceed as in (1a).

If $M \vDash x=a$, then $M \vDash x*a=aa \,\&\, y=ab$. But then $M \vDash Rt_L(a,x*a,y)$, so

$M \vDash \exists z\, Rt_L(z,x*a,y)$, as required.

If $M \vDash x=b$, then $M \vDash x*a=ba \,\&\, y=bb$. Then $M \vDash Rt_L(b,x*a,y)$, and

$M \vDash \exists z\, Rt_L(z,x*a,y)$, as required.

(2b):  $M \vDash \exists z\, zBx$.

If $M \vDash x=y_1$, then $M \vDash Rt_L(y_1,x*a,y)$ and $M \vDash \exists z\, Rt_L(z,x*a,y)$, as required.

If $M \vDash x \neq y_1$, we appeal to $M \vDash I_{Rt_L}(x)$, and we have

  $M \vDash y_1=a \lor y_1=b \lor (aBx \,\&\, bBy_1) \lor (bBx \,\&\, aBy_1) \lor xBy_1 \lor y_1Bx \lor \exists z\, Rt_L(z,x,y_1)$.

(2b1):  $M \vDash y_1=a$.

$\Rightarrow M \vDash \exists z\, zBx \,\&\, y=ab$,

$\Rightarrow$ by (QT5) and (QT2), $M \vDash aBx \lor bBx$,

$\Rightarrow$ from hypothesis $M \vDash \neg(bB(x*a) \,\&\, aBy)$, $M \vDash aBx \,\&\, aBy$,

$\Rightarrow M \vDash \exists x_1\, ax_1=x$,

$\Rightarrow$ by (QT5), $M \vDash aBx_1 \lor a=x_1 \lor b=x_1 \lor bBx_1$.

(2b1a):  $M \vDash aBx_1$.

$\Rightarrow M \vDash \exists x_2\, a(ax_2)=x \,\&\, ab=y$,

$\Rightarrow M \vDash Rt_L(a,x*a,y)$,

$\Rightarrow M \vDash \exists z\, Rt_L(z,x*a,y)$, as required.

(2b1b):  $M \vDash a=x_1$.

$\Rightarrow$ $M \vDash aa=x$ & $ab=y$,

$\Rightarrow$ $M \vDash Rt_L(a,x*a,y)$,

$\Rightarrow$ $M \vDash \exists z\, Rt_L(z,x*a,y)$, as required.

(2b1c):  $M \vDash b=x_1$.

$\Rightarrow$ $M \vDash ab=x$ & $ab=y$,

$\Rightarrow$ $M \vDash yB(x*a)$, as required.

(2b1d):  $M \vDash bBx_1$.

$\Rightarrow$ $M \vDash \exists x_2\, a(bx_2)=x$ & $ab=y$,

$\Rightarrow$ $M \vDash yB(x*a)$, as required.

(2b2):  $M \vDash y_1=b$.

$\Rightarrow$ $M \vDash \exists z\, zBx$ & $y=bb$,

$\Rightarrow$ by (QT5) and (QT2),  $M \vDash aBx \lor bBx$,

$\Rightarrow$ from hypothesis $M \vDash \neg(aB(x*a)$ & $bBy)$,  $M \vDash bBx$ & $bBy$,

and we proceed analogously to (2b1).

(2b3):  $M \vDash aBx$ & $bBy_1$.

Exactly as in (1b3).

(2b4):  $M \vDash bBx$ & $aBy_1$.

Exactly as in (1b4).

(2b5):  $M \vDash xBy_1$.

Exactly as in (1b5).

(2b6): $M \vDash y_1Bx$.

$\Rightarrow M \vDash \exists y_2 \ y_1y_2=x$,

$\Rightarrow$ by (QT5), $M \vDash aBy_2 \lor a=y_2 \lor b=y_2 \lor bBy_2$.

(2b6a): $M \vDash aBy_2$.

Have $M \vDash y=y_1b$.

$\Rightarrow M \vDash Rt_L(y_1,x*a,y)$,

$\Rightarrow M \vDash \exists z \ Rt_L(z,x*a,y)$, as required.

(2b6b): $M \vDash a=y_2$.

$\Rightarrow M \vDash x=y_1a \ \& \ y_1b=y$,

$\Rightarrow M \vDash Rt_L(y_1,x*a,y)$,

$\Rightarrow M \vDash \exists z \ Rt_L(z,x*a,y)$, as required.

(2b6c): $M \vDash b=y_2$.

Have $M \vDash y_1Bx \ \& \ y_1b=x \ \& \ y_1b=y$,

$\Rightarrow M \vDash yB(x*a)$, as required.

(2b6d): $M \vDash bBy_2$.

$\Rightarrow M \vDash y_1bBx \ \& \ y_1b=y$,

$\Rightarrow M \vDash yB(x*a)$, as required.

(2b7): $M \vDash \exists z \ Rt_L(z,x,y_1)$.

Exactly as in (1b7).

This completes the proof of $M \vDash I_{RtL}(x*a)$ in Case (2).

So, by (QT5), we have established that

$$QT \vdash \forall x \ (I_{RtL}(x) \to I_{RtL}(x*a)).$$

Exactly analogously we show that $QT^+ \vdash \forall x\ (I_{RtL}(x) \rightarrow I_{RtL}(x*b))$.

This completes the proof that $I_{RtL}(x)$ is a string concept and the proof of (6.3).

Analogously, we may introduce "z is the (right) root of x and y": let $Rt_R(z,x,y)$ abbreviate

$(((azEx \lor az=x) \& (bzEy \lor bz=y)) \lor ((bzEx \lor bz=x) \& (azEy \lor az=y)))$.

We then have

Right Root Lemma (6.4)  For any string concept $I \subseteq I_0$ there is a string concept $I_{RtR} \subseteq I$ such that

$QT^+ \vdash \forall x \in I_{RtR} (\exists z\, zEx \rightarrow \forall y\, (x \neq y \rightarrow y=a \lor y=b \lor$

$\lor (aEx \& bEy) \lor (bEx \& aEy) \lor xEy \lor yEx \lor \exists z\, Rt_R(z,x,y)))$.

We let $I_{RtR}(x)$ abbreviate

$I(x)\ \&\ (\exists z\, zEx \rightarrow \forall y\, (x \neq y \rightarrow y=a \lor y=b \lor$

$\lor (aEx \& bEy) \lor (bEx \& aEy) \lor xEy \lor yEx \lor \exists z\, Rt_R(z,x,y)))$.

The proof is completely analogous to that of (6.3).

We let  u≪v   abbreviate

  ((u=a v aBu) & (v=b v bBv))  v  uBv  v

$\qquad\qquad$ v  ∃z (Rt$_L$(z,u,v) & ((za=u v zaBu) & (zb=v v zbBv))).

We then say that string  u <u>lexically precedes</u> string v.

(6.5)  For any string concept  I⊆I$_0$  there is a string  concept  J⊆I such that

$\qquad\qquad$ QT$^+$ ⊢ ∀u,v∈J (u≪v  v  u=v  v  v≪u).

Let  J  ≡  I$_{RtL}$.

Assume  M ⊨ u,v∈J & u≠v.

<u>Case 1.</u>  M ⊨ ¬∃z zBu.

⟹ by (QT5),  M ⊨ u=a v u=b.

 (1a)  M ⊨ u=a.

  (1ai)  M ⊨ ¬∃z zBv.

⟹ by (QT5),  M ⊨ v=a v v=b,

⟹ by hypothesis  M ⊨ u≠v,  M ⊨ v=b,

⟹ M ⊨ u≪v.

  (1aii)  M ⊨ ∃z zBv.

⟹  M ⊨ aBv  v  bBv.

   (1aiia)  M ⊨ aBv.

⟹  M ⊨ uBv,

⟹  M ⊨ u≪v.

    (1aiib)  $M \vDash bBv$.

$\Rightarrow M \vDash u \ll v$.

 (1b)  $M \vDash u=b$.

  (1bi)  $M \vDash \neg \exists z\, zBv$.

$\Rightarrow$ by (QT5), $M \vDash v=a \lor v=b$,

$\Rightarrow$ by hypothesis $M \vDash u \neq v$, $M \vDash v=a$,

$\Rightarrow M \vDash v \ll u$.

  (1bii)  $M \vDash \exists z\, zBv$.

$\Rightarrow M \vDash aBv \lor bBv$.

    (1biia)  $M \vDash aBv$.

$\Rightarrow M \vDash v \ll u$.

    (1biib)  $M \vDash bBv$.

$\Rightarrow M \vDash uBv$,

$\Rightarrow M \vDash u \ll v$.

<u>Case 2.</u>  $M \vDash \exists z\, zBu$.

By the LEFT ROOT LEMMA, we have that

$M \vDash v=a \lor v=b \lor (aBu\ \&\ bBv) \lor (bBu\ \&\ aBv) \lor uBv \lor vBu \lor \exists z\, Rt_L(z,u,v)$.

 (2a)  $M \vDash v=a$.

$\Rightarrow$ by (QT5), $M \vDash aBu \lor bBu$.

  (2ai)  $M \vDash aBu$.

$\Rightarrow M \vDash vBu$,

$\Rightarrow M \vDash v \ll u$.

(2aii) $M \vDash bBu$.

$\Rightarrow M \vDash v \ll u$.

(2b) $M \vDash v=b$.

$\Rightarrow$ by (QT5), $M \vDash aBu \lor bBu$.

(2bi) $M \vDash aBu$.

$\Rightarrow M \vDash u \ll v$.

(2bii) $M \vDash bBu$.

$\Rightarrow M \vDash vBu$,

$\Rightarrow M \vDash v \ll u$.

(2c) $M \vDash aBu \ \& \ bBv$.

$\Rightarrow M \vDash u \ll v$.

(2d) $M \vDash bBu \ \& \ aBv$.

$\Rightarrow M \vDash v \ll u$.

(2e) $M \vDash uBv$.

$\Rightarrow M \vDash u \ll v$.

(2f) $M \vDash vBu$.

$\Rightarrow M \vDash v \ll u$.

(2g) $M \vDash \exists z \ Rt_L(z,u,v)$.

$\Rightarrow M \vDash zBu \ \& \ zBv \ \& \ (((zaBu \lor za=u) \ \& \ (zbBv \lor zb=v)) \lor$

$\lor ((zbBu \lor zb=u) \ \& \ (zaBv \lor za=v)))$.

There are 8 possible cases. The desired result, $M \vDash u \ll v \lor v \ll u$, follows in each case immediately from the definition of $\ll$.

This completes the proof of (6.5).

(6.6) For any string concept $I \subseteq I_0$ there is a string concept $J \subseteq I$ such that

$$QT^+ \vdash \forall v \in J \, \forall u,w \, (u \ll v \, \& \, v \ll w \to u \ll w).$$

Let $J \equiv I_{3.7} \, \& \, I_{3.8}$.

Assume that $M \vDash u \ll v \, \& \, v \ll w$ where $M \vDash J(v)$.

We distinguish nine cases based on the defining conditions for $u \ll v$ and $v \ll w$:

(1a) $M \vDash ((u=a \lor aBu) \, \& \, (v=b \lor bBv)) \, \& \, ((v=a \lor aBv) \, \& \, (w=b \lor bBw))$.

We have that $QT^+ \vdash \neg((v=b \lor bBv) \, \& \, (v=a \lor aBv))$, so this case is ruled out.

(1b) $M \vDash ((u=a \lor aBu) \, \& \, (v=b \lor bBv)) \, \& \, vBw$.

$\implies$ from $M \vDash (v=b \lor bBv)) \, \& \, vBw$, $M \vDash bBw$,

$\implies M \vDash (u=a \lor aBu) \, \& \, bBw$,

$\implies M \vDash u \ll w$.

(1c) $M \vDash ((u=a \lor aBu) \, \& \, (v=b \lor bBv)) \, \&$

$\& \, \exists z \, ((za=v \lor zaBv) \, \& \, (zb=w \lor zbBw))$.

$\implies M \vDash (v=b \lor bBv) \, \& \, (za=v \lor zaBv)$.

By (QT2), $M \vDash \neg(v=b \, \& \, za=v) \, \& \, \neg(v=b \, \& \, zaBv)$.

$\implies M \vDash bBv \, \& \, (za=v \lor zaBv)$,

$\implies M \vDash \exists v_1 \, bv_1=v \, \& \, (za=v \lor \exists v_2 \, zav_2=v)$,

$\implies M \vDash bv_1=za \lor bv_1=zav_2$,

$\implies$ by (QT4) and (QT5), either way $M \vDash z=b \lor bBz$,

$\implies M \vDash bB(zb)$,

$\implies$ from $M \vDash zb=w \lor zbBw$, $M \vDash bBw$,

$\Rightarrow$ M ⊨ (u=a ∨ aBu) & bBw,

$\Rightarrow$ M ⊨ u≪w.

(2a)  M ⊨ uBv & ((v=a ∨ aBv) & (w=b ∨ bBw)).

We have that  QT$^+$ ⊢ ¬(uBv & v=a).

$\Rightarrow$ M ⊨ uBv & aBv,

$\Rightarrow$ M ⊨ ∃v$_1$ uv$_1$=v  &  ∃v$_2$ av$_2$=v,

$\Rightarrow$ M ⊨ uv$_1$=av$_2$,

$\Rightarrow$ by (QT4) and (QT5),  M ⊨ u=a ∨ aBu.

So we have  M ⊨ (u=a ∨ aBu) & (w=b ∨ bBw),  whence  M ⊨ u≪w.

(2b)  M ⊨ uBv & vBw.

$\Rightarrow$ M ⊨ uBw,  whence  M ⊨ u≪w.

(2c)  M ⊨ uBv & ∃z ((za=v ∨ zaBv) & (zb=w ∨ zbBw)).

 (2ci)  M ⊨ uBv & za=v.

$\Rightarrow$ M ⊨ ∃v$_1$ uv$_1$=v=za,

$\Rightarrow$ by (QT4) and (QT5),  M ⊨ v$_1$=a ∨ aEv$_1$,

$\Rightarrow$ M ⊨ ua=za  ∨  ∃v$_2$ u(v$_2$a)=za,

$\Rightarrow$ M ⊨ u=z  ∨  uv$_2$=z,

$\Rightarrow$ from M ⊨ zb=w ∨ zbBw,  M ⊨ uBw,

$\Rightarrow$ M ⊨ u≪w.

 (2cii)  M ⊨ uBv & zaBv.

$\Rightarrow$ M ⊨ uBv & zBv,

$\Rightarrow$ by (3.8),  M ⊨ uBz ∨ u=z ∨ zBu.

(2cii1) $M \vDash uBz \lor u=z$.

$\Rightarrow$ from $M \vDash zb=w \lor zbBw$, $M \vDash uBw$,

$\Rightarrow$ $M \vDash u \ll w$.

(2cii2) $M \vDash zBu$.

$\Rightarrow$ $M \vDash \exists u_1\ zu_1=u$,

$\Rightarrow$ by (QT5), $M \vDash u_1=a \lor u_1=b \lor aBu_1 \lor bBu_1$.

(2cii2a) $M \vDash u_1=b \lor bBu_1$.

$\Rightarrow$ $M \vDash zb=u \lor \exists u_2\ z(bu_2)=u$,

$\Rightarrow$ from (2cii), $M \vDash zaBv\ \&\ zbBv$,

$\Rightarrow$ $M \vDash \exists v_1\ zav_1=v\ \&\ \exists v_2\ zbv_2=v$,

$\Rightarrow$ $M \vDash zav_1=zbv_2$, ,

$\Rightarrow$ by (3.7), $M \vDash av_1=bv_2$, a contradiction.

(2cii2b) $M \vDash u_1=a \lor aBu_1$.

$\Rightarrow$ $M \vDash u_1=a \lor \exists u_2\ au_2=u_1$,

$\Rightarrow$ $M \vDash za=u \lor zaBu$,

$\Rightarrow$ $M \vDash (za=u \lor zaBu)\ \&\ (zb=w \lor zbBw)$,

$\Rightarrow$ $M \vDash u \ll w$.

(3a) $M \vDash ((za=u \lor zaBu)\ \&\ (zb=v \lor zbBv))\ \&\ ((v=a \lor aBv)\ \&\ (w=b \lor bBw))$.

We have that $QT^+ \vdash \neg((zb=v \lor zbBv)\ \&\ v=a)$.

$\Rightarrow$ $M \vDash (zb=v \lor zbBv)\ \&\ aBv$,

$\Rightarrow$ $M \vDash \exists v_1\ (av_1=v\ \&\ (av_1=zb \lor \exists v_2\ av_1=v=zbv_2))$,

$\Rightarrow$ by (QT4) and (QT5), $M \vDash z=a \lor aBz$,

$\implies$ from $M \vDash za=u \lor zaBu$, $M \vDash aBu$,

$\implies M \vDash aBu \;\&\; (w=b \lor bBw)$,

$\implies M \vDash u \ll w$.

(3b)  $M \vDash ((za=u \lor zaBu) \;\&\; (zb=v \lor zbBv)) \;\&\; vBw$.

$\implies M \vDash (zb=v \lor zbBv) \;\&\; vBw$,

$\implies M \vDash zbBw$,

$\implies M \vDash (za=u \lor zaBu) \;\&\; zbBw$,

$\implies M \vDash u \ll w$.

(3c)  $M \vDash \exists z_1((z_1 a=u \lor (z_1 a)Bu) \;\&\; (z_1 b=v \lor (z_1 b)Bv)) \;\&\;$

$$\&\; \exists z_2((z_2 a=v \lor (z_1 b)Bv) \;\&\; (z_2 a=v \lor (z_2 a)Bv)).$$

We have that $QT^+ \vdash \neg(z_1 b=v \;\&\; z_2 a=v)$.

We consider each of the remaining three subcases:

  (3ci)  $M \vDash z_1 b=v \;\&\; (z_2 a)Bv$.

$\implies M \vDash \exists v_2 \; z_1 b=v=z_2 a v_2$,

$\implies$ by (QT4) and (QT5), $M \vDash v_2=b \lor bEv_2$,

$\implies M \vDash z_1 b=z_2 ab \lor \exists v_3 \; z_1 b=z_2 a(v_3 b)$,

$\implies M \vDash z_1=z_2 a \lor z_1=z_2 a v_3$.

  (3ci1)  $M \vDash z_1=z_2 a \;\&\; (z_1 a=u \lor (z_1 a)Bu)$.

$\implies M \vDash ((z_2 a)a=u \lor (z_2 a)aBu) \;\&\; (z_2 b=w \lor (z_2 b)Bw)$,

$\implies M \vDash (z_2 a)Bu \;\&\; (z_2 b=w \lor (z_2 b)Bw)$,

$\implies M \vDash u \ll w$.

  (3ci2)  $M \vDash z_1=z_2 a v_3 \;\&\; (z_1 a=u \lor (z_1 a)Bu)$.

$\Rightarrow$ M ⊨ $(z_2av_3)a$=u v $(z_2av_3)aBu$),

$\Rightarrow$ M ⊨ $(z_2a)Bu$ & $(z_2b$=w v $(z_2b)Bw)$,

$\Rightarrow$ M ⊨ u≪w.

(3cii)  M ⊨ $(z_1b)Bv$ & $z_2a$=v.

$\Rightarrow$ M ⊨ $\exists v_1$ $z_1bv_1$=v=$z_2a$,

$\Rightarrow$ by (QT4) and (QT5),  M ⊨ $v_1$=a v $aEv_1$,

$\Rightarrow$ M ⊨ $z_1ba$=$z_2a$ v $\exists v_2$ $z_1b(v_2a)$=$z_2a$,

$\Rightarrow$ M ⊨ $z_1b$=$z_2$ v $z_1bv_2$=$z_2$,

$\Rightarrow$ from M ⊨ $z_2b$=w v $(z_2b)Bw$,

      M ⊨ $(z_1b)b$=w v $(z_1bv_2)b$=w v $(z_1b)bBw$ v $(z_1bv_2)bBw$,

$\Rightarrow$ M ⊨ $(z_1b)Bw$,

$\Rightarrow$ M ⊨ $(z_1a$=u v $(z_1a)Bu)$ & $(z_1b)Bw$,

$\Rightarrow$ M ⊨ u≪w.

(3ciii)  M ⊨ $(z_1b)Bv$ & $(z_2a)Bv$.

$\Rightarrow$ M ⊨ $z_1Bv$ & $z_2Bv$,

$\Rightarrow$ by (3.8),  M ⊨ $z_1Bz_2$ v $z_1$=$z_2$ v $z_2Bz_1$.

(3ciii1)  M ⊨ $z_1Bz_2$.

$\Rightarrow$ M ⊨ $\exists z_3$ $z_1z_3$=$z_2$,

$\Rightarrow$ from M ⊨ $(z_1a$=u v $(z_1a)Bu)$ & $(z_2b$=w v $z_2Bw)$,

      M ⊨ $(z_1a$=u v $(z_1a)Bu)$ & $((z_1z_3)b$=w v $(z_1z_3)bBw)$,

$\Rightarrow$ by (QT5),  M ⊨ $z_3$=a v $z_3$=b v $aBz_3$ v $bBz_3$.

(3ciii1a)  M ⊨ $z_3$=a v $aBz_3$.

$\Rightarrow$ M ⊨ $(z_1a=u \lor (z_1a)Bu)$ &

& $(z_1ab=w \lor (z_1ab)Bw) \lor \exists z_4\, z_1(az_4)b=w \lor \exists z_4\, z_1(az_4)bBw)$,

$\Rightarrow$ M ⊨ uBw,

$\Rightarrow$ M ⊨ u≪w.

(3ciii1b)  M ⊨ $z_3=b \lor bBz_3$.

$\Rightarrow$ M ⊨ $(z_1a=u \lor (z_1a)Bu)$ &

& $(z_1bb=w \lor (z_1bb)Bw) \lor \exists z_4\, z_1(bz_4)b=w \lor \exists z_4\, z_1(bz_4)bBw)$,

$\Rightarrow$ M ⊨ $(z_1a=u \lor (z_1a)Bu)$ & $(z_1b)Bw$,

$\Rightarrow$ M ⊨ u≪w.

(3ciii2)  M ⊨ $z_1=z_2$.

$\Rightarrow$ M ⊨ $(z_1b)Bv$ & $(z_1a)Bv$,

$\Rightarrow$ M ⊨ $\exists v_1\, z_1bv_1=v$ & $\exists v_2\, z_1av_2=v$,

$\Rightarrow$ M ⊨ $z_1bv_1=z_1av_2$,

$\Rightarrow$ by (3.7),  M ⊨ $bv_1=av_2$, a contradiction.

(3ciii3)  M ⊨ $z_2Bz_1$.

$\Rightarrow$ M ⊨ $\exists z_3\, z_2z_3=z_1$,

$\Rightarrow$ M ⊨ $((z_2z_3)a=u \lor (z_2z_3)aBu)$ & $(z_2b=w \lor (z_2b)Bw)$,

$\Rightarrow$ by (QT5),  M ⊨ $z_3=a \lor z_3=b \lor aBz_3 \lor bBz_3$.

(3ciii3a)  M ⊨ $z_3=a \lor aBz_3$.

$\Rightarrow$ M ⊨ $(z_2aa=u \lor (z_2aa)Bu) \lor \exists z_4\, z_2(az_4)a=u \lor \exists z_4\, z_2(az_4)aBu)$ &

& $(z_2b=w \lor (z_2b)Bw)$,

$\Rightarrow$ M ⊨ $(z_2a)Bu$ & $(z_2b=w \lor (z_2b)Bw)$,

$\Rightarrow$ $M \vDash u \ll w$.

   (3ciii3b)  $M \vDash z_3 = b \lor bBz_3$.

$\Rightarrow$ from $M \vDash z_2 z_3 = z_1$ & $(z_1 b)Bv$, $M \vDash (z_2 z_3)bBv$,

$\Rightarrow$ $M \vDash \exists v_1\ (z_2 z_3)bv_1 = v$.

But we also have from (3ciii) that $M \vDash (z_2 a)Bv$.

$\Rightarrow$ $M \vDash \exists v_2\ z_2 a v_2 = v$,

$\Rightarrow$ $M \vDash (z_2 z_3)bv_1 = z_2 a v_2$,

$\Rightarrow$ by (3.7), $M \vDash z_3 b v_1 = a v_2$,

$\Rightarrow$ by (QT4) and (QT5), $M \vDash z_3 = a \lor aBz_3$,

contradicting the hypothesis (3ciii3b).

This completes the proof of (6.6).

(6.7) For any string concept $I\subseteq I_0$ there is a string concept $J\subseteq I$ such that

$$QT^+ \vdash \forall u,v \in J\ (u \ll v \rightarrow \neg(v \ll u)).$$

Let $J \equiv I_{3.7}\ \&\ I_{3.8}$.

Assume that $M \vDash u \ll v$ where $M \vDash J(u)\ \&\ J(v)$.

We distinguish three cases:

Case 1. $M \vDash (u=a \lor aBu)\ \&\ (v=b \lor bBv)$.

 (1a) $M \vDash u=a\ \&\ v=b$.

$\Rightarrow$ by (QT2), $M \vDash \neg vBu$ and $M \vDash \neg \exists z((za=v \lor zaBv)\ \&\ (zb=u \lor zbBu))$.

Suppose, for a reductio, that $M \vDash v \ll u$.

$\Rightarrow M \vDash (v=a \lor aBv)\ \&\ (u=b \lor bBu)$.

   (1ai)   $M \vDash v=a$.

$\Rightarrow M \vDash a=v=b$, contradicting (QT4).

   (1aii)   $M \vDash aBv$.

$\Rightarrow M \vDash aBb$, contradicting (QT2).

 (1b) $M \vDash u=a\ \&\ bBv$.

We reason as in (1a) and obtain $M \vDash (v=a \lor aBv)\ \&\ (u=b \lor bBu)$.

   (1bi)   $M \vDash u=b$.

$\Rightarrow M \vDash a=u=b$, contradicting (QT4).

   (1bii)   $M \vDash bBu$.

$\Rightarrow M \vDash bBa$, contradicting (QT2).

 (1c) $M \vDash aBu\ \&\ v=b$.

We have $M \vDash \neg vBu$ by (QT4), and again by (QT2),

$M \vDash \neg \exists z((za=v \lor zaBv) \,\&\, (zb=u \lor zbBu))$.

Assuming $M \vDash v \ll u$, we then have

$M \vDash (v=a \lor aBv) \,\&\, (u=b \lor bBu)$,

and we then argue exactly as in (1ai) and (1aii).

(1d)  $M \vDash aBu \,\&\, bBv$.

$\Rightarrow M \vDash \exists u_1\, au_1 = u \,\&\, \exists v_1\, bv_1 = v$.

If $M \vDash vBu$, then $M \vDash \exists v_2\, vv_2 = u$. But then $M \vDash bv_1v_2 = u = au_1$, contradicting (QT4). Therefore $M \vDash \neg vBu$.

Suppose, for a reductio, that $M \vDash (za=v \lor zaBv) \,\&\, (zb=u \lor zbBu)$. By (QT5),

$M \vDash z=a \lor z=b \lor aBz \lor bBz$.

(1di)  $M \vDash z=a \lor aBz$.

$\Rightarrow M \vDash aa=v \lor (aa)Bv \lor \exists z_1(az_1 a = v \lor (az_1 a)Bv)$.

In each case, we obtain a contradiction from $M \vDash bBv$ by (QT4).

(1dii)  $M \vDash z=b \lor bBz$.

$\Rightarrow M \vDash bb=u \lor (bb)Bu \lor \exists z_1(bz_1 b = u \lor (bz_1 b)Bu)$.

In each case, we obtain a contradiction from $M \vDash aBu$ by (QT4).

Therefore, $M \vDash \neg \exists z((za=v \lor zaBv) \,\&\, (zb=u \lor zbBu))$.

Thus, assuming $M \vDash v \ll u$, we have $M \vDash (v=a \lor aBv) \,\&\, (u=b \lor bBu)$. But then we immediately derive a contradiction from hypothesis (1d) using (QT2) and (QT4).

Case 2. $M \vDash uBv$.

$\Rightarrow M \vDash \exists u_1\ uu_1=v$.

Then $M \vDash \neg vBu$, for otherwise from $M \vDash uBv$ & $vBu$ we have $M \vDash uBu$, contradicting $M \vDash u \in I_0$.

Suppose, for a reductio, that $M \vDash (za=v \lor zaBv)$ & $(zb=u \lor zbBu)$.

$\Rightarrow M \vDash za=uu_1 \lor \exists v_1\ zav_1=uu_1$.

 (2a) $M \vDash za=uu_1$.

$\Rightarrow M \vDash za=zbu_1 \lor \exists u_2\ za=(zbu_2)u_1$,

$\Rightarrow$ by (3.7), $M \vDash a=bu_1 \lor a=bu_2u_1$, contradicting (QT2).

 (2b) $M \vDash \exists v_1\ zav_1=uu_1$.

$\Rightarrow M \vDash zav_1=zbu_1 \lor \exists u_2\ zav_1=(zbu_2)u_1$,

$\Rightarrow$ by (3.7), $M \vDash av_1=bu_1 \lor av_1=bu_2u_1$, contradicting (QT4).

Therefore, $M \vDash \neg \exists z((za=v \lor zaBv)$ & $(zb=u \lor zbBu))$.

Thus, assuming $M \vDash v \ll u$, we have $M \vDash (v=a \lor aBv)$ & $(u=b \lor bBu)$.

$\Rightarrow$ by (QT2), $M \vDash aBv$,

$\Rightarrow$ from $M \vDash u=b \lor bBu$ and $M \vDash uBv$, $M \vDash bBv$, contradicting (QT4).

Case 3. $M \vDash \exists z\ (Rt_L(z,u,v)$ & $((za=u \lor zaBu)$ & $(zb=v \lor zbBv)))$.

Assume $M \vDash v \ll u$.

Suppose, for a reductio, that $M \vDash (v=a \lor aBv)$ & $(u=b \lor bBu)$.

If $M \vDash (v=a$ & $u=b) \lor (v=a$ & $bBu) \lor (aBv$ & $u=b)$, we obtain a contradiction from the hypothesis by (QT2).

If $M \vDash aBv$ & $bBu$, we argue analogously to (1d).

Therefore, $M \vDash \neg((v=a \lor aBv) \& (u=b \lor bBu))$.

Suppose that $M \vDash vBu$. We then derive a contradiction analogously to Case 2.

Finally, suppose that

$$M \vDash \exists z_1 (Rt_L(z_1,u,v) \& ((z_1a=v \lor z_1aBv) \& (z_1b=u \lor z_1bBu))).$$

$\Longrightarrow$ $M \vDash Rt_L(z,u,v) \& Rt_L(z_1,u,v)$,

$\Longrightarrow$ by (6.2), $M \vDash z=z_1$,

$\Longrightarrow$ $M \vDash (za=u \lor zaBu) \& (zb=u \lor zbBu)$.

There are four subcases:

(3a) $M \vDash za=u \& zb=u$.

$\Longrightarrow$ $M \vDash za=zb$, contradicting (QT4).

(3b) $M \vDash za=u \& zbBu$.

$\Longrightarrow$ $M \vDash \exists u_2\ za=zbu_2$,

$\Longrightarrow$ by (3.7), $M \vDash a=bu_2$, contradicting (QT2).

(3c) $M \vDash zaBu \& zb=u$.

$\Longrightarrow$ $M \vDash \exists u_1\ (za)u_1=zb$,

$\Longrightarrow$ by (3.7), $M \vDash au_1=b$, contradicting (QT2).

(3d) $M \vDash zaBu \& zb=u$.

$\Longrightarrow$ $M \vDash \exists u_1,u_2\ (zau_1=zbu_2)$,

$\Longrightarrow$ by (20), $M \vDash au_1=bu_2$, contradicting (QT4).

This completes Case 3 and the proof of (6.7).

We let

$$\text{MinMax}^+T_b(t,u) \equiv \text{Max}^+T_b(t,u) \ \&\ \forall t'(\text{Max}^+T_b(t',u) \to t \leq t').$$

We say that t is a <u>shortest non-occurrent b-tally</u> in string u.

We then have:

(6.8) For any string concept $I \subseteq I_0$ there is a string concept $J \subseteq I$ such that

$$QT^+ \vdash \forall x \in J\ \exists! t \in J\ \text{MinMax}^+T_b(t,x).$$

(We read "$\exists! x \in J\ (\ldots)$" as "$\exists x\ (J(x)\ \&\ (\ldots)\ \&\ \forall y(J(y)\ \&\ (\ldots) \to y=x))$").

Let $J \equiv I_{3.12}\ \&\ I_{4.6}\ \&\ I_{\text{MaxT}}$.

Assume $M \vDash J(x)$. Suppose $M \vDash \text{Tally}_a(x)$.

If $M \vDash t=b$, then from (QT4), $M \vDash \neg(b \subseteq_p x)$, and also

$$M \vDash \forall t'(\text{Tally}_b(t')\ \&\ t' \subseteq_p x \to t' \subseteq_p b)$$

because $QT^+ \vdash \text{Tally}_a(x) \to \forall v \subseteq_p x\ \neg\text{Tally}_b(v)$. Hence $M \vDash \text{Max}^+T_b(b,x)$.

If $M \vDash \text{Max}^+T_b(t',x)$, then $M \vDash \text{Tally}_b(t')$, so $M \vDash b \subseteq_p t'$. Therefore,

$M \vDash \forall t'(\text{Max}^+T_b(t',u) \to b \leq t')$, that is, $M \vDash \text{MinMax}^+T_b(b,x)$, as required.

So we may assume that $M \vDash \neg\text{Tally}_a(x)$.

$\Rightarrow$ from $M \vDash J(x)$ by (4.12), $M \vDash \exists z\ (\text{MaxT}_b(z,x)\ \&\ (\neg\text{Tally}_a(x) \to z \subseteq_p x))$,

$\Rightarrow M \vDash z \subseteq_p x$,

$\Rightarrow$ since we may assume that J is closed under $\subseteq_p$, $M \vDash J(z)$,

$\Rightarrow M \vDash J(zb)$.

We claim that $M \vDash \text{MinMax}^+T_b(zb,x)$.

Assume that $M \vDash \text{Tally}_b(t') \,\&\, t' \subseteq_p x$.

$\Rightarrow$ from $M \vDash \text{MaxT}_b(z,x)$, $M \vDash t' \subseteq_p z$,

$\Rightarrow M \vDash t' \subseteq_p z \subseteq_p zb$.

From $M \vDash \text{MaxT}_b(z,x)$, we have $M \vDash \text{Tally}_b(z)$, whence $M \vDash \text{Tally}_b(zb)$ by (4.1).

Suppose, for a reductio, that $M \vDash zb \subseteq_p x$. Then $M \vDash zb \subseteq_p z$ follows from $M \vDash \text{MaxT}_b(z,x)$, contradicting (3.12). Therefore, $M \vDash \neg(zb \subseteq_p x)$, and thus we have established that $M \vDash \text{Max}^+\text{T}_b(zb,x)$.

Assume now that $M \vDash \text{Max}^+\text{T}_b(t',x) \,\&\, t' < zb$.

$\Rightarrow M \vDash \text{Tally}_b(t')$,

$\Rightarrow$ by (1.13), $M \vDash t' \leq z$,

$\Rightarrow M \vDash t' \subseteq_p z \subseteq_p x$, contradicting $M \vDash \text{Max}^+\text{T}_b(t',x)$.

Therefore, $M \vDash \text{Max}^+\text{T}_b(t',x) \to \neg(t' < zb)$.

$\Rightarrow$ by (4.6), $M \vDash \neg(t' < zb) \to zb \leq t'$.

Hence $M \vDash \text{Max}^+\text{T}_b(t',x) \to zb \leq t'$, and we have $M \vDash \text{MinMax}^+\text{T}_b(zb,x)$, as claimed.

Finally, suppose that

$\qquad M \vDash \text{MinMax}^+\text{T}_b(t_1,x) \,\&\, \text{MinMax}^+\text{T}_b(t_2,x) \,\&\, J(t_1) \,\&\, J(t_2)$.

$\Rightarrow M \vDash \text{Max}^+\text{T}_b(t_1,x) \,\&\, \forall t(\text{Max}^+\text{T}_b(t,x) \to t_1 \leq t)$

$\qquad$ and $M \vDash \text{Max}^+\text{T}_b(t_2,x) \,\&\, \forall t(\text{Max}^+\text{T}_b(t,x) \to t_2 \leq t)$,

$\Rightarrow M \vDash t_1 \leq t_2 \,\&\, t_2 \leq t_1$,

$\Rightarrow$ since $M \vDash t_1, t_2 \in J \subseteq I_0$, by (2.2), $M \vDash t_1 = t_2$.

This completes the proof of (6.8).

(6.9)     $QT^+ \vdash MinMax^+T_b(b,x) \leftrightarrow Tally_a(x)$.

Suppose  $M \vDash Tally_a(x)$.

$\Rightarrow M \vDash \neg \exists t(Tally_b(t) \& t\subseteq_p x) \& Tally_b(b) \& \forall t(Tally_b(t) \rightarrow b\leq t)$.

Therefore,  $M \vDash \forall t(t\subseteq_p x \& Tally_b(t) \rightarrow t\leq b) \& \neg b\subseteq_p x$, that is, $M \vDash Max^+T_b(b,x)$.

$\Rightarrow M \vDash \forall t\,(Max^+T_b(t,x) \rightarrow b\leq t)$,

$\Rightarrow M \vDash MinMax^+T_b(b,x)$.

Conversely, suppose  $M \vDash MinMax^+T_b(b,x)$.

Then  $M \vDash \neg b\subseteq_p x$.  But  $M \vDash Digit(b)$.

$\Rightarrow M \vDash \neg \exists y(y\subseteq_p x \& Digit(y) \& y=b)$,

$\Rightarrow M \vDash \forall y(y\subseteq_p x \& Digit(y) \rightarrow y=a)$,

$\Rightarrow M \vDash Tally_a(x)$.

This completes the proof of (6.9).

# 7. The Set Adjunction Lemma

We have the following SET ADJUNCTION LEMMA.

(7.1) For any string concept $I \subseteq I_0$ there is a string concept $J \subseteq I$ such that

$QT^+ \vdash \forall x,y \in J$ (Set(x) & $\neg(y \; \varepsilon \; x)$ $\rightarrow$

$\rightarrow \exists z \in J$ (Set(z) & $y \; \varepsilon \; z$ & $\forall w(w \; \varepsilon \; z \leftrightarrow w \; \varepsilon \; x \; \vee \; w=y))$.

Let $J \equiv I_{5.22}$ & $I_{5.46}$ & $I_{6.8}$.

Assume $M \vDash$ Set(x) & $\neg(y \; \varepsilon \; x)$ where $M \vDash J(x)$ & $J(y)$.

$\Rightarrow$ by (6.8), $M \vDash \exists! t_0 \in J$ MinMax$^+$T$_b(t_0,y)$,

$\Rightarrow$ from $M \vDash$ Set(x), $M \vDash x=aa \; \vee \; \exists t$ Env(t,x).

(1) $M \vDash x=aa$.

Let $z = t_0 a y a t_0$.

$\Rightarrow$ since we may assume that J is closed under *, $M \vDash J(z)$,

$\Rightarrow$ by the choice of $t_0$, $M \vDash$ Pref(aya,$t_0$),

$\Rightarrow M \vDash$ Firstf(z,$t_0$,aya,$t_0$),

$\Rightarrow$ by (5.22), $M \vDash$ Env($t_0$,z) & $\forall w(w \; \varepsilon \; z \leftrightarrow w=y)$,

$\Rightarrow$ by (5.18), from $M \vDash x=aa$, $M \vDash \forall w \neg(w \; \varepsilon \; x)$,

$\Rightarrow M \vDash y \; \varepsilon \; z$ & $\forall w(w \; \varepsilon \; z \leftrightarrow w \; \varepsilon \; x \; \vee \; w=y)$, as required.

(2) $M \vDash \exists t$ Env(t,x).

$\Rightarrow$ by (5.11), $M \vDash \exists t',w'($Tally$_b(t')$ & $x=t'w't$ & aBw' & aEw'),

$\Rightarrow$ from the choice of $t_0$, by (4.5), $M \vDash \text{Tally}_b(tt_0)$ & $\text{Max}^+T_b(tt_0,y)$,

$\Rightarrow$ $M \vDash \text{Pref}(aya,tt_0)$,

$\Rightarrow$ $M \vDash \text{Firstf}(y^*,tt_0,aya,tt_0)$ & $\text{Lastf}(y^*,tt_0,aya,tt_0)$, where $y^*=tt_0ayatt_0$,

$\Rightarrow$ by (5.22), $M \vDash \text{Env}(tt_0,y^*)$ & $\forall w(w \ \varepsilon \ y^* \leftrightarrow w=y)$ & $y \ \varepsilon \ y^*$.

Let $z=t'w'tt_0ayatt_0$.

$\Rightarrow$ from hypothesis $M \vDash \neg(y \ \varepsilon \ x)$, $M \vDash \neg\exists w(w \ \varepsilon \ x \ \& \ w \ \varepsilon \ y^*)$,

$\Rightarrow$ by (5.46), $M \vDash J(z)$ and $M \vDash \text{Env}(tt_0,z)$ & $\forall w(w \ \varepsilon \ z \leftrightarrow w \ \varepsilon \ x \ \vee \ w \ \varepsilon \ y^*)$,

$\Rightarrow$ $M \vDash \text{Set}(z)$ & $y \ \varepsilon \ z$ & $\forall w(w \ \varepsilon \ z \leftrightarrow w \ \varepsilon \ x \ \vee \ w=y)$.

This completes the proof of (7.1).

We let σ(x,y,z) abbreviate

$((y \; \varepsilon \; x \; \& \; z=x) \; v \; (\neg(y \; \varepsilon \; x) \; \& \; \exists t^+(MinMax^+T_b(t^+,y) \; \&$

$\& \; ((x=aa \; \& \; z=t^+ayat^+) \; v \; \exists t(Env(t,x) \; \& \; z=xt^+ayatt^+))))$.

# 8. Tally modified lexicographic ordering

We first introduce a comparison according to length of the shortest non-occurrent b-tallies in given strings u, v:

$$u \triangleleft_{Tb} v \equiv \exists t_1, t_2 \, (\text{MinMax}^+ T_b(t_1, u) \, \& \, \text{MinMax}^+ T_b(t_2, v) \, \& \, t_1 < t_2),$$

and $\quad u \approx_{Tb} v \equiv \exists t_1, t_2 \, (\text{MinMax}^+ T_b(t_1, u) \, \& \, \text{MinMax}^+ T_b(t_2, v) \, \& \, t_1 = t_2).$

We then have:

(8.1) For any string concept $I \subseteq I_0$ there is a string concept $J \subseteq I$ such that

$$QT^+ \vdash \forall u, v \in J ((u \triangleleft_{Tb} v \, \vee \, u \approx_{Tb} v \, \vee \, v \triangleleft_{Tb} u) \, \& \, \neg(u \triangleleft_{Tb} v \, \& \, v \triangleleft_{Tb} u)).$$

Let $J \equiv I_{6.8}$.

Assume $M \vDash J(u) \, \& \, J(v)$.

$\Longrightarrow$ by (6.8), $M \vDash \exists t_1 \in J \, \text{MinMax}^+ T_b(t_1, u) \, \& \, \exists t_2 \in J \, \text{MinMax}^+ T_b(t_2, v)$,

$\Longrightarrow$ by (4.6), since $M \vDash \text{Tally}_b(t_1) \, \& \, \text{Tally}_b(t_2)$, $M \vDash t_1 < t_2 \, \vee \, t_1 = t_2 \, \vee \, t_2 < t_1$,

$\Longrightarrow M \vDash u \triangleleft_{Tb} v \, \vee \, u \approx_{Tb} v \, \vee \, v \triangleleft_{Tb} u$, as required.

Suppose, for a reductio, that $M \vDash u \triangleleft_{Tb} v \, \& \, v \triangleleft_{Tb} u$.

$\Longrightarrow M \vDash \exists t_1, t_2 \in J \, (\text{MinMax}^+ T_b(t_1, u) \, \& \, \text{MinMax}^+ T_b(t_2, v) \, \& \, t_1 < t_2) \, \&$

$\& \, \exists t_3, t_4 \in J \, (\text{MinMax}^+ T_b(t_3, v) \, \& \, \text{MinMax}^+ T_b(t_4, u) \, \& \, t_3 < t_4).$

Now, we have that

$M \vDash \text{MinMax}^+ T_b(t_1, u) \, \& \, \text{MinMax}^+ T_b(t_4, u) \, \rightarrow \, t_1 \leq t_4 \, \& \, t_4 \leq t_1.$

By (2.2), since $M \vDash t_1, t_4 \in J \subseteq I_0$, $M \vDash t_1 = t_4$.

Likewise, $M \vDash t_2 = t_3$.

But then $M \vDash t_1 < t_2$ & $t_2 < t_1$, whence $M \vDash t_1 < t_1$, contradicting $M \vDash t_1 \in I_0$.

This completes the proof of (8.1).

We now set $\quad u \prec v \equiv (u \triangleleft_{T_b} v \lor (u \approx_{T_b} v \ \& \ u \ll v))$.

We call $\prec$ the <u>tally modified lexicographic ordering</u>.

We then have:

(8.2) For any string concept $I \subseteq I_0$ there is a string concept $J \subseteq I$ such that

$$QT^+ \vdash \forall u,v \in J((u \prec v \lor u = v \lor v \prec u) \ \& \ \neg(u \prec v \ \& \ v \prec u)).$$

Let $J \equiv I_{6.5} \ \& \ I_{6.7} \ \& \ I_{6.8}$.

Assume $M \vDash J(u) \ \& \ J(v)$.

$\Rightarrow$ by (8.1), $M \vDash u \triangleleft_{T_b} v \lor u \approx_{T_b} v \lor v \triangleleft_{T_b} u$.

If $M \vDash u \triangleleft_{T_b} v$, then $M \vDash u \prec v$.

If $M \vDash v \triangleleft_{T_b} u$, then $M \vDash v \prec u$.

If $M \vDash u \approx_{T_b} v$, then the desired disjunction follows immediately from (6.5).

Suppose, for a reductio, that $M \vDash u \prec v \ \& \ v \prec u$.

$\Rightarrow M \vDash (u \triangleleft_{T_b} v \lor (u \approx_{T_b} v \ \& \ u \ll v)) \ \& \ (v \triangleleft_{T_b} u \lor (v \approx_{T_b} u \ \& \ v \ll u))$.

We distinguish four cases:

<u>Case 1.</u> $M \vDash u \triangleleft_{T_b} v \ \& \ v \triangleleft_{T_b} u$.

This is ruled out by (8.1).

<u>Case 2.</u> $M \vDash u \triangleleft_{T_b} v \ \& \ (v \approx_{T_b} u \ \& \ v \ll u)$.

$\Rightarrow$ by (6.8), $M \vDash \exists t_1, t_2 \in I_{6.8} (\text{MinMax}^+ T_b(t_1,u) \ \& \ \text{MinMax}^+ T_b(t_2,v) \ \& \ t_1 < t_2)$

and $M \vDash \exists t_3, t_4 \in I_{6.8} (\text{MinMax}^+ T_b(t_3,v) \ \& \ \text{MinMax}^+ T_b(t_4,u) \ \& \ t_3 = t_4)$,

$\Rightarrow$ as in the proof of (8.1), $M \vDash t_1 = t_4 \ \& \ t_2 = t_3$,

$\Rightarrow$ M ⊨ $t_1 < t_1$, contradicting M ⊨ $t_1 \in I_0$.

Case 3.  M ⊨ (u$\approx_{T_b}$v & u≪v) & v◁$_{T_b}$u.

Analogously to Case 2.

Case 4.  M ⊨ (u$\approx_{T_b}$v & u≪v) & (v$\approx_{T_b}$u & v≪u).

This contradicts (6.7).

This completes the proof of (8.2).

(8.3) For any string concept $I \subseteq I_0$ there is a string concept $J \subseteq I$ such that

$$QT^+ \vdash \forall v \in J \, \forall u,w \, (u \prec v \,\&\, v \prec w \rightarrow u \prec w).$$

Let $J \equiv I_{6.6} \,\&\, I_{6.8}$.

Assume $M \vDash u \prec v \,\&\, v \prec w$ where $M \vDash J(v)$.

$\Longrightarrow M \vDash (u \triangleleft_{T_b} v \,\vee\, (u \approx_{T_b} v \,\&\, u \ll v)) \,\&\, (v \triangleleft_{T_b} w \,\vee\, (v \approx_{T_b} w \,\&\, v \ll w))$.

Again, we distinguish four cases:

<u>Case 1.</u> $M \vDash u \triangleleft_{T_b} v \,\&\, v \triangleleft_{T_b} w$.

$\Longrightarrow M \vDash \exists t_1, t_2 (\text{MinMax}^+ T_b(t_1, u) \,\&\, \text{MinMax}^+ T_b(t_2, v) \,\&\, t_1 < t_2) \,\&$

$\qquad \,\&\, \exists t_3, t_4 (\text{MinMax}^+ T_b(t_3, v) \,\&\, \text{MinMax}^+ T_b(t_4, w) \,\&\, t_3 < t_4)$,

$\Longrightarrow$ by (6.8), $M \vDash t_2 = t_3$,

$\Longrightarrow M \vDash t_1 < t_2 = t_3 < t_4$ where $M \vDash \text{Tally}_b(t_1) \,\&\, \text{Tally}_b(t_2) \,\&\, \text{Tally}_b(t_4)$,

$\Longrightarrow M \vDash t_1 < t_4$,

$\Longrightarrow M \vDash u \triangleleft_{T_b} w$,

$\Longrightarrow M \vDash u \prec w$.

<u>Case 2.</u> $M \vDash u \triangleleft_{T_b} v \,\&\, (v \approx_{T_b} w \,\&\, v \ll w)$.

$\Longrightarrow M \vDash \exists t_1, t_2 (\text{MinMax}^+ T_b(t_1, u) \,\&\, \text{MinMax}^+ T_b(t_2, v) \,\&\, t_1 < t_2) \,\&$

$\qquad \,\&\, \exists t_3, t_4 (\text{MinMax}^+ T_b(t_3, v) \,\&\, \text{MinMax}^+ T_b(t_4, u) \,\&\, t_3 = t_4)$,

$\Longrightarrow$ by (6.8), $M \vDash t_2 = t_3$,

$\Longrightarrow M \vDash t_1 < t_2 = t_3 < t_4$,

$\Longrightarrow M \vDash u \triangleleft_{T_b} w$,

$\Longrightarrow M \vDash u \prec w$.

Case 3.  $M \vDash (u \approx_{T_b} v \,\&\, u \ll v) \,\&\, v \triangleleft_{T_b} w$.

Exactly analogous to Case 2.

Case 4.  $M \vDash (u \approx_{T_b} v \,\&\, u \ll v) \,\&\, (v \approx_{T_b} w \,\&\, v \ll w)$.

$\Longrightarrow$ $M \vDash \exists t_1, t_2 (\text{MinMax}^+ T_b(t_1, u) \,\&\, \text{MinMax}^+ T_b(t_2, v) \,\&\, t_1 = t_2)$ &

&  $\exists t_3, t_4 (\text{MinMax}^+ T_b(t_3, v) \,\&\, \text{MinMax}^+ T_b(t_4, u) \,\&\, t_3 = t_4)$,

$\Longrightarrow$ by (6.8), $M \vDash t_2 = t_3$,

$\Longrightarrow$ $M \vDash t_1 = t_2 = t_3 = t_4$,

$\Longrightarrow$ $M \vDash u \approx_{T_b} w$.

But we also have that $M \vDash u \ll v \,\&\, v \ll w$.

$\Longrightarrow$ by (6.6), $M \vDash u \ll w$,

$\Longrightarrow$ $M \vDash u \approx_{T_b} w \,\&\, u \ll w$,

$\Longrightarrow$ $M \vDash u \prec w$.

This completes the proof of (8.3).

# 9. Ordering frames in a set code

We let $u<_x v$ abbreviate the formula

$\exists t_1,t_2,t_3,t_4[Fr(x,t_1,aua,t_2)$ & $Fr(x,t_3,ava,t_4)$ &

  & $((Firstf(x,t_1,aua,t_2)$ & $t_1 \neq t_3)$ v $(Lastf(x,t_3,ava,t_4)$ & $t_1 \neq t_3)$ v

  v $(\exists w_1(Intf(x,w_1,t_1,aua,t_2)$ & $\exists w_3(Intf(x,w_3,t_3,ava,t_4)$ & $t_2=t_3))$ v

  v $(\exists w_1(Intf(x,w_1,t_1,aua,t_2)$ & $\exists w_3(Intf(x,w_3,t_3,ava,t_4)$ & $t_2<t_3)))]$.

If x is a set code, we read $u<_x v$ as saying that $u \, \varepsilon \, x$ and $v \, \varepsilon \, x$ and u's frame precedes v's frame in the string x.

(9.1) For any string concept $I \subseteq I_0$ there is a string concept $J \subseteq I$ such that

$QT^+ \vdash \forall x \in J \, \forall v,t,t'$ $(Set(x)$ & $Firstf(x,t,ava,t')) \to \forall u \, \neg(u<_x v))$.

Let $J \equiv I_{4.23b}$ & $I_{5.15}$ & $I_{5.21}$.

Assume $M \vDash Set(x)$ $M \vDash$ where $M \vDash J(x)$.

$\Rightarrow M \vDash \exists t''\, Env(t'',x)$.

Suppose, for a reductio, that $M \vDash u<_x v$.

$\Rightarrow M \vDash \exists t_1,t_2 Fr(x, t_1,aua,t_2)$.

$\Rightarrow$ from $M \vDash Firstf(x,t,ava,t')$, by (5.19), $M \vDash \neg\exists w_1 \, Intf(x,w_1,t,ava,t')$.

Given the hypothesis $M \vDash u<_x v$, we then have only two cases to consider:

<u>Case 1.</u> $M \vDash Firstf(x,t_1,aua,t_2)$ & $t_1 \neq t$.

$\Rightarrow$ from $M \vDash Firstf(x,t,ava,t')$, by (5.15), $M \vDash u=v$,

$\Rightarrow$ M ⊨ (t₁u)Bx & (t₃v)Bx,

$\Rightarrow$ by (4.23ᵇ), M ⊨ t₁=t, a contradiction.

<u>Case 2.</u>  M ⊨ Lastf(x,t,aua,t') & t₁≠t .

$\Rightarrow$  M ⊨ t=t',

$\Rightarrow$ since M ⊨ t ∈ J ⊆ I₀, M ⊨ ¬(t<t'),

$\Rightarrow$ from M ⊨ Firstf(x,t,ava,t'), M ⊨ x=tavat',

$\Rightarrow$ by (5.21), M ⊨ ∀w (w ε x ↔ w=v),

$\Rightarrow$ from M ⊨ Fr(x,t₁,aua,t₂), M ⊨ u ε x,

$\Rightarrow$ M ⊨ u=v,

But then from M ⊨ Fr(x,t₁,ava,t₂) & Fr(x,t,ava,t'), by (d) of M ⊨ Env(t'',x), we have M ⊨ t₁=t, a contradiction.

This completes the proof of (9.1).

Let  $\text{Lex}^+(z) \equiv \forall u,v\, (u<_z v \to u \prec v)$.

We say that z is <u>lexicographically ordered.</u>

Then we have :

(9.2)  For any string concept $I \subseteq I_0$ there is a string concept $J \subseteq I$ such that

   $QT^+ \vdash \forall x,y[\text{Set}(x)\,\&\,\text{Set}(y)\,\&\,x \sim y\,\&\,\text{Lex}^+(x)\,\&\,\text{Lex}^+(y) \to$

   $\to \forall u,v \in J\ \forall t_1,t_2,t_3,t_4\ ((\text{Firstf}(x,t_1,aua,t_2)\,\&\,\text{Firstf}(y,t_3,ava,t_4))\ \vee$

   $\vee\ (\text{Lastf}(x,t_1,aua,t_2)\,\&\,\text{Lastf}(y,t_3,ava,t_4)) \to u=v)]$.

Let  $J \equiv I_{8.2}$.

Assume  $M \vDash \text{Set}(x)\,\&\,\text{Set}(y)$  where  $M \vDash x \sim y\,\&\,\text{Lex}^+(x)\,\&\,\text{Lex}^+(y)$.

Assume also that  $M \vDash \text{Firstf}(x,t_1,aua,t_2)\,\&\,\text{Firstf}(y,t_3,ava,t_4)$  where

$M \vDash J(u)\,\&\,J(v)$.

Then $M \vDash \exists t\ \text{Env}(t,x)\,\&\,\exists t'\ \text{Env}(t',y)$.

Suppose, for a reductio, that  $M \vDash u \neq v$.

By (8.2), we have that  $M \vDash u \prec v\ \vee\ v \prec u$.

Assume that  $M \vDash u \prec v$.

$\Rightarrow$ from hypothesis $M \vDash x \sim y$,  $M \vDash u\ \varepsilon\ y\,\&\,v\ \varepsilon\ x$,

$\Rightarrow$  $M \vDash \exists t_5,t_6\ \text{Fr}(y,t_5,aua,t_6)$,

$\Rightarrow$ from $M \vDash \text{Fr}(y,t_3,ava,t_4)\,\&\,\text{Fr}(y,t_5,aua,t_6)\,\&\,u \neq v$, by (e) of  $M \vDash \text{Env}(t',y)$,

   $M \vDash t_3 \neq t_5$,

$\Rightarrow$ from $M \vDash \text{Fr}(y,t_3,ava,t_4)$ by definition of  $<_y$, $M \vDash v <_y u$,

⟹ from $M \vDash Lex^+(y)$, $M \vDash v \prec u$.

But this, by (8.2), contradicts the hypothesis $M \vDash u \prec v$.

An analogous argument derives a contradiction from the assumption $M \vDash v \prec u$.

Therefore, $M \vDash Firstf(x,t_1,aua,t_2)$ & $Firstf(y,t_3,ava,t_4) \rightarrow u=v$, as required.

Analogous reasoning establishes that

$$M \vDash Lastf(x,t_1,aua,t_2) \text{ \& } Lastf(y,t_3,ava,t_4) \rightarrow u=v.$$

This completes the proof of (9.2).

(9.3)   $QT^+ \vdash \forall x,v,t_3,t_4$ (Set(x) & Lastf(x,t_3,ava,t_4) $\rightarrow \forall u$ (u $\varepsilon$ x & u$\neq$v $\rightarrow$ u$<_x$v)).

Assume M ⊨ Set(x) & Lastf(x,t_3,ava,t_4), and let M ⊨ u $\varepsilon$ x & u$\neq$v.

$\implies$ by (5.18), M ⊨ x$\neq$a,

$\implies$ M ⊨ $\exists$t Env(t,x) & $\exists t_1,t_2$ Fr(x,t_1,aua,t_2),

$\implies$ by (d) of M ⊨ Env(t,x), M ⊨ $t_1 \neq t_3$,

$\implies$ by definition of $<_x$, M ⊨ u$<_x$v.

This completes the proof of (9.3).

(9.4) For any string concept $I \subseteq I_0$, $QT^+ \vdash \forall x \in I$ (Set(x) $\to \forall u \neg(u \leq_x u))$.

Assume $M \vDash$ Set(x) where $M \vDash I(x)$.

$\Rightarrow$ $M \vDash$ x=aa $\vee \exists t$ Env(t,x).

If $M \vDash$ x=aa, then, by (5.18), $M \vDash \forall u, t_1, t_2 \neg Fr(x, t_1, aua, t_2)$. But then $M \vDash \forall u \neg(u \leq_x u)$, as required.

So we may assume that $M \vDash \exists t$ Env(t,x).

Assume, for a reductio, that $M \vDash u \leq_x u$.

Let $M \vDash Fr(x, t_1, aua, t_2)$ & $Fr(x, t_3, aua, t_4)$.

By (d) of $M \vDash$ Env(t,x), $M \vDash t_1 = t_3$. From the definition of $<_x$, we then immediately have that

$$M \vDash \neg Firstf(x, t_1, aua, t_2) \& \neg Firstf(x, t_3, aua, t_4)$$

and $M \vDash \neg Lastf(x, t_1, aua, t_2) \& \neg Lastf(x, t_3, aua, t_4)$.

Again by the definition of $<_x$, there remain two cases:

<u>Case 1.</u> $M \vDash \exists w_1 (Intf(x, w_1, t_1, aua, t_2) \& \exists w_3 (Intf(x, w_3, t_3, ava, t_4) \& t_2 = t_3))$.

$\Rightarrow M \vDash t_1 < t_2 \& t_3 < t_4$,

$\Rightarrow M \vDash t_3 = t_1 < t_2 = t_3$, contradicting $M \vDash t_3 \in I \subseteq I_0$, since we may assume that I is closed under $\subseteq_p$.

<u>Case 2.</u> $M \vDash \exists w_1 (Intf(x, w_1, t_1, aua, t_2) \& \exists w_3 (Intf(x, w_3, t_3, ava, t_4) \& t_2 < t_3))$.

$\Rightarrow M \vDash t_1 < t_2 < t_3 = t_1$ contradicting $M \vDash t_1 \in I \subseteq I_0$.

This completes the proof of (9.4).

(9.5) For any string concept I⊆I₀ there is a string concept J⊆I such that

$$QT^+ \vdash \forall x \in J \; \forall v,t,t' \; (Set(x) \; \& \; Lastf(x,t,ava,t')) \to \forall u \; \neg(v<_x u)).$$

Let $J \equiv I_{5.20} \; \& \; I_{9.1}$.

Assume $M \vDash Set(x) \; \& \; Lastf(x,t,ava,t')$  where $M \vDash J(x)$.

$\Rightarrow M \vDash Env(t,x) \; \& \; t=t'$.

If $M \vDash u=v$, we have by (9.4) that $M \vDash \neg(v<_x u)$.

So we may assume that $M \vDash u \neq v$.

Assume, for a reductio, that $M \vDash v<_x u$.

$\Rightarrow M \vDash \exists t_1,t_2,t_3,t_4 \; (Fr(x,t_1,ava,t_2) \; \& \; Fr(x,t_3,aua,t_4))$.

$\Rightarrow$ from (d) of $M \vDash Env(t,x)$, $M \vDash t_1 \neq t_3$.

We distinguish three cases:

Case 1.  $M \vDash Firstf(x,t_3,aua,t_4)$.

This is ruled out by hypothesis $M \vDash v<_x u$  by (9.1).

Case 2.  $M \vDash Lastf(x,t_3,aua,t_4)$.

$\Rightarrow$ from $M \vDash Env(t,x)$, $M \vDash t_3=t_4=t$,

$\Rightarrow M \vDash Fr(x,t,ava,t) \; \& \; Fr(x,t,aua,t)$.

But then by (e) of $M \vDash Env(t,x)$, $M \vDash u=v$, which again contradicts the assumption that $M \vDash u \neq v$.

Case 3.  $M \vDash \exists w_3 \; Intf(x,w_3,t_3,aua,t_4)$.

$\Rightarrow$ by (5.19),

  $M \vDash \neg Firstf(x,t_3,aua,t_4) \; \& \; \neg Lastf(x,t_3,aua,t_4) \; \& \; \neg \exists w_1 \; Intf(x,w_1,t_1,ava,t_2)$,

$\Rightarrow M \vDash \text{Firstf}(x,t_1,ava,t_2) \lor \text{Lastf}(x,t_1,ava,t_2)$,

$\Rightarrow$ from $M \vDash \text{Fr}(x,t,ava,t)$ by (d) of $M \vDash \text{Env}(t,x)$, $M \vDash t=t_1$,

$\Rightarrow$ by (a) of $M \vDash \text{Env}(t,x)$, $M \vDash \text{MaxT}_b(t,x)$.

(3a) $M \vDash \text{Firstf}(x,t_1,ava,t_2)$.

$\Rightarrow$ from $M \vDash \neg\text{Firstf}(x,t_3,aua,t_4)$ by (5.20), $M \vDash t_1<t_3$,

$\Rightarrow$ from $M \vDash \text{MaxT}_b(t,x)$, $M \vDash t_3 \leq t$,

$\Rightarrow M \vDash t=t_1<t_2\leq t$, which contradicts $M \vDash t \in J \subseteq I_0$.

(3b) $M \vDash \text{Lastf}(x,t_1,ava,t_2)$.

$\Rightarrow M \vDash t_1=t_2$,

$\Rightarrow$ since $M \vDash t_1 \in J \subseteq I_0$, $M \vDash \neg(t_1<t_2)$,

$\Rightarrow$ from $M \vDash v<_x u \ \& \ \neg\exists w_1 \ \text{Intf}(x,w_1,t_1,ava,t_2) \ \& \ \neg\text{Lastf}(x,t_3,aua,t_4)$,

$\qquad M \vDash \text{Firstf}(x,t_1,ava,t_2)$,

$\Rightarrow M \vDash x=t_1 ava t_2$,

$\Rightarrow$ from $M \vDash \text{Fr}(x,t_3,aua,t_4) \ \& \ \text{MaxT}_b(t,x)$, $M \vDash t_3 \leq t$,

$\Rightarrow$ from $M \vDash \text{Intf}(x,w_3,t_3,aua,t_4)$, $M \vDash \exists w_4 \ x=w_3 at_3 aua t_4 a w_4 \ \& \ \text{Max}^+\text{T}_b(t_3,w_3)$,

$\Rightarrow M \vDash t_1 ava t_2 = x = w_3 at_3 aua t_4 a w_4$,

$\Rightarrow$ from $M \vDash \text{Tally}_b(t_1)$ by (4.14$^b$), $M \vDash w_3=t_1 \lor t_1 B w_3$,

$\Rightarrow M \vDash t_1 \subseteq_p w_3$,

$\Rightarrow$ from $M \vDash \text{Max}^+\text{T}_b(t_3,w_3)$, $M \vDash t_1<t_3$,

$\Rightarrow M \vDash t=t_1<t_3 \leq t$, contradicting $M \vDash t \in J \subseteq I_0$.

This completes the proof of (9.5).

(9.6) For any string concept $I \subseteq I_0$ there is a string concept $J \subseteq I$ such that

$$QT^+ \vdash \forall x \in J\ \forall u,v\ (Set(x) \rightarrow \neg(u<_x v\ \&\ v<_x u)).$$

Let $J \equiv I_{9.5}$.

Assume $M \vDash Set(x)$ where $M \vDash J(x)$.

Suppose, for a reductio, that $M \vDash u<_x v\ \&\ v<_x u$.

$\Rightarrow M \vDash x \neq aa$,

$\Rightarrow M \vDash \exists t\ Env(t,x)$,

$\Rightarrow M \vDash \exists t_1,t_2\ Fr(x,t_1,aua,t_2)\ \&\ \exists t_3,t_4\ Fr(x,t_3,ava,t_4)$.

<u>Case 1.</u> $M \vDash Firstf(x,t_1,aua,t_2)$.

Then $M \vDash \neg(v<_x u)$ by (9.1), which contradicts the hypothesis.

<u>Case 2.</u> $M \vDash Lastf(x,t_3,ava,t_4)$.

Then $M \vDash \neg(v<_x u)$ by (9.5), again contradicting the hypothesis.

<u>Case 3.</u> $M \vDash \exists w_1\ Intf(x,w_1,t_1,aua,t_2)\ \&\ \exists w_3\ Intf(x,w_3,t_3,ava,t_4)$.

$\Rightarrow$ from the definition of $u<_x v$, $M \vDash t_1<t_2 \leq t_3<t_4$,

$\Rightarrow$ from the hypothesis $M \vDash v<_x u$, we have

$$M \vDash \exists t_7,t_8,t_5,t_6\ (Fr(x,t_7,ava,t_8)\ \&\ Fr(x,t_5,aua,t_6)),$$

$\Rightarrow M \vDash Fr(x,t_3,ava,t_4)\ \&\ Fr(x,t_7,ava,t_8)$,

$\Rightarrow$ by (d) of $M \vDash Env(t,x)$, $M \vDash t_3=t_7$.

From the definition of $v<_x u$, we distinguish three subcases:

 (3a) $M \vDash Firstf(x,t_3,ava,t_4)\ \&\ t_3 \neq t_7$.

This is ruled out immediately.

(3b) $M \vDash \text{Lastf}(x,t_7,\text{ava},t_8)$ & $t_3 \neq t_7$.

Likewise.

(3c) $M \vDash \exists w' \text{Intf}(x,t_7,\text{ava},t_8)$ & $\exists w'' \text{Intf}(x,w'',t_5,\text{aua},t_6)$.

$\Rightarrow M \vDash \text{Fr}(x,t_1,\text{aua},t_2)$ & $\text{Fr}(x,t_5,\text{aua},t_6)$,

$\Rightarrow$ by (d) of $M \vDash \text{Env}(t,x)$, $M \vDash t_1 = t_5$,

$\Rightarrow M \vDash t_7 < t_8 \leq t_5 < t_6$,

$\Rightarrow M \vDash t_3 = t_7 < t_5 = t_1$.

But then from $M \vDash t_1 < t_3$ we have $M \vDash t_3 < t_3$, which contradicts $M \vDash t_3 \in J \subseteq I_0$.

This completes the proof of (9.6).

(9.7) For any string concept $I \subseteq I_0$ there is a string concept $J \subseteq I$ such that

$$QT^+ \vdash \forall x \in J \, (Set(x) \to \forall u,v \, (u \, \varepsilon \, x \, \& \, v \, \varepsilon \, x \, \& \, u \neq v \to u <_x v \lor v <_x u)).$$

Let $J \equiv I_{4.6} \, \& \, I_{5.27}$.

Assume $M \vDash Set(x)$ where $M \vDash J(x)$.

Let $M \vDash u \, \varepsilon \, x \, \& \, v \, \varepsilon \, x \, \& \, u \neq v$.

From $M \vDash Set(x)$, by (5.18), we have that $M \vDash \exists t \, Env(t,x)$.

We also have $M \vDash \exists t_1, t_2 \, Fr(x, t_1, aua, t_2) \, \& \, \exists t_3, t_4 \, Fr(x, t_3, ava, t_4)$.

$\Rightarrow$ by (d) of $M \vDash Env(t,x)$, $M \vDash t_1 \neq t_3$.

We distinguish three cases:

<u>Case 1.</u> $M \vDash Firstf(x, t_1, aua, t_2)$.

$\Rightarrow$ by definition of $<_x$, $M \vDash u <_x v$.

<u>Case 2.</u> $M \vDash Lastf(x, t_1, aua, t_2)$.

$\Rightarrow$ by definition of $<_x$, $M \vDash v <_x u$.

<u>Case 3.</u> $M \vDash \exists w_1 Intf(x, w_1, t_1, aua, t_2)$.

(3a) $M \vDash Firstf(x, t_3, ava, t_4)$.

$\Rightarrow$ by definition of $<_x$, $M \vDash v <_x u$.

(3b) $M \vDash Lastf(x, t_3, ava, t_4)$.

$\Rightarrow$ by definition of $<_x$, $M \vDash u <_x v$.

(3c) $M \vDash \exists w_3 Intf(x, w_3, t_3, ava, t_4)$.

$\Rightarrow$ by (4.6), $M \vDash t_2 < t_3 \lor t_2 = t_3 \lor t_3 < t_2$.

If $M \vDash t_2 \leq t_3$, then, by definition of $<_x$, $M \vDash u <_x v$.

Otherwise, $M \vDash t_3 < t_2$, and we have, again by (4.6),

$$M \vDash t_4 < t_1 \vee t_4 = t_1 \vee t_1 < t_4.$$

Now, by (5.2), we have that $M \vDash t_2 \leq t_3 \vee t_4 \leq t_1$, that is, $M \vDash \neg(t_3 < t_2 \,\&\, t_1 < t_4)$.

Therefore $M \vDash t_4 \leq t_1$, and, by definition of x, $M \vDash v <_x u$.

This completes the proof of (9.7).

(9.8) For any string concept $I \subseteq I_0$ there is a string concept $J \subseteq I$ such that

$QT^+ \vdash \forall x \in J \ \forall t, t_1, t_2, v \ (Env(t,x) \ \& \ Fr(x, t_1, ava, t_2) \ \& \ \forall u \ (u \ \varepsilon \ x \rightarrow u \leq_x v) \rightarrow$

$\rightarrow t = t_1 = t_2 \ \& \ Lastf(x, t_1, ava, t_2))$.

Let $J \equiv I_{9.5}$.

Assume $M \vDash Env(t,x) \ \& \ Fr(x, t_1, ava, t_2) \ \& \ \forall u \ (u \ \varepsilon \ x \rightarrow u \leq_x v)$ where $M \vDash J(x)$.

From $M \vDash Env(t,x)$ we have $M \vDash \exists v' \ Lastf(x, t, av'a, t)$.

$\Rightarrow$ by (9.3), $M \vDash v <_x v'$,

$\Rightarrow$ from hypothesis, $M \vDash v' \leq_x v$,

$\Rightarrow$ by (9.6), $M \vDash \neg(v' <_x v)$,

$\Rightarrow M \vDash v' = v$,

$\Rightarrow M \vDash Lastf(x, t, ava, t)$,

$\Rightarrow$ by (d) of $M \vDash Env(t,x)$, $M \vDash t_1 = t$,

$\Rightarrow$ from hypothesis $M \vDash Fr(x, t_1, ava, t_2)$, $M \vDash t_1 \leq t_2$.

If $M \vDash t_1 < t_2$, then from $M \vDash MaxT_b(t,x)$, $M \vDash t = t_1 < t_2 \leq t$, contradicting

$M \vDash t \in I \subseteq I_0$.

Therefore, $M \vDash t = t_1 = t_2$. Hence also $M \vDash Lastf(x, t_1, ava, t_2)$.

This completes the proof of (9.8).

(9.9) For any string concept $I \subseteq I_0$ there is a string concept $J \subseteq I$ such that

$QT^+ \vdash \forall x \in J \ \forall x',t,t_1,t_2,u,v,v'$ $(Env(t,x)$ & $Intf(x,w_1,t_1,ava,t_2)$ &

$\qquad\qquad\qquad\qquad\qquad$ & $x'=w_1at_1avat_1$ & $v' \ \varepsilon \ x'$ & $u <_x v' \rightarrow u \ \varepsilon \ x')$.

Let $J \equiv I_{5.4}$ & $I_{5.27}$ & $I_{5.28}$ & $I_{5.53}$ & $I_{9.5}$.

Assume $M \vDash Env(t,x)$ & $Intf(x,w_1,t_1,ava,t_2)$ & $x'=w_1at_1avat_1$

along with $M \vDash v' \ \varepsilon \ x'$ & $u<_x v'$ and $M \vDash J(x)$.

$\Rightarrow M \vDash \exists w_2 \ x=w_1at_1avat_2aw_2$,

$\Rightarrow M \vDash (w_1at_1)Bx$,

$\Rightarrow$ by (5.53), $M \vDash Env(t_1,x')$,

$\Rightarrow$ from the proof of (5.53), $M \vDash Lastf(x',t_1,ava,t_1)$,

$\Rightarrow$ from hypothesis, $M \vDash u<_x v'$,

$\qquad M \vDash \exists t_3,t_4,t_5,t_6 (Fr(x,t_3,aua,t_4)$ & $Fr(x,t_5,av'a,t_6)$ & $t_3 < t_5)$,

$\Rightarrow$ from hypothesis, $M \vDash v' \ \varepsilon \ x'$, $M \vDash \exists t_7,t_8 \ Fr(x',t_7,av'a,t_8)$,

$\Rightarrow$ from $M \vDash (w_1at_1)Bx$, by (5.51), $M \vDash \exists t_9 \ Fr(x,t_7,av'a,t_9)$,

$\Rightarrow$ from $M \vDash Env(t,x)$, $M \vDash t_5=t_7$,

$\Rightarrow$ from $M \vDash Env(t_1,x')$, $M \vDash MaxT_b(t_1,x')$,

$\Rightarrow M \vDash t_3<t_5=t_7 \leq t_1$.

From $M \vDash Fr(x,t_3,aua,t_4)$ we distinguish three cases:

(1) $M \vDash Firstf(x,t_3,aua,t_4)$.

$\Rightarrow$ from $M \vDash Env(t,x)$ & $Env(t_1,x')$ & $x'Bx$, by (5.4), $M \vDash \exists t_8 \ Firstf(x',t_3,aua,t_8)$,

$\Rightarrow M \vDash u \ \varepsilon \ x'$, as required.

(3) $M \vDash Lastf(x,t_3,aua,t_4)$.

$\Rightarrow$ by (9.5), $M \vDash \neg u<_x v'$, contradicting the hypothesis.

(2) $M \vDash \exists w'\ Intf(x,w',t_3,aua,t_4)$.

Now, we have $M \vDash Intf(x,w_1,t_1,ava,t_2)$.

(2a) $M \vDash t_4=t_1$.

$\Rightarrow$ by (5.28), $M \vDash w'at_3auat_4=w_1at_1$,

$\Rightarrow M \vDash (w'at_3auat_4)avat_1=(w_1at_1)avat_1=x'$,

$\Rightarrow$ from $M \vDash Intf(x,w',t_3,aua,t_4)$,

$\qquad M \vDash Pref(aua,t_3)\ \&\ Tally_b(t_4)\ \&\ t_3<t_4\ \&\ Max^+T_b(t_3,w')$,

$\Rightarrow$ since $M \vDash \exists w''\ x'=w'at_3auat_4aw''$, $M \vDash Intf(x',w',t_3,aua,t_4)$,

$\Rightarrow M \vDash u\ \varepsilon\ x'$, as required.

(2b) $M \vDash t_4 \neq t_1$.

Suppose $M \vDash u=v$.

$\Rightarrow M \vDash Fr(x,t_3,aua,t_4)\ \&\ Fr(x,t_1,aua,t_2)$,

$\Rightarrow$ from $M \vDash Env(t,x)$, $M \vDash t_3=t_1$,

$\Rightarrow M \vDash t_3<t_1=t_3$, contradicting $M \vDash t_3 \in I \subseteq I_0$.

Therefore, $M \vDash u \neq v$.

$\Rightarrow$ by (5.27), $M \vDash t_3auat_4 \subseteq_p w_1\ \vee\ t_1avat_2 \subseteq_p w'$.

If $M \vDash t_1avat_2 \subseteq_p w'$, then $M \vDash t_3<t_1 \subseteq_p w'$, contradicting $M \vDash Max^+T_b(t_3,w')$.

Therefore $M \vDash t_3auat_4 \subseteq_p w_1$.

In fact, from the proof of (5.27), part (2iib), we have

$\qquad M \vDash w'at_3auat_4=w_1\ \vee\ w'at_3auat_4a=w_1\ \vee\ \exists w_2\ w'at_3auat_4aw_2=w_1$,

$\implies$ M ⊨ x'=(w'at₃auat₄)at₁avat₁ v x'=(w'at₃auat₄a)at₁avat₁ v

           v ∃w₂ x'=(w'at₃auat₄aw₂)at₁avat₁,

$\implies$ in each case, M ⊨ ∃w" x'=w'at₃auat₄aw".

Along with  M ⊨ Pref(aua,t₃) & Tally_b(t₄) & t₃<t₄ & Max⁺T_b(t₃,w')

we then have  M ⊨ Intf(x',w',t₃,aua,t₄),  whence  M ⊨ u ε x',  as required.

This completes the proof of (9.9).

(9.10) For any string concept $I \subseteq I_0$ there is a string concept $J \subseteq I$ such that

$QT^+ \vdash \forall x \in J\ \forall t, t_1, t_2, v\ (Env(t,x)\ \&\ Fr(x, t_1, ava, t_2)\ \&\ \forall u\ (u\ \varepsilon\ x \to v \leq_x u)\ \to$

$$\to Firstf(x, t_1, ava, t_2)).$$

Let $J \equiv I_{4.14b}\ \&\ I_{9.1}$.

Assume $M \vDash Env(t,x)\ \&\ Fr(x, t_1, ava, t_2)\ \&\ \forall u\ (u\ \varepsilon\ x \to v \leq_x u)$ where $M \vDash J(x)$.

From $M \vDash Env(t,x)$ we have $M \vDash \exists v'\ \exists t_3, t_4\ Firstf(x, t_3, av'a, t_4)$.

By (9.1), $M \vDash \forall u\ (u\ \varepsilon\ x \to \neg(u <_x v'))$.

$\Rightarrow$ from hypothesis, $M \vDash v \leq_x v'$,

$\Rightarrow M \vDash v = v'$,

$\Rightarrow$ by (d) of $M \vDash Env(t,x)$, $M \vDash t_3 = t_1$.

From $M \vDash Firstf(x, t_1, ava, t_4)$ we have that $M \vDash \neg \exists w\ Intf(x, t_1, ava, t_2)$.

$\Rightarrow$ from $M \vDash Fr(x, t_1, ava, t_2)$, $M \vDash Firstf(x, t_1, ava, t_2) \lor Lastf(x, t_1, ava, t_2)$.

Suppose that $M \vDash Lastf(x, t_1, ava, t_2)$.

$\Rightarrow M \vDash Pref(ava, t_1)\ \&\ Tally_b(t_2)\ \&\ t_1 = t_2\ \&$

$\&\ (x = t_1 avat_2 \lor \exists w(x = wat_1 avat_2\ \&\ Max^+T_b(t_1, w)))$.

From $M \vDash Firstf(x, t_3, av'a, t_4)$ we have $M \vDash (t_3 a)Bx$.

$\Rightarrow M \vDash \exists x_1\ t_3 a x_1 = x = wat_1 avat_2$,

$\Rightarrow$ by $(4.14^b)$, $M \vDash w = t_3 \lor t_3 Bw$,

$\Rightarrow M \vDash t_3 \subseteq_p w$,

$\Rightarrow$ from $M \vDash Max^+T_b(t_1, w)$, $M \vDash t_3 < t_1$,

$\Rightarrow M \vDash t_3 < t_1 = t_3$, contradicting $M \vDash t_3 \in I \subseteq I_0$.

Therefore, $M \vDash t_1 = t_2$ & $x = t_1 \text{ava} t_2$, whence $M \vDash \text{Firstf}(x, t_1, \text{ava}, t_2)$.

This completes the proof of (9.10).

(9.11) For any string concept $I \subseteq I_0$ there is a string concept $J \subseteq I$ such that

$QT^+ \vdash \forall z \in J\ \forall x,t,t',t'',w[Env(t,x)\ \&\ xt'wt''=z\ \&\ aBw\ \&\ aEw\ \&\ Tally_b(t')\ \&$
$\&\ Tally_b(t')\ \to\ \forall u,v(u\ \varepsilon\ x\ \&\ v\ \varepsilon\ x \to (u<_x v \leftrightarrow u<_z v))]$.

Let $J \equiv I_{4.5}\ \&\ I_{9.5}\ \&\ I_{9.7}$.

Assume $M \vDash Env(t,x)\ \&\ xt'wt''=z$

where $M \vDash aBw\ \&\ aEw\ \&\ Tally_b(t')\ \&\ Tally_b(t')$ and $M \vDash J(z)$.

Let $M \vDash u\ \varepsilon\ x\ \&\ v\ \varepsilon\ x$, and assume $M \vDash u<_x v$.

$\Rightarrow M \vDash \exists t_1,t_2\ Fr(x,t_1,ava,t_2)$,

$\Rightarrow$ by (5.6), $M \vDash \exists t_3\ Fr(z,t_1,ava,t_3)$.

(1) $M \vDash Firstf(x,t_1,ava,t_2)$.

$\Rightarrow$ by (9.1), $M \vDash \forall u\ \neg(u<_x v)$, contradicting the hypothesis.

(2) $M \vDash \exists w_1\ Intf(x,w_1,t_1,ava,t_2)$.

$\Rightarrow$ from $M \vDash u<_x v$, $M \vDash u\ \varepsilon\ x$,

$\Rightarrow M \vDash \exists t_4,t_5\ Fr(x,t_4,aua,t_5)$,

$\Rightarrow$ by (9.5), $M \vDash \neg Lastf(x, t_4,aua,t_5)$,

$\Rightarrow M \vDash Firstf(x,t_4,aua,t_5)\ \lor\ \exists w_3 Intf(x,w_3,t_4,aua,t_5)$.

(2a) $M \vDash Firstf(x,t_4,aua,t_5)$.

$\Rightarrow$ from $M \vDash u<_x v$, $M \vDash t_4 \neq t_1$,

$\Rightarrow$ from the proof of (5.6),(1), $M \vDash Firstf(z,t_4,aua,t_5)$,

$\Rightarrow$ from $M \vDash Fr(z,t_1,ava,t_3)$, $M \vDash u<_z v$.

(2b) $M \vDash \exists w_3 Intf(x,w_3,t_4,aua,t_5)$.

$\Rightarrow$ from $M \vDash u<_x v$, $M \vDash t_5 \leq t_1$,

$\Rightarrow$ from the proof of (5.6),(2), $M \vDash \text{Intf}(z,w_3,t_4,aua,t_5)$,

$\Rightarrow$ from $M \vDash \text{Fr}(z,t_1,ava,t_3)$, $M \vDash u<_z v$.

(3) $M \vDash \text{Lastf}(x,t_1,ava,t_2)$.

$\Rightarrow$ from $M \vDash u<_x v$, $M \vDash u \, \varepsilon \, x$,

$\Rightarrow$ $M \vDash \exists t_4, t_5 \, \text{Fr}(x,t_4,aua,t_5)$.

We proceed as in (2).

(3a) $M \vDash \text{Firstf}(x,t_4,aua,t_5)$.

Exactly the same as (2a).

(3b) $M \vDash \exists w_3 \text{Intf}(x,w_3,t_4,aua,t_5)$.

$\Rightarrow$ from $M \vDash \text{Env}(t,x)$, $M \vDash t_4 < t_5 \leq t_1 = t_2 = t$.

From the proof of (5.6),(3), we have, from $M \vDash \text{Lastf}(x,t_1,ava,t_2)$,

$$M \vDash \exists w_5, t_6 \, \text{Intf}(z,w_6,t_1,ava,t_6).$$

On the other hand, from the proof of (5.6),(2), we also have

$M \vDash \text{Intf}(z,w_3,t_4,aua,t_5)$, whence $M \vDash u<_z v$, as required.

So we have shown, under the principal hypothesis, that

$$M \vDash \forall u,v (u<_x v \rightarrow u<_z v).$$

Conversely, assume $M \vDash u<_z v$, where $M \vDash u \, \varepsilon \, x \, \& \, v \, \varepsilon \, x$.

Suppose, for a reductio, that $M \vDash \neg(u<_x v)$.

$\Rightarrow$ by (9.7), $M \vDash v \leq_x u$.

But $M \vDash v=u$ contradicts the hypothesis $M \vDash u<_z v$, by (9.4).

$\Rightarrow$ $M \vDash v<_x u$,

$\Rightarrow$ by the first part of the proof, $M \vDash v<_z u$.

But this contradicts the hypothesis $M \vDash u<_z v$, by (9.6).

This completes the proof of (9.11).

(9.12) For any string concept $I \subseteq I_0$ there is a string concept $J \subseteq I$ such that

$QT^+ \vdash \forall x,x' \in J$ $(Set(x)$ & $Set(x')$ & $x \sim x'$ & $\forall v,w$ $(v \leq_x w \leftrightarrow v \leq_{x'} w)$ $\rightarrow$

$\rightarrow \forall t_1,t_2,w$ $((Firstf(x,t_1,awa,t_2) \rightarrow \exists t_3,t_4$ $Firstf(x',t_3,awa,t_4))$ &

& $(Lastf(x,t_1,awa,t_2) \rightarrow \exists t_3$ $Lastf(x',t_3,awa,t_3))))$.

Let $J \equiv I_{9.7}$ & $I_{9.8}$ & $I_{9.10}$.

Assume $M \vDash Set(x)$ & $Set(x')$ & $x \sim x'$ & $\forall v,w$ $(v \leq_x w \leftrightarrow v \leq_{x'} w)$ where

$M \vDash J(x)$ & $J(x')$.

Assume also that $M \vDash Firstf(x,t_1,awa,t_2)$.

$\implies$ by (9.1) and (9.7), $M \vDash \forall v$ $(v \, \varepsilon \, x \rightarrow w \leq_x v)$,

$\implies$ by hypothesis $M \vDash x \sim x'$, $M \vDash \exists t_3,t_4$ $Fr(x', t_3,awa,t_4)$,

$\implies$ from principal hypothesis, $M \vDash \forall v$ $(v \, \varepsilon \, x' \rightarrow w \leq_{x'} v)$,

$\implies M \vDash w \, \varepsilon \, x'$,

$\implies$ from $M \vDash Set(x')$, $M \vDash x' \neq aa$,

$\implies M \vDash \exists t' Env(t',x')$,

$\implies$ by (9.10), $M \vDash Firstf(x', t_3,awa,t_4)$, as required.

On the other hand, assume $M \vDash Lastf(x,t_1,awa,t_2)$.

$\implies$ by (9.3), $M \vDash \forall v$ $(v \, \varepsilon \, x \rightarrow v \leq_x w)$,

$\implies$ by hypothesis $M \vDash x \sim x'$, $M \vDash \exists t_3,t_4$ $Fr(x', t_3,awa,t_4)$,

$\implies$ from principal hypothesis, $M \vDash \forall v$ $(v \, \varepsilon \, x' \rightarrow v \leq_{x'} w)$,

$\implies$ by (9.8), $M \vDash t_3=t'=t_4$ & $Lastf(x', t_3,awa,t_4)$, as required.

This completes the proof of (9.12).

(9.13) For any string concept $I \subseteq I_0$ there is a string concept $J \subseteq I$ such that

$QT^+ \vdash \forall z \in J \ \forall t,t',t_0,u,v \ (z=xt'auat \ \& \ Env(t_0,x) \ \& \ Env(t,z) \ \& \ Tally_b(t') \ \&$

$\& \ Lastf(x,t_0,ava,t_0) \ \& \ Lastf(z,t,aua,t) \ \to \ v<_z u \ \& \ \neg \exists w(v<_z w \ \& \ w<_z u))$.

Let $J \equiv I_{5.22} \ \& \ I_{5.41} \ \& \ I_{5.46} \ \& \ I_{9.11}$.

Assume $M \vDash z=xt'auat \ \& \ Env(t_0,x) \ \& \ Env(t,z)$

where $M \vDash Tally_b(t') \ \& \ Lastf(x,t_0,ava,t_0) \ \& \ Lastf(z,t,aua,t)$ and $M \vDash J(z)$.

$\Rightarrow M \vDash Fr(x,t_0,ava,t_0)$,

$\Rightarrow$ from $M \vDash Env(t,z), \ M \vDash Tally_b(t)$,

$\Rightarrow$ from $M \vDash Tally_b(t')$, by (5.6), $M \vDash \exists t_1 Fr(z,t_0,ava,t_1)$,

$\Rightarrow$ from $M \vDash Lastf(z,t,aua,t)$, by (9.3), $M \vDash v \leq_z u$.

Suppose, for a reductio, that $M \vDash v=u$.

$\Rightarrow$ from $M \vDash Env(t,z), \ M \vDash t_0=t$,

$\Rightarrow$ from $M \vDash Lastf(x,t_0,ava,t_0), \ M \vDash (at_0)Ex$,

$\Rightarrow M \vDash \exists x_1 \ x=x_1 at_0$,

$\Rightarrow$ from $M \vDash Lastf(z,t,aua,t), \ M \vDash z=tauat \ \lor \ \exists w_1(z=w_1 atauat \ \& \ Max^+T_b(t,w_1))$.

Suppose that $M \vDash z=tauat$.

$\Rightarrow M \vDash xt'auat=z=tauat$,

$\Rightarrow$ by (3.6), $M \vDash xt'=t$,

$\Rightarrow M \vDash (x_1 at_0)t'=t$, which contradicts $M \vDash Tally_b(t)$.

Suppose that $M \vDash \exists w_1(z=w_1 atauat \ \& \ Max^+T_b(t,w_1))$.

$\Rightarrow M \vDash w_1 atauat=z=xt'auat$,

$\Rightarrow$ by (3.6), $M \vDash w_1at=xt'=(x_1at_0)t'$,

$\Rightarrow$ from $M \vDash Tally_b(t_0)$ & $Tally_b(t')$, by (4.5), $M \vDash Tally_b(t_0t')$,

$\Rightarrow$ by (4.24$^b$), $M \vDash t=t_0t'$,

$\Rightarrow$ from $M \vDash t=t_0$, $M \vDash tBt$, contradicting $M \vDash t \in I \subseteq I_0$.

Therefore, $M \vDash v \neq u$.

This establishes that $M \vDash v<_z u$.

$\Rightarrow$ from $M \vDash Env(t,z)$, $M \vDash t_0 \neq t$,

$\Rightarrow$ from $M \vDash Env(t,z)$, $M \vDash MaxT_b(t,z)$,

$\Rightarrow$ from $M \vDash t \subseteq_p x \subseteq_p z$ & $M \vDash Tally_b(t_0)$, $M \vDash t_0<t$,

$\Rightarrow$ from $M \vDash Env(t_0,x)$, by (5.11),

$$M \vDash \exists t^*,w'(Tally_b(t^*) \& aBw' \& aEw' \& x=t^*w't_0).$$

Let $u^+=tauat$.

$\Rightarrow$ from $M \vDash Lastf(z,t,aua,t)$, $M \vDash Pref(aua,t)$,

$\Rightarrow$ $M \vDash Firstf(u^+,t,aua,t)$ & $Lastf(u^+,t,aua,t)$,

$\Rightarrow$ by (5.22), $M \vDash Env(t,u^+)$ & $\forall w(w \, \varepsilon \, u^+ \leftrightarrow w=u)$.

We have that $M \vDash z=t^*w't_0t'auat$ & $x=t^*w't_0$.

We have, as above, from $M \vDash Lastf(z,t,aua,t)$, that $M \vDash \exists w_1 z=w_1atauat$,

whence again, just as above, it follows that $M \vDash t_0t'=t$.

Therefore, $M \vDash u^+=tauat=t_0t'auat$.

$\Rightarrow$ by (5.41), $M \vDash \neg \exists w(w \, \varepsilon \, x \, \& \, w \, \varepsilon \, u^+)$,

$\Rightarrow$ by (5.46), $M \vDash \forall w(w \, \varepsilon \, z \leftrightarrow (w \, \varepsilon \, x \, \vee \, w \, \varepsilon \, u^+))$.

Suppose, for a reductio, that $M \vDash \exists w(v<_z w \, \& \, w<_z u)$.

$\Rightarrow$ by (9.4), $M \vDash w \, \varepsilon \, z \, \& \, w{\neq}v \, \& \, w{\neq}u$,

$\Rightarrow M \vDash w \, \varepsilon \, x \, \& \, w{\neq}v$,

$\Rightarrow$ from $M \vDash \text{Lastf}(x,t_0,ava,t_0)$, by (9.3), $M \vDash w<_x v$,

$\Rightarrow$ by (9.11), $M \vDash w<_z v$, which contradicts hypothesis $M \vDash v<_z w$, by (9.6).

Therefore, $M \vDash \neg \exists w(v<_z w \, \& \, w<_z u)$, as claimed.

This completes the proof of (9.13).

(9.14) For any string concept $I \subseteq I_0$ there is a string concept $J \subseteq I$ such that

$QT^+ \vdash \forall x \in J \; \forall t_1, t_2, t', t'' \; (Set(x) \; \& \; Fr(x, t_1, ava, t_2) \; \& \; Fr(x, t', aua, t'')) \to$

$$\to \forall u, v \; (u <_x v \leftrightarrow t' < t_1)).$$

Let $J \equiv I_{5.19} \; \& \; I_{5.20} \; \& \; I_{9.1} \; \& \; I_{9.4} \; \& \; I_{9.5} \; \& \; I_{9.7}$.

Assume $M \vDash Set(x) \; \& \; Fr(x, t_1, ava, t_2)$ where $M \vDash J(x)$.

By (5.18), $M \vDash \exists t \; Env(t, x)$.

We first show that $M \vDash \forall u \; \forall t', t'' (Fr(x, t', aua, t'') \to (u <_x v \to t' < t_1))$.

<u>Case 1.</u> $M \vDash Firstf(x, t_1, ava, t_2)$.

By (9.1), $M \vDash \forall u \; \neg(u <_x v)$. Hence the claim holds immediately.

<u>Case 2.</u> $M \vDash \exists w_1 Intf(x, w_1, t_1, ava, t_2)$.

(2a) $M \vDash u <_x v \; \& \; Firstf(x, t', aua, t'')$.

By (5.19), we have that $M \vDash \neg Firstf(x, t_1, ava, t_2)$. By (5.20), it follows that $M \vDash t' < t_1$.

(2b) $M \vDash u <_x v \; \& \; \exists w' \; Intf(x, w', t', aua, t'')$.

$\Rightarrow$ by definition of $<_x$, $M \vDash t'' \leq t_1$,

$\Rightarrow$ $M \vDash t' < t'' \leq t_1$, as required.

(2c) $M \vDash u <_x v \; \& \; Lastf(x, t', aua, t'')$.

This is ruled out by (9.5).

<u>Case 3.</u> $M \vDash Lastf(x, t_1, ava, t_2)$.

Suppose $M \vDash u <_x v \; \& \; Fr(x, t', aua, t'')$.

From (a) of $M \vDash Env(t, x)$, we have $M \vDash MaxT_b(t, x)$.

$\Rightarrow$ M ⊨ t'≤t=$t_1$=$t_2$.

By (9.4), from M ⊨ u<$_x$v we have M ⊨ u≠v.

$\Rightarrow$ by (d) of M ⊨ Env(t,x), M ⊨ t'≠$t_1$,

$\Rightarrow$ M ⊨ t'<$t_1$, as required.

We have thus shown that

  M ⊨ Set(x) & Fr(x, $t_1$,v,$t_2$) & Fr(x,t',aua,t'') → (u<$_x$v → t'<$t_1$).

Conversely, assume M ⊨ Set(x) & Fr(x, $t_1$,v,$t_2$) & Fr(x,t',u,t'') and that

M ⊨ t'<$t_1$.

Suppose, for a reductio, that M ⊨ ¬(u<$_x$v).

$\Rightarrow$ by (9.7), M ⊨ v≤$_x$u.

If M ⊨ v=u, then from M ⊨ Env(t,x), M ⊨ $t_1$=t'. But then M ⊨ t'<$t_1$=t',

contradicting M ⊨ t'∈I ⊆ $I_0$. Hence M ⊨ v≠u, and M ⊨ v<$_x$u. But then, by the

first part of the proof, M ⊨ $t_1$<t', whence M ⊨ $t_1$<t'<$t_1$, contradicting

M ⊨ $t_1$∈I ⊆ $I_0$.

Therefore, M ⊨ u<$_x$v, as required.

This completes the proof of (9.14).

(9.15) For any string concept $I\subseteq I_0$ there is a string concept $J\subseteq I$ such that

$QT^+ \vdash \forall z\in J\ \forall x,t,t_1,t_2,w_1,v(Env(t,z)\ \&\ Intf(z,w_1,t_1,ava,t_2)\ \&\ x=w_1at_1avat_1 \rightarrow$

$\rightarrow \forall u(u<_z v \leftrightarrow u<_x v))$.

Let $J \equiv I_{9.14}$.

Assume $M \vDash Env(t,z)\ \&\ Intf(z,w_1,t_1,ava,t_2)$.

Let $M \vDash x=w_1at_1avat_1\ \&\ J(z)$.

$\Rightarrow$ by (5.53), $M \vDash Env(t_1,x)\ \&\ Lastf(x,t_1,ava,t_1)$,

$\Rightarrow M \vDash v\ \varepsilon\ x$.

Assume that $M \vDash u<_z v$.

$\Rightarrow M \vDash u\ \varepsilon\ z$,

$\Rightarrow M \vDash \exists t_3,t_4\ Fr(z,t_3,aua,t_4)$,

$\Rightarrow$ from $M \vDash u<_z v$, by (9.14), $M \vDash t_3<t_1$.

We also have

$\qquad M \vDash Env(t,z)\ \&\ Intf(z,w_1,t_1,ava,t_2)\ \&\ x=w_1at_1avat_1\ \&\ v\ \varepsilon\ x\ \&\ u<_z v$.

$\Rightarrow$ by (9.9), $M \vDash u\ \varepsilon\ x$,

$\Rightarrow M \vDash \exists t_5,t_6\ Fr(x,t_5,aua,t_6)$,

$\Rightarrow$ since $M \vDash (w_1at_1)Bx$, by (5.51), $M \vDash \exists t_7\ Fr(z,t_5,aua,t_7)$,

$\Rightarrow$ from $M \vDash Env(t,z)$, $M \vDash t_3=t_5$,

$\Rightarrow M \vDash t_5=t_3<t_1$,

$\Rightarrow$ from $M \vDash Fr(x,t_5,aua,t_6)\ \&\ Fr(x,t_1,ava,t_1)$, by (9.14), $M \vDash u<_x v$, as required.

Conversely, assume $M \vDash u<_x v$.

$\Rightarrow$ by (9.14), $M \vDash \exists t_3,t_4,t_5,t_6\ (Fr(x,t_3,aua,t_4)\ \&\ Fr(x,t_5,ava,t_6)\ \&\ t_3<t_5)$,

$\Rightarrow$ by (5.51), $M \vDash \exists t_7, t_8 \, (Fr(z,t_3,aua,t_7) \,\&\, Fr(z,t_5,ava,t_8) \,\&\, t_3 < t_5)$,

$\Rightarrow$ by (9.14), $M \vDash u <_z v$, as required.

This completes the proof of (9.15).

(9.16) For any string concept $I \subseteq I_0$ there is a string concept $J \subseteq I$ such that

$$QT^+ \vdash \forall x \in J \; \forall u,v,w \; (Set(x) \; \& \; u<_x v \; \& \; v<_x w \; \rightarrow \; u<_x w \,).$$

Let $J \equiv I_{9.14}$.

Assume $M \vDash Set(x) \; \& \; u<_x v \; \& \; v<_x w$ where $M \vDash J(x)$.

$\Rightarrow \; M \vDash u \, \varepsilon \, x \; \& \; v \, \varepsilon \, x \; \& \; w \, \varepsilon \, x$,

$\Rightarrow \; M \vDash \exists t_1, t_2 \; Fr(x, t_1, aua, t_2) \; \& \; \exists t_3, t_4 \; Fr(x, t_3, ava, t_4) \; \& \; \exists t_5, t_6 \; Fr(x, t_5, awa, t_6)$,

$\Rightarrow \;$ by (9.14), $M \vDash t_1 < t_3 \; \& \; t_3 < t_5 \; \& \; Tally_b(t_1) \; \& \; Tally_b(t_3)$,

$\Rightarrow \; M \vDash t_1 B t_3 \; \& \; t_3 B t_5$,

$\Rightarrow \; M \vDash t_1 B t_5$,

$\Rightarrow \; M \vDash t_1 < t_5$,

$\Rightarrow \;$ by (9.14), $M \vDash u<_x w$.

This completes the proof of (9.16).

(9.17) For any string concept I⊆I₀ there is a string concept J⊆I such that

  QT⁺ ⊢ ∀x∈J ∀u,v,t',t'',t₁,t₂ (Set(x) & Free⁺(x,t',ava,t'') & Fr(x,t₁,aua,t₂) & u≤$_x$v

→

                             → Free⁺(x,t₁,aua,t₂)).

Let  J ≡ I$_{9.16}$.

Assume  M ⊨ Set(x) & Free⁺(x,t',ava,t'') & Fr(x,t₁,aua,t₂) & u≤$_x$v

where  M ⊨ J(x).

Assume further that   M ⊨ Fr(x,t₃,au'a,t₄) & u'≤$_x$u.

⟹  by (9.16),  M ⊨ u'≤$_x$v,

⟹  from hypothesis  M ⊨ Free⁺(x,t',ava,t''),

      M ⊨ Firstf(x,t₃,au'a,t₄) ∨ Free(x,t₃,au'a,t₄).

But then we have shown that

  M ⊨ ∀u',t₃,t₄(Fr(x,t₃,au'a,t₄) & u'≤$_x$u  →  Firstf(x,t₃,au'a,t₄) ∨ Free(x,t₃,au'a,t₄)),

that is,  M ⊨ Free⁺(x,t₁,aua,t₂),  as required.

This completes the proof of (9.17).

(9.18) For any string concept $I \subseteq I_0$ there is a string concept $J \subseteq I$ such that

$QT^+ \vdash \forall x \in J \forall w,z,x',t,t',t'',t_3,t_4,v_0[\text{Env}(t,x)$ & $x=t'wz$ & $\text{Env}(t'',x')$ & $x'=t'wt''$ &

& $aBw$ & $aEw$ & $\text{Env}(t,z)$ & $\text{Firstf}(z,t_3,av_0a,t_4)$ & $t''<t_3 \rightarrow$

$\rightarrow \forall u,v(u \, \varepsilon \, z$ & $v \, \varepsilon \, z \rightarrow (u<_z v \leftrightarrow u<_x v))]$.

Let $J \equiv I_{5.25}$ & $I_{5.46}$ & $I_{9.14}$.

Assume

$M \vDash \text{Env}(t,x)$ & $x=t'wz$ & $\text{Env}(t,z)$  where

$M \vDash x'=t'wt''$ & $\text{Env}(t'',x')$ & $aBw$ & $aEw$  and  $M \vDash \text{Firstf}(z, t_3,av_0a,t_4)$ & $t''<t_3$.

Let $M \vDash J(x)$, and assume $M \vDash u \, \varepsilon \, z$ & $v \, \varepsilon \, z$.

Suppose that $M \vDash u<_z v$.

$\Rightarrow M \vDash \exists t_1,t_2 \text{Fr}(z,t_1,aua,t_2)$ & $\exists t_5,t_6 \text{Fr}(z,t_5,ava,t_6)$,

$\Rightarrow$ from $M \vDash \text{Env}(t,z)$ & $u<_z v$, by (9.14), $M \vDash t_1<t_5$,

$\Rightarrow$ as in the proof of (5.46), $M \vDash \text{Max}^+T_b(t_1,t'w)$ & $\text{Max}^+T_b(t_5,t'w)$,

$\Rightarrow$ by (5.25), $M \vDash \text{Fr}(x,t_1,aua,t_2)$ & $\text{Fr}(x,t_5,ava,t_6)$,

$\Rightarrow$ since $M \vDash M \vDash t_1<t_5$, by (9.14), $M \vDash u<_x v$, as required.

This shows, under the principal hypothesis, that $M \vDash \forall u,v(u<_z v \rightarrow u<_x v)$.

Conversely, assume $M \vDash u<_x v$.

Suppose, for a reductio, that $M \vDash \neg(u<_z v)$.

$\Rightarrow$ by (9.7), $M \vDash v \leq_z u$.

Now, $M \vDash v=u$ contradicts hypothesis $M \vDash u<_x v$ by (9.4).

$\Rightarrow M \vDash v<_z u$,

$\Rightarrow$ by the first part of the proof, $M \vDash v<_x u$.

But this contradicts hypothesis $M \vDash u<_x v$ by (9.6).

Therefore, $M \vDash \forall u,v(u<_x v \rightarrow u<_z v)$.

This completes the proof of (9.18).

(9.19) For any string concept $I \subseteq I_0$ there is a string concept $J \subseteq I$ such that

$QT^+ \vdash \forall x \in J \forall t, t_1, t_2, t_3, u, v, x_2, x^-(Env(t,x)\ \&\ x=t_1auat_2ax_2\ \&\ Firstf(x,t_1,aua,t_2)\ \&$

$\&\ x^-=t_2ax_2\ \&\ Env(t,x^-)\ \&\ Firstf(x^-,t_2,ava,t_3)\ \rightarrow\ \forall w(w <_x v \rightarrow w=u))$.

Let $J \equiv I_{5.22}\ \&\ I_{5.41}\ \&\ I_{9.18}$.

Assume

$M \vDash Env(t,x)\ \&\ x=t_1auat_2ax_2\ \&\ Firstf(x,t_1,aua,t_2)\ \&\ x^-=t_2ax_2\ \&\ Env(t,x^-)\ \&$

$\&\ Firstf(x^-,t_2,ava,t_3)$

where $M \vDash J(x)$.

Let $x_0 = t_1auat_1$.

$\Rightarrow$ from $M \vDash Firstf(x,t_1,aua,t_2)$,

$M \vDash Pref(aua,t_1)\ \&\ Tally_b(t_2)\ \&\ ((t_1=t_2\ \&\ x=t_1auat_2)\ v\ (t_1<t_2\ \&\ (t_1auat_2a)Bx))$.

Assume, for a reductio, that $M \vDash t_1=t_2\ \&\ x=t_1auat_2$.

$\Rightarrow M \vDash t_1auat_2=x=t_1auat_2ax_2$,

$\Rightarrow M \vDash xBx_1$, contradicting $M \vDash x \in I \subseteq I_0$.

Therefore $M \vDash t_1 < t_2$.

$\Rightarrow$ from $M \vDash Env(t,x^-)$, by (5.11),

$M \vDash \exists t_0, w_0 (Tally_b(t_0)\ \&\ x^-=t_0w_0t\ \&\ aBw_0\ \&\ aEw_0)$,

$\Rightarrow M \vDash \exists w_1\ t_2ax_2=x^-=t_0aw_1t$,

$\Rightarrow$ by (4.23$^b$), $M \vDash t_2=t_0$,

$\Rightarrow M \vDash x^-=t_2w_0t$,

$\Rightarrow$ from $M \vDash t_1<t_2\ \&\ Tally_b(t_1)\ \&\ Tally_b(t_2)$, $M \vDash \exists t'(Tally_b(t')\ \&\ t_1t'=t_2)$,

$\Rightarrow$ M ⊨ $x = t_1auat_2ax_2 = t_1auax^- = t_1auat_2w_0t = t_1auat_1t'w_0t$,

$\Rightarrow$ from M ⊨ $Pref(aua, t_1)$, M ⊨ $Firstf(x_0, t_1, aua, t_1)$ & $Lastf(x_0, t_1, aua, t_1)$,

$\Rightarrow$ by (5.22), M ⊨ $Env(t_1, x_0)$ & $\forall w(w \, \varepsilon \, x_0 \leftrightarrow w = u)$.

So we have

  M ⊨ $Env(t,x)$ & $x = t_1auat_1t'w_0t$ & $Tally_b(t_1)$ & $aBw_0$ & $aEw_0$ & $x_0 = t_1auat_1$ &

    & $x^- = t_1t'w_0t$ & $Env(t_1, x_0)$ & $Env(t, x^-)$ & $Firstf(x^-, t_1t', ava, t_3)$.

$\Rightarrow$ by (5.41), M ⊨ $\neg \exists w(w \, \varepsilon \, x_0 \, \& \, w \, \varepsilon \, x^-)$,

$\Rightarrow$ by (5.46), M ⊨ $\forall w(w \, \varepsilon \, x \leftrightarrow w \, \varepsilon \, x_0 \vee w \, \varepsilon \, x^-)$.

Assume now that  M ⊨ $w <_x v$.

$\Rightarrow$ M ⊨ $w \, \varepsilon \, x$,

$\Rightarrow$ M ⊨ $w \, \varepsilon \, x_0 \vee w \, \varepsilon \, x^-$.

Suppose, for a reductio, that  M ⊨ $w \, \varepsilon \, x^-$.

$\Rightarrow$ M ⊨ $\exists t_4, t_5 Fr(x^-, t_4, awa, t_5)$,

$\Rightarrow$ from M ⊨ $Firstf(x^-, t_2, ava, t_3)$, by (9.1) and (9.7),  M ⊨ $v \leq_x w$,

$\Rightarrow$ by (9.18), M ⊨ $v \leq_x w$.

But this contradicts the hypothesis  M ⊨ $w <_x v$, by (9.4) and (9.6).

Therefore  M ⊨ $w \, \varepsilon \, x_0$.

$\Rightarrow$ M ⊨ $w = u$.

Hence  M ⊨ $\forall w(w <_x v \rightarrow w = u)$,  as required.

This completes the proof of (9.19).

(9.20) For any string concept $I \subseteq I_0$ there is a string concept $J \subseteq I$ such that

$QT^+ \vdash \forall x,y \in J$ (Set(x) & Set(z) & $Lex^+(x)$ & $Lex^+(y)$ &

     & $\forall v,w$ (v ε x & w ε x & v ε y & w ε y) $\rightarrow$ ($v<_x w \leftrightarrow v<_y w$))).

Let $J \equiv I_{8.2}$ & $I_{9.7}$.

Assume $M \vDash$ Set(x) & Set(y) & $Lex^+(x)$ & $Lex^+(y)$ where $M \vDash J(x)$ & $J(y)$, and

let $M \vDash$ v ε x & w ε x & v ε y & w ε y.

Suppose, for a reductio, that $M \vDash v<_x w$ & $\neg(v<_y w)$.

$\Rightarrow$ by (9.7), $M \vDash w \leq_y v$,

$\Rightarrow$ from $M \vDash Lex^+(y)$, $M \vDash w \lessdot v$,

$\Rightarrow$ from $M \vDash v<_x w$ & $Lex^+(x)$, $M \vDash v \lessdot w$.

But from $M \vDash w \lessdot v$ & $v \lessdot w$ we have a contradiction, by (8.2).

Therefore, $M \vDash v<_x w \rightarrow v<_y w$.

A completely analogous argument shows that $M \vDash v<_y w \rightarrow v<_x w$.

This completes the proof of (9.20).

(9.21) For any string concept $I \subseteq I_0$ there is a string concept $J \subseteq I$ such that

$QT^+ \vdash \forall x,z \in J \; \forall y,t,v (Set(x) \; \& \; Set(z) \; \& \; Lex^+(x) \; \& \; Lex^+(z) \;\&$

$\quad \& \; \forall w (w \, \varepsilon \, z \leftrightarrow (w \, \varepsilon \, x \lor w=y)) \; \& \; Lastf(x,t,ava,t) \; \& \; y \prec v \;\rightarrow$

$\qquad\qquad\qquad\qquad\qquad\qquad\qquad \rightarrow \; \exists t' \, Lastf(z,t',ava,t')).$

Let $J \equiv I_{9.6} \; \& \; I_{9.20}$.

Assume $M \vDash Set(x) \; \& \; Set(z) \; \& \; Lex^+(x) \; \& \; Lex^+(z)$

where $M \vDash \forall w (w \, \varepsilon \, z \leftrightarrow (w \, \varepsilon \, x \lor w=y)) \; \& \; Lastf(x,t,ava,t) \; \& \; y \prec v$

and $M \vDash J(x) \; \& \; J(z)$.

$\Rightarrow$ from $M \vDash Set(z) \; \& \; z \neq aa$, $M \vDash \exists t' \, Env(t',z)$,

$\Rightarrow M \vDash \exists v' \, Lastf(z,t',av'a,t')$.

$\Rightarrow$ since $M \vDash v \, \varepsilon \, x$, $M \vDash v \, \varepsilon \, z$,

$\Rightarrow$ by (9.3), $M \vDash v \leq_z v'$,

$\Rightarrow$ from $M \vDash v' \, \varepsilon \, z$, $M \vDash v' \, \varepsilon \, x \lor v'=y$.

Suppose $M \vDash v' \, \varepsilon \, x$.

$\Rightarrow$ from $M \vDash Lastf(x,t,ava,t)$ by (9.3), $M \vDash v' \leq_x v$,

$\Rightarrow$ by (9.20), $M \vDash v' \leq_z v$.

Suppose $M \vDash v'=y$.

Assume, for a reductio, that $M \vDash v <_z v'=y$.

Then from $M \vDash Lex^+(z)$, $M \vDash v \prec y$. But this contradicts $M \vDash y \prec v$ by (8.2).

Hence $M \vDash \neg(v <_z y)$.

$\Rightarrow$ by (9.7), $M \vDash y \leq_z v$,

$\Rightarrow$ M ⊨ v'$\leq_z$v.

So we have, in either case, that M ⊨ v'$\leq_z$v.

By (9.6), this, together with M ⊨ v$\leq_z$v', yields M ⊨ v'=v.

$\Rightarrow$ M ⊨ Lastf(z,t',ava,t').

This completes the proof of (9.21).

(9.22) For any string concept $I \subseteq I_0$ there is a string concept $J \subseteq I$ such that

$QT^+ \vdash \forall z \in J\ \forall w, t_1, t_2 (Fr(z, t_1, awa, t_2)\ \&\ \neg Firstf(z, t_1, awa, t_2) \rightarrow$

$\rightarrow Free(z, t_1, awa, t_2)\ v\ Bound(z, t_1, awa, t_2))$.

Let $J \equiv I_{9.14}$.

Assume $M \vDash Fr(z, t_1, awa, t_2)\ \&\ \neg Firstf(z, t_1, awa, t_2)$ where $M \vDash J(z)$.

Assume now that $M \vDash Fr(z, t', ava, t'')\ \&\ v <_z w\ \&\ \neg \exists u (u <_z v\ \&\ v <_z w)$.

$\Rightarrow$ from $M \vDash v <_z w$, by (9.14), $M \vDash t' < t_1 \leq t_2$,

$\Rightarrow$ by (1.13), $M \vDash t'b \leq t_1$.

Hence, if $M \vDash t'b = t_1$, then $M \vDash Free(z, t_1, awa, t_2)$, and if $M \vDash t'b < t_1$, then

$M \vDash Bound(z, t_1, awa, t_2)$.

This completes the proof of (9.22).

(9.23) For any string concept $I \subseteq I_0$ there is a string concept $J \subseteq I$ such that

$QT^+ \vdash \forall z \in J \, \forall x,t,t',t_0,u(Env(t,z)$ & $z=xt'auat$ & $Tally_b(t')$ & $Env(t_0,x)$ &

  & $Lastf(z,t,aua,t) \rightarrow$

    $\rightarrow \forall w,t_1(w \, \varepsilon \, x \rightarrow (\exists t_2 Free(x,t_1,awa,t_2) \leftrightarrow \exists t_3 Free(z,t_1,awa,t_3))$ &

      & $(\exists t_2 Bound(x,t_1,awa,t_2) \leftrightarrow \exists t_3 Bound(z,t_1,awa,t_3)))))$.

Let $J \equiv I_{5.13}$ & $I_{.5.22}$ & $I_{5.46}$ & $I_{9.11}$ & $I_{9.16}$.

Assume $M \vDash Env(t,z)$ & $z=xt'auat$ & $Tally_b(t')$ along with

$M \vDash Env(t_0,x)$ & $Lastf(z,t,aua,t)$ and $M \vDash J(x)$.

Let $M \vDash w \, \varepsilon \, x$.

Assume that $M \vDash Free(x,t_1,awa,t_2)$.

$\Rightarrow M \vDash Fr(x,t_1,awa,t_2)$ & $\neg Firstf(x,t_1,awa,t_2)$ &

    & $\forall v,t_3,t_4(Fr(x,t_3,ava,t_4)$ & $v<_x w$ & $\neg \exists u'(v<_x u'$ & $u'<_x w) \rightarrow t_1=t_3 b)$,

$\Rightarrow$ from $M \vDash Env(t_0,x)$ & $z=xt'auat$ & $Tally_b(t')$ & $Tally_b(t)$ & $Fr(x,t_1,awa,t_2)$,

$\Rightarrow$ by (5.6), $M \vDash \exists t_5 Fr(z,t_1,awa,t_5)$.

Suppose, for a reductio, that $M \vDash Firstf(z,t_1,awa,t_5)$.

$\Rightarrow$ from $M \vDash Env(t_0,x)$ & $xBz$, by (5.4), $M \vDash \exists t_6 Firstf(x,t_1,awa,t_6)$,

$\Rightarrow$ from $M \vDash \neg Firstf(x,t_1,awa,t_2)$, by (5.20), $M \vDash t_1<t_1$,

contradicting $M \vDash t_1 \in I \subseteq I_0$.

Therefore, $M \vDash \neg Firstf(z,t_1,awa,t_5)$.

Now, assume that $M \vDash Fr(z,t_3,ava,t_4)$ & $v<_z w$ & $\neg \exists u'(v<_z u'$ & $u'<_z w)$.

$\Rightarrow$ from $M \vDash Env(t_0,x)$ & $z=xt'auat$ & $Tally_b(t')$, by (5.13),

$M \vDash t'<t$ & $\exists t^*,w^*(Tally_b(t^*)$ & $aBw^*$ & $aEw^*$ & $z=t^*w^*tauat)$.

In fact, from the proof of (5.13), $M \vDash t_0 t'=t$.

$\Rightarrow M \vDash xt'auat=z=t^*w^*tauat=t^*w^*t_0t'auat$,

$\Rightarrow$ by (3.6), $M \vDash x=t^*w^*t_0$.

Let $u^+=tauat$.

$\Rightarrow$ from $M \vDash Lastf(z,t,aua,t)$, $M \vDash Pref(aua,t)$,

$\Rightarrow M \vDash Firstf(u^+,t,aua,t)$ & $Lastf(u^+,t,aua,t)$,

$\Rightarrow$ by (5.22), $M \vDash Env(t,u^+)$ & $\forall w'(w' \, \varepsilon \, u^+ \leftrightarrow w'=u)$.

Suppose, for a reductio, that $M \vDash u \, \varepsilon \, x$.

$\Rightarrow M \vDash \exists t_7, t_8 \, Fr(x,t_7,aua,t_8)$,

$\Rightarrow M \vDash Tally_b(t_7)$,

$\Rightarrow$ from $M \vDash Env(t_0,x)$, $M \vDash MaxT_b(t_0,x)$,

$\Rightarrow M \vDash t_7 \leq t_0 < t$,

$\Rightarrow$ from $M \vDash Env(t,z)$ & $Fr(z,t,aua,t)$, $M \vDash t=t_7<t$, contradicting $M \vDash t \in I \subseteq I_0$.

Therefore, $M \vDash \neg(u \, \varepsilon \, x)$ & $\neg \exists w'(w' \, \varepsilon \, x$ & $w' \, \varepsilon \, u^+)$.

So we have

$M \vDash Env(t_0,x)$ & $x=t^*w^*t_0$ & $aBw^*$ & $aEw^*$ & $z=t^*w^*u^+$ & $Env(t,u^+)$ & $t_0<t$ &

& $Firstf(u^+,t,aua,t)$ & $\neg\exists w'(w' \, \varepsilon \, x$ & $w' \, \varepsilon \, u^+)$,

$\Rightarrow$ by (5.46), $M \vDash \forall w'(w' \, \varepsilon \, z \leftrightarrow w' \, \varepsilon \, x \, \vee \, w' \, \varepsilon \, u^+ \leftrightarrow w' \, \varepsilon \, x \, \vee \, w'=u)$.

From $M \vDash Lastf(z,t,aua,t)$, by (9.3), $M \vDash w \leq_z u$.

$\Rightarrow$ from $M \vDash v<_z w$, by (9.16), $M \vDash v<_z u$,

$\Rightarrow$ by (9.4), $M \vDash v \neq u$,

$\Rightarrow$ since $M \vDash v \, \varepsilon \, z$, $M \vDash v \, \varepsilon \, x$.

From $M \vDash Fr(z,t_3,ava,t_4)$, by (5.39),

$\Rightarrow M \vDash \exists t_7 \, Fr(x,t_3,ava,t_7) \lor Fr(u^+,t,ava,t) \lor \exists t_7 \, Fr(u^+,t_3,ava,t_7)$.

But from $M \vDash Fr(u^+,t,ava,t) \lor Fr(u^+,t_3,ava,t_7)$, it follows that $M \vDash v=u$,

a contradiction.

Therefore $M \vDash \exists t_7 \, Fr(x,t_3,ava,t_7)$.

Then from

$\quad M \vDash Env(t_0,x) \, \& \, z=xt'auat \, \& \, Tally_b(t') \, \& \, Tally_b(t) \, \& \, v \, \varepsilon \, x \, \& \, w \, \varepsilon \, x \, \& \, v<_z w$,

by (9.11), $M \vDash v<_x w$.

Suppose, for a reductio, that $M \vDash \exists u'(v<_x u' \, \& \, u'<_x w)$.

$\Rightarrow M \vDash u' \, \varepsilon \, x$,

$\Rightarrow$ by (9.11), $M \vDash v<_z u' \, \& \, u'<_z w$, contradicting hypothesis.

Therefore, $M \vDash \neg \exists u'(v<_x u' \, \& \, u'<_x w)$.

Hence, along with $M \vDash Fr(x,t_3,ava,t_7) \, \& \, v<_x w$, we have, from the hypothesis

$M \vDash Free(x,t_1,awa,t_2)$, that $M \vDash t_3 = t_1 b$.

Therefore, we have established that

$\quad M \vDash Fr(z,t_1,awa,t_5) \, \& \, \neg Firstf(z,t_1,awa,t_5) \, \&$

$\qquad \& \, \forall v,t_3,t_4(Fr(z,t_3,ava,t_4) \, \& \, v<_z w \, \& \, \neg \exists u'(v<_z u' \, \& \, u'<_z w) \rightarrow t_3=t_1 b)$.

This shows that $M \vDash Free(z,t_1,awa,t_5)$.

Conversely, assume that $M \vDash Free(z,t_1,awa,t_3)$.

$\Rightarrow M \vDash Fr(z,t_1,awa,t_3) \, \& \, \neg Firstf(z,t_1,awa,t_3) \, \&$

$\qquad \& \, \forall v,t_4,t_5(Fr(z,t_4,ava,t_5) \, \& \, v<_z w \, \& \, \neg \exists u'(v<_z u' \, \& \, u'<_z w) \rightarrow t_1=t_4 b)$.

From hypothesis  $M \vDash w \, \varepsilon \, x$,  $M \vDash \exists t_6, t_7 \, Fr(x, t_6, awa, t_7)$.

$\implies$ as in the first part of the argument,  $M \vDash \exists t_8 \, Fr(z, t_6, awa, t_8)$,

$\implies$ from  $M \vDash Env(t, z) \, \& \, Fr(z, t_1, awa, t_3)$,  $M \vDash t_6 = t_1$,

$\implies M \vDash Fr(x, t_1, awa, t_7)$.

Suppose, for a reductio, that  $M \vDash Firstf(x, t_1, awa, t_7)$.

$\implies$ by the proof of (5.6),  $M \vDash Firstf(z, t_1, awa, t_7)$,

$\implies$ from  $M \vDash \neg Firstf(z, t_1, awa, t_3)$  and (5.20),  $M \vDash t_1 < t_1$,  contradicting  $M \vDash t_1 \in I \subseteq I_0$.

Therefore,  $M \vDash \neg Firstf(x, t_1, awa, t_7)$.

Now, assume that  $M \vDash Fr(x, t_9, ava, t_{10}) \, \& \, v <_x w \, \& \, \neg \exists u'(v <_x u' \, \& \, u' <_x w)$.

$\implies$ by (5.6),  $M \vDash \exists t_{11} Fr(z, t_9, ava, t_{11})$,

$\implies$ by (9.11),  $M \vDash v <_z w$.

Suppose, for a reductio, that  $M \vDash \exists u'(v <_z u' \, \& \, u' <_z w)$.

Then we have, just as above, that  $M \vDash u' \neq u \, \& \, u' \, \varepsilon \, z$.

$\implies M \vDash u' \, \varepsilon \, x$,

$\implies$ from  $M \vDash v <_z u' \, \& \, u' <_z w$,  by (9.11),  $M \vDash v <_x u' \, \& \, u' <_x w$,  contradicting hypothesis.

Therefore,  $M \vDash \neg \exists u'(v <_z u' \, \& \, u' <_z w)$.

Along with  $M \vDash Fr(z, t_9, ava, t_{11}) \, \& \, v <_z w$,  we have, from hypothesis  $M \vDash Free(z, t_1, awa, t_3)$, that  $M \vDash t_1 = t_9 b$.

But then we have shown that

   $M \vDash Fr(x, t_1, awa, t_7) \, \& \, \neg Firstf(x, t_1, awa, t_7) \, \&$

& $\forall v, t_9, t_{10}(Fr(x, t_9, ava, t_{10})$ & $v <_x w$ & $\neg \exists u'(v <_x u'$ & $u' <_x w) \rightarrow t_1 = t_9 b)$,

which shows that  $M \vDash Free(x, t_1, awa, t_7)$.

Therefore, $M \vDash \exists t_2 Free(x, t_1, awa, t_2) \leftrightarrow \exists t_3 Free(z, t_1, awa, t_3)$.

A completely analogous argument shows that

$M \vDash \exists t_2 Bound(x, t_1, awa, t_2) \leftrightarrow \exists t_3 Bound(z, t_1, awa, t_3)$.

This completes the proof of (9.23).

(9.24) For any string concept $I \subseteq I_0$ there is a string concept $J \subseteq I$ such that

$QT^+ \vdash \forall z \in J \; \forall x,t,t',t_0,u(Env(t,z) \;\&\; z=xt'auat \;\&\; Tally_b(t') \;\&\; Env(t_0,x) \;\&$

   $\&\; Lastf(z,t,aua,t) \to$

      $\to \forall w,t_1(w \;\varepsilon\; x \to (\exists t_2 Free^+(x,t_1,awa,t_2) \leftrightarrow \exists t_3 Free^+(z,t_1,awa,t_3))))$.

Let $J \equiv I_{5.11} \;\&\; I_{5.41} \;\&\; I_{5.42} \;\&\; I_{9.6} \;\&\; I_{9.23}$.

Assume $M \vDash Env(t,z) \;\&\; z=xt'auat \;\&\; Tally_b(t')$ along with

$M \vDash Env(t_0,x) \;\&\; Lastf(z,t,aua,t)$ and $M \vDash J(x)$.

Let $M \vDash w \;\varepsilon\; x$.

Assume $M \vDash Free^+(x,t_1,awa,t_2)$.

$\Longrightarrow M \vDash Fr(x,t_1,awa,t_2) \;\&\; \forall v,t_4,t_5(Fr(x,t_4,ava,t_5) \;\&\; v \leq_x w \to$

                                    $\to Firstf(x,t_4,ava,t_5) \;v\; Free(x,t_4,ava,t_5))$,

$\Longrightarrow$ by (5.6), $M \vDash \exists t_3 Fr(z,t_1,awa,t_3)$.

Assume $M \vDash Fr(z,t_4,ava,t_5) \;\&\; v \leq_z w$.

We have that $M \vDash w \;\varepsilon\; x$. We claim that $M \vDash v \;\varepsilon\; x$.

First, note that, from $M \vDash Env(t_0,x)$, by (5.11),

    $M \vDash \exists t_1,w_1(x=t_1 w_1 t_0 \;\&\; aBw_1 \;\&\; aEw_1 \;\&\; Tally_b(t_1) \;\&\; Tally_b(t_0))$,

$\Longrightarrow M \vDash z=xt'auat=(t_1 w_1 t_0)t'auat$,

$\Longrightarrow$ from $M \vDash Lastf(z,t,aua,t)$, $M \vDash z=tauat \;v\; \exists w' \; z=w'atauat$.

If $M \vDash z=tauat$, then $M \vDash t(auat)=z=(t_1 w_1 t_0)t'auat$. Then, by (3.6),

$M \vDash t=t_1 w_1 t_0 t'$, a contradiction because $M \vDash Tally_b(t) \;\&\; a \subseteq_p w_1$.

Therefore, $M \vDash \exists w' \; z=w'atauat$.

$\Rightarrow$ $M \vDash$ w'at(auat)=z=$t_1 w_1 t_0$t'(auat),

$\Rightarrow$ by (3.6), $M \vDash$ w'at=$t_1 w_1 t_0$t',

$\Rightarrow$ from $M \vDash$ Tally$_b$($t_0$) & Tally$_b$(t'), by (4.5), $M \vDash$ Tally$_b$($t_0$t'),

$\Rightarrow$ by (4.24$^b$), $M \vDash$ t=$t_0$t' since $M \vDash$ aE$w_1$.

Let $u^+$=tauat.

$\Rightarrow$ from $M \vDash$ Lastf(z,t,aua,t), $M \vDash$ Pref(aua,t),

$\Rightarrow$ $M \vDash$ Firstf($u^+$,t,aua,t) & Lastf($u^+$,t,aua,t),

$\Rightarrow$ by (5.22), $M \vDash$ Env(t,$u^+$) & $\forall$w'(w' $\varepsilon$ $u^+$ $\leftrightarrow$ w'=u),

$\Rightarrow$ by (5.41), $M \vDash$ $\neg\exists$w'(w' $\varepsilon$ x & w' $\varepsilon$ $u^+$),

$\Rightarrow$ since $M \vDash$ $t_0$<t, by (5.46), $M \vDash$ $\forall$w'(w' $\varepsilon$ z $\leftrightarrow$ w' $\varepsilon$ x $\vee$ w'=u).

If $M \vDash$ $\neg$(v $\varepsilon$ x), then $M \vDash$ v=u. But from $M \vDash$ Lastf(z,t,aua,t), by (9.3), from $M \vDash$ w $\varepsilon$ z, we have $M \vDash$ w$\leq_z$v. Then from $M \vDash$ v$\leq_z$w, by (9.6), we have $M \vDash$ v=w. So $M \vDash$ v $\varepsilon$ x after all, which proves the claim.

It follows that $M \vDash$ $\exists t_6,t_7$Fr(x,$t_6$,ava,$t_7$).

$\Rightarrow$ from $M \vDash$ v $\varepsilon$ x & w $\varepsilon$ x & v$\leq_z$w, by (9.11), $M \vDash$ v$\leq_x$w,

$\Rightarrow$ from hypothesis, $M \vDash$ Firstf(x,$t_6$,ava,$t_7$) $\vee$ Free(x,$t_6$,ava,$t_7$).

If $M \vDash$ Firstf(x,$t_6$,ava,$t_7$), then, by the proof of (5.6), $M \vDash$ $\exists t_8$Firstf(z,$t_6$,ava,$t_8$),

$\Rightarrow$ from $M \vDash$ Fr(z,$t_4$,ava,$t_5$), by (5.42), $M \vDash$ $t_4$=$t_6$ & $t_5$=$t_8$,

$\Rightarrow$ $M \vDash$ Firstf(z,$t_4$,ava,$t_5$).

If $M \vDash$ Free(x,$t_6$,ava,$t_7$), then, by (9.23), $M \vDash$ $\exists t_8$Free(z,$t_6$,ava,$t_8$).

$\Rightarrow$ $M \vDash$ Fr(z,$t_6$,ava,$t_8$),

$\Rightarrow$ from $M \vDash$ Fr(z,$t_4$,ava,$t_5$), by (5.42), $M \vDash$ $t_4$=$t_6$ & $t_5$=$t_8$,

$\Rightarrow$ M ⊨ Free(z,t_4,ava,t_5).

Hence we proved that

M ⊨ ∀v,t_4,t_5(Fr(z,t_4,ava,t_5) & v≤_z w → Firstf(z,t_4,ava,t_5) v Free(z,t_4,ava,t_5)).

Along with  M ⊨ Fr(z,t_1,awa,t_3)  this yields  M ⊨ Free$^+$(z,t_1,awa,t_3).

Conversely, assume  M ⊨ Free$^+$(z,t_1,awa,t_3).

$\Rightarrow$ M ⊨ Fr(z,t_1,awa,t_3) &

    & ∀v,t_4,t_5(Fr(z,t_4,ava,t_5) & v≤_z w → Firstf(z,t_4,ava,t_5) v Free(z,t_4,ava,t_5)).

From hypothesis  M ⊨ w ε x,  we have   M ⊨ ∃t_6,t_7 Fr(x,t_6,awa,t_7).

$\Rightarrow$ by (5.6),  M ⊨ ∃t_8 Free(z,t_6,awa,t_8),

$\Rightarrow$ from  M ⊨ Env(t,z),   M ⊨ t_1=t_6,

$\Rightarrow$ M ⊨ Fr(x,t_1,awa,t_8).

Assume   M ⊨ Fr(x,t_4,ava,t_5) & v≤_x w.

$\Rightarrow$ since  M ⊨ v ε x & w ε x,  by (9.7),  M ⊨ v≤_z w,

$\Rightarrow$ by (5.6),  M ⊨ ∃t_9 Fr(z,t_4,ava,t_9),

$\Rightarrow$ from hypothesis,  M ⊨ Firstf(z,t_4,ava,t_9) v Free(z,t_4,ava,t_9).

If  M ⊨ Firstf(z,t_4,ava,t_9), then by (5.4),   M ⊨ ∃t_{10} Firstf(x,t_4,ava,t_{10}),

$\Rightarrow$ from  M ⊨ Fr(x,t_4,ava,t_5), by (5.42),   M ⊨ t_{10}=t_5.

If  M ⊨ Free(z,t_4,ava,t_9), then from  M ⊨ v ε x,  by (9.23),

        M ⊨ ∃t_{10} Free(x,t_4,ava,t_{10}),

$\Rightarrow$ M ⊨ Fr(x,t_4,ava,t_{10}),

$\Rightarrow$ by (5.42),  M ⊨ t_{10}=t_5,

$\Rightarrow$ M ⊨ Free(x,t_4,ava,t_5).

So we have proved that

$M \vDash \forall v, t_4, t_5 (Fr(x, t_4, ava, t_5) \& v \leq_x w \rightarrow Firstf(x, t_4, ava, t_5) \vee Free(x, t_4, ava, t_5))$.

Along with $M \vDash Fr(x, t_1, awa, t_7)$ this yields $M \vDash Free^+(x, t_1, awa, t_7)$, as required.

This completes the proof of (9.24).

(9.25) For any string concept $I \subseteq I_0$ there is a string concept $J \subseteq I$ such that

$QT^+ \vdash \forall x,y \in J$ (Set(x) & Set(y) & x~y & Lex$^+$(x) →

$\qquad\qquad\qquad\qquad$ → (Lex$^+$(y) ↔ ∀u,v (u<$_x$v ↔ u<$_y$v))).

Let $J \equiv I_{8.2}$ & $I_{9.4}$ & $I_{9.7}$.

Assume $M \vDash$ Set(x) & Set(y)

where $M \vDash$ x~y & Lex$^+$(x) & Lex$^+$(y) and $M \vDash$ J(x) & J(y).

Assume $M \vDash$ u<$_x$v.

From $M \vDash$ x~y we have $M \vDash$ u ε y & v ε y.

Suppose, for a reductio, that $M \vDash$ ¬(u<$_y$v).

⟹ from $M \vDash$ u<$_x$v, by (9.4), $M \vDash$ u≠v.

⟹ by (9.7), $M \vDash$ v<$_y$u.

⟹ from $M \vDash$ Lex$^+$(x) & u<$_x$v, $M \vDash$ u≺v,

⟹ from $M \vDash$ Lex$^+$(y) & v<$_y$u, $M \vDash$ v≺u,

⟹ $M \vDash$ u≺v & v≺u, which contradicts (8.2).

Therefore $M \vDash$ u<$_x$v → u<$_y$v.

A symmetric argument establishes the converse.

Therefore $M \vDash$ u<$_x$v ↔ u<$_y$v.

On the other hand, assume $M \vDash$ ∀u,v (u<$_x$v ↔ u<$_y$v), and let $M \vDash$ w$_1$<$_y$w$_2$.

⟹ $M \vDash$ w$_1$<$_x$w$_2$,

⟹ from hypothesis $M \vDash$ Lex$^+$(x), $M \vDash$ w$_1$≺ w$_2$.

Hence $M \vDash$ ∀w$_1$,w$_2$ (w$_1$<$_y$w$_2$ → w$_1$≺ w$_2$), that is, $M \vDash$ Lex$^+$(y).

This completes the proof of (9.25).

(9.26) For any string concept $I \subseteq I_0$ there is a string concept $J \subseteq I$ such that

$QT^+ \vdash \forall x,z \in J \forall t,t',t_1,t_2,u,y[\text{Env}(t,x)$ & $x=t_1 auat_2 ax_1$ & $\text{Firstf}(x,t_1,aua,t_2)$ &

& $\neg(y \varepsilon x)$ & $\text{Max}^+T_b(t',aya)$ & $t'<t_2$ & $z=t'ayat_2ax_1$ &

& $\forall v(v \varepsilon x$ & $v \neq u \rightarrow y \prec v)$ & $\text{Lex}^+(x) \rightarrow \text{Lex}^+(z)]$.

Let $J \equiv I_{5.35}$ & $I_{9.1}$ & $I_{9.5}$.

Assume $M \vDash x=t_1 auat_2 ax_1$ & $z=t'ayat_2 ax_1$ where

$M \vDash \text{Env}(t,x)$ & $\text{Firstf}(x,t_1,aua,t_2)$ & $\neg(y \varepsilon x)$ & $\text{Max}^+T_b(t',aya)$ & $t'<t_2$

along with $M \vDash \forall v(v \varepsilon x$ & $v \neq u \rightarrow y \prec v)$ & $\text{Lex}^+(x)$ and $M \vDash J(x)$ & $J(z)$.

Assume also that $M \vDash w,w' \subseteq_p z$ & $w<_z w'$.

We want to show that $M \vDash w \prec w'$.

By (5.35), we have from the principal hypothesis that

$M \vDash \text{Env}(t,z)$ & $\text{Firstf}(z,t',aya,t_2)$ & $\forall w(w \varepsilon z \leftrightarrow (w \varepsilon x$ & $w \neq u) \vee w=y)$.

$\Rightarrow$ from $M \vDash w<_z w'$, $M \vDash w \varepsilon z$ & $w' \varepsilon z$,

$\Rightarrow M \vDash \exists t_3,t_4 \text{Fr}(z,t_3,awa,t_4)$ & $\exists t_5,t_6 \text{Fr}(z,t_5,aw'a,t_6)$.

<u>Case 1</u>. $M \vDash w'=y$.

$\Rightarrow$ by (9.1), $M \vDash \neg(w<_z w')$, contradicting the hypothesis.

<u>Case 2</u>. $M \vDash w' \varepsilon x$ & $w' \neq u$.

$\Rightarrow$ since $M \vDash \neg(y \varepsilon x)$, $M \vDash w' \neq y$,

$\Rightarrow$ by (5.15), $M \vDash \neg \text{Firstf}(z,t_5,aw'a,t_6)$,

$\Rightarrow M \vDash \exists w_2 \text{Intf}(z,w_2,t_5,aw'a,t_6) \vee \text{Lastf}(z,t_5,aw'a,t_6)$.

(2a) $M \vDash w=y$.

⟹ from the principal hypothesis, $M \vDash w=y\prec w'$, as required.

(2b) $M \vDash w \neq y$.

We claim that $M \vDash w<_x w'$.

From $M \vDash w \neq y$ we have $M \vDash w \, \varepsilon \, x \, \& \, w \neq u$.

⟹ since $M \vDash \neg(y \, \varepsilon \, x)$, $M \vDash w \neq y$,

⟹ by (5.15), $M \vDash \neg \text{Firstf}(z,t_3,awa,t_4)$,

⟹ from hypothesis $M \vDash w<_z w'$, by (9.5), $M \vDash \neg \text{Lastf}(z,t_3,awa,t_4)$,

⟹ $M \vDash \exists w_1 \, \text{Intf}(z,w_1,t_3,awa,t_4)$.

We now verify that the hypothesis of (5.31) holds, with the roles of x and z reversed.

From $M \vDash \text{Firstf}(x,t_1,aua,t_2) \, \& \, x=t_1auat_2ax_1$, we have that $M \vDash x \neq t_1auat_2$, for otherwise $M \vDash xBx$, contradicting $M \vDash x \in I \subseteq I_0$.

⟹ $M \vDash t_1 < t_2$.

On the other hand, from $M \vDash \text{Firstf}(z,t',aya,t_2) \, \& \, \text{Intf}(z,w_1,t_3,awa,t_4)$, we have, by (5.34), that $M \vDash t' < t_2 \leq t_3$.

Likewise, if $M \vDash \text{Intf}(z,w_2,t_5,aw'a,t_6)$, we have $M \vDash t' < t_2 \leq t_5$.

If $M \vDash \text{Lastf}(z,t_5,aw'a,t_6)$, we have from $M \vDash \text{Env}(t,z)$, that $M \vDash t' < t_2 \leq t = t_6 = t_5$.

Thus, $M \vDash t_1 < t_2 \, \& \, t' < t_3 \, \& \, t' < t_5$.

Applying (5.31) with the roles of x and z reversed, we obtain

$\qquad M \vDash \text{Fr}(x,t_3,awa,t_4) \, \& \, \text{Fr}(x,t_5,aw'a,t_6)$.

It then follows immediately from $M \vDash w<_z w'$ that $M \vDash w<_x w'$, as claimed.

But then $M \vDash w \prec w'$ follows from the hypothesis $M \vDash \text{Lex}^+(x)$.

We therefore have  $M \vDash \forall w, w'(w<_z w' \to w \prec w')$, that is  $M \vDash \text{Lex}^+(z)$.

This completes the proof of (9.26).

(9.27) For any string concept $I \subseteq I_0$ there is a string concept $J \subseteq I$ such that

$QT^+ \vdash \forall x \in J \forall t, t_1, t_2, w_1, v, x'(Env(t,x)$ & $Intf(x,w_1,t_1,ava,t_2)$ & $x'=w_1at_1avat_1$ &

& $Lex^+(x) \rightarrow Env(t_1,x')$ & $Lex^+(x'))$.

Let $J \equiv I_{5.55}$ & $I_{9.11}$.

Assume $M \vDash Env(t,x)$ & $Intf(x,w_1,t_1,ava,t_2)$ & $Lex^+(x)$ and let

$M \vDash x'=w_1at_1avat_1$ where $M \vDash J(x)$.

$\Rightarrow$ by (5.55), $M \vDash Env(t_1,x')$.

Assume $M \vDash u <_{x'} w$.

$\Rightarrow$ from $M \vDash Intf(x,w_1,t_1,ava,t_2)$, $M \vDash \exists w_2\, x = w_1at_1avat_2aw_2$ & $t_1 < t_2$,

$\Rightarrow$ $M \vDash \exists t_3\, (Tally_b(t_3)$ & $t_2 = t_1t_3)$,

$\Rightarrow$ $M \vDash x = w_1at_1avat_1t_3aw_2 = x't_3aw_2$,

$\Rightarrow$ by (5.55), $M \vDash \exists w'(x = x't_3w't$ & $aBw'$ & $aEw')$ & $Tally_b(t_3)$ & $Tally_b(t)$,

$\Rightarrow$ by (9.11), $M \vDash u <_x w$,

$\Rightarrow$ from $M \vDash Lex^+(x)$, $M \vDash u \not< w$.

Hence we have $M \vDash \forall u,w\, (u <_{x'} w \rightarrow u \not< w)$, and so $M \vDash Lex^+(x')$.

This completes the proof of (9.27).

(9.28) For any string concept $I \subseteq I_0$ there is a string concept $J \subseteq I$ such that

$QT^+ \vdash \forall x \in J \forall x',z,u,t,t',t'',t_3,t_4,v_0[Env(t,x)$ & $x=t'uz$ & $Env(t'',x')$ & $x'=t'ut''$ &

& $aBu$ & $aEu$ & $Env(t,z)$ & $Firstf(z, t_3,av_0a,t_4)$ & $t''<t_3$ & $Lex^+(x) \rightarrow$

$\rightarrow Lex^+(x')$ & $Lex^+(z)]$.

Let $J \equiv I_{9.11}$ & $I_{9.18}$.

Assume $M \vDash Env(t,x)$ & $x=t'uz$ & $Lex^+(x)$ & $J(x)$ where

$M \vDash x'=t'ut''$ & $Env(t'',x')$ & $aBu$ & $aEu$ & $Env(t,z)$ along with

$M \vDash Firstf(z, t_3,av_0a,t_4)$ & $t''<t_3$.

$\Rightarrow$ from hypothesis $M \vDash Lex^+(x)$, $M \vDash \forall v,w(v<_x w \rightarrow v \prec w)$.

We want to show that

$M \vDash \forall v,w(v<_{x'} w \rightarrow v \prec w)$ & $\forall v,w(v<_z w \rightarrow v \prec w)$.

It will thus suffice to prove that

$M \vDash \forall v,w(v<_{x'} w \lor v<_z w \rightarrow v<_x w)$.

Assume $M \vDash v<_{x'} w$.

$\Rightarrow$ from $M \vDash Env(t,z)$, by (5.11),

$M \vDash \exists t_5,z'(Tally_b(t_5)$ & $Tally_b(t)$ & $z=t_5z't$ & $aBz'$ & $aEz')$,

$\Rightarrow M \vDash \exists z_2\, az_2=z'$,

$\Rightarrow M \vDash z=t_5az_2t$,

$\Rightarrow$ from $M \vDash Firstf(z,t_3,av_0a,t_4)$, $M \vDash (t_3a)Bz$,

$\Rightarrow M \vDash \exists z_1\, t_3az_1=z=t_5az_2t$,

$\Rightarrow$ since $M \vDash Tally_b(t_3)$ & $Tally_b(t_5)$, by (4.23$^b$), $M \vDash t_3=t_5$,

$\Rightarrow$ $M \vDash z=t_3z't$,

$\Rightarrow$ from $M \vDash t''<t_3$, $M \vDash \exists t_0\, t''t_0=t_3$,

$\Rightarrow$ $M \vDash x=t'uz=t'ut_3z't=t'ut''t_0z't=x't_0z't$,

$\Rightarrow$ $M \vDash \text{Env}(t'',x')$ & $x't_0z't=x$ & $aBz'$ & $aEz'$ & $\text{Tally}_b(t_0)$ & $\text{Tally}_b(t)$,

$\Rightarrow$ by (9.11), $M \vDash v<_x w$, as required.

Assume $M \vDash v<_z w$.

$\Rightarrow$ by (9.18), $M \vDash v<_x w$, as required.

This completes the proof of (9.28).

(9.29) For any string concept $I \subseteq I_0$ there is a string concept $J \subseteq I$ such that

$$QT^+ \vdash \forall x \in J \forall x',z,t,t',t'',t_3,t_4,v_0,w[\text{Env}(t,x) \;\&\; x=t'wz \;\&\; \text{Env}(t'',x') \;\&\; x'=t'wt'' \;\&$$
$$\&\; aBw \;\&\; aEw \;\&\; \text{Env}(t,z) \;\&\; \text{Firstf}(z, t_3, av_0a, t_4) \;\&\; t''<t_3 \;\&\; \neg\exists u(u \,\varepsilon\, x' \;\&\; u \,\varepsilon\, z) \;\&$$
$$\&\; \forall u,v(u \,\varepsilon\, x' \;\&\; v \,\varepsilon\, z \to u \prec v) \;\&\; \text{Lex}^+(x') \;\&\; \text{Lex}^+(z) \to \text{Lex}^+(x)].$$

Let $J \equiv I_{9.28}$.

Assume

$M \vDash \text{Env}(t,x) \;\&\; \text{Env}(t'',x') \;\&\; \text{Env}(t,z) \;\&\; \text{Lex}^+(x') \;\&\; \text{Lex}^+(z) \;\&$

$$\&\; \forall u,v(u \,\varepsilon\, x' \;\&\; v \,\varepsilon\, z \to u \prec v),$$

where

$M \vDash x=t'wz \;\&\; x'=t'wt'' \;\&\; aBw \;\&\; aEw \;\&\; \text{Firstf}(z, t_3, av_0a, t_4) \;\&\; t''<t_3 \;\&$

$$\&\; \neg\exists u(u \,\varepsilon\, x' \;\&\; u \,\varepsilon\, z)$$

and $M \vDash J(x)$.

We want to show that $M \vDash \forall u,v(u<_x v \to u \prec v)$.

From the hypothesis $M \vDash \text{Lex}^+(x') \;\&\; \text{Lex}^+(z)$ we have that

$$M \vDash \forall u,v(u<_{x'} v \to u \prec v) \;\&\; \forall u,v(u<_z v \to u \prec v).$$

Hence it suffices to show that

$$M \vDash \forall u,v(u<_x v \to u<_{x'} v \;\vee\; u<_z v \;\vee\; u \prec v).$$

Assume $M \vDash u<_x v$.

$\Rightarrow$ by (5.11), as in the proof of (9.28),

$$M \vDash \exists z'(z=t_3 z't \;\&\; aBz' \;\&\; aEz' \;\&\; \text{Tally}_b(t_3) \;\&\; \text{Tally}_b(t)),$$

$\Rightarrow$ from hypothesis $M \vDash u<_x v$, $M \vDash u \,\varepsilon\, x \;\&\; v \,\varepsilon\, x$,

$\Rightarrow$ by (5.46), $M \vDash (u \, \varepsilon \, x' \, \& \, u \, \varepsilon \, z) \, \& \, (v \, \varepsilon \, x' \, \& \, v \, \varepsilon \, z)$.

We distinguish four scenarios:

(1) $M \vDash u \, \varepsilon \, x' \, \& \, v \, \varepsilon \, x'$.

$\Rightarrow$ from $M \vDash u <_x v$, by (9.11), $M \vDash u <_{x'} v$, as required.

(2) $M \vDash u \, \varepsilon \, x' \, \& \, v \, \varepsilon \, z$.

$\Rightarrow$ from hypothesis $M \vDash \forall u,v (u \, \varepsilon \, x' \, \& \, v \, \varepsilon \, z \rightarrow u \prec v)$, $M \vDash u \prec v$, as required.

(3) $M \vDash u \, \varepsilon \, z \, \& \, v \, \varepsilon \, x'$.

$\Rightarrow M \vDash \exists t_1, t_2 \, Fr(z, t_1, aua, t_2) \, \& \, \exists t_5, t_6 \, Fr(x', t_5, ava, t_6)$,

$\Rightarrow$ from $M \vDash Env(t'', x)$, $M \vDash Max^+T_b(t'', x')$,

$\Rightarrow M \vDash t_5 \leq t''$,

$\Rightarrow$ from $M \vDash Firstf(z, t_3, av_0 a, t_4)$, by (5.20), $M \vDash t_3 \leq t_1$,

$\Rightarrow$ from hypothesis $M \vDash t'' < t_3$, $M \vDash t_5 \leq t'' < t_3 \leq t_1$,

$\Rightarrow$ from $M \vDash Max^+T_b(t'', x') \, \& \, t'' < t_3 \leq t_1 \, \& \, t'w \subseteq_p x'$, $M \vDash Max^+T_b(t_1, t'w)$,

$\Rightarrow$ by (5.25), $M \vDash Fr(x, t_1, aua, t_2)$,

$\Rightarrow$ from $M \vDash x = t'wz = t'wt_3z't \, \& \, aBz' \, \& \, aEz' \, \& \, Tally_b(t_3) \, \& \, Tally_b(t) \, \&$

$\& \, Fr(x', t_5, ava, t_6)$,

by (5.6), $M \vDash \exists t_7 \, Fr(x', t_5, ava, t_7)$,

$\Rightarrow$ from hypothesis $M \vDash u <_x v$, by (9.14), $M \vDash t_1 < t_5$,

$\Rightarrow M \vDash t_1 < t_5 < t_1$, contradicting $M \vDash t_1 \in I \subseteq I_0$.

(4) $M \vDash u \, \varepsilon \, z \, \& \, v \, \varepsilon \, z$.

$\Rightarrow$ from $M \vDash u <_x v$, by (9.18), $M \vDash u <_z v$, as required.

This completes the proof of (9.29).

(9.30) For any string concept I⊆I₀ there is a string concept J⊆I such that

$QT^+ \vdash \forall x,y,t',t'' \in J \forall z(Max^+T_b(t',x)$ & $Max^+T_b(t'',y)$ & $z=t'axat''ayat''$ & $t'<t''$ &

& $x \prec y \rightarrow Set(z)$ & $\forall w(w \varepsilon z \leftrightarrow w=x \vee w=y))$ & $Lex^+(z))$.

Let $J \equiv I_{5.22}$ & $I_{8.2}$ & $I_{9.29}$.

Assume $M \vDash Max^+T_b(t',x)$ & $Max^+T_b(t'',y)$ & $z=t'axat''ayat''$

where $M \vDash J(x)$ & $J(y)$ & $J(t')$ & $J(t'')$ & $t'<t''$ & $x \prec y$.

Since we may assume that J is closed under *, we have that $M \vDash J(z)$.

Let $x'=t'axat'$ and $x''=t''ayat''$.

$\Rightarrow$ from $M \vDash Max^+T_b(t',x)$ & $Max^+T_b(t'',y)$,

$M \vDash Firstf(x',t',axa,t')$ & $Lastf(x',t',axa,t')$ &

& $Firstf(x'',t'',aya,t'')$ & $Lastf(x'',t'',aya,t'')$,

$\Rightarrow$ by (5.22), $M \vDash Env(t',x')$ & $Env(t'',x'')$ & $x \varepsilon x'$ & $y \varepsilon x''$ &

& $\forall w(w \varepsilon x' \leftrightarrow w=x)$ & $\forall w(w \varepsilon x'' \leftrightarrow w=y)$,

$\Rightarrow$ trivially, by (9.4), $M \vDash Lex^+(x')$ & $Lex^+(x'')$,

$\Rightarrow$ from $M \vDash x \prec y$, by (8.2), $M \vDash x \neq y$,

$\Rightarrow$ from $M \vDash Max^+T_b(t',x)$ & $Max^+T_b(t'',y)$, $M \vDash Pref(axa,t')$ & $Pref(aya,t'')$,

$\Rightarrow$ from hypothesis $M \vDash t'<t''$, by (5.58),

$M \vDash Env(t'',z)$ & $\forall w(w \varepsilon z \leftrightarrow w=x \vee w=y)$,

$\Rightarrow M \vDash Env(t'',z)$ & $z=t'ayax''$ & $Env(t',x')$ & $x'=t'axat'$ & $Env(t'',x'')$ &

& $Firstf(x'',t'',aya,t'')$ & $t'<t''$ & $\neg \exists w(w \varepsilon x'$ & $w \varepsilon x'')$ &

& $\forall u,v(u \varepsilon x'$ & $v \varepsilon x'' \rightarrow u \prec v)$ & $Lex^+(x')$ & $Lex^+(x'')$,

$\Rightarrow$ by (9.29), $M \vDash \text{Lex}^+(z)$, as required.

This completes the proof of (9.30).

## 10. Minimal, special and canonical set codes

Let $Occ(w_1,z,w_2,x,t',v,t'')$ abbreviate

$\quad w_1zw_2=x$ & $Fr(x,t',v,t'')$ &

& $[(Firstf(x,t',v,t'')$ & $(t'=w_1 \lor (w_1z)B(t'v) \lor w_1z=t'v)) \lor$

$\lor \exists w'((Intf(x,w',t'v,t'') \lor (Lastf(x,t'vt'')$ & $w'at'vt''=x))$ &

& $(w'at'=w_1 \lor \exists v_1(v_1Bv$ & $w'at'v_1=w_1z) \lor w'at'v=w_1z))]$,

meaning that the occurrence of z sandwiched between $w_1,w_2$ in x appears within the frame $t'vt''$.

We then let $MinSet(x)$ abbreviate

$\quad Set(x)$ & $\forall w_1,w_2 \subseteq_p x$ $(w_1aw_2=x \to \exists v \subseteq_p x \exists t',t'' \subseteq_p x\ Occ(w_1,a,w_2,x,t',v,t''))$,

meaning "every occurrence of the digit a in the set code x appears within some frame in x".

Let $Max^+(t,v,x)$ abbreviate

$Set(x)$ & $v \, \varepsilon \, x$ & $Tally_b(t)$ &

$\quad$ & $\forall u,t_1,t_2 \subseteq_p x\ (Fr(x, t_1,aua,t_2)$ & $u<_x v \to t_1 < t)$.

We then set $MMax^+T_b(t,v,x)$ to abbreviate

$\quad Max^+(t,v,x)$ & $\forall t'(Max^+(t',v,x)$ & $Max^+T_b(t',v) \to t \le t')$.

Let Special(x) abbreviate

   Set(x) & $\forall v, t_1, t_2$ (Fr(x, $t_1$,ava,$t_2$) → MMax$^+$T$_b$($t_1$,v,x)).

Let Set*(x) abbreviate

   MinSet(x) & Lex$^+$(x) & Special(x).

In that case we say that x is a <u>canonical set code</u>.

THE SUBTRACTION LEMMA. (10.1) For any string concept $I \subseteq I_0$ there is a string concept $I_{SUB} \subseteq I$ such that

$QT^+ \vdash \forall x \in I_{SUB} [MinSet(x) \rightarrow \forall y \exists z \in I_{SUB} (Set(z) \& \neg(y \varepsilon z) \&$

$\& \forall w(w \varepsilon z \leftrightarrow w \varepsilon x \& w \neq y) \&$

$\& \forall w, t_1, t_2 (Fr(x,t_1,awa,t_2) \& w \neq y \rightarrow \exists t_3 Fr(z,t_1,awa,t_3) \& (Lex^+(x) \rightarrow Lex^+(z))))]$

Let $I_{SUB} \equiv I_{5.22} \& I_{5.38} \& I_{5.44} \& I_{5.45} \& I_{5.46} \& I_{5.48} \& I_{5.50} \& I_{5.53} \& I_{9.14}$.

We may assume that $I_{SUB}(x)$ is closed under * and downward closed under $\subseteq_p$.

Assume $M \vDash MinSet(x)$ where $M \vDash I_{SUB}(x)$.

$\Rightarrow M \vDash Set(x)$, $M \vDash x=aa \lor \exists t Env(t,x)$.

(1) $M \vDash x=aa$.

Let $z=aa$.

$\Rightarrow M \vDash Set(z) \& I_{SUB}(z)$,

$\Rightarrow$ by (5.17), $M \vDash \forall y \neg(y \varepsilon z)$,

$\Rightarrow$ since $M \vDash \forall w \neg(w \varepsilon x)$, $M \vDash \forall y, w(w \varepsilon z \leftrightarrow w \varepsilon x \& w \neq y)$,

and trivially $M \vDash \forall y, w, t_1, t_2(Fr(x,t_1,awa,t_2) \& w \neq y \rightarrow \exists t_3 Fr(z,t_1,awa,t_3))$.

(2) $M \vDash \exists t Env(t,x)$.

Fix y.

Then $M \vDash y \varepsilon x \lor \neg(y \varepsilon x)$.

(2a) $M \vDash \neg(y \varepsilon x)$.

Let $z=x$.

$\Rightarrow M \vDash Set(z) \& \neg(y \varepsilon z)$.

Assume $M \vDash w \varepsilon z$.

Then $M \vDash w \, \varepsilon \, x$. If $M \vDash w=y$, then $M \vDash y \, \varepsilon \, x$, contradicting the hypothesis.

Therefore, $M \vDash \forall w(w \, \varepsilon \, z \rightarrow w \, \varepsilon \, x \, \& \, w \neq y)$.

The converse follows immediately from the hypothesis $M \vDash z=x$.

So $M \vDash \forall w(w \, \varepsilon \, z \leftrightarrow w \, \varepsilon \, x \, \& \, w \neq y)$.

Also, $M \vDash \forall w, t_1, t_2(Fr(x,t_1,awa,t_2) \, \& \, w \neq y \rightarrow \exists t_3 Fr(z,t_1,awa,t_3))$ follows trivially.

(2b) $M \vDash y \, \varepsilon \, x$.

$\Longrightarrow$ $M \vDash \exists t_1, t_2 Fr(x,t_1,aya,t_2)$.

There are three subcases:

(2bi) $M \vDash Firstf(x,t_1,aya,t_2)$.

$\Longrightarrow M \vDash Pref(aya,t_1) \, \& \, Tally_b(t_2) \, \&$

$\& \, ((t_1=t_2 \, \& \, t_1 ayat_2 = x) \vee (t_1 < t_2 \, \& \, (t_1 ayat_2 a)Bx))$.

(2bi1) $M \vDash t_1=t_2 \, \& \, t_1 ayat_2 = x$.

Let $z=aa$.

As in (1) we have that $M \vDash Set(z) \, \& \, J(z) \, \& \, \forall w \, \neg(w \, \varepsilon \, z)$.

$\Longrightarrow$ by (5.21), $M \vDash \forall w(w \, \varepsilon \, x \leftrightarrow w=y)$,

$\Longrightarrow$ $M \vDash \forall w(w \, \varepsilon \, z \leftrightarrow w \, \varepsilon \, x \, \& \, w \neq y)$.

Assume $M \vDash Fr(x,t',awa,t'') \, \& \, w \neq y$.

$\Longrightarrow$ $M \vDash w \, \varepsilon \, x$,

$\Longrightarrow$ $M \vDash w=y$, a contradiction.

Therefore, trivially, $M \vDash \forall w,t',t''(Fr(x,t',awa,t'') \, \& \, w \neq y \rightarrow \exists t_3 Fr(z,t',awa,t_3))$.

(2bi2) $M \vDash t_1 < t_2 \, \& \, (t_1 ayat_2 a)Bx$.

$\Longrightarrow$ $M \vDash \exists x_1 (t_1 ayat_2 a)x_1 = x$,

$\Rightarrow$ from $M \vDash \text{MinSet}(x)$, $M \vDash \exists v,t',t''\ \text{Occ}(t_1 ayat_2,a,x_1,x,t',v,t'')$.

    (2bi2a) $M \vDash \text{Firstf}(x,t',v,t'')$.

      (2bi2a1) $M \vDash t' = t_1 ayat_2$.

This is a contradiction because $M \vDash \text{Tally}_b(t')$.

      (2bi2a2) $M \vDash (t_1 ayat_2 a)B(t'v)$.

$\Rightarrow M \vDash \exists x_2 (t_1 ayat_2 a)x_2 = t'v$,

$\Rightarrow$ from $M \vDash \text{Firstf}(x,t_1,aya,t_2)\ \&\ \text{Firstf}(x,t',v,t'')$, by (5.15), $M \vDash aya = v$,

$\Rightarrow M \vDash t_1 vt_2 ax_2 = t'v$,

$\Rightarrow$ by (4.23$^b$), $M \vDash t_1 = t'$,

$\Rightarrow$ by (3.7), $M \vDash vt_2 ax_2 = v$,

$\Rightarrow M \vDash vBv$, contradicting $M \vDash v \in I \subseteq I_0$.

      (2bi2a3) $M \vDash t_1 ayat_2 a = t'v$.

Same as (2bi2a2) except that $x_2$ is omitted throughout.

    (2bi2b) $M \vDash \exists w'\text{Intf}(x,w',t',v,t'')$.

$\Rightarrow M \vDash \text{Pref}(v,t')\ \&\ t'<t''\ \&\ \text{Tally}_b(t'')\ \&\ \exists w''\ x = w'at'vt''aw''\ \&\ \text{Max}^+ T_b(t',w')$.

      (2bi2b1) $M \vDash w'at' = t_1 ayat_2$.

$\Rightarrow$ by (4.24$^b$), $M \vDash t' = t_2$,

$\Rightarrow$ by (3.6), $M \vDash w'a = t_1 aya$.

Let $z = t'vt''aw''$. Since $I_{\text{SUB}}$ is downward closed under $\subseteq_p$, we have $M \vDash I_{\text{SUB}}(z)$.

$\Rightarrow$ from $M \vDash \text{Env}(t,x)$, by (5.44), $M \vDash \text{Env}(t,z)$,

$\Rightarrow M \vDash \text{Set}(z)$.

Note that $M \vDash t_1 ayaz = x = t_1 ayat_2 ax_1$.

$\Rightarrow$ by (3.7),  $M \vDash z = t_2 a x_1$.

Suppose that  $M \vDash y \,\varepsilon\, z$.

$\Rightarrow$  $M \vDash \exists v',t_3,t_4 \ (Fr(z,t_3,v',t_4) \ \& \ v'=aya)$,

$\Rightarrow$ by hypothesis (2bi2),  $M \vDash t_1 a y a t_2 a x_1 = x = t_1 a y a z \ \& \ t_1 < t_2$.

From  $M \vDash Pref(v,t') \ \& \ t'<t'' \ \& \ Tally_b(t'') \ \& \ (t'vt''a)Bz$

we have that  $M \vDash Firstf(z,t',v,t'')$.

$\Rightarrow$ by (5.20),  $M \vDash t' \leq t_3$,

$\Rightarrow$  $M \vDash t_1 < t_2 = t' \leq t_3$,

$\Rightarrow$ $M \vDash t_1 < t_3 = t_1$,  contradicting  $M \vDash t_1 \in I \subseteq I_0$.

Therefore  $M \vDash \neg(y \,\varepsilon\, z)$.

We now proceed to show that  $M \vDash \forall w(w \,\varepsilon\, z \leftrightarrow w \,\varepsilon\, x \ \& \ w \neq y)$.

Assume  $M \vDash w \,\varepsilon\, z$.

$\Rightarrow$  $M \vDash \exists v',t_3,t_4 \ (Fr(z,t_3,v',t_4) \ \& \ v'=awa)$.

Then, just as in the proof of  $M \vDash Env(t,z)$  part (d) in (5.44),  we have that

$\qquad M \vDash Fr(x,t_3,awa,t_4)$.

Therefore  $M \vDash w \,\varepsilon\, x$, and we have just shown that  $M \vDash w \neq y$.

We now prove the converse,  $M \vDash \forall w(w \,\varepsilon\, x \ \& \ w \neq y \rightarrow w \,\varepsilon\, z)$.

Assume  $M \vDash w \,\varepsilon\, x \ \& \ w \neq y$.

$\Rightarrow$  $M \vDash \exists v',t_3,t_4 \ (Fr(x,t_3,v',t_4) \ \& \ v'=awa)$,

$\Rightarrow$ from $M \vDash Firstf(x,t_1,aya,t_2)$ by (5.15), $M \vDash \neg Firstf(x,t_3,v',t_4)$,

$\Rightarrow$ $M \vDash \exists w_1 \ Intf(x,w_1,t_3,v',t_4) \ \vee \ Lastf(x,t_3,v',t_4)$.

Suppose that  $M \vDash Lastf(x,t_3,v',t_4)$.

⇒ from  M ⊨ Env(t,z),  M ⊨ $t_3=t_4=t$ & $t'<t''\leq t=t_3$.

Thus

M ⊨ Firstf(x,$t_1$,aya,$t_2$) & Lastf(x,$t_3$,v',$t_4$) & Intf(x,w',t',v,t'') & x=w'at'vt''aw'' &

& x=w'at'vt''aw'' & z=t'vt''aw'' & $t'\leq t_3$,

⇒ by (5.43),   M ⊨ Lastf(z,$t_3$,v',$t_4$),

⇒  M ⊨ w ε z,   as required.

Suppose that   M ⊨ ∃$w_1$ Intf(x,$w_1$,$t_3$,v',$t_4$).

Here we apply (5.45) to derive  M ⊨ ∃$t_5$ Fr(z,$t_3$,v',$t_5$),  whence M ⊨ w ε z.

Thus we also have shown that   M ⊨ ∀w(w ε x & w≠y  →  w ε z).

Finally, assume that  M ⊨ Fr(x,$t_3$,awa,$t_4$) & w≠y. Then  M ⊨ w ε x & w≠y. But then from the proof just given we have that  M ⊨ ∃$t_5$ Fr(z,$t_3$,awa,$t_5$),  as required.

   (2bi2b2)  M ⊨ ∃$v_1$ ($v_1$Bv & w'at'$v_1$=$t_1$ayat$_2$a).

⇒ M ⊨ $v_1$=a ∨ aE$v_1$.

If  M ⊨ $v_1$=a, then M ⊨ w'at'a=$t_1$ayat$_2$a,  whence  M ⊨ w'at'=$t_1$ayat$_2$,  and we proceed exactly as in (2bi2b1).

So we may assume that  M ⊨ aE$v_1$.

⇒ M ⊨ ∃$v_2$ $v_2$=$v_1$,

⇒ M ⊨ w'at'($v_2$a)=$t_1$ayat$_2$a,

⇒ M ⊨ w'at'$v_2$=$t_1$ayat$_2$,

⇒ M ⊨ Tally$_b$(t'$v_2$) ∨ ¬ Tally$_b$(t'$v_2$).

   (2bi2b2a) $M \vDash \text{Tally}_b(t'v_2)$.

$\Rightarrow$ from $M \vDash w'at'v_2 = t_1ayat_2$, by (4.24$^b$), $M \vDash t'v_2 = t_2$,

$\Rightarrow$ $M \vDash t' < t_2$.

But from $M \vDash \text{Firstf}(x, t_1, aya, t_2)$ & $\text{Intf}(x, w', t', v, t'')$ by (5.34), $M \vDash t_2 \leq t'$,

$\Rightarrow$ $M \vDash t' < t_2 \leq t'$, contradicting $M \vDash t' \in I \subseteq I_0$.

   (2bi2b2b) $M \vDash \neg \text{Tally}_b(t'v_2)$.

$\Rightarrow$ $M \vDash a \subseteq_p t'v_2$,

$\Rightarrow$ since $M \vDash \text{Tally}_b(t_2)$, $M \vDash \neg(v_2 = a \vee aEv_2)$,

$\Rightarrow$ $M \vDash aBv_2 \vee \exists v_3, v_4\; v_3av_4 = v_2$,

$\Rightarrow$ $M \vDash \exists v_5\; (av_5 = v_2\; \&\; w'at'(av_5) = t_1ayat_2) \vee w'at'(v_3av_4) = t_1ayat_2$,

$\Rightarrow$ by (4.16), $M \vDash t' \subseteq_p t'v_3 \subseteq_p u$.

But by (5.34), $M \vDash t_1 < t_2 \leq t' \subseteq_p u$, contradicting $M \vDash \text{Pref}(aya, t_1)$.

   (2bi2b3) $M \vDash w'at'v = t_1ayat_2a$.

$\Rightarrow$ from $M \vDash \text{Pref}(v, t')$, $M \vDash \exists v_0\; w'at'(av_0a) = t_1ayat_2a$,

$\Rightarrow$ $M \vDash w'at'av_0 = t_1ayat_2$,

$\Rightarrow$ by (4.16), $M \vDash t' \subseteq_p y$.

But by (5.34), $M \vDash t_1 < t_2 \leq t' \subseteq_p y$, again contradicting $M \vDash \text{Pref}(aya, t_1)$.

  (2bi2c) $M \vDash \text{Lastf}(x, t', v, t'')$.

$\Rightarrow$ $M \vDash \text{Pref}(v, t')\; \&\; t' = t''\; \&\; \exists w'\; (x = w'at'vt''\; \&\; \text{Max}^+T_b(t', w'))$.

   (2bi2c1) $M \vDash w'at' = t_1ayat_2$.

$\Rightarrow$ by (4.24$^b$), $M \vDash t' = t_2$,

$\Rightarrow$ by (3.6), $M \vDash w'a = t_1aya$.

Let z=t'vt''. Since we may assume that $I_{SUB}$ is downward closed under $\subseteq_p$, we have that $M \vDash I_{SUB}(z)$.

Then $M \vDash \text{Pref}(v,t')$ & $t'=t''$ & $z=t'vt''$.

$\Rightarrow M \vDash \text{Firstf}(z,t',v,t'')$ & $\text{Lastf}(z,t',v,t'')$,

$\Rightarrow$ by (5.22), $M \vDash \text{Set}(z)$ & $\exists w \forall u'(u' \varepsilon\ z \leftrightarrow u'=w)$,

$\Rightarrow$ from $M \vDash \text{Pref}(v,t')$, $M \vDash \exists v_0\ v=av_0a$,

$\Rightarrow$ from $M \vDash \text{Firstf}(z,t',v,t'')$, $M \vDash v_0\ \varepsilon\ z$,

$\Rightarrow M \vDash \forall u'(u' \varepsilon\ z \leftrightarrow u'=v_0)$.

We argue that $M \vDash \neg(y\ \varepsilon\ z)$.

Assume, for a reductio, that $M \vDash y\ \varepsilon\ z$.

$\Rightarrow M \vDash y=v_0$,

$\Rightarrow$ from $M \vDash \text{Fr}(x,t_1,aya,t_2)$ & $\text{Fr}(x,t',av_0a,t'')$ & $\text{Set}(x)$, $M \vDash t_1=t'$,

$\Rightarrow$ from $M \vDash \text{Env}(t,x)$ & $x=t_1ayat_2ax_1$, $M \vDash t_2 \leq t''=t'=t$,

$\Rightarrow M \vDash t_1 < t_2 \leq t'=t_1a$, contradicting $M \vDash t_1 \in I \subseteq I_0$.

Therefore, $M \vDash \neg(y\ \varepsilon\ z)$, as required.

We now show that $M \vDash \forall w(w\ \varepsilon\ z \leftrightarrow w\ \varepsilon\ x$ & $w \neq y)$.

Assume $M \vDash w\ \varepsilon\ z$.

$\Rightarrow M \vDash w=v_0$,

$\Rightarrow$ from hypothesis $M \vDash \text{Lastf}(x,t',v,t'')$ & $v=av_0a$, $M \vDash v_0\ \varepsilon\ x$,

$\Rightarrow M \vDash w\ \varepsilon\ x$.

That $M \vDash w \neq y$ follows from $M \vDash \neg(y\ \varepsilon\ z)$ and the hypothesis $M \vDash w\ \varepsilon\ z$.

Hence $M \vDash \forall w(w\ \varepsilon\ z \rightarrow w\ \varepsilon\ x$ & $w \neq y)$.

Conversely, assume that $M \vDash w \varepsilon x$ & $w \neq y$.

We have that $M \vDash w'a = t_1aya$ & $x = w'at'av_0at''$.

$\Rightarrow M \vDash Set(x)$ & $x = t_1ayat'av_0at''$ & $Firstf(x,t_1,aya,t')$ & $Lastf(x,t',av_0a,t'')$,

$\Rightarrow$ by (5.58), $M \vDash \forall w(w \varepsilon x \leftrightarrow w=y \vee w=v_0)$,

$\Rightarrow$ from hypothesis $M \vDash w \varepsilon x$ & $w \neq y$, $M \vDash w=v_0$,

$\Rightarrow$ from $M \vDash v_0 \varepsilon z$, $M \vDash w \varepsilon z$.

Therefore, we also have $M \vDash \forall w(w \varepsilon x$ & $w \neq y \rightarrow w \varepsilon z)$.

Finally, assume $M \vDash Fr(x,t_3,awa,t_4)$ & $w \neq y$.

$\Rightarrow M \vDash w \varepsilon x$ & $w \neq y$,

$\Rightarrow M \vDash w \varepsilon z$,

$\Rightarrow M \vDash w=v_0$,

$\Rightarrow$ from $M \vDash Fr(z,t',av_0a,t'')$, $M \vDash Fr(z,t',awa,t'')$,

$\Rightarrow$ from hypothesis $M \vDash Fr(z,t',av_0a,t'')$ and $M \vDash Env(t,x)$, $M \vDash t'=t_3$,

$\Rightarrow M \vDash \exists t_5 Fr(z,t_3,awa,t_5)$, as required.

So we have proved that

$M \vDash \forall w,t_3,t_4(Fr(x,t_3,awa,t_4)$ & $w \neq y \rightarrow \exists t_5 Fr(z,t_3,awa,t_5))$.

   (2bi2c2) $M \vDash \exists v_1(v_1 Bv$ & $w'at'v_1 = t_1ayat_2a)$.

$\Rightarrow M \vDash v_1 = a \vee aEv_1$.

If $M \vDash v_1 = a$, then $M \vDash w'at'a = t_1ayat_2a$, whence $M \vDash w'at' = t_1ayat_2$, and we proceed exactly as in (2bi2c1).

So we may assume that $M \vDash aEv_1$. In this case we derive a contradiction by following the pattern of the argument in (2bi2b2) with the sole exception that, instead of appealing to (5.34), we obtain $M \vDash t_2 \leq t'$ from the hypothesis

$M \vDash \text{Set}(x) \ \& \ \text{Lastf}(x,t',v,t'')$.

   (2bi2c3)  $M \vDash w'at'v = t_1ayat_2a$.

Again, we derive a contradiction as in (2bi2b3) by modifying the argument as just described.

This completes the subcase (2bi).

   (2bii)  $M \vDash \exists w_1 \text{Intf}(x,w_1,t_1,aya,t_2)$.

$\Rightarrow M \vDash \text{Pref}(aya,t_1) \ \& \ \text{Tally}_b(t_2) \ \& \ t_1 < t_2 \ \&$

$$\& \ \exists w_2 \ x = w_1at_1ayat_2aw_2 \ \& \ \text{Max}^+(t_1,w_1).$$

Let $z = w_1at_2aw_2$.

Since we may assume that J is closed under $\subseteq_p$ and *, we have that $M \vDash I_{\text{SUB}}(z)$.

$\Rightarrow$ by (5.48), $M \vDash \text{Env}(t,z)$,

$\Rightarrow M \vDash \text{Set}(z)$.

Suppose, for a reductio, that $M \vDash y \ \varepsilon \ z$.

$\Rightarrow M \vDash \exists t',t'' \text{Fr}(z,t',aya,t'')$,

$\Rightarrow$ by (5.49), $M \vDash t' < t_1 \ \text{v} \ t_1 < t'$,

$\Rightarrow$ by (5.47), $M \vDash \text{Fr}(x,t',aya,t'') \ \text{v} \ \text{Fr}(x,t'aya,t_1)$,

$\Rightarrow$ from $M \vDash \text{Fr}(x,t_1,aya,t_2)$ by (d) of $M \vDash \text{Env}(t,x)$, $M \vDash t_1 = t' < t_1$,

a contradiction because $M \vDash t_1 \in I \subseteq I_0$.

Therefore, $M \vDash \neg(y \ \varepsilon \ z)$.

We now proceed to show that $M \vDash \forall w(w \ \varepsilon \ z \leftrightarrow w \ \varepsilon \ x \ \& \ w \neq y)$.

Assume $M \vDash w \ \varepsilon \ z$.

$\Rightarrow$ M ⊨ ∃t',t''Fr(z,t',awa,t''),

$\Rightarrow$ by (5.47), M ⊨ Fr(x,t',awa,t'') v Fr(x,t'awa,$t_1$),

$\Rightarrow$ either way, M ⊨ w ε x.

That M ⊨ w≠y then follows from the preceding argument and the hypothesis M ⊨ w ε z.

Assume now that M ⊨ w ε x & w≠y.

$\Rightarrow$ M ⊨ ∃t',t''Fr(x,t',awa,t'').

We distinguish three cases:

    (2bii1) M ⊨ Firstf(x,t',awa,t'').

$\Rightarrow$ by (5.36), M ⊨ Firstf(z,t',ava,$t_2$) v Firstf(z,t',ava,t''),

$\Rightarrow$ either way, M ⊨ w ε z.

    (2bii2) M ⊨ ∃w'Intf(x,w',t',awa,t'').

$\Rightarrow$ M ⊨ ∃w''(x=w'at'awat''aw'' & Max$^+$T$_b$(t,w')),

$\Rightarrow$ from hypothesis M ⊨ Intf(x,$w_1$,$t_1$,aya,$t_2$), M ⊨ x=$w_1$a$t_1$aya$t_2$a$w_2$.

Suppose, for a reductio, that M ⊨ $t_1$=t'.

$\Rightarrow$ from M ⊨ Firstf(x,t',awa,t'') & Fr(x,$t_1$,aya,$t_2$) & Env(t,x), M ⊨ w=y.

But this contradicts the hypothesis M ⊨ w≠y.

Therefore, M ⊨ $t_1$≠t'.

$\Rightarrow$ by (5.50), M ⊨ Fr(z,t',awa,t'') v Fr(z,t',awa,$t_2$),

$\Rightarrow$ either way, M ⊨ w ε z.

    (2bii3) M ⊨ Lastf(x,t',awa,t'').

$\Rightarrow$ M ⊨ t'=t'',

$\Rightarrow M \vDash Lastf(x,t',awa,t')$,

$\Rightarrow$ by (5.37), $M \vDash Lastf(z,t',awa,t')$,

$\Rightarrow M \vDash w \, \varepsilon \, z$.

This completes the argument that $M \vDash \forall w(w \, \varepsilon \, x \, \& \, w \neq y \rightarrow w \, \varepsilon \, z)$.

Finally, assume that $M \vDash Fr(x, t_3, awa, t_4) \, \& \, w \neq y$.

Then the arguments of (2bii1)-(2bii3) show that $M \vDash \exists t_5 \, Fr(z, t_3, awa, t_5)$.

So we also have that

$$M \vDash \forall w, t_3, t_4 (Fr(x,t_3,awa,t_4) \, \& \, w \neq y \rightarrow \exists t_5 Fr(z,t_3,awa,t_5)).$$

This completes the subcase (2bii).

(2biii) $M \vDash Lastf(x,t_1,aya,t_2)$.

$\Rightarrow M \vDash Pref(aya,t_1) \, \& \, Tally_b(t_2) \, \& \, t_1 = t_2 \, \&$

$\qquad \qquad \& \, (x = t_1 ayat_2 \, \vee \, \exists w_1 (x = w_1 at_1 ayat_2 \, \& \, Max^+T_b(t_1,w_1))$.

(2biii1) $M \vDash x = t_1 ayat_2$.

Let $z = aa$. Then the argument is the same as in (2bi1).

(2biii2) $M \vDash \exists w_1 (x = w_1 at_1 ayat_2 \, \& \, Max^+T_b(t_1,w_1))$.

$\Rightarrow$ from $M \vDash MinSet(x)$, $M \vDash \exists v, t', t'' \, Occ(w_1, a, t_1 ayat_2, x, t', v, t'')$.

We distinguish three cases:

(2biii2a) $M \vDash Firstf(x,t',v,t'')$.

$\Rightarrow M \vDash Pref(v,t') \, \& \, Tally_b(t'') \, \& \, ((t' = t'' \, \& \, x = t'vt'') \, \vee \, (t' < t'' \, \& \, (t'vt''a)Bx))$.

(2biii2a1) $M \vDash t' = t'' \, \& \, x = t'vt''$.

$\Rightarrow M \vDash w_1 at_1 ayat_2 = x = t'vt''$,

$\Rightarrow$ from $M \vDash Pref(v,t')$, $M \vDash \exists v_0 (v = av_0 a \, \& \, Max^+T_b(t',av_0 a) \, \& \, Tally_b(t''))$,

$\Rightarrow$ from $M \vDash Env(t,x)$, $M \vDash t_1=t_2=t$ & $MaxT_b(t,x)$.

But this contradicts (4.21).

  (2biii2a2) $M \vDash t'<t''$ & $(t'vt''a)Bx$.

$\Rightarrow M \vDash \exists x_1\ t'vt''ax_1=x=w_1at_1ayat_2$,

$\Rightarrow$ from $M \vDash Env(t,x)$ & $t_1=t_2=t$, $M \vDash t'<t''\leq t_1$.

From $M \vDash Occ(w_1,a,t_1ayat_2,x,t',v,t'')$, we distinguish three subcases:

  (2biii2a2a) $M \vDash t'=w_1$.

$\Rightarrow M \vDash t'vt''ax_1=x=t'at_1ayat_2$,

$\Rightarrow$ from $M \vDash Pref(v,t')$, $M \vDash \exists v_0\ v=av_0a$,

$\Rightarrow M \vDash t'av_0at''ax_1=x=t'at_1ayat_2$,

$\Rightarrow$ by (3.7), $M \vDash v_0at''ax_1=t_1ayat_2$,

$\Rightarrow$ by (4.14$^b$), $M \vDash v_0=t_1$ v $t_1Bv_0$,

$\Rightarrow M \vDash t_1\subseteq_p v_0\subseteq_p v$,

$\Rightarrow M \vDash t'<t_1\subseteq_p v$, contradicting $M \vDash Pref(v,t')$.

  (2biii2a2b) $M \vDash (w_1a)B(t'v)$.

$\Rightarrow M \vDash \exists w_2\ w_1aw_2=t'v$,

$\Rightarrow M \vDash (w_1aw_2)t''ax_1=w_1at_1ayat_2$,

$\Rightarrow$ by (3.7), $M \vDash w_2t''ax_1=t_1ayat_2$,

$\Rightarrow$ from $M \vDash Pref(v,t')$, $M \vDash \exists v_0\ v=av_0a$,

$\Rightarrow$ from $M \vDash w_1aw_2=t'v=av_0a$, $M \vDash w_2=a$ v $aEw_2$.

Now, we cannot have $M \vDash w_2=a$ because $M \vDash Tally_b(t_1)$.

$\Rightarrow M \vDash aEw_2$,

$\Rightarrow M \vDash \exists w_3\ w_2 = w_3 a$,

$\Rightarrow M \vDash (w_3 a)t''ax_1 = t_1 ayat_2$,

$\Rightarrow$ by (4.14$^b$), $M \vDash w_3 = t_1\ \vee\ t_1 B w_3$,

$\Rightarrow M \vDash t_1 \subseteq_p w_3 \subseteq_p w_2$,

$\Rightarrow M \vDash t' < t'' \leq t_1 \subseteq_p w_2$,

$\Rightarrow$ from $M \vDash w_1 a w_2 = t'v$, by (4.14$^b$), $M \vDash w_1 = t'\ \vee\ t'Bw_1$,

$\Rightarrow M \vDash t'aw_2 = t'v\ \vee\ \exists w_4 (t'w_4 = w_1\ \&\ (t'w_4)aw_2 = t'v)$,

$\Rightarrow$ by (3.7), $M \vDash aw_2 = v\ \vee\ w_4 aw_2 = v$,

$\Rightarrow M \vDash w_2 \subseteq_p v$,

$\Rightarrow M \vDash t' \subseteq_p w_2 \subseteq_p v$, contradicting $M \vDash \text{Pref}(v, t')$.

    (2biii2a2c) $M \vDash w_1 a = t'v$.

$\Rightarrow M \vDash x = (t'v)t_1 ayat_2$.

Let $z = t'vt'$. We then apply an argument analogous to that of (2bi2c1).

    (2biii2b) $M \vDash \exists w'\text{Intf}(x, w', t', v, t'')$.

$\Rightarrow M \vDash \text{Pref}(v, t')\ \&\ \text{Tally}_b(t'')\ \&\ t' < t''\ \&\ \exists w''(x = w'at'vt''aw''\ \&\ \text{Max}^+T_b(t', w'))$,

$\Rightarrow M \vDash w'at'vt''aw'' = x = w_1 at_1 ayat_2$.

Again, from $M \vDash \text{Env}(t, x)\ \&\ t_1 = t_2 = t$, we have $M \vDash t' < t'' \leq t_1$.

From $M \vDash \text{Occ}(w_1, a, t_1 ayat_2, x, t', v, t'')$, we distinguish three subcases:

    (2biii2b1) $M \vDash w'at' = w_1$.

$\Rightarrow$ from $M \vDash \text{Pref}(v, t')$, $M \vDash \exists v_0\ v = av_0 a$,

$\Rightarrow M \vDash w'at'av_0 at''aw'' = w_1 at_1 ayat_2$,

$\Rightarrow$ by (3.7), $M \vDash v_0 at''aw'' = t_1 ayat_2$,

$\Rightarrow$ by (4.14$^b$), M $\vDash$ v$_0$=t$_1$ v t$_1$Bv$_0$,

$\Rightarrow$ M $\vDash$ t$_1 \subseteq_p$v$_0 \subseteq_p$v,

$\Rightarrow$ M $\vDash$ t'<t$_1 \subseteq_p$v,  contradicting  M $\vDash$ Pref(v,t').

(2biii2b2)  M $\vDash \exists$v$_1$ (v$_1$Bv & w'at'v$_1$=w$_1$a).

$\Rightarrow$  M $\vDash \exists$v$_2$ v$_1$v$_2$=v,

$\Rightarrow$ M $\vDash$ w'at'(v$_1$v$_2$)t''aw''=x=w$_1$at$_1$ayat$_2$,

$\Rightarrow$ by (3.7), M $\vDash$v$_2$t''aw''=t$_1$ayat$_2$,

$\Rightarrow$ from M $\vDash$ Pref(v,t'), M $\vDash \exists$v$_0$ av$_0$a=v= v$_1$v$_2$,

$\Rightarrow$  M $\vDash$ v$_2$=a v aEv$_2$.

Now, we cannot have  M $\vDash$ v$_2$=a  because  M $\vDash$ Tally$_b$(t$_1$).

$\Rightarrow$ M $\vDash$ aEv$_2$,

$\Rightarrow$ M $\vDash \exists$v$_3$ v$_2$=v$_3$a,

$\Rightarrow$  M $\vDash$ (v$_3$a)t''aw''=t$_1$ayat$_2$,

$\Rightarrow$ by (4.14$^b$), M $\vDash$ v$_3$=t$_1$ v t$_1$Bv$_3$,

$\Rightarrow$ M $\vDash$ t$_1 \subseteq_p$v$_3 \subseteq_p$v$_2 \subseteq_p$v,

$\Rightarrow$ M $\vDash$ t'<t$_1 \subseteq_p$v,  contradicting  M $\vDash$ Pref(v,t').

(2biii2b3)  M $\vDash$ w'at'v=w$_1$a.

$\Rightarrow$ from  M $\vDash$ w'at'vt''aw''=x=w$_1$at$_1$ayat$_2$, by (3.7),  M $\vDash$ t''aw''=t$_1$ayat$_2$,

$\Rightarrow$ by (4.14$^b$), M $\vDash$ t''=t$_1$=t$_2$=t.

Let  z=w'at'vt'.

Since  M $\vDash$ zBx, we have that  M $\vDash$ J(z).

$\Rightarrow$ by (5.53),  M $\vDash$ Env(t',z).

We now argue that $M \vDash \neg(y \, \varepsilon \, z)$.

Assume, for a reductio, that $M \vDash y \, \varepsilon \, z$.

$\Rightarrow M \vDash \exists t_3, t_4 \, Fr(z, t_3, aya, t_4)$,

$\Rightarrow$ from $M \vDash Env(t', z)$, $M \vDash t_3 \leq t'$,

$\Rightarrow$ by (5.51), $M \vDash Fr(x, t_3, aya, t_4) \lor Fr(x, t_3, aya, t'')$,

$\Rightarrow$ from $M \vDash Fr(x, t_1, aya, t_2)$ & $Env(t, x)$, $M \vDash t_3 = t_1 = t$,

$\Rightarrow M \vDash t' < t'' \leq t = t_3 \leq t'$, contradicting $M \vDash t' \in I \subseteq I_0$.

Therefore $M \vDash \neg(y \, \varepsilon \, z)$.

Next, we show that $M \vDash \forall w (w \, \varepsilon \, z \leftrightarrow w \, \varepsilon \, x \, \& \, w \neq y)$.

Assume $M \vDash w \, \varepsilon \, z$.

$\Rightarrow M \vDash \exists t_3, t_4 \, Fr(z, t_3, awa, t_4)$,

$\Rightarrow$ by (5.51), $M \vDash Fr(x, t_3, awa, t_4) \lor Fr(x, t_3, awa, t'')$,

$\Rightarrow M \vDash w \, \varepsilon \, x$.

That $M \vDash w \neq y$ follows from $M \vDash w \, \varepsilon \, z \, \& \, \neg(y \, \varepsilon \, z)$.

Conversely, assume $M \vDash w \, \varepsilon \, x \, \& \, w \neq y$.

$\Rightarrow M \vDash \exists t_3, t_4 \, Fr(x, t_3, awa, t_4)$,

$\Rightarrow$ from $M \vDash Lastf(x, t_1, aya, t_2)$, by (5.15),

$\quad \quad \quad M \vDash Firstf(x, t_3, awa, t_4) \lor \exists w_3 \, Intf(x, w_3, t_3, awa, t_4)$,

$\Rightarrow$ by (5.38), $M \vDash \exists t_5 \, Fr(z, t_3, awa, t_5)$,

$\Rightarrow M \vDash w \, \varepsilon \, z$.

Finally, assume $M \vDash Fr(x, t_3, awa, t_4) \, \& \, w \neq y$.

Then, by the argument just given, $M \vDash \exists t_5 \, Fr(z, t_3, awa, t_5)$, as required.

(2biii2c)  $M \vDash Lastf(x,t',v,t'')$.

$\Rightarrow M \vDash Pref(v,t')$ & $Tally_b(t'')$ & $t'=t''$ &

& $(x=t'vt'' \lor \exists w'(x=w'at'vt''$ & $Max^+T_b(t',w')))$.

If $M \vDash x=t'vt''$, we derive a contradiction the same way as in (2biii2a1).

$\Rightarrow M \vDash \exists w'(x=w'at'vt''$ & $Max^+T_b(t',w'))$,

$\Rightarrow M \vDash w'at'vt''=x=w_1at_1ayat_2$,

$\Rightarrow$ from $M \vDash Pref(v,t')$, $M \vDash \exists v_0\ av_0a=v$,

$\Rightarrow$ from $M \vDash Lastf(x,t',av_0a,t'')$ & $Lastf(x,t_1,aya,t_2)$, by (5.15), $M \vDash v_0=y$,

$\Rightarrow$ from $M \vDash w'at'av_0at''=x=w_1at_1ayat_2$, by $(4.24^b)$, $M \vDash t''=t_2$,

$\Rightarrow$ by (3.6),  $M \vDash w'at'=w_1at_1$,

$\Rightarrow$ by $(4.24^b)$,  $M \vDash t'=t_1$,

$\Rightarrow$ by (3.6),  $M \vDash w'=w_1$.

From  $M \vDash Occ(w_1,a,t_1ayat_2,x,t',v,t'')$, we distinguish three subcases:

(2biii2c1)  $M \vDash w'at'=w_1$.

$\Rightarrow M \vDash w_1at'=w_1$,

$\Rightarrow M \vDash w_1Bw_1$, contradicting  $M \vDash w_1 \in I \subseteq I_0$.

(2biii2c2)  $M \vDash \exists v_1(v_1Bv$ & $w'at'v_1=w_1a)$.

$\Rightarrow M \vDash w_1at'v_1=w_1a$,

$\Rightarrow M \vDash (w_1a)B(w_1a)$, contradicting  $M \vDash (w_1a) \in I \subseteq I_0$.

(2biii2c3)  $M \vDash w'at'v=w_1a$.

$\Rightarrow M \vDash w_1at'v=w_1a$,

$\Rightarrow M \vDash (w_1a)B(w_1a)$, contradicting  $M \vDash (w_1a) \in I \subseteq I_0$.

Finally, we show that $M \vDash \text{Lex}^+(x) \rightarrow \text{Lex}^+(z)$.

Assume $M \vDash \text{Lex}^+(x)$.

$\Rightarrow$ $M \vDash \forall u,v \, (u <_x v \rightarrow u \prec v)$.

Assume now that $M \vDash u <_z v$.

$\Rightarrow$ $M \vDash \exists t_1, t_2 \, Fr(z, t_1, aua, t_2) \, \& \, \exists t_3, t_4 \, Fr(z, t_3, ava, t_4)$,

$\Rightarrow$ by (9.14), $M \vDash t_1 < t_3$,

$\Rightarrow$ from $M \vDash Env(t,x) \, \& \, u \, \varepsilon \, x \, \& \, v \, \varepsilon \, x$ and earlier argument,

$\quad\quad M \vDash \exists t_5 Fr(x, t_1, aua, t_5) \, \& \, \exists t_6 Fr(x, t_3, ava, t_5) \, \& \, t_1 < t_3$,

$\Rightarrow$ by (9.14), $M \vDash u <_x v$,

$\Rightarrow$ $M \vDash u \prec v$.

Thus, $M \vDash \forall u,v \, (u <_z v \rightarrow u \prec v)$, that is, $M \vDash \text{Lex}^+(z)$.

This completes the proof of THE SUBTRACTION LEMMA.

SPECIAL SET CODES LEMMA. (10.2) For any string concept $I \subseteq I_0$ there is a string concept $J \subseteq I$ such that

$QT^+ \vdash \forall x,y \in J$ ($Lex^+(x)$ & $Lex^+(y)$ & $Special(x)$ & $Special(y)$ & $x \sim y$ &

& $Fr(x,t_1,aua,t_2)$ & $Fr(y,t_3,aua,t_4) \rightarrow t_1=t_3$).

Let $I' \equiv I_{9.1}$ & $I_{9.10}$ & $I_{9.14}$ & $I_{9.25}$.

Let $J(t)$ abbreviate

$I'(t)$ & $\forall x,y \in I', t_2, t_3, t_4, u$ ($Lex^+(x)$ & $Lex^+(y)$ & $Special(x)$ & $Special(y)$ & $x \sim y$ &

& $Env(t,x)$ → $\forall t' \leq t$ ($Fr(x,t',aua,t_2)$ & $Fr(y,t_3,aua,t_4) \rightarrow t'=t_3$)).

Claim 1. $J(t)$ is a string concept.

Let $t=b$.

Assume $M \vDash J(x)$ & $J(y) \in I'$ along with

$M \vDash Lex^+(x)$ & $Lex^+(y)$ & $Special(x)$ & $Special(y)$ & $x \sim y$ & $Env(t,x)$.

Let $M \vDash t' \leq t$ & $Fr(x,t',aua,t_2)$ & $Fr(y,t_3,aua,t_4)$. Then $M \vDash Tally_b(t')$. Since

$QT^+ \vdash Tally_b(t')$ & $t' \leq b \rightarrow t'=b$,

we have that $M \vDash Fr(x,b,aua,t_2)$ & $Fr(y,t_3,aua,t_4)$.

⇒ from hypothesis $M \vDash Special(x)$ & $Special(y)$, we have that

$M \vDash Set(x)$ & $Set(y)$.

Now, $QT^+ \vdash Tally_b(t_0) \rightarrow \neg(t_0<b)$.

Assume $M \vDash w \: \varepsilon \: x$.

Suppose, for a reductio, that $M \vDash w <_x u$.

$\Rightarrow$ by (9.14), $M \vDash \exists t_5, t_6 (Fr(x,t_5,awa,t_6)\ \&\ t_5<b)$,

$\Rightarrow M \vDash Tally_b(t_5)\ \&\ t_5<b$, a contradiction.

Therefore, $M \vDash \forall w\ (w\ \varepsilon\ x \to \neg(w<_x u))$,

$\Rightarrow$ by (9.7), $M \vDash \forall w\ (w\ \varepsilon\ x \to u\leq_x w)$.

So we have

$\quad M \vDash Env(t,x)\ \&\ Fr(x,t',aua,t_2)\ \&\ \forall w\ (w\ \varepsilon\ x \to u\leq_x w)$,

$\Rightarrow$ by (9.10), $M \vDash Firstf(x,t',aua,t_2)$,

$\Rightarrow$ from $M \vDash Set(x)\ \&\ Set(y)\ \&\ x\sim y\ \&\ Lex^+(x)\ \&\ Lex^+(y)$, by (9.25),

$\quad M \vDash \forall w\ (w\ \varepsilon\ y \to u\leq_y w)$,

$\Rightarrow$ from $M \vDash Set(y)\ \&\ Fr(y,t_3,aua,t_4)$, by (9.10), $M \vDash Firstf(y,t_3,aua,t_4)$,

$\Rightarrow$ by (9.1), $M \vDash \forall w\ \neg(w<_y u)$,

$\Rightarrow$ since $M \vDash u\ \varepsilon\ x\ \&\ Set(x)\ \&\ u\ \varepsilon\ y\ \&\ Set(y)$,

$\quad M \vDash \forall t_0 (Tally_b(t_0) \to Max^+(t_0,u,y))$.

Hence, in particular, $M \vDash Max^+(t_3,u,x)\ \&\ Max^+(t',u,y)$.

From $M \vDash Special(x)$, we have that $M \vDash MMax^+T_b(t',u,x)$, whence $M \vDash t' \in I$.

Likewise, from $M \vDash Special(y)$, we obtain $M \vDash MMax^+T_b(t_3,u,y)$.

On the other hand, from $M \vDash Fr(x,t',aua,t_2)\ \&\ Fr(y,t_3,aua,t_4)$, we have that

$\quad M \vDash Max^+T_b(t',u)\ \&\ Max^+T_b(t_3,u)$.

Hence from $M \vDash MMax^+T_b(t',u,x)\ \&\ Max^+(t_3,u,x)\ \&\ Max^+T_b(t_3,u)$ we have that

$\quad M \vDash t'\leq t_3$,

whereas from $M \vDash MMax^+T_b(t_3,u,y)\ \&\ Max^+(t',u,y)\ \&\ Max^+T_b(t',u)$ we have

$\quad M \vDash t_3 \leq t'$.

Since $M \vDash t' \in I \subseteq I_0$, it follows by (2.2) that $M \vDash t'=t_3$.

Therefore $M \vDash J(b)$.

If $t=a$, we have that $M \vDash J(a)$ holds trivially because then $M \vDash \neg Env(t,x)$.

Assume now that $M \vDash J(t)$.

Again, $M \vDash J(ta)$ holds trivially as for $t=a$.

To show that $M \vDash J(tb)$, assume that $M \vDash I'(x)\ \&\ I'(y)$ and

$$M \vDash Lex^+(x)\ \&\ Lex^+(y)\ \&\ Special(x)\ \&\ Special(y)\ \&\ x \sim y\ \&\ Env(tb,x).$$

Let $M \vDash t' \leq tb\ \&\ Fr(x,t',aua,t_2)\ \&\ Fr(y,t_3,aua,t_4)$.

$\Rightarrow\ M \vDash Tally_b(tb)\ \&\ (t' \leq t\ v\ t'=tb)$.

If $M \vDash t' \leq t$, then $M \vDash t'=t_3$ follows from hypothesis $M \vDash J(t)$.

So we may assume that $M \vDash t'=tb$.

From $M \vDash Special(x)\ \&\ Special(y)$ we have

$$M \vDash MMax^+T_b(tb,u,x)\ \&\ MMax^+T_b(t_3,u,y).$$

$\Rightarrow\ M \vDash Max^+(tb,u,x)\ \&\ Max^+(t_3,u,y)$.

Then we claim:

(1.1) $M \vDash Max^+(tb,u,y)$.

To see this, note that we have $M \vDash Set(y)\ \&\ u\ \varepsilon\ y\ \&\ Tally_b(tb)$.

Assume that $M \vDash Fr(y,t_5,awa,t_6)\ \&\ w <_y u$.

$\Rightarrow$ from $M \vDash x \sim y$, $M \vDash w\ \varepsilon\ x$,

$\Rightarrow$ by (9.25), from $M \vDash Set(x)\ \&\ Set(y)\ \&\ x \sim y\ \&\ Lex^+(x)\ \&\ Lex^+(y)$,

$$M \vDash w <_x u,$$

$\Rightarrow$ by (9.14), $M \vDash \exists t_7, t_8\ (Fr(x,t_7,awa,t_8)\ \&\ t_7 < t'=tb)$,

$\Rightarrow$ $M \vDash t_7 \leq t$.

So we have that $\quad M \vDash Fr(x,t_7,awa,t_8)$ & $Fr(y,t_5,awa,t_6)$.

Hence from hypothesis $M \vDash J(t)$, $M \vDash t_5 = t_7$, whence $M \vDash t_5 < tb$.

Therefore, $M \vDash Fr(y,t_5,awa,t_6)$ & $w <_y u \rightarrow t_5 < tb$ which suffices to prove (1.1).

Next, we claim:

(1.2) $M \vDash Max^+(t_3,u,x)$.

Again, we have that $M \vDash Set(x)$ & $u \, \varepsilon \, x$ & $Tally_b(t_3)$.

Assume that $M \vDash Fr(x,t_5,awa,t_6)$ & $w <_x u$,

$\Rightarrow$ from $M \vDash x \sim y$, $M \vDash w \, \varepsilon \, y$,

$\Rightarrow$ by (9.25), $M \vDash w <_y u$,

$\Rightarrow$ by (9.14), $M \vDash \exists t_7, t_8 \, (Fr(y,t_7,awa,t_8)$ & $t_7 < t_3)$.

But from $M \vDash Fr(x,t_5,awa,t_6)$ & $w <_x u$, by (9.14), we have that $M \vDash t_5 < tb$.

$\Rightarrow$ $M \vDash t_5 \leq t$.

But we have that $\quad M \vDash Fr(x,t_5,awa,t_6)$ & $Fr(y,t_7,awa,t_8)$,

$\Rightarrow$ from $M \vDash J(t)$, $M \vDash t_5 = t_7$,

$\Rightarrow$ $M \vDash t_5 < t_3$.

We have thus shown that $M \vDash Fr(x,t_5,awa,t_6)$ & $w <_x u \rightarrow t_5 < t_3$,

which suffices to establish (1.2).

Now, from hypothesis $M \vDash Fr(x,tb,aua,t_2)$ & $Fr(y,t_3,aua,t_4)$ we have that

$\qquad M \vDash Max^+T_b(tb,u)$ & $Max^+T_b(t_3,u)$.

Then from $M \vDash MMax^+T_b(tb,u,x)$ & $MMax^+T_b(t_3,u,y)$ and (1.1)-(1.2) we have that

$\qquad M \vDash tb \leq t_3$ & $t_3 \leq tb$.

Therefore $M \vDash t'=tb=t_3$, as required.

Hence $M \vDash J(tb)$, which completes the proof of Claim 1.

We now prove the main claim.

Assume $M \vDash \text{Lex}^+(x)$ & $\text{Lex}^+(y)$ & $\text{Special}(x)$ & $\text{Special}(y)$ & $x \sim y$

along with $M \vDash \text{Fr}(x,t_1,aua,t_2)$ & $\text{Fr}(y,t_3,aua,t_4)$ where $M \vDash J(x)$ & $J(y) \subseteq I'$.

$\Rightarrow M \vDash \text{Set}(x)$.

Let $M \vDash \text{Env}(t,x)$.

Since we may assume that J is closed under $\subseteq_p$, we have from $M \vDash J(x)$ that $M \vDash J(t)$. But from $M \vDash \text{Env}(t,x)$, we have $M \vDash \text{MaxT}_b(tb,x)$. Hence $M \vDash t_1 \leq t$.

So from $M \vDash J(t)$ it follows that $M \vDash t_1=t_3$, as required.

This completes the proof of the SPECIAL SET CODES LEMMA.

(10.3) For any string concept $I \subseteq I_0$ there is a string concept $J \subseteq I$ such that

$QT^+ \vdash \forall x \in J \ \forall t, t_1, t_2, t_3, t_4, u, v \ [MinSet(x) \ \& \ Fr(x, t_1, aua, t_2) \ \& \ Fr(x, t_3, ava, t_4) \ \&$

$\& \ u <_x v \ \& \ \neg \exists w(u <_x w \ \& \ w <_x v) \ \rightarrow \ (t_2 = t_4 \ \& \ (x = t_1 auat_2 avat_4 \ v$

$v \ \exists w_1(x = w_1 at_1 auat_2 avat_4 \ \& \ Max^+ T_b(t_1, w_1)))) \ v \ (t_2 < t_4 \ \& \ Pref(ava, t_2) \ \&$

$\& \ \exists w_1, w_2((x = w_1 at_1 auat_2 avat_4 aw_2 \ \& \ Max^+ T_b(t_1, w_1)) \ v \ x = t_1 auat_2 avat_4 aw_2))]$.

Let $J \equiv I_{3.10} \ \& \ I_{4.5} \ \& \ I_{4.20} \ \& \ I_{5.22} \ \& \ I_{5.28} \ \& \ I_{5.34} \ \& \ I_{9.7}$.

Assume $M \vDash MinSet(x) \ \& \ Fr(x, t_1, aua, t_2) \ \& \ Fr(x, t_3, ava, t_4)$

where $M \vDash u <_x v \ \& \ \neg \exists w(u <_x w \ \& \ w <_x v)$ and $M \vDash J(x)$.

$\Longrightarrow M \vDash Set(x) \ \& \ x \neq aa$,

$\Longrightarrow$ by (5.18), $M \vDash \exists t \ Env(t, x)$.

We distinguish cases based on $M \vDash Fr(x, t_3, ava, t_4)$.

(1) $M \vDash Lastf(x, t_3, ava, t_4)$.

$\Longrightarrow M \vDash Pref(ava, t_3) \ \& \ Tally_b(t_4) \ \& \ t_3 = t_4 = t \ \&$

$\& \ (x = t_3 avat_4 \ v \ \exists w(x = wat_3 avat_4 \ \& \ Max^+ T_b(t_3, w)))$.

Suppose, for a reductio, that $M \vDash x = t_3 avat_4$.

$\Longrightarrow$ from $M \vDash Pref(ava, t_3) \ \& \ Tally_b(t_4) \ \& \ t_3 = t_4$,

$M \vDash Firstf(x, t_3, ava, t_4) \ \& \ Lastf(x, t_3, ava, t_4)$,

$\Longrightarrow$ by (5.22), $M \vDash \exists y \forall w(w \ \varepsilon \ x \leftrightarrow w = y)$,

$\Longrightarrow$ from $M \vDash u \ \varepsilon \ x \ \& \ v \ \varepsilon \ x$, $M \vDash u = v$, which contradicts $M \vDash u <_x v$ by (9.4).

Therefore, $M \vDash \neg(x = t_3 avat_4)$.

$\Rightarrow$ M ⊨ ∃w(x=wat$_3$avat$_4$ & Max$^+$T$_b$(t$_3$,w)).

$\Rightarrow$ from M ⊨ Fr(x,t$_1$,aua,t$_2$), by (5.15) and (9.4), M ⊨ ¬Lastf(x,t$_1$,aua,t$_2$),

$\Rightarrow$ M ⊨ Firstf(x,t$_1$,aua,t$_2$) v ∃w$_1$Intf(x,w$_1$,t$_1$,aua,t$_2$).

(1a)  M ⊨ Firstf(x,t$_1$,aua,t$_2$).

$\Rightarrow$ M ⊨ Pref(aua,t$_1$) & Tally$_b$(t$_2$) & ((t$_1$=t$_2$ & x=t$_1$auat$_2$) v

$\qquad\qquad\qquad\qquad\qquad\qquad$ v (t$_1$<t$_2$ & (t$_1$auat$_2$a)Bx)),

$\Rightarrow$ as in (1) above, M ⊨ ¬(x=t$_1$auat$_2$),

$\Rightarrow$ M ⊨ t$_1$<t$_2$ & (t$_1$auat$_2$a)Bx,

$\Rightarrow$ M ⊨ ∃x$_1$ t$_1$auat$_2$ax$_1$=x=wat$_3$avat$_4$,

$\Rightarrow$ M ⊨ (t$_1$auat$_2$)Bx & wBx,

$\Rightarrow$ by (3.8), M ⊨ wB(t$_1$auat$_2$) v w=t$_1$auat$_2$ v (t$_1$auat$_2$)Bw.

(1ai)  M ⊨ wB(t$_1$auat$_2$).

$\Rightarrow$ M ⊨ ∃w' ww'=t$_1$auat$_2$,

$\Rightarrow$ from M ⊨ (t$_1$auat$_2$a)Bx, M ⊨ ∃x$_1$ ww'ax$_1$=x=wat$_3$avat$_4$,

$\Rightarrow$ by (3.7), M ⊨ w'ax$_1$=at$_3$avat$_4$,

$\Rightarrow$ M ⊨ w'=a v aBw'.

Now, M ⊨ w'=a is ruled out because M ⊨ Tally$_b$(t$_2$).

$\Rightarrow$ M ⊨ aBw',

$\Rightarrow$ M ⊨ ∃w'' aw''=w',

$\Rightarrow$ M ⊨ (aw'')ax$_1$=at$_3$avat$_4$,

$\Rightarrow$ M ⊨ w''ax$_1$=t$_3$avat$_4$,

$\Rightarrow$ by (4.14$^b$), M ⊨ w''=t$_3$ v t$_3$Bw'',

$\Rightarrow$ $M \vDash t_3 \subseteq_p w''$,

$\Rightarrow$ $M \vDash t_3 \subseteq_p w'' \subseteq_p w' \subseteq_p t_1 a u a t_2$,

$\Rightarrow$ by (4.17$^b$), $M \vDash t_3 \subseteq_p t_1$ v $t_3 \subseteq_p u a t_2$,

$\Rightarrow$ from $M \vDash t_3 = t_4 = t$ & $Env(t,x)$, $M \vDash MaxT_b(t_3,x)$,

$\Rightarrow$ $M \vDash t_1 < t_2 \leq t_3$.

Hence, if $M \vDash t_3 \subseteq_p t_1$, then $M \vDash t_1 < t_3 \leq t_1$, contradicting $M \vDash t_1 \in I \subseteq I_0$.

So $M \vDash t_3 \subseteq_p u a t_2$.

$\Rightarrow$ by (4.17$^b$), $M \vDash t_3 \subseteq_p u$ v $t_3 \subseteq_p t_2$.

If $M \vDash t_3 \subseteq_p u$, then from $M \vDash Pref(aua, t_1)$ we have $M \vDash Max^+T_b(t_1, aua)$, so

$M \vDash t_3 < t_1$, a contradiction again. Therefore $M \vDash \neg(t_3 \subseteq_p u)$.

Hence, we may assume that $M \vDash t_3 \subseteq_p t_2$.

But then we have $M \vDash t_2 \leq t_3$ & $t_3 \leq t_2$, whence, since $M \vDash t_2 \in I_0$ & $t_3 \in I_0$, we have,

by (2.2), that $M \vDash t_2 = t_3 = t_4$.

$\Rightarrow$ $M \vDash t_1 a u a t_2 a x_1 = x = w a t_2 a v a t_4$,

$\Rightarrow$ $M \vDash (t_1 a u a) B x$ & $(wa) B x$,

$\Rightarrow$ by (3.8), $M \vDash (t_1 a u a) B (wa)$ v $wa = t_1 a u a$ v $(wa) B (t_1 a u a)$.

(1ai1) $M \vDash (t_1 a u a) B (wa)$.

$\Rightarrow$ $M \vDash \exists w_3\ t_1 a u a w_3 = wa$,

$\Rightarrow$ $M \vDash t_1 a u a t_2 a x_1 = x = (t_1 a u a w_3) t_2 a v a t_4$,

$\Rightarrow$ by (3.7), $M \vDash t_2 a x_1 = w_3 t_2 a v a t_4$,

$\Rightarrow$ from $M \vDash t_1 a u a w_3 = wa$, $M \vDash w_3 = a$ v $a E w_3$,

$\Rightarrow$ since $M \vDash Tally_b(t_3)$, $M \vDash \neg(w_3 = a)$,

$\Rightarrow M \vDash aEw_3$,

$\Rightarrow M \vDash \exists w_4\ w_3=w_4a$,

$\Rightarrow M \vDash t_2ax_1=(w_4a)t_2avat_4$,

$\Rightarrow$ by (4.14$^b$), $M \vDash w_4=t_2\ v\ t_2Bw_4$,

$\Rightarrow M \vDash t_2\subseteq_p w_4$,

$\Rightarrow$ from $M \vDash t_1auaw_3=wa$, $M \vDash t_1aua(w_4a)=wa$,

$\Rightarrow M \vDash t_1auaw_4=w$,

$\Rightarrow M \vDash t_2\subseteq_p w_4\subseteq_p w$,

$\Rightarrow M \vDash t_3=t_2\subseteq_p w$,

$\Rightarrow$ from $M \vDash \text{Max}^+T_b(t_3,w)$, $M \vDash t_3<t_3$, contradicting $M \vDash t_3\in I\subseteq I_0$.

(1ai2) $M \vDash (wa)B(t_1aua)$.

$\Rightarrow M \vDash \exists w_3\ waw_3=t_1aua$,

$\Rightarrow M \vDash (waw_3)t_2ax_1=x=wat_2avat_4$,

$\Rightarrow$ by (3.7), $M \vDash w_3t_2ax_1=t_2avat_4$,

$\Rightarrow$ from $M \vDash waw_3=t_1aua$, $M \vDash w_3=a\ v\ aEw_3$,

$\Rightarrow$ since $M \vDash \text{Tally}_b(t_3)$, $M \vDash \neg(w_3=a)$,

$\Rightarrow M \vDash aEw_3$,

$\Rightarrow M \vDash \exists w_4\ t_2avat_4=(w_4a)t_2ax_1$,

$\Rightarrow$ by (4.16), $M \vDash t_2\subseteq_p ava$,

$\Rightarrow$ from $M \vDash \text{Pref}(ava,t_3)\ \&\ t_2=t_3$, $M \vDash \text{Max}^+T_b(t_2,ava)$,

$\Rightarrow M \vDash t_2<t_2$, contradicting $M \vDash t_2\in I\subseteq I_0$.

(1ai3) $M \vDash wa=t_1aua$.

$\Rightarrow$ M ⊨ $t_1auat_2ax_1=x=wat_2avat_4=t_1auat_2avat_4$ & $t_2=t_4$, as required.

   (1aii)  M ⊨ $t_1auat_2=w$.

$\Rightarrow$ M ⊨ $t_1auat_2ax_1=x=t_1auat_2at_3avat_4$,

$\Rightarrow$ from M ⊨ MinSet(x), M ⊨ $\exists v',t',t''$ Occ($t_1auat_2,a,t_3avat_4,x,t',v',t''$),

$\Rightarrow$ M ⊨ Fr(x,t',v',t''),

$\Rightarrow$ M ⊨ Pref(v',t') & Tally$_b$(t''),

$\Rightarrow$ M ⊨ $\exists v''$ $v'=av''a$ & Max$^+$T$_b$(t',v').

   (1aii1)  M ⊨ Firstf(x,t',av''a,t'').

$\Rightarrow$ M ⊨ ((t'=t'' & x=t'av''at'') v (t'<t'' & (t'av''at''a)Bx)),

$\Rightarrow$ as in (1) above, M ⊨ ¬(x=t'av''at''),

$\Rightarrow$ M ⊨ $\exists x_2$ t'av''at''a$x_2$=x=$t_1auat_2ax_1$=wat$_3$avat$_4$,

$\Rightarrow$ from M ⊨ Firstf(x,$t_1$,aua,$t_2$), by (5.15), M ⊨ u=v''.

From M ⊨ Occ($t_1auat_2,a,t_3avat_4,x,t',v',t''$) we distinguish three scenarios:

    (1aii1a)  M ⊨ t'=$t_1auat_2$.

$\Rightarrow$ M ⊨ a⊆$_p$t', which contradicts M ⊨ Tally$_b$(t').

    (1aii1b)  M ⊨ ($t_1auat_2$a)B(t'av''a).

$\Rightarrow$ M ⊨ $\exists z_1$ $t_1auat_2az_1$=t'av''a,

$\Rightarrow$ by (4.24$^b$), M ⊨ $t_1$=t',

$\Rightarrow$ from M ⊨ u=v'', M ⊨ t'auat$_2$az$_1$=t'aua,

$\Rightarrow$ M ⊨ (t'aua)B(t'aua), contradicting M ⊨ (t'aua)∈I⊆I$_0$.

    (1aii1c)  M ⊨ $t_1auat_2a$=t'av''a.

$\Rightarrow$ just as in (1aii1b), M ⊨ t'auat$_2$a=t'aua,

$\Rightarrow$ M ⊨ (t'aua)B(t'aua), again contradicting M ⊨ (t'aua)∈I⊆I$_0$.

   (1aii2) M ⊨ ∃w'Intf(x,w',t',av''a,t'').

$\Rightarrow$ by (5.19), M ⊨ ¬Firstf(x,t',av''a,t''),

$\Rightarrow$ from M ⊨ Firstf(x,t$_1$,aua,t$_2$), by (5.20), M ⊨ t$_1$<t',

$\Rightarrow$ by definition of <$_x$, M ⊨ u<$_x$v'',

$\Rightarrow$ from M ⊨ MaxT$_b$(t$_4$,x), M ⊨ t'<t''≤t$_4$,

$\Rightarrow$ by definition of <$_x$, M ⊨ v''<$_x$v,

$\Rightarrow$ M ⊨ u<$_x$v'' & v''<$_x$v, contradicting the principal hypothesis.

   (1aii3) M ⊨ Lastf(x,t',av''a,t'').

$\Rightarrow$ from hypothesis M ⊨ Lastf(x,t$_3$,ava,t$_4$), by (5.15), M ⊨ t'=t''=t$_4$=t$_3$ & v''=v,

$\Rightarrow$ M ⊨ ∃w' w'at'av''at''=x=wat$_3$avat$_4$=wat'av''at'',

$\Rightarrow$ by (3.6), M ⊨ w'=w.

From M ⊨ Occ(t$_1$auat$_2$,a,t$_3$avat$_4$,x,t',v',t'') we then have:

   (1aii3a) M ⊨ w'at'=t$_1$auat$_2$.

$\Rightarrow$ M ⊨ wat$_3$=t$_1$auat$_2$,

$\Rightarrow$ from hypothesis (1aii), M ⊨ wat$_3$=w,

$\Rightarrow$ M ⊨ wBw, contradicting M ⊨ w∈I⊆I$_0$.

   (1aii3b) M ⊨ ∃u'(u'Bv' & w'at'u'=t$_1$auat$_2$a).

$\Rightarrow$ from hypothesis (1aii), M ⊨ w'at'u'=wa,

$\Rightarrow$ M ⊨ wat'u'=wa,

$\Rightarrow$ M ⊨ (wa)B(wa), contradicting M ⊨ (wa)∈I⊆I$_0$.

   (1aii3c) M ⊨ w'at'av''a=t$_1$auat$_2$a.

$\Rightarrow$ from hypothesis (1aii), $M \vDash w'at'av''a=wa$,

$\Rightarrow$ $M \vDash wat'av''a=wa$,

$\Rightarrow$ $M \vDash (wa)B(wa)$, again contradicting $M \vDash (wa) \in I \subseteq I_0$.

   (1aiii) $M \vDash (t_1auat_2)Bw$.

$\Rightarrow$ $M \vDash \exists w'\ t_1auat_2w'=w$,

$\Rightarrow$ $M \vDash t_1auat_2ax_1=x=(t_1auat_2w')at_3avat_4$,

$\Rightarrow$ from $M \vDash MinSet(x)$, $M \vDash \exists v',t',t''Occ(t_1auat_2w',a,t_3avat_4,x,t',v',t'')$,

$\Rightarrow$ $M \vDash Fr(x,t',v',t'')$,

$\Rightarrow$ $M \vDash Pref(v',t')\ \&\ Tally_b(t'')$,

$\Rightarrow$ $M \vDash \exists v''\ v'=av''a\ \&\ Max^+T_b(t',v')$.

   (1aiii1) $M \vDash Firstf(x,t',av''a,t'')$.

$\Rightarrow$ from hypothesis (1a), by (5.15), $M \vDash t_1=t'\ \&\ u=v''$.

    (1aiii1a) $M \vDash t'=t_1auat_2w'$.

$\Rightarrow$ $M \vDash a \subseteq_p t'$, which contradicts $M \vDash Tally_b(t')$.

    (1aiii1b) $M \vDash (t_1auat_2w'a)B(t'av''a)$.

$\Rightarrow$ $M \vDash \exists z_1\ t_1auat_2w'az_1=t'av''a$,

$\Rightarrow$ from $M \vDash t_1=t'\ \&\ u=v''$, $M \vDash t_1auat_2w'az_1=t_1aua$,

$\Rightarrow$ $M \vDash (t_1aua)B(t_1aua)$, contradicting $M \vDash (t_1aua) \in I \subseteq I_0$.

    (1aii1c) $M \vDash t_1auat_2w'a=t'av''a$.

A contradiction as in (1aiii1b), omitting $z_1$.

   (1aiii2) $M \vDash \exists w_1 Intf(x,w_1,t',av''a,t'')$.

This is ruled out by the same argument as in (1aii2).

(1aiii3)  $M \vDash Lastf(x,t',av''a,t'')$.

$\Rightarrow$ from hypothesis  $M \vDash Lastf(x,t_3,ava,t_4)$, by (5.15), $M \vDash t'=t''=t_4=t_3$ & $v''=v$,

$\Rightarrow M \vDash \exists w_1\ w_1at'av''at''=wat_3avat_4=wat'av''at''$,

$\Rightarrow$ by (3.7),  $M \vDash w_1=w$.

From  $M \vDash Occ(t_1auat_2w',a,t_3avat_4,x,t',v',t'')$  we then have:

(1aiii3a)  $M \vDash w_1at'=t_1auat_2w'$.

A contradiction as in (1aii3a).

(1aii3b)  $M \vDash \exists u'(u'Bv'\ \&\ w_1at'u'=t_1auat_2w'a)$.

A contradiction as in (1aii3b).

(1aii3c)  $M \vDash w_1at'av''a=t_1auat_2w'a$.

A contradiction as in (1aii3c).

This completes the subcase (1a).

(1b)  $M \vDash \exists w_1\ Intf(x,w_1,t_1,aua,t_2)$.

$\Rightarrow M \vDash Pref(aua,t_1)\ \&\ Tally_b(t_2)\ \&$

$\qquad\qquad\qquad \&\ \exists w_2(x=w_1at_1auat_2aw_2\ \&\ t_1<t_2\ \&\ Max^+T_b(t_1,w_1))$,

$\Rightarrow M \vDash w_1at_1auat_2aw_2=x=wat_3avat_4$,

$\Rightarrow M \vDash (w_1at_1auat_2)Bx\ \&\ wBx$,

$\Rightarrow$ by (3.8), $M \vDash (w_1at_1auat_2)Bw\ \vee\ w_1at_1auat_2=w\ \vee\ wB(w_1at_1auat_2)$.

(1bi)  $M \vDash (w_1at_1auat_2)Bw$.

$\Rightarrow M \vDash \exists w'\ w_1at_1auat_2w'=w$,

$\Rightarrow M \vDash w_1at_1auat_2aw_2=x=(w_1at_1auat_2w')at_3avat_4$,

$\Rightarrow$ from $M \vDash MinSet(x)$, $M \vDash \exists v',t',t''Occ(w_1at_1auat_2w',a,t_3avat_4,x,t',v',t'')$,

$\Rightarrow$ M ⊨ Fr(x,t',v',t''),

$\Rightarrow$ M ⊨ Pref(v',t') & Tally$_b$(t''),

$\Rightarrow$ M ⊨ ∃v'' v'=av''a & Max$^+$T$_b$(t',v').

(1bi1)  M ⊨ Firstf(x,t',av''a,t'').

$\Rightarrow$ as in (1) above, M ⊨ ¬(x=t'av''at''),

$\Rightarrow$ M ⊨ t'<t'' & (t'av''at''a)Bx,

$\Rightarrow$ M ⊨ ∃x$_1$ t'av''at''ax$_1$=x=w$_1$at$_1$auat$_2$w'at$_3$avat$_4$.

We distinguish three scenarios from  M ⊨ Occ(w$_1$at$_1$auat$_2$w',a,t$_3$avat$_4$,x,t',v',t''):

(1bi1a)  M ⊨ t'=w$_1$at$_1$auat$_2$w'.

This is ruled out because M ⊨ Tally$_b$(t').

(1bi1b)  M ⊨ (w$_1$at$_1$auat$_2$w'a)B(t'av''a).

$\Rightarrow$ M ⊨ ∃w'' w$_1$at$_1$auat$_2$w'aw''=t'av''a,

$\Rightarrow$ by (4.14$^b$),  M ⊨ w$_1$=t' ∨ t'Bw$_1$,

$\Rightarrow$  M ⊨ t'at$_1$auat$_2$w'aw''=t'av''a ∨ ∃w$_3$ (t'w$_3$)at$_1$auat$_2$w'aw''=t'av''a,

$\Rightarrow$  by (3.7), M ⊨ t$_1$auat$_2$w'aw''=v''a ∨ w$_3$at$_1$auat$_2$w'aw''=av''a,

$\Rightarrow$ M ⊨ t$_2$⊆$_p$v''.

But from hypotheses (1bi1) and (1b) we have, by (5.20), that   M ⊨ t'≤t$_1$<t$_2$.

$\Rightarrow$ from  M ⊨ Max$^+$T$_b$(t',av''a),  M ⊨ t'<t$_2$<t', contradicting  M ⊨ t'∈I⊆I$_0$.

(1bi1c)  M ⊨ w$_1$at$_1$auat$_2$w'a=t'av''a.

The same argument as in (1bi1b), omitting w''.

(1bi2)  $M \vDash \exists w_3 \text{Intf}(x,w_3,t',av''a,t'')$.

$\Rightarrow M \vDash \text{Pref}(av''a,t')$ & $\text{Tally}_b(t'')$ & $t'<t''$ &

$\qquad\qquad\qquad\qquad$ & $\exists w_4\, x=w_3at'av''at''aw_4$ & $\text{Max}^+T_b(t',w_3)$,

$\Rightarrow M \vDash w_3at'av''at''aw_4 = x = w_1at_1auat_2aw_2 = (w_1at_1auat_2w')at_3avat_4$.

Three subcases from  $M \vDash \text{Occ}(w_1at_1auat_2w', a, t_3avat_4, x, t', v', t'')$:

(1bi2a)  $M \vDash w_3at' = w_1at_1auat_2w'$.

$\Rightarrow$ from  $M \vDash w_1at_1auat_2aw_2 = (w_1at_1auat_2w')at_3avat_4$, by (3.7),

$\qquad\qquad M \vDash aw_2 = w'at_3avat_4$,

$\Rightarrow M \vDash w' = a \lor aBw'$,

$\Rightarrow$ from $M \vDash \text{Tally}_b(t')$,  $M \vDash \neg(w'=a)$,

$\Rightarrow M \vDash aBw'$,

$\Rightarrow M \vDash \exists w''(w_3at' = w_1at_1auat_2(aw'') \& w' = aw'')$,

$\Rightarrow$ by (4.15$^b$), $M \vDash w''=t' \lor t'Ew''$,

$\Rightarrow M \vDash w_3at' = w_1at_1auat_2at' \lor \exists w_5(w_3at' = w_1at_1auat_2a(w_5t') \& w_5t' = w'')$,

$\Rightarrow$ by (3.6), $M \vDash w_3 = w_1at_1auat_2 \lor w_3a = w_1at_1auat_2aw_5$.

If  $M \vDash w_3a = w_1at_1auat_2aw_5$, then $M \vDash w_5 = a \lor aEw_5$, whence

$\qquad M \vDash w_3a = w_1at_1auat_2aa \lor \exists w_6(w_3a = w_1at_1auat_2a(w_6a) \& w_6a = w_5)$.

Therefore  $M \vDash w_3 = w_1at_1auat_2 \lor w_3 = w_1at_1auat_2a \lor w_3 = w_1at_1auat_2aw_6$.

$\Rightarrow M \vDash t_2 \subseteq_p w_3$,

$\Rightarrow$ from  $M \vDash \text{Max}^+T_b(t',w_3)$,  $M \vDash t_2 < t'$.

But then we have, from hypotheses (1b), (1bi2) and (1), respectively, that

$\qquad M \vDash \text{Intf}(x,w_1,t_1,aua,t_2) \& \text{Intf}(x,w_3,t',av''a,t'') \& \text{Lastf}(x,t_3,ava,t_4)$.

$\Rightarrow$ from $M \vDash MaxT_b(t_4,x)$, $M \vDash t_2<t'$ & $t'<t''=t_4=t_3$,

$\Rightarrow$ by definition of $<_x$, $M \vDash u<_xv''$ & $v''<_xv$, contradicting the principal hypothesis.

(1bi2b) $M \vDash \exists u''(u''B(av''a) \& w_3at'u''=w_1at_1auat_2w'a)$.

$\Rightarrow$ $M \vDash \exists v_3\ u''v_3=av''a$,

$\Rightarrow$ $M \vDash w_3at'av''at''aw_4=x=w_3at'(u''v_3)t''aw_4=(w_1at_1auat_2w'a)v_3t''aw_4=$
$=w_1at_1auat_2aw_2=(w_1at_1auat_2w')at_3avat_4$,

$\Rightarrow$ by (3.7), $M \vDash v_3t''aw_4=t_3avat_4$,

$\Rightarrow$ from $M \vDash u''v_3=av''a$, $M \vDash v_3=a \lor aEv_3$,

$\Rightarrow$ from $M \vDash Tally_b(t_3)$, $M \vDash \neg(v_3=a)$,

$\Rightarrow$ $M \vDash \exists v_4\ t_3avat_4=(v_4a)t''aw_4$,

$\Rightarrow$ by (4.14[b]), $M \vDash v_4=t_3 \lor t_3Bv_4$,

$\Rightarrow$ $M \vDash t_3\subseteq_p v_4\subseteq_p v_3\subseteq_p av''a$,

$\Rightarrow$ from $M \vDash Max^+T_b(t',v')$, $M \vDash t_3<t'$,

$\Rightarrow$ from hypothesis (1), $M \vDash t'\leq t_3$,

$\Rightarrow$ $M \vDash t'\leq t_3<t'$, contradicting $M \vDash t'\in I\subseteq I_0$.

(1bi2c) $M \vDash w_3at'av''a=w_1at_1auat_2w'a$.

$\Rightarrow$ $M \vDash w_3at'av''at''aw_4=x=(w_1at_1auat_2w'a)t''aw_4$,

$\Rightarrow$ $M \vDash (w_1at_1auat_2w'a)t_3avat_4=x=(w_3at'av''a)t_3avat_4$,

$\Rightarrow$ $M \vDash w_3at'av''at''aw_4=x=(w_3at'av''a)t_3avat_4$,

$\Rightarrow$ by (3.7), $M \vDash t''aw_4=t_3avat_4$,

$\Rightarrow$ by (4.24[b]), $M \vDash t''=t_3$,

$\Rightarrow$ $M \vDash t'<t''=t_3=t_4$,

$\Rightarrow$ from $M \vDash Fr(x,t',av''a,t'')$ & $Lastf(x,t_3,ava,t_4)$,  $M \vDash v''<_x v$,

$\Rightarrow$ from $M \vDash w_3at'av''at''aw_4=x=w_1at_1auat_2w'at''aw_4$, by (3.6),

$\qquad\qquad M \vDash w_3at'av''=w_1at_1auat_2w'$,

$\Rightarrow$ $M \vDash (t'av'')E(w_3at'av'')$ & $(t_2w')E(w_3at'av'')$,

$\Rightarrow$ by (3.10), $M \vDash (t'av'')E(t_2w')$ v $t'av''=t_2w'$ v $(t_2w')E(t'av'')$.

$\quad$ (1bi2ci)  $M \vDash (t'av'')E(t_2w')$.

$\Rightarrow$ $M \vDash \exists w_5\ t_2w'=w_5t'av''$,

$\Rightarrow$ $M \vDash w_3at'av''=w_1at_1aua(w_5t'av'')$,

$\Rightarrow$ by (3.6), $M \vDash w_3a=w_1at_1auaw_5$,

$\Rightarrow$ $M \vDash w_5=a$ v $aEw_5$,

$\Rightarrow$ from $M \vDash Tally_b(t_2)$, $M \vDash \neg(w_5=a)$,

$\Rightarrow$ $M \vDash \exists w_6(t_2w'=(w_6a)t'av''$ & $w_5=w_6a)$,

$\Rightarrow$ by (4.14[b]), $M \vDash w_6=t_2$ v $t_2Bw_6$,

$\Rightarrow$ $M \vDash t_2 \subseteq_p w_6$,

$\Rightarrow$ from $M \vDash w_1at_1aua(t_2w')=w_1at_1aua(w_6at'av'')=w_3at'av''$, by (3.6),

$\qquad\qquad M \vDash w_1at_1auaw_6=w_3$,

$\Rightarrow$ $M \vDash t_2 \subseteq_p w_6 \subseteq_p w_3$,

$\Rightarrow$ from $M \vDash Max^+T_b(t',w_3)$, $M \vDash t_2<t'$,

$\Rightarrow$ from hypotheses (1b) and (1bi2), $M \vDash u<_x v''$,

which along with $M \vDash v''<_x v$ contradicts the principal hypothesis.

   (1bi2cii)  $M \vDash t'av''=t_2w'$.

$\Rightarrow$ from $M \vDash (w_1at_1auat_2)aw_2=x=(w_1at_1auat_2)w'at''aw_4$, by (3.7),

     $M \vDash aw_2=w'at''aw_4$,

$\Rightarrow$ $M \vDash w'=a \ v \ aBw'$,

$\Rightarrow$ $M \vDash t'av''=t_2a \ v \ \exists w_5 \ t'av''=t_2(aw_5)$,

$\Rightarrow$ either way, by (4.23$^b$), $M \vDash t_2=t'$,

$\Rightarrow$ from hypotheses (1b) and (1bi2), $M \vDash u<_x v''$,

which along with $M \vDash v''<_x v$ contradicts the principal hypothesis.

  (1bi2ciii)  $M \vDash (t_2w')E(t'av'')$.

$\Rightarrow$ $M \vDash \exists w_5(t_2w'w_5=t'av'')$,

$\Rightarrow$ just as in (1bi2cii), $M \vDash w'=a \ v \ aBw'$, and by (4.23$^b$), $M \vDash t_2=t'$, whence $M \vDash u<_x v''$, a contradiction again.

  (1bi3)  $M \vDash Lastf(x,t',av''a,t'')$.

$\Rightarrow$ $M \vDash Pref(av''a,t')$ & $Tally_b(t'')$ & $t'=t''$ &

    & $(x=t'av''at'' \ v \ \exists w''(x=w''at'av''at''$ & $Max^+T_b(t',w'')))$,

$\Rightarrow$ as in (1), $M \vDash \neg(x=t'av''at'')$,

$\Rightarrow$ from (1) by (5.15) and (4.24$^b$), $M \vDash v''=v$ & $t''=t_4$,

$\Rightarrow$ $M \vDash \exists w'' \ (w_1at_1auat_2w')at_3avat_4=wat_3avat_4= x=w''at'av''at''$.

Three subcases from $M \vDash Occ(w_1at_1auat_2w',a,t_3avat_4,x,t',v',t'')$:

  (1bi3a)  $M \vDash w''at'=w_1at_1auat_2w'$.

$\Rightarrow$ from $M \vDash (w_1at_1auat_2w')at_3avat_4=x=wat_3avat_4$, by (3.6),

    $M \vDash w_1at_1auat_2w'=w$,

$\Rightarrow$ $M \vDash w=w''at'$,

$\Rightarrow$ $M \vDash t' \subseteq_p w$,

$\Rightarrow$ from $M \vDash t'=t''=t_4=t_3$ & $Max^+T_b(t_3,w)$, $M \vDash t'<t_3=t'$,

contradicting $M \vDash t' \in I \subseteq I_0$.

    (1bi3b) $M \vDash \exists u''(u''B(av''a)$ & $w''at'u''=w_1at_1auat_2w'a)$.

$\Rightarrow$ $M \vDash u''=a \lor aEu''$,

$\Rightarrow$ $M \vDash w''at'a=w_1at_1auat_2w'a \lor \exists v_3\ w''at'(v_3a)=w_1at_1auat_2w'a$,

$\Rightarrow$ $M \vDash w''at'=w_1at_1auat_2w' \lor w''at'v_3=w_1at_1auat_2w'$,

$\Rightarrow$ just as in (1bi3a), $M \vDash w_1at_1auat_2w'=w$,

$\Rightarrow$ $M \vDash w''at'=w \lor w''at'v_3=w$,

$\Rightarrow$ $M \vDash t' \subseteq_p w$.

Then a contradiction follows as in (1bi3a).

    (1bi3c) $M \vDash w''at'av''a=w_1at_1auat_2w'a$.

$\Rightarrow$ from $M \vDash w_1at_1auat_2w'=w$, $M \vDash w''at'av''a=wa$,

$\Rightarrow$ $M \vDash w''at'av''=w$,

$\Rightarrow$ $M \vDash t' \subseteq_p w$, again a contradiction as in (1bi3a).

  (1bii) $M \vDash w_1at_1auat_2=w$.

$\Rightarrow$ $M \vDash w_1at_1auat_2aw_2=x=(w_1at_1auat_2)at_3avat_4$,

$\Rightarrow$ from $M \vDash MinSet(x)$, $M \vDash \exists v',t',t''Occ(w_1at_1auat_2,a,t_3avat_4,x,t',v',t'')$,

$\Rightarrow$ $M \vDash Fr(x,t',v',t'')$,

$\Rightarrow$ $M \vDash \exists v''\ v'=av''a$ & $Max^+T_b(t',v')$.

(1bii1)  $M \vDash \text{Firstf}(x,t',av''a,t'')$.

$\Rightarrow$ as in (1) above, $M \vDash \neg(x=t'av''at'')$,

$\Rightarrow M \vDash t'<t''$ & $(t'av''at''a)Bx$,

$\Rightarrow M \vDash \exists x_1\, t'av''at''ax_1=x=w_1at_1auat_2at_3avat_4$.

We distinguish three scenarios from $M \vDash \text{Occ}(w_1at_1auat_2,a,t_3avat_4,x,t',v',t'')$:

(1bii1a)  $M \vDash t'=w_1at_1auat_2$.

This is ruled out because $M \vDash \text{Tally}_b(t')$.

(1bii1b)  $M \vDash (w_1at_1auat_2a)B(t'av''a)$.

We derive a contradiction exactly as in (1bi1b), omitting w'.

(1bii1c)  $M \vDash w_1at_1auat_2a=t'av''a$.

A contradiction exactly as in (1bi1c).

(1bii2)  $M \vDash \exists w_3 \text{Intf}(x,w_3,t',av''a,t'')$.

$\Rightarrow M \vDash \text{Pref}(av''a,t')$ & $\text{Tally}_b(t'')$ & $t'<t''$ &

$\qquad\qquad\qquad$ & $\exists w_4\, x=w_3at'av''at''aw_4$ & $\text{Max}^+T_b(t',w_3)$,

$\Rightarrow M \vDash w_3at'av''at''aw_4=x=w_1at_1auat_2aw_2=w_1at_1auat_2at_3avat_4$.

Again, three subcases from $M \vDash \text{Occ}(w_1at_1auat_2w',a,t_3avat_4,x,t',v',t'')$:

(1bii2a)  $M \vDash w_3at'=w_1at_1auat_2$.

$\Rightarrow M \vDash w_3at'=w_1at_1auat_2$,

$\Rightarrow$ by (4.24$^b$), $M \vDash t_2=t'$,

$\Rightarrow$ from hypotheses (1b) and (1bii2),

$\qquad M \vDash \text{Intf}(x,w_1,t_1,aua,t_2)$ & $\text{Intf}(x,w_3,t',av''a,t'')$,

$\Rightarrow$ M ⊨ u<$_x$v'',

$\Rightarrow$ from hypothesis (1), M ⊨ Lastf(x,t$_3$,ava,t$_4$),

$\Rightarrow$ from M ⊨ Env(t,x), M ⊨ MaxT$_b$(t$_4$,x),

$\Rightarrow$ from M ⊨ t'<t''=t$_4$=t$_3$.

We then proceed just as in (1bi2a) to derive M ⊨ v''<$_x$v and a contradiction.

(1bii2b) M ⊨ ∃u''(u''B(av''a) & w$_3$at'u''=w$_1$at$_1$auat$_2$a).

We derive a contradiction exactly as in (1bi2b), omitting w'.

(1bii2c) M ⊨ w$_3$at'av''a=w$_1$at$_1$auat$_2$a.

We proceed to derive, exactly as in (1bi2c), that M ⊨ v''<$_x$v.

On the other hand, we also have M ⊨ w$_3$at'av''at''aw$_4$=w$_1$at$_1$auat$_2$at''aw$_4$,

$\Rightarrow$ by (3.6), M ⊨ w$_3$at'av''=w$_1$at$_1$auat$_2$,

$\Rightarrow$ by (4.15$^b$), M ⊨ v''=t$_2$ v t$_2$Ev'',

$\Rightarrow$ M ⊨ t$_2$⊆$_p$v''⊆$_p$av''a,

$\Rightarrow$ from M ⊨ Pref(av''a,t'), M ⊨ Max$^+$T$_b$(t',v''),

$\Rightarrow$ M ⊨ t$_2$<t',

$\Rightarrow$ from M ⊨ Intf(x,w$_1$,t$_1$,aua,t$_2$) & Intf(x,w$_3$,t',av''a,t''), M ⊨ u<$_x$v'',

which, along with M ⊨ v''<$_x$v contradicts the principal hypothesis.

(1bii3) M ⊨ Lastf(x,t',av''a,t'').

The argument is exactly the same as in (1bi3), omitting w'.

(1biii) M ⊨ wB(w$_1$at$_1$auat$_2$).

$\Rightarrow$ M ⊨ ∃w' ww'=w$_1$at$_1$auat$_2$,

$\Rightarrow$ M ⊨ ww'aw$_2$=w$_1$at$_1$auat$_2$aw$_2$=x= wat$_3$avat$_4$,

$\Rightarrow$ by (3.7),  $M \vDash w'aw_2 = at_3avat_4$,

$\Rightarrow$  $M \vDash w' = a \lor aBw'$,

$\Rightarrow$ from $M \vDash Tally_b(t_2)$, $M \vDash \neg(w'=a)$,

$\Rightarrow$  $M \vDash aBw'$,

$\Rightarrow$  $M \vDash \exists w''\ aw''=w'$,

$\Rightarrow$  $M \vDash (aw'')aw_2 = at_3avat_4$,

$\Rightarrow$  $M \vDash w''aw_2 = t_3avat_4$,

$\Rightarrow$ by (4.17$^b$), $M \vDash t_3 \subseteq_p w_1 \lor t_3 \subseteq_p t_1auat_2$.

Now, if  $M \vDash t_3 \subseteq_p w_1$, then from  $M \vDash Max^+T_b(t_1,w_1)$  we have  $M \vDash t_3 < t_1$.

But from $M \vDash MaxT_b(t_4,x)$ we have  $M \vDash t_1 \leq t_4 = t_3$. But then $M \vDash t_3 < t_1 \leq t_3$,

contradicting  $M \vDash t_3 \in I \subseteq I_0$.

Therefore,  $M \vDash t_3 \subseteq_p t_1auat_2$.

$\Rightarrow$ by (4.17$^b$),  $M \vDash t_3 \subseteq_p t_1 \lor t_3 \subseteq_p uat_2$.

If $M \vDash t_3 \subseteq_p t_1$, then from  $M \vDash MaxT_b(t_4,x)$ we have  $M \vDash t_1 \leq t_4 = t_3 \leq t_1$,

so  $M \vDash t_1 = t_3$. But then from hypothesis (1b)  we have  $M \vDash t_1 < t_2$, whence

$M \vDash t_1 < t_2 \leq t_4 = t_3 \leq t_1$,  again contradicting  $M \vDash t_1 \in I \subseteq I_0$.

Therefore,  $M \vDash t_3 \subseteq_p uat_2$.

$\Rightarrow$ by (4.17$^b$),  $M \vDash t_3 \subseteq_p u \lor t_3 \subseteq_p t_2$.

If  $M \vDash t_3 \subseteq_p u \subseteq_p aua$, then from $M \vDash Pref(aua,t_1)$  we have  $M \vDash Max^+T_b(t_1,aua)$,

so $M \vDash t_3 < t_1$ and a contradiction follows just as above.

Finally, if $M \vDash t_3 \subseteq_p t_2$, then from  $M \vDash MaxT_b(t_4,x)$  we have  $M \vDash t_2 \leq t_4 = t_3 \leq t_2$,

so  $M \vDash t_2 = t_3$.

$\Rightarrow$ $M \vDash w_1at_1auat_3aw_2=x= wat_3avat_4$,

$\Rightarrow$ $M \vDash (w_1at_1aua)Bx$ & $(wa)Bx$,

$\Rightarrow$ by (3.8), $M \vDash (w_1at_1aua)B(wa)$ v $w_1at_1aua=wa$ v $(wa)B(w_1at_1aua)$.

(1biii1) $M \vDash (w_1at_1aua)B(wa)$.

$\Rightarrow$ $M \vDash \exists w_3\ w_1at_1auaw_3=wa$,

$\Rightarrow$ $M \vDash w_1at_1auat_2aw_2=x=(w_1at_1auaw_3)t_3avat_4$,

$\Rightarrow$ by (3.7), $M \vDash t_2aw_2=w_3t_3avat_4$,

$\Rightarrow$ from $M \vDash t_2=t_3$, $M \vDash t_2aw_2=w_3t_2avat_4$.

We then derive a contradiction exactly as in (1ai2a) with $w_2$ in place of $x_1$.

(1biii2) $M \vDash w_1at_1aua=wa$.

$\Rightarrow$ from $M \vDash t_2=t_3$, $M \vDash w_1at_1auat_2aw_2=x=(w_1at_1aua)t_2avat_4$ & $t_2=t_4$,

where $M \vDash Max^+T_b(t_1,w_1)$, as required.

(1biii3) $M \vDash (wa)B(w_1at_1aua)$.

$\Rightarrow$ $M \vDash \exists w_3\ waw_3=w_1at_1aua=w_1at_1aua$,

$\Rightarrow$ $M \vDash (waw_3)t_2aw_2=w_1at_1auat_2aw_2=x=wat_3avat_4$,

$\Rightarrow$ by (3.7), $M \vDash w_3t_2aw_2=t_3avat_4$,

$\Rightarrow$ from $M \vDash t_2=t_3$, $M \vDash w_3t_2aw_2=t_2avat_4$.

We then derive a contradiction exactly as in (1ai2b) with $w_2$ in place of $x_1$.

This completes the subcase (1b).

(1c) $M \vDash Lastf(x,t_1,aua,t_2)$.

$\Rightarrow$ from hypothesis (1), by (5.15), $M \vDash u=v$.

By (9.4) this contradicts the hypothesis $M \vDash u<_x v$.

(2)  $M \vDash \exists w_3 Intf(x,w_3,t_3,ava,t_4)$.

$\Rightarrow M \vDash Pref(ava,t_3) \, \& \, Tally_b(t_4) \, \& \, t_3<t_4 \, \&$

$\& \, \exists w_4 \, x=w_3at_3avat_4aw_4 \, \& \, Max^+T_b(t_3,w_3)$.

(2a)  $M \vDash Firstf(x,t_1,aua,t_2)$.

By the same argument as in (1), $M \vDash \neg(x=t_1auat_2)$.

$\Rightarrow M \vDash Pref(aua,t_1) \, \& \, Tally_b(t_2) \, \& \, (t_1<t_2 \, \& \, (t_1auat_2a)Bx))$,

$\Rightarrow M \vDash \exists x_1 \, x=t_1auat_2ax_1$,

$\Rightarrow$ by (5.34), $M \vDash t_2 \leq t_3$,

$\Rightarrow M \vDash t_1<t_2 \leq t_3$,

$\Rightarrow M \vDash (t_1auat_2a)Bx \, \& \, w_3Bx$,

$\Rightarrow$ by (3.8), $M \vDash w_3B(t_1auat_2) \vee w_3=t_1auat_2 \vee (t_1auat_2)Bw_3$.

(2ai)  $M \vDash w_3B(t_1auat_2)$.

$\Rightarrow M \vDash \exists w' \, w_3w'=t_1auat_2$.

We distinguish the subcases:

(2ai1) $M \vDash w'=a$    and    (2ai2) $M \vDash aBw'$,

and proceed as in (1ai1) and (1ai2) replacing $wat_3avat_4$ with $w_3at_3avat_4aw_4$.

We obtain $M \vDash t_2=t_3$.

We then have $M \vDash (t_1aua)B(w_3a) \vee t_1aua=w_3a \vee (w_3a)B(t_1aua)$.

(2ai2a)  $M \vDash (t_1aua)B(w_3a)$.

Exactly analogous to (1ai2a).

(2ai2b)  $M \vDash (w_3a)B(t_1aua)$.

Arguing analogously to (1ai2b) we obtain $M \vDash \exists w_6 \, t_2avat_4aw_4=(w_6a)t_2ax_1$,

where $M \vDash w_6a=w_5$ & $waw_5=t_1aua$.

$\Rightarrow$ by (4.14$^b$), $M \vDash w_6=t_2$ v $t_2Bw_6$,

$\Rightarrow M \vDash t_2 \subseteq_p w_6$,

$\Rightarrow$ from $M \vDash wa(w_6a)=t_1aua$, $M \vDash waw_6=t_1au$,

$\Rightarrow M \vDash t_2 \subseteq_p w_6 \subseteq_p t_1au$,

$\Rightarrow$ by (4.17$^b$), $M \vDash t_2 \subseteq_p t_1$ v $t_2 \subseteq_p u$.

If $M \vDash t_2 \subseteq_p u$, then from $M \vDash Pref(aua,t_1)$ we have $M \vDash t_2 < t_1$. Hence, either way $M \vDash t_2 \leq t_1$, and so we have $M \vDash t_1 < t_2 \leq t_1$, contradicting $M \vDash t_1 \in I \subseteq I_0$.

(2ai2c) $M \vDash w_3a=t_1aua$.

$\Rightarrow M \vDash x=t_1auat_2avat_4$, as required.

(2aii) $M \vDash w_3=t_1auat_2$.

$\Rightarrow M \vDash t_1auat_2ax_1=x=t_1auat_2at_3avat_4aw_4$.

We now proceed as in (1aii): from hypothesis $M \vDash MinSet(x)$, we have

$M \vDash \exists v',t',t''\ Occ(t_1auat_2,a,t_3avat_4aw_4,x,t',v',t'')$.

$\Rightarrow M \vDash Fr(x,t',v',t'')$,

$\Rightarrow M \vDash \exists v''\ v'=av''a$.

(2aii1) $M \vDash Firstf(x,t',av''a,t'')$.

Exactly the same argument as in (1aii1).

(2aii2) $M \vDash \exists w'Intf(x,w',t',av''a,t'')$.

$\Rightarrow M \vDash Pref(av''a,t')$ & $Tally_b(t'')$ & $t'<t''$ &

& $\exists w''\ x=w'at'av''at''aw''$ & $Max^+T_b(t',w')$.

(2aii2a) $M \vDash w'at'=t_1auat_2$.

$\Rightarrow$ by (4.24$^b$), $M \vDash t'=t_2$,

$\Rightarrow M \vDash t_1 < t'$,

$\Rightarrow$ from $M \vDash \text{Firstf}(x,t_1,aua,t_2)$, $M \vDash u <_x v''$,

$\Rightarrow$ from the principal hypothesis, $M \vDash \neg(v'' <_x v)$,

$\Rightarrow$ by (9.7), $M \vDash v=v'' \lor v <_x v''$.

    (2aii2ai) $M \vDash v=v''$.

$\Rightarrow$ from hypothesis (2aii2a),

       $M \vDash (t_1auat_2)av''at''aw''=(w'at')av''at''aw''=x=(w_3at_3)avat_4aw_4$,

$\Rightarrow$ from $M \vDash \text{Env}(t,x)$ & $\text{Fr}(x,t_3,ava,t_4)$ & $\text{Fr}(x,t',av''a,t'')$ & $v=v''$, $M \vDash t_3=t'=t_2$,

$\Rightarrow$ from $M \vDash t_3 < t_4$, $M \vDash t_2 < t_4$,

$\Rightarrow$ from $M \vDash \text{Pref}(aua,t_1)$ & $t_1 < t_2$, $M \vDash \text{Max}^+T_b(t_2,t_1au)$,

$\Rightarrow$ from $M \vDash \text{Max}^+T_b(t_3,w_3)$ & $t_2=t_3$, $M \vDash \text{Max}^+T_b(t_2,w_3)$ & $\text{Max}^+T_b(t_3,t_1au)$,

$\Rightarrow$ by (4.20), $M \vDash t_1au=w_3$,

$\Rightarrow M \vDash (t_1au)at_2ava=w_3at_3ava$,

$\Rightarrow$ by (3.7), $M \vDash t''aw''=t_4aw_4$,

$\Rightarrow$ by (4.23$^b$), $M \vDash t''=t_4$,

$\Rightarrow M \vDash x=t_1auat_2avat_4aw_4$,

$\Rightarrow M \vDash \exists w_2\; x=t_1auat_2avat_4aw_2$, as required.

    (2aii2aii) $M \vDash v <_x v''$.

$\Rightarrow$ from $M \vDash \text{Intf}(x,w_3,t_3,ava,t_4)$ & $\text{Intf}(x,w',t',av''a,t'')$, $M \vDash t_4 \leq t'$,

$\Rightarrow M \vDash t_3 < t_4 \leq t'$,

$\Rightarrow M \vDash t_2 \leq t_3 < t'=t_2$, contradicting $M \vDash t_2 \in I \subseteq I_0$.

(2aii2b) $M \vDash \exists u''(u''B(av''a) \mathbin{\&} w'at'u''=t_1auat_2a)$.

$\implies M \vDash \exists v_3\ u''v_3=av''a$,

$\implies M \vDash w'at'av''at''aw''=x=w'at'(u''v_3)t''aw''=(t_1auat_2a)v_3t''aw''=$

$$=t_1auat_2ax_1=t_1auat_2at_3avat_4aw_4,$$

$\implies$ by (3.7), $M \vDash v_3t''aw''=t_3avat_4aw_4$.

We then derive $M \vDash t_3\subseteq_p av''a$ as in (1bi2b), and further that $M \vDash t_3<t'$.

On the other hand, we also have from $M \vDash w'at'u''=t_1auat_2a$ that

$$M \vDash u''=a \vee aEu''.$$

If $M \vDash u''=a$, then $M \vDash w'at'a=t_1auat_2a$, hence $M \vDash w'at'=t_1auat_2$, and we proceed to reason as in (2aii2a).

If $M \vDash aEu''$, then $M \vDash \exists u_3\ u''=u_3a$.

$\implies M \vDash w'at'(u_3a)=t_1auat_2a$,

$\implies M \vDash w'at'u_3=t_1auat_2$.

Assume that $M \vDash \mathrm{Tally}_b(u_3)$.

$\implies$ by (4.5), $M \vDash \mathrm{Tally}_b(t'u_3)$,

$\implies$ from $M \vDash \mathrm{Tally}_b(t_2)$, by (4.24$^b$), $M \vDash t'u_3=t_2$,

$\implies M \vDash w'at_2=t_1auat_2$,

$\implies$ from $M \vDash w'at'u''=t_1auat_2a$, $M \vDash w'at'u''=w'at_2a$,

$\implies$ by (3.7), $M \vDash t'u''=t_2a$,

$\implies$ from $M \vDash \mathrm{Tally}_b(t')$, by (4.14$^b$), $M \vDash t_2=t' \vee t'Bt_2$,

$\implies M \vDash t'\subseteq_p t_2$,

$\implies M \vDash t'\leq t_2\leq t_3<t'$, contradicting $M \vDash t'\in I\subseteq I_0$.

Assume that $M \vDash \neg Tally_b(u_3)$.

$\Rightarrow M \vDash a \subseteq_p u_3$,

$\Rightarrow M \vDash u_3 = a \lor aBu_3 \lor aEu_3 \lor \exists u_4, u_5 \; u_3 = u_4 a u_5$,

$\Rightarrow$ from $M \vDash Tally_b(t_2)$ & $t'u_3 = t_2$, $M \vDash \neg(u_3 = a \lor aEu_3)$,

$\Rightarrow M \vDash aEu_3 \lor \exists u_4, u_5 \; u_3 = u_4 a u_5$,

$\Rightarrow M \vDash \exists u_6 \; w'at'(au_6) = t_1 auat_2 \lor \exists u_4, u_5 \; w'at'(u_4 au_5) = t_1 auat_2$,

$\Rightarrow$ either way, by (4.16), $M \vDash t' \subseteq_p t'u_4 \subseteq_p u \subseteq_p aua$,

$\Rightarrow$ from $M \vDash Pref(aua, t_1)$, $M \vDash Max^+T_b(t_1, aua)$,

$\Rightarrow M \vDash t' < t_1 < t_3 < t'$, contradicting $M \vDash t' \in I \subseteq I_0$.

    (2aii2c) $M \vDash w'at'av''a = t_1 auat_2 a$.

$\Rightarrow M \vDash w'at'av''at''aw'' = x = t_1 auat_2 at_3 avat_4 aw_4 = w'at'av''at_3 avat_4 aw_4$,

$\Rightarrow$ by (3.7), $M \vDash t''aw'' = t_3 avat_4 aw_4$,

$\Rightarrow$ from $M \vDash Tally_b(t'')$ & $Tally_b(t_3)$, by (4.24$^b$), $M \vDash t'' = t_3$,

$\Rightarrow$ from $M \vDash Intf(x, w', t', av''a, t'')$ & $Intf(x, w_3, t_3, ava, t_4)$, $M \vDash v'' <_x v$,

$\Rightarrow$ from $M \vDash Firstf(x, t_1, aua, t_2)$, as in (2aii2a), $M \vDash u <_x v''$.

But this contradicts the principal hypothesis.

    (2aii3) $M \vDash Lastf(x, t', av''a, t'')$.

$\Rightarrow M \vDash Pref(av''a, t')$ & $Tally_b(t'')$ & $t' = t'' = t$ &

                       & $\exists w' (x = w'at'av''at''$ & $Max^+T_b(t', w'))$.

As above, from $M \vDash Occ(t_1 auat_2, a, t_3 avat_4 aw_4, x, t', v', t'')$ we have three scenarios:

    (2aii3a) $M \vDash w'at' = t_1 auat_2$.

Then we have, as in (2aii2a), that $M \vDash t'=t_2$.

But we also have $M \vDash t_1<t_2\leq t_3<t_4$, whereas from $M \vDash Env(t,x)$, $M \vDash t_4\leq t=t'$.

But then $M \vDash t_2<t_4=t'=t_2$, contradicting $M \vDash t_2 \in I \subseteq I_0$.

    (2aii3b) $M \vDash \exists u''(u''B(av''a)$ & $w'at'u''=t_1auat_2a)$.

$\Rightarrow M \vDash \exists v_3\ u''v_3=av''a$,

$\Rightarrow M \vDash w'at'av''at''=x=w'at'(u''v_3)t''=(t_1auat_2a)v_3t''=t_1auat_2at_3avat_4aw_4$,

$\Rightarrow$ by (3.7), $M \vDash v_3t''=t_3avat_4aw_4$,

$\Rightarrow$ as in (1bi2b), $M \vDash t_3\subseteq_p av''a$ and further that $M \vDash t_3<t'$.

Then the rest of the argument proceeds exactly as in (2aii2b).

    (2aii3c) $M \vDash w'at'av''a=t_1auat_2a$.

$\Rightarrow M \vDash w'at'av''at''=x=t_1auat_2at_3avat_4aw_4=w'at'av''at_3avat_4aw_4$,

$\Rightarrow$ by (3.7), $M \vDash t''=t_3avat_4aw_4$, a contradiction because $M \vDash Tally_b(t'')$.

  (2aiii) $M \vDash (t_1auat_2)Bw_3$.

$\Rightarrow M \vDash \exists w_5\ t_1auat_2w_5=w_3$,

$\Rightarrow M \vDash w_3at_3avat_4aw_4=x=(t_1auat_2w_5)at_3avat_4aw_4=t_1auat_2ax_1$,

$\Rightarrow$ from $M \vDash MinSet(x)$, we have

    $M \vDash \exists v',t',t''\ Occ(t_1auat_2w_5,a,t_3avat_4aw_4,x,t',v',t'')$.

$\Rightarrow M \vDash Fr(x,t',v',t'')$,

$\Rightarrow M \vDash \exists v''\ v'=av''a$.

  (2aiii1) $M \vDash Firstf(x,t',av''a,t'')$.

Same argument as in (1aiii1), with $w_5$ in place of w'.

  (2aiii2) $M \vDash \exists w'Intf(x,w',t',av''a,t'')$.

$\Rightarrow$ M ⊨ Pref(av''a,t') & Tally$_b$(t'') & t'<t'' &

& ∃w'' x=w'at'av''at''aw'' & Max$^+$T$_b$(t',w').

(2aiii2a)  M ⊨ w'at'=t$_1$auat$_2$w$_5$.

$\Rightarrow$ M ⊨ w'at'=w$_3$,

$\Rightarrow$ M ⊨ t'⊆$_p$w$_3$,

$\Rightarrow$ from  M ⊨ Max$^+$T$_b$(t$_3$,w$_3$),  M ⊨ t'<t$_3$.

But also from  M ⊨ w'at'=w$_3$,  we have

M ⊨ w'at'av''at''aw''=(w'at')at$_3$avat$_4$aw$_4$,

$\Rightarrow$ by (3.7),  M ⊨ v''at''aw''=t$_3$avat$_4$aw$_4$,

$\Rightarrow$ by (4.14$^b$),  M ⊨ v''=t$_3$ v t$_3$Bv'',

$\Rightarrow$ M ⊨ t$_3$⊆$_p$v'',

$\Rightarrow$ from  M ⊨ Pref(av''a,t'),  M ⊨ Max$^+$T$_b$(t',v''),

$\Rightarrow$ M ⊨ t$_3$<t'<t$_3$,  contradicting  M ⊨ t$_3$∈I⊆I$_0$.

(2aiii2b)  M ⊨ ∃u''(u''B(av''a) & w'at'u''=t$_1$auat$_2$w$_5$a).

$\Rightarrow$ M ⊨ ∃v$_3$ u''v$_3$=av''a,

$\Rightarrow$ M ⊨ w'at'av''at''aw''=x=w'at'(u''v$_3$)t''aw''=(t$_1$auat$_2$w$_5$a)v$_3$t''aw''=

=(t$_1$auat$_2$w$_5$)at$_3$avat$_4$aw$_4$,

$\Rightarrow$ by (3.6),  M ⊨ v$_3$t''aw''=t$_3$avat$_4$aw$_4$.

We then derive  M ⊨ t$_3$⊆$_p$v''  as in (1bi2b), and further that  M ⊨ t$_3$<t'.

But we also have from  M ⊨ w'at'u''=t$_1$auat$_2$w$_5$a  that  M ⊨ w'at'u''=w$_3$a.

$\Rightarrow$ M ⊨ u''=a v aEu'',

$\Rightarrow$ M ⊨ w'at'a=w$_3$a v ∃u$_3$ w'at'(u$_3$a)=w$_3$a,

$\Rightarrow M \vDash w'at'=w_3 \lor w'at'u_3=w_3, \Rightarrow M \vDash t' \subseteq_p w_3,$

$\Rightarrow$ from $M \vDash \text{Max}^+T_b(t_3,w_3)$, $M \vDash t'<t_3,$

$\Rightarrow M \vDash t_3<t'<t_3,$ contradicting $M \vDash t_3 \in I \subseteq I_0.$

(2aiii2c) $M \vDash w'at'av''a=t_1auat_2w_5a.$

$\Rightarrow M \vDash w'at'av''at''aw''=x=(t_1auat_2w_5)at_3avat_4aw_4=w'at'av''at_3avat_4aw_4,$

$\Rightarrow$ by (3.7), $M \vDash t''aw''=t_3avat_4aw_4.$

We then proceed to derive a contradiction exactly as in (2aii2c).

(2aiii3) $M \vDash \text{Lastf}(x,t',av''a,t'').$

$\Rightarrow M \vDash \text{Pref}(av''a,t') \,\&\, \text{Tally}_b(t'') \,\&\, t'=t''=t \,\&$

$\&\, \exists w'(\, x=w'at'av''at''aw'' \,\&\, \text{Max}^+T_b(t',w')).$

(2aiii3a) $M \vDash w'at'=t_1auat_2w_5.$

Exactly the same as (2aiii2a), omitting $aw''$ throughout.

(2aiii3b) $M \vDash \exists u''(u''B(av''a) \,\&\, w'at'u''=t_1auat_2w_5a).$

We proceed as in (2aii2b) and derive $M \vDash v_3t''=t_3avat_4aw_4.$

We then derive $M \vDash t_3 \subseteq_p v''$ as in (1bi2b), and further that $M \vDash t_3<t'.$ We then obtain a contradiction as in (2aii2b).

(2aiii2c) $M \vDash w'at'av''a=t_1auat_2w_5a.$

$\Rightarrow M \vDash w'at'av''at''aw''=x=t_1auat_2w_5at_3avat_4aw_4=w'at'av''at_3avat_4aw_4,$

$\Rightarrow$ by (3.7), $M \vDash t''=t_3avat_4aw_4,$ a contradiction because $M \vDash \text{Tally}_b(t'').$

(2b) $M \vDash \exists w_1 \, \text{Intf}(x,w_1,t_1,aua,t_2).$

$\Rightarrow M \vDash \text{Pref}(aua,t_1) \,\&\, \text{Tally}_b(t_2) \,\&\, t_1<t_2 \,\&$

$\&\, \exists w_2 \, x=w_1at_1auat_2aw_2 \,\&\, \text{Max}^+T_b(t_1,w_1).$

$\Rightarrow$ from hypothesis (2) and $M \vDash u<_x v$, $M \vDash t_2 \leq t_3$,

$\Rightarrow$ by (9.4), $M \vDash u \neq v$.

(2bi) $M \vDash t_2<t_3$.

$\Rightarrow M \vDash t_1<t_2<t_3<t_4$,

$\Rightarrow$ since $M \vDash t_1 \in I \subseteq I_0$, $M \vDash t_1 \neq t_4$,

$\Rightarrow$ by (5.27), $M \vDash t_1 a u a t_2 \subseteq_p w_3 \lor t_3 a v a t_4 \subseteq_p w_1$.

If $M \vDash t_3 \subseteq_p t_3 a v a t_4 \subseteq_p w_1$, then from $M \vDash \text{Max}^+ T_b(t_1, w_1)$ we have $M \vDash t_3<t_1$. But then $M \vDash t_1<t_3<t_1$, contradicting $M \vDash t_1 \in I \subseteq I_0$.

Hence $M \vDash \neg(t_3 a v a t_4 \subseteq_p w_1)$.

$\Rightarrow M \vDash t_1 a u a t_2 \subseteq_p w_3$.

From the proof of (5.27), case (2iib), we in fact have

$M \vDash w_1 a t_1 a u a t_2 = w_3 \lor w_1 a t_1 a u a t_2 a = w_3 \lor \exists w_5\, w_1 a t_1 a u a t_2 a w_5 = w_3$.

(2bi1) $M \vDash w_1 a t_1 a u a t_2 = w_3$.

$\Rightarrow M \vDash w_1 a t_1 a u a t_2 a w_2 = x = (w_1 a t_1 a u a t_2) a t_3 a v a t_4 a w_4$.

From $M \vDash \text{MinSet}(x)$ we have

$M \vDash \exists v', t', t''\, \text{Occ}(w_1 a t_1 a u a t_2, a, t_3 a v a t_4 a w_4, x, t', v', t'')$.

(2bi1a) $M \vDash \text{Firstf}(x, t', v', t'')$.

$\Rightarrow M \vDash \text{Pref}(v', t')\ \&\ \text{Tally}_b(t'')\ \&\ (t'<t''\ \&\ (t'v't''a)Bx)$,

$\Rightarrow M \vDash \exists v''\, v' = a v'' a$.

From $M \vDash \text{Occ}(w_1 a t_1 a u a t_2, a, t_3 a v a t_4 a w_4, x, t', v', t'')$, we distinguish three scenarios:

(2bi1ai) $M \vDash t' = w_1 a t_1 a u a t_2$.

This is ruled out because $M \vDash \text{Tally}_b(t')$.

   (2bi1aii) $M \vDash (w_1at_1auat_2a)B(t'av''a)$.

$\Rightarrow M \vDash \exists w_5\ w_1at_1auat_2aw_5=t'av''a$,

$\Rightarrow$ by (4.14$^b$), $M \vDash w_1=t' \lor t'Bw_1$,

$\Rightarrow M \vDash t'at_1auat_2aw_5=t'av''a \lor \exists w_6\ (t'w_6)at_1auat_2aw_5=t'av''a$,

$\Rightarrow$ by (3.7), $M \vDash t_1auat_2aw_5=v''a \lor w_6at_1auat_2aw_5=av''a$,

$\Rightarrow M \vDash t_1\subseteq_p av''a$,

$\Rightarrow$ from $M \vDash \text{Pref}(av''a,t')$, $M \vDash \text{Max}^+T_b(t',av''a)$,

$\Rightarrow M \vDash t_1<t'$,

$\Rightarrow$ from $M \vDash \text{Intf}(x,w_1,t_1,aua,t_2)$ by (5.19), $M \vDash \neg\text{Firstf}(x,t_1,aua,t_2)$,

$\Rightarrow$ by (5.20), $M \vDash t'<t_1$,

$\Rightarrow M \vDash t_1<t'<t_1$, contradicting $M \vDash t_1\in I\subseteq I_0$.

   (2bi1aiii) $M \vDash w_1at_1auat_2a=t'av''a$.

Same argument as in (2bi1aii).

   (2bi1b) $M \vDash \exists w'\text{Intf}(x,w',t',av''a,t'')$.

$\Rightarrow M \vDash \text{Pref}(av''a,t') \& \text{Tally}_b(t'') \& t'<t'' \&$

$\qquad\qquad\qquad\qquad \& \exists w''\ x=w'at'av''at''aw'' \& \text{Max}^+T_b(t',w')$.

Again, from $M \vDash \text{Occ}(w_1at_1auat_2,a,t_3avat_4aw_4,x,t',v',t'')$, we distinguish three scenarios:

   (2bi1bi) $M \vDash w'at'=w_1at_1auat_2$.

$\Rightarrow$ by (2bi1), $M \vDash w'at'=w_3$,

$\Rightarrow M \vDash t'\subseteq_p w_3$,

$\Rightarrow$ from $M \vDash \text{Max}^+T_b(t_3,w_3)$, $M \vDash t'<t_3$,

$\Rightarrow M \vDash w'at'av''at''aw''=x=w_3at_3avat_4aw_4=(w'at')at_3avat_4aw_4$,

$\Rightarrow$ by (3.7), $M \vDash v''at''aw''=t_3avat_4aw_4$,

$\Rightarrow$ by (4.14$^b$), $M \vDash v''=t_3 \lor t_3Bv''$,

$\Rightarrow M \vDash t_3 \subseteq_p v''$,

$\Rightarrow$ from $M \vDash \text{Pref}(av''a,t')$, $M \vDash \text{Max}^+T_b(t',av''a)$,

$\Rightarrow M \vDash t_3<t'$,

$\Rightarrow M \vDash t'<t_3<t'$, contradicting $M \vDash t' \in I \subseteq I_0$.

(2bi1bii) $M \vDash \exists u''(u''B(av''a) \& w'at'u''=w_1at_1auat_2a)$.

$\Rightarrow$ by (2bi1), $M \vDash w'at'u''=w_3a$,

$\Rightarrow M \vDash u''=a \lor aEu''$,

$\Rightarrow M \vDash w'at'a=w_3a \lor \exists w_5\ w'at'(w_5a)=w_3a$,

$\Rightarrow M \vDash w'at'=w_3 \lor w'at'w_5=w_3$,

$\Rightarrow M \vDash t' \subseteq_p w_3$,

$\Rightarrow$ from $M \vDash \text{Max}^+T_b(t_3,w_3)$, $M \vDash t'<t_3$.

On the other hand, $M \vDash \exists v_3\ u''v_3=av''a$.

$\Rightarrow M \vDash w'at'(u''v_3)t''aw''=x=w_3at_3avat_4aw_4=(w'at'u'')t_3avat_4aw_4$,

$\Rightarrow$ by (3.7), $M \vDash v_3t''aw''=t_3avat_4aw_4$.

We then derive as in (1bi2b) that $M \vDash t_3 \subseteq_p v''$, and further that $M \vDash t_3<t'$.

But then $M \vDash t_3<t'<t_3$, contradicting $M \vDash t_3 \in I \subseteq I_0$.

(2bi1biii) $M \vDash w'at'av''a=w_1at_1auat_2a$.

$\Rightarrow M \vDash w'at'av''=w_1at_1auat_2$,

$\Rightarrow$ by (4.15$^b$), $M \vDash v''=t_2 \vee t_2 E v''$,

$\Rightarrow M \vDash t_2 \subseteq_p v''$,

$\Rightarrow$ from $M \vDash \text{Pref}(av''a,t')$, $M \vDash \text{Max}^+ T_b(t',av''a)$,

$\Rightarrow M \vDash t_2 < t'$,

$\Rightarrow$ from $M \vDash \text{Intf}(x,w_1,t_1,aua,t_2)$ & $\text{Intf}(x,w',t',av''a,t'')$, $M \vDash u <_x v''$.

On the other hand, we also have $M \vDash w'at'u''=w_3 a$, whence

$\quad M \vDash w'at'av''at''aw''=x=w_3 at_3 avat_4 aw_4=(w'at'av''a)t_3 avat_4 aw_4$,

$\Rightarrow$ by (3.7), $M \vDash t''aw''=t_3 avat_4 aw_4$,

$\Rightarrow$ by (4.23$^b$), $M \vDash t''=t_3$,

$\Rightarrow$ from $M \vDash \text{Intf}(x,w',t',av''a,t'')$ & $\text{Intf}(x,w_3,t_3,ava,t_4)$, $M \vDash v'' <_x v$.

But this along with $M \vDash u <_x v''$ contradicts the principal hypothesis.

$\quad$ (2bi1c) $M \vDash \text{Lastf}(x,t',av''a,t'')$.

$\Rightarrow M \vDash \text{Pref}(av''a,t')$ & $\text{Tally}_b(t'')$ & $t'=t''=t$ &

$\quad\quad\quad\quad\quad$ & $\exists w'(x=w'at'av''at''$ & $\text{Max}^+ T_b(t',w'))$.

From $M \vDash \text{Occ}(w_1 at_1 auat_2, a, t_3 avat_4 aw_4, x, t', v', t'')$, we have:

$\quad$ (2bi1ci) $M \vDash w'at'=w_1 at_1 auat_2$.

$\Rightarrow M \vDash w'at'=w_3$,

$\Rightarrow M \vDash t' \subseteq_p w_3$,

$\Rightarrow$ from $M \vDash \text{Max}^+ T_b(t_3,w_3)$, $M \vDash t' < t_3$.

From $M \vDash \text{Env}(t,x)$, $M \vDash t_3 < t_4 \leq t = t'$.

But then $M \vDash t_3 < t' < t_3$, contradicting $M \vDash t_3 \in I \subseteq I_0$.

(2bi1cii)  $M \vDash \exists u''(u''B(av''a) \& w'at'u''=w_1at_1auat_2a)$.

$\Longrightarrow M \vDash w'at'u''=w_3a$.

We then derive, as in (2bi1bii), that $M \vDash t' \subseteq_p w_3$, whence a contradiction follows as in (2bi1ci).

(2bi1ciii)  $M \vDash w'at'av''a=w_1at_1auat_2a$.

$\Longrightarrow M \vDash w'at'av''=w_1at_1auat_2=w_3$,

$\Longrightarrow M \vDash t' \subseteq_p w_3$, and we obtain a contradiction just as in (2bi1ci).

(2bi2)  $M \vDash w_1at_1auat_2a=w_3$.

$\Longrightarrow M \vDash w_1at_1auat_2aw_2=x=(w_1at_1auat_2a)at_3avat_4aw_4$.

From  $M \vDash MinSet(x)$,    $M \vDash \exists v',t',t''\ Occ(w_1at_1auat_2a,a,t_3avat_4aw_4,x,t',v',t'')$.

(2bi2a)  $M \vDash Firstf(x,t',v',t'')$.

$\Longrightarrow M \vDash Pref(v',t') \& Tally_b(t'') \& (t'<t'' \& (t'v't''a)Bx)$,

$\Longrightarrow M \vDash \exists v''\ v'=av''a$.

From  $M \vDash Occ(w_1at_1auat_2a,a,t_3avat_4aw_4,x,t',v',t'')$, there are three scenarios:

(2bi2ai)  $M \vDash t'=w_1at_1auat_2$.

A contradiction because $M \vDash Tally_b(t')$.

(2bi2aii)  $M \vDash (w_1at_1auat_2aa)B(t'av''a)$.

The same argument applies as in (2bi1aii), with $w_1at_1auat_2aa$ in place of $w_1at_1auat_2a$.

(2bi2aiii)  $M \vDash w_1at_1auat_2aa=t'av''a$.

The same argument as in (2bi2aii).

(2bi2b)  $M \vDash \exists w'\text{Intf}(x,w',t',av''a,t'')$.

$\Rightarrow M \vDash \text{Pref}(av''a,t') \,\&\, \text{Tally}_b(t'') \,\&\, t'<t'' \,\&$

$\&\, \exists w''\, x=w'at'av''at''aw'' \,\&\, \text{Max}^+T_b(t',w')$.

From $M \vDash \text{Occ}(w_1at_1auat_2a,a,t_3avat_4aw_4,x,t',v',t'')$, we have:

(2bi2bi)  $M \vDash w'at'=w_1at_1auat_2a$.

A contradiction because $M \vDash \text{Tally}_b(t')$.

(2bi2bii)  $M \vDash \exists u''(u''B(av''a) \,\&\, w'at'u''=w_1at_1auat_2aa)$.

Exactly as in (2bi1bii).

(2bi2biii)  $M \vDash w'at'av''a=w_1at_1auat_2aa$.

$\Rightarrow M \vDash w'at'av''=w_1at_1auat_2a$,

$\Rightarrow M \vDash v''=a \,\vee\, aEv''$,

$\Rightarrow M \vDash w'at'aa=w_1at_1auat_2a \,\vee\, \exists v_3\, w'at'a(v_3a)=w_1at_1auat_2a$,

$\Rightarrow M \vDash w'at'a=w_1at_1auat_2 \,\vee\, w'at'av_3=w_1at_1auat_2$.

Now, $M \vDash w'at'a \neq w_1at_1auat_2$ because $M \vDash \text{Tally}_b(t_2)$.

$\Rightarrow M \vDash w'at'av_3=w_1at_1auat_2$,

$\Rightarrow$ by (4.15$^b$),  $M \vDash v_3=t_2 \,\vee\, t_2Ev_3$,

$\Rightarrow M \vDash t_2 \subseteq_p v_3 \subseteq_p v''$,

$\Rightarrow$ from  $M \vDash \text{Pref}(av''a,t')$, $M \vDash \text{Max}^+T_b(t',av''a)$,

$\Rightarrow M \vDash t_2<t'$,

$\Rightarrow$ from  $M \vDash \text{Intf}(x,w_1,t_1,aua,t_2) \,\&\, \text{Intf}(x,w',t',av''a,t'')$, $M \vDash u<_x v''$.

We then argue, exactly as in (2bi1biii), that also $M \vDash v''<_x v$, which contradicts the principal hypothesis.

(2bi2c)  $M \vDash \text{Lastf}(x,t',av''a,t'')$.

$\Rightarrow M \vDash \text{Pref}(av''a,t') \& \text{Tally}_b(t'') \& t'=t''=t \&$

$\& \exists w'( x=w'at'av''at'' \& \text{Max}^+T_b(t',w'))$.

Again, from  $M \vDash \text{Occ}(w_1at_1auat_2a,a,t_3avat_4aw_4,x,t',v',t'')$, we have:

(2bi2ci)  $M \vDash w'at'=w_1at_1auat_2a$.

Same as (2bi2bi).

(2bi2cii)  $M \vDash \exists u''(u''B(av''a) \& w'at'u''=w_1at_1auat_2aa)$.

We derive, as in (2bi1bii), that  $M \vDash t'<t_3$.  By the same argument as in (1bi2b) we obtain $M \vDash t'<t_3<t'$,  a contradiction since  $M \vDash t' \in I \subseteq I_0$.

(2bi2ciii)  $M \vDash w'at'av''a=w_1at_1auat_2aa$.

Exactly as (2bi1ciii).

(2bi3)  $M \vDash \exists w_5\ w_1at_1auat_2aw_5=w_3$.

$\Rightarrow M \vDash w_1at_1auat_2aw_2=x=(w_1at_1auat_2aw_5)at_3avat_4aw_4$.

From  $M \vDash \text{MinSet}(x)$,   $M \vDash \exists v',t',t''\ \text{Occ}(w_1at_1auat_2aw_5,a,t_3avat_4aw_4,x,t',v',t'')$.

(2bi3a)  $M \vDash \text{Firstf}(x,t',v',t'')$.

$\Rightarrow M \vDash \text{Pref}(v',t') \& \text{Tally}_b(t'') \& (t'<t'' \& (t'v't''a)Bx)$,

$\Rightarrow M \vDash \exists v''\ v'=av''a$.

From  $M \vDash \text{Occ}(w_1at_1auat_2aw_5,a,t_3avat_4aw_4,x,t',v',t'')$, there are three scenarios:

(2bi3ai)  $M \vDash t'=w_1at_1auat_2aw_5$.

A contradiction because $M \vDash \text{Tally}_b(t')$.

(2bi3aii) $M \vDash (w_1at_1auat_2aw_5a)B(t'av''a)$.

The same argument applies as in (2bi1aii), with $w_1at_1auat_2aw_5a$ in place of $w_1at_1auat_2a$, with appropriate change of variables throughout.

(2bi3aiii) $M \vDash w_1at_1auat_2aw_5a = t'av''a$.

Same argument as in (2bi3aii).

(2bi3b) $M \vDash \exists w' \mathrm{Intf}(x,w',t',av''a,t'')$.

$\Rightarrow M \vDash \mathrm{Pref}(av''a,t') \,\&\, \mathrm{Tally}_b(t'') \,\&\, t'<t'' \,\&\,$

$\&\, \exists w''\, x = w'at'av''at''aw'' \,\&\, \mathrm{Max}^+T_b(t',w')$.

From $M \vDash \mathrm{Occ}(w_1at_1auat_2aw_5,a,t_3avat_4aw_4,x,t',v',t'')$, we have:

(2bi3bi) $M \vDash w'at' = w_1at_1auat_2aw_5$.

$\Rightarrow$ by (2bi3), $M \vDash w'at' = w_3a$.

We then proceed exactly as in (2bi1bi) to derive a contradiction.

(2bi3bii) $M \vDash \exists u''(u''B(av''a) \,\&\, w'at'u'' = w_1at_1auat_2aw_5a)$.

$\Rightarrow$ by (2bi3), $M \vDash w'at'u'' = w_3a$.

Exactly the same as (2bi1bii), with appropriate change of variables.

(2bi3biii) $M \vDash w'at'av''a = w_1at_1auat_2aw_5a$.

$\Rightarrow$ from (2bi3), $M \vDash w'at'av''a = w_3a$,

$\Rightarrow M \vDash w'at'av''at''aw'' = x = (w'at'av''a)t_3avat_4aw_4$,

$\Rightarrow$ by (3.7), $M \vDash t''aw'' = t_3avat_4aw_4$,

$\Rightarrow$ by (4.23[b]), $M \vDash t'' = t_3$,

$\Rightarrow$ from $M \vDash \mathrm{Intf}(x,w',t',av''a,t'') \,\&\, \mathrm{Intf}(x,w_3,t_3,ava,t_4)$, $M \vDash v'' <_x v$.

Now, we have that $M \vDash w'at'av''at''aw'' = x = w_1at_1auat_2aw_2$.

Suppose that $M \vDash v''=u$.

$\Rightarrow$ by (d) of $M \vDash Env(t,x)$, $M \vDash t'=t_1$,

$\Rightarrow$ from $M \vDash w'at'(av''at''aw'')=w_1at_1(auat_2aw_2)$ & $Max^+T_b(t',w')$ &

& $Max^+T_b(t_1,w_1)$, by (4.20), $M \vDash w'=w_1$,

$\Rightarrow$ from hypothesis (2bi3biii), by (3.7), $M \vDash v''=uat_2aw_5$,

$\Rightarrow$ from $M \vDash v''=u$, $M \vDash uBu$, contradicting $M \vDash u\in I\subseteq I_0$.

Therefore $M \vDash v''\neq u$.

Suppose that $M \vDash t''=t_1$.

$\Rightarrow$ from hypothesis (2bi), $M \vDash t''=t_1<t_2<t_3=t''$, contradicting $M \vDash t''\in I\subseteq I_0$.

Therefore, $M \vDash t''\neq t_1$.

Suppose that $M \vDash t_2=t'$.

Then from $M \vDash Intf(x,w_1,t_1,aua,t_2)$ & $Intf(x,w',t',av''a,t'')$, $M \vDash u<_x v''$, which along with $M \vDash v''<_x v$ contradicts the principal hypothesis.

Therefore, $M \vDash t_2 \neq t'$.

So we have $M \vDash v''\neq u$ & $t''\neq t_1$ & $t_2\neq t'$.

$\Rightarrow$ from $M \vDash Intf(x,w_1,t_1,aua,t_2)$ & $Intf(x,w',t',av''a,t'')$, by (5.27),

$M \vDash t_1auat_2\subseteq_p w'$ v $t'av''at''\subseteq_p w_1$.

If $M \vDash t_1auat_2\subseteq_p w'$, then $M \vDash t_2\subseteq_p t_1auat_2\subseteq_p w'$, hence $M \vDash t_2<t'$ from $M \vDash Max^+T_b(t',w')$. Then, just as above, $M \vDash u<_x v''$, contradicting the principal hypothesis.

If $M \vDash t'av''at''\subseteq_p w_1$, then from hypothesis (2bi3), $M \vDash t''\subseteq_p t'av''at''\subseteq_p w_1\subseteq_p w_3$. From $M \vDash Max^+T_b(t_3,w_3)$ we then have $M \vDash t''<t_3=t''$, contradicting

$M \vDash t'' \in I \subseteq I_0$.

 (2bi3c)  $M \vDash \text{Lastf}(x,t',av''a,t'')$.

$\Rightarrow M \vDash \text{Pref}(av''a,t')$ & $\text{Tally}_b(t'')$ & $t'=t''=t$ &

       & $\exists w'(x=w'at'av''at''$ & $\text{Max}^+T_b(t',w'))$.

From $M \vDash \text{Occ}(w_1at_1auat_2aw_5,a,t_3avat_4aw_4,x,t',v',t'')$, three subcases:

 (2bi3ci)  $M \vDash w'at'=w_1at_1auat_2aw_5$.

The same argument as in (2bi1ci).

 (2bi3cii)  $M \vDash \exists u''(u''B(av''a)$ & $w'at'u''=w_1at_1auat_2aw_5a)$.

The same argument as in (2bi1cii).

 (2bi3ciii)  $M \vDash w'at'av''a=w_1at_1auat_2aw_5a$.

The same argument as in (2bi1ci).

 (2bii)  $M \vDash t_2=t_3$.

$\Rightarrow$ from $M \vDash \text{Intf}(x,w_1,t_1,aua,t_2)$ & $\text{Intf}(x,w_3,t_3,ava,t_4)$, by (5.28),

    $M \vDash w_1at_1auat_2=w_3at_3$,

$\Rightarrow$ from $M \vDash w_1at_1auat_2aw_2=x=w_3at_3avat_4aw_4$,

 $M \vDash t_2<t_4$ & $\text{Pref}(ava,t_2)$ & $\text{Pref}(aua,t_1)$ & $x=w_1at_1auat_2avat_4aw_4$ &

       & $\text{Max}^+T_b(t_1,w_1)$,

as required.

(3)  $M \vDash \text{Firstf}(x,t_3,ava,t_4)$.

By (9.1) this contradicts the principal hypothesis $M \vDash u<_x v$.

This completes the proof of (10.3).

The following is a converse to (5.58).

(10.4) For any string concept $I \subseteq I_0$ there is a string concept $J \subseteq I$ such that

$QT^+ \vdash \forall x \in J \; \forall y,y'$ [MinSet(x) & $\forall w(w \; \varepsilon \; x \leftrightarrow w=y \lor w=y')$ & $y \neq y' \rightarrow$

$\rightarrow \exists t_1, t_2 (Tally_b(t_1)$ & $Tally_b(t_2)$ & $((y <_x y'$ & $x = t_1 a y a t_2 a y' a t_2) \lor$

$\lor (y' <_x y$ & $x = t_1 a y' a t_2 a y a t_2)))$].

Let $J \equiv I_{9.10}$ & $I_{10.3}$.

Assume $M \vDash MinSet(x)$ & $\forall w(w \; \varepsilon \; x \leftrightarrow w=y \lor w=y')$    where $M \vDash y \neq y'$ & $J(x)$.

$\Rightarrow M \vDash Set(x)$ & $y \; \varepsilon \; x$ & $y' \; \varepsilon \; x$,

$\Rightarrow$ by (9.7), $M \vDash y <_x y' \lor y' <_x y$.

Suppose that $M \vDash y <_x y'$.

We claim   (i) $M \vDash \neg \exists v(y <_x v$ & $v <_x y')$.

Otherwise, let $M \vDash y <_x v$ & $v <_x y'$.

$\Rightarrow$ by (9.4), $M \vDash y \neq v$ & $v \neq y'$, $M \vDash v \; \varepsilon \; x$ & $v \neq y$ & $v \neq y'$, contradicting the principal hypothesis. This proves (1).

$\Rightarrow$ from $M \vDash y \; \varepsilon \; x$ & $y' \; \varepsilon \; x$,   $M \vDash \exists t_1, t_2 \; Fr(x, t_1, aya, t_2)$ & $\exists t_3, t_4 \; Fr(x, t_3, ay'a, t_4)$.

We claim that   (ii) $M \vDash Firstf(x, t_1, aya, t_2)$ & $Lastf(x, t_3, ay'a, t_4)$.

From $M \vDash Set(x)$, by (5.18), $M \vDash \exists t \; Env(t,x)$.

Assume $M \vDash w \; \varepsilon \; x$.

$\Rightarrow M \vDash w=y \lor w=y'$,

$\Rightarrow$ from $M \vDash y <_x y'$, $M \vDash y \leq_x w$.

Hence $M \vDash \forall w(w \; \varepsilon \; x \rightarrow y \leq_x w)$.

$\Rightarrow$ by (9.10), $M \vDash Firstf(x, t_1, aya, t_2)$.

On the other hand, again assuming $M \vDash w \, \varepsilon \, x$, we have $M \vDash w=y \lor w=y'$,

whence from $M \vDash y<_x y'$ we obtain $M \vDash w \leq_x y'$.

Hence $M \vDash \forall w(w \, \varepsilon \, x \rightarrow w \leq_x y')$.

$\Rightarrow$ by (9.8), $M \vDash \text{Lastf}(x,t_3,ay'a,t_4) \,\&\, t_3=t_4=t$.

This proves (ii).

$\Rightarrow$ from $M \vDash \text{MinSet}(x) \,\&\, \text{Fr}(x,t_1,aya,t_2) \,\&\, \text{Fr}(x,t_3,ay'a,t_4) \,\&\, y<_x y' \,\&$

$\&\, \neg \exists v(y<_x v \,\&\, v<_x y')$

we have, by (10.3), that

$M \vDash (t_2=t_4 \,\&\, (x=t_1 ayat_2 ay'at_4 \lor$

$\lor \exists w_1(x=w_1 at_1 ayat_2 ay'at_4 \,\&\, \text{Max}^+T_b(t_1,w_1)))) \lor (t_2<t_4 \,\&\, \text{Pref}(ay'a,t_2) \,\&$

$\&\, \exists w_1,w_2((x=w_1 at_1 ayat_2 ay'at_4 aw_2 \,\&\, \text{Max}^+T_b(t_1,w_1)) \lor x=t_1 ayat_2 ay'at_4 aw_2))$.

Assume, for a reductio, that $M \vDash x=w_1 at_1 ayat_2 ay'at_4 \,\&\, \text{Max}^+T_b(t_1,w_1)$.

$\Rightarrow$ from $M \vDash \text{Firstf}(x,t_1,aya,t_2)$, $M \vDash (t_1 a)Bx$,

$\Rightarrow M \vDash \exists x_1 \, t_1 a x_1 = x = w_1 at_1 ayat_2 ay'at_4$,

$\Rightarrow$ by (4.14$^b$), $M \vDash w_1=t_1 \lor t_1 B w_1$,

$\Rightarrow M \vDash t_1 \subseteq_p w_1$,

$\Rightarrow$ from $M \vDash \text{Max}^+T_b(t_1,w_1)$, $M \vDash t_1<t_1$, contradicting $M \vDash t_1 \in I \subseteq I_0$.

Exactly the same argument applies if

$M \vDash x=w_1 at_1 ayat_2 ay'at_4 aw_2 \,\&\, \text{Max}^+T_b(t_1,w_1)$.

Assume, again for a reductio, that

$M \vDash t_2<t_4 \,\&\, \text{Pref}(ay'a,t_2) \,\&\, x=t_1 ayat_2 ay'at_4 aw_2$.

$\Rightarrow$ from $M \vDash \text{Lastf}(x,t_3,ay'a,t_4)$,

$M \vDash \text{Pref}(ay'a, t_3) \,\&\, (x = t_3 ay' at_4 \,\vee\, \exists w'(x = w' at_3 ay' at_4 \,\&\, \text{Max}^+ T_b(t_3, w')))$.

Suppose that $M \vDash x = t_3 ay' at_4$. Then from $M \vDash \text{Pref}(ay', t_3) \,\&\, t_3 = t_4$, we also have $M \vDash \text{Firstf}(x, t_3, ay'a, t_4)$. Hence from $M \vDash \text{Set}(x)$, by (5.21),

$$M \vDash \forall w(w \,\varepsilon\, x \leftrightarrow w = y').$$

But then from $M \vDash y \,\varepsilon\, x$, $M \vDash y = y'$, contradicting the principal hypothesis.

Suppose, on the other hand, that $M \vDash \exists w'(x = w' at_3 ay' at_4 \,\&\, \text{Max}^+ T_b(t_3, w'))$.

$\implies M \vDash w' at_3 ay' at_4 = x = t_1 ay at_2 ay' at_4 aw_2$,

$\implies$ by (4.15$^b$), $M \vDash w_2 = t_4 \,\vee\, t_4 E w_2$.

(1) $M \vDash w_2 = t_4$.

$\implies M \vDash w' at_3 ay' at_4 = x = t_1 ay at_2 ay' at_4 at_4$,

$\implies$ by (3.6), $M \vDash w' at_3 ay' = x = t_1 ay at_2 ay' at_4$,

$\implies$ by (4.15$^b$), $M \vDash y' = t_4 \,\vee\, t_4 E y'$,

$\implies M \vDash t_4 \subseteq_p y'$,

$\implies$ from $M \vDash \text{Pref}(ay', t_3)$, $M \vDash \text{Max}^+ T_b(t_3, ay'a)$,

$\implies M \vDash t_4 < t_3$,

$\implies M \vDash t_4 < t_3 = t_4$, contradicting $M \vDash t_4 \in I \subseteq I_0$.

(2) $M \vDash t_4 E w_2$.

$\implies M \vDash \exists w_3 \, w_2 = w_3 t_4$,

$\implies M \vDash w' at_3 ay' at_4 = x = t_1 ay at_2 ay' at_4 aw_3 t_4$,

$\implies$ by (3.6), $M \vDash w' at_3 (ay'a) = t_1 ay at_2 ay'a(t_4 aw_3)$.

We claim that $M \vDash \text{Max}^+ T_b(t_3, t_1 ay at_2 ay')$.

For, suppose $M \vDash \text{Tally}_b(t') \,\&\, t' \subseteq_p t_1 ay at_2 ay'$.

$\Rightarrow$ by (4.17$^b$), $M \vDash t' \subseteq_p t_1ay \lor t' \subseteq_p t_2ay'$,

$\Rightarrow$ by (4.17$^b$), $M \vDash t' \subseteq_p t_1 \lor t' \subseteq_p y \lor t' \subseteq_p t_2 \lor t' \subseteq_p y'$,

$\Rightarrow$ by an argument analogous to one above, $M \vDash x \neq t_1ayat_2$,

$\Rightarrow$ from $M \vDash \text{Firstf}(x,t_1,aya,t_2)$, $M \vDash t_1 < t_2$ & $\text{Pref}(aya,t_1)$,

$\Rightarrow$ from $M \vDash \text{Pref}(ay'a,t_2)$, $M \vDash t' \leq t_1 \lor t' \leq t_2$,

$\Rightarrow M \vDash t' \leq t_2 < t_4 = t_3$.

Hence, $M \vDash \text{Max}^+T_b(t_3, t_1ayat_2ay')$ as claimed.

$\Rightarrow$ from $M \vDash w'at_3(ay'a) = t_1ayat_2ay'a(t_4aw_3)$, by (4.19),

$\qquad M \vDash (t_1ayat_2ay'a)B(w'a) \lor w' = t_1ayat_2ay'$.

(2a) $M \vDash w' = t_1ayat_2ay'$.

$\Rightarrow M \vDash (w'at_3a)y'at_4 = (w'at_4a)w_3t_4$,

$\Rightarrow$ by (3.7), $M \vDash y'at_4 = w_3t_4$,

$\Rightarrow$ by (3.6), $M \vDash y'a = w_3$,

$\Rightarrow M \vDash w'at_3ay'at_4 = x = t_1ayat_2ay'at_4a(y'a)t_4$,

$\Rightarrow M \vDash \text{Pref}(ay'a,t_2)$ & $t_2 < t_4$ & $\exists w_4\ x = (t_1ay)at_2ay'at_4aw_4$ & $\text{Max}^+T_b(t_2,t_1ay)$,

$\Rightarrow M \vDash \text{Intf}(x,t_1ay,t_2,ay'a,t_4)$,

$\Rightarrow M \vDash \text{Env}(t,x)$ & $\text{Lastf}(x,t_3,ay'a,t_4)$ & $\text{Intf}(x,t_1ay,t_2,ay'a,t_4)$, contradicting (5.19).

(2b) $M \vDash (t_1ayat_2ay'a)B(w'a)$.

$\Rightarrow M \vDash \exists w_4\ (t_1ayat_2ay'a)w_4 = w'a$,

$\Rightarrow M \vDash (t_1ayat_2ay'aw_4)at_3ay'at_4 = x = t_1ayat_2ay'at_4aw_3t_4$,

$\Rightarrow$ by (3.7), $M \vDash w_4at_3ay'at_4 = t_4aw_3t_4$,

$\Rightarrow$ by (3.6), $M \vDash w_4at_3ay'a=t_4aw_3$,

$\Rightarrow$ by (4.14[b]), $M \vDash w_4=t_4 \vee t_4Bw_4$,

$\Rightarrow$ from $M \vDash t_1ayat_2ay'aw_4=w'a$ & $Tally_b(t_4)$, $M \vDash \neg(w_4=t_4)$,

$\Rightarrow M \vDash t_4Bw_4$,

$\Rightarrow M \vDash \exists w_5\ t_4w_5=w_4$,

$\Rightarrow M \vDash t_1ayat_2ay'a(t_4w_5)=w'a$,

$\Rightarrow M \vDash w_5=a \vee aEw_5$,

$\Rightarrow M \vDash t_1ayat_2ay'at_4a=w'a \vee \exists w_6\ t_1ayat_2ay'at_4(w_6a)=w'a$,

$\Rightarrow M \vDash t_1ayat_2ay'at_4=w' \vee t_1ayat_2ay'at_4w_6=w'$,

$\Rightarrow M \vDash t_4\subseteq_p w'$,

$\Rightarrow$ from $M \vDash Max^+T_b(t_3,w')$,

$\Rightarrow M \vDash t_4<t_3=t_4$, contradicting $M \vDash t_4 \in I \subseteq I_0$.

Therefore, $M \vDash \neg\exists w'(x=w'at_3ay'at_4$ & $Max^+T_b(t_3,w'))$,

whence $M \vDash \neg(t_2<t_4$ & $Pref(ay'a,t_2)$ & $x=t_1ayat_2ay'at_4aw_2)$.

The only remaining possibility is that

$$M \vDash t_2=t_4 \text{ & } x=t_1ayat_2ay'at_4.$$

Hence $M \vDash x=t_1ayat_2ay'at_2$, as required.

An analogous argument shows that

$$M \vDash \exists t_1,t_2(Tally_b(t_1) \text{ & } Tally_b(t_2) \text{ & } x=t_1ay'at_2ayat_2)$$

if $M \vDash y'<_x y$.

This completes the proof of (10.4).

(10.5) For any string concept $I \subseteq I_0$ there is a string concept $J \subseteq I$ such that

$$QT^+ \vdash \forall x \in J \, \forall t,u \, (Pref(aua,t) \,\&\, x = tauat \rightarrow MinSet(x)).$$

Let $J \equiv I_{4.14b} \,\&\, I_{4.15b} \,\&\, I_{5.22}$.

Assume $M \vDash Pref(aua,t) \,\&\, x=tauat$ where $M \vDash J(x)$.

$\Rightarrow M \vDash Max^+T_b(t,u)$,

$\Rightarrow M \vDash Firstf(x,t,aua,t) \,\&\, Lastf(x,t,aua,t)$,

$\Rightarrow$ by (5.22), $M \vDash Env(t,x)$,

$\Rightarrow M \vDash Set(x)$.

Assume now that $M \vDash w_1 a w_2 = z$.

$\Rightarrow M \vDash w_1 a w_2 = tauat \,\&\, Tally_b(t)$,

$\Rightarrow$ by (4.14$^b$) and (4.15$^b$), $M \vDash (w_1=t \,\vee\, tBw_1) \,\&\, (w_2=t \,\vee\, tEw_2)$.

We distinguish four scenarios:

(1) $M \vDash w_1=t=w_2$.

$\Rightarrow M \vDash tat = tauat$,

$\Rightarrow$ by (3.7), $M \vDash t = uat$, a contradiction because $M \vDash Tally_b(t)$.

(2) $M \vDash w_1=t \,\&\, tEw_2$.

$\Rightarrow M \vDash w_1 a w_2 = z \,\&\, Firstf(z,t,aua,t) \,\&\, t=w_1$,

$\Rightarrow M \vDash Occ(w_1,a,w_2,z,t,aua,t)$.

(3) $M \vDash tBw_1 \,\&\, w_2=t$.

$\Rightarrow M \vDash \exists w_3 \, tw_3 = w_1$,

$\Rightarrow M \vDash w_1 a w_2 = (tw_3)at$,

$\Rightarrow$ $M \vDash w_1at=(tw_3)at=tauat$,

$\Rightarrow$ by (3.6), $M \vDash w_1a=taua$,

$\Rightarrow$ $M \vDash w_1aw_2=z$ & $Firstf(z,t,aua,t)$ & $w_1a=taua$,

$\Rightarrow$ $M \vDash Occ(w_1,a,w_2,z,t,aua,t)$.

(4) $M \vDash tBw_1$ & $tEw_2$.

$\Rightarrow$ $M \vDash \exists w_4\ w_1a(w_4t)=w_1aw_2=tauat$,

$\Rightarrow$ by (3.6), $M \vDash w_1aw_4=taua$,

$\Rightarrow$ $M \vDash (w_1a)B(taua)$,

$\Rightarrow$ $M \vDash w_1aw_2=z$ & $Firstf(z,t,aua,t)$ & $(w_1a)B(taua)$,

$\Rightarrow$ $M \vDash Occ(w_1,a,w_2,z,t,aua,t)$.

Hence we have shown that

$M \vDash \forall w_1,w_2\ (w_1aw_2=z \rightarrow \exists v,t_1,t_2\ Occ(w_1,a,w_2,z,t_1,ava,t_2))$.

It follows that $M \vDash MinSet(z)$.

This completes the proof of (10.5).

(10.6) For any string concept $I\subseteq I_0$ there is a string concept $J\subseteq I$ such that

$QT^+ \vdash \forall x\in J \forall t,t_1,t_2,t_3,v,w,x',z[Env(t_2,x')$ & $x'=t_1wt_2$ & $aBw$ & $aEw$ & $x=t_1wz$ &

& $Env(t,z)$ & $t_2<t_3$ & $Firstf(z,t_3,ava,t_4)$ & $\neg\exists u(u\ \varepsilon\ x'\ \&\ u\ \varepsilon\ z)$ &

& $MinSet(x')$ & $MinSet(z)\ \to\ Env(t,x)$ & $MinSet(x)]$.

Let $J \equiv I_{5.46}$.

Assume $M \vDash Env(t_2,x')$ & $Env(t,z)$ & $MinSet(x')$ & $MinSet(z)$

where $M \vDash x'=t_1wt_2$ & $aBw$ & $aEw$ & $x=t_1wz$ & $t_2<t_3$ & $Firstf(z,t_3,ava,t_4)$ &

& $\neg\exists u(u\ \varepsilon\ x'\ \&\ u\ \varepsilon\ z)$

and $M \vDash J(x)$.

$\Rightarrow$ by (5.46), $M \vDash Env(t,x)$,

$\Rightarrow M \vDash Set(x)$.

Assume $M \vDash w_1aw_2=x$.

$\Rightarrow M \vDash w_1aw_2=x=t_1wz$,

$\Rightarrow M \vDash (w_1a)Bx$ & $(t_1w)Bx$,

$\Rightarrow$ by (3.8), $M \vDash (w_1a)B(t_1w)$ ∨ $w_1a=t_1w$ ∨ $(t_1w)B(w_1a)$.

(1) $M \vDash (w_1a)B(t_1w)$.

$\Rightarrow M \vDash \exists w_3\ w_1aw_3=t_1w$,

$\Rightarrow M \vDash (w_1aw_3)t_2=t_1wt_2=x'$ & $w_1aw_3z=t_1wz=x=w_1aw_2$,

$\Rightarrow$ from $M \vDash MinSet(x')$, $M \vDash \exists v',t',t''Occ(w_1,a,w_3t_2,x',t',v',t'')$.

 (1a) $M \vDash Firstf(x',t',v',t'')$ & $(t'=w_1$ ∨ $(w_1a)B(t'v')$ ∨ $w_1a=t'v')$.

$\Rightarrow$ from $M \vDash Env(t,z)$, by (5.11), $M \vDash \exists t',z'(z=t'z't$ & $aBz'$ & $aEz')$.

Hence from $M \vDash \text{Firstf}(z,t_3,av_0a,t_4)$ we have, by the proof of (5.6)(1), that

$$M \vDash \text{Firstf}(x,t',v',t'') \ \&\ (t'=w_1 \vee (w_1a)B(t'v') \vee w_1a=t'v').$$

$\Longrightarrow M \vDash \text{Occ}(w_1,a,w_2,x,t',v',t'')$.

(1b) $M \vDash \exists w'(\text{Intf}(x',w',t',v',t'') \ \&\ (w'at'=w_1 \vee \exists v_1(v_1Bv' \ \&\ w'at'v_1=w_1a) \vee$

$$\vee\ w'at'v'=w_1a)).$$

Analogous to (1a), except that from the proof of (5.6)(2) we obtain

$M \vDash \text{Intf}(x,w',t',v',t'')$.

(1c) $M \vDash \exists w'((\text{Lastf}(x',t',v',t'') \ \&\ w'at'v't''=x') \ \&$

$$\&\ (w'at'=w_1 \vee \exists v_1(v_1Bv' \ \&\ w'at'v_1=w_1a) \vee w'at'v'=w_1a)).$$

$\Longrightarrow M \vDash w'at'v't''=x'=t_1wt_2$,

$\Longrightarrow M \vDash \text{Pref}(v',t') \ \&\ t'=t''$,

$\Longrightarrow M \vDash \exists v_0\ v=av_0a \ \&\ \text{Tally}_b(t')$,

$\Longrightarrow$ from $M \vDash aEw$, $M \vDash \exists w_4\ w=w_4a$,

$\Longrightarrow M \vDash w'at'(av_0a)t''=x'=t_1(w_4a)t_2$,

$\Longrightarrow$ from $M \vDash \text{Env}(t_2,x')$, $M \vDash \text{Tally}_b(t_2)$,

$\Longrightarrow$ by ($4.24^b$), $M \vDash t''=t_2$,

$\Longrightarrow$ from hypothesis $M \vDash t_2<t_3$, $M \vDash t''<t_3$,

$\Longrightarrow M \vDash \exists t_0\ t''t_0=t_3$,

$\Longrightarrow$ as in the proof of (4.24)(3), $M \vDash \text{Intf}(x,w',t',v',t''t_0)$,

$\Longrightarrow M \vDash \exists w'(\text{Intf}(x,w',t',v',t_3) \ \&\ (w'at'=w_1 \vee \exists v_1(v_1Bv' \ \&\ w'at'v_1=w_1a) \vee$

$$\vee\ w'at'v'=w_1a)),$$

$\Longrightarrow M \vDash \text{Occ}(w_1,a,w_2,x,t',v',t_3)$.

We have thus shown that

$$M \vDash w_1aw_2 = x \ \& \ (w_1a)B(t_1w) \ \rightarrow \ \exists v',t_5,t_6 \ Occ(w_1,a,w_2,x,t_5,v',t_6).$$

(2) $M \vDash w_1a = t_1w$.

$\Rightarrow$ from $M \vDash Env(t_2,x')$, $M \vDash \exists v_0 \ Lastf(x',t_2,av_0a,t_2)$,

$\Rightarrow M \vDash Pref(av_0a,t_2) \ \& \ (x' = t_2av_0at_2 \ \lor \ \exists w'(x' = w'at_2av_0at_2 \ \& \ Max^+T_b(t_2,w')))$.

 (2a) $M \vDash x' = t_2av_0at_2$.

$\Rightarrow M \vDash t_1wt_2 = x' = t_2av_0at_2$,

$\Rightarrow$ by (3.6), $M \vDash t_1w = t_2av_0a$,

$\Rightarrow M \vDash w_1a = t_1w = t_2av_0a$,

$\Rightarrow$ from $M \vDash Firstf(z,t_3,ava,t_4)$, $M \vDash (t_3a)Bz \ \& \ Tally_b(t_3)$,

$\Rightarrow M \vDash \exists z_1 \ t_3az_1 = z$,

$\Rightarrow M \vDash x = t_1wz = t_2av_0a(t_3az_1) = t_2v't_3az_1 \quad \text{for } v' = av_0a$,

$\Rightarrow M \vDash Pref(v',t_2) \ \& \ Tally_b(t_3) \ \& \ t_2 < t_3 \ \& \ (t_2v't_3a)Bx$,

$\Rightarrow M \vDash Firstf(x,t_2,v',t_3) \ \& \ w_1a = t_2v' \ \& \ w_1aw_2 = x$,

$\Rightarrow M \vDash Occ(w_1,a,w_2,x,t_2,v',t_3)$.

 (2b) $M \vDash \exists w'(x' = w'at_2av_0at_2 \ \& \ Max^+T_b(t_2,w'))$.

$\Rightarrow M \vDash t_1wt_2 = x' = w'at_2av_0at_2$,

$\Rightarrow$ by (3.6), $M \vDash t_1w = w'at_2av_0a$,

$\Rightarrow M \vDash x = t_1wz = w'at_2av_0az$,

$\Rightarrow$ as in (2a), $M \vDash w'at_2v'z = w'at_2v't_3az_1$,

$\Rightarrow M \vDash Pref(v',t_2) \ \& \ Tally_b(t_3) \ \& \ t_2 < t_3 \ \& \ \exists w'' \ x = w'at_2v't_3aw'' \ \& \ Max^+T_b(t_2,w')$,

$\Rightarrow M \vDash Intf(x,w',t_2,v',t_3) \ \& \ w'at_2v' = w_1a \ \& \ w_1aw_2 = x$,

$\Rightarrow$ $M \vDash Occ(w_1,a,w_2,x,t_2,v',t_3)$.

Hence we have shown that

$$M \vDash w_1aw_2 = x \,\&\, w_1a = t_1w \rightarrow \exists v',t_5,t_6\, Occ(w_1,a,w_2,x,t_5,v',t_6).$$

(3) $M \vDash (t_1w)B(w_1a)$.

$\Rightarrow$ $M \vDash \exists w_3\, t_1ww_3 = w_1a$,

$\Rightarrow$ $M \vDash t_1ww_3w_2 = w_1aw_2 = x = t_1wz$,

$\Rightarrow$ by (3.7), $M \vDash w_3w_2 = z$,

$\Rightarrow$ from $M \vDash t_1ww_3 = w_1a$, $M \vDash w_3 = a \lor aEw_3$,

$\Rightarrow$ from $M \vDash Firstf(z,t_3,ava,t_4)$, $M \vDash (t_3a)Bz \,\&\, Tally_b(t_3)$,

$\Rightarrow$ $M \vDash \exists z_1\, t_3az_1 = z$,

$\Rightarrow$ $M \vDash w_3w_2 = t_3az_1$.

Now, we cannot have $M \vDash w_3 = a$ because $M \vDash Tally_b(t_3)$.

$\Rightarrow$ $M \vDash aEw_3$,

$\Rightarrow$ $M \vDash \exists w_4\, w_3 = w_4a$,

$\Rightarrow$ $M \vDash t_1w(w_4a) = w_1$,

$\Rightarrow$ $M \vDash w_4aw_2 = t_3az_1 = z$,

$\Rightarrow$ by (4.14[b]), $M \vDash w_4 = t_3 \lor t_3Bw_4$.

(3a) $M \vDash w_4 = t_3$.

$\Rightarrow$ $M \vDash t_1wt_3 = w_1$.

We claim that $M \vDash Max^+T_b(t_3,t_1w)$.

Assume $M \vDash t_0 \subseteq_p t_1w \,\&\, Tally_b(t_0)$.

$\Rightarrow$ from $M \vDash Env(t_2,x')$, $M \vDash MaxT_b(t_2,x')$,

$\implies$ from $M \vDash t_1w \subseteq_p x'$, $M \vDash t_0 \leq t_2 < t_3$, as required.

$\implies M \vDash \text{Firstf}(z,t_3,ava,t_4)$ & $x=t_1wz$ & $aBw$ & $aEw$ & $\text{Max}^+T_b(t_3,t_1w)$.

Then, from the proof of (5.25) cases (1) and (3), we have that

$\quad M \vDash (\text{Lastf}(x,t_3,ava,t_4)$ & $x=w'at_3avat_4) \lor \text{Intf}(x,w',t_3,ava,t_4)$

where $M \vDash w'=t_1w_5$ & $w_5a=w$. This, along with $M \vDash w_1aw_2=x$ and

$M \vDash w'at_3=t_1w_5at_3=t_1wt_3=w_1$ yields $M \vDash \text{Occ}(w_1,a,w_2,x,t_3,ava,t_4)$.

(3b) $M \vDash t_3Bw_4$.

$\implies M \vDash \exists w_6\ t_3w_6=w_4$,

$\implies M \vDash t_1wt_3w_6=w_1$,

$\implies M \vDash w_1aw_2=x=(t_1wt_3w_6)aw_2=t_1wz$,

$\implies$ by (3.7), $M \vDash t_3w_6aw_2=z$.

From hypothesis $M \vDash \text{MinSet}(z)$ we have $M \vDash \exists v',t',t''\text{Occ}(t_3,w_6,a,w_2,z,t',v',t'')$.

(3bi) $M \vDash \text{Firstf}(z,t',v',t'')$ & $(t'=t_3w_6 \lor (t_3w_6)B(t'v') \lor t_3w_6a=t'v')$.

$\implies$ from $M \vDash \text{Firstf}(z,t_3,ava,t_4)$, by (5.15), $M \vDash v'=ava$ & $t'=t_3$ & $t''=t_4$,

$\implies$ as in (3a), $M \vDash (\text{Lastf}(x,t_3,ava,t_4)$ & $x=w'at_3avat_4) \lor \text{Intf}(x,w',t_3,ava,t_4)$

where $M \vDash w'=t_1w_5$ & $w_5a=w$.

If $M \vDash t_3=t_3w_6$, then $M \vDash t_3Bt_3$, contradicting the hypothesis $M \vDash t_3 \in I \subseteq I_0$.

$\implies M \vDash (t_3w_6)B(t_3ava) \lor t_3w_6a=t_3ava$,

$\implies M \vDash \exists w_7\ t_3w_6w_7=t_3ava \lor t_3w_6a=t_3ava$,

$\implies M \vDash t_1wt_3w_6w_7=t_1wt_3ava \lor t_1wt_3w_6a=t_1wt_3ava$,

$\implies$ from hypothesis $M \vDash \text{Firstf}(z,t_3,ava,t_4)$, $M \vDash (t_3ava)Bz$,

$\implies M \vDash \exists z_2\ t_3avaz_2=z$.

(3bi1)  $M \vDash t_1wt_3w_6a=t_1wt_3ava$.

$\Rightarrow$ $M \vDash w_1a=t_1wt_3ava$.

(3bi2)  $M \vDash t_1wt_3w_6w_7=t_1wt_3ava$.

$\Rightarrow$ $M \vDash t_1wt_3w_6w_7z_2=t_1w(t_3avaz_2)=t_1wz=x=w_1aw_2=(t_1wt_3w_6)aw_2$,

$\Rightarrow$ by (3.7), $M \vDash w_7z_2=aw_2$,

$\Rightarrow$ $M \vDash w_7=a \lor aBw_7$,

$\Rightarrow$ $M \vDash t_1wt_3w_6a=t_1wt_3ava \lor \exists w_8 (w_7=aw_8 \ \& \ t_1wt_3w_6aw_8=t_1wt_3ava)$.

(3bi2a)  $M \vDash t_1wt_3w_6a=t_1wt_3ava$.

$\Rightarrow$ $M \vDash w_1a=t_1wt_3ava$.

(3bi2b)  $M \vDash t_1wt_3w_6aw_8=t_1wt_3ava$.

$\Rightarrow$ by (3.7), $M \vDash w_6aw_8=ava$,

$\Rightarrow$ $M \vDash \exists v_1 (v_1w_8=ava \ \& \ v_1B(ava))$,

$\Rightarrow$ $M \vDash w_1aw_8=(t_1wt_3w_6)aw_8=t_1wt_3v_1w_8$,

$\Rightarrow$ by (3.6), $M \vDash w_1a=t_1wt_3v_1$.

Hence from (3bi1)-(3bi2) we have

$\quad M \vDash w_1a=t_1wt_3ava \lor \exists v_1(v_1B(ava) \ \& \ w_1a=t_1wt_3v_1)$.

$\Rightarrow M \vDash t_1w_5at_3ava=w_1a \lor \exists v_1(v_1B(ava) \ \& \ w_1a=t_1w_5at_3v_1)$,

$\Rightarrow M \vDash w_1aw_2=x \ \& \ ((Lastf(x,t_3,ava,t_4) \ \& \ x=w'at_3avat_4) \lor Intf(x,w',t_3,ava,t_4)) \ \&$

$\quad\quad\quad\quad\quad \& \ (w'at_3ava=w_1a \lor \exists v_1(v_1B(ava) \ \& \ w'at_3v_1=w_1a))$,

$\Rightarrow M \vDash Occ(w_1,a,w_2,x,t_3,ava,t_4)$.

(3bii)  $M \vDash \exists w'((Intf(z,w',t',v',t'') \lor (Lastf(z,t',v',t'') \ \& \ x=w'at'v't'')) \ \&$

$\quad \& \ (w'at'=t_3w_6 \lor \exists v_1(v_1Bv' \ \& \ w'at'v_1=t_3w_6a) \lor w'at'v'=t_3w_6a))$.

As in the proof of (5.25), parts (4) and (2), we have that

$$M \vDash \text{Intf}(x,w'',t',v',t'') \lor (\text{Lastf}(x,t',v',t'') \,\&\, x=w''at'v't'')$$

where $M \vDash w''=t_1w_4aw' \,\&\, w_4a=w$.

Now, if $M \vDash w'at'=t_3w_6$, then $M \vDash w''at'=(t_1w)w'at'=(t_1w)t_3w_6=w_1$;

and if $M \vDash \exists v_1(v_1Bv' \,\&\, w'at'v_1=t_3w_6a)$, then

$$M \vDash w'at'v_1=(t_1w)w'at'v_1=(t_1w)t_3w_6a=w_1a;$$

and if $M \vDash w'at'v'=t_3w_6a$, then $M \vDash w''at'v'=(t_1w)w'at'v'=(t_1w)t_3w_6a=w_1a$.

Therefore

$$M \vDash w_1aw_2=x \,\&\, \exists w''((\text{Intf}(x,w'',t',v',t'') \lor (\text{Lastf}(x,t',v',t'') \,\&\, x=w''at'v't'')) \,\&$$
$$\&\, (w''at'=w_1 \lor \exists v_1(v_1Bv' \,\&\, w''at'v_1=w_1a) \lor w''at'v'=w_1a),$$
$$\implies M \vDash \text{Occ}(w_1,a,w_2,x,t',v',t'').$$

So we have established that

$$M \vDash w_1aw_2=x \,\&\, (t_1w)B(w_1a) \;\to\; \exists v',t_5,t_6 \, \text{Occ}(w_1,a,w_2,x,t_5,v',t_6).$$

From (1)-(3) we thus have $M \vDash w_1aw_2=x \;\to\; \exists v',t_5,t_6 \, \text{Occ}(w_1,a,w_2,x,t_5,v',t_6)$.

Along with $M \vDash \text{Set}(x)$, this suffices to obtain $M \vDash \text{MinSet}(x)$.

This completes the proof of (10.6).

(10.7) For any string concept $I\subseteq I_0$ there is a string concept $J\subseteq I$ such that

$QT^+ \vdash \forall x \in J\ \forall t_1,t_2,t_3,u,v(Pref(aua,t_1)\ \&\ Pref(ava,t_2)\ \&\ t_1<t_2\ \&\ t_2=t_3\ \&$

$\&\ u\neq v\ \&\ x=t_1auat_2avat_3\ \rightarrow\ MinSet(x))$.

Let $J \equiv I_{5.46}\ \&\ I_{10.5}$.

Assume $M \vDash Pref(aua,t_1)\ \&\ Pref(ava,t_2)\ \&\ t_1<t_2\ \&\ t_2=t_3\ \&\ u\neq v$ and $M \vDash J(x)$.

Let $M \vDash x=t_1auat_2avat_3$.

$\Longrightarrow$ by (5.58), $M \vDash Set(x)$,

$\Longrightarrow$ by (10.5), $M \vDash MinSet(x')\ \&\ MinSet(z)$ where $x'=t_1auat_1$ and $z=t_2avat_2$,

$\Longrightarrow M \vDash Firstf(x',t_1,aua,t_1)\ \&\ Lastf(x',t_1,aua,t_1)$,

$\Longrightarrow$ by (5.22), $M \vDash \forall w(w\ \varepsilon\ x' \leftrightarrow w=u)$.

Likewise, $M \vDash \forall w(w\ \varepsilon\ z \leftrightarrow w=v)$.

$\Longrightarrow$ from $M \vDash u\neq v$, $M \vDash \neg\exists w(w\ \varepsilon\ x'\ \&\ w\ \varepsilon\ z)$,

$\Longrightarrow$ by (10.6), $M \vDash MinSet(x)$.

This completes the proof of (10.7).

(10.8) For any string concept $I \subseteq I_0$ there is a string concept $J \subseteq I$ such that

$QT^+ \vdash \forall y \in J \forall x,t,t',t'',u(Env(t,y)$ & $MinSet(y)$ & $Fr(y,t',aua,t'')$ & $Env(t',x)$ &

& $Lastf(x,t',aua,t')$ & $xBy \rightarrow MinSet(x))$.

Let $J \equiv I_{4.6}$ & $I_{5.24}$ & $I_{5.27}$ & $I_{5.33}$ & $I_{5.57}$.

Assume $M \vDash Env(t,y)$ & $MinSet(y)$ & $Fr(y,t',aua,t'')$ and $M \vDash xBy$ where

$M \vDash Env(t',x)$ & $Lastf(x,t',aua,t')$ and $M \vDash J(y)$.

$\Rightarrow M \vDash Set(x)$.

Assume $M \vDash w_1aw_2=x$.

$\Rightarrow$ from hypothesis $M \vDash xBy$, $M \vDash \exists z \, xz=y$,

$\Rightarrow M \vDash (w_1aw_2)z=y$,

$\Rightarrow$ from hypothesis $M \vDash MinSet(y)$, $M \vDash \exists v,t_1,t_2 \, Occ(w_1,a,w_2z,y,t_1,ava,t_2)$,

$\Rightarrow M \vDash Fr(y,t_1,ava,t_2)$.

We distinguish three cases:

(1) $M \vDash Firstf(y,t_1,ava,t_2)$.

$\Rightarrow$ from $M \vDash Env(t',x)$ & $xBy$ & $Firstf(y,t_1,ava,t_2)$, by (5.4),

$M \vDash \exists t_3 Firstf(x,t_1,ava,t_3)$.

$\Rightarrow$ from $M \vDash Occ(w_1,a,w_2z,y,t_1,ava,t_2)$, $M \vDash t_1=w_1 \vee (w_1a)B(t_1ava) \vee w_1=t_1ava$,

$\Rightarrow$ from $M \vDash w_1aw_2=x$ & $Firstf(x,t_1,ava,t_3)$, $M \vDash Occ(w_1,a,w_2,x,t_1,ava,t_3)$.

(2) $M \vDash \exists w' \, Intf(y,w',t_1,ava,t_2)$.

$\Rightarrow$ from $M \vDash Occ(w_1,a,w_2z,y,t_1,ava,t_2)$,

$M \vDash w'at_1=w_1 \vee \exists v_1(v_1B(ava)$ & $w'at_1v_1=w_1a) \vee w'at_1ava=w_1a$.

Suppose that

(2a)  $M \vDash t' < t_1$.

If  $M \vDash w'at_1v_1 = w_1a$,  then  $M \vDash v_1 = a \vee aEv_1$,  hence

$$M \vDash w'at_1a = w_1a \vee \exists v_2\ w'at_1(v_2a) = w_1a.$$

Likewise, if  $M \vDash w'at_1ava = w_1a$,  then  $M \vDash w'at_1a = w_1a \vee \exists v'\ w'at_1v'a = w_1a$.

Hence we have

  $M \vDash w'at_1 = w_1 \vee \exists v_1, v_2\ (v_1 B(ava)\ \&\ w'at_1v_2 = w_1) \vee \exists v'\ w'at_1v' = w_1$.

$\Longrightarrow M \vDash t_1 \subseteq_p w_1 \subseteq_p x$,

$\Longrightarrow$ from $M \vDash \text{Env}(t',x)$, $M \vDash \text{MaxT}_b(t',x)$,

$\Longrightarrow$ from  $M \vDash t_1 \subseteq_p x$,  $M \vDash t_1 \leq t'$,

$\Longrightarrow M \vDash t_1 \leq t' < t_1$,  contradicting  $M \vDash t_1 \in I \subseteq I_0$.

(2b)  $M \vDash t' = t_1$.

$\Longrightarrow$ from  $M \vDash \text{Fr}(y,t',aua,t'')\ \&\ \text{Fr}(y,t_1,ava,t_2)\ \&\ \text{Env}(t,y)$, $M \vDash u = v$,

$\Longrightarrow M \vDash \text{Lastf}(x,t',ava,t')$.

We show that  $M \vDash x \neq t'avat'$.

Note from $M \vDash \text{Env}(t,y)$  we have  $M \vDash \exists v_0, t_3, t_4\ \text{Firstf}(y,t_3,av_0a,t_4)$,

$\Longrightarrow$ by (5.4), $M \vDash \exists t_5\ \text{Firstf}(x,t_3,av_0a,t_5)$,

$\Longrightarrow$ from  $M \vDash \text{Intf}(y,w',t_1,ava,t_2)$, by (5.19), $M \vDash \neg\exists t_6\ \text{Firstf}(y,t',ava,t_6)$,

$\Longrightarrow$ from  $M \vDash \text{Env}(t,y)\ \&\ \text{Env}(t',x)\ \&\ y'By$, by (5.15), $M \vDash \neg\text{Firstf}(x,t',ava,t')$,

$\Longrightarrow$ from  $M \vDash \text{Lastf}(x,t',ava,t')$, $M \vDash \text{Pref}(ava,t')$,

$\Longrightarrow$ from  $M \vDash \neg\text{Firstf}(x,t',ava,t')$, $M \vDash x \neq t'avat'$,  as required.

Therefore, from  $M \vDash \text{Lastf}(x,t',ava,t')$,

$$M \vDash \exists w_3 \ (x = w_3 at'avat' \ \& \ Max^+T_b(t',w_3)).$$

$\Longrightarrow$ from hypotheses (2) and (2b), $M \vDash \exists w'' \ y = w'at'avat_2aw'' \ \& \ Max^+T_b(t',w')$,

$\Longrightarrow$ from $M \vDash xz = y$, $M \vDash (w_3at'avat')z = w'at'avat_2aw''$,

$\Longrightarrow$ from $M \vDash Max^+T_b(t',w_3) \ \& \ Max^+T_b(t',w')$, by (4.20), $M \vDash w' = w_3$.

Therefore, $M \vDash Lastf(x,t',ava,t') \ \& \ x = w'at'avat'$.

Hence from

$$M \vDash w'at' = w_1 \ \lor \ \exists v_1,v_2 \ (v_1B(ava) \ \& \ w'at'v_1 = w_1a) \ \lor \ w'at'ava = w_1a$$

we obtain $M \vDash Occ(w_1,a,w_2,x,t',ava,t')$, as required.

(2c) $M \vDash t_1 < t'$.

Then we have

$M \vDash Env(t,y) \ \& \ Fr(y,t',aua,t'') \ \& \ xBy \ \& \ Env(t',x) \ \& \ Lastf(x,t',aua,t') \ \&$

$\& \ Intf(y,w',t_1,ava,t_2) \ \& \ t_1 < t'$.

$\Longrightarrow$ by the proof of (5.56), part (2b), $M \vDash \exists t_3 Intf(x,w',t_1,ava,t_3)$.

Hence from

$$M \vDash w'at' = w_1 \ \lor \ \exists v_1,v_2 \ (v_1B(ava) \ \& \ w'at'v_1 = w_1a) \ \lor \ w'at'ava = w_1a$$

we have $M \vDash Occ(w_1,a,w_2,x,t_1,ava,t_3)$, as required.

(3) $M \vDash Lastf(y,t_1,ava,t_2) \ \& \ \exists w' \ y = w'at_1avat_2$.

$\Longrightarrow$ from hypothesis $M \vDash Env(t,y)$, $M \vDash t_1 = t_2 = t$.

(3a) $M \vDash t' < t_1$.

$\Longrightarrow$ from hypothesis

$$M \vDash w'at_1 = w_1 \ \lor \ \exists v_1(v_1B(ava) \ \& \ w'at_1v_1 = w_1a) \ \lor \ w'at_1ava = w_1a,$$

we have $M \vDash t_1 \subseteq_p w_1 \subseteq_p x$,

$\Rightarrow$ from  $M \vDash Env(t',x)$,  $M \vDash MaxT_b(t',x)$,

$\Rightarrow$  $M \vDash t_1 \leq t' < t_1$,  contradicting  $M \vDash t_1 \in I \subseteq I_0$.

Suppose that

   (3b)   $M \vDash t'=t_1$.

$\Rightarrow$ from  $M \vDash Fr(y,t',aua,t'')$ & $Fr(y,t_1,ava,t_2)$ & $Env(t,y)$,  $M \vDash u=v$,

$\Rightarrow$ from  $M \vDash t_1=t$,  $M \vDash t'=t_1=t$,

$\Rightarrow$  $M \vDash Env(t',y)$ & $Fr(y,t',aua,t'')$ & $Env(t'x)$ & $Lastf(y,t,ava,t)$ &

$\qquad\qquad\qquad\qquad\qquad\qquad$ & $xBy$ & $Lastf(x,t',ava,t')$

contradicting (5.57).

   (3c)   $M \vDash t_1<t'$.

$\Rightarrow$  $M \vDash t=t_1<t'$,

$\Rightarrow$ from  $M \vDash Env(t',x)$,  $M \vDash MaxT_b(t',x)$,

$\Rightarrow$ from  $M \vDash Env(t,y)$,  $M \vDash MaxT_b(t,y)$.

Now, if  $M \vDash Tally_b(t_0)$ & $t_0 \subseteq_p x$, then $M \vDash t_0 \subseteq_p y$, since $M \vDash xBy$. Hence  $M \vDash t' \leq t$.

But then  $M \vDash t' \leq t = t_1 < t'$,  contradicting  $M \vDash t' \in I \subseteq I_0$.

Hence we have, by (4.6),  that

$\qquad$ $M \vDash w_1aw_2=x \rightarrow \exists v,t',t''\ Occ(w_1,a,w_2,x,t',ava,t'')$,

which suffices to show that in Case (3) also  $M \vDash MinSet(x)$,  as required.

This completes the proof of (10.8).

RESOLUTION LEMMA. (10.9)  For any string concept $I \subseteq I_0$ there is a string concept $I_{RES} \subseteq I$ such that

$QT^+ \vdash \forall y \in I_{RES} \forall t, t_1, t_2, u, w'[Env(t,y) \& MinSet(y) \& Inf(y,w',t_1,aua,t_2) \rightarrow$

$\rightarrow \exists y', y'', t_0, t', w^*, w''(Env(t',y') \& Env(t,y'') \& y'=t_0w^*t' \&$

$\& y''=t_1auat_2aw'' \& y=t_0w^*y'' \& Firstf(y'',t_1,aua,t_2) \&$

$\& t' < t_1 \& aBw^* \& aEw^* \& \neg \exists w(w \;\varepsilon\; y' \& w \;\varepsilon\; y''))].$

Let $I_{RES} \equiv (I_{10.8})_{SUB}$.

We may assume that $I_{RES}$ is closed under * and downward closed under $\subseteq_p$.

Assume  $M \vDash Env(t,y) \& MinSet(y) \& Inf(y,w',t_1,aua,t_2)$  where  $M \vDash I_{RES}(y)$.

$\Rightarrow M \vDash Pref(aua,t_1) \& Tally_b(t_2) \& t_1 < t_2 \& \exists w'' \; y=w'at_1auat_2aw'' \&$

$\& Max^+T_b(t_1,w')$.

$\Rightarrow$ from $M \vDash I_{RES}(y) \& w' \subseteq_p y \& t_1 \subseteq_p y \& u \subseteq_p y \& t_2 \subseteq_p y \& w'' \subseteq_p y$,

$M \vDash I_{RES}(w') \& I_{RES}(t_1) \& I_{RES}(u) \& I_{RES}(t_2) \& I_{RES}(w'')$.

Letting  $y''=t_1auat_2aw''$, we have that  $M \vDash I_{RES}(y'')$.

$\Rightarrow$ by (5.44),  $M \vDash Env(t,y'') \& Firstf(y'',t_1,aua,t_2)$.

$\Rightarrow$ by (5.11),  $M \vDash \exists y_1, t_0(Tally_b(y_1) \& y=t_0y_1t \& aBy_1 \& aEy_1)$,

$\Rightarrow$ from $M \vDash I_{RES}(y) \& t_0 \subseteq_p y \& y_1 \subseteq_p y \& t \subseteq_p y$,  $M \vDash I_{RES}(t_0) \& I_{RES}(y_1) \& I_{RES}(t)$,

$\Rightarrow$ from  $M \vDash Env(t,y)$,  $M \vDash \exists v_0, t_3, t_4 \; Firstf(y,t_3,av_0a,t_4)$,

$\Rightarrow M \vDash Pref(av_0a,t_3) \& Tally_b(t_4) \& ((t_3=t_4 \& y=t_3av_0at_4) \lor$

$\lor (t_3 < t_4 \& (t_3av_0at_4a)By))$.

If  $M \vDash t_3=t_4 \& y=t_3av_0at_4$, then from  $M \vDash Pref(av_0a,t_3)$ we have

$M \vDash Firstf(y,t_3,av_0a,t_4) \& Lastf(y,t_3,av_0a,t_4)$.

$\Rightarrow$ by (5.22), $M \vDash Env(t_3,y)$ & $\forall w(w \, \varepsilon \, y \leftrightarrow w=v_0)$,

$\Rightarrow M \vDash u=v_0$,

$\Rightarrow M \vDash Firstf(y,t_3,av_0a,t_4)$ & $Inf(y,w',t_1,aua,t_2)$, contradicting (5.19).

Therefore $M \vDash t_3 < t_4$ & $(t_3av_0at_4a)By$.

$\Rightarrow$ from $M \vDash aBy_1$, $M \vDash \exists y_2 \, ay_2=y_1$,

$\Rightarrow$ from $M \vDash (t_3av_0at_4a)By$,

$\qquad M \vDash \exists y_2 \, (t_3av_0at_4a)y_3=y=w'at_1auat_2aw''=t_0y_1t=t_0(ay_2)t$,

$\Rightarrow$ since $M \vDash Tally_b(t_3)$ & $Tally_b(t_0)$, by $(4.23^b)$, $M \vDash t_3=t_0$,

$\Rightarrow$ from $M \vDash Firstf(y,t_3,av_0a,t_4)$ & $Inf(y,w',t_1,aua,t_2)$, by (5.34), $M \vDash t_4 \leq t_1$,

$\Rightarrow$ from $M \vDash Pref(av_0a,t_0)$, $M \vDash Max^+T_b(t_0,av_0a)$,

$\Rightarrow$ since $M \vDash t_0 < t_4$, $M \vDash Max^+T_b(t_4,av_0a)$,

$\Rightarrow$ by $(4.17^b)$, $M \vDash Max^+T_b(t_4,t_0av_0a)$,

$\Rightarrow$ from $M \vDash t_4 \leq t_1$, $M \vDash Max^+T_b(t_1,t_0av_0a)$,

$\Rightarrow$ from $M \vDash t_0av_0at_4ay_3=y=w'at_1auat_2aw''$, by (4.19),

$\qquad M \vDash (t_0av_0a)B(w'a) \vee t_0av_0a=w'a$,

$\Rightarrow M \vDash \exists w_5 \, t_0av_0aw_5=w'a \vee t_0av_0a=w'a$,

$\Rightarrow$ either way, $M \vDash \exists w^*(t_0w^*=w'a \, \& \, aBw^* \, \& \, aEw^*)$,

$\Rightarrow M \vDash y=w'at_1auat_2aw''=t_0w^*(t_1auat_2aw'')=t_0w^*y''$,

where $M \vDash aBw^*$ & $aEw^*$.

Let $x'=w'at_1auat_1$.

$\Rightarrow$ from $M \vDash I_{RES}(w')$ & $I_{RES}(t_1)$ & $I_{RES}(u)$, $M \vDash I_{RES}(x')$,

$\Rightarrow M \vDash Env(t,y)$ & $(w'at_1)By$ & $Inf(y,w',t_1,aua,t_2)$,

$\Rightarrow$ by (5.53), $M \vDash Env(t_1,x') \& Lastf(x',t_1,aua,t_1)$,

$\Rightarrow M \vDash Env(t,y) \& MinSet(y) \& Firstf(y,t_1,aua,t_2) \& Env(t_1,x') \&$

$\& Lastf(x',t_1,aua,t_1) \& x'By$,

$\Rightarrow$ by (10.8), $M \vDash MinSet(x')$,

$\Rightarrow$ by SUBTRACTION LEMMA,

$M \vDash \exists x^- \in I_{RES}(Set(x^-) \& \neg(u \varepsilon x^-) \& \forall w(w \varepsilon x^- \leftrightarrow w \varepsilon x' \& w \neq u))$.

In fact, by the proof of SUBTRACTION LEMMA, parts (2biii2a2c) and (2biii2b3), we have

$M \vDash \exists t',v_0(Env(t',x^-) \& ((x'=w'at_1auat_1=t'av_0at_1auat_1 \& t'<t_1 \& x^-=t'av_0at') \vee$

$\vee \exists w_0(x'=w'at_1auat_1=w_0t'av_0at_1auat_1 \& t'<t_1 \& x^-=w_0at'av_0at')))$,

$\Rightarrow$ from $M \vDash t'<t_1 \& Tally_b(t')$, $M \vDash \exists t'' \; t_1=t't''$,

$\Rightarrow$ from $M \vDash I_{RES}(t_1) \& t' \subseteq_p t_1 \& t'' \subseteq_p t_1$, $M \vDash I_{RES}(t') \& I_{RES}(t'')$,

$\Rightarrow M \vDash w'at'(t''auat_1)=t'av_0at'(t''auat_1) \vee$

$\vee \exists w_0 \; w'at'(t''auat_1)=w_0t'av_0at'(t''auat_1)$,

$\Rightarrow$ from $M \vDash I_{RES}(t'') \& I_{RES}(u) \& I_{RES}(t_1)$, $M \vDash I_{RES}(t''auat_1)$,

$\Rightarrow$ by (3.6), $M \vDash x^-=t'av_0at'=w'at' \vee x^-=w_0at'av_0at'=w'at'$,

$\Rightarrow M \vDash x^-=w'at'$,

$\Rightarrow M \vDash x^-t''auat_1=(w'at')t''auat_1=w'at_1auat_1=x'$,

$\Rightarrow$ from $M \vDash t_1<t_2 \& Tally_b(t_1)$, $M \vDash \exists t^+ \; t_1t^+=t_2$,

$\Rightarrow M \vDash (x^-t''auat_1)t^+aw''=(w'at_1auat_1)t^+aw''=w'at_1auat_2aw''=y=$

$=t_0w*t_1auat_2aw''=t_0w*t'(t''auat_1t^+aw'')$,

$\Rightarrow$ from $M \vDash I_{RES}(t_2) \& t^+ \subseteq_p t_2,$ , $M \vDash I_{RES}(t^+)$,

$\Rightarrow$ from $M \vDash I_{RES}(t^+)$ & $I_{RES}(w")$, $M \vDash I_{RES}(t^+aw")$,

$\Rightarrow$ by (3.6), $M \vDash x^- = t_0 w^* t'$.

So we may let $y' = x^-$.

Suppose, for a reductio, that $M \vDash w \, \varepsilon \, y'$ & $w \, \varepsilon \, y"$.

$\Rightarrow M \vDash \exists t_3, t_4, t_5, t_6 (Fr(y', t_3, awa, t_4) \, \& \, Fr(y", t_5, awa, t_6))$,

$\Rightarrow$ from $M \vDash Env(t', y')$, $M \vDash MaxT_b(t', y')$,

$\Rightarrow M \vDash t_3 \leq t'$,

$\Rightarrow$ from $M \vDash Firstf(y", t_1, aua, t_2)$ & $Fr(y", t_5, awa, t_6)$, by (5.20), $M \vDash t_1 \leq t_5$,

$\Rightarrow M \vDash t_3 \leq t' < t_1 \leq t_5$.

Applying (5.11) to $y"$ we have from $M \vDash Env(t, y")$ that

$$M \vDash \exists t_6, w^{**}(Tally_b(t_6) \, \& \, y" = t_6 w^{**} t \, \& \, aBw^{**} \, \& \, aEw^{**}),$$

$\Rightarrow M \vDash \exists w_6 \, t_1 auat_2 aw" = y" = t_6(aw_6)t \quad$ where $M \vDash w^{**} = aw_6$,

$\Rightarrow$ by (4.23$^b$), $M \vDash t_1 = t_6$,

$\Rightarrow M \vDash y = t_0 w^* t_1 auat_2 aw" = t_0 w^* t_1 w^{**} t = t_0 w^* t' t" w^{**} t = y' t" w^{**} t$.

Hence we have

$$M \vDash Env(t, y) \, \& \, y = y' t" w^{**} t \, \& \, aBw^{**} \, \& \, aEw^{**} \, \& \, Tally_b(t") \, \& \, Tally_b(t),$$

$\Rightarrow$ from $M \vDash Fr(y', t_3, awa, t_4)$, by (5.6), $M \vDash \exists t_7 \, Fr(y, t_3, awa, t_7)$,

$\Rightarrow$ from $M \vDash MaxT_b(t', y')$ & $t' < t_5$ & $t_0 w^* \subseteq_p y'$, $M \vDash Max^+T_b(t_5, t_0 w^*)$,

$\Rightarrow M \vDash y = t_0 w^* y"$ & $aBw^*$ & $aEw^*$ & $Max^+T_b(t_5, t_0 w^*)$,

$\Rightarrow$ from $M \vDash Fr(y", t_5, awa, t_6)$, by (5.25), $M \vDash Fr(y, t_5, awa, t_6)$,

$\Rightarrow$ from $M \vDash Env(t, y)$, $M \vDash t_3 = t_5$,

$\Rightarrow M \vDash t_3 < t_5 = t_3$, contradicting $M \vDash t_3 \in I \subseteq I_0$.

Hence  M ⊨ ¬∃w(w ε y' & w ε y'').

This completes the proof of THE RESOLUTION LEMMA.

(10.10) For any string concept $I\subseteq I_0$ there is a string concept $J\subseteq I$ such that

$QT^+ \vdash \forall y \in J \forall t,t_1,t_2,w',w'',u,x[Env(t,y)$ & $MinSet(y)$ & $Intf(y,w',t_1,aua,t_2)$ &

& $x=t_1auat_2aw''$ & $y=w'at_1auat_2aw''$ & $Firstf(x,t_1,aua,t_2) \rightarrow MinSet(x)]$.

Let $J \equiv I_{RES}$,

where $I_{RES}$ is obtained from I as in THE RESOLUTION LEMMA.

Assume $M \vDash Env(t,y)$ & $MinSet(y)$ & $y=w'at_1auat_2aw''$ where

$M \vDash Intf(y,w',t_1,aua,t_2)$ and $M \vDash J(y)$.

Assume $M \vDash x=t_1auat_2aw''$ & $Firstf(x,t_1,aua,t_2)$.

$\Rightarrow$ by (5.44), $M \vDash Env(t,x)$.

To show that $M \vDash MinSet(x)$, assume that $M \vDash w_1aw_2=x$.

$\Rightarrow$ from $M \vDash y=w'at_1auat_2aw''=w'ax$, $M \vDash (w'aw_1)aw_2=y$,

$\Rightarrow$ from hypothesis $M \vDash MinSet(y)$, $M \vDash \exists v,t',t''\ Occ(w'aw_1,a,w_2,y,t',ava,t'')$,

$\Rightarrow M \vDash Fr(y,t',ava,t'')$.

(1) $M \vDash Firstf(y,t',ava,t'')$ & $(t'=w'aw_1 \lor (w'aw_1)B(t'ava) \lor w'aw_1=t'ava)$.

  (1a) $M \vDash t'=w'aw_1$.

This contradicts $M \vDash Tally_b(t')$.

  (1b) $M \vDash (w'aw_1)B(t'ava) \lor w'aw_1=t'ava$.

$\Rightarrow$ from $M \vDash w_1aw_2=x$ & $Firstf(x,t_1,aua,t_2)$, $M \vDash (t_1a)Bx$,

$\Rightarrow M \vDash \exists x_1\ t_1ax_1=x=w_1aw_2$,

$\Rightarrow$ by (4.14$^b$), $M \vDash w_1=t_1 \lor t_1Bw_1$,

$\Rightarrow M \vDash t_1\subseteq_p w_1$,

$\Rightarrow$ M ⊨ $\exists w_3$ w'aw$_1$aw$_3$=t'ava v w'aw$_1$a=t'ava,

$\Rightarrow$ by (4.14$^b$), M ⊨ w'=t' v t'Bw'.

(1bi) M ⊨ w'=t'.

$\Rightarrow$ M ⊨ (t'a)w$_1$aw$_3$=t'ava v (t'a)w$_1$a=(t'a)va,

$\Rightarrow$ M ⊨ w$_1$aw$_3$=va v w$_1$a=va,

$\Rightarrow$ M ⊨ $t_1 \subseteq_p$va,

$\Rightarrow$ from M ⊨ Firstf(y,t',ava,t''), M ⊨ Pref(ava,t'),

$\Rightarrow$ M ⊨ $t_1$<t',

$\Rightarrow$ from M ⊨ Firstf(y,t',ava,t'') & Intf(y,w',$t_1$,aua,$t_2$), by (5.20), M ⊨ t'≤$t_1$,

$\Rightarrow$ M ⊨ $t_1$<t'≤$t_1$, contradicting M ⊨ $t_1 \in I \subseteq I_0$.

(1bii) M ⊨ t'Bw'.

$\Rightarrow$ M ⊨ $\exists w_4$ t'w$_4$=w',

$\Rightarrow$ M ⊨ (t'w$_4$)aw$_1$aw$_3$=t'ava v t'w$_4$aw$_1$a=t'ava,

$\Rightarrow$ by (3.7), M ⊨ w$_4$aw$_1$aw$_3$=ava v w$_4$aw$_1$a=ava,

$\Rightarrow$ M ⊨ $w_1 \subseteq_p$ava,

$\Rightarrow$ M ⊨ $t_1 \subseteq_p$ava, whence we derive a contradiction as in (1bi).

(2) M ⊨ $\exists w_3$ (Intf(y,w$_3$,t',ava,t'') & (w$_3$at'=w'aw$_1$ v

v $\exists v_1$($v_1$B(ava) & w$_3$at'v$_1$=w'aw$_1$a) v w$_3$at'ava=w'aw$_1$a)).

(2a) M ⊨ t'<$t_1$.

$\Rightarrow$ M ⊨ $w_1 \subseteq_p$w$_3$at' v $w_1 \subseteq_p$w$_3$at'v$_1$ v $w_1 \subseteq_p$w$_3$at'(ava),

$\Rightarrow$ M ⊨ $t_1 \subseteq_p w_1 \subseteq_p$w$_3$at'ava,

$\Rightarrow$ by (4.17$^b$), M ⊨ $t_1 \subseteq_p$w$_3$ v $t_1 \subseteq_p$t'ava,

$\Longrightarrow$ by (4.17$^b$), $M \vDash t_1 \subseteq_p w_3 \lor t_1 \subseteq_p t' \lor t_1 \subseteq_p va$,

$\Longrightarrow$ from $M \vDash \text{Intf}(y,w_3,t',ava,t'')$, $M \vDash \text{Pref}(ava,t') \& \text{Max}^+T_b(t',w_3)$,

$\Longrightarrow M \vDash t_1 \leq t'$,

$\Longrightarrow M \vDash t_1 \leq t' < t_1$, contradicting $M \vDash t_1 \in I \subseteq I_0$.

(2b) $M \vDash t' = t_1$.

$\Longrightarrow$ from $M \vDash \text{Fr}(y,t',ava,t'') \& \text{Fr}(y,t_1,aua,t_2) \& \text{Env}(t,y)$, $M \vDash v = u$,

$\Longrightarrow$ from $M \vDash \text{Intf}(y,w',t_1,aua,t_2) \& \text{Intf}(y,w_3,t',ava,t'')$,

$\quad M \vDash \exists w_4\, w'at_1auat_2aw''=y=w_3at'avat''aw_4 \& \text{Max}^+T_b(t_1,w') \&$

$\quad\quad\quad\quad\quad\quad\quad\quad \& \text{Max}^+T_b(t',w_3) \& \text{Tally}_b(t_2) \& \text{Tally}_b(t'')$,

$\Longrightarrow$ by (4.20), $M \vDash w' = w_3$,

$\Longrightarrow$ by (3.7), $M \vDash t_2 aw'' = t'' aw_4$,

$\Longrightarrow$ by (4.23$^b$), $M \vDash t_2 = t''$.

(2bi) $M \vDash w_3 at' = w' a w_1$.

$\Longrightarrow$ by (3.7), $M \vDash t' = w_1$,

$\Longrightarrow M \vDash \text{Firstf}(x,t',ava,t'') \& t' = w_1$,

$\Longrightarrow M \vDash \text{Occ}(w_1,a,w_2,x,t',ava,t'')$,

$\Longrightarrow$ since $M \vDash t' = t_1 \& t_2 = t''$, $M \vDash \text{Occ}(w_1,a,w_2,x,t_1,ava,t_2)$.

(2bii) $M \vDash \exists v_1(v_1 B(ava) \& (w_3 a)t'v_1 = (w'a)w_1 a)$.

$\Longrightarrow$ since $M \vDash w' = w_3$, by (3.7), $M \vDash t'v_1 = w_1 a$,

$\Longrightarrow M \vDash \text{Firstf}(x,t',ava,t'') \& \exists v_1(v_1 B(ava) \& t'v_1 = w_1 a)$,

$\Longrightarrow M \vDash \text{Occ}(w_1,a,w_2,x,t',ava,t'')$,

$\Longrightarrow$ since $M \vDash t' = t_1 \& t_2 = t'' \& v = u$, $M \vDash \text{Occ}(w_1,a,w_2,x,t_1,aua,t_2)$.

(2biii)  $M \vDash w_3at'ava=w'aw_1a$.

$\Rightarrow$ from  $M \vDash w'=w_3$,  by (3.7),  $M \vDash t'ava=w_1a$,

$\Rightarrow$ $M \vDash \text{Firstf}(x,t',ava,t'') \ \& \ t'ava=w_1a$,

$\Rightarrow$ $M \vDash \text{Occ}(w_1,a,w_2,x,t',ava,t'')$,

$\Rightarrow$ since  $M \vDash t'=t_1 \ \& \ t_2=t''$,  $M \vDash \text{Occ}(w_1,a,w_2,x,t_1,ava,t_2)$.

(2c)  $M \vDash t_1<t'$.

$\Rightarrow$ since  $M \vDash t_1 \in I \subseteq I_0$,  $M \vDash t' \neq t_1$,

$\Rightarrow$ from  $M \vDash \text{Fr}(y,t_1,aua,t_2) \ \& \ \text{Fr}(y,t',ava,t'') \ \& \ \text{Env}(t,y)$,  $M \vDash u \neq v$,

$\Rightarrow$ from  $M \vDash \text{Env}(t,y) \ \& \ \text{MinSet}(y) \ \& \ \text{Intf}(y,w',t_1,aua,t_2)$,

  by RESOLUTION LEMMA,

$M \vDash \exists y',t_0,t^+,w^*(y=t_0w^*x \ \& \ y'=t_0w^*t^+ \ \& \ t^+<t_1 \ \& \ aBw^* \ \& \ aEw^* \ \& \ \text{Env}(t^+,y'))$.

We have that  $M \vDash v \ \varepsilon \ y$.

$\Rightarrow$ from  $M \vDash \text{Env}(t,x)$,  by (5.11),

   $M \vDash \exists t^*,w^+(x=t^*w^+t \ \& \ aBw^+ \ \& \ aEw^+ \ \& \ \text{Tally}_b(t^*))$,

$\Rightarrow$ $M \vDash \exists x_1 \ w^+=ax_1$,

$\Rightarrow$ $M \vDash t^*ax_1t=t^*w^+t=x$,

$\Rightarrow$ from  $M \vDash \text{Firstf}(x,t_1,aua,t_2)$,  $M \vDash (t_1a)Bx \ \& \ \text{Tally}_b(t_1)$,

$\Rightarrow$ $M \vDash \exists x_2 \ t_1ax_2=t$,

$\Rightarrow$ $M \vDash t^*ax_1t=x=t_1ax_2$,

$\Rightarrow$ by (4.23$^b$),  $M \vDash t^*=t_1$,

$\Rightarrow$ from  $M \vDash t^+<t_1 \ \& \ \text{Tally}_b(t^+)$,  $M \vDash \exists t^{++}(\text{Tally}_b(t^{++}) \ \& \ t_1=t^+t^{++})$,

$\Rightarrow$ $M \vDash y=t_0w^*x=t_0w^*(t^*w^+t)=t_0w^*t_1w^+t=t_0w^*(t^+t^{++})w^+t=y't^{++}w^+t$,

$\implies$ M ⊨ Env(t,y) & y=$t_0$w*$t^+t^{++}$w$^+$t & Tally$_b$($t_0$) & aBw* & aEw* & aBw$^+$ &

& aEw$^+$ & y'=$t_0$w*$t^+$ & x=$t^+t^{++}$w$^+$t & Env($t^+$,y') &

& Env(t,x) & Firstf(x,$t^+t^{++}$,aua,$t_2$),

$\implies$ by (5.41), M ⊨ ¬∃w(w ε y' & w ε x),

$\implies$ by (5.46), M ⊨ ∀w(w ε y ↔ w ε y' v w ε x),

$\implies$ from  M ⊨ v ε y,  M ⊨ v ε y' v v ε x.

Suppose that   M ⊨ v ε y'.

$\implies$ M ⊨ ∃$t_3$,$t_4$ Fr(y',$t_3$,ava,$t_4$),

$\implies$ by (5.6), M ⊨ ∃$t_5$ Fr(y,$t_3$,ava,$t_5$),

$\implies$ from  M ⊨ Fr(y,t',ava,t'') & Env(t,y),  M ⊨ $t_3$=t',

$\implies$ from  M ⊨ Env($t^+$,y'),  M ⊨ MaxT$_b$($t^+$,y'),

$\implies$ M ⊨ $t_3$≤$t^+$<$t_1$<t'=$t_3$,  contradicting  M ⊨ $t_3$∈I⊆$I_0$.

Hence  M ⊨ ¬(v ε y').

$\implies$ M ⊨ v ε x,

$\implies$ since  M ⊨ v≠u & Firstf(x,$t_1$,aua,$t_2$),  by (5.15),  M ⊨ ¬∃$t_3$,$t_4$ Firstf(x,$t_3$,ava,$t_4$),

$\implies$ from  M ⊨ v ε x,  M ⊨ ∃$t_3$,$t_4$ Fr(x,$t_3$,ava,$t_4$),

$\implies$ M ⊨ ∃$w_5$ Intf(x,$w_5$,$t_3$,ava,$t_4$) v Lastf(x,$t_3$,ava,$t_4$).

Suppose, for a reductio, that  M ⊨ Lastf(x,$t_3$,ava,$t_4$).

$\implies$ from the proof of (5.25), parts (1) and (2),  M ⊨ Lastf(y,$t_3$,ava,$t_4$).

But, by (5.19), this contradicts the hypothesis  M ⊨ Intf(y,$w_3$,t',ava,t'').

Hence   M ⊨ ¬Lastf(x,$t_3$,ava,$t_4$).

$\implies$ M ⊨ ∃$w_5$ Intf(x,$w_5$,$t_3$,ava,$t_4$),

$\implies$ by (5.25), $M \vDash Fr(y,t_3,ava,t_4)$,

$\implies$ from $M \vDash Fr(y,t',ava,t'') \& Env(t,y)$, $M \vDash t_3=t'$,

$\implies$ $M \vDash \exists w_6\ x=w_5at'avat_4aw_6 \& Max^+T_b(t',w_5)$,

$\implies$ $M \vDash y=w'ax=w'aw_5at'avat_4aw_6$,

$\implies$ from the proof of (5.25), part (4), $M \vDash Max^+T_b(t',w'aw_5)$.

But we also have $M \vDash y=w_3at'avat''aw_4 \& Max^+T_b(t',w_3)$.

$\implies$ by (4.15), $M \vDash w'aw_5=w_3$,

$\implies$ from hypothesis (2),

$M \vDash (w'aw_5)at'=w'aw_1\ \lor\ \exists v_1(v_1B(ava) \& (w'aw_5)at'v_1=w'aw_1a)\ \lor$

$\lor\ (w'aw_5)at'ava=w'aw_1a$,

$\implies$ by (3.7), $M \vDash w_5at'=w'\ \lor\ \exists v_1(v_1B(ava) \& w_5at'v_1=w_1a)\ \lor$

$\lor\ w_5at'(ava)=w_1a)\ \&\ Intf(x,w_5,t',ava,t_4)$,

$\implies$ $M \vDash Occ(w_1,a,w_2,x,t',ava,t_2)$.

Hence from hypothesis (2) by (4.6) we have that

$\qquad M \vDash w_1aw_2=x\ \rightarrow\ \exists v,t_3,t_4\ Occ(w_1,a,w_2,x,t_3,ava,t_4)$.

(3) $M \vDash Lastf(y,t',ava,t'') \& y=w_3at'avat'' \& (w_3at'=w'aw_1\ \lor$

$\lor\ \exists v_1(v_1B(ava) \& w_3at'v_1=w'aw_1a)\ \lor w_3at'ava=w'aw_1a)$.

$\implies$ from $M \vDash Env(t,y)$,

$M \vDash \exists v_0,t_3,t_4\ Firstf(y,t_3,av_0a,t_4) \& \exists v'Lastf(y,t,av'a,t) \& MaxT_b(t,y)$,

$\implies$ from hypothesis $M \vDash Intf(y,w',t_1,aua,t_2)$,

$M \vDash y=w'at_1auat_2aw'' \& x=t_1auat_2aw'' \& t_1<t_2\leq t$,

$\implies$ by (5.43), $M \vDash Lastf(x,t,av'a,t)$,

$\Rightarrow$ from hypothesis (3), by (5.15), $M \vDash t'=t''=t$ & $v=v'$,

$\Rightarrow$ from $M \vDash Lastf(x,t,av'a,t)$,

$M \vDash Pref(av'a,t)$ & $(x=tav'at \lor \exists w_4(x=w_4atav'at$ & $Max^+T_b(t',w_4)))$.

Suppose that $M \vDash x=tav'at$.

$\Rightarrow M \vDash tav'at=t_1auat_2aw''$,

$\Rightarrow$ by (4.23$^b$), $M \vDash t=t_1<t_2\leq t$, contradicting $M \vDash t\in I\subseteq I_0$.

Therefore, $M \vDash x\neq tav'at$.

$\Rightarrow M \vDash \exists w_4(x=w_4atav'at$ & $Max^+T_b(t',w_4))$,

$\Rightarrow M \vDash y=w'ax=w'aw_4atav'at$, while also $M \vDash y=w_3at'avat''=w_3atav'at$,

$\Rightarrow$ by (3.6), $M \vDash w'aw_4=w_3$,

$\Rightarrow$ from hypothesis (3),

$M \vDash (w'aw_4)at'=w'aw_1 \lor \exists v_1(v_1B(ava)$ & $(w'aw_4)at'v_1=w'aw_1a) \lor$

$\lor (w'aw_4)at'ava=w'aw_1a$,

$\Rightarrow$ by (3.7), $M \vDash (w_4at=w' \lor \exists v_1(v_1B(ava)$ & $w_4at'v_1=w_1a) \lor w_4at'ava=w_1a)$ &

& $Lastf(x,t',ava,t'')$ & $x=w_4atav'at$,

$\Rightarrow M \vDash Occ(w_1,a,w_2,x,t',ava,t'')$.

From (1)-(3) we therefore have

$M \vDash w_1aw_2=x \rightarrow \exists v,t_3,t_4\, Occ(w_1,a,w_2,x,t_3,ava,t_4)$,

which suffices to establish $M \vDash MinSet(x)$.

This completes the proof of (10.10).

(10.11) For any string concept $I \subseteq I_0$ there is a string concept $J \subseteq I$ such that

$QT^+ \vdash \forall x \in J \forall t, t_1, t_2, w_1, w_2, u, x''(Env(t,x) \& Intf(x,w_1,t_1,aua,t_2) \&$

$\& x = w_1 a t_1 a u a t_2 a w_2 \& x'' = t_2 a w_2 \& Env(t,x'') \& MinSet(x) \rightarrow MinSet(x''))$.

Let $J \equiv I_{10.5} \& I_{10.10}$.

Assume $M \vDash Env(t,x) \& Intf(x,w_1,t_1,aua,t_2) \& x = w_1 a t_1 a u a t_2 a w_2$ and

$M \vDash x'' = t_2 a w_2 \& Env(t,x'')$ where $M \vDash J(x)$.

Assume also that $M \vDash MinSet(x)$.

Let $x' = w_1 a t_1 a u a t_1$.

$\implies$ by (5.53), $M \vDash Env(t_1, x')$,

$\implies$ by (5.11), $M \vDash \exists t_0, w^* (Tally_b(t_0) \& x' = t_0 w^* t_1 \& aBw^* \& aEw^*)$,

$\implies$ from $M \vDash Intf(x,w_1,t_1,aua,t_2)$,

$M \vDash Pref(aua,t_1) \& Tally_b(t_2) \& t_1 < t_2 \& \exists w_2\, x = w_1 a t_1 a u a t_2 a w_2 \& Max^+T_b(t_1,w_1)$,

$\implies$ since $M \vDash Tally_b(t_1)$, $M \vDash \exists t'(Tally_b(t') \& t_1 t' = t_2)$,

$\implies M \vDash x = w_1 a t_1 a u a t_2 a w_2 = w_1 a t_1 a u a t_1 t' a w_2 = x' t' a w_2 = t_0 w^* t_1 t' a w_2 = t_0 w^* t_2 a w_2 =$

$= t_0 w^* x''$.

We claim that $M \vDash Max^+T_b(t_2, t_0 w^*)$.

Assume $M \vDash t'' \subseteq_p t_0 w^* \& Tally_b(t'')$.

$\implies M \vDash t'' \subseteq_p t_0 w^* \subseteq_p x' = w_1 a t_1 a u a t_1$,

$\implies$ by (4.17$^b$), $M \vDash t'' \subseteq_p w_1 \lor t'' \subseteq_p x' = t_1 a u a t_1$,

$\implies$ by (4.17$^b$), $M \vDash t'' \subseteq_p w_1 \lor t'' \subseteq_p t_1 \lor t'' \subseteq_p u a t_1$,

$\implies$ by (4.17$^b$), $M \vDash t'' \subseteq_p w_1 \lor t'' \subseteq_p t_1 \lor t'' \subseteq_p u \subseteq_p aua$,

$\Rightarrow$ from $M \vDash Max^+T_b(t_1,w_1)$ & $Pref(aua,t_1)$, $M \vDash t''\leq t_1<t_2$.

Hence $M \vDash Max^+T_b(t_2,t_0w^*)$, as claimed.

From $M \vDash Env(t,x'')$, $M \vDash \exists v,t_3,t_4\ Firstf(x'',t_3,ava,t_4)$.

$\Rightarrow M \vDash (t_3a)Bx''$ & $Tally_b(t_3)$,

$\Rightarrow M \vDash \exists x_2\ t_3ax_2=x''=t_2aw_2$,

$\Rightarrow$ by (4.23$^b$), $M \vDash t_3=t_2$.

Then we have

$\quad M \vDash Fr(x'',t_2,ava,t_4)$ & $x=t_0w^*x''$ & $aBw^*$ & $aEw^*$ & $Max^+T_b(t_2,t_0w^*)$.

$\Rightarrow$ by (5.25), $M \vDash Fr(x,t_2,ava,t_4)$.

We distinguish three cases.

(1) $M \vDash Firstf(x,t_2,ava,t_4)$.

$\Rightarrow$ from $M \vDash Intf(x,w_1,t_1,aua,t_2)$, by (5.20), $M \vDash t_2\leq t_1$,

$\Rightarrow M \vDash t_2\leq t_1<t_2$, contradicting $M \vDash t_2 \in I \subseteq I_0$.

(2) $M \vDash \exists w_3 Intf(x,w_3,t_2,ava,t_4)$.

$\Rightarrow M \vDash \exists w_4\ x=w_3at_2avat_4aw_4$ & $Max^+T_b(t_2,w_3)$,

$\Rightarrow M \vDash w_1at_1auat_2aw_2=x=w_3at_2avat_4aw_4$.

We also have $M \vDash w_1at_1auat_2aw_2=x=t_0w^*t_2aw_2$.

$\Rightarrow$ by (3.6), $M \vDash w_1at_1aua=t_0w^*$,

$\Rightarrow M \vDash w_1at_1auat_2aw_2=x=w_3at_2avat_4aw_4$ & $Max^+T_b(t_2,w_1at_1au)$,

$\Rightarrow$ by (4.20), $M \vDash w_1at_1au=w_3$,

$\Rightarrow$ by (3.7), $M \vDash x''=t_2aw_2=t_2avat_4aw_4$.

So we have,

$M \vDash \text{Env}(t,x) \ \& \ \text{MinSet}(x) \ \& \ \text{Intf}(x,w_3,t_2,\text{ava},t_4) \ \& \ x=w_3at_2\text{ava}t_4aw_4 \ \&$

$\& \ x''=t_2\text{ava}t_4aw_4 \ \& \ \text{Firstf}(x'',t_2,\text{ava},t_4),$

$\Rightarrow$ by (10.10), $M \vDash \text{MinSet}(x'')$, as required.

(3) $M \vDash \text{Lastf}(x,t_2,\text{ava},t_4)$.

$\Rightarrow M \vDash \text{Pref}(\text{ava},t_2) \ \& \ t_2=t_4 \ \& \ (x=t_2\text{ava}t_4 \ \vee \ \exists w(x=wat_2\text{ava}t_4 \ \& \ \text{Max}^+T_b(t_2,w)))$.

(3a) $M \vDash x=t_2\text{ava}t_4$.

$\Rightarrow M \vDash w_1at_1auat_2aw_2=t_2\text{ava}t_4$,

$\Rightarrow$ by (4.14$^b$), $M \vDash w_1=t_2 \ \vee \ t_2Bw_1$,

$\Rightarrow M \vDash t_2\subseteq_p w_1$,

$\Rightarrow$ from $M \vDash \text{Max}^+T_b(t_1,w_1)$, $M \vDash t_2<t_1$,

$\Rightarrow M \vDash t_2<t_1<t_2$, contradicting $M \vDash t_2\in I\subseteq I_0$.

(3b) $M \vDash \exists w(x=wat_2\text{ava}t_4 \ \& \ \text{Max}^+T_b(t_2,w))$.

$\Rightarrow M \vDash w_1at_1auat_2aw_2=wat_2\text{ava}t_4$,

$\Rightarrow$ as in (2), $M \vDash \text{Max}^+T_b(t_2,w_1at_1au)$,

$\Rightarrow$ by (4.20), $M \vDash w_1at_1au=w$,

$\Rightarrow$ by (3.7), $M \vDash x''=t_2aw_2=t_2\text{ava}t_4$,

$\Rightarrow$ from $M \vDash \text{Pref}(\text{ava},t_2) \ \& \ t_2=t_4$, by (10.5), $M \vDash \text{MinSet}(x'')$, as required.

This completes the proof of (10.11).

(10.12) For any string concept $I \subseteq I_0$ there is a string concept $J \subseteq I$ such that

$$QT^+ \vdash \forall x \in J \forall s,t_1,t_2,t_3,t_4,u,w[\text{Set}(x) \,\&\, \text{Firstf}(x,t_1,aua,t_2) \,\&$$

$$\&\, ((\text{Lastf}(x,t_3,awa,t_4) \,\&\, x=t_1auat_2st_3awat_4) \,\vee$$

$$\vee\, \exists w_3,w_4 (\text{Intf}(x,w_3,t_3,awa,t_4) \,\&\, x=t_1auat_2st_3awat_4)) \,\&\, (s=a \,\vee\, s=aa) \,\rightarrow$$

$$\rightarrow \neg\text{MinSet}(x)].$$

Let $J \equiv I_{5.15} \,\&\, I_{5.20}$.

(i) Assume $M \vDash x=t_1auat_2st_3awat_4$

where $M \vDash \text{Set}(x) \,\&\, \text{Firstf}(x,t_1,aua,t_2) \,\&\, \text{Lastf}(x,t_3,awa,t_4)$ and $M \vDash s=a \,\&\, J(x)$.

Assume for a reductio that $M \vDash \text{MinSet}(x)$.

$\Rightarrow$ from $M \vDash \text{Set}(x)$, $M \vDash \exists t \text{Env}(t,x)$,

$\Rightarrow$ from $M \vDash \text{Lastf}(x,t_3,awa,t_4)$,

$M \vDash \text{Pref}(awa,t_3) \,\&\, \text{Tally}_b(t_4) \,\&\, t_3 < t_4 \,\&$

$\&\, (x=t_3awat_4 \,\vee\, \exists w_3(x=w_3at_3awat_4 \,\&\, \text{Max}^+T_b(t_3,w_3))$.

If $M \vDash x=t_3awat_4$, then $M \vDash t_1auat_2st_3awat_4=x=t_3awat_4$, whence, by (3.6),

$M \vDash t_1auat_2st_3=t_3$, a contradiction because $M \vDash \text{Tally}_b(t_3)$.

$\Rightarrow M \vDash \exists w_3(x=w_3at_3awat_4 \,\&\, \text{Max}^+T_b(t_3,w_3))$.

Let $w_1=t_1auat_2$ and $w_2=t_3awat_4$.

$\Rightarrow$ from $M \vDash \text{MinSet}(x)$, $M \vDash \exists v,t',t''\text{Occ}(w_1,a,w_2,x,t',ava,t'')$.

  (i1) $M \vDash \text{Firstf}(x,t',ava,t'')$.

$\Rightarrow$ from $M \vDash \text{Firstf}(x,t_1,aua,t_2)$, by (5.15), $M \vDash v=u$.

We also have $M \vDash (t'a)Bx \,\&\, (t_1a)Bx$.

$\Rightarrow$ $M \vDash \exists x_1,x_2\ t'ax_1=x=t_1ax_2$,

$\Rightarrow$ since $M \vDash Tally_b(t')$ & $Tally_b(t_1)$, $M \vDash t'=t_1$.

   (i1a) $M \vDash t'=w_1$.

$\Rightarrow$ $M \vDash t'=t_1auat_2$, a contradiction because $M \vDash Tally_b(t')$.

   (i1b) $M \vDash (w_1a)B(t'v)$.

$\Rightarrow$ $M \vDash (t_1auat_2)B(t'v)$,

$\Rightarrow$ $M \vDash \exists y\ (t_1auat_2)y=t'v=t_1aua$,

$\Rightarrow$ by (3.7), $M \vDash at_2y=a$, a contradiction.

   (i1c) $M \vDash w_1a=t'v$.

$\Rightarrow$ $M \vDash t_1auat_2a=t'v=t_1aua$,

$\Rightarrow$ by (3.7), $M \vDash at_2a=a$, a contradiction.

  (i2) $M \vDash \exists w'Intf(x,w',t',ava,t'')$.

$\Rightarrow$ $M \vDash Pref(ava,t')$ & $Tally_b(t'')$ & $t'<t''$ &

                             & $\exists w''\ x=w'at'avat''aw''$ & $Max^+T_b(t',w')$,

$\Rightarrow$ $M \vDash Env(t,x)$ & $Lastf(x,t_3,awa,t_4)$, $M \vDash t''\leq t=t_3=t_4$,

$\Rightarrow$ $M \vDash t'<t''\leq t_3$.

   (i2a) $M \vDash w'at'=w_1$.

$\Rightarrow$ $M \vDash w'at'=t_1auat_2$,

$\Rightarrow$ since $M \vDash Tally_b(t')$ & $Tally_b(t_2)$, by (4.24[b]), $M \vDash t'=t_2$.

We also have $M \vDash t_1auat_2at_3awat_4=x=w'at'avat''aw''$.

$\Rightarrow$ by (3.7), $M \vDash t_3awat_4=vat''aw''$,

$\Rightarrow$ by (4.14[b]), $M \vDash v=t_3v\ t_3Bv$,

$\implies M \vDash t_3 \subseteq_p v,$

$\implies$ from $M \vDash \text{Pref}(ava,t')$, $M \vDash \text{Max}^+T_b(t',v),$

$\implies M \vDash t_3 < t = t_2,$

$\implies M \vDash t' < t_3 < t'$, contradicting $M \vDash t' \in I \subseteq I_0.$

(i2b) $M \vDash \exists v_1(v_1 B(ava) \,\&\, w'at'v_1 = w_1 a).$

$\implies M \vDash w'at'v_1 = t_1 auat_2 a,$

$\implies$ from $M \vDash v_1 B(ava),$

$\implies M \vDash \exists v_2 \, v_1 v_2 = ava,$

$\implies M \vDash t_1 auat_2 at_3 awat_4 = x = w'at'(v_1 v_2)t''aw'',$

$\implies$ by (3.7), $M \vDash t_3 awat_4 = v_2 t''aw'',$

$\implies$ from $M \vDash v_1 v_2 = ava$, $M \vDash v_2 = a \lor aEv_2.$

But $M \vDash v_2 \neq a$ because $M \vDash \text{Tally}_b(t_3).$

$\implies M \vDash aEv_2,$

$\implies M \vDash \exists v_3 \, v_3 a = v_2,$

$\implies M \vDash t_3 awat_4 = (v_3 a)t''aw'',$

$\implies$ by (4.14$^b$), $M \vDash v_3 = t_3 \lor t_3 Bv_3,$

$\implies M \vDash t_3 \subseteq_p v_3 \subseteq_p v_2 \subseteq_p ava,$

$\implies$ from $M \vDash \text{Pref}(ava,t')$, $M \vDash \text{Max}^+T_b(t',v),$

$\implies M \vDash t_3 < t'$, whence a contradiction follows as in (i2a).

(i2c) $M \vDash w'at'ava = w_1 a.$

$\implies M \vDash w'at'ava = t_1 auat_2 a,$

$\implies M \vDash w'at'av = t_1 auat_2,$

$\Rightarrow$ by (4.16), $M \vDash t' \subseteq_p aua$,

$\Rightarrow$ from $M \vDash Firstf(x,t_1,aua,t_2)$, $M \vDash Pref(aua,t_1)$,

$\Rightarrow$ $M \vDash Max^+T_b(t_1,aua)$,

$\Rightarrow$ $M \vDash t'<t_1$,

$\Rightarrow$ from $M \vDash Firstf(x,t_1,aua,t_2)$, by (5.20), $M \vDash t_1 \leq t'$,

$\Rightarrow$ $M \vDash t'<t_1 \leq t'$, contradicting $M \vDash t' \in I \subseteq I_0$.

  (i3) $M \vDash Lastf(x,t',ava,t'')$.

$\Rightarrow$ $M \vDash Pref(ava,t')$ & $Tally_b(t'')$ & $t'=t''$ &

                 & $(x=t'avat'' \vee \exists w'(x=w'at'avat''$ & $Max^+T_b(t',w'))$,

$\Rightarrow$ from $M \vDash Lastf(x,t_3,awa,t_4)$, $M \vDash v=w$ & $t'=t''=t_4=t_3$.

  (i3i) $M \vDash t'avat''=x$.

$\Rightarrow$ $M \vDash t'avat''=x=t_1auat_2at_3awat_4$,

$\Rightarrow$ by (4.16), $M \vDash t_3 \subseteq_p ava$,

$\Rightarrow$ from $M \vDash Pref(ava,t')$, $M \vDash t_3<t'=t_3$, contradicting $M \vDash t_3 \in I \subseteq I_0$.

  (i3ii) $M \vDash \exists w'(x=w'at'avat''$ & $Max^+T_b(t',w'))$.

Exactly analogous to (i2a)-(i2c) with aw" omitted throughout the argument.

We now consider

(iaa) $M \vDash x=t_1auat_2st_3awat_4$

where $M \vDash Set(x)$ & $Firstf(x,t_1,aua,t_2)$ & $Lastf(x,t_3,awa,t_4)$ and $M \vDash s=aa$ along with $M \vDash J(x)$.

Assume that $M \vDash MinSet(x)$.

We proceed as in (i), deriving that $M \vDash \exists w_3(x=w_3at_3awat_4$ & $Max^+T_b(t_3,w_3))$.

Let $w_1 = t_1 auat_2$ and $w_2 = at_3 awat_4$.

$\Rightarrow$ from $M \vDash \text{MinSet}(x)$, $M \vDash \exists v, t', t''\ \text{Occ}(w_1, a, w_2, x, t', ava, t'')$.

The arguments of (i1a)-(i1c) remain unchanged.

  (iaa2) $M \vDash \exists w'\ \text{Intf}(x, w', t', ava, t'')$.

We argue as in (i2) that $M \vDash t' < t_3$.

  (iaa2a) $M \vDash w'at' = w_1$.

$\Rightarrow M \vDash t_1 auat_2 aat_3 awat_4 = x = w'at'avat''aw''$.

$\Rightarrow$ by (3.7), $M \vDash aat_3 awat_4 = avat''aw''$,

$\Rightarrow M \vDash at_3 awat_4 = vat''aw''$,

$\Rightarrow M \vDash v = a \lor aBv$.

If $M \vDash v = a$, then $M \vDash at_3 awat_4 = aat''aw''$, whence $M \vDash t_3 awat_4 = at''aw''$,

a contradiction because $M \vDash \text{Tally}_b(t_3)$.

Therefore $M \vDash aBv'$.

$\Rightarrow M \vDash \exists v'(av' = v\ \&\ at_3 awat_4 = (av')at''aw'')$,

$\Rightarrow M \vDash t_3 awat_4 = v'at''aw''$.

We then continue to derive a contradiction as in (i2a) with $v'$ in place of $v$.

  (iaa2b) $M \vDash \exists v_1(v_1 B(ava)\ \&\ w'at'v_1 = w_1 a)$.

$\Rightarrow$ for $M \vDash v_1 v_2 = ava$, $M \vDash t_1 auat_2 aat_3 awat_4 = w'at'(v_1 v_2)t''aw''$,

$\Rightarrow$ by (3.7), $M \vDash at_3 awat_4 = v_2 t''aw''$,

$\Rightarrow M \vDash v_2 = a \lor aBv_2$.

  (iaa2b1) $M \vDash v_2 = a$.

$\Rightarrow M \vDash at_3 awat_4 = at''aw''$,

$\Rightarrow$ M ⊨ $t_3awat_4=t"aw"$,

$\Rightarrow$ by (3.6), M ⊨ $t_1auat_2aa=w'at'ava$,

$\Rightarrow$ M ⊨ $t_1auat_2a=w'at'av$,

$\Rightarrow$ M ⊨ $v=a \vee aEv$.

    (iaa2b1a) M ⊨ $v=a$.

$\Rightarrow$ M ⊨ $v_1a=v_1v_2=ava=aaa$,

$\Rightarrow$ M ⊨ $v_1=aa$.

But M ⊨ $w'at'v_1=t_1auat_2a$.

$\Rightarrow$ M ⊨ $w'at'aa=t_1auat_2a$,

$\Rightarrow$ M ⊨ $w'at'a=t_1auat_2$, a contradiction.

    (iaa2b1b) M ⊨ $aEv$.

$\Rightarrow$ M ⊨ $\exists v'\ v'a=v$,

$\Rightarrow$ M ⊨ $v_1a=v_1v_2=ava=a(v'a)a$,

$\Rightarrow$ M ⊨ $t_1auat_2a=w'at'a(v'a)$,

$\Rightarrow$ M ⊨ $t_1auat_2=w'at'av'$,

$\Rightarrow$ by (4.16), M ⊨ $t'\subseteq_p aua$,

$\Rightarrow$ from M ⊨ $Pref(aua,t_1)$, M ⊨ $t'<t_1$,

$\Rightarrow$ from M ⊨ $Firstf(x,t_1,aua,t_2)$, by (5.20), M ⊨ $t_1\leq t'$,

$\Rightarrow$ M ⊨ $t'<t_1\leq t'$, contradicting M ⊨ $t'\in I\subseteq I_0$.

    (iaa2b2) M ⊨ $aBv_2$.

$\Rightarrow$ M ⊨ $\exists v_3\ (av_3=v_2\ \&\ at_3awat_4=(av_3)t"aw")$.

We have M ⊨ $t_1auat_2a=w'at'v_1$.

$\Rightarrow$ $M \vDash v_1=a \lor aEv_1$.

(iaa2b2a) $M \vDash v_1=a$.

$\Rightarrow$ $M \vDash w'at'a=w'at'v_1=t_1auat_2a$,

$\Rightarrow$ $M \vDash w'at'=t_1auat_2$,

$\Rightarrow$ from $M \vDash t_1auat_2aat_3awat_4=w'at'avat''aw''$, by (3.7),

$M \vDash aat_3awat_4=avat''aw''$,

$\Rightarrow$ $M \vDash at_3awat_4=vat''aw''$,

$\Rightarrow$ $M \vDash v=a \lor aBv$.

If $M \vDash v=a$, then $M \vDash at_3awat_4=aat''aw''$, whence $M \vDash t_3awat_4=at''aw''$, a contradiction.

$\Rightarrow$ $M \vDash aBv$,

$\Rightarrow$ $M \vDash \exists v'(av'=v \ \& \ at_3awat_4=(av')at''aw'')$,

$\Rightarrow$ $M \vDash t_3awat_4=v'at''aw''$.

We then derive a contradiction as in (i2a) with v' in place of v.

(iaa2b2b) $M \vDash aEv_1$.

$\Rightarrow$ $M \vDash \exists v_2 \ v_3a=v_1$,

$\Rightarrow$ $M \vDash w'at'(v_2a)=t_1auat_2a$,

$\Rightarrow$ $M \vDash w'at'v_2=t_1auat_2$.

But $M \vDash v_2a=v_1 \ \& \ v_1av_3=ava$.

$\Rightarrow$ $M \vDash (v_2a)av_3=ava$,

$\Rightarrow$ $M \vDash v_2=a \lor aBv_2$.

Now, $M \vDash v_2 \neq a$ because $M \vDash Tally_b(t_2)$.

$\Rightarrow$ $M \vDash aBv_2$,

$\Rightarrow$ $M \vDash \exists v_4\ av_4 = v_2$,

$\Rightarrow$ $M \vDash w'at'(av_4) = t_1auat_2$.

We now derive a contradiction as in (iaa2b1b).

   (iaa2c)  $M \vDash w'at'v = w_1a$.

Exactly as in (i2c).

  (iaa3)  $M \vDash Lastf(x,w',t',ava,t'')$.

The same argument as in (i3).

This completes the proof that $M \vDash \neg MinSet(x)$ under (i).

(ii) Assume $M \vDash x = t_1auat_2st_3awat_4$

where $M \vDash Set(x)\ \&\ Firstf(x,t_1,aua,t_2)\ \&\ Intf(x,w_3,t_3,awa,t_4)$ and

$M \vDash s = a\ \&\ J(x)$.

Assume for a reductio that $M \vDash MinSet(x)$.

$\Rightarrow$ $M \vDash Pref(awa,t_3)\ \&\ Tally_b(t_4)\ \&\ t_3 < t_4\ \&$

                                     $\&\ x = w_3at_3awat_4aw_4\ \&\ Max^+T_b(t_3,w_3)$.

Let $w_1 = t_1auat_2$ and $w_2 = t_3awat_4aw_4$.

$\Rightarrow$ from $M \vDash MinSet(x)$, $M \vDash \exists v,t',t''\ Occ(w_1,a,w_2,x,t',ava,t'')$.

  (ii1)  $M \vDash Firstf(x,t',ava,t'')$.

The same argument as in (i1).

  (ii2)  $M \vDash \exists w'Intf(x,w',t',ava,t'')$.

The same argument applies as in (i2) with $aw_4$ appended throughout.

  (ii3)  $M \vDash Lastf(x,t',ava,t'')$.

$\Rightarrow$ $M \vDash \text{Pref}(ava,t') \& \text{Tally}_b(t'') \& t'=t'' \&$

$\& (x=t'avat'' \vee \exists w'(x=w'at'avat'' \& \text{Max}^+T_b(t',w')))$.

 (ii3i)  $M \vDash t'avat''=x$.

$\Rightarrow$ $M \vDash \text{Firstf}(x,t',ava,t'')$,

$\Rightarrow$ from $M \vDash \text{Firstf}(x,t_1,aua,t_2)$, by (5.15), $M \vDash v=u$,

$\Rightarrow$ from $M \vDash (t'a)Bx \& (t_1a)Bx \& \text{Tally}_b(t') \& \text{Tally}_b(t'')$, by (4.23$^b$), $M \vDash t'=t_1$,

$\Rightarrow$ $M \vDash t'avat''=x=t_1auat_2at_3awat_4aw_4=t'avat_2at_3awat_4aw_4$,

$\Rightarrow$ by (3.7), $M \vDash t''=t_2at_3awat_4aw_4$, a contradiction because $M \vDash \text{Tally}_b(t'')$.

 (ii3ii)  $M \vDash \exists w'(x=w'at'avat'' \& \text{Max}^+T_b(t',w'))$.

The same argument applies as in (i2a)-(i2c) with $aw''$ omitted from $t'avat''aw''$ and $aw_4$ appended to $t_1auat_2at_3awat_4$.

We now consider

(iiaa)  $M \vDash t_1auat_2st_3awat_4aw_4=x$

where $M \vDash \text{Set}(x) \& \text{Firstf}(x,t_1,aua,t_2) \& \text{Intf}(x,w_3,t_3,awa,t_4)$ along with $M \vDash s=aa \& J(x)$.

Assume that $M \vDash \text{MinSet}(x)$.

We proceed as in (ii), deriving $M \vDash x=w_3at_3awat_4aw_4 \& \text{Max}^+T_b(t_3,w_3)$.

Let $w_1=t_1auat_2$ and $w_2=t_3awat_4aw_4$.

$\Rightarrow$ from $M \vDash \text{MinSet}(x)$, $M \vDash \exists v,t',t''\ \text{Occ}(w_1,a,w_2,x,t',ava,t'')$.

We adapt the arguments of (ii) just as we adapted the arguments of (i) for (iaa).

This completes the proof of $M \vDash \neg\text{MinSet}(x)$ under hypothesis (ii), and the proof of (10.12).

(10.13) For any string concept $I\subseteq I_0$ there is a string concept $J\subseteq I$ such that

$QT^+ \vdash \forall x,z \in J \forall t',t_1,t_2,u,v,x_1[x=t_1auat_2ax_1$ & $Firstf(x,t_1,aua,t_2)$ & $z=t'avat_2ax_1$ &

& $t'<t_2$ & $\neg(v \,\varepsilon\, x)$ & $Pref(ava,t')$ & $MinSet(x) \rightarrow MinSet(z)]$.

Let $J \equiv I_{5.35}$ & $I_{10.12}$.

Assume that $M \vDash x=t_1auat_2ax_1$ & $z=t'avat_2ax_1$

where $M \vDash Firstf(x,t_1,aua,t_2)$ & $t'<t_2$ & $\neg(v \,\varepsilon\, x)$ & $Pref(ava,t')$.

Assume also that $M \vDash MinSet(x)$.

$\Rightarrow$ by (5.35), $M \vDash \exists t\, Env(t,z)$ & $Firstf(z,t',ava,t_2)$.

Assume that $M \vDash w_1aw_2=z=t'avat_2ax_1$.

$\Rightarrow$ from the proof of (5.35), $M \vDash Env(t,x)$.

We have that $M \vDash w_1Bz$ & $(t'av)Bz$.

$\Rightarrow$ by (3.8), $M \vDash w_1B(t'av) \vee w_1=t'av \vee (t'av)Bw_1$.

(1) $M \vDash w_1B(t'av)$.

$\Rightarrow M \vDash \exists w\, w_1w=t'av$,

$\Rightarrow M \vDash (w_1w)at_2ax_1=t'avat_2ax_1=z=w_1aw_2$,

$\Rightarrow$ by (3.7), $M \vDash wat_2ax_1=aw_2$,

$\Rightarrow M \vDash w=a \vee aBw$,

$\Rightarrow M \vDash aat_2ax_1=aw_2 \vee \exists w_3\, (aw_3)at_2ax_1=aw_2$,

$\Rightarrow M \vDash at_2ax_1=w_2 \vee w_3at_2ax_1=w_2$,

$\Rightarrow$ from $M \vDash w_1aw_2=z$ by (3.6), $M \vDash w_1a=t'av \vee w_1aw_3=t'av$,

$\Rightarrow M \vDash (w_1a)B(t'ava)$ & $Firstf(z,t',ava,t_2)$ & $w_1aw_2=z$,

$\implies$ M ⊨ Occ($w_1$,a,$w_2$,z,t',ava,$t_2$).

(2) M ⊨ $w_1$=t'av.

$\implies$ M ⊨ $w_1$a=t'(ava) & Firstf(z,t',ava,$t_2$) & $w_1$a$w_2$=z,

$\implies$ M ⊨ Occ($w_1$,a,$w_2$,z,t',ava,$t_2$).

(3) M ⊨ (t'av)B$w_1$.

$\implies$ M ⊨ ∃w t'avw=$w_1$,

$\implies$ M ⊨ z=$w_1$a$w_2$=(t'avw)a$w_2$=t'av$t_2$a$x_1$,

$\implies$ by (3.7), M ⊨ wa$w_2$=a$t_2$a$x_1$,

$\implies$ M ⊨ w=a v aBw.

(3a) M ⊨ w=a.

$\implies$ M ⊨ aa$w_2$=a$t_2$a$x_1$,

$\implies$ M ⊨ a$w_2$=$t_2$a$x_1$, a contradiction because M ⊨ Tally$_b$($t_2$).

(3b) M ⊨ aBw.

$\implies$ M ⊨ ∃$w_0$ (a$w_0$=w & (a$w_0$)a$w_2$=a$t_2$a$x_1$),

$\implies$ M ⊨ $w_0$a$w_2$=$t_2$a$x_1$,

$\implies$ by (4.14$^b$), M ⊨ $w_0$=$t_2$ v $t_2$B$w_0$.

(3bi) M ⊨ $w_0$=$t_2$.

$\implies$ M ⊨ $t_2$a$w_2$=$t_2$a$x_1$,

$\implies$ M ⊨ $w_1$=t'av(a$w_0$)=t'av$t_2$,

$\implies$ M ⊨ $t_1$aua$t_2$a$x_1$=($t_1$aua$t_2$)a$w_2$.

By hypothesis, M ⊨ MinSet(x).

Let w'=$t_1$aua$t_2$.

$\Rightarrow M \vDash x = w'aw_2$,

$\Rightarrow$ from $M \vDash \text{MinSet}(x)$, $M \vDash \exists w, t_3, t_4\ \text{Occ}(w', a, w_2, x, t_3, awa, t_4)$.

   (3bi1)  $M \vDash \text{Firstf}(x, t_3, awa, t_4)$.

$\Rightarrow$ from $M \vDash \text{Firstf}(x, t_1, aua, t_2)$, by (5.15) and (4.23$^b$), $M \vDash u=w\ \&\ t_3=t_1$,

$\Rightarrow M \vDash \text{Pref}(awa, t_3)\ \&\ \text{Tally}_b(t_4)\ \&$

$$\&\ ((t_3=t_4\ \&\ x=t_3awat_4)\ \vee\ (t_3<t_4\ \&\ (t_3awat_4a)Bx)).$$

Suppose that $M \vDash t_3=t_4\ \&\ x=t_3awat_4$.

$\Rightarrow M \vDash t_3awat_4 = x = t_1auat_2ax_1$,

$\Rightarrow$ from $M \vDash \text{Env}(t, x)$, $M \vDash \text{MaxT}_b(t, x)$,

$\Rightarrow M \vDash t_2 \leq t = t_3$,

$\Rightarrow M \vDash t_1 < t_2 \leq t_3 = t_1$, contradicting $M \vDash t_1 \in I \subseteq I_0$.

Therefore $M \vDash \neg(t_3=t_4\ \&\ x=t_3awat_4)$.

$\Rightarrow M \vDash t_3<t_4\ \&\ (t_3awat_4a)Bx$,

$\Rightarrow M \vDash \exists x_2\ t_3awat_4ax_2 = x = t_1auat_2ax_1$,

$\Rightarrow$ since $M \vDash u=w\ \&\ t_3=t_1$, by (3.7), $M \vDash t_4ax_2 = t_2ax_1$,

$\Rightarrow$ since $M \vDash \text{Tally}_b(t_4)\ \&\ \text{Tally}_b(t_2)$, by (4.23$^b$), $M \vDash t_4=t_2$,

$\Rightarrow M \vDash t_3awat_4 = t_1auat_2$.

   (3bi1a)  $M \vDash t_3 = w'$.

$\Rightarrow M \vDash t_3 = t_1auat_2$, a contradiction because $M \vDash \text{Tally}_b(t_3)$.

   (3bi1b)  $M \vDash (w'a)B(t_3w)$.

$\Rightarrow M \vDash (t_1auat_2a)B(t_3aua)$,

$\Rightarrow M \vDash \exists w_3\ t_1auat_2aw_3 = t_3aua$,

$\Rightarrow$ since $M \vDash t_1=t_3$, by (3.7), $M \vDash at_2aw_3=a$, a contradiction.

(3bi1c) $M \vDash w'a=t_3w$.

$\Rightarrow M \vDash t_1auat_2a=t_3aua$, whence a contradiction follows as in (3bi1b).

(3bi2) $M \vDash \exists w_3\, Intf(x,w_3,t_3,awa,t_4)$.

$\Rightarrow M \vDash Pref(awa,t_3)\ \&\ Tally_b(t_4)\ \&\ t_3<t_4\ \&$
$$\&\ \exists w_4\ (x=w_3at_3awat_4aw_4\ \&\ Max^+T_b(t_3,w_3)).$$

(3bi2a) $M \vDash w_3at_3=w'$.

$\Rightarrow M \vDash w_3at_3=w'=t_1auat_2$,

$\Rightarrow$ by (4.24$^b$), $M \vDash t_3=t_2$,

$\Rightarrow M \vDash w_3at_3awat_4aw_4=x=t_1auat_2ax_1=w'ax_1=w'aw_2$,

$\Rightarrow M \vDash w'(awat_4aw_4)=x=w'aw_2$,

$\Rightarrow$ by (3.7), $M \vDash wat_4aw_4=w_2$,

$\Rightarrow M \vDash z=w_1aw_2=w_1awat_4aw_4$.

But $M \vDash w_1=t'avat_2$.

$\Rightarrow M \vDash w_1aw_2=z=(t'avat_2)awat_4aw_4=t'avat_3awat_4aw_4$.

Let $w_5=t'av$.

Then, as in the proof of (5.31), (1b), with the roles of x and z reversed, we have that $M \vDash Max^+T_b(t_3,w_5)$.

$\Rightarrow M \vDash Pref(awa,t_3)\ \&\ Tally_b(t_4)\ \&\ t_3<t_4\ \&\ z=w_5at_3awat_4aw_4\ \&\ Max^+T_b(t_3,w_5)$,

$\Rightarrow M \vDash Intf(z,w_5,t_3,awa,t_4)\ \&\ w_5at_3=w_1\ \&\ w_1aw_2=z$,

$\Rightarrow M \vDash Occ(w_1,a,w_2,z,t_3,awa,t_4)$.

(3bi2b) $M \vDash \exists v_1(v_1Bw \ \& \ w_3at_3v_1=w'a)$.

$\Longrightarrow M \vDash v_1=a \ v \ aEv_1$.

(3bi2bi) $M \vDash v_1=a$.

$\Longrightarrow M \vDash w_3at_3a=w'a$,

$\Longrightarrow M \vDash w_3at_3=w'=t_1auat_2$.

We then proceed exactly as in (3bi2a).

(3bi2bii) $M \vDash aEv_1$.

$\Longrightarrow M \vDash \exists v_2 \ v_2a=v_1$,

$\Longrightarrow M \vDash w_3at_3(v_2a)=w'a$,

$\Longrightarrow M \vDash w_3at_3v_2=w'=t_1auat_2$.

We have $M \vDash Tally_b(v_2) \ v \ \neg Tally_b(v_2)$.

(3bi2bii1) $M \vDash Tally_b(v_2)$.

$\Longrightarrow$ by (4.5), $M \vDash Tally_b(t_3v_2)$,

$\Longrightarrow$ by (4.24$^b$), $M \vDash t_3v_2=t_2$,

$\Longrightarrow M \vDash t_3<t_2$,

$\Longrightarrow$ by (5.34), $M \vDash t_2 \leq t_3$,

$\Longrightarrow M \vDash t_3<t_2 \leq t_3$, contradicting $M \vDash t_3 \in I \subseteq I_0$.

(3bi2bii2) $M \vDash \neg Tally_b(v_2)$.

$\Longrightarrow M \vDash a \subseteq_p v_2$,

$\Longrightarrow$ since $M \vDash Tally_b(t_2)$, $M \vDash \neg(v_2=a \ v \ aEv_2)$,

$\Longrightarrow M \vDash aBv_2 \ v \ \exists v_3,v_4 \ v_3av_4=v_2$.

(3bi2biia) $M \vDash aBv_2$.

$\Rightarrow M \vDash \exists v_3\, av_3 = v_2$,

$\Rightarrow M \vDash w_3 at_3(av_3) = w' = t_1 auat_2$,

$\Rightarrow$ by (4.16), $M \vDash t_3 \subseteq_p aua$.

But from $M \vDash \text{Pref}(aua, t_1)$, $M \vDash \text{Max}^+ T_b(t_1, aua)$.

$\Rightarrow M \vDash t_3 < t_1$,

$\Rightarrow$ from (5.34), $M \vDash t_2 \leq t_3$,

$\Rightarrow M \vDash t_3 < t_1 < t_2 \leq t_3$, again contradicting $M \vDash t_3 \in I \subseteq I_0$.

(3bi2biib) $M \vDash \exists v_3, v_4\, v_3 av_4 = v_2$.

$\Rightarrow M \vDash w_3 at_3(v_3 av_4) = w' = t_1 auat_2$,

$\Rightarrow$ by (4.16), $M \vDash t_3 \subseteq_p t_3 v_3 \subseteq_p aua$, and we proceed as in (3bi2bii2a).

(3bi2c) $M \vDash w_3 at_3 awa = w'a$.

$\Rightarrow M \vDash w_3 at_3 aw = w' = t_1 auat_2$,

$\Rightarrow$ by (4.16), $M \vDash t_3 \subseteq_p aua$, and we proceed exactly as in (3bi2bii2a).

(3bi3) $M \vDash \text{Lastf}(x, t_3, awa, t_4)$.

$\Rightarrow M \vDash \text{Pref}(awa, t_3)\ \&\ \text{Tally}_b(t_4)\ \&\ t_3 = t_4\ \&\ (x = t_3 awat_4\ v$

$v\ \exists w_3\ (x = w_3 at_3 awat_4\ \&\ \text{Max}^+ T_b(t_3, w_3)))$.

We have that $M \vDash \neg(x = t_3 awat_4)$ as in (3bi1).

$\Rightarrow M \vDash \exists w_3\ (x = w_3 at_3 awat_4\ \&\ \text{Max}^+ T_b(t_3, w_3))$.

From this point on the argument is exactly analogous to (3bi2) with $aw_4$ omitted throughout.

(3bii)  $M \vDash t_2Bw_0$.

$\implies$ $M \vDash \exists w_4\ t_2w_4 = w_0$,

$\implies$ $M \vDash t'ava(t_2w_4)aw_2 = w_1aw_2 = z = t'avat_2ax_1$,

$\implies$ by (3.6), $M \vDash w_1 = t'avat_2w_4$,

$\implies$ from $M \vDash w_0aw_2 = t_2ax_1$, $M \vDash (t_2w_4)aw_2 = t_2ax_1$,

$\implies$ by (3.7), $M \vDash w_4aw_2 = ax_1$,

$\implies$ $M \vDash (t_1auat_2)ax_1 = x = (t_1auat_2)w_4aw_2$.

By hypothesis, we have that $M \vDash \mathrm{MinSet}(x)$.

Let $w' = t_1auat_2w_4$.

$\implies$ $M \vDash x = w'aw_2$,

$\implies$ $M \vDash \exists w, t_3, t_4\ \mathrm{Occ}(w', a, w_2, x, t_3, awa, t_4)$.

(3bii1)  $M \vDash \mathrm{Firstf}(x, t_3, awa, t_4)$.

We proceed as in (3bi1) with $t_2w_4$ replacing $t_2$ throughout.

(3bii2)  $M \vDash \exists w_5\ \mathrm{Intf}(x, w_5, t_3, awa, t_4)$.

$\implies$ $M \vDash \mathrm{Pref}(awa, t_5)\ \&\ \mathrm{Tally}_b(t_4)\ \&\ t_3 < t_4\ \&$

$\quad\quad\quad\quad\quad\quad \&\ \exists w_6\ (x = w_5at_3awat_4aw_6\ \&\ \mathrm{Max}^+T_b(t_3, w_5))$.

(3bii2a)  $M \vDash w_5at_3 = w'$.

$\implies$ $M \vDash w_5at_3 = w' = t_1auat_2w_4$,

$\implies$ $M \vDash w'aw_2 = (t_1auat_2w_4)aw_2 = x = t_1auat_2ax_1 = w_5at_3awat_4aw_6$,

$\implies$ $M \vDash x = (t_1auat_2w_4)awat_4aw_6$,

$\implies$ by (3.7), $M \vDash w_4awat_4aw_6 = ax_1\ \&\ w_4aw_2 = ax_1$ and $M \vDash aw_2 = awat_4aw_6$,

$\implies$ $M \vDash w_4 = a \lor aBw_4$,

$\Rightarrow$ from $M \vDash \text{Tally}_b(t_3)$, $M \vDash \neg(w_4=a)$,

$\Rightarrow M \vDash aBw_4$,

$\Rightarrow M \vDash \exists w_7\ aw_7=w_4$,

$\Rightarrow M \vDash w_1=t'avat_2(aw_7)$,

$\Rightarrow M \vDash (aw_7)awat_4aw_6=ax_1$,

$\Rightarrow M \vDash w_5at_3=w'=t_1auat_2aw_7$,

$\Rightarrow$ by (4.15$^b$), $M \vDash w_7=t_3\ \lor\ t_3Ew_7$,

$\Rightarrow M \vDash w_1=t'avat_2at_3\ \lor\ \exists w_8(w_8t_3=w_7\ \&\ w_1=t'avat_2aw_8t_3)$,

$\Rightarrow M \vDash w_5at_3=w'=t_1auat_2at_3\ \lor\ \exists w_8(w_5at_3=w'=t_1auat_2aw_8t_3\ \&\ w_8t_3=w_7)$,

$\Rightarrow$ by (3.6), $M \vDash w_5=t_1auat_2\ \lor\ w_5a=t_1auat_2aw_8$.

If $M \vDash w_5a=t_1auat_2aw_8$, then $M \vDash w_8=a\ \lor\ aEw_8$.

$\Rightarrow M \vDash z=w_1aw_2=t'avat_2(aw_7)awat_4aw_6$,

$\Rightarrow M \vDash z=t'avat_2at_3awat_4aw_6\ \lor\ z=t'avat_2a(w_8t_3)awat_4aw_6$,

$\Rightarrow$ since $M \vDash w_8=a\ \lor\ aEw_8$,

$\quad M \vDash (w_1=t'avat_2at_3\ \&\ z=t'avat_2at_3awat_4aw_6)\ \lor$

$\qquad\qquad \lor\ (w_1=t'avat_2aat_3\ \&\ z=t'avat_2aat_3awat_4aw_6)\ \lor$

$\qquad\qquad \lor\ \exists w_9(w_8=w_9a\ \&\ w_1=t'avat_2aw_9at_3\ \&\ z=t'avat_2aw_9at_3awat_4aw_6)$.

We now claim that $M \vDash \text{Max}^+T_b(t_3,t'avat_2aw_9)$.

Assume $M \vDash \text{Tally}_b(t'')\ \&\ t''\subseteq_p t'avat_2aw_9$.

$\Rightarrow$ by (4.17$^b$), $M \vDash t''\subseteq_p t'avat_2\ \lor\ t''\subseteq_p w_9$.

If $M \vDash t''\subseteq_p t'avat_2$, then, by (4.17$^b$), $M \vDash t''\subseteq_p t'av\ \lor\ t''\subseteq_p t_2$.

If $M \vDash t''\subseteq_p t'av$, then $M \vDash t''\subseteq_p t'\ \lor\ t''\subseteq_p v\subseteq_p ava$, and from $M \vDash \text{Pref}(ava,t')$

we have $M \vDash t''\leq t'<t_2$.

Hence, if $M \vDash t''\subseteq_p t'avat_2$, then $M \vDash t''\leq t_2$.

If $M \vDash t''\subseteq_p w_9$, note that from $M \vDash w_5=t_1auat_2 \vee w_5a=t_1auat_2a(w_9a)$, we have

$$M \vDash w_5=t_1auat_2 \vee w_5=t_1auat_2aw_9.$$

Hence either way $M \vDash t_2\subseteq_p w_5$, and if $M \vDash w_5=t_1auat_2aw_9$, then $M \vDash w_9\subseteq_p w_5$.

It follows that if $M \vDash t''\subseteq_p t'avat_2aw_9$, then $M \vDash t''\subseteq_p w_5$, so $M \vDash t''<t_3$ from

$M \vDash \text{Max}^+T_b(t_3,w_5)$.

Therefore, $M \vDash \text{Max}^+T_b(t_3,t'avat_2aw_9)$, and a fortiori, $M \vDash \text{Max}^+T_b(t_3,z_1)$,

where $M \vDash z_1=t'avat_2 \vee z_1=t'avat_2a \vee z_1=t_1auat_2aw_9$.

Hence we have

$$M \vDash \text{Pref}(awa,t_3) \mathbin{\&} \text{Tally}_b(t_4) \mathbin{\&} t_3<t_4 \mathbin{\&}$$
$$\mathbin{\&} \exists z_1 (z=z_1at_3awat_4aw_6 \mathbin{\&} \text{Max}^+T_b(t_3,z_1)),$$

that is, $M \vDash \text{Intf}(z,z_1,t_3,awa,t_4)$, along with $M \vDash w_1=z_1at_3 \mathbin{\&} w_1aw_2=z$.

But then $M \vDash \text{Occ}(w_1,a,w_2,z,t_3,awa,t_4)$, as required.

  (3bii2b) $M \vDash \exists w''(w''B(awa) \mathbin{\&} w_5at_3w''=w'a)$.

$\implies M \vDash w_5at_3w''=w'a=t_1auat_2w_4a$,

$\implies M \vDash w'aw_2=x=t_1auat_2w_4aw_2=t_1auat_2ax_1=w_5at_3awat_4aw_6$,

$\implies$ from $M \vDash w''B(awa)$, $M \vDash \exists x' \, w''x'=awa$,

$\implies M \vDash x=w_5at_3(w''x')t_4aw_6=(t_1auat_2w_4a)x't_4aw_6$,

$\implies$ by (3.6), $M \vDash w_5at_3w''=t_1auat_2w_4a$,

$\implies M \vDash w''=a \vee aEw''$.

(3bii2bi) $M \vDash w''=a$.

$\implies M \vDash w_5at_3a=w'a$,

$\implies M \vDash w_5at_3=w'=t_1auat_2w_4$.

We then proceed exactly as in (3bii2a).

(3bii2bii) $M \vDash aEw''$.

$\implies M \vDash \exists x_3\ x_3a=w''$,

$\implies M \vDash w_5at_3(x_3a)=t_1auat_2w_4a$,

$\implies M \vDash w_5at_3x_3=t_1auat_2w_4$,

$\implies$ from $M \vDash w_4aw_2=ax_1$, $M \vDash w_4=a \lor aBw_4$,

$\implies M \vDash (w_4=a\ \&\ w_5at_3x_3=t_1auat_2a) \lor \exists w_7\ (w_4=aw_7\ \&\ t_1auat_2aw_7=w_5at_3x_3)$.

(3bii2bii1) $M \vDash w_4=a\ \&\ w_5at_3x_3=t_1auat_2a$.

$\implies M \vDash x_3=a \lor aEx_3$.

(3bii2bii1a) $M \vDash x_3=a$.

$\implies M \vDash w_5at_3a=t_1auat_2a$,

$\implies M \vDash w_5at_3=w'=t_1auat_2$,

$\implies$ by (4.24$^b$), $M \vDash t_2=t_3$,

$\implies M \vDash w_5at_2=w'=t_1auat_2$,

$\implies$ by (3.6), $M \vDash w_5=t_1au$,

$\implies$ by (3.7), $M \vDash at_2ax_1=at_3awat_4aw_6$,

$\implies M \vDash z=(t'av)at_2ax_1=(t'av)at_3awat_4aw_6$.

We then reason as in (3bi2a) that $M \vDash Max^+T_b(t_3,t'av)$, and derive that

$$M \vDash Intf(z,t'av,t_3,awa,t_4),$$

along with $\ M \vDash (t'av)at_3aa=(t'av)at_2w_4a=w_1a\ $ and $\ M \vDash w_1aw_2=z$,

so that $\ M \vDash (t'av)at_3v_1=w_1a\ $ for $\ v_1=aa$.

But $\ M \vDash w''=x_3a=aa$, and by hypothesis $\ M \vDash w''Bw$.

Hence also $\ M \vDash v_1Bw$.

Therefore $\ M \vDash Occ(w_1,a,w_2,z,t_3,awa,t_4)$, as required.

   (3bii2bii1b) $\ M \vDash aEx_3$.

$\implies M \vDash \exists x_4\ x_4a=x_3$,

$\implies M \vDash w_5at_3x_4a=w'=t_1auat_2a$,

$\implies M \vDash w_5at_3x_4=t_1auat_2$,

$\implies$ by (4.15$^b$), $M \vDash t_3x_4=t_2\ v\ t_2E(t_3x_4)$.

By an argument analogous to (3bi2bii1) we have that $\ M \vDash \neg Tally_b(x_4)$.

$\implies M \vDash \neg Tally_b(t_3x_4)$,

$\implies M \vDash t_3x_4 \neq t_2$,

$\implies M \vDash t_2E(t_3x_4)$,

$\implies M \vDash \exists x_5\ t_3x_4=x_5t_2$,

$\implies M \vDash w_5ax_5t_2=t_1auat_2$,

$\implies$ by (3.6), $M \vDash w_5ax_5=t_1aua$,

$\implies M \vDash x_5=a\ v\ aEx_5$,

$\implies$ since $\ M \vDash \neg Tally_b(t_3)$, $\ M \vDash x_5 \neq a$,

$\implies M \vDash aEx_5$,

$\implies M \vDash \exists x_6\ x_6a=x_5$,

$\implies M \vDash w_5a(x_6a)=t_1aua$,

$\Rightarrow$ $M \vDash w_5ax_6=t_1au$,

$\Rightarrow$ from $M \vDash t_1auat_2ax_1=x=w_5at_3awat_4aw_6$,

$\qquad\qquad M \vDash (w_5ax_6a)t_2ax_1=x=w_5at_3awat_4aw_6$,

$\Rightarrow$ by (3.7), $M \vDash x_6at_2ax_1=t_3awat_4aw_6$,

$\Rightarrow$ by (4.14$^b$), $M \vDash x_6=t_3 \vee t_3Bx_6$,

$\Rightarrow$ $M \vDash t_3 \subseteq_p x_6$,

$\Rightarrow$ from $M \vDash w_5ax_6=t_1au$, by (4.14$^b$), $M \vDash w_5=t_1 \vee t_1Bw_5$,

$\Rightarrow$ $M \vDash t_1ax_6=t_1au \vee \exists x_7(w_5=t_1x_7 \,\&\, t_1au=(t_1x_7)ax_6)$,

$\Rightarrow$ by (3.7), $M \vDash x_6=u \vee au=x_7ax_6$,

$\Rightarrow$ $M \vDash x_6=u \vee au=aax_6 \vee \exists x_8\, au=(ax_8)ax_6$,

$\Rightarrow$ $M \vDash x_6=u \vee u=ax_6 \vee u=x_8ax_6$,

$\Rightarrow$ from $M \vDash t_3\subseteq_p x_6$, $M \vDash t_3\subseteq_p u$,

$\Rightarrow$ from $M \vDash \text{Intf}(x,w_5,t_3,awa,t_4)$, by (5.34), $M \vDash t_1<t_2\leq t_3$,

contradicting $M \vDash \text{Pref}(aua,t_1)$.

$\qquad$ (3bii2bii2) $M \vDash \exists w_7(w_4=aw_7 \,\&\, t_1auat_2aw_7=w_5at_3x_3)$.

$\Rightarrow$ $M \vDash (t_1auat_2)B(w_5at_3x_3) \,\&\, (w_5at_3x_3)B(w_5at_3x_3)$,

$\Rightarrow$ by (3.8), $M \vDash (t_1auat_2)B(w_5at_3) \vee t_1auat_2=w_5at_3 \vee (w_5at_3)B(t_1auat_2)$.

$\qquad$ (3bii2bii2a) $M \vDash (w_5at_3)B(t_1auat_2)$.

$\Rightarrow$ $M \vDash \exists w_8\, t_1auat_2=w_5at_3w_8$,

$\Rightarrow$ by (4.15$^b$), $M \vDash t_3w_8=t_2 \vee t_2E(t_3w_8)$.

If $M \vDash t_3w_8=t_2$, then $M \vDash t_3<t_2$, whereas by (5.34), $M \vDash t_2\leq t_3$, hence $M \vDash t_3<t_2\leq t_3$, contradicting $M \vDash t_3 \in I \subseteq I_0$.

Therefore $M \vDash t_2E(t_3w_8)$.

$\implies M \vDash \exists w_9\ w_9t_2=t_3w_8$,

$\implies M \vDash t_1auat_2=w_5aw_9t_2$,

$\implies$ by (3.6), $M \vDash t_1aua=w_5aw_9$.

We then derive a contradiction analogously to (3bii2bii1b) with $w_9$ replacing $x_5$ throughout the argument.

Therefore $M \vDash \neg(w_5at_3)B(t_1auat_2)$.

   (3bii2bii2b)   $M \vDash t_1auat_2=w_5at_3$.

$\implies$ by (4.15$^b$), $M \vDash t_2=t_3$,

$\implies M \vDash t_1auat_2=w_5at_2$,

$\implies$ by (3.6), $M \vDash t_1au=w_5$.

We then derive, as in (3bii2bii1a), that $M \vDash z=t'avat_2ax_1=(t'av)at_3awat_4aw_6$, and that $M \vDash \text{Intf}(z,t'av,t_3,awa,t_4)$.

Since also $M \vDash (t'av)at_3aw_7a=(t'av)at_2w_4a=w_1a\ \&\ w_1aw_2=z$,

we have that $M \vDash t'avat_3v_1a=w_1a$ for $v_1=aw_7a=w_4a$.

$\implies$ from $M \vDash t_1auat_2w_4a=w_5at_3w''$, $M \vDash w''=w_4a$,

$\implies$ since by $M \vDash w''B(awa)$ by hypothesis, $M \vDash v_1B(awa)$.

Therefore $M \vDash \text{Occ}(w_1,a,w_2,z,t_3,awa,t_4)$, as required.

   (3bii2bii2c)   $M \vDash (t_1auat_2)B(w_5at_3)$.

$\implies M \vDash \exists x_7\ t_1auat_2x_7=w_5at_3$,

$\implies$ from $M \vDash t_1auat_2aw_7=w_5at_3x_3$, $M \vDash t_1auat_2aw_7=t_1auat_2x_7x_3$,

$\implies$ by (3.7), $M \vDash aw_7=x_7x_3$, $M \vDash x_7=a \lor aBx_7$.

Now, $M \vDash x_7 \neq a$ because $M \vDash \text{Tally}_b(t_3)$.

$\implies M \vDash aBx_7$,

$\implies M \vDash \exists x_9\, ax_9 = x_7$,

$\implies M \vDash t_1 auat_2(ax_9) = w_5 at_3$,

$\implies$ by (4.15$^b$), $M \vDash x_9 = t_3\ \vee\ t_3 E x_9$,

$\implies M \vDash t_1 auat_2 at_3 = w_5 at_3\ \vee\ \exists y_1(y_1 t_3 = x_9\ \&\ t_1 auat_2 ay_1 t_3 = w_5 at_3)$,

$\implies$ by (3.6), $M \vDash t_1 auat_2 = w_5\ \vee\ t_1 auat_2 ay_1 = w_5 a$,

$\implies M \vDash t_1 auat_2 = w_5\ \vee\ t_1 auat_2 aa = w_5 a\ \vee\ \exists y_2\, t_1 auat_2 a(y_2 a) = w_5 a$,

$\implies M \vDash t_1 auat_2 = w_5\ \vee\ t_1 auat_2 a = w_5\ \vee\ t_1 auat_2 ay_2 = w_5$.

Now, if $M \vDash t_1 auat_2 = w_5\ \vee\ t_1 auat_2 a = w_5$, then

$\quad M \vDash (t_1 auat_2)at_3 awat_4 aw_6 = w_5 at_3 awat_4 aw_6 = x\ \vee$

$\qquad\qquad\qquad \vee\ (t_1 auat_2 a)at_3 awat_4 aw_6 = w_5 at_3 awat_4 aw_6 = x$.

$\implies$ by (10.12), $M \vDash \neg \text{MinSet}(x)$, contradicting the hypothesis.

$\implies M \vDash \neg(t_1 auat_2 = w_5\ \vee\ t_1 auat_2 a = w_5)$,

$\implies M \vDash t_1 auat_2 ay_2 = w_5$,

$\implies M \vDash t_1 auat_2 ay_2 (at_3 x_3) = w_5 at_3 x_3 = t_1 auat_2 aw_7$,

$\implies$ by (3.7), $M \vDash ay_2 at_3 x_3 = aw_7 = w_4$,

$\implies M \vDash w_1 = t'avat_2 aw_4 = t'avat_2(ay_2 at_3 x_3)$ whereas $M \vDash w_2 = x't_4 aw_6$,

$\implies M \vDash z = w_1 aw_2 = (t'avat_2 ay_2 at_3 x_3)a(x't_4 aw_6)$, where $M \vDash x_3 ax' = awa$.

We claim that $M \vDash \text{Max}^+ T_b(t_3, t'avat_2 ay_2)$.

Assume $M \vDash \text{Tally}_b(t'')\ \&\ t'' \subseteq_p t'avat_2 ay_2$.

$\implies$ by (4.17$^b$), $M \vDash t'' \subseteq_p t'av\ \vee\ t'' \subseteq_p t_2 ay_2$.

From the proof of (3bii2a) we have $M \vDash t'' < t_2$.

If $M \vDash t'' \subseteq_p t_2 a y_2 \subseteq_p w_5$, then $M \vDash t'' < t_3$ from $M \vDash Max^+T_b(t_3, w_5)$, hence also $M \vDash t_2 < t_3$.

It follows that $M \vDash Max^+T_b(t_3, t'avat_2ay_2)$.

We then have

$\quad M \vDash Pref(awa, t_3) \, \& \, Tally_b(t_4) \, \& \, t_3 < t_4 \, \&$

$\quad\quad\quad\quad \& \, z = (t'avat_2ay_2)at_3awat_4aw_6 \, \& \, Max^+T_b(t_3, t'avat_2ay_2)$,

hence $M \vDash Intf(z, t'avat_2ay_2, t_3, awa, t_4)$.

We also have

$\quad M \vDash (t'avat_2ay_2)at_3(x_3a) = w_1a \, \& \, (x_3a)Bw \, \& \, w_1aw_2 = z$.

Therefore $M \vDash Occ(w_1, a, w_2, z, t_3, awa, t_4)$, as required.

$\quad$ (3bii2c) $M \vDash w_5 a t_3 awa = w'a$.

$\Longrightarrow M \vDash w_5 a t_3 awa = w'a = t_1 auat_2 w_4 a$,

$\Longrightarrow M \vDash w' = t_1 auat_2 w_4$,

$\Longrightarrow M \vDash w'aw_2 = x = t_1 auat_2 w_4 aw_2 = t_1 auat_2 ax_1 = w_5 at_3 awat_4 aw_6$,

$\Longrightarrow$ by (3.7), $M \vDash w_2 = t_4 aw_6 \, \& \, w_4 aw_2 = ax_1$,

$\Longrightarrow M \vDash z = w_1 aw_2 = (t'avat_2 w_4)a(t_4 aw_6)$.

We have that $M \vDash (w_5 at_3)B(w'a) \, \& \, (t_1 auat_2)B(w'a)$.

$\Longrightarrow$ by (3.8), $M \vDash (w_5 at_3)B(t_1 auat_2) \lor w_5 at_3 = t_1 auat_2 \lor (t_1 auat_2)B(w_5 at_3)$.

If $M \vDash (w_5 at_3)B(t_1 auat_2)$, the argument is the same as in (3bii2bii2a).

If $M \vDash w_5 at_3 = t_1 auat_2$, we have, by $(4.24^b)$, that $M \vDash t_2 = t_3$, whence, by $(4.15^b)$, $M \vDash w_5 = t_1 au$.

$\Rightarrow$ by (3.7), $M \vDash t_3awa=t_2w_4a$,

$\Rightarrow$ by (3.7), $M \vDash awa=w_4a$.

We then derive, as in (3bii2bii1a), that

$$M \vDash z=t'avat_2ax_1=t'avat_3awat_4aw_6,$$

and further, that $M \vDash \text{Max}^+T_b(t_3,t'av)$ and $M \vDash \text{Intf}(z,t'av,t_3,awa,t_4)$.

Now, $M \vDash w_1=t'avat_2w_4=t'avat_3w_4$.

So we have $M \vDash (t'ava)t_3w_4a=(t'av)at_3awa=w_1a$ & $w_1aw_2=z$.

Therefore $M \vDash \text{Occ}(w_1,a,w_2,z,t_3,awa,t_4)$.

So we may assume that $M \vDash (t_1auat_2)B(w_5at_3)$.

$\Rightarrow$ $M \vDash \exists x_7\ t_1auat_2x_7=w_5at_3$.

We proceed as in (3bii2bii2c) and derive

$$M \vDash x_7=a \vee aBx_7,$$

and further, using (10.12), that $M \vDash t_1auat_2ay_2=w_5$.

$\Rightarrow$ $M \vDash t_1auat_2ay_2(at_3aw)=w_5at_3aw=w'=t_1auat_2w_4$,

$\Rightarrow$ by (3.7), $M \vDash w_4=ay_2at_3aw$,

$\Rightarrow$ $M \vDash w_1=t'avat_2w_4=t'avat_2(ay_2at_3aw)$ whereas $M \vDash w_2=t_4aw_6$.

We then prove that $M \vDash \text{Max}^+T_b(t_3,t'avat_2ay_2)$ exactly as in (3bii2bii2c),

whence we derive $M \vDash \text{Intf}(z,t'avat_2ay_2,t_3,awa,t_4)$.

We then have $M \vDash (t'avat_2ay_2)at_3awa=w_1a$ & $w_1aw_2=z$.

Therefore $M \vDash \text{Occ}(w_1,a,w_2,z,t_3,awa,t_4)$, as required.

 (3bii3) $M \vDash \text{Lastf}(x,t_3,awa,t_4)$.

$\Rightarrow M \vDash \text{Pref}(awa,t_3)$ & $\text{Tally}_b(t_4)$ & $t_3=t_4$ &

$$\& \ (x=t_3awat_4 \ \lor \ \exists w_5 \ (x=w_5at_3awat_4 \ \& \ Max^+T_b(t_3,w_5))).$$

Now, $M \vDash x=t_3awat_4$ is ruled out as in (3bi1).

So we may assume that $M \vDash \exists w_5 \ (x=w_5at_3awat_4 \ \& \ Max^+T_b(t_3,w_5))$.

Here an argument exactly analogous to that in (3bii2) applies with $aw_6$ omitted throughout.

This completes the argument that, under the principal hypothesis,

$\quad M \vDash \forall w_1,w_2(w_1aw_2=z \ \rightarrow \ Occ(w_1,a,w_2,z,t_3,awa,t_4))$,

that is, $M \vDash MinSet(z)$.

This completes the proof of (10.13).

(10.14) For any string concept $I \subseteq I_0$ there is a string concept $J \subseteq I$ such that

$QT^+ \vdash \forall z \in J \forall t',t'',y,v[\text{MinSet}(z)\ \&\ z=t'\text{ayat}''\text{avat}''\ \&\ \text{Firstf}(z,t',\text{aya},t'')\ \&$

$\&\ \text{Lastf}(z,t'',\text{ava},t'')\ \to\ \exists!z' \in J(\text{MinSet}(z')\ \&\ z \sim z'\ \&$

$\&\ \forall w(w\ \varepsilon\ z\ \to\ \forall w'(w'<_z w \leftrightarrow w'<_{z'} w))\ \&$

$\&\ \forall w, t_1, t_2\ (\text{Free}^+(z,t_1,\text{awa},t_2)\ \to\ \exists t_3 \text{Fr}(z',t_1 b,\text{awa},t_3))\ \&$

$\&\ \forall w, t_1, t_2\ (\text{Bound}(z,t_1,\text{awa},t_2)\ \vee\ \text{Free}^-(z,t_1,\text{awa},t_2)\ \to\ \exists t_3 \text{Fr}(z',t_1,\text{awa},t_3)))]$.

Let $J \equiv I_{4.8}\ \&\ I_{5.58}\ \&\ I_{9.14}\ \&\ I_{10.4}\ \&\ I_{10.7}$.

Assume $M \vDash \text{MinSet}(z)\ \&\ \text{Firstf}(z,t',\text{aya},t'')\ \&\ \text{Lastf}(z,t'',\text{ava},t'')$

where $M \vDash z=t'\text{ayat}''\text{avat}''$ and $M \vDash J(z)$.

$\Rightarrow$ since we may assume that J is downward closed under $\subseteq_p$, from $M \vDash J(z)$,

$M \vDash J(t')\ \&\ J(y)\ \&\ J(t'')\ \&\ J(v)$,

$\Rightarrow$ from $M \vDash \text{MinSet}(z)$, $M \vDash \text{Set}(z)$,

$\Rightarrow$ from $M \vDash z \neq aa$, $M \vDash \exists t^* \text{Env}(t^*,z)$,

$\Rightarrow$ from $M \vDash \text{Env}(t^*,z)$, $M \vDash t^* \subseteq_p z$,

$\Rightarrow$ from $M \vDash J(z)$, $M \vDash J(t^*)$,

$\Rightarrow$ from $M \vDash \text{Lastf}(z,t'',\text{ava},t'')$, $M \vDash \text{Env}(t'',z)$,

$\Rightarrow$ from $M \vDash \text{Firstf}(z,t',\text{aya},t'')\ \&\ (t'\text{ayat}''a)Bz$, by (5.1), $M \vDash t'<t''$,

$\Rightarrow$ from $M \vDash t' \in I \subseteq I_0$, $M \vDash t' \neq t''$,

$\Rightarrow$ from $M \vDash \text{Env}(t'',z)\ \&\ \text{Fr}(z,t',\text{aya},t'')\ \&\ \text{Fr}(z,t'',\text{ava},t'')$, $M \vDash y \neq v$,

$\Rightarrow$ from (5.58), $M \vDash \forall w(w\ \varepsilon\ z \leftrightarrow w=y \vee w=v)$,

$\Rightarrow$ from $M \vDash t'<t''$ & $\text{Tally}_b(t'')$, by (4.7), $M \vDash t'b \leq t''$.

(i) $M \vDash t'b<t''$.

$\Rightarrow$ from $M \vDash \text{Tally}_b(t')$, by (4.8), $M \vDash \text{Tally}_b(t'b)$,

$\Rightarrow$ from $M \vDash \text{Tally}_b(t'')$, $M \vDash \exists t^*(\text{Tally}_b(t^*) \& t'bt^*=t'')$.

Let $z'=t'bayat'bt^*avat''$.

$\Rightarrow$ since we may assume that the string concept J is closed with respect to *,

from $M \vDash J(t') \& J(y) \& J(t^*) \& J(v) \& J(t'')$,  $M \vDash J(z')$.

We first show that  (i1)  $M \vDash \text{MinSet}(z')$.

From  $M \vDash \text{Firstf}(z,t',aya,t'') \& \text{Lastf}(z,t'',ava,t'')$,  $M \vDash \text{Pref}(aya,t') \& \text{Pref}(ava,t'')$.

$\Rightarrow$ from $M \vDash \text{Tally}_b(t'b)$, $M \vDash \text{Pref}(aya,t'b)$,

$\Rightarrow$ $M \vDash \text{Pref}(aya,t'b) \& \text{Pref}(ava,t'bt^*) \& t'b<t'bt^* \& t'bt^*=t'' \&$

$\qquad\qquad\qquad\qquad\qquad$ & $y \neq v \& z'=t'bayat'bt^*avat''$,

$\Rightarrow$ by (10.7), $M \vDash \text{MinSet}(z')$,  as claimed.

Next, we have that  (i2)  $M \vDash z \sim z'$.

This is because, by (5.58),

   (i2a)  $M \vDash \forall w(w\, \varepsilon\, z \leftrightarrow w=y \lor w=v)$,

while, on the other hand, again by (5.58),

   (i2b)  $M \vDash \forall w(w\, \varepsilon\, z' \leftrightarrow w=y \lor w=v)$.

Thus (i2) follows as claimed.

Next, we show that  (i3)  $M \vDash \forall w(w\, \varepsilon\, z \rightarrow \forall u(u<_z w \leftrightarrow u<_{z'} w))$.

We have that  $M \vDash z'=t'bayat'bt^*avat''=t'bayat''avat''$.

⇒ as in the proof of (i1),  M ⊨ Pref(aya,t'b),

⇒ from M ⊨ (t'bayat''a)Bz' & t'b<t'',  M ⊨ Firstf(z',t'b,aya,t''),

⇒ from M ⊨ Pref(aya,t'),  M ⊨ Max$^+$T$_b$(t',aya),

⇒ M ⊨ Max$^+$T$_b$(t'',aya),

⇒ from  M ⊨ t'b<t'' by (4.17$^b$),  M ⊨ Max$^+$T$_b$(t'',t'bay),

⇒ from M ⊨ Pref(ava,t'') & Tally$_b$(t'') &

& z'=(t'bay)at''avat'' & Max$^+$T$_b$(t'',t'baya),

  M ⊨ Lastf(z',t'',ava,t'').

⇒ from M ⊨ Firstf(z,t',aya,t''), by (9.1), M ⊨ ∀u ¬(u<$_z$y),

⇒ from M ⊨ Firstf(z',t'b,aya,t''), by (9.1), M ⊨ ∀u ¬(u<$_{z'}$y).

Hence we have  (i3a)   M ⊨ ∀u(u<$_z$y ↔ u<$_{z'}$y).

⇒ from M ⊨ Lastf(z,t'',ava,t''),  M ⊨ ∀u(u ε z → u≤$_z$v),

⇒ from M ⊨ Lastf(z',t'',ava,t'') & z~z',  M ⊨ ∀u(u ε z → u≤$_{z'}$v),

⇒ M ⊨ ∀u(u ε z → (u<$_z$v ↔ u<$_{z'}$v)),

⇒ from  M ⊨ ∀u(¬(u ε z) → (¬(u<$_z$v) & ¬(u<$_{z'}$v)),

M ⊨ ∀u(¬(u ε z) → (u<$_z$v ↔ u<$_{z'}$v)).

Therefore,  we have  (i3b)   M ⊨ ∀u(u<$_z$v ↔ u<$_{z'}$v).

⇒ from (i2a), (i3a) and (i3b),  M ⊨ ∀w(w ε z → ∀u(u<$_z$w ↔ u<$_{z'}$w)), as claimed.

Next, we show that

  (i4)   M ⊨ ∀w,t$_1$,t$_2$ (Free$^+$(z,t$_1$,awa,t$_2$) → ∃t$_3$Fr(z',t$_1$b,awa,t$_3$)).

Assume   M ⊨ Free$^+$(z,t$_1$,awa,t$_2$).

We show that $M \vDash \forall w, t_1, t_2 \, \neg Free(z, t_1, awa, t_2)$.

Suppose, for a reductio, that $M \vDash Free(z, t_1, awa, t_2)$.

$\Rightarrow$ by (5.59), $M \vDash w=v \, \& \, t_1=t''$,

$\Rightarrow$ from $M \vDash Firstf(z, t', aya, t'')$, $M \vDash Fr(z, t', aya, t'')$,

$\Rightarrow$ from hypothesis $M \vDash Free(z, t_1, awa, t_2)$, $M \vDash \neg Firstf(z, t_1, awa, t_2)$,

$\Rightarrow$ by (5.15), $M \vDash y \neq w$,

$\Rightarrow$ since $M \vDash t' < t'' \, \& \, t' \neq t''$, $M \vDash y <_z v$.

Suppose, again for a reductio, that $M \vDash \exists u(y <_z u \, \& \, u <_z v)$.

$\Rightarrow$ by (9.4), $M \vDash u \neq y \, \& \, u \neq v \, \& \, u \, \varepsilon \, z$.

But this contradicts (5.58). Hence $M \vDash \neg \exists u(y <_z u \, \& \, u <_z v)$.

But then $M \vDash Fr(z, t', aya, t'') \, \& \, y <_z w \, \& \, \neg \exists u(y <_z u \, \& \, u <_z w)$ whereas $M \vDash t'b \neq t_1$

by hypothesis (i) since $M \vDash t'b \in I \subseteq I_0$.

This contradicts $M \vDash Free(z, t_1, awa, t_2)$.

Therefore $M \vDash \forall w, t_1, t_2 \, \neg Free(z, t_1, awa, t_2)$, as claimed.

Hence from hypothesis $M \vDash Free^+(z, t_1, awa, t_2)$ we have $M \vDash Firstf(z, t_1, awa, t_2)$.

$\Rightarrow$ from $M \vDash Firstf(z, t', aya, t'')$, by (5.15), $M \vDash t_1 = t' \, \& \, w = y$,

$\Rightarrow$ from $M \vDash Firstf(z', t'b, aya, t'')$, $M \vDash \exists t_3 Fr(z', t_1 b, awa, t_3)$, as required.

Finally, we show that

  (i5) $M \vDash \forall w, t_1, t_2 \, (Bound(z, t_1, awa, t_2) \lor Free^-(z, t_1, awa, t_2) \rightarrow$

$\rightarrow \exists t_3 Fr(z', t_1, awa, t_3))$.

$\Rightarrow$ from $M \vDash \neg Free(z, t_1, awa, t_2)$, $M \vDash Bound(z, t_1, awa, t_2)$,

$\Rightarrow M \vDash Fr(z, t_1, awa, t_2) \, \& \, \neg Firstf(z, t_1, awa, t_2)$,

$\Rightarrow$ from $M \vDash \text{Firstf}(z,t',aya,t'')$, by (5.15), $M \vDash w \neq y$,

$\Rightarrow$ $M \vDash w=v$,

$\Rightarrow$ from $M \vDash \text{Env}(t^*,z)$ & $\text{Fr}(z,t_1,ava,t_2)$ & $\text{Fr}(z,t'',ava,t'')$, $M \vDash t_1=t''$,

$\Rightarrow$ from $M \vDash \text{Last}(z,t'',ava,t'')$, $M \vDash \exists t_3 \text{Fr}(z',t_1,awa,t_3)$, as required.

To establish the uniqueness of z', assume that

$M \vDash \text{MinSet}(z'')$ & $z \sim z''$ & $\forall w(w \,\varepsilon\, z \to \forall u(u <_z w \leftrightarrow u <_{z''} w))$ &

        & $\forall w,t_1,t_2 \,(\text{Free}^+(z,t_1,awa,t_2) \to \exists t_3 \text{Fr}(z'',t_1 b,awa,t_3))$ &

& $\forall w,t_1,t_2 \,(\text{Bound}(z,t_1,awa,t_2) \vee \text{Free}^-(z,t_1,awa,t_2) \to \exists t_3 \text{Fr}(z'',t_1,awa,t_3))$.

where $M \vDash J(z'')$.

$\Rightarrow$ from $M \vDash \forall w(w \,\varepsilon\, z \leftrightarrow w=y \vee w=v)$ & $z \sim z''$, $M \vDash \forall w(w \,\varepsilon\, z'' \leftrightarrow w=y \vee w=v)$,

$\Rightarrow$ from hypothesis about z'' and $M \vDash \text{Firstf}(z,t',aya,t'')$, $M \vDash \exists t_3 \text{Fr}(z'',t'b,aya,t_3)$,

$\Rightarrow$ from $M \vDash \text{Firstf}(z,t',aya,t'')$, by (9.1) and (9.7), $M \vDash \forall w(w \,\varepsilon\, z \to y \leq_z w)$,

$\Rightarrow$ $M \vDash \forall w(w \,\varepsilon\, z'' \to y \leq_{z''} w)$,

$\Rightarrow$ from $M \vDash \text{MinSet}(z'')$ & $y \,\varepsilon\, z''$, $M \vDash \text{Set}(z'')$ & $z'' \neq aa$,

$\Rightarrow$ $M \vDash \exists t^{**} \text{Env}(t^{**},z'')$ & $\text{Fr}(z'',t'b,aya,t_3)$ & $\forall w(w \,\varepsilon\, z'' \to y \leq_{z''} w)$,

$\Rightarrow$ by (9.10), $M \vDash \text{Firstf}(z'',t'b,aya,t_3)$,

$\Rightarrow$ from $M \vDash \text{Last}(z,t'',ava,t'')$, by (9.3), $M \vDash \forall w(w \,\varepsilon\, z \to w \leq_z v)$,

$\Rightarrow$ $M \vDash \forall w(w \,\varepsilon\, z'' \to w \leq_{z''} v)$,

$\Rightarrow$ from $M \vDash v \,\varepsilon\, z$, $M \vDash \exists t_4,t_5 \,\text{Fr}(z'',t_4,ava,t_5)$,

$\Rightarrow$ by (9.8), $M \vDash t_4 = t^{**} = t_5$ & $\text{Lastf}(z'',t_4,ava,t_5)$,

$\Rightarrow$ $M \vDash \text{Lastf}(z'',t^{**},ava,t^{**})$,

$\Rightarrow$ from $M \vDash y \neq v$ & $\text{Firstf}(z,t',aya,t'')$, by (5.15), $M \vDash \neg \text{Firstf}(z,t'',ava,t'')$,

$\Rightarrow$ by (9.22), $M \vDash \text{Free}(z,t'',ava,t'') \lor \text{Bound}(z,t'',ava,t'')$,

$\Rightarrow$ from $M \vDash \forall w,t_1,t_2 \neg \text{Free}(z,t_1,awa,t_2)$, $M \vDash \text{Bound}(z,t'',ava,t'')$,

$\Rightarrow$ from hypothesis about $z''$, $M \vDash \exists t_6 \text{Fr}(z'',t'',ava,t_6)$,

$\Rightarrow$ from $M \vDash \text{MinSet}(z'')$ & $\forall w(w \: \varepsilon \: z'' \leftrightarrow w=y \lor w=v)$ & $y \neq v$, by (10.4),

$M \vDash \exists t_7,t_8 (\text{Tally}_b(t_7)$ & $\text{Tally}_b(t_8))$ & $((y<_{z''}v$ & $z''=t_7 ayat_8 avat_8) \lor$

$\lor ((v<_{z''}y$ & $z''=t_7 avat_8 ayat_8)))$,

$\Rightarrow$ from $M \vDash y<_z v$, by (9.6), $M \vDash \neg(v<_z y)$,

$\Rightarrow$ from hypothesis about $z''$, $M \vDash \neg(v<_{z''}y)$,

$\Rightarrow M \vDash z''=t_7 ayat_8 avat_8$,

$\Rightarrow$ from $M \vDash \text{Firstf}(z'',t'b,aya,t_3)$, $M \vDash \text{Tally}_b(t'b)$ & $(t'ba)Bz''$,

$\Rightarrow$ by (4.23$^b$), $M \vDash t_7 = t'b$,

$\Rightarrow$ from $M \vDash \text{Env}(t^{**},z'')$ & $\text{Fr}(z'',t'',ava,t'')$ & $\text{Lastf}(z'',t^{**},ava,t^{**})$,

$M \vDash t^{**}=t''$ & $\text{Tally}_b(t^{**})$ & $(t''a)Ez''$,

$\Rightarrow$ by (4.24$^b$), $M \vDash t_8 = t''$,

$\Rightarrow M \vDash z''=t_7 ayat_8 avat_8 = t'bayat''avat''=z$, as required.

(ii) $M \vDash t'b=t''$.

Let $z'=t'bayat''bavat''b$.

$\Rightarrow$ from $M \vDash J(t')$ & $J(y)$ & $J(t'')$ & $J(v)$, $M \vDash J(z')$.

We first show that (ii1) $M \vDash \text{MinSet}(z')$.

We proceed analogously to (i1), observing that, since by (4.8),

$M \vDash \text{Tally}_b(t''b)$, we have $M \vDash \text{Pref}(ava,t''b)$. Then, since $M \vDash t'b<t''b$,

by (10.7), $M \vDash \text{MinSet}(z')$.

We show  (ii2)  $M \vDash z \sim z'$  exactly analogously to (i2).

Also, the proof of  (ii3)  $M \vDash \forall w(w \: \varepsilon \: z \rightarrow \forall u(u<_z w \leftrightarrow u<_{z'} w))$

is exactly analogous to (i3), with t'' replaced throughout by t''b.

To show that  (ii4)  $M \vDash \forall w, t_1, t_2 \: (Free^+(z,t_1,awa,t_2) \rightarrow \exists t_3 Fr(z',t_1 b,awa,t_3))$

assume  $M \vDash Free^+(z,t_1,awa,t_2)$.

$\Rightarrow M \vDash Fr(z,t_1,awa,t_2) \: \& \: (Firstf(z,t_1,awa,t_2) \: v \: Free(z,t_1,awa,t_2)$.

Suppose that  $M \vDash Firstf(z,t_1,awa,t_2)$.

We then obtain $M \vDash \exists t_3 Fr(z',t_1 b,awa,t_3)$  exactly as in the proof of (i4).

Suppose, on the other hand, that  $M \vDash Free(z,t_1,awa,t_2)$.

$\Rightarrow M \vDash \neg Firstf(z,t_1,awa,t_2)$,

$\Rightarrow$ from $M \vDash Firstf(z,t',aya,t'')$, by (5.15),  $M \vDash w \neq y$,

$\Rightarrow$ from $M \vDash \forall w(w \: \varepsilon \: z \leftrightarrow w=y \: v \: w=v)$,  $M \vDash w=v$.

But we have, analogously to the reasoning in the proof of (i3), that

$\qquad M \vDash Lastf(z,t''b,ava,t''b)$.

$\Rightarrow$ from $M \vDash Env(t^*,z) \: \& \: Fr(z,t_1,awa,t_2) \: \& \: Fr(z,t'',awa,t'')$,  $M \vDash t_1 = t''$,

$\Rightarrow$ from $M \vDash Lastf(z',t''b,ava,t''b)$,  $M \vDash \exists t_3 Fr(z',t_1 b,awa,t_3)$,  as required.

This suffices to establish (ii4).

Finally, we show that

  (ii5)  $M \vDash \forall w, t_1, t_2 \: (Bound(z,t_1,awa,t_2) \: v \: Free^-(z,t_1,awa,t_2) \rightarrow$

$\qquad\qquad\qquad\qquad\qquad\qquad\qquad \rightarrow \exists t_3 Fr(z',t_1,awa,t_3))$.

Suppose that  $M \vDash Bound(z,t_1,awa,t_2)$,

$\Rightarrow M \vDash Fr(z,t_1,awa,t_2) \: \& \: \neg Firstf(z,t_1,awa,t_2)$,

$\Rightarrow$ as in the proof of (ii4), $M \vDash w=v$,

$\Rightarrow$ as in the proof of (i4), $M \vDash \neg \exists u(y<_z u \ \& \ u<_z v)$,

$\Rightarrow$ from $M \vDash \text{Bound}(z,t_1,awa,t_2) \ \& \ \text{Fr}(z,t',aya,t'') \ \& \ y<_z v)$, $M \vDash t'b<t_1$,

$\Rightarrow$ from $M \vDash \text{Env}(t^*,z) \ \& \ \text{Fr}(z,t'',ava,t'') \ \& \ \text{Fr}(z,t_1,ava,t_2)$, $M \vDash t_1=t''$,

$\Rightarrow$ from hypothesis (ii), $M \vDash t''=t'b<t''$, contradicting $M \vDash t'' \in I \subseteq I_0$.

Therefore, $M \vDash \neg \text{Bound}(z,t_1,awa,t_2)$.

Hence we have $M \vDash \text{Free}^-(z,t_1,awa,t_2)$.

$\Rightarrow M \vDash \text{Fr}(z,t_1,awa,t_2) \ \& \ \neg \text{Firstf}(z,t_1,awa,t_2)$,

$\Rightarrow$ as in the proof of (ii4), $M \vDash w=v$.

We are going to show that $M \vDash \text{Free}^+(z,t_1,awa,t_2)$.

Assume that $M \vDash \text{Fr}(z,t_4,aua,t_5) \ \& \ u \leq_z w$.

$\Rightarrow M \vDash u \ \varepsilon \ z$,

$\Rightarrow M \vDash u=y \ \lor \ u=v=w$.

Now, we have $M \vDash \text{Firstf}(z,t',aya,t'')$.

Hence, if $M \vDash u=y$, from $M \vDash \text{Env}(t,z) \ \& \ \text{Fr}(z,t_4,aua,t_5) \ \& \ \text{Fr}(z,t',aua,t'')$,

$M \vDash t'=t_4$. So, if $M \vDash u=y$, then $M \vDash \text{Firstf}(z,t_4,aua,t'')$, whence, by (5.15),

$M \vDash t''=t_5$. Hence, if $M \vDash u=y$, then $M \vDash \text{Firstf}(z,t_4,aua,t_5)$.

On the other hand, if $M \vDash u=v=w$, from

$$M \vDash \text{Env}(t,z) \ \& \ \text{Fr}(z,t_4,aua,t_5) \ \& \ \text{Fr}(z,t',aua,t''),$$

we have $M \vDash t_4=t''$.

Also, if $M \vDash u=v=w$, then $M \vDash \neg \text{Firstf}(z,t_4,aua,t_5)$.

For, if $M \vDash \text{Firstf}(z,t_4,aua,t_5)$, then $M \vDash \text{Firstf}(z,t_4,awa,t_5)$, whence

$M \vDash v=w=y$ from $M \vDash \text{Firstf}(z,t',aya,t'')$ by (5.15), a contradiction.

Assume that $M \vDash \text{Fr}(z,t_6,au'a,t_7)$ & $u'<_z u$ & $\neg\exists w'(u'<_z w'$ & $w'<_z u)$.

$\Rightarrow$ by (9.4), $M \vDash u' \neq u$,

$\Rightarrow M \vDash u' \neq v$,

$\Rightarrow$ from $M \vDash u' \varepsilon z$, $M \vDash u'=y$,

$\Rightarrow$ from $M \vDash \text{Env}(t,z)$ & $\text{Fr}(z,t',aya,t'')$ & $\text{Fr}(z,t_6,aya,t_7)$, $M \vDash t_6=t'$,

$\Rightarrow$ from hypothesis (ii), $M \vDash t_6 b = t'b = t''$,

$\Rightarrow M \vDash t_4 = t'' = t_6 b$.

So we have that, if $M \vDash u=v=w$, then

$\quad M \vDash \forall u',t_6,t_7(\text{Fr}(z,t_6,au'a,t_7)$ & $u'<_z u$ & $\neg\exists w'(u'<_z w'$ & $w'<_z u) \rightarrow t_4=t_6 b)$

along with $M \vDash \text{Fr}(z,t_4,aua,t_5)$ & $\neg\text{Firstf}(z,t_4,aua,t_5)$.

In other words, if $M \vDash u=v=w$, then $M \vDash \text{Free}(z,t_4,aua,t_5)$.

So we have shown that

$\quad M \vDash \forall u,t_4,t_5(\text{Fr}(z,t_4,aua,t_5)$ & $u \leq_z w \rightarrow \text{Firstf}(z,t_4,aua,t_5) \vee \text{Free}(z,t_4,aua,t_5))$.

Along with $M \vDash \text{Fr}(z,t_1,awa,t_2)$, this suffices to establish $M \vDash \text{Free}^+(z,t_1,awa,t_2)$.

But then $M \vDash \neg\text{Free}^-(z,t_1,awa,t_2)$, contrary to hypothesis.

Therefore, $M \vDash \forall w,t_1,t_2(\neg\text{Bound}(z,t_1,awa,t_2)$ & $\neg\text{Free}^-(z,t_1,awa,t_2))$,

and (ii5) holds trivially.

It remains to establish the uniqueness of z' under hypothesis (ii).

Assume that

$M \vDash \text{MinSet}(z'')$ & $z \sim z''$ & $\forall w(w \varepsilon z \rightarrow \forall u(u<_z w \leftrightarrow u<_{z''} w))$ &

$\quad\quad$ & $\forall w,t_1,t_2 (\text{Free}^+(z,t_1,awa,t_2) \rightarrow \exists t_3 \text{Fr}(z'',t_1 b,awa,t_3))$ &

& $\forall w,t_1,t_2$ (Bound($z,t_1,awa,t_2$) v Free$^-$($z,t_1,awa,t_2$) → $\exists t_3$Fr($z'',t_1,awa,t_3$))

where  M ⊨ J($z''$).

We proceed as in (i) to obtain    M ⊨ $z''=t_7ayat_8avat_8$

and further that    M ⊨ $t_7=t'b$.

We have, as in the proof of (ii5), that if  M ⊨ $t'b=t''$, then

  M ⊨ $\forall w,t_1,t_2$(Fr($z,t_1,awa,t_2$) & ¬Firstf($z,t_1,awa,t_2$) → Free$^+$($z,t_1,awa,t_2$)).

Now, we do have that  M ⊨ Fr($z,t'',ava,t''$).

Also, from  M ⊨ Firstf($z,t',aya,t''$) & $y \neq v$,  by (5.15),   M ⊨ ¬Firstf($z,t'',ava,t''$).

Hence  M ⊨ Free$^+$($z,t'',ava,t''$).

⟹ by the choice of $z''$,  M ⊨ $\exists t_3$Fr($z'',t''b,ava,t_3$),

⟹ from M ⊨ Env($t^{**},z''$) & Lastf($z'',t^{**},ava,t^{**}$),  M ⊨ $t^{**}=t''b$,

⟹ M ⊨ ($at^{**}$)E$z''$ & ($at_8$)E$z''$,

⟹  by (4.24$^b$),  M ⊨ $t''b=t^{**}=t_8$,

⟹ M ⊨ $z''= t'bayat''bavat''b=z'$,  as required.

This completes the proof of (10.14).

(10.15) For any string concept $I \subseteq I_0$ there is a string concept $J \subseteq I$ such that

$QT^+ \vdash \forall y \in J \exists t_0, z \in J(MinMax^+T_b(t_0,y)\ \&\ z=t_0ayat_0\ \&\ Env(t_0,z)\ \&\ Set^*(z)\ \&$

$\&\ \forall w(w\ \varepsilon\ z \leftrightarrow w=y))$.

Let $I' \equiv I_{6.8}\ \&\ I_{9.1}\ \&\ I_{9.4}\ \&\ I_{10.5}$, and let J be obtained from I' as in (6.8).

Let $M \vDash J(y)$.

$\Rightarrow$ by (6.8), $M \vDash \exists t_0, z \in J\ MinMax^+T_b(t_0,y)$.

Since we may assume that J is closed under *, we have that $M \vDash J(z)$ where $z = t_0ayat_0$.

We then have $M \vDash Tally_b(t_0)\ \&\ Max^+T_b(t_0,aya)$.

$\Rightarrow M \vDash Pref(aya,t_0)$,

$\Rightarrow$ by (10.5), $M \vDash MinSet(z)$,

$\Rightarrow$ since $M \vDash Firstf(z,t_0,aya,t_0)\ \&\ Lastf(z,t_0,aya,t_0)$, by (5.22),

$M \vDash Env(t_0,z)\ \&\ y\ \varepsilon\ z\ \&\ \forall w(w\ \varepsilon\ z \leftrightarrow w=y)$.

To show $M \vDash Lex^+(z)$, we need $M \vDash \forall v,w(v<_z w \rightarrow v \prec w)$.

But we have $M \vDash Firstf(z,t_0,aya,t_0)$.

Assume $M \vDash v <_z w$.

$\Rightarrow M \vDash v\ \varepsilon\ z\ \&\ w\ \varepsilon\ z$,

$\Rightarrow M \vDash v=y=w$, $M \vDash y <_z y$, contradicting (9.4).

Therefore, $M \vDash \forall v,w\ \neg(v <_z w)$.

Hence, trivially, $M \vDash \forall v,w(v <_z w \rightarrow v \prec w)$, and so $M \vDash Lex^+(z)$.

Finally, to show $M \vDash Special(z)$, assume $M \vDash Fr(z,t_1,ava,t_2)$.

$\Rightarrow$ M ⊨ v ε z,

$\Rightarrow$ M ⊨ v=y,

$\Rightarrow$ M ⊨ Fr(z,$t_0$,aya,$t_0$) & Fr(z,$t_1$,aya,$t_2$),

$\Rightarrow$ from M ⊨ Env($t_0$,z), M ⊨ $t_0$=$t_1$.

Hence, to establish M ⊨ Special(z), it suffices to show that

$$M ⊨ MMax^+T_b(t_0,y,z).$$

We have that M ⊨ y ε z & Set(z) & Tally$_b$($t_0$).

$\Rightarrow$ since M ⊨ Fr(z,$t_0$,aya,$t_0$), by (9.1), M ⊨ $\forall$w ¬(w<$_z$y),

$\Rightarrow$ M ⊨ $\forall$t(Tally$_b$(t) $\to$ Max$^+$(t,y,z)),

$\Rightarrow$ M ⊨ Max$^+$($t_0$,y,z).

Assume that M ⊨ Max$^+$(t',y,z) & Max$^+$T$_b$(t',y).

$\Rightarrow$ since M ⊨ MinMax$^+$T$_b$($t_0$,y), M ⊨ $t_0$≤t'.

So we have that

$\quad$ M ⊨ Max$^+$($t_0$,y,z) & $\forall$t'(Max$^+$(t',y,z) & Max$^+$T$_b$(t',y) $\to$ $t_0$≤t'),

that is, M ⊨ MMax$^+$T$_b$($t_0$,y,z).

This completes the proof that M ⊨ Special(z), and the proof of (10.15).

(10.16) For any string concept $I \subseteq I_0$ there is a string concept $J \subseteq I$ such that

$QT^+ \vdash \forall x \in J \forall t_1,t_2,v[\text{Special}(x) \ \& \ \text{Firstf}(x,t_1,ava,t_2) \rightarrow \text{MinMax}^+T_b(t_1,v)]$.

Let $J \equiv I_{9.1}$.

Assume $M \vDash \text{Special}(x) \ \& \ \text{Firstf}(x,t_1,ava,t_2)$ where $M \vDash J(x)$.

$\Rightarrow$ from $M \vDash \text{Firstf}(x,t_1,ava,t_2)$, $M \vDash \text{Pref}(ava,t_1)$,

$\Rightarrow M \vDash \text{Max}^+T_b(t_1,v)$.

Assume now that $M \vDash \text{Max}^+T_b(t,v)$.

$\Rightarrow$ from $M \vDash \text{Firstf}(x,t_1,ava,t_2)$, by (9.1), $M \vDash \forall u \neg(u <_x v)$,

$\Rightarrow M \vDash \forall u,t',t''(\text{Fr}(x,t',aua,t'') \ \& \ u <_x v \rightarrow t' < t)$.

But we have $M \vDash v \ \varepsilon \ x \ \& \ \text{Set}(x) \ \& \ \text{Tally}_b(t)$.

$\Rightarrow M \vDash \text{Max}^+(t,v,x)$,

$\Rightarrow M \vDash \text{Max}^+(t,v,x) \ \& \ \text{Max}^+T_b(t,v)$,

$\Rightarrow$ from hypothesis $M \vDash \text{Special}(x)$, $M \vDash \text{MMax}^+T_b(t_1,v,x)$,

$\Rightarrow M \vDash t_1 \leq t$.

So we have shown that $M \vDash \text{Max}^+T_b(t,v) \rightarrow t_1 \leq t$.

Along with $M \vDash \text{Max}^+T_b(t_1,v)$, this yields $M \vDash \text{MinMax}^+T_b(t_1,v)$, as required.

This completes the proof of (10.16).

(10.17) For any string concept I⊆I₀ there is a string concept J⊆I such that

$QT^+ \vdash \forall x \in J \forall x',z,t,t',t'',t_3,t_4,v_0,w[Env(t,x)$ & $x=t'wz$ & $Env(t'',x')$ & $x'=t'wt''$ &

& $aBw$ & $aEw$ & $Env(t,z)$ & $Firstf(z,t_3,av_0a,t_4)$ & $t''<t_3$ & $\neg\exists u(u \; \varepsilon \; x'$ & $u \; \varepsilon \; z)$ &

& $\forall u,v(u \; \varepsilon \; x'$ & $v \; \varepsilon \; z \rightarrow u \prec v)$ & $Special(x')$ & $Special(z) \rightarrow Special(x)]$.

Let $J \equiv I_{5.39}$ & $I_{9.11}$ & $I_{9.18}$.

Assume

$M \vDash Env(t,x)$ & $x=t'wz$ & $Env(t'',x')$ & $x'=t'wt''$ & $aBw$ & $aEw$ & $Env(t,z)$ &

& $Firstf(z,t_3,av_0a,t_4)$ & $t''<t_3$,

where $M \vDash \neg\exists u(u \; \varepsilon \; x'$ & $u \; \varepsilon \; z)$ & $\forall u,v(u \; \varepsilon \; x'$ & $v \; \varepsilon \; z \rightarrow u \prec v)$ and $M \vDash J(x)$.

Assume now that $M \vDash Fr(x,t_1,ava,t_2)$.

$\Rightarrow$ by (5.39), $M \vDash \exists t_5 Fr(x',t_1,ava,t_5) \vee Fr(z,t_3,ava,t_4) \vee \exists t_5 Fr(z,t_1,ava,t_5)$.

We distinguish the three cases:

(1) $M \vDash \exists t_5 Fr(x',t_1,ava,t_5)$.

We first show that (1a) $M \vDash Max^+(t_1,v,x)$.

We already have that $M \vDash v \; \varepsilon \; x$ & $Set(x)$ & $Tally_b(t_1)$.

Assume $M \vDash Fr(x,t_6,aua,t_7)$ & $u<_x v$. Then from $M \vDash Fr(x,t_1,ava,t_2)$, by (9.14),

$M \vDash t_6<t_1$. Hence (1a) holds.

Assume now that $M \vDash Max^+(t_0,v,x)$ & $Max^+T_b(t_0,v)$.

We claim that $M \vDash Max^+(t_0,v,x')$.

For, assume $M \vDash Fr(x',t_6,aua,t_7)$ & $u<_{x'} v$.

$\Rightarrow$ by (5.6) and (9.11), $M \vDash \exists t_8 Fr(x,t_6,aua,t_8)$ & $u<_x v$,

$\implies$ from $M \vDash Max^+(t_0,v,x)$, $M \vDash t_6 < t_0$.

Hence $M \vDash Max^+(t_0,v,x')$ as claimed, since $M \vDash v \,\varepsilon\, x'$ & $Set(x')$ & $Tally_b(t_0)$.

Now, from hypotheses $M \vDash Special(x')$ and (1) we have

$$M \vDash MMax^+T_b(t_1,v,x').$$

$\implies$ from $M \vDash Max^+(t_0,v,x')$ & $Max^+T_b(t_0,v)$, $M \vDash t_1 \leq t_0$.

So we have shown that

(1b) $M \vDash Max^+(t_0,v,x)$ & $Max^+T_b(t_0,v) \rightarrow t_1 \leq t_0$.

Then (1a) and (1b) establish that $M \vDash MMax^+T_b(t_1,v,x)$.

(2) $M \vDash Fr(z,t_3,ava,t_4)$.

$\implies$ from hypothesis $M \vDash Firstf(z,t_3,av_0a,t_4)$,

$$M \vDash Fr(z,t_3,av_0a,t_4) \,\&\, Fr(z,t_3,ava,t_4),$$

$\implies$ from hypothesis $M \vDash Env(t,z)$, $M \vDash v=v_0$.

We first show that (2a) $M \vDash Max^+(t_3,v,x)$.

We already have that $M \vDash v \,\varepsilon\, x$ & $Set(x)$ & $Tally_b(t_3)$.

Assume that $M \vDash Fr(x,t_5,aua,t_6)$ & $u <_x v$.

$\implies$ by (9.14), $M \vDash t_5 < t_1$.

Let $M \vDash Tally_b(t_0)$ & $t_0 \subseteq_p t'w$.

$\implies M \vDash t_0 \subseteq_p t'w \subseteq_p x'$,

$\implies$ from $M \vDash Env(t'',x')$, $M \vDash Max^+T_b(t'',x')$, $M \vDash t_0 \leq t'' < t_3$.

Hence we have that $M \vDash Max^+T_b(t_3,t'w)$.

$\implies M \vDash Fr(z,t_3,ava,t_4)$ & $x=t'wz$ & $aBw$ & $aEw$ & $Max^+T_b(t_3,t'w)$,

$\implies$ by (5.25), $M \vDash Fr(x,t_3,ava,t_4)$,

$\Rightarrow$ from hypothesis $M \vDash Fr(x,t_1,ava,t_2)$ since $M \vDash Env(t,x)$, $M \vDash t_1=t_3$,

$\Rightarrow M \vDash t_5<t_3$.

Thus we established that $M \vDash (Fr(x,t_5,aua,t_6) \& u<_x v \to t_5<t_3)$,

which suffices to prove $M \vDash Max^+(t_3,v,x)$.

Next, assume that $M \vDash Max^+(t_0,v,x) \& Max^+T_b(t_0,v)$.

$\Rightarrow$ from hypothesis $M \vDash Firstf(x,t_3,ava,t_4)$, by (9.1), $M \vDash \forall w' \neg(w'<_z v)$,

$\Rightarrow$ since $M \vDash v \varepsilon z \& Set(z)$, $M \vDash \forall t_5(Tally_b(t_5) \to Max^+(t_5,v,z))$,

$\Rightarrow$ in particular, $M \vDash Max^+(t_0,v,z)$,

$\Rightarrow$ from hypothesis $M \vDash Special(z)$, $M \vDash MMax^+T_b(t_3,v,z)$,

$\Rightarrow M \vDash Max^+(t_0,v,z) \& Max^+T_b(t_0,v)$, $M \vDash t_3 \leq t_0$.

Thus we established that

(2b) $M \vDash (Max^+(t_0,v,x) \& Max^+T_b(t_0,v) \to t_3 \leq t_0)$.

Along with (2a), this suffices to show that $M \vDash MMax^+T_b(t_3,v,x)$.

(3) $M \vDash \exists t_5 Fr(z,t_1,ava,t_5)$.

Again, we first show that (3a) $M \vDash Max^+(t_1,v,x)$.

This is proved exactly as (1a).

Assume now that $M \vDash Max^+(t_0,v,x) \& Max^+T_b(t_0,v)$.

We claim that $M \vDash Max^+(t_0,v,z)$.

For this, assume that $M \vDash Fr(z,t_6,aua,t_7) \& u<_z v$.

$\Rightarrow$ from $M \vDash Firstf(z,t_3,av_0a,t_4)$, by (5.20), $M \vDash t_3 \leq t_6$.

But we have, as in the proof of (2a), that $M \vDash Max^+T_b(t_3,t'w)$.

$\Rightarrow M \vDash Max^+T_b(t_6,t'w)$,

⟹ M ⊨ Fr(z,$t_6$,aua,$t_7$) & x=t'wz & aBw & aEw & Max$^+$T$_b$($t_6$,t'w),

⟹ by (5.25), M ⊨ Fr(x,$t_6$,aua,$t_7$),

⟹ M ⊨ Env(t,x) & x=t'wz & aBw & aEw & Env(t'',x') & x'=t'wt'' & Env(t,z) &

& Firstf(z,$t_3$,a$v_0$a,$t_4$) & t''<$t_3$,

⟹ from M ⊨ u<$_z$v, by (9.18), M ⊨ u<$_x$v,

⟹ from M ⊨ Max$^+$($t_0$,v,x), M ⊨ $t_6$<$t_0$,

⟹ since M ⊨ v ε z & Set(z) & Tally$_b$($t_0$), M ⊨ Max$^+$($t_0$,v,z), as claimed.

Now, from hypothesis M ⊨ Special(z) and (3), we have that

$$M \vDash MMax^+T_b(t_1,v,z).$$

⟹ from M ⊨ Max$^+$($t_0$,v,z) & Max$^+$T$_b$($t_0$,v), M ⊨ $t_1$≤$t_0$.

So we have shown that

(3b) M ⊨ (Max$^+$($t_0$,v,x) & Max$^+$T$_b$($t_0$,v) → $t_1$≤$t_0$).

But (3a) and (3b) suffice to establish M ⊨ MMax$^+$T$_b$($t_1$,v,x).

In (1)-(3) we therefore proved

$$M \vDash (Fr(x,t_1,ava,t_2) \to MMax^+T_b(t_1,v,x)).$$

Along with M ⊨ Set(x), this yields M ⊨ Special(x), as required.

This complete the proof of (10.17).

(10.18) For any string concept I⊆I₀ there is a string concept J⊆I such that

$$QT^+ \vdash \forall x \in J \forall t,t',x',y,w,z[Env(tb,x) \& x=t'wz \& Env(t,x') \& x'=t'wt \&$$

$$\& aBw \& aEw \& z=tbayatb \& Max^+T_b(tb,y) \& \neg(y \varepsilon x') \&$$

$$\& \forall u(u \varepsilon x' \rightarrow u \prec y) \& Special(x') \rightarrow Special(x)].$$

Let $J \equiv I_{5.22} \& I_{5.46} \& I_{9.14}$.

Assume  M ⊨ Env(tb,x) & x=t'wz & Env(t,x') & x'=t'wt & aBw & aEw

where  M ⊨ z=tbayatb & Max⁺T_b(tb,y) & ¬(y ε x') & ∀u(u ε x' → u≺y)

and  M ⊨ J(x).

Suppose that  M ⊨ Special(x').

⟹ from  M ⊨ z=tbayatb & Max⁺T_b(tb,y),  M ⊨ Pref(aya,tb),

⟹ M ⊨ Firstf(z,tb,aya,tb) & Lastf(z,tb,aya,tb),

⟹ by (5.22),  M ⊨ Env(tb,z) & ∀u(u ε z ↔ u=y).

Assume now that  M ⊨ Fr(x,t₁,ava,t₂).

⟹ since  M ⊨ Firstf(z,tb,aya,tb), by (5.39),

  M ⊨ ∃t₃Fr(x',t₁,ava,t₃) v Fr(z,tb,ava,tb) v ∃t₃Fr(z,t₁,ava,t₃),

⟹ from  M ⊨ Env(tb,z),  M ⊨ MaxT_b(tb,z),

⟹ by (5.5),  M ⊨ (Fr(z,t₁,ava,t₃) → t₁=tb=t₃ & v=y),

⟹ M ⊨ ∃t₃Fr(x',t₁,ava,t₃) v (Fr(z,t₁,ava,t₁) & t₁=tb & v=y).

We distinguish two cases:

(1)  M ⊨ ∃t₃Fr(x',t₁,ava,t₃).

We argue exactly as in (10.17), part (1), to show that  M ⊨ MMax⁺T_b(t₁,v,x).

(2)  $M \vDash Fr(z,t_1,ava,t_1)$ & $t_1=tb$ & $v=y$.

We show that

  (2a)  $M \vDash Max^+(tb,v,x)$.

We have that  $M \vDash v \, \varepsilon \, x$ & $Set(x)$ & $Tally_b(tb)$.

Assume  $M \vDash Fr(x,t_4,aua,t_5)$ & $u<_x v$.

$\Rightarrow M \vDash u \, \varepsilon \, x$,

$\Rightarrow$ by the proof of (5.46),  $M \vDash u \, \varepsilon \, x'$ v $u \, \varepsilon \, z$.

If  $M \vDash u \, \varepsilon \, z$, then $M \vDash u=y=v$, whence $M \vDash u<_x u$, contradicting (9.4).

$\Rightarrow M \vDash \neg(u \, \varepsilon \, z)$,

$\Rightarrow M \vDash u \, \varepsilon \, x'$,

$\Rightarrow$ from (5.39),  $M \vDash \exists t_6 Fr(x',t_4,aua,t_6)$,

$\Rightarrow$ from  $M \vDash Env(t,x')$, $M \vDash t_4 \leq t$,

$\Rightarrow M \vDash t_4 \leq t < tb$.

This establishes (2a).

Assume now that  $M \vDash Max^+(t',v,x)$ & $Max^+T_b(t',v)$.

We want to show that  $M \vDash tb \leq t'$.

Suppose, for a reductio, that  $M \vDash t' < tb$.

$\Rightarrow M \vDash t' \leq t$,

$\Rightarrow$ from  $M \vDash Env(t,x')$, $M \vDash \exists v' Lastf(x',t,av'a,t)$,

$\Rightarrow M \vDash Fr(x',t,av'a,t)$,

$\Rightarrow$ since  $M \vDash z=t'wz=t'wtbayatb$, by (5.6), $M \vDash \exists t_6 Fr(x,t,av'a,t_6)$,

$\Rightarrow$ since  $M \vDash Fr(x,tb,ava,tb)$ & $t<tb$, by (9.14), $M \vDash v'<_x v$,

⟹ from $M \vDash Max^+(t',v,x)$, $M \vDash t<t'$,

⟹ $M \vDash t<t'\leq t$, contradicting $M \vDash t\in I\subseteq I_0$.

Therefore $M \vDash \neg(t'<tb)$.

⟹ by (4.6), $M \vDash tb\leq t'$.

So we established that

  (2b) $M \vDash (Max^+(t',v,x) \,\&\, Max^+T_b(t',v) \rightarrow tb\leq t')$.

Along with (2a) this shows that $M \vDash MMax^+T_b(tb,v,x)$, that is,

$$M \vDash MMax^+T_b(t_1,v,x).$$

From (1) and (2) we have that $M \vDash (Fr(x,t_1,ava,t_2) \rightarrow MMax^+T_b(t_1,v,x))$.

Along with $M \vDash Set(x)$ this shows that $M \vDash Special(x)$.

This completes the proof of (10.18).

(10.19) For any string concept $I \subseteq I_0$ there is a string concept $J \subseteq I$ such that

$QT^+ \vdash \forall u,v \in J(u \neq v \rightarrow \exists t',t'',z(Tally_b(t') \& Tally_b(t'') \& z=t'auat''avat'' \&$

$\& Set^*(z) \& \forall w(w \ \varepsilon \ z \leftrightarrow w=u \ v \ w=v)))$.

Let $J \equiv I_{9.30} \& I_{10.7} \& I_{10.15} \& I_{10.18}$.

Assume $M \vDash J(u) \& J(v) \& u \neq v$.

$\Rightarrow$ by (6.8), $M \vDash \exists!t_0 \ MinMax^+T_b(t_0,u) \& \exists!t_1 \ MinMax^+T_b(t_1,v)$,

$\Rightarrow$ by (8.2), $M \vDash u \prec v \ v \ v \prec u$.

Without loss of generality, we may assume that $M \vDash u \prec v$.

$\Rightarrow M \vDash u \triangleleft_{Tb} v \ v \ (u \approx_{Tb} v \ \& \ u \ll v)$.

(1) $M \vDash \vDash u \triangleleft_{Tb} v$.

$\Rightarrow M \vDash t_0 < t_1$.

Let $t'=t_0$ and $t''=t_1$.

Then $z=t_0 auat_1 avat_1$.

$\Rightarrow$ from $M \vDash Max^+T_b(t_0,u) \& Max^+T_b(t_1,v)$, $M \vDash Pref(aua,t_0) \& Pref(ava,t_1)$,

$\Rightarrow$ by (5.58), $M \vDash Env(t_1,z) \& \forall w(w \ \varepsilon \ z \leftrightarrow w=u \ v \ w=v)$,

$\Rightarrow$ by (10.7), $M \vDash MinSet(z)$,

$\Rightarrow$ from hypothesis $M \vDash u \prec v$, by (9.30), $M \vDash Lex^+(z)$.

Let $x=t_0 auat_0$ and $y=t_1 avat_1$.

$\Rightarrow$ by (10.15),

$M \vDash \exists t_2,x'(MinMax^+T_b(t_2,u) \& x'=t_2 auat_2 \& Special(x') \& \forall w(w \ \varepsilon \ x' \leftrightarrow w=u))$

and

$M \vDash \exists t_3, y'(\text{MinMax}^+T_b(t_3,v) \ \& \ y'=t_3 a v a t_3 \ \& \ \text{Special}(y') \ \& \ \forall w(w \ \varepsilon \ y' \leftrightarrow w=v))$,

$\Rightarrow$ from (6.8), $M \vDash t_0=t_2 \ \& \ t_1=t_3$,

$\Rightarrow M \vDash x=t_0 a u a t_0 = t_2 a u a t_2 = x' \ \& \ y=t_1 a v a t_1 = t_3 a v a t_3 = y'$,

$\Rightarrow M \vDash \text{Special}(x) \ \& \ \text{Special}(y)$,

$\Rightarrow M \vDash \text{Env}(t_1,z) \ \& \ z=t_0 a u a y \ \& \ \text{Env}(t_0,x) \ \& \ x=t_0 a u a t_0 \ \& \ \text{Env}(t_1,y) \ \&$

$\& \ \text{Firstf}(y,t_1,a v a,t_1) \ \& \ t_0 < t_1 \ \& \ \neg \exists w(w \ \varepsilon \ x \ \& \ w \ \varepsilon \ y) \ \&$

$\& \ \forall w_1, w_2(w_1 \ \varepsilon \ x \ \& \ w_2 \ \varepsilon \ y \to w_1 \prec w_2) \ \& \ \text{Special}(x) \ \& \ \text{Special}(y)$,

$\Rightarrow$ by (10.17), $M \vDash \text{Special}(z)$,

$\Rightarrow M \vDash \text{Set}^*(z)$, as required.

(2) $M \vDash u \approx_{T_b} v \ \& \ u \ll v$.

$\Rightarrow M \vDash t_0 = t_1$.

Let $t'=t_0$ and $t''=t_1 b$.

$\Rightarrow M \vDash z = t_0 a u a t_1 b a v a t_1 b$.

We obtain $M \vDash \text{Env}(t_1,z) \ \& \ \forall w(w \ \varepsilon \ z \leftrightarrow w=u \ \lor \ w=v)$ as in (1), and further that

$$M \vDash \text{MinSet}(z) \ \& \ \text{Lex}^+(z).$$

Now, for $x = t_0 a u a t_0$ we have $M \vDash \text{Special}(x)$ as in (1).

Let $y^+ = t_1 b a v a t_1 b$.

$\Rightarrow M \vDash \text{Env}(t_1 b, z) \ \& \ z = t_0 a u a y^+ \ \& \ \text{Env}(t_0,x) \ \& \ x=t_0 a u a t_0 \ \& \ \text{Max}^+T_b(t_1 b, v) \ \&$

$\& \ \neg(v \ \varepsilon \ x) \ \& \ \forall w(w \ \varepsilon \ x \to w \prec v) \ \& \ \text{Special}(x)$,

$\Rightarrow$ by (10.18), $M \vDash \text{Special}(z)$,

$\Rightarrow M \vDash \text{Set}^*(z)$, as required.

This completes the proof of (10.19).

(10.20) For any string concept $I \subseteq I_0$ there is a string concept $J \subseteq I$ such that

$QT^+ \vdash \forall x \in J \forall t, t_0, t'', y, v_0, z(Env(t,x)$ & $x = t_0 ayaz$ & $MinMax^+T_b(t_0,y)$ & $\neg(y \, \varepsilon \, z)$ &

& $Firstf(z, t_0 b, av_0 a, t'')$ & $\forall u(u \, \varepsilon \, z \to y \prec u)$ & $Special(z) \to Special(x))$.

Let $J \equiv I_{5.26}$ & $I_{10.15}$ & $I_{10.17}$ & $I_{10.18}$.

Assume $M \vDash Env(t,x)$ & $x = t_0 ayaz$ & $MinMax^+T_b(t_0,y)$

where $M \vDash \neg(y \, \varepsilon \, z)$ & $Firstf(z, t_0 b, av_0 a, t'')$ & $\forall u(u \, \varepsilon \, z \to y \prec u)$ & $Special(z)$

and $M \vDash J(x)$.

From $M \vDash MinMax^+T_b(t_0,y)$ we have, as in the proof of (10.18), that

$M \vDash Env(t_0, x')$ & $\forall w(w \, \varepsilon \, x' \leftrightarrow w = y)$,

where $x' = t_0 ayat_0$.

$\Rightarrow$ as in the proof of (10.15), $M \vDash Special(x')$.

We now claim that $M \vDash Env(t,z)$.

From $M \vDash Env(t_0, x')$ we have that $M \vDash MaxT_b(t_0, x')$, and from hypothesis

$M \vDash Special(z)$, $M \vDash Set(z)$.

$\Rightarrow M \vDash Env(t,x)$ & $x = t_0 ayaz$ & $x' = t_0 ayat_0$ & $MaxT_b(t_0, x')$ &

& $Firstf(z, t_0 b, av_0 a, t'')$ & $t_0 \prec t_0 b$ & $Set(z)$,

$\Rightarrow$ by (5.26), $M \vDash Env(t,z)$, as claimed.

We now have

$M \vDash Env(t,x)$ & $x = t_0 ayaz$ & $x' = t_0 ayat_0$ & $Env(t_0, x')$ & $Env(t,z)$ &

& $Firstf(z, t_0 b, av_0 a, t'')$ & $t_0 \prec t_0 b$ & $\neg \exists u(u \, \varepsilon \, x'$ & $u \, \varepsilon \, z)$ &

& $\forall u, v(u \, \varepsilon \, x'$ & $v \, \varepsilon \, z \to u \prec v)$ & $Special(x')$ & $Special(z)$.

$\implies$ by (10.17), $M \vDash \text{Special}(x)$, as required.

This completes the proof of (10.20).

(10.21) For any string concept $I \subseteq I_0$ there is a string concept $J \subseteq I$ such that

$QT^+ \vdash \forall z \in J \forall w,x,x',t,t_1,t_2,t_3,t_4,v_0(Env(t,z)\ \&\ z=t_1wx\ \&\ aBw\ \&\ aEw\ \&\ Env(t_2,x')\ \&$

$\&\ x'=t_1wt_2\ \&\ Env(t,x)\ \&\ Firstf(x,t_3,av_0a,t_4)\ \&\ t_2<t_3\ \rightarrow$

$\rightarrow \forall t',v(v\ \varepsilon\ x\ \&\ Max^+(t',v,z) \rightarrow Max^+(t',v,x)))$.

Let $J(x)$ be as in (9.18).

Assume $M \vDash Env(t,z)\ \&\ z=t_1wx\ \&\ Env(t,x)\ \&\ Firstf(x,t_3,av_0a,t_4)$ where

$M \vDash aBw\ \&\ aEw\ \&\ Env(t_2,x')\ \&\ x'=t_1wt_2\ \&\ t_2<t_3$ and $M \vDash J(z)$.

Suppose that $M \vDash v\ \varepsilon\ x\ \&\ Max^+(t',v,z)$.

Assume $M \vDash Fr(x,t_5,aua,t_6)\ \&\ u<_x v$.

From $M \vDash Env(t_2,x')$ we have $M \vDash MaxT_b(t_2,x')$.

If $M \vDash Tally_b(t_0)\ \&\ t_0 \subseteq_p t_1w_1$, then $M \vDash t_0 \subseteq_p x'$, whence $M \vDash t_0 \leq t_2$.

$\Rightarrow$ from $M \vDash Firstf(x,t_3,av_0a,t_4)\ \&\ Fr(x,t_5,aua,t_6)$, by (5.20), $M \vDash t_3 \leq t_5$,

$\Rightarrow M \vDash t_0 \leq t_2 < t_3 \leq t_5$.

So we established that $M \vDash Max^+T_b(t_5,t_1w)$.

Hence we have

$M \vDash Fr(x,t_5,aua,t_6)\ \&\ z=t_1wx\ \&\ aBw\ \&\ aEw\ \&\ Max^+T_b(t_5,t_1w)$.

$\Rightarrow$ by (5.25), $M \vDash Fr(z,t_5,aua,t_6)$,

$\Rightarrow$ from $M \vDash u<_x v$, by (9.18), $M \vDash u<_z v$,

$\Rightarrow$ from hypothesis $M \vDash Max^+(t',v,z)$, $M \vDash t_5 < t'$.

Hence we have shown that $M \vDash (Fr(x,t_5,aua,t_6)\ \&\ u<_x v \rightarrow t_5 < t')$.

Along with $M \vDash v\ \varepsilon\ x\ \&\ Set(x)\ \&\ Tally_b(t')$, which we have, this suffices for

$$M \vDash \text{Max}^+(t',v,x),$$

as required.

This completes the proof of (10.21).

(10.22) For any string concept $I \subseteq I_0$ there is a string concept $J \subseteq I$ such that

$QT^+ \vdash \forall x \in J \forall t, t_1, t_2, v, x'(Env(t,x) \& Intf(x,w_1,t_1,ava,t_2) \& x'=w_1at_1avat_1 \rightarrow$

$\rightarrow Env(t_1,x') \& \forall v', t'(v' \varepsilon x' \rightarrow (Max^+(t',v',x) \leftrightarrow Max^+(t',v',x)))$.

Let $J \equiv I_{5.6} \& I_{5.32} \& I_{9.9} \& I_{9.11}$.

Assume $M \vDash Env(t,x) \& Intf(x,w_1,t_1,ava,t_2)$ along with $M \vDash J(x)$.

Suppose that $M \vDash x'=w_1at_1avat_1$.

$\Rightarrow$ by the proof of (5.53), $M \vDash Env(t_1,x') \& Lastf(x',t_1,ava,t_1)$.

Assume $M \vDash v' \varepsilon x'$.

$\Rightarrow$ by (5.32), $M \vDash \exists t^+, w^*(Tally_b(t^+) \& x=x't^+w^*t \& aBw^* \& aEw^*)$.

Now, assume that $M \vDash Max^+(t',v',x)$.

$\Rightarrow M \vDash v' \varepsilon x \& Tally_b(t')$.

Assume also that $M \vDash Fr(x',t_3,aua,t_4) \& u<_{x'}v'$.

$\Rightarrow$ by (5.6) and (9.11), $M \vDash \exists t_5 \, Fr(x,t_3,aua,t_5) \& u<_x v'$,

$\Rightarrow$ from hypothesis $M \vDash Max^+(t',v',x)$, $M \vDash t_3<t'$.

Along with $M \vDash Set(x') \& v' \varepsilon x' \& Tally_b(t')$, what we just derived, namely

$\quad M \vDash (Fr(x',t_3,aua,t_4) \& u<_{x'}v' \rightarrow t_3<t')$,

is sufficient to establish $M \vDash Max^+(t',v',x')$.

Conversely, assume $M \vDash Max^+(t',v',x')$.

$\Rightarrow M \vDash Tally_b(t')$.

Assume that $M \vDash Fr(x,t_3,aua,t_4) \& u<_x v'$.

Now, we have that

$M \vDash Env(t,x)$ & $Intf(x,w_1,t_1,ava,t_2)$ & $x'=w_1at_1avat_1$ & $v' \varepsilon x'$ & $u<_x v'$.

$\Rightarrow$ by (9.9), $M \vDash u \varepsilon x'$,

$\Rightarrow M \vDash \exists t_6,t_7 \, Fr(x',t_6,aua,t_7)$,

$\Rightarrow$ by (9.11), $M \vDash u<_{x'} v'$,

$\Rightarrow$ from hypothesis $M \vDash Max^+(t',v',x')$, $M \vDash t_6<t'$,

$\Rightarrow$ by (5.6), $M \vDash \exists t_8 \, Fr(x,t_6,aua,t_8)$,

$\Rightarrow$ from $M \vDash Env(t,x)$, $M \vDash t_3=t_6$,

$\Rightarrow M \vDash t_3<t'$.

So we proved that $M \vDash (Fr(x,t_3,aua,t_4)$ & $u<_x v' \to t_3<t')$.

Along with $M \vDash Set(x)$ & $v' \varepsilon x$, where $M \vDash v' \varepsilon x$ follows from $M \vDash v' \varepsilon x'$ by (5.6), it follows that $M \vDash Max^+(t',v',x)$.

This completes the proof of (10.22).

(10.23) For any string concept $I \subseteq I_0$ there is a string concept $J \subseteq I$ such that

$QT^+ \vdash \forall x \in J \forall t_1,t_2,w_1,w_2,u,x'(Special(x)$ & $Intf(x,w_1,t_1,aua,t_2)$ &

& $x=w_1at_1auat_2aw_2$ & $x'=w_1at_1auat_1 \rightarrow Special(x'))$.

Let $J \equiv I_{5.51}$ & $I_{10.22}$.

Assume $M \vDash Special(x)$ & $Intf(x,w_1,t_1,aua,t_2)$ where

$M \vDash x=w_1at_1auat_2aw_2$ & $J(x)$.

Let $M \vDash x'=w_1at_1auat_1$.

$\Rightarrow$ from hypothesis $M \vDash Special(x)$, $M \vDash Set(x)$.

We have that $M \vDash x \neq aa$.

$\Rightarrow M \vDash \exists t\ Env(t,x)$,

$\Rightarrow$ by (10.22), $M \vDash Env(t_1,x')$,

$\Rightarrow M \vDash Set(x')$.

Assume $M \vDash Fr(x',t_3,ava,t_4)$.

$\Rightarrow$ by (5.51), $M \vDash \exists t_5\ Fr(x,t_3,ava,t_5)$,

$\Rightarrow$ from hypothesis $M \vDash Special(x)$, $M \vDash MMax^+T_b(t_3,v,x)$,

$\Rightarrow M \vDash Max^+(t_3,v,x)$,

$\Rightarrow$ by (10.22), $M \vDash Max^+(t_3,v,x')$.

Assume that $M \vDash Max^+(t',v,x')$ & $Max^+T_b(t',v)$.

$\Rightarrow$ by (10.22), $M \vDash Max^+(t',v,x)$ & $Max^+T_b(t',v)$,

$\Rightarrow$ from $M \vDash MMax^+T_b(t_3,v,x)$, $M \vDash t_3 \leq t'$.

Therefore, we derived $M \vDash (Max^+(t',v,x')$ & $Max^+T_b(t',v) \rightarrow t_3 \leq t')$.

Along with $M \vDash Max^+(t_3,v,x')$ this shows that $M \vDash MMax^+T_b(t_3,v,x')$.

So we established that $M \vDash (Fr(x',t_3,ava,t_4) \rightarrow MMax^+T_b(t_3,v,x'))$.

Along with $M \vDash Set(x')$ this yields $M \vDash Special(x')$.

This completes the proof of (10.23).

(10.24) For any string concept $I \subseteq I_0$ there is a string concept $J \subseteq I$ such that

$QT^+ \vdash \forall x \in J \forall t, t_0, t_1, t_2, t_3, x', x'', v_0, w^* (Env(t,x)$ & $x = t_0 w^* x''$ & $aBw^*$ & $aEw^*$ &

& $Env(t_1, x')$ & $x' = t_0 w^* t_1$ & $Env(t, x'')$ & $Firstf(x'', t_2, av_0 a, t_3)$ & $t_1 < t_2$ &

& $MinMax^+ T_b(t_2, v_0)$ & $Lex^+(x)$ & $Special(x)$ $\to$ $Lex^+(x'')$ & $Special(x'')$).

Let $J \equiv I_{5.41}$ & $I_{9.28}$.

Assume $M \vDash Env(t, x)$ & $x = t_0 w^* x''$ & $aBw^*$ & $aEw^*$

along with $M \vDash Env(t_1, x')$ & $x' = t_0 w^* t_1$ & $Env(t, x'')$ & $Firstf(x'', t_2, av_0 a, t_3)$ & $t_1 < t_2$

where $M \vDash J(x)$.

Suppose that $M \vDash MinMax^+ T_b(t_2, v_0)$ & $Lex^+(x)$ & $Special(x)$.

$\Rightarrow$ by (9.28), $M \vDash Lex^+(x'')$.

Assume $M \vDash Fr(x'', t_4, ava, t_5)$.

We want to show that $M \vDash MMax^+ T_b(t_4, v, x'')$.

First, we claim that

(1) $M \vDash Max^+ T_b(t_4, t_0 w^*)$.

Assume that $M \vDash Tally_b(t')$ & $t' \subseteq_p t_0 w^*$.

$\Rightarrow M \vDash t' \subseteq_p t_0 w^* \subseteq_p x'$,

$\Rightarrow$ from $M \vDash Env(t_1, x')$, $M \vDash MaxT_b(t_1, x')$,

$\Rightarrow M \vDash t' \leq t_1 < t_2$,

$\Rightarrow$ from $M \vDash Firstf(x'', t_2, av_0 a, t_3)$ & $Fr(x'', t_4, ava, t_5)$, by (5.20), $M \vDash t_2 \leq t_4$,

$\Rightarrow M \vDash t' < t_2 \leq t_4$.

Hence $M \vDash Max^+ T_b(t_4, t_0 w^*)$, as claimed.

From (1) along with $M \vDash Fr(x'',t_4,ava,t_5)$ & $x=t_0w*x''$ & $aBw*$ & $aEw*$

we have, by (5.25), that $M \vDash Fr(x,t_4,ava,t_5)$.

$\Rightarrow$ from hypothesis $M \vDash Special(x)$, $M \vDash MMax^+T_b(t_4,v,x)$.

Now, since $M \vDash t_1<t_2$ & $Tally_b(t_2)$, we have that $M \vDash \exists t''(t_1t''=t_2$ & $Tally_b(t''))$.

$\Rightarrow$ from $M \vDash Env(t,x'')$, by (5.11),

$\qquad M \vDash \exists t^*,w_1(x''=t^*w_1t$ & $aBw_1$ & $aEw_1$ & $Tally_b(t^*))$,

$\Rightarrow$ from $M \vDash Firstf(x'',t_2,av_0a,t_3)$, $M \vDash (t_2a)Bx''$,

$\Rightarrow M \vDash \exists x_1,w_2\ t_2ax_1=x''=t^*(aw_2)t$,

$\Rightarrow$ by (4.23$^b$), $M \vDash t_2=t^*$,

$\Rightarrow M \vDash x''=t^*w_1t=t_2w_1t=t_1t''w_1t$,

$\Rightarrow M \vDash x=t_0w*x''=t_0w*t_1t''w_1t$,

$\Rightarrow$ by (5.41), $M \vDash \neg\exists w(w\ \varepsilon\ x'$ & $w\ \varepsilon\ x'')$,

$\Rightarrow$ by (5.46), $M \vDash \forall w(w\ \varepsilon\ x \leftrightarrow w\ \varepsilon\ x'$ & $w\ \varepsilon\ x'')$.

Now, from $M \vDash MMax^+T_b(t_4,v,x)$, $M \vDash Max^+(t_4,v,x)$.

$\Rightarrow$ since $M \vDash v\ \varepsilon\ x''$, by (10.21), we have that

(2) $M \vDash Max^+(t_4,v,x'')$.

Assume $M \vDash Max^+(t',v,x'')$ & $Max^+T_b(t',v)$.

We want to show that $M \vDash t_4 \leq t'$.

We claim that

(3) $M \vDash Max^+(t',v,x)$.

We already have that $M \vDash v\ \varepsilon\ x$ & $Set(x)$ & $Tally_b(t')$.

Assume $M \vDash Fr(x,t_6,aua,t_7)$ & $u<_xv$.

$\Rightarrow$ M ⊨ u ε x,

$\Rightarrow$ M ⊨ u ε x' v u ε x".

Suppose that (3a)  M ⊨ u ε x'.

$\Rightarrow$ M ⊨ $\exists t_8,t_9$ Fr(x', $t_8$,aua,$t_9$),

$\Rightarrow$ since  M ⊨ x=$t_0$w*$t_1$t"$w_1$t=x't"$w_1$t, by (5.6),  M ⊨ $\exists t_{10}$ Fr(x,$t_8$,aua,$t_{10}$),

$\Rightarrow$ from  M ⊨ Env(t,x),  M ⊨ $t_8$=$t_6$,

$\Rightarrow$ from  M ⊨ MaxT$_b$($t_1$,x'),  M ⊨ $t_6$=$t_8 \leq t_1 < t_2$,

$\Rightarrow$ from  M ⊨ Firstf(x",$t_2$,a$v_0$a,$t_3$) & Fr(x",$t_4$,ava,$t_5$),  M ⊨ $v_0 \leq_{x"} v$,

$\Rightarrow$ from  M ⊨ Lex$^+$(x"),  M ⊨ $v_0 \leqslant v$,

$\Rightarrow$ from hypothesis  M ⊨ MinMax$^+$T$_b$($t_2$,$v_0$) & Max$^+$T$_b$(t',v),  M ⊨ $t_2 \leq t'$,

$\Rightarrow$ M ⊨ $t_6 < t_2 \leq t'$.

Suppose that (3b)  M ⊨ u ε x".

$\Rightarrow$ M ⊨ $\exists t_8,t_9$ Fr(x", $t_8$,aua,$t_9$),

$\Rightarrow$ from  M ⊨ u<$_x$v, by (9.18),  M ⊨ u$\leq_{x"} v$,

$\Rightarrow$ from  M ⊨ Max$^+$(t',v,x"),  M ⊨ $t_6$=$t_8 < t'$.

In (3a)-(3b)  we thus proved   M ⊨ (Fr(x,$t_6$,aua,$t_7$) & u<$_x$v  →  $t_6$<t').

This suffices to prove (3) under the hypothesis

$\quad\quad$ M ⊨ Max$^+$(t',v,x") & Max$^+$T$_b$(t',v).

Since we have that  M ⊨ MMax$^+$T$_b$($t_4$,v,x),  from (3) and  M ⊨ Max$^+$T$_b$(t',v) it follows that  M ⊨ $t_4 \leq t'$.

Therefore, we proved that   M ⊨ (Max$^+$(t',v,x") & Max$^+$T$_b$(t',v) → $t_4 \leq t'$).

Along with (2) this establishes that   M ⊨ MMax$^+$T$_b$($t_4$,v,x).

$\implies$ M ⊨ Lex$^+$(x″) & Special(x″).

This completes the proof of (10.24).

(10.25) For any string concept $I \subseteq I_0$ there is a string concept $J \subseteq I$ such that

$QT^+ \vdash \forall x \in J \forall t, t_0, t_1, t_2, t_3, t^{++}, u_0, v_0, w^*, x', x''($Env$(t,x)$ & $x = t_0 w^* x''$ & $aBw^*$ & $aEw^*$ &

& Env$(t_1, x')$ & $x' = t_0 w^* t_1$ & Env$(t, x'')$ & Lastf$(x', t_1, au_0 a, t_1)$ &

& Firstf$(x'', t_2, av_0\, a, t_3)$ & MinMax$^+$T$_b(t^{++}, v_0)$ & $t_1 < t^{++}$ & Special$(x)$ →

→ $t_2 = t^{++}$).

Let $J \equiv I_{5.41}$ & $I_{5.46}$ & $I_{9.11}$ & $I_{9.14}$.

Assume  $M \vDash$ Env$(t,x)$ & $x = t_0 w^* x''$ & $aBw^*$ & $aEw^*$ & Env$(t_1, x')$ & $x' = t_0 w^* t_1$ &

& Env$(t, x'')$ & Lastf$(x', t_1, au_0 a, t_1)$ & Firstf$(x'', t_2, av_0\, a, t_3)$

along with  $M \vDash J(x)$.

Suppose further that  $M \vDash$ MinMax$^+$T$_b(t^{++}, v_0)$ & $t_1 < t^{++}$ & Special$(x)$.

⇒ from  $M \vDash$ Env$(t, x'')$, by (5.11),

$M \vDash \exists t'', w^{**}(x'' = t'' w^{**} t$ & Tally$_b(t'')$ & $aBw^{**}$ & $aEw^{**})$,

⇒ from  $M \vDash$ Firstf$(x'', t_2, av_0\, a, t_3)$,  $M \vDash (t_2 a) B x''$,

⇒ $M \vDash \exists x_1, w_1\, t_2 a x_1 = x'' = t'' a w_1 t$,

⇒ by (4.23$^b$),  $M \vDash t_2 = t''$,

⇒ from  $M \vDash$ Firstf$(x'', t_2, av_0\, a, t_3)$,  $M \vDash$ Max$^+$T$_b(t_2, v_0)$,

⇒ from  $M \vDash$ MinMax$^+$T$_b(t^{++}, v_0)$,  $M \vDash t^{++} \leq t_2$,

⇒ from hypothesis  $M \vDash t_1 < t^{++}$,  $M \vDash t_1 < t_2$,

⇒  $M \vDash \exists t^*($Tally$_b(t^*)$ & $t_1 t^* = t_2)$,

⇒ $M \vDash x = t_0 w^* x'' = t_0 w^* t_2 w^{**} t = t_0 w^* t_1 t^* w^{**} t_2$,

⇒ by (5.41),  $M \vDash \neg \exists w(w\, \varepsilon\, x'$ & $w\, \varepsilon\, x'')$,

$\Rightarrow$ by (5.46), $M \vDash \forall w(w \, \varepsilon \, x \leftrightarrow w \, \varepsilon \, x' \, \& \, w \, \varepsilon \, x'')$.

<u>Claim 1.</u>  $M \vDash \forall u(u <_x v_0 \rightarrow u \leq_x u_0)$.

Assume  $M \vDash u <_x v_0$.

$\Rightarrow M \vDash u \, \varepsilon \, x$,

$\Rightarrow M \vDash u \, \varepsilon \, x' \lor u \, \varepsilon \, x''$.

Suppose, for a reductio, that  $M \vDash u \, \varepsilon \, x''$.

$\Rightarrow$ from  $M \vDash \text{Firstf}(x'', t_2, av_0 \, a, t_3)$,  by (9.1) and (9.7),  $M \vDash v_0 \leq_{x''} u$,

$\Rightarrow$ by (9.18),  $M \vDash v_0 \leq_x u$,  which, by (9.4) and (9.6), contradicts hypothesis

$M \vDash u <_x v_0$.

Therefore,  $M \vDash \neg u \, \varepsilon \, x''$.

$\Rightarrow M \vDash u \, \varepsilon \, x'$,

$\Rightarrow$ from  $M \vDash \text{Lastf}(x', t_1, au_0 a, t_1)$, by (9.3),  $M \vDash u \leq_{x'} u_0$,

$\Rightarrow$ by (9.11),  $M \vDash u \leq_x u_0$,  as required.

<u>Claim 2.</u>  $M \vDash \text{Max}^+ T_b(t_2, t_0 w^*)$.

Suppose  $M \vDash \text{Tally}_b(t') \, \& \, t' \subseteq_p t_0 w^*$.

$\Rightarrow$ from  $M \vDash \text{Env}(t_1, x')$,  $M \vDash \text{MaxT}_b(t_1, x')$,

$\Rightarrow M \vDash t' \subseteq_p t_0 w^* \subseteq_p x'$,

$\Rightarrow M \vDash t' \leq t_1 < t_2$.

Therefore,  $M \vDash \text{Max}^+ T_b(t_2, t_0 w^*)$.

<u>Claim 3.</u>  $M \vDash \forall t'(t_1 < t' < t_2 \rightarrow \text{Max}^+(t', v_0, x))$.

Assume  $M \vDash t_1 < t' < t_2$.

We have that  $M \vDash \text{Fr}(x'', t_2, av_0 \, a, t_3)$.

$\Rightarrow$ from Claim 2, by (5.25), $M \vDash Fr(x,t_2,av_0a,t_3)$,

$\Rightarrow$ $M \vDash v_0 \, \varepsilon \, x$.

We have $M \vDash Set(x)$ by hypothesis.

$\Rightarrow$ from $M \vDash Tally_b(t_1)$ & $Tally_b(t_2)$ & $t_1 < t' < t_2$, $M \vDash Tally_b(t')$.

Assume $M \vDash Fr(x,t_4,aua,t_5)$ & $u <_x v_0$.

$\Rightarrow$ by Claim 1, $M \vDash u \leq_x u_0$,

$\Rightarrow$ by (9.14), $M \vDash t_4 \leq t_1 < t'$.

Therefore, $M \vDash (Fr(x,t_4,aua,t_5)$ & $u <_x v_0 \rightarrow t_4 < t')$, which suffices to show $M \vDash Max^+(t',v_0,x)$.

Suppose now, for a reductio, that $M \vDash t_2 \neq t^{++}$.

$\Rightarrow$ $M \vDash t^{++} < t_2$,

$\Rightarrow$ $M \vDash t_1 < t^{++} < t_2$,

$\Rightarrow$ by Claim 3, $M \vDash Max^+(t^{++},v_0,x)$,

$\Rightarrow$ $M \vDash Max^+(t^{++},v_0,x)$ & $Max^+T_b(t^{++},v_0)$,

$\Rightarrow$ from hypothesis $M \vDash Special(x)$ and $M \vDash Fr(x,t_2,av_0a,t_3)$,

$\qquad M \vDash MMax^+T_b(t_2,v_0,x)$,

$\Rightarrow$ $M \vDash t_2 \leq t^{++}$,

$\Rightarrow$ $M \vDash t_2 \leq t^{++} < t_2$, contradicting $M \vDash t_2 \in I \subseteq I_0$.

Therefore $M \vDash t_2 = t^{++}$.

This completes the proof of (10.25).

(10.26) For any string concept $I \subseteq I_0$ there is a string concept $J \subseteq I$ such that

$$QT^+ \vdash \forall x \in J \forall t_1,t_2,x(Set(x) \ \& \ Fr(x,t_1,ava,t_2) \ \rightarrow \ Max^+(t_1,v,x)).$$

Let $J \equiv I_{9.14}$.

Assume $M \vDash Set(x) \ \& \ Fr(x,t_1,ava,t_2)$ where $M \vDash J(x)$.

$\Rightarrow M \vDash v \ \varepsilon \ x \ \& \ Tally_b(t_1)$.

Assume that $M \vDash Fr(x,t_3,aua,t_4) \ \& \ u<_x v$.

$\Rightarrow$ by (9.14), $M \vDash t_3 < t_1$.

This suffice to show that $M \vDash Max^+(t_1,v,x)$, and completes the proof of (10.26).

(10.27) For any string concept $I \subseteq I_0$ there is a string concept $J \subseteq I$ such that

$QT^+ \vdash \forall x, z \in J \forall v,t,t',t''[\text{MinSet}(x) \& \text{Set}(z) \&$

$\quad \& \forall w, t_1, t_2 (\text{Free}^+(x, t_1, awa, t_2) \rightarrow \exists t_3 \text{Fr}(z, t_1 b, awa, t_3)) \&$

$\& \forall w, t_1, t_2 (\text{Bound}(x, t_1, awa, t_2) \vee \text{Free}^-(x, t_1, awa, t_2) \rightarrow \exists t_3 \text{Fr}(z, t_1, awa, t_3)) \&$

$\quad \& \forall u, w (u \, \varepsilon \, x \& v \, \varepsilon \, x \rightarrow (u <_x w \rightarrow u <_z w)) \& \text{Max}^+(t', v, z) \&$

$\quad \& \text{Free}^+(x, t_6, ava, t_7) \& \neg \text{Firstf}(x, t_6, ava, t_7) \rightarrow tb \leq t']$.

Let $J \equiv (I_{5.32} \& I_{9.13} \& I_{9.17} \& I_{10.8})_{SUB}$.

As usual, we may assume that J is downward closed with respect to under $\subseteq_p$.

Assume that $M \vDash \text{MinSet}(x) \& \text{Set}(z) \& \text{Max}^+(t', v, z)$

and

$\quad M \vDash \forall w, t_1, t_2 (\text{Free}^+(x, t_1, awa, t_2) \rightarrow \exists t_3 \text{Fr}(z, t_1 b, awa, t_3)) \&$

$\& \forall w, t_1, t_2 (\text{Bound}(x, t_1, awa, t_2) \vee \text{Free}^-(x, t_1, awa, t_2) \rightarrow \exists t_3 \text{Fr}(z, t_1, awa, t_3)) \&$

$\quad \& \forall u, w (u \, \varepsilon \, x \& v \, \varepsilon \, x \rightarrow (u <_x w \rightarrow u <_z w))$

along with $M \vDash \text{Free}^+(x, t_6, ava, t_7) \& \neg \text{Firstf}(x, t_6, ava, t_7) \& J(x) \& J(z)$.

$\Rightarrow M \vDash \exists t_3 \text{Fr}(z, tb, ava, t'')$,

$\Rightarrow M \vDash \exists w_1 (\text{Intf}(x, w_1, t, ava, t'') \vee (\text{Lastf}(x, t, ava, t'') \& t = t'' \& x = w_1 atavat''))$.

If $M \vDash \text{Intf}(x, w_1, t, ava, t'')$, then, by the proof of (5.53),

$\quad M \vDash \exists x_1 (\text{Env}(t, x_1) \& x_1 = w_1 atavat \& \text{Lastf}(x_1, t, ava, t) \& x_1 Bx)$.

$\Rightarrow$ from $M \vDash \text{MinSet}(x)$, by (10.8), $M \vDash \text{MinSet}(x_1)$,

$\Rightarrow M \vDash \exists x' (\text{MinSet}(x') \& \text{Env}(t, x') \& \text{Lastf}(x', t, ava, t) \& (x'Bx \vee x' = x))$,

$\Rightarrow$ by the proof of SUBTRACTION LEMMA,

$M \vDash \exists x^-, t^- (Env(t^-,x^-) \& \forall w(w \varepsilon x^- \leftrightarrow w \varepsilon x' \vee w \neq v) \& t^- < t \& x^- Bx' \&$

$\& \exists v' Lastf(x^-,t^-,av'a,t^-)),$

$\implies M \vDash Tally_b(t^-) \& Tally_b(t) \& t^- < t,$

$\implies M \vDash \exists t_0 (Tally_b(t_0) \& t = t^- t_0 \& x' = x^- t_0 avat^- t_0),$

$\implies$ from $M \vDash Env(t^-,x^-) \& Env(t,x') \& Tally_b(t^-) \& Lastf(x^-,t^-,av'a,t^-) \&$

$\& Lastf(x',t,ava,t),$

by (9.13), $M \vDash v' <_{x'} v \& \neg \exists w(v' <_{x'} w \& w <_{x'} v),$

$\implies$ from $M \vDash MinSet(x)$, $M \vDash Set(x),$

$\implies$ from $M \vDash v \varepsilon x$, $M \vDash x \neq aa,$

$\implies$ by (5.18), $M \vDash \exists t^* Env(t^*,x),$

$\implies$ from $M \vDash Intf(x,w_1,t,ava,t'') \& x' = w_1 atavat$, by (5.32),

$M \vDash \exists t',w'(Tally_b(t') \& x = x't'w't^* \& aBw' \& aEw'),$

$\implies$ from $M \vDash Env(t^*,x) \& Tally_b(t^*)$, by (9.11), $M \vDash v' <_x v.$

If $M \vDash v' <_x w \& w <_x v$, then $M \vDash v' <_{x'} w \& w <_{x'} v$, a contradiction.

Hence also $M \vDash \neg \exists w(v' <_x w \& w <_x v).$

Now, we have that $M \vDash x = x' \vee x = x't'w't^*$, that is,

$M \vDash x = x' = x^- t_0 avat^- t_0 = x^- t_0 avat \vee x = x't'w't^* = (x^- t_0 avat^- t_0)t'w't^*.$

$\implies$ from $M \vDash Env(t^-,x^-) \& Fr(x^-,t^-,av'a,t^-) \& Tally_b(t_0) \& Tally_b(t^*) \& Tally_b(t),$

by (5.6), $M \vDash \exists t_4 Fr(x,t^-,av'a,t_4),$

$\implies$ from $M \vDash Free^+(x,t,ava,t'') \& \neg Firstf(x,t,ava,t'')$, $M \vDash Free(x,t,ava,t''),$

$\implies M \vDash t = t^- b,$

$\implies$ from $M \vDash Free^+(x,t,ava,t'') \& Fr(x,t^-,av'a,t_4) \& v' <_x v$, by (9.17),

$$M \vDash \text{Free}^+(x,t^-,av'a,t_4),$$

$\Longrightarrow$ by principal hypothesis,  $M \vDash \exists t_5\, \text{Fr}(z,t^-b,av'a,t_5)$,

$\Longrightarrow$ from  $M \vDash t=t^-b$,  $M \vDash \text{Fr}(z,t,av'a,t_5)$,

$\Longrightarrow$ from principal hypothesis and $M \vDash v'<_x v$,  $M \vDash v'<_z v$,

$\Longrightarrow$ from hypothesis  $M \vDash \text{Max}^+(t',v,z)$,  $M \vDash t<t'$,

$\Longrightarrow$ from  $M \vDash \text{Tally}_b(t')$, by (4.7),  $M \vDash tb \leq t'$,  as required.

This completes the proof of (10.27).

(10.28) For any string concept $I \subseteq I_0$ there is a string concept $J \subseteq I$ such that

$QT^+ \vdash \forall z, x \in J \forall x', x^*, t, t^*, t'', t_0, t^+, w, v_0 [Env(t,z) \ \& \ z = t_0 w x^* \ \& \ aBw \ \& \ aEw \ \&$

$\& \ Set(x^*) \ \& \ x' = t_0 w t^+ \ \& \ t^+ = t^* \ \& \ Env(t^+, x') \ \& \ Firstf(x^*, t^*b, av_0a, t'') \ \&$

$\& \ \forall u,v (u \ \varepsilon \ x' \ \& \ v \ \varepsilon \ x^* \rightarrow u \prec v) \ \& \ \neg \exists u (u \ \varepsilon \ x' \ \& \ u \ \varepsilon \ x^*) \ \&$

$\& \ Set(x) \ \& \ x \sim x^* \ \& \ \forall u,v (u <_x v \leftrightarrow u <_{x^*} v) \ \&$

$\& \ \forall v, t_1, t_2 \ (Free^+(x, t_1, ava, t_2) \rightarrow \exists t_3 Fr(x^*, t_1 b, ava, t_3)) \ \&$

$\& \ \forall v, t_1, t_2 \ (Bound(x, t_1, ava, t_2) \lor Free^-(x, t_1, ava, t_2) \rightarrow \exists t_3 Fr(x^*, t_1, ava, t_3)) \ \&$

$\& \ Special(x) \ \& \ Special(x') \rightarrow Special(z)].$

Let $J \equiv I_{5.26} \ \& \ I_{10.17} \ \& \ I_{10.27}$.

Assume $M \vDash Env(t,z) \ \& \ z = t_0 w x^* \ \& \ aBw \ \& \ aEw \ \& \ Set(x^*) \ \& \ t^+ = t^*$

along with $M \vDash Env(t^+, x') \ \& \ x' = t_0 w t^+ \ \& \ Firstf(x^*, t^*b, av_0a, t'') \ \&$

$\& \ \forall u,v (u \ \varepsilon \ x' \ \& \ v \ \varepsilon \ x^* \rightarrow u \prec v) \ \& \ \neg \exists u (u \ \varepsilon \ x' \ \& \ u \ \varepsilon \ x^*)$

and

$M \vDash Set(x) \ \& \ x \sim x^* \ \& \ \forall u,v (u <_x v \leftrightarrow u <_{x^*} v) \ \&$

$\& \ \forall v, t_1, t_2 \ (Free^+(x, t_1, ava, t_2) \rightarrow \exists t_3 Fr(x^*, t_1 b, ava, t_3)) \ \&$

$\& \ \forall v, t_1, t_2 \ (Bound(x, t_1, ava, t_2) \lor Free^-(x, t_1, ava, t_2) \rightarrow \exists t_3 Fr(x^*, t_1, ava, t_3))$.

Assume that $M \vDash Special(x) \ \& \ Special(x')$.

Assume $M \vDash Fr(z, t_1, ava, t_2)$.

$\Rightarrow$ from $M \vDash Env(t^+, x')$, $M \vDash Max^+T_b(t^+, x')$,

$\Rightarrow M \vDash Env(t,z) \ \& \ z = t_0 w x^* \ \& \ aBw \ \& \ aEw \ \& \ x' = t_0 w t^+ \ \& \ Max^+T_b(t^+, x') \ \&$

$\& \ Firstf(x^*, t^*b, av_0a, t'') \ \& \ t^+ < t^*b \ \& \ Set(x^*),$

⟹ by (5.26), $M \vDash Env(t,x^*)$,

⟹ from $M \vDash Fr(z,t_1,ava,t_2)$, by (5.39),

$M \vDash \exists t_3 Fr(x',t_1,ava,t_3) \lor Fr(x^*,t^*b,ava,t'') \lor \exists t_3 Fr(x^*,t_1,ava,t_3)$.

(1) $M \vDash \exists t_3 Fr(x',t_1,ava,t_3)$.

We argue exactly as in the proof of (10.17),(1), that $M \vDash MMax^+T_b(t_1,v,z)$.

(2) $M \vDash Fr(x^*,t^*b,ava,t'')$.

We first show (2a) $M \vDash Max^+(t^*b,v,z)$.

We already have that $M \vDash v \, \varepsilon \, z \, \& \, Set(z) \, \& \, Tally_b(t^*b)$.

Assume $M \vDash Fr(z,t_4,aua,t_5) \, \& \, u<_z v$.

Let $M \vDash Tally_b(t_6) \, \& \, t_6 \subseteq_p t_0 w$.

⟹ $M \vDash t_6 \subseteq_p x'$,

⟹ from $M \vDash Max^+T_b(t^+,x')$, $M \vDash t_6 \leq t^+ < t^*b$,

⟹ $M \vDash Max^+T_b(t^*b,t_0 w)$,

⟹ $M \vDash Fr(x^*,t^*b,ava,t'') \, \& \, z=t_0 w x^* \, \& \, aBw \, \& \, aEw \, \& \, Max^+T_b(t^*b,t_0 w)$,

⟹ by (5.25), $M \vDash Fr(z,t^*b,ava,t'')$,

⟹ by (10.26), $M \vDash Max^+(t^*b,v,z)$.

Hence (2a) holds.

Assume now $M \vDash Max^+(t',v,z) \, \& \, Max^+T_b(t',v)$.

⟹ from $M \vDash Env(t^+,x')$, $M \vDash \exists v' \, Lastf(x',t^+,av'a,t^+)$.

We claim that $M \vDash \exists t_3 \, Fr(z,t^+,av'a,t_3)$.

First note that from $M \vDash Firstf(x^*,t^*b,av_0 a,t'')$, we have $M \vDash (t^*ba)Bx^*$, whence $M \vDash \exists x_1 \, t^*bax_1 = x^*$.

On the other hand, from $M \models Env(t,x^*)$, by (5.11),

$$M \models \exists x_2, t_6 \, (Tally_b(t_6) \,\&\, x^*=t_6 x_2 t \,\&\, aBx_2 \,\&\, aEx_2),$$

whence $M \models \exists x_3 \, x^* = t_6(ax_3)t$.

$\Rightarrow M \models t_6(ax_3)t = x^* = t^*bax_1$,

$\Rightarrow$ since $M \models Tally_b(t_6) \,\&\, Tally_b(t^*b)$, by (4.23$^b$), $M \models t_6 = t^*b = t^+b$,

$\Rightarrow M \models z = t_0 w x^* = t_0 w(t_6 x_2 t) = (t_0 w \, t^+)b \, x_2 t = x'b \, x_2 t$,

$\Rightarrow M \models Env(t^+,x') \,\&\, Fr(x',t^+,av'a,t^+) \,\&\, x'b \, x_2 t = z \,\&\, aBx_2 \,\&\, aEx_2 \,\&\, Tally_b(b) \,\&$

$\&\, Tally_b(t)$,

$\Rightarrow$ by (5.6), $M \models \exists t_3 \, Fr(z,t^+,av'a,t_3)$, as claimed.

Now, from $M \models Set(z) \,\&\, Fr(z,t^+b,ava,t'')$, by (9.14), $M \models v' <_z v$.

$\Rightarrow$ from hypothesis $M \models Max^+(t',v,z)$, $M \models t^+ < t'$,

$\Rightarrow$ since $M \models Tally_b(t')$, by (4.7), $M \models t^*b \leq t'$.

So we have established (2b) $M \models (Max^+(t',v,z) \,\&\, Max^+T_b(t',v) \rightarrow t^*b \leq t')$.

Along with (2a) this gives $M \models MMax^+T_b(t^*b,v,z)$.

(3) $M \models \exists t_3 Fr(x^*,t_1,ava,t_3)$.

To show (3a) $M \models Max^+(t_1,v,z)$, we use (10.26) and the hypothesis

$M \models Fr(z,t_1,ava,t_2)$.

Assume now that $M \models Max^+(t',v,z) \,\&\, Max^+T_b(t',v)$.

$\Rightarrow$ from $M \models Firstf(x^*,t^*b,av_0 a,t'') \,\&\, Fr(x^*,t_1,ava,t_3)$, by (5.20), $M \models t^*b \leq t_1$.

If $M \models v_0 = v$, then $M \models t_1 = t^*b$.

Then we proceed exactly as in (2b) to derive $M \models t_1 = t^*b \leq t'$.

Suppose that $M \models v_0 \neq v$.

$\Rightarrow$ from $M \vDash \text{Firstf}(x^*,t^*b,av_0a,t'') \ \& \ \text{Fr}(x^*,t_1,ava,t_3)$, by (9.1) and (9.7),

$$M \vDash v_0<_{x^*}v,$$

$\Rightarrow$ from principal hypothesis, $M \vDash v_0 \ \varepsilon \ x \ \& \ v \ \varepsilon \ x \ \& \ v_0<_xv,$

$\Rightarrow M \vDash \exists t_6,t_7 \text{Fr}(x,t_6,ava,t_7),$

$\Rightarrow$ by (9.1), $M \vDash \neg \text{Firstf}(x,t_6,ava,t_7),$

$\Rightarrow$ by (9.22), $M \vDash \text{Free}(x,t_6,ava,t_7) \ v \ \text{Bound}(x,t_6,ava,t_7).$

Suppose that

  (3bi)  $M \vDash \text{Free}(x,t_6,ava,t_7) \ \& \ \text{Free}^+(x,t_6,ava,t_7).$

To establish (3b), we proceed under the hypothesis

$$M \vDash \text{Max}^+(t',v,z) \ \& \ \text{Max}^+\text{T}_b(t',v).$$

$\Rightarrow M \vDash \text{Tally}_b(t'),$

$\Rightarrow$ by (4.2), $M \vDash t'=b \ v \ \exists t^- \ (\text{Tally}_b(t^-) \ \& \ t'=t^-b).$

Suppose, for a reductio, $M \vDash t'=b.$

$\Rightarrow$ from $M \vDash \text{Max}^+(t',v,z)$, $M \vDash \forall u,t_8,t_9 \ (\text{Fr}(z,t_8,aua,t_9) \ \& \ u<_zv \ \rightarrow \ t_8<t').$

Now, we have that $M \vDash \text{Fr}(z,t^*b,av_0a,t'')$ from $M \vDash \text{Firstf}(x^*,t^*b,av_0a,t'')$

as in (1). Since we also have

$M \vDash \text{Env}(t,z) \ \& \ z=t_0wx^* \ \& \ aBw \ \& \ aEw \ \& \ \text{Env}(t^+,x') \ \& \ x'=t_0wt^+ \ \& \ \text{Env}(t,x^*) \ \&$

$$\& \ \text{Firstf}(x^*,t^*b,av_0a,t'') \ \& \ t^+<t^*b \ \& \ v_0 \ \varepsilon \ x^* \ \& \ v \ \varepsilon \ x^*,$$

by (9.18), we have that $M \vDash v_0<_{x^*}v \leftrightarrow v_0<_zv.$

$\Rightarrow$ from $M \vDash v_0<_{x^*}v$, $M \vDash v_0<_zv,$

$\Rightarrow M \vDash t^*b<t'=b,$ a contradiction.

Therefore, $M \vDash t'\neq b,$ and we have that $M \vDash \exists t^-(\text{Tally}_b(t^-) \ \& \ t'=t^-b).$

$\Rightarrow$ from hypothesis (3bi), $M \vDash \exists t_8\ Fr(x^*,t_6b,ava,t_8)$,

$\Rightarrow$ from $M \vDash Env(t,x^*)\ \&\ Fr(x^*,t_1,ava,t_3)$, $M \vDash t_1=t_6b$.

Also, from hypothesis $M \vDash Special(x)$, $M \vDash MMax^+T_b(t_6,v,x)$.

$\Rightarrow M \vDash \forall t^{**}(Max^+(t^{**},v,x)\ \&\ Max^+T_b(t^{**},v) \to t_6 \leq t^{**})$.

We claim that to show $M \vDash t_1 \leq t'$, it is sufficient to establish that

$M \vDash Max^+(t^-,v,x)\ \&\ Max^+T_b(t^-,v)$.

To show that $M \vDash Max^+(t^-,v,x)$ note that we already have that

$$M \vDash v\ \varepsilon\ x\ \&\ Tally_b(t^-).$$

Assume that $M \vDash Fr(x,t_9,aua,t_{10})\ \&\ u<_x v$.

$\Rightarrow$ from hypothesis $M \vDash Free^+(x,t_6,ava,t_7)$, by (9.17), $M \vDash Free^+(x,t_9,aua,t_{10})$,

$\Rightarrow M \vDash \exists t_{11}\ Fr(x^*,t_9b,aua,t_{11})\ \&\ u<_{x^*}v$,

$\Rightarrow$ as in (1), $M \vDash Fr(z,t_9b,aua,t_{11})$,

$\Rightarrow$ using (9.18), $M \vDash u<_z v$,

$\Rightarrow$ from hypothesis $M \vDash Max^+(t',v,z)$, $M \vDash t_9b<t'=t^-b$,

$\Rightarrow$ by (4.9), $M \vDash t_9<t^-$.

This suffices to establish $M \vDash Max^+(t^-,v,x)$.

To see that $M \vDash Max^+T_b(t^-,v)$, note that we have

$M \vDash MinSet(x)\ \&\ \forall w,t_1,t_2\ (Free^+(x,t_1,awa,t_2) \to \exists t_3 Fr(z,t_1b,awa,t_3))\ \&$

$\&\ \forall w,t_1,t_2\ (Bound(x,t_1,awa,t_2)\ v\ Free^-(x,t_1,awa,t_2) \to \exists t_3 Fr(z,t_1,awa,t_3))\ \&$

$\&\ \forall u,w(u\ \varepsilon\ x\ \&\ v\ \varepsilon\ x \to (u<_x w \leftrightarrow u<_z w))\ \&\ Max^+(t',v,z)\ \&$

$\&\ Free^+(x,t_6,ava,t_7)\ \&\ \neg Firstf(x,t_6,ava,t_7)$.

$\Rightarrow$ by (10.27), $M \vDash t_6b \leq t'=t^-b$,

$\Rightarrow$ by (4.9), $M \vDash t_6 \leq t^-$.

Since $M \vDash Fr(x,t_6,ava,t_7)$, this suffices to establish $M \vDash Max^+T_b(t^-,v)$.

This completes the proof of (3b) under the hypothesis (3bi).

Suppose that

(3bii) $M \vDash Bound(x,t_6,ava,t_7)$ & $Free^-(x,t_6,ava,t_7)$.

Again, to show (3b), assume that

$$M \vDash Max^+(t',v,z) \ \& \ Max^+T_b(t',v).$$

From hypothesis $M \vDash Special(x)$, $M \vDash MMax^+T_b(t_6,v,x)$.

$\Rightarrow M \vDash \forall t''(Max^+(t'',v,x) \ \& \ Max^+T_b(t'',v) \to t_6 \leq t'')$,

$\Rightarrow$ from $M \vDash Fr(x,t_6,ava,t_7)$ and principal hypothesis, $M \vDash \exists t_8 Fr(x^*,t_6,ava,t_8)$,

$\Rightarrow$ from $M \vDash Env(t,x^*)$ & $Fr(x^*,t_1,ava,t_3)$, $M \vDash t_1 = t_6$.

Hence to show that $M \vDash t_1 \leq t'$, it suffices to establish that $M \vDash Max^+(t',v,x)$.

Now, we have that $M \vDash Set(x)$ & $v \varepsilon x$ & $Tally_b(t')$.

Assume $M \vDash Fr(x,t_4,aua,t_5)$ & $u <_x v$.

Suppose that

(3bii1) $M \vDash Bound(x,t_4,aua,t_5) \lor Free^-(x,t_4,aua,t_5)$.

$\Rightarrow$ from principal hypothesis, $M \vDash \exists t_9 Fr(x^*,t_4,aua,t_9) \ \& \ u <_{x^*} v$,

$\Rightarrow$ from $M \vDash Firstf(x^*,t^*b,av_0a,t'')$, by (5.20), $M \vDash t^*b \leq t_4$.

We have that $M \vDash Max^+T_b(t^*b,t_0w)$, as in (2).

$\Rightarrow M \vDash Max^+T_b(t_4,t_0w)$,

$\Rightarrow M \vDash Fr(x^*,t_4,aua,t_9) \ \& \ z = t_0wx^* \ \& \ aBw \ \& \ aEw \ \& \ Max^+T_b(t_4,t_0w)$,

$\Rightarrow M \vDash Fr(z,t_4,aua,t_9)$,

$\Rightarrow$ $M \vDash Env(t,z)$ & $z=t_0wx^*$ & $aBw$ & $aEw$ & $Env(t^+,x')$ & $x'=t_0wt^+$ &

& $Env(t,x^*)$ & $Firstf(x^*,t^*b,av_0a,t'')$ & $t^+ < t^*b$,

$\Rightarrow$ since $M \vDash u \, \varepsilon \, x^*$ & $v \, \varepsilon \, x^*$, by (9.18), $M \vDash u <_{x^*} v \leftrightarrow u <_z v$,

$\Rightarrow$ from $M \vDash u <_{x^*} v$, $M \vDash u <_z v$,

$\Rightarrow$ from $M \vDash Max^+(t',v,z)$ & $Fr(z,t_4,aua,t_9)$, $M \vDash t_4 < t'$.

But then we have shown that $M \vDash Fr(x,t_4,aua,t_5)$ & $u <_x v \to t_4 < t'$,

that is, that $M \vDash Max^+(t',v,x)$.

This completes the proof of (3b) under hypothesis (3bii1).

(3bii2) $M \vDash Free^+(x,t_4,aua,t_5)$.

Then exactly the same argument as in (3bii1) with $t_4$ replaced throughout by $t_4b$ shows that $M \vDash t_4b < t'$ from hypothesis $M \vDash Fr(x,t_4,aua,t_5)$ & $u <_x v$.

Hence again $M \vDash Max^+(t',v,x)$.

This completes the proof of (3b) under the hypothesis (3bii).

Thus $M \vDash MMax^+T_b(t_1,v,z)$.

With (1), (2) and (3) we have thus shown that

$$M \vDash Fr(z,t_1,ava,t_2) \to MMax^+T_b(t_1,v,z).$$

This completes the proof that $M \vDash Special(z)$, and the proof of (10.28).

(10.29) For any string concept $I \subseteq I_0$ there are string concepts $I^* \subseteq I_0$ and $J \subseteq I^*$ such that

$QT^+ \vdash \forall t \in J \; \forall z \in I^*[Env(t,z) \; \& \; MinSet(z) \rightarrow \exists!z' \in I^* \exists t' \in J \; (Env(t',z') \; \& \; z \sim z' \; \&$

$\& \; MinSet(z') \; \& \; \forall v,w \; (v <_z w \leftrightarrow v <_{z'} w) \; \&$

$\& \; \forall w,t_1,t_2 \; (Free^+(z,t_1,awa,t_2) \rightarrow \exists t_3 Fr(z',t_1 b,awa,t_3)) \; \&$

$\& \; \forall w,t_1,t_2 \; (Bound(z,t_1,awa,t_2) \lor Free^-(z,t_1,awa,t_2) \rightarrow \exists t_3 Fr(z',t_1,awa,t_3)))]$.

Roughly, any set of strings with a minimal set code z has a minimal set code z' ordered in the same way as z in which every Free⁺ framed element has the length of its initial marker increased by 1 while all Bound and Free⁻ framed elements have their initial markers unchanged.

Let

$I^*(x) \equiv (I_{4.10} \; \& \; I_{5.7} \; \& \; I_{9.12} \; \& \; I_{9.13} \; \& \; I_{9.24} \; \& \; I_{10.8} \; \& \; I_{10.14})_{SUB}$.

Then we set

$J(x) \equiv I^*(x) \; \& \; \forall t \leq x \; \forall z \in I^*[Env(t,z) \; \& \; MinSet(z) \rightarrow$

$\rightarrow \exists!z' \in I^* \exists t'(Env(t',z') \; \& \; t' \leq tb \; \& \; z \sim z' \; \& \; MinSet(z') \; \& \; \forall v,w \; (v <_z w \leftrightarrow v <_{z'} w) \; \&$

$\& \; \forall w,t_1,t_2 \; (Free^+(z,t_1,awa,t_2) \rightarrow \exists t_3 Fr(z',t_1 b,awa,t_3)) \; \&$

$\& \; \forall w,t_1,t_2 \; (Bound(z,t_1,awa,t_2) \lor Free^-(z,t_1,awa,t_2) \rightarrow \exists t_3 Fr(z',t_1,awa,t_3)))]$.

We claim that J(x) is a string concept.

If x=a, we have that $M \vDash \neg Tally_b(x)$. Since $M \vDash Env(x,z) \rightarrow Tally_b(x)$, we have that $M \vDash \neg Env(x,z)$, hence $M \vDash J(x)$ holds trivially.

Suppose x=b.

To show that $M \vDash J(b)$, assume that $M \vDash Env(t,z)$ & $MinSet(z)$ & $t \leq b$ where $M \vDash I^*(z)$.

$\Rightarrow M \vDash Tally_b(t)$ & $t \leq b$,

$\Rightarrow$ by (1.3), $M \vDash t=b$,

$\Rightarrow$ from $M \vDash Env(t,z)$ & $t=b$ by (5.12), $M \vDash \exists u(z=bauab$ & $Pref(aua,b))$,

$\Rightarrow M \vDash Firstf(z,b,aua,b)$ & $Lastf(z,b,aua,b)$.

$\Rightarrow$ since we may assume that $I^*$ is downward closed with respect to $\subseteq_p$ and $M \vDash u \subseteq_p z$, $M \vDash I^*(u)$.

Let z'=bbauabb.

$\Rightarrow$ since we may assume that $I^*$ is closed with respect to $*$, $M \vDash I^*(z')$,

$\Rightarrow$ from $M \vDash Pref(aua,b)$, $M \vDash Pref(aua,bb)$,

$\Rightarrow$ since $M \vDash Tally_b(bb)$,

$\qquad M \vDash Firstf(z',bb,aua,bb)$ and also $M \vDash Lastf(z',bb,aua,bb)$,

$\Rightarrow$ by (5.22), $M \vDash Env(bb,z')$ & $\forall w(w \: \varepsilon \: z' \leftrightarrow w=u)$.

But, also by (5.22), $M \vDash \forall w(w \: \varepsilon \: z \leftrightarrow w=u)$. Hence $M \vDash z \sim z'$.

$\Rightarrow$ by (10.5), $M \vDash MinSet(z')$.

By (9.4), $M \vDash \forall v,w \: \neg(v<_z w)$ & $\forall v,w \: \neg(v<_{z'} w)$.

Hence $M \vDash \forall v,w \: (v<_z w \leftrightarrow v<_{z'} w)$.

Assume now that $M \vDash Fr(z,t_1,awa,t_2)$.

$\Rightarrow\ M \vDash \text{Tally}_b(t_1)\ \&\ \text{Tally}_b(t_2),$

$\Rightarrow$ from $M \vDash \text{Env}(b,z),\ M \vDash \text{MaxT}_b(b,z),$

$\Rightarrow\ M \vDash t_1 \leq b\ \&\ t_2 \leq b,$

$\Rightarrow$ by (1.3), $M \vDash t_1 = b = t_2,$

$\Rightarrow$ from $M \vDash \text{Env}(b,z),\ M \vDash w = u,$

$\Rightarrow\ M \vDash \text{Firstf}(z, t_1, awa, t_2).$

Therefore, $M \vDash \forall w, t_1, t_2\ (\text{Fr}(z, t_1, awa, t_2) \rightarrow \text{Firstf}(z, t_1, awa, t_2)).$

Assume that $M \vDash \text{Free}^+(z, t_1, awa, t_2).$

$\Rightarrow M \vDash \text{Fr}(z, t_1, awa, t_2),$

$\Rightarrow M \vDash \text{Firstf}(z, t_1, awa, t_2),$

$\Rightarrow$ from $M \vDash \text{Firstf}(z, b, aua, b)$ by (5.15), $M \vDash t_1 = b\ \&\ w = u.$

But we have $M \vDash \text{Fr}(z, t_1 b, awa, t_1 b).$

Hence $M \vDash \forall w, t_1, t_2\ (\text{Free}^+(z, t_1, awa, t_2) \rightarrow \exists t_3 \text{Fr}(z', t_1 b, awa, t_3)).$

Assume that $M \vDash \text{Bound}(z, t_1, awa, t_2)\ \vee\ \text{Free}^-(z, t_1, awa, t_2).$

$\Rightarrow M \vDash \text{Fr}(z, t_1, awa, t_2)\ \&\ \neg \text{Firstf}(z, t_1, awa, t_2).$

But then $M \vDash \forall w, t_1, t_2\ \neg \text{Fr}(z, t_1, awa, t_2).$

$\Rightarrow M \vDash \forall w\ \neg(w\ \varepsilon\ z),$ which contradicts $M \vDash u\ \varepsilon\ z.$

Hence $M \vDash \forall w, t_1, t_2\ \neg(\text{Bound}(z, t_1, awa, t_2)\ \vee\ \text{Free}^-(z, t_1, awa, t_2)),$ and so, trivially,

$M \vDash \forall w, t_1, t_2\ (\text{Bound}(z, t_1, awa, t_2)\ \vee\ \text{Free}^-(z, t_1, awa, t_2) \rightarrow \exists t_3 \text{Fr}(z', t_1, awa, t_3)).$

To establish uniqueness of z', assume that

  $M \vDash Env(t'',z'')$ & $t'' \leq tb$ & $z \sim z''$ &

  & $\forall w, t_1, t_2 \ (Free^+(z,t_1,awa,t_2) \rightarrow \exists t_3 Fr(z'',t_1 b,awa,t_3))$

where $M \vDash I^*(z'')$.

From $M \vDash Firstf(z,b,aua,b)$, $M \vDash Free^+(z,b,aua,b)$.

$\Rightarrow M \vDash \exists t_3 Fr(z'',bb,awa,t_3)$.

From $M \vDash z \sim z''$, we have $M \vDash \forall w(w \ \varepsilon \ z'' \leftrightarrow w = u)$.

$\Rightarrow$ by (5.22), $M \vDash z'' = t''aut''$ & $Firstf(z'',t'',aua,t'')$,

$\Rightarrow$ from $M \vDash Env(t'',z'')$ & $Fr(z'',bb,aua,t_3)$ & $Fr(z'',t'',aua,t'')$, $M \vDash t'' = bb$,

$\Rightarrow M \vDash z'' = t''aut'' = bbauabb = z'$, as required.

This completes the argument that $M \vDash J(b)$.

Assume now that $M \vDash J(t)$.

$\Rightarrow M \vDash I^*(t)$,

$\Rightarrow$ since $I^*$ is a string concept, $M \vDash I^*(tb)$.

Suppose that $M \vDash Env(x,z)$ & $MinSet(z)$ where $M \vDash x \leq tb$ & $I^*(z)$.

$\Rightarrow M \vDash Tally_b(x)$,

$\Rightarrow$ by (1.13), $M \vDash x \leq t \lor x = tb$.

If $M \vDash x \leq t$, the desired claim follows from $M \vDash J(t)$.

So we may assume that $M \vDash x = tb$.

$\Rightarrow M \vDash \exists u Lastf(z,tb,aua,tb)$,

$\Rightarrow M \vDash Pref(aua,tb)$ & $(z = tbauatb \lor \exists w_1(z = w_1 atbauatb$ & $Max^+ T_b(tb,w_1)))$.

(1) $M \vDash z = tbauatb$.

Let z'=tbbauatbb.

Then $M \vDash J(tb)$ by an argument exactly analogous to that for $M \vDash J(b)$.

(2) $M \vDash \exists w_1(z=w_1 atbauatb\ \&\ Max^+T_b(tb,w_1))$.

$\Rightarrow$ from $M \vDash I^*(z)$, just as in (1), $M \vDash I^*(u)$,

$\Rightarrow$ by the SUBTRACTION LEMMA,

$M \vDash \exists z^- \in I^* \exists v',t'[Env(t',z^-)\ \&\ Lastf(z^-,t',av'a,t')\ \&\ ((z^-=t'av'at'\ \&\ w_1=t'av')\ v$

$v\ \exists w',w'',t''(z=w'at'av'at''aw''\ \&\ t''aw''=tbauatb\ \&$

$\&\ z^-=w'at'av'at'\ \&\ t'<t''=tb))\ \&\ \forall w(w\ \varepsilon\ z^- \leftrightarrow w\ \varepsilon\ z\ \&\ w \neq u)]$.

(2a) $M \vDash z^-=t'av'at'\ \&\ w_1=t'av'$.

$\Rightarrow$ from $M \vDash Lastf(z^-,t',av'a,t')$, $M \vDash Tally_b(t')\ \&\ t' \subseteq_p w_1\ \&\ Pref(av'a,t')$,

$\Rightarrow$ from $M \vDash Max^+T_b(tb,w_1)$, $M \vDash t'<tb$,

$\Rightarrow$ by (1.13), $M \vDash t' \leq t$,

$\Rightarrow$ since we may assume that $I^*$ is downward closed under $\leq$, from

hypothesis $M \vDash J(t)$, $M \vDash J(t')$,

$\Rightarrow M \vDash \exists!z^+ \in I^* \exists t^+(Env(t^+,z^+)\ \&\ t^+ \leq t'b\ \&\ z^- \sim z^+\ \&\ MinSet(z^+)\ \&$

$\&\ \forall v,w\ (v<_{z^-}w \leftrightarrow v<_{z^+}w)\ \&$

$\&\ \forall w,t_1,t_2\ (Free^+(z^-,t_1,awa,t_2) \to \exists t_3 Fr(z^+,t_1 b,awa,t_3))\ \&$

$\&\ \forall w,t_1,t_2\ (Bound(z^-,t_1,awa,t_2)\ v\ Free^-(z^-,t_1,awa,t_2) \to \exists t_3 Fr(z^+,t_1,awa,t_3)))$,

$\Rightarrow$ from $M \vDash Pref(av'a,t')$, $M \vDash Firstf(z^-,t',av'a,t')$,

$\Rightarrow$ by (5.22), $M \vDash \forall w(w\ \varepsilon\ z^- \leftrightarrow w=v')$,

$\Rightarrow$ from $M \vDash z^- \sim z^+$, $M \vDash \forall w(w\ \varepsilon\ z^+ \leftrightarrow w=v')$,

$\Rightarrow M \vDash \exists t_3 Fr(z^+,t'b,av'a,t_3)$,

$\Rightarrow$ from $M \vDash Env(t^+,z^+)$, by (5.22),

$M \vDash z^+ = t^+av'at^+$ & $Firstf(z^+,t^+,av'a,t^+)$ & $Lastf(z^+,t^+,av'a,t^+)$,

$\Rightarrow M \vDash t'b = t^+$,

$\Rightarrow M \vDash z^+ = t'bav'at'b$,

$\Rightarrow$ from $M \vDash t' < tb$ & $Tally_b(tb)$, by (4.7), $M \vDash t'b \leq tb$.

(2ai) $M \vDash t'b = tb$.

Let $z' = z^+bauatbb$.

$\Rightarrow$ from $M \vDash I^*(z^+b)$ & $I^*(u)$, since we may assume that $I^*$ is closed with respect to $*$, $M \vDash I^*(z')$.

(2aii) $M \vDash t'b < tb$.

$\Rightarrow M \vDash \exists t^*(Tally_b(t^*)$ & $t'bt^* = tb)$.

$\Rightarrow$ by (4.10), $M \vDash (t'b)t^* = t^*(t'b)$,

$\Rightarrow M \vDash t^*(t'b) = tb$,

$\Rightarrow M \vDash t^* < tb$,

$\Rightarrow$ from $M \vDash I^*(tb)$, $M \vDash I^*(t^*)$.

Let $z' = z^+t^*auatb$.

$\Rightarrow$ from $M \vDash I^*(z^+)$ & $I^*(t^*)$ & $I^*(u)$ & $I^*(tb)$, $M \vDash I^*(z')$.

Note that from $M \vDash z = w_1atbauatb$ & $w_1 = t'av'$ we have $M \vDash z = t'av'atbauatb$.

$\Rightarrow$ from $M \vDash Pref(av'a,t')$ & $Tally_b(tb)$ & $t' < tb$ & $(t'av'atba)Bz$,

$M \vDash Firstf(z,t',av'a,tb)$.

Hence in either case from $M \vDash Lastf(z,tb,aua,tb)$, by repeating the argument of the uniqueness proof in (10.14), which is not affected by the additional conditions in I*, we have that

  $M \vDash \exists! z' \in I^*$ (Env(tb,z) & MinSet(z') & z∼z' & $\forall v,w$ (v$<_z$w ↔ v$<_{z'}$w) &

    & $\forall w, t_1, t_2$ (Free$^+$(z,$t_1$,awa,$t_2$) → $\exists t_3$Fr(z',$t_1$b,awa,$t_3$)) &

& $\forall w, t_1, t_2$ (Bound(z,$t_1$,awa,$t_2$) v Free$^-$(z,$t_1$,awa,$t_2$) → $\exists t_3$Fr(z',$t_1$,awa,$t_3$))).

 (2b)  $M \vDash \exists w', w'', t''$ (z=w'at'av'at''aw'' & t''aw''=tbauatb &

    & z$^-$=w'at'av'at' & t'<t''=tb).

$\Longrightarrow$ from $M \vDash$ t'<tb, by (1.13),  $M \vDash$ t'≤t,

$\Longrightarrow$ from hypothesis $M \vDash$ J(t) & t'≤t, just as in (2a),  $M \vDash$ J(t').

We first claim that   $M \vDash \exists t_3$Fr(z,t',ava,$t_3$).

From  $M \vDash$ t'≤t & Tally$_b$(t'),  we have as in (2a) that

    $M \vDash \exists t^*$(Tally$_b$(t*) & t't*=t'').

$\Longrightarrow$ from $M \vDash$ t''=tb & t''aw''=tbauatb & z=w'at'av'at''aw'' & z$^-$=w'at'av'at',

 $M \vDash$ z=w'at'av'atbauatb=w'at'av'at''auatb=(w'at'av'at')t*auatb=z$^-$t*auatb,

$\Longrightarrow$ from $M \vDash$ Env(t',z$^-$) & Fr(z$^-$,t',av'a,t') & Tally$_b$(t*) & Tally$_b$(tb),

    $M \vDash \exists t_3$Fr(z,t',ava,$t_3$).

So we have

$M \vDash$ Env(tb,z) & MinSet(z) & Fr(z,t',ava,$t_3$) & Env(t',z$^-$) &

    & Lastf(z$^-$,t',av'a,t') & z$^-$Bz,

$\Longrightarrow$ by (10.8),  $M \vDash$ MinSet(z$^-$).

Then from $M \vDash$ Env(t',z$^-$) & MinSet(z$^-$) we obtain from $M \vDash$ J(t') a unique z$^+$ in I* as in (2a).

$\Longrightarrow$ from $M \vDash$ MinSet(z$^+$),  $M \vDash$ Set(z$^+$),

$\Rightarrow$ from $M \vDash z^-\sim z^+$ & $v' \varepsilon z^-$, $M \vDash v' \varepsilon z^+$,

$\Rightarrow$ by (5.18), $M \vDash z^+ \neq aa$,

$\Rightarrow$ $M \vDash \exists t^+ \, Env(t^+,z^+)$.

We have that $M \vDash z^- = w'at'av'at'$ & $Lastf(z^-,t',av'a,t')$.

Suppose, for a reductio, that $M \vDash Firstf(z^-,t',av'a,t')$.

$\Rightarrow$ by (5.7), $M \vDash z^- = t'av'at'$,

$\Rightarrow$ $M \vDash t'av'at' = z^- = w'at'av'at'$,

$\Rightarrow$ $M \vDash z^- E z^-$, contradicting $M \vDash z^- \in I$ by (3.4).

Therefore, $M \vDash \neg Firstf(z^-,t',av'a,t')$.

We claim that $M \vDash Lastf(z^+,t^+,av'a,t^+)$.

From $M \vDash Lastf(z^-,t',av'a,t')$ and the choice of $z^+$ we have, by (9.12), that

$$M \vDash \exists t_3, t_4 \, Lastf(z^+,t_3,av'a,t_4).$$

$\Rightarrow$ from $M \vDash Env(t^+,z^+)$, $M \vDash t_3 = t^+ = t_4$,

$\Rightarrow$ $M \vDash Lastf(z^+,t^+,av'a,t^+)$,

$\Rightarrow$ from $M \vDash Lastf(z^-,t',av'a,t')$ & $\neg Firstf(z^-,t',av'a,t')$, by (9.22),

$$M \vDash Free(z^-,t',av'a,t') \lor Bound(z^-,t',av'a,t').$$

(2bi) $M \vDash Free^+(z^-,t',av'a,t')$.

$\Rightarrow$ by choice of $z^+$, $M \vDash \exists t_3 \, Fr(z^+,t'b,av'a,t_3)$,

$\Rightarrow$ from $M \vDash Env(t^+,z^+)$ & $Lastf(z^+,t^+,av'a,t^+)$, $M \vDash t^+ = t'b$,

$\Rightarrow$ $M \vDash Lastf(z^+,t'b,av'a,t'b)$,

$\Rightarrow$ from $M \vDash t' < t'' = tb$ & $Tally_b(t'')$, by (4.7), $M \vDash t'b \leq t'' = tb$.

(2bia) $M \vDash t'b = t'' = tb$.

Let $z'=z^+bauatbb$ as in (2ai).

   (2bib)   $M \vDash t'b<t''=tb$.

$\Rightarrow$ from $M \vDash Tally_b(t'')$, $M \vDash \exists t^*(Tally_b(t^*)\ \&\ t'bt^*=t'')$.

Let $z'=z^+t^*auatb$ as in (2aii).

   (2bii)   $M \vDash Free^-(z^-,t',av'a,t') \vee Bound(z^-,t',av'a,t')$.

$\Rightarrow$ by choice of $z^+$, $M \vDash \exists t_3 Fr(z^+,t',av'a,t_3)$,

$\Rightarrow$ from $M \vDash Env(t^+,z^+)\ \&\ Lastf(z^+,t^+,av'a,t^+)$, $M \vDash t^+=t'$,

$\Rightarrow$ $M \vDash Lastf(z^+,t',av'a,t')$.

Again, from $M \vDash t'b \leq t''=tb$, we distinguish subcases:

   (2biia)   $M \vDash t'b=t''=tb$.

$\Rightarrow$ $M \vDash tbb=t'bb=t^+bb$.

Let $z'=z^+bbauatbb$.

   (2biib)   $M \vDash t'b<t''=tb$.

$\Rightarrow$ $M \vDash t'bt^*=t''=tb$,

$\Rightarrow$ $M \vDash t^+<tb$.

Let $z'=z^+bt^*auatb$.

$\Rightarrow$ by choice of $z^+$, $M \vDash \exists t_3 Firstf(z^+,t_1,av_0a,t_3)$,

$\Rightarrow$ $M \vDash (t_1a)Bz^+\ \&\ Tally_b(t_1)$,

$\Rightarrow$ from $M \vDash Env(t^+,z^+)$, by (5.11),

       $M \vDash \exists t_4, z_1(Tally_b(t_4)\ \&\ z^+=t_4z_1t^+\ \&\ aBz_1\ \&\ aEz_1)$,

$\Rightarrow$ $M \vDash (t_4a)Bz^+$,

$\Rightarrow$ by (4.23[b]), $M \vDash t_1=t_4$.

Let $u_1 = tbauatb$ and $u_2 = tbbauatbb$.

$\Rightarrow$ from $M \vDash I^*(t)$ & $I^*(u)$, $M \vDash I^*(u_1)$ & $I^*(u_2)$,

$\Rightarrow$ from $M \vDash Pref(aua, tb)$ & $Tally_b(tbb)$, $M \vDash Pref(aua, tbb)$.

Now, in (2bia) and (2biia) we have that

$$M \vDash z' = t_1 z_1 tbbauatbb = t_1 z_1 u_2,$$

and in (2bib) and (2biib) we have

$$M \vDash z' = t_1 z_1 t'bt^*auatb = t_1 z_1 tbauatb = t_1 z_1 u_1.$$

$\Rightarrow$ by (10.5),

$$M \vDash MinSet(u_1) \text{ \& } Env(tb, u_1) \text{ and } M \vDash MinSet(u_1) \text{ \& } Env(tbb, u_2),$$

$\Rightarrow$ from $M \vDash Pref(aua, tb)$ & $Pref(aua, tbb)$,

$$M \vDash Firstf(u_1, tb, aua, tb) \text{ \& } Firstf(u_2, tbb, aua, tbb),$$

$\Rightarrow$ by (5.21), $M \vDash \forall w(w \; \varepsilon \; u_1 \leftrightarrow w = u \leftrightarrow w \; \varepsilon \; u_2)$,

$\Rightarrow$ from $M \vDash \neg(u \; \varepsilon \; z^-)$ & $z^- \sim z^+$, $M \vDash \neg(u \; \varepsilon \; z^+)$,

$\Rightarrow$ $M \vDash \neg \exists w(w \; \varepsilon \; z^+ \text{ \& } w \; \varepsilon \; u_1)$ and $M \vDash \neg \exists w(w \; \varepsilon \; z^+ \text{ \& } w \; \varepsilon \; u_2)$.

So we have, in (2bia) and (2biia), that

$M \vDash Env(t^+, z^+)$ & $z^+ = t_1 z_1 t^+$ & $aBz_1$ & $aEz_1$ & $z' = t_1 z_1 u_2$ & $Env(tbb, u_2)$ &

& $t^+ < tbb$ & $Firstf(u_2, tbb, aua, tbb)$ & $\neg \exists w(w \; \varepsilon \; z^+ \text{ \& } w \; \varepsilon \; u_2)$ & $MinSet(z^+)$ &

& $MinSet(u_2)$,

and in (2bib) and (2biib) that

$M \vDash Env(t^+, z^+)$ & $z^+ = t_1 z_1 t^+$ & $aBz_1$ & $aEz_1$ & $z' = t_1 z_1 u_1$ & $Env(tb, u_1)$ &

& $t^+ < tb$ & $Firstf(u_1, tb, aua, tb)$ & $\neg \exists w(w \; \varepsilon \; z^+ \text{ \& } w \; \varepsilon \; u_1)$ & $MinSet(z^+)$ &

& $MinSet(u_1)$.

$\Rightarrow$ by (10.6),

$M \vDash \text{MinSet}(z') \ \& \ \text{Env}(tbb,z')$ and $M \vDash \text{MinSet}(z') \ \& \ \text{Env}(tb,z')$,

respectively.

Since $M \vDash \neg(u \ \varepsilon \ z^+)$, by the proof of the SET ADJUNCTION LEMMA, (7.1)(2),

$$M \vDash \forall w(w \ \varepsilon \ z' \leftrightarrow w \ \varepsilon \ z^+ \vee w=u).$$

$\Rightarrow$ from $M \vDash z^- \sim z^+$, $M \vDash \forall w(w \ \varepsilon \ z' \leftrightarrow (w \ \varepsilon \ z^- \vee w=u) \leftrightarrow w \ \varepsilon \ z)$.

Therefore, $M \vDash z' \sim z$.

Next, we show that $M \vDash \forall w(w \ \varepsilon \ z \rightarrow \forall v \ (v <_z w \leftrightarrow v <_{z'} w))$.

We have $M \vDash \text{Env}(t', z^-) \ \& \ z = z^- t^* auatb \ \& \ \text{Tally}_b(t^*) \ \& \ \text{Tally}_b(tb)$.

$\Rightarrow$ by (9.11), $M \vDash \forall v,w(v \ \varepsilon \ z^- \ \& \ w \ \varepsilon \ z^- \rightarrow (v <_{z^-} w \leftrightarrow v <_z w))$.

Assume $M \vDash w \ \varepsilon \ z$. Then $M \vDash w \ \varepsilon \ z^- \vee w=u$.

Assume $M \vDash w \ \varepsilon \ z^-$.

We need only consider the case $M \vDash \neg(v \ \varepsilon \ z^-)$.

Then $M \vDash \neg(v <_{z^-} w)$. If $M \vDash \neg(v \ \varepsilon \ z)$, then also $M \vDash \neg(v <_z w)$, and so

$$M \vDash v <_{z^-} w \leftrightarrow v <_z w.$$

So we may assume $M \vDash v \ \varepsilon \ z$.

$\Rightarrow$ from $M \vDash \neg(v \ \varepsilon \ z^-)$, $M \vDash v=u$,

$\Rightarrow$ from $M \vDash w \ \varepsilon \ z^- \ \& \ \neg(v \ \varepsilon \ z^-)$, $M \vDash w \neq v$,

$\Rightarrow$ from $M \vDash \text{Lastf}(z,tb,aua,tb)$, by (9.3), $M \vDash w <_z v$,

$\Rightarrow$ by (9.6), $M \vDash \neg(v <_z w)$.

Hence, again, $M \vDash v <_{z^-} w \leftrightarrow v <_z w$.

Therefore, we have shown that $M \vDash \forall w(w \ \varepsilon \ z^- \rightarrow (v <_{z^-} w \leftrightarrow v <_z w))$.

Likewise, we have

$$M \vDash Env(t^+,z^+) \ \& \ z'=z^+bauatbb \ \& \ Tally_b(b) \ \& \ Tally_b(tbb),$$

$$M \vDash Env(t^+,z^+) \ \& \ z'=z^+t^*auatb \ \& \ Tally_b(t^*) \ \& \ Tally_b(tb),$$

$$M \vDash Env(t^+,z^+) \ \& \ z'=z^+bbauatbb \ \& \ Tally_b(bb) \ \& \ Tally_b(tbb),$$

$$M \vDash Env(t^+,z^+) \ \& \ z'=z^+bt^*auatb \ \& \ Tally_b(bt^*) \ \& \ Tally_b(tb),$$

respectively, whence we obtain, by (9.11),

$$M \vDash \forall w(w \ \varepsilon \ z^+ \to \forall v(v<_{z^+}w \leftrightarrow v<_{z'}w)).$$

$\Rightarrow$ from the choice of $z^+$ and $M \vDash z^- \sim z^+$,

$$M \vDash \forall w(w \ \varepsilon \ z^- \to \forall v(v<_z w \leftrightarrow v<_{z'}w)).$$

If $M \vDash w=u$, we first observe that, from the proof of (10.6), we have, in (2bia) and (2biia) from $M \vDash Env(tbb,z')$, that $M \vDash Lastf(z',tbb,aua,tbb)$, and in (2bib) and (2biib) from $M \vDash Env(tb,z')$, that $M \vDash Lastf(z',tb,aua,tb)$. Then, in (2bia) and (2biia), from

$$M \vDash Lastf(z,tb,aua,tb) \ \& \ Lastf(z',tbb,aua,tbb) \ \& \ z \sim z',$$

for $M \vDash v \ \varepsilon \ z$, by (9.3), that $M \vDash v \leq_z u \leftrightarrow v \leq_{z'} u$. For $M \vDash \neg(v \ \varepsilon \ z)$, from $M \vDash z \sim z'$, we have $M \vDash \neg(v \leq_z u) \ \& \ \neg(v \leq_{z'} u)$, so again $M \vDash v \leq_z u \leftrightarrow v \leq_{z'} u$. Therefore, $M \vDash \forall w(w \ \varepsilon \ z \to \forall v \ (v \leq_z w \leftrightarrow v \leq_{z'} w))$, as required.

We reason the same way in (2bib) and (2biib).

Next, we show that $M \vDash \forall w,t_1,t_2(Free^+(z,t_1,awa,t_2) \to \exists t_3 Fr(z',t_1b,awa,t_3))$.

Assume $M \vDash Free^+(z,t_1,awa,t_2)$.

$\Rightarrow M \vDash w \ \varepsilon \ z$,

$\Rightarrow M \vDash w \ \varepsilon \ z^- \lor w=u.$

Assume that $M \vDash w=u$.

$\implies$ from $M \vDash Fr(z,tb,aua,tb)$ & $Env(tb,z)$, $M \vDash t_1=tb$,

Now, we have

$M \vDash Env(t', z^-)$ & $z=z^-t^*auatb$ & $Tally_b(t^*)$ & $Env(t',z^-)$ & $Env(tb,z)$ &

& $Lastf(z^-,t',av'a,t')$ & $Lastf(z,tb,aua,tb)$.

$\implies$ by (9.13), $M \vDash v'<_z u$ & $\neg \exists w'(v'<_z w'$ & $w'<_z u)$,

$\implies$ from $M \vDash z= z^-t^*auatb$ & $Tally_b(tb)$ & $Fr(z^-,t',av'a,t')$, by (5.6),

$M \vDash \exists t_3 Fr(z,t',av'a,t_3)$,

$\implies$ from $M \vDash Set(z)$ & $v'<_z u$, by (9.1), $M \vDash \neg Firstf(z,t_1,aua,t_2)$,

$\implies$ from hypothesis $M \vDash Free^+(z,t_1,awa,t_2)$, $M \vDash t_1=t'b$,

$\implies$ from $M \vDash t_1=tb$, $M \vDash tb=t_1=t'b$,

$\implies$ by (2bia) and (2biia),

$M \vDash (z'= z^+bauatbb$ v $z'=z^+bbauatbb)$ & $Lastf(z',tbb,aua,tbb)$,

$\implies M \vDash \exists t_3 Fr(z',t_1b,aua,t_3)$, as required.

Assume that $M \vDash w \varepsilon z^-$.

$\implies$ from $M \vDash Env(tb, z)$ & $z=z^-t^*auatb$ & $Tally_b(t^*)$ & $Env(t',z^-)$ &

& $Lastf(z,tb,aua,t')$ & $Free^+(z,t_1,awa,t_2)$, by (9.24),

$M \vDash \exists t_4 Fr(z^-,t_1,awa,t_4)$,

$\implies$ by choice of $z^+$, $M \vDash \exists t_5 Fr(z^+,t_1b,awa,t_5)$,

$\implies$ by (5.6), $M \vDash \exists t_3 Fr(z',t_1b,awa,t_3)$, as required.

Finally, we show that

$M \vDash \forall w,t_1,t_2 (Bound(z,t_1,awa,t_2)$ v $Free^-(z,t_1,awa,t_2) \rightarrow \exists t_3 Fr(z',t_1,awa,t_3))$.

Assume $M \models \text{Bound}(z,t_1,awa,t_2) \lor \text{Free}^-(z,t_1,awa,t_2)$.

In particular, assume that $M \models \text{Bound}(z,t_1,awa,t_2)$.

$\Rightarrow M \models w \; \varepsilon \; z$,

$\Rightarrow M \models w \; \varepsilon \; z^- \lor w=u$.

If $M \models w \; \varepsilon \; z^-$, then an argument completely analogous to the one just given applies, using (9.23) in place of (9.24).

If $M \models w=u$, then we proceed just as under the hypothesis

$M \models \text{Free}(z,t_1,awa,t_2)$ until we obtain $M \models \exists t_3 \text{Fr}(z,t',av'a,t_3)$.

We have that $M \models t_1=tb$.

From $M \models \text{Bound}(z,t_1,awa,t_2) \; \& \; w=u$, by the choice of v', $M \models t'b<t_1$.

$\Rightarrow M \models t'b<tb$.

Then, by (2bib) and (2biib), we have that

$$M \models (z'= z^+t^*auatb \lor z'=z^+bt^*bauatb) \; \& \; \text{Lastf}(z',tb,aua,tb).$$

Hence $M \models \text{Fr}(z',t_1,aua,t_1)$, as required.

Assume, on the other hand, that $M \models \text{Free}^-(z,t_1,awa,t_2)$.

$\Rightarrow M \models \text{Free}(z,t_1,awa,t_2) \; \& \; \neg \text{Free}^+(z,t_1,awa,t_2)$,

$\Rightarrow$ by (9.23), $M \models \exists t_3 \text{Free}(z^-,t_1,awa,t_3)$.

Suppose, for a reductio, that $M \models \text{Free}^+(z^-,t_1,awa,t_3)$.

$\Rightarrow$ by (9.24), $M \models \exists t_4 \text{Free}^+(z,t_1,awa,t_4)$,

$\Rightarrow M \models \text{Fr}(z,t_1,awa,t_4)$.

From $M \models \text{Free}(z,t_1,awa,t_2)$, $M \models \text{Fr}(z,t_1,awa,t_2)$,

$\Rightarrow$ by (5.42), $M \models t_2=t_4$.

But then  $M \vDash \text{Free}^+(z,t_1,awa,t_4)$ & $\neg\text{Free}^+(z,t_1,awa,t_2)$  is a contradiction.

Therefore,  $M \vDash \neg\text{Free}^+(z^-,t_1,awa,t_3)$.

$\implies$ since $M \vDash \text{Free}(z^-,t_1,awa,t_3)$,  $M \vDash \text{Free}^-(z^-,t_1,awa,t_3)$,

$\implies$ by choice of $z^+$,  $M \vDash \exists t_4 \text{Fr}(z^+,t_1,awa,t_4)$,

$\implies$ by (5.6),  $M \vDash \exists t_5 \text{Fr}(z',t_1,awa,t_5)$,  as required.

It now remains to establish the uniqueness of $z'$.

Assume that

$\quad M \vDash \text{MinSet}(z'')$ & $z \sim z''$ & $\forall v,w\ (v<_z w \leftrightarrow v<_{z''} w)$ &

$\qquad$ & $\forall w,t_1,t_2\ (\text{Free}^+(z,t_1,awa,t_2) \rightarrow \exists t_3 \text{Fr}(z'',t_1 b,awa,t_3))$ &

$\quad$ & $\forall w,t_1,t_2\ (\text{Bound}(z,t_1,awa,t_2) \vee \text{Free}^-(z,t_1,awa,t_2) \rightarrow \exists t_3 \text{Fr}(z'',t_1,awa,t_3))$

where  $M \vDash z'' \in I^*$.

$\implies$ from  $M \vDash \text{MinSet}(z'')$,  $M \vDash \text{Set}(z'')$,

$\implies$ from  $M \vDash \text{Env}(tb,z)$,  $M \vDash \exists t_1,t_2,v_0\ \text{Firstf}(z,t_1,av_0 a,t_2)$.

Just as for $z^-$ earlier, we show that  $M \vDash t_1 < t_2$.

We have that $M \vDash \text{Lastf}(z,tb,aua,tb)$.

$\implies$ by (9.12),  $M \vDash \exists t''(\text{Lastf}(z'',t'',aua,t'')$ & $\text{Env}(t'',z''))$,

$\implies$ from  $M \vDash \exists t_1,t_2,v_0\ \text{Firstf}(z,t_1,av_0 a,t_2)$,  $M \vDash (t_1 a)Bz$ & $\text{Tally}_b(t_1)$,

$\implies$ from  $M \vDash z=w_1 atbauatb$, by (4.14[b]),  $M \vDash w_1=t_1 \vee t_1 B w_1$,

$\implies$ $M \vDash t_1 \subseteq_p w_1$,

$\implies$ from  $M \vDash \text{Max}^+ T_b(tb,w_1)$,  $M \vDash t_1 < tb$,

$\implies$ from  $M \vDash tb \in I \subseteq I_0$, $M \vDash t_1 \neq tb$,

$\implies$ by (5.15), $M \vDash \neg Firstf(z,tb,aua,tb)$,

$\implies$ from $M \vDash Lastf(z,tb,aua,tb)$, by (9.22),

$M \vDash Free(z,tb,aua,tb) \lor Bound(z,tb,aua,tb)$,

$\implies$ from hypothesis about $z''$,

$M \vDash (Free^+(z,tb,aua,tb) \& \exists t_3 Fr(z'',tbb,awa,t_3)) \lor$

$\lor ((Bound(z,tb,aua,tb) \lor Free^-(z,tb,aua,tb)) \& \exists t_4 Fr(z'',tb,aua,t_4))$,

$\implies$ from $M \vDash Env(t'',z'') \& Fr(z'',t'',aua,t'')$,

$M \vDash (Free^+(z,tb,aua,tb) \& t''=tbb) \lor$

$\lor ((Bound(z,tb,aua,tb) \lor Free^-(z,tb,aua,tb)) \& t''=tb)$.

(2bA) Suppose that, as in (2bia),

$M \vDash Free^+(z,tb,aua,tb) \& t''=tbb$.

$\implies M \vDash Lastf(z'',tbb,aua,tbb)$,

$\implies M \vDash z''=tbbauatbb \lor \exists w_2 (z''=w_2atbbauatbb \& Max^+T_b(tbb,w_2))$.

Suppose $M \vDash z''=tbbauatbb$.

$\implies M \vDash z''=u_2$,

$\implies$ by earlier argument, $M \vDash \forall w(w \, \varepsilon \, z'' \leftrightarrow w=u)$.

But we have $M \vDash Fr(z,t',av'a,t'') \& Fr(z,tb,aua,tb) \& t'<tb$.

$\implies$ from $M \vDash Env(tb,z) \& t' \in I \subseteq I_0$, $M \vDash v' \neq u \, \& \, u \, \varepsilon \, z \, \& \, v' \, \varepsilon \, z$,

$\implies$ from $M \vDash z \sim z''$, $M \vDash u \, \varepsilon \, z'' \, \& \, v' \, \varepsilon \, z''$,

$\implies M \vDash u=v'$, a contradiction.

Therefore, $M \vDash \neg(z''=tbbauatbb)$.

$\implies M \vDash \exists w_2 (z''=w_2atbbauatbb \& Max^+T_b(tbb,w_2))$,

$\Rightarrow$ from $M \vDash \text{Env}(tb,z)$ & $\text{MinSet}(z'')$, by the SUBTRACTION LEMMA,

$M \vDash \exists z^{**} \in I^* \exists v'', t^{**}[\text{Env}(t^{**},z^{**})$ & $\text{Lastf}(z^{**},t^{**},av''a,t^{**})$ &

& $((z^{**}=t^{**}av''at^{**}$ & $w_2=t^{**}av'')$ v

v $\exists w_3, w_4, t_3 (z''=w_3at^{**}av''at_3aw_4$ & $t_3aw_4=tbbauatbb$ & $z^{**}=w_3at^{**}av''at^{**}$ &

& $t^{**}<t_3=tbb))$ &

& $\forall w(w \; \varepsilon \; z^{**} \leftrightarrow w \; \varepsilon \; z''$ & $w \neq u)]$,

$\Rightarrow$ from $M \vDash z' \sim z \sim z''$, $M \vDash \forall w(w \; \varepsilon \; z^{**} \leftrightarrow w \; \varepsilon \; z'$ & $w \neq u \leftrightarrow w \; \varepsilon \; z^+)$.

So $M \vDash z^{**} \sim z^+$.

We claim that $M \vDash z^{**}=z^+$.

There are two scenarios:

(2bAi) $M \vDash z^{**}=t^{**}av''at^{**}$ & $w_2=t^{**}av''$.

$\Rightarrow$ from $M \vDash \text{Lastf}(z^{**},t^{**},av''a,t^{**})$, $M \vDash \text{Pref}(av''a,t^{**})$,

$\Rightarrow M \vDash \text{Firstf}(z^{**},t^{**},av''a,t^{**})$,

$\Rightarrow$ by (5.21), $M \vDash \forall w(w \; \varepsilon \; z^{**} \leftrightarrow w=v'')$,

$\Rightarrow$ from $M \vDash z^{**} \sim z^+$, $M \vDash \forall w(w \; \varepsilon \; z^+ \leftrightarrow w=v'')$,

$\Rightarrow$ from $M \vDash v' \; \varepsilon \; z$ & $z \sim z'$, $M \vDash v' \; \varepsilon \; z'$,

$\Rightarrow$ from $M \vDash v' \neq u$, $M \vDash v' \; \varepsilon \; z^+$,

$\Rightarrow M \vDash v''=v'$,

$\Rightarrow M \vDash z^{**}=t^{**}av'at^{**}$.

$\Rightarrow M \vDash w_2=t^{**}av''$ & $z''=w_2atbbauatbb$, $M \vDash z''=t^{**}av''atbbauatbb$,

$\Rightarrow M \vDash z''=t^{**}av'atbbauatbb$,

$\Rightarrow$ from $M \vDash \text{Firstf}(z,t_1,av_0a,t_2)$, by hypothesis about $z''$ and (9.12),

$$M \vDash \exists t_5 \text{ Firstf}(z'',t_1b,av_0a,t_5),$$

$\implies$ from $M \vDash t^{**}<tbb$ & $\text{Tally}_b(t^{**})$, $M \vDash \exists t^+(\text{Tally}_b(t^+)$ & $t^{**}t^+=tbb)$,

$\implies M \vDash z''=t^{**}av'atbbauatbb=(t^{**}av'at^{**})t^+auatbb=z^{**}t^+auatbb$,

$\implies M \vDash z^{**}Bz''$,

$\implies$ from $M \vDash \text{Env}(t^{**},z^{**})$ & $z^{**}Bz''$ & $\text{Firstf}(z'',t_1b,av_0a,t_5)$, by (5.4),

$$M \vDash \exists t_6 \text{ Firstf}(z^{**},t_1b,av_0a,t_6),$$

$\implies$ from $M \vDash \text{Firstf}(z^{**},t^{**},av'a,t^{**})$, by (5.15), $M \vDash t_1b=t^{**}$ & $v_0=v'$,

$\implies M \vDash z^{**}=t_1bav_0at_1b$.

On the other hand, from $M \vDash \text{Env}(t',z^-)$ & $z^-Bz$ & $\text{Firstf}(z,t_1,av_0a,t_2)$, by (5.4),

we have $M \vDash \exists t_7 \text{ Firstf}(z^-,t_1,av_0a,t_7)$.

$\implies$ by choice of $z^+$, $M \vDash \exists t_8 \text{Fr}(z^+,t_1b,av_0a,t_8)$,

$\implies$ from $M \vDash \text{Env}(t^{**},z^{**})$ & $z''=z^{**}t^+auatbb$, by (9.11),

$$M \vDash \forall v,w(v \, \varepsilon \, z^{**} \, \& \, w \, \varepsilon \, z^{**} \to (v<_{z^{**}}w \leftrightarrow v<_{z''}w)).$$

Then from $M \vDash z^{**}\sim z^+\sim z^-$ and $M \vDash \forall v,w(v<_{z''}w \leftrightarrow v<_z w)$,

$M \vDash \forall w(w \, \varepsilon \, z^- \to \forall v(v<_z w \leftrightarrow v<_{z'}w))$, and

$M \vDash \forall w(w \, \varepsilon \, z^+ \to \forall v(v<_{z'}w \leftrightarrow v<_{z^+}w))$ we have that

$$M \vDash \forall v,w(v \, \varepsilon \, z^+ \, \& \, w \, \varepsilon \, z^+ \to (v<_{z^{**}}w \leftrightarrow v<_{z^+}w)).$$

On the other hand, since

$$M \vDash \forall v,w(\neg v \, \varepsilon \, z^+ \, \& \, \neg w \, \varepsilon \, z^+ \to \neg v<_{z^{**}}w \, \& \neg v<_{z^+}w)),$$

we have $M \vDash \forall v,w(\neg v \, \varepsilon \, z^+ \, \& \, \neg w \, \varepsilon \, z^+ \to (v<_{z^{**}}w \leftrightarrow v<_{z^+}w))$.

Therefore, $M \vDash \forall v,w \, (v<_{z^{**}}w \leftrightarrow v<_{z^+}w)$.

Assume now $M \vDash \text{Firstf}(z^{**},t_9,awa,t_{10})$.

$\Rightarrow$ from $M \vDash \text{Firstf}(z^{**},t_1b,av'a,t_1b)$, by (5.15), $M \vDash t_9=t_1b \ \& \ w=v'=v_0$.

Now, we have that

$$M \vDash \text{Set}(z^{**}) \ \& \ \text{Set}(z^+) \ \& \ z^{**} \sim z^+ \ \& \ \forall v,w \ (v<_{z^{**}}w \leftrightarrow v<_{z^+}w).$$

$\Rightarrow$ by (9.12), $M \vDash \exists t_{11},t_{12},v_0 \ \text{Firstf}(z^+,t_{11},av'a,t_{12})$,

$\Rightarrow$ from $M \vDash \text{Env}(t^+,z^+) \ \& \ \text{Fr}(z^+,t_1b,av_0a,t_8)$, by (5.42), $M \vDash t_{11}=t_1b \ \& \ t_{12}=t_8$,

$\Rightarrow M \vDash \text{Firstf}(z^+,t_9,awa,t_8)$.

So we have shown that

$$M \vDash \forall w,t_9,t_{10}(\text{Firstf}(z^{**},t_9,awa,t_{10}) \to \exists t_8 \ \text{Firstf}(z^+,t_9,awa,t_8)).$$

Along with

$$M \vDash \text{Set}(z^{**}) \ \& \ z^{**}=t^{**}av'at^{**} \ \& \ \text{Pref}(av'a,t^{**}) \ \& \ \text{Set}(z^+) \ \& \ z^{**} \sim z^+,$$

this implies, by (5.23), that $M \vDash z^{**}=z^+$, as required.

(2bAii) $M \vDash \exists w_3,w_4,t_3 \ (z''=w_3at^{**}av''at_3aw_4 \ \& \ t_3aw_4=tbbauatbb \ \&$
$$\& \ z^{**}=w_3at^{**}av''at^{**} \ \& \ t^{**}<t_3=tbb).$$

$\Rightarrow$ from $M \vDash \text{Tally}_b(t^{**}) \ \& \ \text{Tally}_b(t_3)$, $M \vDash \exists t^{++}(\text{Tally}_b(t^{++}) \ \& \ t^{**}t^{++}=t_3)$,

$\Rightarrow M \vDash z''=(w_3at^{**}av''at^{**})t^{++}auatbb=z^{**}t^{++}auatbb$,

$\Rightarrow M \vDash z^{**}Bz''$.

Also, from

$$M \vDash \text{Env}(t^{**},z^{**}) \ \& \ \text{Fr}(z^{**},t^{**},av''a,t^{**}) \ \& \ \text{Tally}_b(t^{++}) \ \& \ \text{Tally}_b(tbb),$$

by (5.6), $M \vDash \exists t_4 \ \text{Fr}(z'',t^{**},av''a,t_4)$.

So we have

$$M \vDash \text{Env}(tbb,z'') \ \& \ \text{MinSet}(z'') \ \& \ \text{Fr}(z'',t^{**},av''a,t_4) \ \& \ \text{Env}(t^{**},z^{**}) \ \&$$
$$\& \ \text{Lastf}(z^{**},t^{**},av''a,t^{**}) \ \& \ z^{**}Bz''.$$

$\implies$ by (10.8), $M \vDash \text{MinSet}(z^{**})$.

Next, we show that $M \vDash \forall v,w \ (v<_{z^-}w \leftrightarrow v<_{z^{**}}w)$.

We have, by earlier argument, that $M \vDash \forall w(w \ \varepsilon \ z^- \rightarrow \forall v(v<_{z^-}w \leftrightarrow v<_z w))$.

By hypothesis, $M \vDash \forall v,w \ (v<_z w \leftrightarrow v<_{z''}w)$.

We also have

$\quad M \vDash \text{Env}(t^{**},z^{**}) \ \& \ z'' = z^{**}t^{++}\text{auatbb} \ \& \ \text{Tally}_b(t^{++}) \ \& \ \text{Tally}_b(\text{tbb})$.

$\implies$ from $M \vDash z^{**} \sim z^+ \sim z^-$, $M \vDash \forall v,w(v \ \varepsilon \ z^- \ \& \ w \ \varepsilon \ z^- \rightarrow (v<_{z^-}w \leftrightarrow v<_{z^{**}}w))$.

Analogously to (2bAi), we then obtain $M \vDash \forall v,w(v<_{z^-}w \leftrightarrow v<_{z^{**}}w)$, as required.

We now show that

$\quad M \vDash \forall w,t_1,t_2 \ (\text{Free}^+(z^-,t_1,\text{awa},t_2) \rightarrow \exists t_3 \text{Fr}(z'',t_1b,\text{awa},t_3))$.

Assume $M \vDash \text{Free}^+(z^-,t_1,\text{awa},t_2)$.

Then from

$M \vDash \text{Env}(\text{tb},z) \ \& \ z = z^-t^*\text{auatb} \ \& \ \text{Tally}_b(t^*) \ \& \ \text{Env}(t',z^-) \ \& \ \text{Lastf}(z,\text{tb},\text{aua},\text{tb}) \ \&$

$\quad\quad\quad\quad\quad\quad\quad\quad\quad\quad\quad\quad\quad\quad\quad \& \ \text{Fr}(z^-,t_1,\text{awa},t_2)$,

by (9.24), $M \vDash \exists t_4 \ \text{Free}^+(z,t_1,\text{awa},t_4)$.

$\implies$ by choice of $z''$, $M \vDash \exists t_5 \ \text{Fr}(z'',t_1b,\text{awa},t_5)$,

$\implies$ from hypothesis $M \vDash \text{Free}^+(z^-,t_1,\text{awa},t_2)$, $M \vDash w \ \varepsilon \ z^-$,

$\implies M \vDash w \ \varepsilon \ z^{**}$,

$\implies M \vDash \exists t_6,t_7 \ \text{Fr}(z^{**},t_6,\text{awa},t_7)$.

Now, from $M \vDash \text{Env}(t^{**},z^{**}) \ \& \ z'' = z^{**}t^{++}\text{auatbb} \ \& \ \text{Tally}_b(t^{++}) \ \& \ \text{Tally}_b(\text{tbb})$,

by (5.6), $M \vDash \exists t_8 \ \text{Fr}(z'',t_6,\text{awa},t_8)$.

$\Rightarrow$ from $M \vDash Env(tbb,z'')$, $M \vDash t_6=t_1b$,

$\Rightarrow$ $M \vDash \exists t_7\ Fr(z^{**},t_1b,awa,t_7)$, as required.

Finally, we show that

$M \vDash \forall w,t_1,t_2\ (Bound(z^-,t_1,awa,t_2) \lor Free^-(z^-,t_1,awa,t_2) \to$

$\to \exists t_3 Fr(z^{**},t_1,awa,t_3))$.

Assume that $M \vDash Bound(z^-,t_1,awa,t_2)$.

Then an argument exactly analogous to the one given under the hypothesis

$M \vDash Free(z^-,t_1,awa,t_2)$ establishes the desired conclusion.

Assume that $M \vDash Free^-(z^-,t_1,awa,t_2)$.

$\Rightarrow$ $M \vDash Free^-(z^-,t_1,awa,t_2)\ \&\ \neg Free^+(z^-,t_1,awa,t_2)$,

$\Rightarrow$ by (9.23), $M \vDash \exists t_3\ Free(z,t_1,awa,t_3)$.

Suppose, for a reductio, that $M \vDash Free^+(z,t_1,awa,t_3)$.

$\Rightarrow$ by (9.24), $M \vDash \exists t_4\ Free^+(z^-,t_1,awa,t_4)$,

$\Rightarrow$ $M \vDash Fr(z^-,t_1,awa,t_4)$,

$\Rightarrow$ from $M \vDash Free^-(z^-,t_1,awa,t_2)$, $M \vDash Fr(z^-,t_1,awa,t_2)$,

$\Rightarrow$ by (5.42), $M \vDash t_2=t_4$,

$\Rightarrow$ $M \vDash Free^+(z^-,t_1,awa,t_4)\ \&\ \neg Free^+(z^-,t_1,awa,t_2)$, a contradiction.

Therefore, $M \vDash \neg Free^+(z,t_1,awa,t_3)$.

$\Rightarrow$ since $M \vDash Free(z,t_1,awa,t_3)$, $M \vDash Free^-(z,t_1,awa,t_3)$,

$\Rightarrow$ by hypothesis about $z''$, $M \vDash \exists t_5\ Fr(z'',t_1,awa,t_5)$,

$\Rightarrow$ from $M \vDash z^{**}\sim z^-\ \&\ w\ \varepsilon\ z^-$, $M \vDash \exists t_6,t_7\ Fr(z^{**},t_6,awa,t_7)$,

$\Rightarrow$ by (5.6), $M \vDash \exists t_8\ Fr(z'',t_6,awa,t_8)$,

$\Rightarrow$ from $M \vDash Env(tbb,z'')$, $M \vDash t_6=t_1$,

$\Rightarrow M \vDash \exists t_7\, Fr(z^{**},t_1,awa,t_7)$, as claimed.

Now, we have that $M \vDash Env(t^+,z^+)\, \&\, Env(t^{**},z^{**})$.

$\Rightarrow M \vDash (at^+)Ez^+\, \&\, Tally_b(t^+)\, \&\, (at^*)Ez^{**}\, \&\, Tally_b(t^{**})$,

$\Rightarrow M \vDash \exists z_1\, z^+=z_1at^+\, \&\, \exists z_2\, z^{**}=z_2at^{**}$,

$\Rightarrow$ from $M \vDash z^{**} \sim z^+$, $M \vDash z_1at^+=z_2at^{**}$,

$\Rightarrow$ from $M \vDash Tally_b(t^+)\, \&\, Tally_b(t^{**})$, by (4.24$^b$), $M \vDash t^+=t^{**}$.

We proceed to distinguish the cases:

   (2bAi)   $M \vDash z^{**}=t^{**}av''at^{**}\, \&\, v''=v'$ and $M \vDash z''=t^{**}av'atbbauatbb$.

    (2bia)   $M \vDash z'=z^+bauatbb\, \&\, t^+=tb$.

$\Rightarrow M \vDash z'=z^+bauatbb=z^{**}bauatbb=(t^{**}av'at^{**})bauatbb=t^{**}av'at^+bauatbb=$

$$=t^{**}av'atbbauatbb=z'',$$

as required.

    (2biia)   $M \vDash z'=z^+bbauatbb\, \&\, t^+b=tb$.

$\Rightarrow M \vDash z'=z^+bbauatbb=z^{**}bbauatbb=(t^{**}av'at^{**})bbauatbb=$

$$=t^{**}av'a(t^+b)bauatbb=t^{**}av'a(tb)bauatbb=z''.$$

   (2bAii)   $M \vDash z^{**}=w_3at^{**}av''at^{**}\, \&\, t^{**}t^{++}=tbb$ and

$$M \vDash z''=t^{**}\, t^{++}auatbb.$$

    (2bia)   $M \vDash z'=z^+bauatbb\, \&\, t^+=tb$.

$\Rightarrow$ from $M \vDash t^{**}=t^+=tb$, $M \vDash t^{**}b=tbb$, $\Rightarrow M \vDash t^{**}t^{++}=tbb=t^{**}b$,

$\Rightarrow$ by (3.7), $M \vDash t^{++}=b$,

$\Rightarrow M \vDash z'=z^+bauatbb=z^{**}bauatbb=z^{**}bauatbb=z^{**}t^{++}auatbb=z''$,

as required.

  (2biia)  $M \vDash z'=z^+bbauatbb\ \&\ t^+=tb$.

$\Rightarrow$ from $M \vDash t^{**}=t^+$, $M \vDash t^{**}bb=t^+bb=tbb$,

$\Rightarrow M \vDash t^{**}t^{++}=tbb=t^{**}bb$,

$\Rightarrow$ by (3.7), $M \vDash t^{++}=bb$,

$\Rightarrow M \vDash z'=z^+bbauatbb=z^{**}t^+auatbb=z''$, as required.

 (2bB)  Suppose that, as in (2bib) and (2biia) and (2biib),

$$M \vDash (Bound(z,tb,aua,tb) \vee Free^-(z,tb,aua,tb))\ \&\ t''=tb.$$

$\Rightarrow\ M \vDash Lastf(z'',tb,aua,tb)$,

$\Rightarrow\ M \vDash z''=tbauatb \vee \exists w_2\,(z''=w_2atbauatb\ \&\ Max^+T_b(tb,w_2))$.

We then proceed exactly analogously to (2bA) replacing tbb with tb to obtain

$$M \vDash z^{**}=z^+,$$

and further $M \vDash t^+=t^{**}$.

We then distinguish the cases:

  (2bBi)  $M \vDash z^{**}=t^{**}av''at^{**}\ \&\ v''=v'$ and $M \vDash z''=t^{**}av'atbauatb$.

   (2bib)  $M \vDash z'=z^+t^*auatb\ \&\ t^+=t'b$ where $M \vDash t'bt^*=tb$.

$\Rightarrow M \vDash z'=z^+t^*auatb=z^{**}t^*auatb=(t^{**}av''at^{**})t^*auatb=t^{**}av'at^+t^*auatb=$

    $=t^{**}av'at'bt^*auatb=t^{**}av'atbauatb=z''$, as required.

   (2biib)  $M \vDash z'=z^+bt^*auatb\ \&\ t^+=t'$ where $M \vDash t'bt^*=tb$.

$\Rightarrow M \vDash z'=z^+bt^*auatb=z^{**}bt^*auatb=(t^{**}av'at^{**})bt^*auatb=$

  $=t^{**}av'at^+bt^*bauatb=t^{**}av'at'bt^*auatb=t^{**}av'atbauatb=z''.$

  (2bBii)  $M \vDash z^{**}=w_3at^{**}av''at^{**}\ \&\ t^{**}t^{++}=tb$ and $M \vDash z''=z^{**}\,t^{++}auatb$

where $M \vDash t'bt^*=tb$.

(2bib) $M \vDash z'=z^+t^*auatb$ & $t^+=t'b$.

$\Rightarrow M \vDash t^+t^{++}=t^{**}t^{++}=tb=t'bt^*=t^+t^*$,

$\Rightarrow$ by (3.7), $M \vDash t^{++}=t^*$,

$\Rightarrow M \vDash z'=z^+t^*auatb=z^{**}t^*auatb=z^{**}t^{++}auatbb=z''$, as required.

(2biib) $M \vDash z'=z^+bt^*auatb$ & $t^+=t'$.

$\Rightarrow M \vDash t^+t^{++}=t^{**}t^{++}=tb=t'bt^*=t^+bt^*$,

$\Rightarrow$ by (3.7), $M \vDash t^{++}=bt^*$,

$\Rightarrow M \vDash z'=z^+bt^*auatb=z^{**}bt^*auatb=z^{**}t^{++}auatb=z''$, as required.

This completes the proof of the uniqueness of z' and the proof of (2b).

Thus we have shown that J(x) is a string concept.

To complete the proof of (10.29), assume that

$M \vDash Env(t,z)$ & $MinSet(z)$    where $M \vDash J(t)$ & $I^*(z)$.

Then

$M \vDash \exists!z' \in I^* \exists t'(Env(t',z')$ & $t' \leq tb$ & $z \sim z'$ &

& $MinSet(z')$ & $\forall v,w\ (v<_z w \leftrightarrow v<_{z'} w)$ &

& $\forall w,t_1,t_2\ (Free^+(z,t_1,awa,t_2) \to \exists t_3 Fr(z',t_1 b,awa,t_3))$ &

& $\forall w,t_1,t_2\ (Bound(z,t_1,awa,t_2) \lor Free^-(z,t_1,awa,t_2) \to \exists t_3 Fr(z',t_1 b,awa,t_3)))$.

Since J has been proved to be a string concept, we have that

$M \vDash \exists!z' \in I^* \exists t' \in J\ (Env(t',z')$ & $z \sim z'$ &

& $MinSet(z')$ & $\forall v,w\ (v<_z w \leftrightarrow v<_{z'} w)$ &

& $\forall w,t_1,t_2\ (Free^+(z,t_1,awa,t_2) \to \exists t_3 Fr(z',t_1 b,awa,t_3))$ &

& $\forall w, t_1, t_2$ (Bound($z, t_1, awa, t_2$) v Free$^-$($z, t_1, awa, t_2$) → $\exists t_3$Fr($z', t_1 b, awa, t_3$)))

as required.

This completes the proof of (10.29).

# 11. The Uniqueness Lemma

(11.1) For any string concept $I \subseteq I_0$ there is a string concept $J \subseteq I$ such that

$QT^+ \vdash \forall x,y \in J$ (Lex$^+$(x) & Lex$^+$(y) & Special(x) & Special(y) & x~y $\rightarrow$

$\rightarrow \neg xBy$ & $\neg yBx$).

Let $J \equiv I_{4.20}$ & $I_{9.2}$ & $I_{10.2}$.

Assume $M \vDash x{\sim}y$

where $M \vDash$ Lex$^+$(x) & Lex$^+$(y) & Special(x) & Special(y) and $M \vDash J(x)$ & $J(y)$.

Then $M \vDash$ Set(x). So $M \vDash x=aa \lor \exists t$ Env(t,x).

If $M \vDash x=aa$, then, by (5.18), $M \vDash \forall u \neg(u \varepsilon x)$. Then $M \vDash \forall u \neg(u \varepsilon y)$ from

$M \vDash x{\sim}y$, whence, again by (5.18), $M \vDash y=aa$. But then $M \vDash \neg xBy$ & $\neg yBx$

follows from $M \vDash I(x)$ & $I(y)$ and $I \subseteq I_0$, as required.

So we may assume that $M \vDash \exists t$ Env(t,x).

$\Rightarrow$ from $M \vDash x{\sim}y$, by (5.18), $M \vDash y \neq aa$,

$\Rightarrow M \vDash \exists t'$ Env(t',y),

$\Rightarrow M \vDash \exists u \subseteq_p x$ Lastf(x,t,u,t) and $M \vDash \exists v \subseteq_p y$ Lastf(y,t',v,t'),

$\Rightarrow$ by (9.2), $M \vDash u=v$,

$\Rightarrow M \vDash$ Fr(x,t,u,t) & Fr(y,t',u,t'),

$\Rightarrow$ by (10.2), $M \vDash t=t'$,

$\Rightarrow M \vDash$ Lastf(x,t,u,t) & Lastf(y,t,u,t),

$\Rightarrow M \vDash$ (x=tut $\lor \exists w_1(w_1$atut=x & Max$^+T_b(t,w_1))$) &

$$\& \ (y=tut \ \lor \ \exists w_2(w_2atut=x \ \& \ Max^+T_b(t,w_2))).$$

We distinguish four cases:

<u>Case 1.</u>  $M \vDash x=tut \ \& \ y=tut$.

$\Rightarrow M \vDash x=y$,

$\Rightarrow M \vDash \neg xBy \ \& \ \neg yBx$   since $M \vDash x,y \in I \subseteq I_0$.

<u>Case 2.</u>  $M \vDash \exists w_1(w_1atut=x \ \& \ Max^+T_b(t,w_1)) \ \& \ y=tut$.

From $M \vDash y=tut$, we have $M \vDash Firstf(y,t,u,t)$.

From $M \vDash Env(t,x)$ we have that $M \vDash \exists v \subseteq_p x \ \exists t_1,t_2 \subseteq_p x \ Firstf(x,t_1,v,t_2)$.

$\Rightarrow$ by (9.2),  $M \vDash v=u$.

$\Rightarrow$ from $M \vDash Fr(x,t_1,u,t_2) \ \& \ Fr(y,t,u,t)$, by (10.2),

$$M \vDash t_1=t,$$

$\Rightarrow M \vDash Firstf(x,t,u,t_2)$,

$\Rightarrow M \vDash (tu)Bx$,

$\Rightarrow M \vDash \exists x_2 \ tux_2=x=w_1atut$,

$\Rightarrow$ by (4.14$^b$),  $M \vDash w_1=t \ \lor \ tBw_1$,

$\Rightarrow M \vDash t \subseteq_p w_1$, which contradicts $M \vDash Max^+T_b(t,w_1)$.

<u>Case 3.</u>  $M \vDash x=tut \ \& \ \exists w_2(w_2atut=y \ \& \ Max^+T_b(t,w_2))$.

Exactly analogous to Case 2.

<u>Case 4.</u>  $M \vDash \exists w_1(w_1atut=x \ \& \ Max^+T_b(t,w_1)) \ \&$

$$\& \ \exists w_2(w_2atut=y \ \& \ Max^+T_b(t,w_2)).$$

Suppose, for a reductio, that $M \vDash xBy$.

$\Rightarrow M \vDash \exists x_1 \ xx_1=y$,

$\implies$ M ⊨ $(w_1 atut)x_1 = y = w_2 atut$,

$\implies$ M ⊨ $(w_1 at)By$ & $(w_2 at)By$ & $Max^+T_b(t,w_1)$ & $Max^+T_b(t,w_2)$,

$\implies$ by (4.20), M ⊨ $w_1 = w_2$,

$\implies$ M ⊨ $(w_1 atut)x_1 = y = w_1 atut$,

$\implies$ M ⊨ $yBy$, which contradicts M ⊨ $y \in I \subseteq I_0$.

Hence M ⊨ $\neg xBy$.

Analogously, M ⊨ $\neg yBx$.

This completes the proof of (11.1).

(11.2) For any string concept $I \subseteq I_0$ there is a string concept $J \subseteq I$ such that

$$QT^+ \vdash \forall x,y \in J \ (Set^*(x) \ \& \ Set^*(y) \ \& \ x \sim y \ \rightarrow \ \neg \exists z \ (zaBx \ \& \ zbBy)).$$

Let $J \equiv I_{11.1}$.

Assume $M \vDash Set^*(x) \ \& \ Set^*(y)$ where $M \vDash x \sim y$ and $M \vDash J(x) \ \& \ J(y)$.

Suppose, for a reductio, that $M \vDash zaBx \ \& \ zbBy$.

$\Rightarrow M \vDash \exists x_1 \ zax_1 = x \ \& \ \exists y_1 \ zby_1 = y$.

From $M \vDash Set^*(x)$ we have $M \vDash MinSet(x)$.

$\Rightarrow M \vDash \exists u \subseteq_p x \ \exists t_1,t_2 \subseteq_p x \ Occ(z,a,x_1,x,t_1,u,t_2)$,

$\Rightarrow M \vDash Fr(x,t_1,u,t_2)$.

We distinguish three cases:

<u>Case 1.</u> $M \vDash Firstf(x,t_1,u,t_2)$.

From $M \vDash Set(y)$ we have, since $M \vDash x \sim y$, by (5.18), that $M \vDash y \neq aa$.

$\Rightarrow M \vDash \exists v \subseteq_p y \ \exists t_3,t_4 \subseteq_p y \ Firstf(y,t_3,v,t_4)$.

From $M \vDash Set^*(x) \ \& \ Set^*(y)$ we have $M \vDash Lex^+(x) \ \& \ Lex^+(y)$.

$\Rightarrow$ by (9.2), $M \vDash v = u$.

Since $M \vDash Special(x) \ \& \ Special(y) \ \& \ Fr(x,t_1,u,t_2) \ \& \ Fr(y,t_3,v,t_4)$, it follows by

(10.2) that $M \vDash t_3 = t_1$.

From $M \vDash Occ(z,a,x_1,x,t_1,u,t_2)$, we have three subcases:

(1ai) $M \vDash t_1 = z$.

By hypothesis (1) we have that $M \vDash zax_1 = x \ \& \ zby_1 = y$.

But from $M \vDash Firstf(y,t_1,u,t_4)$ we also have that $M \vDash (t_1a)By$.

$\Rightarrow$ M ⊨ $\exists y_2\ t_1ay_2=y$,

$\Rightarrow$ from M ⊨ $t_1=z$, M ⊨ $t_1by_1=y=t_1ay_2$,

$\Rightarrow$ by (3.7), M ⊨ $by_1=ay_2$, a contradiction.

(1aii)  M ⊨ $(za)B(t_1u)$.

$\Rightarrow$ M ⊨ $\exists z_1\ zaz_1=t_1u$.

From M ⊨ $Firstf(y,t_1,u,t_4)$ we have M ⊨ $(t_1u)By$.

$\Rightarrow$ M ⊨ $\exists y_2\ t_1uy_2=y$,

$\Rightarrow$ M ⊨ $t_1uy_2=y=(zaz_1)y_2$.

But then from M ⊨ $zby_1=y$ it follows that M ⊨ $zaz_1y_2=y=zby_1$, whence, by

(3.7), M ⊨ $az_1y_2=by_1$, a contradiction.

(1aiii)  M ⊨ $za=t_1u$.

Same as (1aii) except that $z_1$ is omitted throughout.

Case 2.  M ⊨ $\exists w_1 Intf(x,w_1,t_1,u,t_2)$.

$\Rightarrow$ M ⊨ $Pref(u,t_1)$ & $Max^+T_b(t_1,w_1)$.

From M ⊨ $Pref(u,t_1)$ we have that M ⊨ $\exists u_0\ u=au_0a$. So M ⊨ $u_0\ \varepsilon\ x$.

$\Rightarrow$ from the hypothesis M ⊨ $x\sim y$, M ⊨ $u_0\ \varepsilon\ y$.

$\Rightarrow$ M ⊨ $\exists t_3,t_4\subseteq_p y\ Fr(y,t_3,u,t_4)$.

Suppose that M ⊨ $Firstf(y,t_3,u,t_4)$.

From M ⊨ $Set(x)$ we have, by (5.18), that M ⊨ $\exists v\subseteq_p x\ \exists t',t''\subseteq_p x\ Firstf(x,t',v,t'')$.

Then, by (9.2), M ⊨ $v=u$, so M ⊨ $Firstf(x,t',u,t'')$. But this contradicts

M ⊨ $Intf(x,w_1,t_1,u,t_2)$ by (5.19). Hence M ⊨ $\neg Firstf(y,t_3,u,t_4)$.

Likewise, M ⊨ $\neg Lastf(y,t_3,u,t_4)$.

Therefore, $M \vDash \exists w_3 \text{Intf}(x,w_3,t_3,u,t_4)$.

$\Rightarrow M \vDash \exists w_4 \ w_3 a t_3 u t_4 a w_4 = y \ \& \ \text{Max}^+ T_b(t_3,w_3)$.

From $M \vDash \text{Lex}^+(x) \ \& \ \text{Lex}^+(y) \ \& \ \text{Special}(x) \ \& \ \text{Special}(y)$ it follows that

$M \vDash t_1 = t_3$ by (10.2).

Again, from $M \vDash \text{Occ}(z,a,x_1,x,t_1,u,t_2)$ we distinguish three subcases:

(2i)   $M \vDash w_1 a t_1 = z$.

By hypothesis we have that $M \vDash z a x_1 = x \ \& \ z b y_1 = y$.

$\Rightarrow M \vDash w_3 a t_1 u t_4 a w_4 = y = w_1 a t_1 b y_1$,

$\Rightarrow M \vDash (w_3 a t_1) B y \ \& \ (w_1 a t_1) B y \ \& \ \text{Max}^+ T_b(t_1,w_3) \ \& \ \text{Max}^+ T_b(t_1,w_1)$,

$\Rightarrow$ by (4.20), $M \vDash w_3 = w_1$,

$\Rightarrow$ from $M \vDash w_3 a t_1 u t_4 a w_4 = y = z b y_1 = w_1 a t_1 b y_1$, $M \vDash w_1 a t_1 u t_4 a w_4 = y = w_1 a t_1 b y_1$,

$\Rightarrow M \vDash w_1 a t_1 (a u_0 a) t_4 a w_4 = w_1 a t_1 b y_1$,

$\Rightarrow$ by (3.7), $M \vDash a u_0 a t_4 a w_4 = b y_1$, a contradiction.

(2ii)   $M \vDash \exists u_1 (u_1 B u \ \& \ w_1 a t_1 u_1 = z a)$.

We have from the hypothesis of Case 2 that

$\qquad M \vDash \exists w_2 \ w_1 a t_1 u t_2 a w_2 = x$   and   $M \vDash z a x_1 = x \ \& \ z b y_1 = y$.

$\Rightarrow M \vDash (w_1 a t_1 u_1) x_1 = x$.

From $M \vDash \text{Pref}(u,t_1)$ we have $M \vDash \exists u_0 \ u = a u_0 a$,

$\Rightarrow$ from $M \vDash u_1 B u$, $M \vDash \exists u' \ u_1 u' = u = a u_0 a$,

$\Rightarrow M \vDash u_1 = a \lor a B u_1$.

If $M \vDash u_1 = a$, then $M \vDash w_1 a t_1 a = w_1 a t_1 u_1 = z a$, whence $M \vDash w_1 a t_1 = z$. We then proceed exactly as in (2i).

So we may assume that $M \vDash aBu_1$.

$\implies M \vDash \exists u_2\, au_2 = u_1$,

$\implies M \vDash w_1 at_1 au_2 = w_1 at_1 u_1 = za$,

$\implies M \vDash u_2 = a \lor aEu_2$.

(2iia) $M \vDash u_2 = a$.

$\implies M \vDash u_1 = aa$,

$\implies M \vDash w_1 at_1 aa = w_1 at_1 u_1 = za$,

$\implies M \vDash w_1 at_1 a = z$,

$\implies M \vDash x = zax_1 = (w_1 at_1 a)a\, x_1$ and $M \vDash y = zby_1 = (w_1 at_1 a)by_1$.

So we have $M \vDash w_3 at_1 ut_4 aw_4 = y = w_1 at_1 by_1$,

$\implies M \vDash (w_3 at_1)By\ \&\ (w_1 at_1)By\ \&\ Max^+T_b(t_1,w_3)\ \&\ Max^+T_b(t_1,w_1)$,

$\implies$ by (4.20), $M \vDash w_1 = w_3$,

$\implies$ from $M \vDash w_3 at_1 ut_4 aw_4 = y = w_1 at_1 by_1$,

$\qquad M \vDash w_1 at_1 ut_4 aw_4 = y = (w_1 at_1 a)by_1 = zby_1$.

On the other hand, we have

$\qquad M \vDash w_1 at_1 ut_2 aw_2 = x = (w_1 at_1 a)a\, x_1 = zax_1$.

So we have that

$\qquad M \vDash w_1 at_1 (au_0 a)t_2 aw_2 = x = (w_1 at_1 a)a\, x_1$

and $\qquad M \vDash w_1 at_1 (au_0 a)t_4 aw_4 = y = (w_1 at_1 a)by_1$,

$\implies$ by (3.7), $M \vDash u_0 at_2 aw_2 = a\, x_1$ and $M \vDash u_0 at_4 aw_4 = by_1$,

$\implies M \vDash (u_0 = b \lor bBu_0)\ \&\ (u_0 = a \lor aBu_0)$, a contradiction in each case.

(2iib)  $M \vDash aEu_2$.

$\Rightarrow M \vDash \exists u_3\ u_2=u_3a$,

$\Rightarrow M \vDash w_1at_1a(u_3a)=w_1at_1u_1=za$,

$\Rightarrow M \vDash w_1at_1au_3=z$,

$\Rightarrow M \vDash x=zax_1=(w_1at_1au_3)a\,x_1$  and  $M \vDash y=zby_1=(w_1at_1au_3)by_1$,

$\Rightarrow M \vDash w_3at_1ut_4aw_4=y=w_1at_1au_3by_1$,

$\Rightarrow M \vDash (w_3at_1)By\ \&\ (w_1at_1)By\ \&\ Max^+T_b(t_1,w_3)\ \&\ Max^+T_b(t_1,w_1)$,

$\Rightarrow$ by (4.20), $M \vDash w_1=w_3$,

$\Rightarrow M \vDash w_1at_1ut_4aw_4=y=w_1at_1au_3by_1$.

From hypothesis (2ii), $M \vDash u_1Bu$, whence $M \vDash \exists u'\ u_1u'=u$.

We also have that $M \vDash u_3a=u_2\ \&\ au_2=u_1$, so $M \vDash u_1=au_3a$.

$\Rightarrow M \vDash w_1at_1(au_3a)u't_4aw_4=y=zby_1=w_1at_1au_3by_1$.

But then, by (3.7), $M \vDash au't_4aw_4=by_1$, a contradiction.

(2iii)  $M \vDash w_1at_1u=za$.

$\Rightarrow M \vDash \exists u_0\ u=au_0a$,

$\Rightarrow M \vDash w_1at_1(au_0a)=za$,

$\Rightarrow M \vDash w_1at_1au_0=z$.

Therefore,

  $M \vDash w_1at_1ut_2aw_2=x=zax_1=(w_1at_1au_0)ax_1$ and  $M \vDash y=zby_1=(w_1at_1au_0)by_1$,

$\Rightarrow M \vDash w_3at_1ut_4aw_4=y=w_1at_1au_0by_1$,

$\Rightarrow M \vDash (w_3at_1)By\ \&\ (w_1at_1)By\ \&\ Max^+T_b(t_1,w_3)\ \&\ Max^+T_b(t_1,w_1)$,

$\Rightarrow$ by (4.20), $M \vDash w_1=w_3$.

In particular, $M \vDash w_1at_1ut_4aw_4=y=(w_1at_1au_0)by_1$,

$\Longrightarrow M \vDash w_1at_1(au_0a)t_4aw_4=y=(w_1at_1au_0)by_1$,

$\Longrightarrow$ by (3.7), $M \vDash at_4aw_4=by_1$, a contradiction.

<u>Case 3.</u> $M \vDash \exists w_1(Lastf(x,t_1,u,t_2) \& w_1at_1ut_2=x)$.

From $M \vDash Set(y)$ by (5.18) we have that $M \vDash \exists v\subseteq_p y\ \exists t_3,t_4\subseteq_p y\ Lastf(y,t_3,v,t_4)$,

$\Longrightarrow$ by (9.2), $M \vDash v=u$.

Since

$M \vDash Lex^+(x) \& Lex^+(y) \& Special(x) \& Special(y) \& Fr(x,t_1,u,t_2) \& Fr(y,t_3,u,t_4)$,

it follows by (10.2) that $M \vDash t_3=t_1$.

Again, from $M \vDash Occ(z,a,x_1,x,t_1,u,t_2)$, we have three subcases:

(3i) $M \vDash w_1at_1=z$.

From $M \vDash Lastf(y,t_1,u,t_4)$ we have

$$M \vDash y=t_1ut_4 \vee \exists w_3 (w_3at_1ut_4=y \& Max^+T_b(t_1,w_3)).$$

If $M \vDash y=t_1ut_4$, then $M \vDash w_1at_1by_1=zby_1=y=t_1ut_4$, whence $M \vDash w_1=t_1 \vee t_1Bw_1$.

But then $M \vDash t_1 \subseteq_p w_1$, which contradicts $M \vDash Max^+T_b(t_1,w_1)$. So we have that

$$M \vDash \exists w_3 (w_3at_1ut_4=y=zby_1=w_1at_1by_1),$$

$\Longrightarrow M \vDash (w_3at_1)By \& (w_1at_1)By \& Max^+T_b(t_1,w_3) \& Max^+T_b(t_1,w_1)$,

$\Longrightarrow$ by (4.20), $M \vDash w_1=w_3$,

$\Longrightarrow M \vDash w_1at_1ut_4=y=w_1at_1by_1$.

From $M \vDash Pref(u,t_1)$ we have $M \vDash \exists u_0\ u=au_0a$.

$\Longrightarrow M \vDash w_1at_1(au_0a)t_4=y=w_1at_1by_1$,

$\Longrightarrow$ by (3.7), $M \vDash au_0at_4=by_1$, a contradiction.

(3ii)   $M \vDash \exists u_1(u_1Bu \ \& \ w_1at_1u_1=za)$.

We have that

$M \vDash \exists w_1 \ x=w_1at_1ut_2 \ \& \ \exists w_3 \ y=w_3at_1ut_4$  and  $M \vDash zax_1=x \ \& \ zby_1=y$.

We derive a contradiction exactly as in (2ii), omitting $aw_2$ and $aw_4$ throughout the argument.

(3iii)   $M \vDash w_1at_1u=za$.

We reason exactly as in (2iii), omitting $aw_2$ and $aw_4$ throughout.

This completes the proof of (11.2).

(11.3) For any string concept $I \subseteq I_0$ there is a string concept $J \subseteq I$ such that

$$QT \vdash \forall x,y \in J\ (Set^*(x)\ \&\ Set^*(y)\ \&\ x \sim y\ \to\ \neg \exists z\ (zaBx\ \&\ zb=y)).$$

Let $J \equiv I_{4.5}\ \&\ I_{11.2}$.

Assume $M \vDash Set^*(x)\ \&\ Set^*(y)$ where $M \vDash x \sim y$ and $M \vDash J(x)\ \&\ J(y)$.

$\Rightarrow M \vDash Lex^+(x)\ \&\ Lex^+(y)\ \&\ Special(x)\ \&\ Special(y)$.

Suppose, for a reductio, that $M \vDash zaBx\ \&\ zb=y$.

$\Rightarrow M \vDash \exists x_1\ zax_1=x$.

From $M \vDash zb=y$, we have that $M \vDash zb \neq aa$. Hence from $M \vDash Set(y)$ it follows that $M \vDash \exists t'\ Env(t',y)$.

$\Rightarrow M \vDash \exists u \subseteq_p y\ Lastf(y,t',u,t')$,

$\Rightarrow$ from $M \vDash Pref(u,t')$,

$\Rightarrow M \vDash \exists u_0\ u=au_0a$.

From $M \vDash u_0\ \varepsilon\ y\ \&\ x \sim y$, we have $M \vDash u_0\ \varepsilon\ x$.

From $M \vDash Set(x)$, by (5.18), $M \vDash \exists t\ Env(t,x)$.

$\Rightarrow M \vDash \exists v \subseteq_p x\ Lastf(x,t,v,t)$,

$\Rightarrow$ by (9.2), $M \vDash v=u$.

Since $M \vDash Lex^+(x)\ \&\ Lex^+(y)\ \&\ Special(x)\ \&\ Special(y)$, it follows by (10.2) that $M \vDash t=t'$.

From $M \vDash Lastf(x,t,u,t)\ \&\ Lastf(y,t,u,t)$ we have that

$M \vDash (x=tut\ \vee\ \exists w_1(w_1 atut=x\ \&\ Max^+T_b(t,w_1)))\ \&$

$\&\ (y=tut\ \vee\ \exists w_3(w_3 atut=y\ \&\ Max^+T_b(t,w_3)))$

We distinguish four cases:

(a) $M \vDash x=tut$ & $y=tut$.

$\Rightarrow M \vDash zax_1=x=tut=y=zb$,

$\Rightarrow$ by (3.7), $M \vDash ax_1=b$, a contradiction.

(b) $M \vDash x=tut$ & $\exists w_3\, y=w_3 atut$.

$\Rightarrow M \vDash zax_1=x=tut$ & $zb=y=w_3 atut$,

$\Rightarrow$ since $M \vDash Tally_b(t)$, by $(4.14^b)$, $M \vDash z=t \lor tBz$.

If $M \vDash z=t$, then $M \vDash tb=zb=y=w_3 atut$, a contradiction because, by (4.5),

$M \vDash Tally_b(tb)$.

If $M \vDash tBz$, then $M \vDash \exists z_1\, tz_1=z$, so $M \vDash (tz_1)b=zb=y=w_3 atut$,

$\Rightarrow$ by $(4.14^b)$, $M \vDash w_3=t \lor tBw_3$,

$\Rightarrow M \vDash t\subseteq_p w_3$, which contradicts $M \vDash Max^+T_b(t,w_3)$.

(c) $M \vDash \exists w_1\, x=w_1 atut$ & $y=tut$.

From $M \vDash y=tut$ we have $M \vDash Firstf(y,t,u,t)$. From $M \vDash Env(t,x)$, we have that

$\quad\quad M \vDash \exists v\subseteq_p x\, \exists t_1,t_2\subseteq_p x\, Firstf(x,t_1,v,t_2)$,

$\Rightarrow$ by (9.2), $M \vDash v=u$.

From

$M \vDash Lex^+(x)$ & $Lex^+(y)$ & $Special(x)$ & $Special(y)$ & $Fr(x,t_1,u,t_2)$ & $Fr(y,t,u,t)$

it follows that $M \vDash t_1=t$ by (10.2).

$\Rightarrow M \vDash Firstf(x,t,u,t_2)$,

$\Rightarrow M \vDash (tu)Bx$,

$\Rightarrow M \vDash \exists x_2\, tux_2=x=w_1 atut$,

$\Rightarrow$ by (4.14$^b$), $M \vDash w_1=t \lor tBw_1$,

$\Rightarrow M \vDash t\subseteq_p w_1$, which contradicts $M \vDash \text{Max}^+T_b(t,w_1)$.

(d) $M \vDash \exists w_1 x=w_1 atut \ \& \ \exists w_3 \, y=w_3 atut$.

$\Rightarrow M \vDash w_1 atut=x=zax_1 \ \& \ w_3 atut=y=zb$.

From $M \vDash \text{Tally}_b(t)$, we have, by (4.2), that $M \vDash t=b \lor \exists t^- (\text{Tally}_b(t^-) \ \& \ t^-b=t)$,

$\Rightarrow M \vDash (t=b \ \& \ w_3 abub=zb) \lor (t=t^-b \ \& \ w_3 atu(t^-b)=zb)$,

$\Rightarrow M \vDash (t=b \ \& \ z=w_3 abu) \lor (t=t^-b \ \& \ z=w_3 atut^-)$,

$\Rightarrow$ from $M \vDash w_1 atut=x=zax_1$,

$\quad M \vDash (t=b \ \& \ w_1 abub=x=(w_3 abu)ax_1) \lor (t=t^-b \ \& \ w_1 atut=x=(w_3 atut^-)ax_1)$.

Either way, $M \vDash (w_1 at)Bx \ \& \ (w_3 at)Bx$. By exactly the same argument as in

(11.2)(2i) we derive $M \vDash w_1=w_3$.

$\Rightarrow M \vDash w_1 abub=w_1 abuax_1 \lor w_1 atu(t^-b)=w_1 atut^-ax_1$,

$\Rightarrow$ either way by (3.7), $M \vDash b=ax_1$, a contradiction.

This completes the proof of (11.3).

THE UNIQUENESS LEMMA. (11.4) For any string concept $I \subseteq I_0$ there is a string concept $J \subseteq I$ such that

$$QT^+ \vdash \forall x,y \in J \ (Set^*(x) \ \& \ Set^*(y) \ \& \ x \sim y \rightarrow x=y).$$

Let $J \equiv I_{11.3}$.

Assume $M \vDash Set^*(x) \ \& \ Set^*(y)$ where $M \vDash x \sim y$ and $M \vDash J(x) \ \& \ J(y)$.

Suppose, for a reductio, that $M \vDash x \neq y$.

Since from $M \vDash Set(x)$ we have $M \vDash \exists z \ zBx$, we may apply the LEFT ROOT LEMMA to derive

$M \vDash y=a \lor y=b \lor (aBx \ \& \ bBy) \lor (bBx \ \& \ aBy) \lor xBy \lor yBx \lor \exists z \ Rt_L(z,x,y)$.

From $M \vDash Set(x) \ \& \ Set(y)$ we immediately have that

$$M \vDash y \neq a \ \& \ y \neq b \ \& \ \neg(aBx \ \& \ bBy) \ \& \ \neg(bBx \ \& \ aBy).$$

By (11.3), from $M \vDash Set(x) \ \& \ Set(y)$ we have

$$M \vDash \neg xBy \ \& \ \neg yBx.$$

So the remaining possibility is $M \vDash \exists z \ Rt_L(z,x,y)$, that is,

$M \vDash \exists z \ (zBx \ \& \ zBy \ \& \ (((zaBx \lor za=x) \ \& \ (zbBy \lor zb=y)) \lor$

$\lor ((zbBx \lor zb=x) \ \& \ (zaBy \lor za=y)))$.

There are eight cases to consider:

(1) $M \vDash zaBx \ \& \ zbBy$,      (5) $M \vDash zbBx \ \& \ zaBy$,

(2) $M \vDash za=x \ \& \ zbBy$,      (6) $M \vDash zb=x \ \& \ zaBy$,

(3) $M \vDash zaBx \ \& \ zb=y$,      (7) $M \vDash zbBx \ \& \ za=y$,

(4) $M \vDash za=x \ \& \ zb=y$,      (8) $M \vDash zb=x \ \& \ za=y$.

By symmetry, we need only consider (1)-(4).

From the hypothesis $M \vDash x \neq y$ & $x \sim y$, we have that $M \vDash x \neq aa \neq y$.

For, if $M \vDash x = aa$, then $M \vDash \forall z \neg(z \; \varepsilon \; x)$ by (5.18).

$\Rightarrow$ from $M \vDash x \sim y$, $M \vDash \forall z \neg(z \; \varepsilon \; y)$,

$\Rightarrow$ again by (5.18), $M \vDash y = aa$, contradicting the hypothesis $M \vDash x \neq y$.

Similarly if $M \vDash y = aa$.

Hence $M \vDash \text{Set}(x)$ & $\text{Set}(y)$ we have $M \vDash \exists t \; \text{Env}(t,x)$ & $\exists t' \; \text{Env}(t',y)$, which rules out (2) and (4).

So the only cases to consider are (1) and (3). But case (1) is ruled out by (11.2) and case (3) by (11.3).

This completes the proof of the UNIQUENESS LEMMA.

Let $\text{Max}_{\leqslant}(u,x,y)$ abbreviate

$u \,\varepsilon\, x \;\&\; \forall z\, (z \leq_x u \to z \prec y) \;\&\; \forall z\, (z \,\varepsilon\, x \;\&\; \forall v\, (v \leq_x z \to v \prec y) \to z \leq_x u)$.

(11.5)  For any string concept $I \subseteq I_0$ there is a string concept $J \subseteq I$ such that

   $QT^+ \vdash \forall x,y \in J \;\forall u\, (\text{Lex}^+(x) \;\&\; \neg(y \,\varepsilon\, x) \;\&\; \text{Max}_{\leqslant}(u,x,y) \to \forall v(u \leq_x v \to y \prec v))$.

Let  $J(x) \equiv I_{8.2} \;\&\; I_{8.3} \;\&\; I_{9.6}$.

Assume  $M \vDash \text{Lex}^+(x) \;\&\; \neg(y \,\varepsilon\, x) \;\&\; \text{Max}_{\leqslant}(u,x,y)$  where  $M \vDash J(x)$.

Let  $M \vDash u <_x v$.

Assume, for a reductio, that  $M \vDash v \prec y$.

$\Rightarrow$ from  $M \vDash \text{Lex}^+(x)$,  $M \vDash \forall w(w <_x v \to w \prec v)$,

$\Rightarrow$ by (8.3),  $M \vDash \forall w(w <_x v \to w \prec y)$,

$\Rightarrow$ from  $M \vDash \text{Max}_{\leqslant}(u,x,y)$,  $M \vDash v \leq_x u$.

But this contradicts  $M \vDash u <_x v$  by (9.4) and (9.6).

Therefore  $M \vDash \neg(v \prec y)$.

$\Rightarrow$ since  $M \vDash u <_x v$,  $M \vDash v \,\varepsilon\, x$,

$\Rightarrow$  from hypothesis  $M \vDash \neg(y \,\varepsilon\, x)$,  $M \vDash y \neq v$,

$\Rightarrow$ from  $M \vDash \neg(v \prec y) \;\&\; y \neq v$,  by (8.2),  $M \vDash y \prec v$,  as required.

This completes the proof of (11.5).

## 12. Set concepts and the Strong Set Adjunction Lemma

A formula V(x) is said to be a <u>set concept</u> if

(i)     QT ⊢ V(aa)

(ii)    QT ⊢ ∀x,y,z (Set*(x) & Set*(z) & ∀w (w ε z ↔ w ε x ∨ w=y) &

                                                                       & V(x) & V(y) → V(z)).

(12.1) For any string concept $I \subseteq I_0$ there is a set concept $V \subseteq I$ such that

$$QT^+ \vdash \forall x \in V \, \forall y \, (Set(x) \,\&\, \exists u \, (u \, \varepsilon \, x \,\&\, u \prec y) \to \exists! u \, Max_\preccurlyeq(u,x,y)).$$

Let $I'(x) \equiv I_{8.3} \,\&\, I_{9.3} \,\&\, I_{9.21} \,\&\, I_{11.3}$.

We may assume that $I'$ is downward closed under $\subseteq_p$.

Then let

$$V(x) \equiv I'(x) \,\&\, [Set(x) \,\&\, \forall y \, (\exists u(u \, \varepsilon \, x \,\&\, u \prec y) \to \exists! u \, Max_\preccurlyeq(u,x,y)))].$$

We show that $V(x)$ is a set concept.

That $QT \vdash V(aa)$ follows trivially from (5.18).

Assume that

$$M \vDash Set^*(x) \,\&\, Set^*(z) \,\&\, \forall w \, (w \, \varepsilon \, z \leftrightarrow w \, \varepsilon \, x \lor w = y) \,\&\, V(x) \,\&\, V(y).$$

We want to show that $M \vDash V(z)$.

We may assume that $M \vDash \neg(y \, \varepsilon \, x)$, for otherwise $M \vDash z \sim x$, and by the UNIQUENESS LEMMA $M \vDash z = x$, in which case there is nothing to prove.

Assume now that $M \vDash \exists u(u \, \varepsilon \, x \,\&\, u \prec y')$ for some fixed $y'$.

$\Rightarrow M \vDash u \, \varepsilon \, x \lor u = y.$

(1) $M \vDash u \, \varepsilon \, x$.

$\Rightarrow$ from $M \vDash u \prec y' \,\&\, V(x)$, $M \vDash \exists u' \, Max_\preccurlyeq(u',x,y')$.

Hence

(*) $M \vDash u' \, \varepsilon \, x \,\&\, \forall z'(z' \preccurlyeq_x u' \to z' \prec y') \,\&$

$\& \, \forall z''(z'' \, \varepsilon \, x \,\&\, \forall v(v \preccurlyeq_x z'' \to v \prec y') \to z'' \preccurlyeq_x u'),$

⇒ by (5.18), $M \vDash x \neq aa$,

⇒ from $M \vDash Set(x)$, $M \vDash \exists t\, Env(t,x)$.

(1a) $M \vDash Lastf(x,t,au'a,t)$.

(1ai) $M \vDash y \prec y'$.

⇒ since $M \vDash \neg(y\, \varepsilon\, x)$, $M \vDash y \neq u'$,

⇒ since $M \vDash u' \subseteq_p x\, \&\, V(x)\, \&\, V(y)$ and $M \vDash V \subseteq I' \subseteq I_{8.2}$, $M \vDash I_{8.2}(u')\, \&\, I_{8.2}(y)$,

⇒ $M \vDash y \prec u'\, \vee\, u' \prec y$.

(1aia) $M \vDash y \prec u'$.

We show that $M \vDash Max_{\preccurlyeq}(u',z,y')$.

By (9.21), we have that $M \vDash \exists t_1 Lastf(z,t_1,au'a,t_1)$.

Towards proving $M \vDash Max_{\preccurlyeq}(u',z,y')$, assume $M \vDash z' \leq_z u'$.

⇒ from $M \vDash Lex^+(z)$, $M \vDash z' \preccurlyeq u'$,

⇒ from (*), $M \vDash u' \prec y'$,

⇒ by (8.3), $M \vDash z' \prec y'$.

Thus, $M \vDash \forall z'(z' \leq_z u' \to z' \prec y')$.

Assume that $M \vDash z''\, \varepsilon\, z\, \&\, \forall v\, (v \leq_z z'' \to v \prec y')$.

⇒ from $M \vDash Lastf(z,t_1,au'a,t_1)$ by (9.3), $M \vDash z'' \leq_z u'$.

Hence $M \vDash \forall z''(z''\, \varepsilon\, z\, \&\, \forall v\, (v \leq_z z'' \to v \prec y') \to z'' \leq_z u')$.

This suffices to establish that $M \vDash Max_{\preccurlyeq}(u',z,y')$.

So we have $M \vDash \exists u''\, Max_{\preccurlyeq}(u'',z,y')$.

(1aib) $M \vDash u' \prec y$.

We claim that $M \vDash \exists t_1 Lastf(z,t_1,aya,t_1)$.

Suppose not, i.e., that $M \vDash \forall t_1 \neg \text{Lastf}(z,t_1,aya,t_1)$.

$\Longrightarrow$ from $M \vDash \text{Set}(z)$, $M \vDash \exists t',v\, \text{Lastf}(z,t',ava,t')$,

$\Longrightarrow M \vDash v \neq y$,

$\Longrightarrow$ from $M \vDash v\, \varepsilon\, z$, $M \vDash v\, \varepsilon\, x$,

$\Longrightarrow$ from hypothesis (1a), by (9.3), $M \vDash v \leq_x u'$,

$\Longrightarrow$ from $M \vDash \text{Lex}^+(x)$, $M \vDash v \leqslant u'$,

$\Longrightarrow$ from hypothesis (1aib), by (8.3), $M \vDash v \prec y$.

But from $M \vDash \text{Lastf}(z,t',ava,t')$ & $v \neq y$, $M \vDash y <_z v$,

$\Longrightarrow$ from $M \vDash \text{Lex}^+(z)$, $M \vDash y \prec v$, contradicting $M \vDash v \prec y$ by (8.2).

Therefore, $M \vDash \exists t_1 \text{Lastf}(z,t_1,aya,t_1)$, as claimed.

Assume that $M \vDash z' \leq_z y$.

$\Longrightarrow$ from $M \vDash \text{Lex}^+(z)$, $M \vDash z' \leqslant y$,

$\Longrightarrow$ from hypothesis (1ai) by (8.3), $M \vDash z' \prec y'$.

Hence $M \vDash \forall z'\, (z' \leq_z y \to z' \prec y')$.

Assume that $M \vDash z''\, \varepsilon\, z$ & $\forall v\, (v \leq_z z'' \to v \prec y')$.

$\Longrightarrow$ from $M \vDash \text{Lastf}(z,t_1,aya,t_1)$, by (9.3), $M \vDash z'' \leq_z y$.

Hence $M \vDash \forall z''\, (z''\, \varepsilon\, z\, \&\, \forall v\, (v \leq_z z'' \to v \prec y')\, \to\, z'' \leq_z y)$.

This suffices to establish that $M \vDash \text{Max}_\leqslant(y,z,y')$.

So we have again that $M \vDash \exists u''\, \text{Max}_\leqslant(u'',z,y')$.

(1aii) $M \vDash \neg(y \prec y')$.

$\Longrightarrow$ by (8.2), $M \vDash y' \leqslant y$.

We claim that $M \vDash Max_{\preccurlyeq}(u',z,y')$.

Note first that, from (*), we have $M \vDash u' \prec y'$.

$\implies$ by (8.3), $M \vDash u' \prec y$.

If $M \vDash y \leq_z u'$, then from $M \vDash Lex^+(z)$, $M \vDash y \prec u'$, which contradicts $M \vDash u' \prec y$

by (8.2). Therefore $M \vDash \neg(y \leq_z u')$.

Assume now that $M \vDash z' \leq_z u'$.

$\implies M \vDash z' \neq y$,

$\implies$ since $M \vDash z' \varepsilon z$, $M \vDash z' \varepsilon x$,

$\implies$ from hypothesis (1a) by (9.3), $M \vDash M \vDash z' \leq_x u'$,

$\implies$ by (*), $M \vDash z' \prec y'$.

Therefore, $M \vDash \forall z' (z' \leq_z u' \to z' \prec y')$.

Assume now that $M \vDash z'' \varepsilon z$ & $\forall v (v \leq_z z'' \to v \prec y')$.

Suppose that $M \vDash z''=y$. Then $M \vDash y \prec y'$, which contradicts hypothesis (1a).

Hence $M \vDash z'' \neq y$.

$\implies M \vDash z'' \varepsilon x$,

$\implies$ from hypothesis (1a) by (9.3), $M \vDash z'' \leq_x u'$,

$\implies$ by (9.20), $M \vDash z'' \leq_z u'$.

Therefore, $M \vDash \forall z'' (z'' \varepsilon z$ & $\forall v (v \leq_z z'' \to v \prec y') \to z'' \leq_z u')$.

This establishes that $M \vDash Max_{\preccurlyeq}(u',z,y')$. Thus $M \vDash \exists u'' Max_{\preccurlyeq}(u'',z,y')$.

(1b) $M \vDash \neg Lastf(x,t,au'a,t)$.

$\implies$ from $M \vDash Env(t,x)$, $M \vDash \exists v' Lastf(x,t,av'a,t)$, $M \vDash v' \neq u'$.

(1bi) $M \vDash y \prec y'$.

From $M \vDash u' \varepsilon z$ & $y \varepsilon z$ & $u' \neq y$ by (9.7), $M \vDash u' <_z y \vee y <_z u'$.

(1bia) $M \vDash u' <_z y$.

Two subcases:

(1bia1) $M \vDash \exists v''(u' <_z v'' \;\&\; v'' <_z y)$.

$\Rightarrow$ by (8.2), $M \vDash v'' \neq y$,

$\Rightarrow M \vDash v'' \varepsilon x$, $\Rightarrow$ from $M \vDash \text{Lex}^+(z)$, $M \vDash v'' \prec y$.

Assume now that $M \vDash z' \leq_z u'$.

$\Rightarrow$ from $M \vDash \text{Lex}^+(z)$, $M \vDash z' \prec u'$,

$\Rightarrow$ from hypothesis (1bia) by $M \vDash \text{Lex}^+(z)$, $M \vDash u' \prec y$,

$\Rightarrow$ by (8.3), $M \vDash z' \prec y$.

Thus, $M \vDash \forall z' (z' \leq_z u' \to z' \prec y')$.

Assume now that $M \vDash z'' \varepsilon z$ & $\forall v (v \leq_z z'' \to v \prec y')$.

Suppose, for a reductio, that $M \vDash \neg(z'' \leq_z u')$.

$\Rightarrow$ by (9.7), $M \vDash u' <_z z''$.

Suppose, again for a reductio, that $M \vDash z'' = y$.

$\Rightarrow$ from $M \vDash v'' <_z y$, by (9.16), $M \vDash \forall v (v \leq_z v'' \to v \prec y')$.

Assume that $M \vDash v \leq_x v''$.

$\Rightarrow M \vDash v \varepsilon x$ & $v'' \varepsilon x$,

$\Rightarrow M \vDash v \varepsilon z$ & $v'' \varepsilon z$,

$\Rightarrow$ by (9.20), $M \vDash v \leq_z v''$.

Hence $M \vDash \forall v (v \leq_x v'' \to v \leq_z z'')$, whence further $M \vDash \forall v (v \leq_x v'' \to v \prec y')$.

From $M \vDash u' \varepsilon x$ & $v'' \varepsilon x$ & $u' <_z v''$, by (9.20), $M \vDash u' <_x v''$. But then

M ⊨ v''ε x  &  ∀v (v≤$_x$v'' → v≺y')  &  u'<$_x$v''   contradicts (*) by (9.4) and (9.6).

Therefore  M ⊨ z''≠y.

⟹  M ⊨ z''ε x.

We then derive, analogously to above, that

$$M \vDash z''ε\ x\ \&\ \forall v\ (v\leq_x z'' \to v\prec y')\ \&\ u'<_x z'',$$

again contradicting (*) by (9.4) and (9.6).

So we proved that

$$M \vDash \forall z''(z''\ \varepsilon\ z\ \&\ \forall v\ (v\leq_z z'' \to v\prec y')\ \to\ z''\leq_z u').$$

This suffices to show that  M ⊨ Max$_≼$(u',z,y'), and so  M ⊨ ∃u'' Max$_≼$(u'',z,y').

    (1bia2)  M ⊨ ¬∃v''(u'<$_z$v'' & v''<$_z$y).

We show that  M ⊨ Max$_≼$(y,z,y').

Assume  M ⊨ z'≤$_z$y.

If  M ⊨ z'=y,  then  M ⊨ z'≺y'  by hypothesis (1bi).

If  M ⊨ z'≠y,  then  M ⊨ z'ε x.

Claim:  M ⊨ z'≤$_x$u'.

Suppose otherwise, i.e. that  M ⊨ ¬(z'≤$_x$u').

⟹ by (9.7),  M ⊨ u'<$_z$z',

⟹ from hypothesis (1bia2),  M ⊨ ¬(z'<$_z$y),

⟹ from hypothesis M ⊨ z'≤$_z$y,  M ⊨ z'=y,  contradicting hypothesis.

This proves the claim.

⟹ from (*), M ⊨ z'≺y' .

So we have that   M ⊨ ∀z'(z'≤$_z$y → z'≺y').

Assume now that $M \vDash z'' \varepsilon z$ & $\forall v (v \leq_z z'' \rightarrow v \prec y')$.

$\Rightarrow$ from (9.20), $M \vDash \forall v (v \leq_x z'' \rightarrow v \prec y')$.

$\Rightarrow$ by (*), $M \vDash z'' \leq_x u'$,

$\Rightarrow$ by (9.20), $M \vDash z'' \leq_z u'$,

$\Rightarrow$ by hypothesis (1bia), by (9.16), $M \vDash z'' <_z y$,

$\Rightarrow$ $M \vDash z'' \leq_z y$.

So we proved that $M \vDash \forall z''(z'' \varepsilon z$ & $\forall v (v \leq_z z'' \rightarrow v \prec y') \rightarrow z'' \leq_z y)$,

which suffices for $M \vDash \text{Max}_\preccurlyeq(y,z,y')$. Then again $M \vDash \exists u'' \text{Max}_\preccurlyeq(u'',z,y')$.

   (1bib) $M \vDash y <_z u'$.

We claim that $M \vDash \text{Max}_\preccurlyeq(u',z,y')$.

Assume $M \vDash z' \leq_z u'$.

If $M \vDash z'=y$, then $M \vDash z' \prec y'$ by hypothesis (1bi).

If $M \vDash z' \neq y$, then $M \vDash z' \varepsilon x$.

$\Rightarrow$ since $M \vDash u' \varepsilon x$, by (9.20) from $M \vDash z' \leq_z u'$, $M \vDash z' \leq_x u'$,

$\Rightarrow$ by (*), $M \vDash z' \prec y'$.

So we have $M \vDash \forall z'(z' \leq_z u' \rightarrow z' \prec y')$.

Assume $M \vDash z'' \varepsilon z$ & $\forall v (v \leq_z z'' \rightarrow v \prec y')$.

If $M \vDash z''=y$, then $M \vDash z' \leq_z u'$ by hypothesis (1bib).

If $M \vDash z'' \neq y$, then $M \vDash z'' \varepsilon x$.

$\Rightarrow$ from (9.20), $M \vDash \forall v (v \leq_x z'' \rightarrow v \prec y')$,

$\Rightarrow$ by (*), $M \vDash z'' \leq_x u'$,

$\Rightarrow$ by (9.20), $M \vDash z'' \leq_z u'$.

Thus $M \vDash \forall z''(z'' \varepsilon z \;\&\; \forall v (v \leq_z z'' \to v \prec y') \to z'' \leq_z u')$,

which establishes $M \vDash \text{Max}_\preccurlyeq(u',z,y')$. Hence $M \vDash \exists u'' \text{Max}_\preccurlyeq(u'',z,y')$.

(1bii) $M \vDash \neg(y \prec y')$.

$\implies$ by (8.2), $M \vDash y' \preccurlyeq y$.

We claim that $M \vDash \text{Max}_\preccurlyeq(u',z,y')$.

From (*) we have $M \vDash u' \prec y'$.

Assume that $M \vDash z' \leq_z u'$.

If $M \vDash z'=y$, then from $M \vDash \text{Lex}^+(z)$, $M \vDash y \preccurlyeq u'$.

But then from $M \vDash y' \preccurlyeq y$ by (8.3), $M \vDash y' \preccurlyeq u'$, which contradicts $M \vDash u' \prec y'$ by (8.2).

$\implies M \vDash z' \neq y$,

$\implies$ since $M \vDash z' \varepsilon z$, $M \vDash z' \varepsilon x$,

$\implies$ from $M \vDash z' \leq_z u'$ by (9.20), $M \vDash z' \leq_x u'$,

$\implies$ by (*), $M \vDash z' \prec y'$.

Therefore $M \vDash \forall z'(z' \leq_z u' \to z' \prec y')$.

Assume now that $M \vDash z'' \varepsilon z \;\&\; \forall v (v \leq_z z'' \to v \prec y')$.

Suppose $M \vDash y'=y$.

$\implies M \vDash z'' \prec y'=y$,

$\implies$ by (8.2), $M \vDash z'' \neq y$,

$\implies M \vDash z'' \varepsilon x$,

$\implies$ from (9.20), $M \vDash \forall v (v \leq_x z'' \to v \prec y')$,

$\implies$ from (*), $M \vDash z'' \leq_x u'$,

$\Rightarrow$ by (9.20), $M \vDash z'' \leq_z u'$.

Suppose $M \vDash y' \neq y$.

$\Rightarrow M \vDash y' \prec y$,

$\Rightarrow$ from $M \vDash z'' \prec y'$ by (8.3), $M \vDash z'' \prec y$,

$\Rightarrow$ by (8.2), $M \vDash z' \neq y$,

$\Rightarrow M \vDash z'' \varepsilon x$,

$\Rightarrow$ just as above, $M \vDash z'' \leq_z u'$.

Therefore, $M \vDash \forall z''(z'' \varepsilon z \ \& \ \forall v (v \leq_z z'' \rightarrow v \prec y') \rightarrow z'' \leq_z u')$,

which suffices for $M \vDash \text{Max}_\leq(u',z,y')$. Thus again $M \vDash \exists u'' \text{Max}_\leq(u'',z,y')$.

(2) $M \vDash u = y$.

Recall that $M \vDash u \varepsilon z \ \& \ u \prec y'$.

$\Rightarrow M \vDash y \prec y'$

We then proceed exactly as in (1ai) and (1bi).

This completes the proof that $M \vDash \exists u'' \text{Max}_\leq(u'',z,y')$.

To establish uniqueness, suppose that $M \vDash \text{Max}_\leq(u_1,z,y') \ \& \ \text{Max}_\leq(u_2,z,y')$.

$\Rightarrow M \vDash u_1 \varepsilon z \ \& \ u_2 \varepsilon z \ \& \ u_1 \prec y' \ \& \ u_2 \prec y'$,

$\Rightarrow$ by (9.7), $M \vDash u_1 <_z u_2 \ \vee \ u_2 <_z u_1 \ \vee \ u_1 = u_2$.

Suppose, for a reductio, that $M \vDash u_1 <_z u_2$.

$\Rightarrow$ from $M \vDash \text{Max}_\leq(u_2,z,y')$, $M \vDash \forall v (v \leq_z u_2 \rightarrow v \prec y')$,

$\Rightarrow$ from $M \vDash \text{Max}_\leq(u_1,z,y')$, $M \vDash u_2 \leq_z u_1$.

But this contradicts $M \vDash u_1 <_z u_2$ by (9.4) and (9.6).

Therefore, $M \vDash \neg(u_1 <_z u_2)$.

Likewise, $M \vDash \neg(u_2 \leq_z u_1)$.

Hence, $M \vDash u_1 = u_2$, as required.

This completes the proof that $M \vDash \forall y'(\exists u(u \: \varepsilon \: z \: \& \: u \prec y) \rightarrow \exists! u \: \text{Max}_{\preccurlyeq}(u,z,y'))$, and hence $M \vDash V(z)$.

Thus $V(x)$ is a set concept.

This completes the proof of (12.1).

We let $I^{**}(x)$ abbreviate the formula

$(I_{9.15}$ & $I_{9.19}$ & $I_{9.26}$ & $I_{10.11}$ & $I_{10.13}$ & $I_{10.19}$ & $I_{10.20}$ & $I_{10.23}$ & $I_{10.24}$ & $I_{10.25}$ & $I_{10.28}$ &

& $I_{10.29}$ & $I_{11.4})_{SUB}$.

Let $V^{**}(x) \equiv Set^*(x)$ & $I^{**}(x)$ & $V^*(x)$ & $[\forall t\ (Env(t,x) \to J^{**}(t))]$,

where $V^*$ is the set concept obtained from $I \equiv I_0$ as in (12.1), and $J^{**}$ is the string concept obtained from $I^{**}$ as $J$ is from $I^*$ in (10.29).

STRONG SET ADJUNCTION LEMMA. (12.2)

$QT^+ \vdash \forall x,y\ (V^{**}(x)$ & $V^{**}(y) \to \exists!z\ (V^{**}(z)$ & $\forall w\ (w\ \varepsilon\ z \leftrightarrow (w\ \varepsilon\ x \lor w=y)))$.

Assume $M \vDash V^{**}(x)$ & $V^{**}(y)$.

Note that, since $V^*$ is a set concept, as proved in (12.1), once we show the existence of z such that $M \vDash Set^*(z)$ & $\forall w\ (w\ \varepsilon\ z \leftrightarrow (w\ \varepsilon\ x \lor w=y))$, that $M \vDash V^*(z)$ follows at once from the hypothesis $M \vDash V^{**}(x)$ & $V^{**}(y)$. The uniqueness of z in $I^{**}$, and hence in $V^{**}$, follows immediately from the UNIQUENESS LEMMA (11.4).

From $M \vDash V^{**}(x)$ & $V^{**}(y)$ we have $M \vDash Set^*(x)$ & $Set^*(y)$.

$\implies M \vDash MinSet(x)$ & $Lex^+(x)$ & $Special(x)$ & $I^{**}(x)$ & $I^{**}(y)$,

$\implies M \vDash Set(x)$,

$\implies$ since $I^{**} \subseteq I_{6.8}$, by (6.8), $M \vDash \exists!t_0 \in I^{**}\ MinMax^+T_b(t_0,y)$.

From $M \vDash Set^*(y)$ we have that $M \vDash Set(y)$, whence $M \vDash y=aa \lor \exists t_y Env(t_y,y)$.

If $M \vDash y=aa$, then $M \vDash Tally_a(y)$, and from (6.9) we have that

$M \vDash \text{MinMax}^+T_b(b,y)$.

$\Rightarrow$ as in the proof of (6.8), from $M \vDash \text{MinMax}^+T_b(t_0,y)$ we have that $M \vDash t_0 = b$, whence $M \vDash J^{**}(t_0)$.

Suppose, on the other hand, that $M \vDash \exists t_y \text{Env}(t_y, y)$.

$\Rightarrow$ from $M \vDash V^{**}(y)$, $M \vDash J^{**}(t_y)$,

$\Rightarrow$ $M \vDash \exists v \text{ Lastf}(y, t_y, \text{ava}, t_y)$,

$\Rightarrow$ from the proof of (6.8), $M \vDash t_0 = t_y b$,

$\Rightarrow$ $M \vDash t_0 \in J^{**} \subseteq I^{**}$.

From $M \vDash \text{Set}^*(x)$ we have $M \vDash \text{Set}(x)$, whence $M \vDash x = \text{aa} \vee \exists t \text{Env}(t, x)$.

Suppose

(*A): $M \vDash x = \text{aa}$.

Let $z = t_0 \text{ay} a t_0$.

Then, since we may assume that $I^{**}$ is closed with respect to *, and $M \vDash I^{**}(t_0) \mathbin{\&} I^{**}(y)$, we have $M \vDash I^{**}(z)$.

$\Rightarrow$ by (10.15), $M \vDash \text{Env}(t_0, z) \mathbin{\&} \text{Set}^*(z) \mathbin{\&} y \, \varepsilon \, z \mathbin{\&} \forall w (w \, \varepsilon \, z \leftrightarrow w = y)$.

$\Rightarrow$ by (5.18), $M \vDash \forall w \neg (w \, \varepsilon \, x)$,

$\Rightarrow$ $M \vDash \forall w (w \, \varepsilon \, z \leftrightarrow w \, \varepsilon \, x \vee w = y)$, whereas also $M \vDash V^{**}(z)$.

(*B): $M \vDash x \neq \text{aa}$.

$\Rightarrow$ $M \vDash \exists t \text{Env}(t, x)$.

If $M \vDash y \, \varepsilon \, x$, we may let $z = x$. So we may assume that $M \vDash \neg (y \, \varepsilon \, x)$.

$\Rightarrow$ from $M \vDash I^{**}(x) \mathbin{\&} I^{**}(y)$, since we may assume that $I^{**}$ is downward

closed with respect to $\subseteq_p$, and $I^{**} \subseteq I_{8.2}$, $M \vDash \forall u(u \, \varepsilon \, x \to u \prec y \lor y \prec u)$.

(*Bi): $M \vDash \forall u(u \, \varepsilon \, x \to u \prec y)$.

(*Bia): $M \vDash t \prec t_0$.

$\Rightarrow M \vDash \exists t^* \, tt^* = t_0$.

Let $z = xt^*ayatt^*$.

$\Rightarrow$ from $M \vDash I^{**}(x) \,\&\, I^{**}(y) \,\&\, I^{**}(t_0)$, since we may assume that $I^{**}$ is closed with respect to *, $M \vDash I^{**}(z)$.

Let $x^+ = tt^*ayatt^*$.

Then $M \vDash I^{**}(x^+)$.

$\Rightarrow$ since $M \vDash Max^+T_b(tt^*, y)$,

$\qquad M \vDash Pref(aya, tt^*) \,\&\, Firstf(x^+, tt^*, aya, tt^*) \,\&\, Lastf(x^+, tt^*, aya, tt^*)$,

$\Rightarrow$ by (5.22), $M \vDash Env(tt^*, x^+) \,\&\, \forall w(w \, \varepsilon \, x^+ \leftrightarrow w = y)$,

$\Rightarrow$ from $M \vDash Env(t, x)$, by (5.11),

$\qquad M \vDash \exists x_1, t_1(Tally_b(t_1) \,\&\, x = t_1 x_1 t \,\&\, abx_1 \,\&\, aEx_1)$,

$\Rightarrow M \vDash Env(t, x) \,\&\, x = t_1 x_1 t \,\&\, abx_1 \,\&\, aEx_1 \,\&\, Env(tt^*, x^+) \,\&\, t \prec tt^* \,\&$

$\qquad\qquad\qquad \&\, Firstf(x^+, tt^*, aya, tt^*) \,\&\, \neg\exists w(w \, \varepsilon \, x \,\&\, w \, \varepsilon \, x^+)$,

$\Rightarrow$ by (5.46), for $z' = t_1 x_1 x^+$, $M \vDash Env(tt^*, z') \,\&\, \forall w(w \, \varepsilon \, z' \leftrightarrow w \, \varepsilon \, x \lor w \, \varepsilon \, x^+)$,

$\Rightarrow M \vDash z = xt^*ayatt^* = (t_1 x_1 t)t^*ayatt^* = t_1 x_1 x^+ = z'$,

$\Rightarrow M \vDash Env(tt^*, z) \,\&\, \forall w(w \, \varepsilon \, z \leftrightarrow w \, \varepsilon \, x \lor w = y)$.

$\Rightarrow$ from $M \vDash Lastf(x^+, tt^*, aya, tt^*)$, by the proof of (5.46),(1c),

$\qquad M \vDash Lastf(z, tt^*, aya, tt^*)$,

$\Rightarrow M \vDash Fr(z, tt^*, aya, tt^*)$,

$\Rightarrow$ $M \vDash y \, \varepsilon \, z$,

$\Rightarrow$ from hypothesis $M \vDash Set^*(x)$, $M \vDash MinSet(x)$,

$\Rightarrow$ $M \vDash Pref(aya, tt^*)$ by (10.5), $M \vDash MinSet(x^+)$,

$\Rightarrow$ $M \vDash Env(t,x)$ & $x = t_1 x_1 t$ & $aBx_1$ & $aEx_1$ & $z = t_1 x_1 x^+$ & $Env(tt^*, x^+)$ & $t < tt^*$ &

& $Firstf(x^+, tt^*, aya, tt^*)$ & $\neg \exists w (w \, \varepsilon \, x \, \& \, w \, \varepsilon \, x^+)$ & $MinSet(x)$ & $MinSet(x^+)$,

$\Rightarrow$ by (10.6), $M \vDash MinSet(z)$.

That $M \vDash Lex^+(x^+)$ follows immediately from (9.4).

$\Rightarrow$ $M \vDash Env(tt^*, z)$ & $z = t_1 x_1 x^+$ & $Env(t,x)$ & $x = t_1 x_1 t$ & $aBx_1$ & $aEx_1$ &

& $Env(tt^*, x^+)$ & $t < tt^*$ & $Firstf(x^+, tt^*, aya, tt^*)$ & $\neg \exists w (w \, \varepsilon \, x \, \& \, w \, \varepsilon \, x^+)$ &

& $\forall w_1, w_2 (w_1 \, \varepsilon \, x \, \& \, w_2 \, \varepsilon \, x^+ \to w_1 \prec w_2)$ & $Lex^+(x)$ & $Lex^+(x^+)$,

$\Rightarrow$ by (9.29), $M \vDash Lex^+(z)$.

Finally, from $M \vDash Special(x)$, by (10.15) and the choice of $t_0$, $M \vDash Special(x^+)$.

$\Rightarrow$ by (10.17), $M \vDash Special(z)$.

Hence $M \vDash Set^*(z)$, as required.

Since also $M \vDash Env(t_0, z) \, \& \, J^{**}(t_0)$, it follows that $M \vDash V^{**}(z)$.

(*Bib): $M \vDash t_0 \leq t$.

Let $z = xbayatb$.

$\Rightarrow$ from $M \vDash I^{**}(x) \, \& \, I^{**}(y)$, $M \vDash I^{**}(z)$.

We proceed as in (*Bia). Then with $t^* = b$, the proof of $M \vDash MinSet(z)$ and

$M \vDash Lex^+(z)$ remains the same as in (*Bia).

We likewise prove that $M \vDash \forall w (w \, \varepsilon \, z \leftrightarrow w \, \varepsilon \, x \, w = y)$.

For $M \vDash \text{Special}(z)$, let $x^* = tbayatb$.

Then $M \vDash I^{**}(x^*)$.

We have, from the choice of $t_0$, that $M \vDash \text{Max}^+T_b(tb,y)$.

$\implies M \vDash \text{Pref}(aya, tb)$,

$\implies M \vDash \text{Firstf}(x^*, tb, aya, tb) \, \& \, \text{Lastf}(x^*, tb, aya, tb)$,

$\implies$ by (5.22), $M \vDash \text{Env}(tb, x^*) \, \& \, \forall w(w \, \varepsilon \, x^* \leftrightarrow w = y)$,

$\implies$ from $M \vDash \text{Env}(t,x)$, by (5.11), $M \vDash \exists t_0, z'(\text{Tally}_b(t_0) \, \& \, x = t_0 x't \, \& \, aBx' \, \& \, aEx')$,

$\implies M \vDash \text{Env}(tb, z) \, \& \, z = t_0 x' x^* \, \& \, \text{Env}(t,x) \, \& \, x = t_0 x't \, \& \, aBx' \, \& \, aEx' \, \& \, x^* = tbayatb \, \&$

$\& \, \text{Max}^+T_b(tb,y) \, \& \, \neg(y \, \varepsilon \, x) \, \& \, \forall u(u \, \varepsilon \, x \to u \prec y) \, \& \, \text{Special}(x)$,

$\implies$ by (10.18), $M \vDash \text{Special}(z)$,

$\implies M \vDash \text{Set}^*(z)$.

Since from $M \vDash \text{Env}(t,x) \, \& \, V^{**}(x)$, $M \vDash J^{**}(t)$, we also have that $M \vDash J^{**}(tb)$.

It follows that $M \vDash V^{**}(z)$.

(*Bii): $M \vDash \neg \forall u(u \, \varepsilon \, x \to u \prec y)$.

$\implies M \vDash \exists u(u \, \varepsilon \, x \, \& \, \neg u \prec y)$,

$\implies$ since $M \vDash \neg(y \, \varepsilon \, x)$, by (8.2), $M \vDash \exists u(u \, \varepsilon \, x \, \& \, y \prec u)$.

(*Biia): $M \vDash \forall u(u \, \varepsilon \, x \to y \prec u)$.

$\implies$ from $M \vDash \text{Env}(t,x)$, $M \vDash \exists v_0, t_1, t_2 \text{Firstf}(x, t_1, av_0 a, t_2)$,

$\implies M \vDash y \prec v_0$,

$\implies$ by choice of $t_0$, $M \vDash \text{MinMax}^+T_b(t_0, y)$,

$\implies$ from $M \vDash I^{**}(x) \, \& \, v_0 \subseteq_p x$, $M \vDash I^{**}(v_0)$,

$\Rightarrow$ by (6.8), $M \vDash \exists t^* \text{MinMax}^+T_b(t^*,v_0)$,

$\Rightarrow$ from $M \vDash y \prec v_0$, $M \vDash y \triangleleft_{T_b} v_0 \lor (y \approx_{T_b} v_0 \,\&\, y \lll v_0)$.

(*Biia1): $M \vDash y \triangleleft_{T_b} v_0$.

$\Rightarrow M \vDash t_0 < t^*$.

Let $z = t_0 a y a x$.

$\Rightarrow$ from $M \vDash I^{**}(t_0) \,\&\, I^{**}(y) \,\&\, I^{**}(x)$, $M \vDash I^{**}(z)$,

$\Rightarrow$ from $M \vDash \text{Firstf}(x,t_1,av_0a,t_2) \,\&\, \text{Special}(x)$, by (10.16),

$$M \vDash \text{MinMax}^+T_b(t_1,v_0),$$

$\Rightarrow$ by (6.8), $M \vDash t_1 = t^*$,

$\Rightarrow$ from $M \vDash \text{Max}^+T_b(t_0,y)$, $M \vDash \text{Pref}(aya,t_0)$,

$\Rightarrow$ just as in (*Bia), $M \vDash \text{Env}(t_0,x^+) \,\&\, \forall w(w \,\varepsilon\, x^+ \leftrightarrow w = y)$ where $x^+ = t_0 a y a t_0$,

$\Rightarrow M \vDash \text{Env}(t_0,x^+) \,\&\, x^+ = t_0 a y a t_0 \,\&\, \text{Env}(t,x) \,\&\, \text{Firstf}(x,t_1,av_0a,t_2) \,\&\, t_0 < t_1 \,\&$

$\&\, \neg \exists w(w \,\varepsilon\, x^+ \,\&\, w \,\varepsilon\, x)$,

$\Rightarrow$ by (5.46), $M \vDash \text{Env}(t,z) \,\&\, \forall w(w \,\varepsilon\, z \leftrightarrow w \,\varepsilon\, x^+ \lor w \,\varepsilon\, x)$,

$\Rightarrow M \vDash \forall w(w \,\varepsilon\, z \leftrightarrow w = y \lor w \,\varepsilon\, x)$.

Now, with $M \vDash (t_1 a) B x \,\&\, z = t_0 a y a x$, we have

$M \vDash \text{Pref}(aya,t_0) \,\&\, t_0 < t_1 \,\&\, \text{Tally}_b(t_1) \,\&\, (t_0 a y a t_1 a) B z$.

$\Rightarrow \text{Firstf}(z,t_0,av_0a,t_1)$,

$\Rightarrow M \vDash y \,\varepsilon\, z$,

$\Rightarrow$ by (10.5) and choice of $x^+$, $M \vDash \text{MinSet}(x^+)$,

$\Rightarrow M \vDash \text{Env}(t_0,x^+) \,\&\, x^+ = t_0 a y a t_0 \,\&\, z = t_0 a y a x \,\&\, \text{Env}(t,x) \,\&\, \text{Firstf}(x,t_1,av_0a,t_2) \,\&$

    & $t_0 < t_1$ & $\neg\exists w(w\ \varepsilon\ x^+\ \&\ w\ \varepsilon\ x)$ & $\text{MinSet}(x^+)$ & $\text{MinSet}(x)$,

$\Rightarrow$ by (10.6),  $M \models \text{MinSet}(z)$.

Next, from hypothesis  $M \models \text{Set}^*(x)$  we have  $M \models \text{Lex}^+(x)$,  and  $M \models \text{Lex}^+(x^+)$ follows from (9.1).

We then have, taking into account hypothesis (*Biia), that

 $M \models \text{Env}(t,z)$ & $z = t_0 a y a x$ & $\text{Env}(t_0, x^+)$ & $x^+ = t_0 a y a t_0$ & $\text{Env}(t,x)$ &

   & $\text{Firstf}(x, t_1, a v_0 a, t_2)$ & $t_0 < t_1$ & $\neg\exists w(w\ \varepsilon\ x^+\ \&\ w\ \varepsilon\ x)$ &

  & $\forall w_1, w_2(w_1\ \varepsilon\ x\ \&\ w_2\ \varepsilon\ x^+ \rightarrow w_1 \prec w_2)$ & $\text{Lex}^+(x)$ & $\text{Lex}^+(x^+)$,

$\Rightarrow$ by (9.29),  $M \models \text{Lex}^+(z)$.

Finally, by (10.15) and the choice of $t_0$, $M \models \text{Special}(x^+)$.

$\Rightarrow$ from hypothesis  $M \models \text{Special}(x)$, by (10.17),  $M \models \text{Special}(z)$.

Hence  $M \models \text{Set}^*(z)$.

From  $M \models \text{Env}(t,x)$ & $V^{**}(x)$,  we have  $M \models J^{**}(t)$. Taking into account $M \models \text{Env}(t,z)$, it follows that  $M \models V^{**}(z)$.

 (*Biia2):  $M \models y \approx_{Tb} v_0$  &  $y \ll v_0$.

$\Rightarrow M \models t_0 = t^*$,

$\Rightarrow$ from hypothesis  $M \models \text{Special}(x)$, $M \models \text{MinSet}(x)$,

$\Rightarrow$ since  $M \models I^{**}(x)$ & $J^{**}(t)$,  by (10.29),

 $M \models \exists! x^* \in I^{**} \exists t' \in J^{**}(\text{Env}(t', x^*)\ \&\ x \sim x^*\ \&$

   & $\text{MinSet}(x^*)$ & $\forall u, v(u <_x v \leftrightarrow u <_{x^*} v)$ &

    & $\forall v,t_1,t_2$ (Free$^+$(x,t$_1$,ava,t$_2$) → $\exists t_3$Fr(x*,t$_1$b,ava,t$_3$)) &

& $\forall v,t_1,t_2$ (Bound(x,t$_1$,ava,t$_2$) v Free$^-$(x,t$_1$,ava,t$_2$) → $\exists t_3$Fr(x*,t$_1$,ava,t$_3$))).

Again from M ⊨ Firstf(x,t$_1$,av$_0$a,t$_2$) & Special(x), by (10.16), we have

$$M \vDash \text{MinMax}^+T_b(t_1,v_0).$$

⟹ M ⊨ t$_1$=t*,

⟹ M ⊨ Firstf(x,t*,av$_0$a,t$_2$)

⟹ M ⊨ Free$^+$(x,t*,av$_0$a,t$_2$),

⟹ M ⊨ $\exists t_3$Fr(x*,t*b,av$_0$a,t$_3$),

⟹ by (9.12), M ⊨ $\exists t_4,t_5$ Firstf(x*,t$_4$,av$_0$a,t$_5$),

⟹ by (5.15), M ⊨ t$_4$=t*b & t$_5$=t$_3$,

⟹ M ⊨ Firstf(x*,t*b,av$_0$a,t$_3$),

⟹ M ⊨ Pref(av$_0$a,t*b) & Tally$_b$(t$_3$) & ((t*b=t$_3$ & x*=t*bav$_0$at$_3$) v

                                                v (t*b<t$_3$ & (t*bav$_0$at$_3$a)Bx*)).

    (*Biia2i): M ⊨ t*b=t$_3$ & x*=t*bav$_0$at$_3$.

Let z=t$_0$ayat*bav$_0$at*b.

⟹ from M ⊨ I**(t$_0$) & I**(y) & I**(t*) & I**(v$_0$), since I** may be assumed to be closed w.r. to *, M ⊨ I**(z),

⟹ from M ⊨ I**(y) & I**(v$_0$), from the proof of (10.19)(2),

        M ⊨ Set*(z) & $\forall w$(w ε z ↔ w=y v w=v$_0$),

⟹ M ⊨ Pref(av$_0$a,t*b), M ⊨ Lastf(x*,t*b,av$_0$a,t*b),

⟹ from M ⊨ Env(t',x*), M ⊨ t'=t*b,

$\Rightarrow$ $M \vDash J^{**}(t^*b)$.

Hence $M \vDash V^{**}(z)$.

Now, from $M \vDash \text{Firstf}(x^*,t^*b,av_0a,t^*b)$ & $\text{Lastf}(x^*,t^*b,av_0a,t^*b)$.

$\Rightarrow$ from $M \vDash I^{**}(x^*)$, by (5.22), $M \vDash \forall w(w \, \varepsilon \, x^* \leftrightarrow w=v_0)$,

$\Rightarrow$ $M \vDash \forall w(w \, \varepsilon \, z \leftrightarrow w=y \vee w \, \varepsilon \, x^*)$,

$\Rightarrow$ from $M \vDash x \sim x^*$, $M \vDash \forall w(w \, \varepsilon \, z \leftrightarrow w=y \vee w \, \varepsilon \, x^*)$, as required.

  (*Biia2ii): $M \vDash t^*b < t_3$ & $(t^*bav_0at_3a)Bx^*$.

$\Rightarrow$ $M \vDash \exists x^{**} \, t^*bav_0at_3ax^{**}=x^*$.

Let $z = t_0ayax^*$.

$\Rightarrow$ from $M \vDash I^{**}(t_0)$ & $I^{**}(y)$ & $I^{**}(x^*)$, since $I^{**}$ may be assumed to be closed w.r. to *,    $M \vDash I^{**}(z)$,

$\Rightarrow$ from $M \vDash t_0=t^*$, $M \vDash t_0 < t_0b = t^*b$.

We have, as in (*Biia1), that $M \vDash \text{MinSet}(x^+)$ & $\forall w(w \, \varepsilon \, x^+ \leftrightarrow w=y)$

where $x^+ = t_0ayat_0$.

$\Rightarrow$ $M \vDash \text{Env}(t_0, x^+)$ & $z=t_0ayax^*$ & $\text{Env}(t',x^*)$ & $\text{Firstf}(x^*,t^*b,av_0a,t_3)$ & $t_0 < t^*b$ &

& $\neg \exists w(w \, \varepsilon \, x^+ \& w \, \varepsilon \, x^*)$ & $\text{MinSet}(x^+)$ & $\text{MinSet}(x^*)$,

$\Rightarrow$ from $M \vDash I^{**}(z)$, by (10.6), $M \vDash \text{MinSet}(z)$ & $\text{Env}(t',z)$.

Note that $M \vDash J^{**}(t')$.

We now argue that $M \vDash \text{Lex}^+(z)$.

We have that $M \vDash \text{Lex}^+(x^+)$ as in (*Biia1).

$\Rightarrow$ from

$M \vDash I^{**}(x)$ & $I^{**}(x^*)$ & $\text{Set}(x)$ & $\text{Set}(x^*)$ & $\text{Lex}^+(x)$ & $\forall u,v(u <_x v \leftrightarrow u <_{x^*} v)$,

by (9.25), $\quad M \vDash \text{Lex}^+(x^*)$.

Taking into account the hypothesis (*Biia), we then have

$\Rightarrow M \vDash \text{Env}(t',z)\ \&\ z=t_0 a y a x^*\ \&\ \text{Env}(t_0, x^+)\ \&\ \text{Firstf}(x^*, t^*b, a v_0 a, t_3)\ \&\ t_0 < t^* b\ \&$

$\&\ \neg \exists w(w\ \varepsilon\ x^+\ \&\ w\ \varepsilon\ x^*)\ \&\ \forall w_1, w_2(w_1\ \varepsilon\ x\ \&\ w_2\ \varepsilon\ x^+ \rightarrow w_1 \prec w_2)\ \&$

$\&\ \text{Lex}^+(x)\ \&\ \text{Lex}^+(x^+),$

$\Rightarrow$ from $M \vDash I^{**}(z)$, by (9.29), $M \vDash \text{Lex}^+(z)$.

Finally, we show that $M \vDash \text{Special}(z)$.

From $M \vDash I^{**}(y)$, by (10.15), $M \vDash \text{Special}(x^+)$.

$\Rightarrow M \vDash \text{Env}(t',z)\ \&\ z=t_0 a y a x^*\ \&\ \text{Set}(x^*)\ \&\ x^+=t_0 a y a t_0\ \&\ t_0=t^*\ \&\ \text{Env}(t_0, x^+)\ \&$

$\&\ \text{Firstf}(x^*, t^* b, a v_0 a, t_3)\ \&\ \forall u,v(u\ \varepsilon\ x^+\ \&\ v\ \varepsilon\ x^* \rightarrow u \prec v)\ \&\ \neg \exists u(u\ \varepsilon\ x^+\ \&\ u\ \varepsilon\ x^*)\ \&$

$\&\ \text{Set}(x)\ \&\ x \sim x^*\ \&\ \forall u,v(u <_x v \leftrightarrow u <_{x^*} v)\ \&$

$\&\ \forall v, t_1, t_2\ (\text{Free}^+(x, t_1, a v a, t_2) \rightarrow \exists t_3 \text{Fr}(x^*, t_1 b, a v a, t_3))\ \&$

$\&\ \forall v, t_1, t_2\ (\text{Bound}(x, t_1, a v a, t_2)\ \lor\ \text{Free}^-(x, t_1, a v a, t_2) \rightarrow \exists t_3 \text{Fr}(x^*, t_1, a v a, t_3))\ \&$

$\&\ \text{Special}(x^+)\ \&\ \text{Special}(x),$

$\Rightarrow$ from $M \vDash I^{**}(x)\ \&\ I^{**}(z)$, by (10.28), $M \vDash \text{Special}(z)$.

Hence $M \vDash \text{Set}^*(z)$, and indeed $M \vDash V^{**}(z)$.

It only remains to show that $M \vDash \forall u(u\ \varepsilon\ z \leftrightarrow u\ \varepsilon\ x \lor u=y)$.

Now, we have that

$\Rightarrow M \vDash \text{Env}(t_0, x^+)\ \&\ x^+=t_0 a y a t_0\ \&\ \text{Env}(t', x^*)\ \&\ \text{Firstf}(x^*, t^* b, a v_0 a, t_3)\ \&\ t_0 < t^* b\ \&$

$\&\ \neg \exists u(u\ \varepsilon\ x^+\ \&\ u\ \varepsilon\ x^*),$

$\Rightarrow$ from $M \vDash I^{**}(x^+)\ \&\ I^{**}(x^*)$, by (5.46), $M \vDash \forall u(u\ \varepsilon\ z \leftrightarrow u\ \varepsilon\ x^+\ \&\ u\ \varepsilon\ x^*)$,

$\Rightarrow$ since $M \vDash x \sim x^*$, $M \vDash \forall u(u \, \varepsilon \, z \leftrightarrow u \, \varepsilon \, x \vee u = y)$, as required.

(*Biib): $M \vDash \neg \forall u(u \, \varepsilon \, x \rightarrow y \prec u)$.

$\Rightarrow M \vDash \exists u(u \, \varepsilon \, x \,\&\, \neg(y \prec u))$,

$\Rightarrow$ since $M \vDash \neg(y \, \varepsilon \, x)$, by (8.2), $M \vDash \exists u(u \, \varepsilon \, x \,\&\, u \prec y)$,

$\Rightarrow$ since $M \vDash V^*(x)$, by (12.1), $M \vDash \exists ! u_0 \, \mathrm{Max}_{\preccurlyeq}(u_0, x, y)$,

$\Rightarrow M \vDash u_0 \, \varepsilon \, x \,\&\, \forall w(w \preccurlyeq_x u_0 \rightarrow w \prec y) \,\&\, \forall w(w \, \varepsilon \, x \,\&\, \forall v(v \preccurlyeq_x w \rightarrow v \prec y) \rightarrow w \preccurlyeq_x u_0))$,

$\Rightarrow M \vDash \exists t_1, t_2 \mathrm{Fr}(x, t_1, au_0a, t_2)$.

We now distinguish three principal subcases:

(*Biib1): $M \vDash \mathrm{Firstf}(x, t_1, au_0a, t_2)$,

(*Biib2): $M \vDash \exists w_1 \mathrm{Intf}(x, w_1, t_1, au_0a, t_2)$,

(*Biib3): $M \vDash \mathrm{Lastf}(x, t_1, au_0a, t_2)$.

We first deal with (*Biib2).

(*Biib2): $M \vDash \exists w_1 \mathrm{Intf}(x, w_1, t_1, au_0a, t_2)$.

$\Rightarrow$ from hypothesis $M \vDash \mathrm{Set}^*(x)$, $M \vDash \mathrm{MinSet}(x)$,

$\Rightarrow$ from $M \vDash \mathrm{Env}(t, x)$, by the RESOLUTION LEMMA (10.9),

$M \vDash \exists x', x'', t^*, t', w_0^*, w''(\mathrm{Env}(t', x') \,\&\, \mathrm{Env}(t, x'') \,\&\, x' = t_0 w_0^* t' \,\&$

$\&\, x'' = t_1 au_0 a t_2 a w'' \,\&\, x = t_0 w_0^* x'' \,\&\, \mathrm{Firstf}(x'', t_1, au_0a, t_2) \,\&$

$\&\, t' \prec t_1 \,\&\, aBw_0^* \,\&\, aEw_0^* \,\&\, \neg \exists w(w \, \varepsilon \, x' \,\&\, w \, \varepsilon \, x''))]$,

$\Rightarrow$ from $M \vDash I^{**}(x)$, since $I^{**}$ may be assumed to be downward closed w.r. to $\subseteq_p$, $M \vDash I^{**}(w_0^*) \,\&\, I^{**}(x'') \,\&\, I^{**}(t_1) \,\&\, I^{**}(u_0) \,\&\, I^{**}(t_2) \,\&\, I^{**}(w_2)$,

$\Rightarrow$ from $M \vDash t' < t_1$, $M \vDash I^{**}(t')$,

$\Rightarrow$ since $I^{**}$ may be assumed to be closed w.r. to $*$, from $M \vDash I^{**}(t_0)$ & $I^{**}(w_0*)$,

$$M \vDash I^{**}(x'),$$

$\Rightarrow$ by (5.46), $M \vDash \forall w(w \, \varepsilon \, x \leftrightarrow w \, \varepsilon \, x' \vee w \, \varepsilon \, x'')$,

$\Rightarrow$ by (10.10), $M \vDash MinSet(x'')$.

Note that

$$M \vDash w_1 a t_1 a u_0 a t_2 a w_2 = w_1 a x'' = t * w * x''.$$

$\Rightarrow$ by (3.6), $M \vDash w_1 a = t * w *$.

By the SUBTRACTION LEMMA, for $x^- = t_2 a w_2$, we have

$M \vDash Env(t, x^-)$ & $\forall w(w \, \varepsilon \, x^- \leftrightarrow w \, \varepsilon \, x''$ & $w \neq u_0)$ & $(Lex^+(x'') \rightarrow Lex^+(x^-))$.

$\Rightarrow$ since $I^{**}$ may be assumed to be closed w.r. to $*$, from $M \vDash I^{**}(t_2)$ & $I^{**}(w_2)$,

$$M \vDash I^{**}(x^-),$$

$\Rightarrow$ $M \vDash \exists v_0, t_2^+, t_3 \, Firstf(x^-, t_2^+, a v_0 a, t_3)$,

$\Rightarrow$ $M \vDash (t_2^+ a) B x^-$ & $Tally_b(t_2^+)$,

$\Rightarrow$ from $M \vDash I^{**}(x^-)$, $M \vDash I^{**}(t_2^+)$ & $I^{**}(v_0)$ & $I^{**}(w_2)$,

$\Rightarrow$ $M \vDash \exists x_1 \, t_2^+ a x_1 = x^- = t_2 a w_2$,

$\Rightarrow$ from $M \vDash I^{**}(t_2^+)$ & $I^{**}(t_2)$, by (4.23$^b$), $M \vDash t_2^+ = t_2$,

$\Rightarrow$ $M \vDash Firstf(x^-, t_2, a v_0 a, t_3)$,

$\Rightarrow$ as in the proof of the SUBTRACTION LEMMA, $M \vDash Fr(x'', t_2, a v_0 a, t_3)$,

$\Rightarrow$ as in the proof of (5.46), $M \vDash Fr(x, t_2, a v_0 a, t_3)$,

$\Rightarrow$ from $M \vDash Intf(x, w_1, t_1, a u_0 a, t_2)$ & $t_1 < t_2$, by (9.14), $M \vDash u_0 <_x v_0$.

By choice of $u_0$, we have $M \vDash Max_\leqslant(u_0, x, y)$.

Also, $M \vDash \neg(y \, \varepsilon \, x)$ & $Lex^+(x)$.

$\Rightarrow$ from $M \vDash I^{**}(x)$ & $I^{**}(y)$, by (11.5), $M \vDash y \prec v_0$,

$\Rightarrow M \vDash u_0 \prec y \prec v_0$,

$\Rightarrow$ by choice of $t_0$, $M \vDash MinMax^+T_b(t_0,y)$,

$\Rightarrow$ from $M \vDash I^{**}(u_0)$ & $I^{**}(v_0)$, by (6.8),

$\qquad M \vDash \exists!t^+ \, MinMax^+T_b(t^+,u_0)$ & $\exists!t^{++} \, MinMax^+T_b(t^{++},v_0)$.

We distinguish four basic cases, based on the definition of $\prec$:

$\qquad$ (*Biib2.1): $M \vDash u_0 \triangleleft_{Tb} y$ & $y \triangleleft_{Tb} v_0$.

$\Rightarrow M \vDash t^+ < t_0 < t^{++}$,

$\Rightarrow$ from $M \vDash Pref(av_0a,t_2)$, $M \vDash Max^+T_b(t_2,v_0)$,

$\Rightarrow M \vDash t^{++} \leq t_2$,

$\Rightarrow M \vDash t_0 < t^{++} \leq t_2$,

$\Rightarrow$ from $M \vDash MinMax^+T_b(t_0,y)$, $M \vDash Pref(aya,t_0)$,

$\Rightarrow$ from $M \vDash \neg(y \, \varepsilon \, x)$, $M \vDash \neg(y \, \varepsilon \, x")$.

So we have

$\quad M \vDash x"=t_1au_0at_2aw"$ & $Firstf(x",t_1,au_0a,t_2)$ & $t_0 < t_2$ & $\neg(y \, \varepsilon \, x")$ & $Pref(aya,t_0)$ &

$\hfill$ & $MinSet(x")$.

Let $x^+ = t_0ayat_2aw_2$.

$\Rightarrow$ from $M \vDash I^{**}(t_0)$ & $I^{**}(y)$ & $I^{**}(t_2)$ & $I^{**}(w_2)$,

since $I^{**}$ may be assumed to be closed w.r. to *, $M \vDash I^{**}(x^+)$,

$\Rightarrow$ from $M \vDash I^{**}(x")$, by (10.13), $M \vDash MinSet(x^+)$,

$\Rightarrow$ by (5.35), $M \vDash Env(t,x^+)$ & $Firstf(x^+,t_0,aya,t_2)$ &

& $\forall w(w \, \varepsilon \, x^+ \leftrightarrow (w \, \varepsilon \, x" \, \& \, w \neq u_0) \vee w=y)$.

Let $x_0 = t_1 a u_0 a t_1$.

$\Rightarrow$ from $M \vDash I^{**}(t_1) \, \& \, I^{**}(u_0)$, since $I^{**}$ may be assumed to be closed w.r. to $*$,

$$M \vDash I^{**}(x_0),$$

$\Rightarrow$ from $M \vDash \text{Firstf}(x",t_1,au_0a,t_2), \, M \vDash \text{Pref}(au_0a,t_1)$,

$\Rightarrow M \vDash \text{Firstf}(x_0,t_1,au_0a,t_1) \, \& \, \text{Lastf}(x_0,t_1,au_0a,t_1)$,

$\Rightarrow$ by (5.22), $M \vDash \text{Env}(t_1,x_0) \, \& \, \forall w(w \, \varepsilon \, x_0 \leftrightarrow w=u_0)$,

$\Rightarrow$ by (10.5), $M \vDash \text{MinSet}(x_0)$,

$\Rightarrow$ from $M \vDash \text{Env}(t',x'), \, M \vDash \exists v' \, \text{Lastf}(x,t',av'a,t')$,

$\Rightarrow$ by (5.6), $M \vDash \exists t_4 \text{Fr}(x,t',av'a,t_4)$.

Hence we have

$M \vDash \text{Env}(t,x) \, \& \, \text{Fr}(x,t',av'a,t_4) \, \& \, \text{Lastf}(x,t',av'a,t') \, \& \, \text{Env}(t',x') \, \&$

& $x'Bx \, \& \, \text{MinSet}(x)$,

$\Rightarrow$ by (10.8), $M \vDash \text{MinSet}(x')$.

So we have

$M \vDash \text{Env}(t',x') \, \& \, x'=t^*w_0^*t' \, \& \, aBw_0^* \, \& \, aEw_0^* \, \& \, x^*=t^*w_0^*x_0 \, \& \, \text{Env}(t_1,x_0) \, \&$

& $t' < t_1 \, \& \, \text{Firstf}(x_0,t_1,au_0a,t_1) \, \& \, \neg \exists w(w \, \varepsilon \, x' \, \& \, w \, \varepsilon \, x_0) \, \&$

& $\text{MinSet}(x') \, \& \, \text{MinSet}(x_0)$,

$\Rightarrow$ from $M \vDash I^{**}(x')$, since $I^{**}$ may be assumed to be downward closed

w.r. to $\subseteq_p$, $M \vDash I^{**}(t^*)$,

$\Rightarrow$ from $M \vDash I^{**}(w_0^*) \, \& \, I^{**}(x_0)$, since $I^{**}$ may be assumed to be closed w.r. to $*$,

$$M \vDash I^{**}(x^*),$$

$\Rightarrow$ by (10.6), $M \models \text{MinSet}(x^*)$ & $\text{Env}(t_1,x^*)$,

whereas, by (5.46), $M \models \forall w(w \, \varepsilon \, x^* \leftrightarrow w \, \varepsilon \, x' \, \vee \, w \, \varepsilon \, x_0)$.

(1a)   $M \models t_1 < t_0 < t^{++}$.

Let $z = w_1 a t_1 a u_0 a t_0 a y a t_2 a w_2$.

$\Rightarrow M \models z = (w_1 a) t_1 a u_0 a (t_0 a y a t_2 a w_2) = t^*(w_0^* t_1 a u_0 a) x^+$.

$\Rightarrow$ from $M \models I^{**}(t^*)$ & $I^{**}(w_0^*)$ & $I^{**}(t_1)$ & $I^{**}(u_0)$ & $I^{**}(x^+)$, since $I^{**}$ may be assumed to be closed w.r. to $*$,    $M \models I^{**}(z)$,

$\Rightarrow$ from $M \models I^{**}(z)$,  since $I^{**}$ may be assumed to be downward closed w.r. to $\subseteq_p$,    $M \models I^{**}(w_1)$.

Then we have

$M \models \text{Env}(t_1,x^*)$ & $x^* = t^*(w_0^* \, t_1 a u_0 a) t_1$ & $a B w_0^*$ & $z = t^* w_0^* t_1 a u_0 a x^+$ &

   & $\text{Env}(t,x^+)$ & $t_1 < t_0$ & $\text{Firstf}(x^+,t_0,a u_0 a,t_2)$ & $\neg \exists w(w \, \varepsilon \, x^* \, \& \, w \, \varepsilon \, x^+)$ &

                                                          & $\text{MinSet}(x^*)$ & $\text{MinSet}(x^+)$,

$\Rightarrow$ by (10.6),  $M \models \text{MinSet}(z)$.

From the principal hypothesis $M \models \text{Set}^*(x)$ we have $M \models \text{Lex}^+(x)$.

Also, by (5.46),  $M \models \forall w(w \, \varepsilon \, z \leftrightarrow w \, \varepsilon \, x^* \, \vee \, w \, \varepsilon \, x^+)$.

Applying the RESOLUTION LEMMA we then have

   $M \models \text{Lex}^+(x')$ & $\text{Lex}^+(x'')$.

We claim that   $M \models \forall w(w \, \varepsilon \, x'' \, \& \, w \neq u_0 \rightarrow y \prec w)$.

Assume that   $M \models w \, \varepsilon \, x'' \, \& \, w \neq u_0$.

$\Rightarrow$ since $M \models I^{**}(x'')$ & $\text{Firstf}(x'',t_1,a u_0 a,t_2)$, by (9.1) and (9.7), $M \models u_0 <_{x''} w$,

$\Rightarrow$ by (9.18), $M \models \models u_0 <_x w$,

$\Rightarrow$ by (11.5), $M \vDash y \prec w$, as claimed.

Then, from $M \vDash I^{**}(x") \& I^{**}(x^+)$, by (9.26), $M \vDash Lex^+(x^+)$.

Trivially, $M \vDash Lex^+(x_0)$.

We claim that $M \vDash \forall u,v\ (u\ \varepsilon\ x'\ \&\ v\ \varepsilon\ x_0 \to u \prec v)$.

Assume $M \vDash u\ \varepsilon\ x'\ \&\ v = u_0$.

$\Rightarrow$ $M \vDash \exists t_6, t_7\ Fr(x', t_6, aua, t_7)$,

$\Rightarrow$ from $M \vDash Env(t', x')$, $M \vDash MaxT_b(t', x')$,

$\Rightarrow$ $M \vDash t_6 \leq t' < t_1$,

$\Rightarrow$ by (5.6), $M \vDash \exists t_8\ Fr(x, t_6, aua, t_8)$,

$\Rightarrow$ from $M \vDash Firstf(x", t_1, au_0a, t_2)$, by (5.25), $M \vDash Fr(x, t_1, au_0a, t_2)$,

$\Rightarrow$ by (9.14), $M \vDash u <_x u_0$,

$\Rightarrow$ since $M \vDash Lex^+(x)$, $M \vDash u \prec u_0$,

$\Rightarrow$ $M \vDash u \prec v$, as claimed.

$\Rightarrow$ since $M \vDash I^{**}(x^*)\ \&\ \forall w(w\ \varepsilon\ x^* \leftrightarrow w\ \varepsilon\ x'\ v\ w\ \varepsilon\ x_0)\ \&\ \neg\exists w(w\ \varepsilon\ x'\ \&\ w\ \varepsilon\ x_0)$,

by (9.29), $M \vDash Lex^+(x^*)$.

$\Rightarrow$ since $M \vDash I^{**}(z)\ \&\ \forall w(w\ \varepsilon\ z \leftrightarrow w\ \varepsilon\ x^*\ v\ w\ \varepsilon\ x^+)\ \&\ \neg\exists w(w\ \varepsilon\ x^*\ \&\ w\ \varepsilon\ x^+)$,

by (9.29), $M \vDash Lex^+(z)$.

Let $x^{++} = w_1 a t_1 a u_0 a t_1$.

$\Rightarrow$ from $M \vDash I^{**}(x)\ \&\ Env(t,x)\ \&\ Intf(x, w_1, t_1, au_0a, t_2)$, by (5.53),

$\qquad M \vDash Env(t_1, x^{++})\ \&\ Lastf(x^{++}, t_1, au_0a, t_1)$,

$\Rightarrow$ $M \vDash x^{++} = w_1 a t_1 a u_0 a t_1 = t^* w_0^*\ t_1 a u_0 a t_1 = x^*$,

$\Rightarrow$ $M \vDash x^{++} = t^* w_1^*\ t_1\ \&\ Tally_b(t^*)\ \&\ aBw_1^*\ \&\ aEw_1^*$ where $w_1^* = w_0^* t_1 au_0 a$.

By (10.23), from $M \vDash Special(x)$, we have $M \vDash Special(x^{++})$.

We also have $M \vDash MinMax^+T_b(t_0,y)$.

Let $x_0^* = t_0ayat_0$.

$\Rightarrow$ from $M \vDash I^{**}(t_0)$ & $I^{**}(y)$, since $I^{**}$ may be assumed to be closed w.r. to $*$,

$$M \vDash I^{**}(x_0^*),$$

$\Rightarrow$ $M \vDash Pref(aya,t_0)$,

$\Rightarrow$ $M \vDash Firstf(x_0^*,t_0,aya,t_0)$ & $Lastf(x_0^*,t_0,aya,t_0)$,

$\Rightarrow$ by (5.22), $M \vDash Env(t_0,x_0^*)$ & $\forall w(w \varepsilon x_0^* \leftrightarrow w=y)$,

$\Rightarrow$ by (10.15), $M \vDash Special(x_0^*)$.

Consider $x^{**} = w_1at_1au_0at_0ayat_0$.

$\Rightarrow$ from $M \vDash I^{**}(w_1)$ & $I^{**}(t_1)$ & $I^{**}(u_0)$ & $I^{**}(t_0)$ & $I^{**}(y)$, since $I^{**}$ may be assumed to be closed w.r. to $*$,    $M \vDash I^{**}(x^{**})$,

$\Rightarrow$ $M \vDash x^{**} = w_1at_1au_0a(t_1t'')ayat_0 = x^{++}t''ayat_0$, where $M \vDash t_1t'' = t_0$,

$\Rightarrow$ since $M \vDash Tally_b(t_1)$ & $Tally_b(t_0)$ & $t_1 < t_0$, $M \vDash Tally_b(t'')$,

$\Rightarrow$ $M \vDash Env(t_0,x_0^*)$ & $x^{**} = t^*w_0^*t_1t''ayat_0$ & $Tally_b(t_0)$ & $aBw_0^*$ & $aEw_0^*$ &

& $x^{++} = t^*w_0^* t_1$ & $x_0^* = t_1t''ayat_0$ & $Env(t_1,x^{++})$ & $Firstf(x_0^*,t_1t'',aya,t_0)$ & $t_1 < t_0$,

$\Rightarrow$ from $M \vDash I^{**}(x^{**})$, by (5.41), $M \vDash \neg \exists w(w \varepsilon x^{++}$ & $w \varepsilon x_0^*)$,

$\Rightarrow$ by (5.46), $M \vDash \forall w(w \varepsilon x^{**} \leftrightarrow w \varepsilon x^{++} \vee w \varepsilon x_0^*)$.

We claim that $M \vDash \forall u,v(u \varepsilon x^{++}$ & $v \varepsilon x_0^* \rightarrow u \prec v)$.

Assume $M \vDash u \varepsilon x^{++}$ & $v \varepsilon x_0^*$.

$\Rightarrow$ $M \vDash v = y$,

$\Rightarrow$ by choice of $u_0$, $M \vDash \forall w(w \leq_x u_0 \rightarrow w \prec y)$.

We have that $M \vDash \text{Lastf}(x^{++},t_1,au_0a,t_1)$.

$\implies$ by (9.3), $M \vDash u \leq_{x^{++}} u_0$,

$\implies$ since $M \vDash \text{Env}(t,x)$ & $\text{Intf}(x,w_1,t_1,au_0a,t_2)$ & $x^{++}=w_1at_1au_0at_1$, by (9.15),

$$M \vDash u \leq_x u_0,$$

$\implies M \vDash u \prec y$,

$\implies M \vDash u \prec v$, as claimed.

So we have

$M \vDash \text{Env}(t_0,x^{**})$ & $x^{**}= x^{++}t"ayat_0$ & $\text{Env}(t_1,x^{++})$ & $x^{++}=t^*w_0^* t_1$ & $aBw_0^*$ &

& $aEw_0^*$ & $\text{Env}(t_0,x_0^*)$ & $\text{Firstf}(x_0^*,t_0,aya,t_0)$ & $t_1<t_0$ &

& $\neg\exists w(w \,\varepsilon\, x^{++}$ & $w \,\varepsilon\, x_0^*)$ & $\forall u,v(u \,\varepsilon\, x^{++}$ & $v \,\varepsilon\, x_0^* \to u \prec v)$ &

& $\text{Special}(x^{++})$ & $\text{Special}(x_0^*)$,

$\implies$ from $M \vDash I^{**}(x^{**})$, by (10.17), $M \vDash \text{Special}(x^{**})$.

We claim that $M \vDash \neg\exists w(w \,\varepsilon\, x^{**}$ & $w \,\varepsilon\, x^-)$.

Suppose, for a reductio, that $M \vDash w \,\varepsilon\, x^{**}$ & $w \,\varepsilon\, x^-$.

$\implies M \vDash w \,\varepsilon\, x^{++} \vee w \,\varepsilon\, x_0^*$.

Assume that $M \vDash w \,\varepsilon\, x^{++}$.

$\implies$ from $M \vDash I^{**}(x^-)$ & $\text{Env}(t,x^-)$, by (5.11),

$M \vDash \exists t_7,w^-(\text{Tally}_b(t_7)$ & $x^-=t_7w^-t$ & $aB\, w^-$ & $aEw^-)$,

$\implies M \vDash \exists w_8\; t_2aw_2=x^-=t_7w^-t=t_7aw_8t$,

$\implies$ from $M \vDash I^{**}(x^-)$, since $I^{**}$ may be assumed to be downward closed

w.r. to $\subseteq_p$,    $M \vDash I^{**}(t_7)$,

$\implies$ from $M \vDash I^{**}(t_2)$, by (4.23$^b$), $M \vDash t_2=t_7$.

We have that  $M \vDash t_1 < t_0 < t^{++} \leq t_2$.

$\Rightarrow$ since $M \vDash \text{Tally}_b(t_0)$ & $\text{Tally}_b(t_2)$, $M \vDash \exists t_2'(\text{Tally}_b(t_2') \& t_0 t_2' = t_2)$,

$\Rightarrow M \vDash t_1 t'' = t_0 \& t_1(t'' t_2') = t_0 t_2' = t_2$.

Therefore,

$\Rightarrow M \vDash x^{++} t_2' a w_2 = w_1 a t_1 a u_0 a t_1 t'' t_2' a w_2 = w_1 a t_1 a u_0 a t_2 a w_2 = w_1 a t_1 a u_0 a\ x^- =$

$\qquad\qquad\qquad = w_1 a t_1 a u_0 a t_2 w^- t = w_1 a t_1 a u_0 a t_1 t'' t_2' w^- t = x^{++} t'' t_2' w^- t$,

$\Rightarrow$ from hypothesis $M \vDash w\ \varepsilon\ x^{++}$, $M \vDash \exists t_9, t_{10}\ \text{Fr}(x^{++}, t_9, awa, t_{10})$,

$\Rightarrow$ from $M \vDash \text{Env}(t_1, x^{++})$, $M \vDash \text{MaxT}_b(t_1, x^{++})$,

$\Rightarrow M \vDash t_9 \leq t_1$,

$\Rightarrow$ by (5.6), $M \vDash \exists t_{11}\ \text{Fr}(x, t_9, awa, t_{11})$.

On the other hand, we also have that

$\qquad M \vDash x = x^{++} t'' t_2' w^- t = (t^* w_0^* \ t_1) t'' t_2' w^- t = t^* w_0^* t_2 w^- t = t^* w_0^* x^-$.

From $M \vDash w\ \varepsilon\ x^-$, we have $M \vDash \exists t_{12}, t_{13}\ \text{Fr}(x^-, t_{12}, awa, t_{13})$.

$\Rightarrow$ from $M \vDash I^{**}(x^-)$ & $\text{Firstf}(x^-, t_2^+, av_0 a, t_3)$, by (5.20), $M \vDash t_2 \leq t_{12}$.

We claim that $M \vDash \text{Max}^+ T_b(t_{12}, t^* w_0^*)$.

Assume $M \vDash \text{Tally}_b(t_{14})\ \&\ t_{14} \subseteq_p t^* w_0^*$.

$\Rightarrow M \vDash t_{14} \subseteq_p t^* w_0^* \subseteq_p t^* w_0^* t_1 = x^{++}$,

$\Rightarrow M \vDash t_{14} \leq t_1 < t_2 \leq t_{12}$.

Therefore, $M \vDash \text{Max}^+ T_b(t_{12}, t^* w_0^*)$.

So we have

$\quad M \vDash \text{Fr}(x^-, t_{12}, awa, t_{13})\ \&\ x = t^* w_0^* x^-\ \&\ aBw_0^*\ \&\ aEw_0^*\ \&\ \text{Max}^+ T_b(t_{12}, t^* w_0^*)$.

$\Rightarrow$ by (5.25),  $M \vDash Fr(x,t_{12},awa,t_{13})$,

$\Rightarrow$ $M \vDash t_9 \leq t_1 < t_2 \leq t_{12}$,

$\Rightarrow$ from  $M \vDash Env(t,x)$,  $M \vDash t_9 = t_{12}$,

$\Rightarrow$ $M \vDash t_9 \leq t_1 < t_2 \leq t_{12} = t_9$, contradicting  $M \vDash t_9 \in I \subseteq I_0$.

Therefore,  $M \vDash \neg w \,\varepsilon\, x^{++}$.

Assume   $M \vDash w \,\varepsilon\, x_0^*$.

$\Rightarrow$ $M \vDash w = y$,

$\Rightarrow$ from hypothesis  $M \vDash w \,\varepsilon\, x^-$,  $M \vDash w \,\varepsilon\, x$,

$\Rightarrow$ $M \vDash y \,\varepsilon\, x$, contradicting hypothesis.

This completes the proof of the claim that   $M \vDash \neg \exists w(w \,\varepsilon\, x^{**} \,\&\, w \,\varepsilon\, x^-)$.

We then have

$\quad M \vDash Env(t_0,x^{**}) \,\&\, x^{**} = t^*w^{**}t_0 \,\&\, aBw^{**} \,\&\, aEw^{**} \,\&\, Env(t,x^-) \,\&$

$\qquad\qquad \,\&\, Firstf(x^-,t_2,av_0a,t_3) \,\&\, t_0 < t_2 \,\&\, \neg \exists w(w \,\varepsilon\, x^{**} \,\&\, w \,\varepsilon\, x^-)$.

We also have

$\quad M \vDash z = w_1 at_1 au_0 at_0 ayat_2 aw_2 = (w_1 at_1 au_0 at_1 \text{''} ayat_1 \text{''}) t_2\text{'} aw_2 = x^{**} t_2\text{'} aw_2 =$

$\qquad\qquad\qquad\qquad\qquad\qquad\qquad = t^*w^{**}t_0 t_2\text{'} aw_2 = t^*w^{**}x^-$,

$\Rightarrow$ from  $M \vDash I^{**}(z)$, by (5.46),  $M \vDash \forall w(w \,\varepsilon\, z \leftrightarrow w \,\varepsilon\, x^{**} \vee w \,\varepsilon\, x^-) \,\&\, Env(t,z)$.

Now we have

$\quad M \vDash Env(t,x) \,\&\, x = t^*w_0^*x^- \,\&\, aBw_0^* \,\&\, aEw_0^* \,\&\, Env(t_1,x^{++}) \,\&\, x^{++} = t^*w_0^*t_1 \,\&$

$\qquad \,\&\, Env(t,x^-) \,\&\, Lastf(x^{++},t_1,au_0a,t_1) \,\&\, Firstf(x^-,t_2,av_0a,t_3) \,\&$

$\qquad\qquad\qquad \,\&\, MinMax^+T_b(t^{++},v_0) \,\&\, t_1 < t^{++} \,\&\, Special(x)$,

$\Rightarrow$ by (10.25),  $M \vDash t_2 = t^{++}$,

$\Rightarrow$ by (10.24), $M \vDash Special(x^-)$.

Note that we already established that

$$M \vDash \forall w(w \, \varepsilon \, x^{++} \to w \prec y) \, \& \, \forall w(w \, \varepsilon \, x^- \to y \prec w),$$

$\Rightarrow$ from $M \vDash I^{**}(y)$, by (8.3),

$$M \vDash \forall w',w''((w' \, \varepsilon \, x^{++} \vee w'=y) \, \& \, w'' \, \varepsilon \, x^- \to w' \prec w''),$$

$\Rightarrow$ since $M \vDash \forall w(w \, \varepsilon \, x^{**} \leftrightarrow w \, \varepsilon \, x^{++} \vee w=y)$,

$$M \vDash \forall w',w''(w' \, \varepsilon \, x^{**} \, \& \, w'' \, \varepsilon \, x^- \to w' \prec w'').$$

We then have

$M \vDash Env(t,z) \, \& \, z=t*w^{**}x^- \, \& \, Env(t_0,x^{**}) \, \& \, x^{**}= t*w^{**}t_0 \, \& \, aBw^{**} \, \& \, aEw^{**} \, \&$

$\quad \& \, Env(t,x^-) \, \& \, Firstf(x^-,t_2,av_0a,t_5) \, \& \, t_0 \prec t_2 \, \& \, \neg \exists w(w \, \varepsilon \, x^{**} \, \& \, w \, \varepsilon \, x^-) \, \&$

$\quad \& \, \forall w',w''(w' \, \varepsilon \, x^{**} \, \& \, w'' \, \varepsilon \, x^- \to w' \prec w'') \, \& \, Special(x^{**}) \, \& \, Special(x^-),$

$\Rightarrow$ from $M \vDash I^{**}(z)$, by (10.17), $M \vDash Special(z)$.

We now show that

$$M \vDash \forall w(w \, \varepsilon \, z \leftrightarrow w \, \varepsilon \, x \vee w=y).$$

Recall that we have $M \vDash \forall w(w \, \varepsilon \, x^* \leftrightarrow w \, \varepsilon \, x' \vee w \, \varepsilon \, x_0)$.

$\Rightarrow$ since $M \vDash \forall w(w \, \varepsilon \, x_0 \leftrightarrow w=u_0)$, $M \vDash \forall w(w \, \varepsilon \, x^* \leftrightarrow w \, \varepsilon \, x' \vee w=u_0)$,

$\Rightarrow$ from $M \vDash u_0 \, \varepsilon \, x''$ and $M \vDash \forall w(w \, \varepsilon \, x \leftrightarrow w \, \varepsilon \, x' \vee w \, \varepsilon \, x'')$,

$$M \vDash \forall w(w \, \varepsilon \, x^* \leftrightarrow w \, \varepsilon \, x' \vee w \, \varepsilon \, x''),$$

$\Rightarrow$ $M \vDash \forall w(w \, \varepsilon \, x^* \to w \, \varepsilon \, x)$.

On the other hand, we have $M \vDash \forall w(w \, \varepsilon \, x^- \leftrightarrow w \, \varepsilon \, x'' \, \& \, w \neq u_0)$.

$\Rightarrow$ $M \vDash \forall w(w \, \varepsilon \, x^- \to w \, \varepsilon \, x'')$,

$\Rightarrow$ $M \vDash \forall w(w \, \varepsilon \, x^- \to w \, \varepsilon \, x)$.

We also have $M \vDash \forall w(w \, \varepsilon \, z \leftrightarrow w \, \varepsilon \, x^{**} \vee w \, \varepsilon \, x^-)$, where

$$M \vDash \forall w(w \, \varepsilon \, x^{**} \leftrightarrow w \, \varepsilon \, x^{++} \vee w \, \varepsilon \, x_0^*).$$

But $M \vDash \forall w(w \, \varepsilon \, x_0^* \leftrightarrow w=y)$, whereas $M \vDash x^{++}=x^*$.

Suppose now that $M \vDash w \, \varepsilon \, z$.

$\Rightarrow M \vDash w \, \varepsilon \, x^{**} \vee w \, \varepsilon \, x^-$,

$\Rightarrow M \vDash (w \, \varepsilon \, x^{++} \vee w \, \varepsilon \, x_0^*) \vee w \, \varepsilon \, x$,

$\Rightarrow M \vDash w \, \varepsilon \, x^* \vee w=y \vee w \, \varepsilon \, x$,

$\Rightarrow M \vDash w \, \varepsilon \, x \vee w=y$, as required.

Conversely, suppose that $M \vDash w \, \varepsilon \, x \vee w=y$.

If $M \vDash w=y$, then $M \vDash w \, \varepsilon \, x_0^*$, whence $M \vDash w \, \varepsilon \, x^{**}$. So $M \vDash w \, \varepsilon \, z$.

Thus, we may assume that $M \vDash w \, \varepsilon \, x$.

$\Rightarrow M \vDash w \, \varepsilon \, x' \vee w \, \varepsilon \, x''$.

Assume that $M \vDash w \, \varepsilon \, x'$.

$\Rightarrow M \vDash w \, \varepsilon \, x^*$,

$\Rightarrow M \vDash w \, \varepsilon \, z$, as required.

Assume that $M \vDash w \, \varepsilon \, x''$.

If $M \vDash w \, \varepsilon \, x''$ & $w=u_0$, then $M \vDash w \, \varepsilon \, x_0$.

$\Rightarrow M \vDash w \, \varepsilon \, x^*$,

$\Rightarrow M \vDash w \, \varepsilon \, z$, as required.

If $M \vDash w \, \varepsilon \, x''$ & $w \neq u_0$, recall that $M \vDash \forall w(w \, \varepsilon \, x^+ \leftrightarrow ((w \, \varepsilon \, x'' \, \& \, w \neq u_0) \vee w=y))$.

$\Rightarrow M \vDash w \, \varepsilon \, x^+$,

$\Rightarrow M \vDash w \, \varepsilon \, z$, as required.

Hence  $M \vDash w \varepsilon x \to w \varepsilon z$.

This completes the argument that  $M \vDash \forall w(w \varepsilon z \leftrightarrow w \varepsilon x \vee w = y)$.

(1b)  $M \vDash t_1 = t_0$.

  (1bi)  $M \vDash t_1 = t_0$ & $t_0 b = t^{++}$.

We have that       $M \vDash MinSet(x^+)$ & $Env(t, x^+)$

where  $M \vDash x^+ = t_0 ayax^- = t_0 ayat_2 aw_2$.

$\Rightarrow$ since $M \vDash J^{**}(t)$, by (10.29),

  $M \vDash \exists! x^{+*} \in I^{**} \exists t' \in J^{**} (Env(t', x^{+*})$ & $x^+ \sim x^{+*}$ &

          & $MinSet(x^{+*})$ & $\forall u, v (u \leq_{x+} v \leftrightarrow u \leq_{x+*} v)$ &

        & $\forall v, t_3, t_4 (Free^+(x^+, t_3, ava, t_4) \to \exists t_5 Fr(x^{+*}, t_3 b, ava, t_5))$ &

& $\forall v, t_3, t_4 (Bound(x^+, t_3, ava, t_4) \vee Free^-(x^+, t_3, ava, t_4) \to \exists t_5 Fr(x^{+*}, t_3, ava, t_5)))$.

Let  $z = w_1 a t_1 a u_0 a x^{+*}$.

$\Rightarrow$ from $M \vDash I^{**}(w_1)$ & $I^{**}(t_1)$ & $I^{**}(u_0)$ & $I^{**}(x^{+*})$, since $I^{**}$ may be

assumed to be closed w.r. to $*$,     $M \vDash I^{**}(z)$.

We have that  $M \vDash Firstf(x^+, t_0, aya, t_2)$.

$\Rightarrow$  $M \vDash Free^+(x^+, t_0, aya, t_2)$,

$\Rightarrow$  $M \vDash \exists t_6 Fr(x^{+*}, t_0 b, aya, t_6)$,

$\Rightarrow$  $M \vDash (t_0 ba) B x^{+*}$,

$\Rightarrow$ from  $M \vDash I^{**}(x^{+*})$ & $Env(t', x^{+*})$, by (5.11),

    $M \vDash \exists t^{**}, w_4 (Tally_b(t^{**})$ & $aBw_4$ & $aEw_4$ & $x^{+*} = t^{**} w_4 t')$,

$\Rightarrow$  $M \vDash \exists x_1, w_5\ t_0 bax_1 = x^{+*} = t^{**}(aw_5) t'$,

$\Rightarrow$ from $M \vDash I^{**}(x^{+*})$ & $t^{**} \subseteq_p x^{+*}$, since $I^{**}$ may be assumed to be downward

closed w.r. to $\subseteq_p$,     $M \vDash I^{**}(t^{**})$,

$\implies$ from $M \vDash I^{**}(t_0 b)$, by (4.23$^b$), $M \vDash t_0 b = t^{**}$,

$\implies M \vDash x^{+*} = t_0 b w_4 t'$,

$\implies M \vDash z = (w_1 a t_1 a u_0 a t_1) b w_4 t' = x^{++} b w_4 t'$.

By (10.8), from $M \vDash \text{MinSet}(x)$ we have $M \vDash \text{MinSet}(x^{++})$.

We have that

$$M \vDash t^* w_0^* t_1 = x^{++} = w_1 a t_1 a u_0 a t_1.$$

$\implies$ by (3.6), $M \vDash t^* w_0^* = w_1 a t_1 a u_0 a$,

$\implies M \vDash z = (w_1 a t_1 a u_0 a) x^{+*} = t^* w_0^* x^{+*}$.

Suppose, for a reductio, that $M \vDash \exists w(w \, \varepsilon \, x^{++} \, \& \, w \, \varepsilon \, x^{+*})$.

We have that

$$M \vDash w \, \varepsilon \, x^{+*} \leftrightarrow w \, \varepsilon \, x^+ \leftrightarrow (w \, \varepsilon \, x'' \, \& \, w \neq u_0) \lor w = y \leftrightarrow w \, \varepsilon \, x^- \lor w = y.$$

Suppose that (a) $M \vDash w \, \varepsilon \, x^{++} \, \& \, w \, \varepsilon \, x^-$.

$\implies M \vDash w \, \varepsilon \, x^{**} \, \& \, w \, \varepsilon \, x^-$, contradicting earlier proof.

Suppose that (b) $M \vDash w \, \varepsilon \, x^{++} \, \& \, w = y$.

$\implies M \vDash w \, \varepsilon \, x^{**} \, \& \, w \, \varepsilon \, x_0^*$, again contradicting earlier proof.

Therefore, $M \vDash \neg \exists w(w \, \varepsilon \, x^{++} \, \& \, w \, \varepsilon \, x^{+*})$.

$\implies$ from $M \vDash I^{**}(z)$, by (10.6), $M \vDash \text{MinSet}(z) \, \& \, \text{Env}(t',z)$.

Now, we have $M \vDash x = (w_1 a t_1 a u_0 a) t_2 a w_2 = t^* w_0^* x^-$.

$\implies M \vDash \text{Env}(t,x) \, \& \, x = t^* w_0^* x^- \, \& \, x^{++} = t^* w_0^* t_1 \, \& \, \text{Env}(t_1, x^{++}) \, \&$

  $\& \, aBw_0^* \, \& \, aEw_0^* \, \& \, \text{Env}(t, x^-) \, \& \, \text{Firstf}(x^-, t_2, a v_0 a, t_3) \, \& \, t_1 < t_2 \, \& \, \text{Lex}^+(x)$,

$\implies$ by (9.28), $M \vDash \text{Lex}^+(x^{++}) \, \& \, \text{Lex}^+(x^-)$.

We have that, trivially, $M \vDash \text{Lex}^+(x_0{}^*)$.

Since $M \vDash \forall w(w \: \varepsilon \: x_0{}^* \leftrightarrow w=y)$, we then have

$\quad M \vDash \text{Env}(t,x^+) \: \& \: x^+ = t_0 a y a x^- \: \& \: \text{Env}(t_0,x_0{}^*) \: \& \: x_0{}^* = t_0 a y a t_0 \: \& \: \text{Env}(t,x^-) \: \&$

$\quad\quad \& \: \text{Firstf}(x^-,t_2,av_0a,t_3) \: \& \: t_0 < t_2 \: \& \: \neg \exists w(w \: \varepsilon \: x_0{}^* \: \& \: w \: \varepsilon \: x^-) \: \&$

$\quad\quad\quad \& \: \forall u,v(u \: \varepsilon \: x_0{}^* \: \& \: v \: \varepsilon \: x^- \to u \prec v) \: \& \: \text{Lex}^+(x_0{}^*) \: \& \: \text{Lex}^+(x^-)$.

$\Longrightarrow$ by (9.29), $M \vDash \text{Lex}^+(x^+)$,

$\Longrightarrow M \vDash \forall u,v \: (u <_{x^+} v \to u \prec v)$,

$\Longrightarrow M \vDash \forall u,v \: (u <_{x^{+*}} v \to u \prec v)$.

Hence $M \vDash \text{Lex}^+(x^{+*})$.

We claim that $M \vDash \forall u,v \: (u \: \varepsilon \: x^{++} \: \& \: v \: \varepsilon \: x^{+*} \to u \prec v)$.

Assume $M \vDash u \: \varepsilon \: x^{++} \: \& \: v \: \varepsilon \: x^{+*}$.

We have $M \vDash \forall w \: (w \: \varepsilon \: x^{++} \to w \prec y)$.

Also, $M \vDash \forall w \: (w \: \varepsilon \: x^{+*} \leftrightarrow w=y \lor w \: \varepsilon \: x^-)$ and $M \vDash \forall w(w \: \varepsilon \: x^- \leftrightarrow w \: \varepsilon \: x'' \: \& \: w \neq u_0)$.

Hence, by earlier proof, $M \vDash \forall w(w \: \varepsilon \: x^- \to y \prec w)$.

But then $M \vDash u \prec y \: \& \: (v=y \lor y \prec v)$, whence $M \vDash u \prec v$, as claimed.

Then from $M \vDash I^{**}(z) \: \& \: \text{Lex}^+(x^{++}) \: \& \: \text{Lex}^+(x^{+*})$, by (9.29), $M \vDash \text{Lex}^+(z)$.

We have that $M \vDash t^{++} = t_2$.

$\Longrightarrow$ from hypothesis $M \vDash \text{Special}(x)$, by (10.24), $M \vDash \text{Special}(x^-)$,

$\Longrightarrow$ from $M \vDash I^{**}(x^+)$, by (10.20), $M \vDash \text{Special}(x^+)$,

$\Longrightarrow$ by earlier proof, $M \vDash \text{Special}(x^{++})$.

Now, we have that

$M \vDash \text{Env}(t',z) \ \& \ z=t*w_0* \ x^{+*} \ \& \ aBw_0* \ \& \ aEw_0* \ \& \ \text{Set}(x^{++}) \ \& \ x^{++}=t*w_0*t_1 \ \&$

$\& \ t_1=t_0 \ \& \ \text{Env}(t_1,x^{++}) \ \& \ \text{Firstf}(x^{+*},t_0b,aya,t_6) \ \&$

$\& \ \forall u,v \ (u \ \varepsilon \ x^{++} \ \& \ v \ \varepsilon \ x^{+*} \to u < v) \ \& \ \neg\exists u(u \ \varepsilon \ x^{++} \ \& \ u \ \varepsilon \ x^{+*}) \ \& \ \text{Set}(x^+) \ \&$

$\& \ x^+ \sim x^{+*} \ \& \ \forall u,v \ (u <_{x^+} v \leftrightarrow u <_{x^{+*}} v) \ \&$

$\& \ \forall v,t_3,t_4 \ (\text{Free}^+(x^+,t_3,ava,t_4) \to \exists t_5 \text{Fr}(x^{+*},t_3 b,ava,t_5)) \ \&$

$\& \ \forall v,t_3,t_4 \ (\text{Bound}(x^+,t_3,ava,t_4) \ \lor \ \text{Free}^-(x^+,t_3,ava,t_4) \to \exists t_5 \text{Fr}(x^{+*},t_3,ava,t_5)) \ \&$

$\& \ \text{Special}(x^{++}) \ \& \ \text{Special}(x^+),$

$\implies$ from $M \vDash I^{**}(z) \ \& \ I^{**}(x^+)$, by (10.28), $M \vDash \text{Special}(z)$.

Hence $M \vDash V^{**}(z)$.

It remains to show that $M \vDash \forall w(w \ \varepsilon \ z \leftrightarrow w \ \varepsilon \ x \ \lor \ w=y)$.

Now, we have that

$M \vDash \text{Env}(t_1,x^{++}) \ \& \ x^{++}=t*w_0*t_1 \ \& \ aBw_0* \ \& \ aEw_0* \ \& \ \text{Env}(t',x^{+*}) \ \&$

$\& \ \text{Firstf}(x^{+*},t_0b,aya,t_6) \ \& \ t_1 < t_0 b \ \& \ \neg\exists u(u \ \varepsilon \ x^{++} \ \& \ u \ \varepsilon \ x^{+*}),$

$\implies$ from $M \vDash I^{**}(x^{++}) \ \& \ I^{**}(x^{+*})$, by (5.46),

$M \vDash \forall w(w \ \varepsilon \ z \leftrightarrow w \ \varepsilon \ x^{++} \ \lor \ w \ \varepsilon \ x^{+*})$.

Assume $M \vDash w \ \varepsilon \ z$.

$\implies M \vDash w \ \varepsilon \ x^{++} \ \lor \ w \ \varepsilon \ x^{+*}$.

Suppose $M \vDash w \ \varepsilon \ x^{++}$.

$\implies$ from $M \vDash x^{++}=x^* \ \& \ \forall u(u \ \varepsilon \ x^* \to u \ \varepsilon \ x), \ M \vDash w \ \varepsilon \ x$.

Suppose $M \vDash w \ \varepsilon \ x^{+*}$.

$\implies M \vDash w \ \varepsilon \ x^+,$

$\implies M \vDash (w \ \varepsilon \ x" \ \& \ w \neq u_0) \ \lor \ w=y,$

$\Rightarrow$ $M \vDash w \, \varepsilon \, x \lor w=y$.

Therefore, $M \vDash \forall w(w \, \varepsilon \, z \to w \, \varepsilon \, x \lor w=y)$.

Conversely, suppose that $M \vDash w \, \varepsilon \, x \lor w=y$.

Assume $M \vDash w=y$.

$\Rightarrow$ from $M \vDash \text{Firstf}(x^+, t_0, aya, t_2)$, $M \vDash w \, \varepsilon \, x^+$,

$\Rightarrow$ $M \vDash w \, \varepsilon \, x^{+*}$,

$\Rightarrow$ $M \vDash w \, \varepsilon \, z$.

The argument that $M \vDash w \, \varepsilon \, x \to w \, \varepsilon \, z$ is exactly analogous to that in (1a).

Therefore also $M \vDash \forall w(w \, \varepsilon \, x \lor w=y \to w \, \varepsilon \, z)$.

(1bii) $M \vDash t_1 = t_0 \,\&\, t_0 b < t^{++}$.

Let $z = w_1 a t_1 a u_0 a t_0 b a y a t_2 a w_2$.

$\Rightarrow$ from $M \vDash I^{**}(w_1) \,\&\, I^{**}(t_1) \,\&\, I^{**}(u_0) \,\&\, I^{**}(t_0) \,\&\, I^{**}(y) \,\&\, I^{**}(t_2) \,\&\, I^{**}(w_2)$,

since $I^{**}$ may be assumed to be closed w.r. to *, $\quad M \vDash I^{**}(z)$.

Let $x_1 = t_0 b a y a t_0 b$.

$\Rightarrow$ from $M \vDash I^{**}(t_0) \,\&\, I^{**}(y)$, $M \vDash I^{**}(x_1)$,

$\Rightarrow$ from $M \vDash \text{MinMax}^+ T_b(t_0, aya)$, $M \vDash \text{Pref}(aya, t_0) \,\&\, \text{Tally}_b(t_0)$,

$\Rightarrow$ by (4.1), $M \vDash \text{Tally}_b(t_0 b)$,

$\Rightarrow$ $M \vDash \text{Pref}(aya, t_0 b)$,

$\Rightarrow$ $M \vDash \text{Firstf}(x_1, t_0 b, aya, t_0 b) \,\&\, \text{Lastf}(x_1, t_0 b, aya, t_0 b)$,

$\Rightarrow$ from $M \vDash I^{**}(x_1)$, by (5.22), $M \vDash \text{Env}(t_0 b, x_1) \,\&\, \forall w(w \, \varepsilon \, x_1 \leftrightarrow w=y)$.

We have that $M \vDash t_0 t_2' = t_2$.

Let $x_0^{**} = t_0 b a y a t_2 a w_2$.

$\Rightarrow M \vDash x_0^{**} = t_0 bayax^-$,

$\Rightarrow$ from $M \vDash I^{**}(t_0)$ & $I^{**}(y)$ & $I^{**}(x^-)$, since $I^{**}$ may be assumed to be closed w.r. to *, $M \vDash I^{**}(x_0^{**})$,

$\Rightarrow$ from $M \vDash \forall w(w \varepsilon x_1 \leftrightarrow w=y \leftrightarrow w \varepsilon x_0^*)$ and $M \vDash \neg \exists w(w \varepsilon x_0^* \& w \varepsilon x^-)$,

$$M \vDash \neg \exists w(w \varepsilon x_1 \& w \varepsilon x^-).$$

So we have

$$M \vDash Env(t_0 b, x_1) \& x_1 = t_0 bayat_0 b \& Env(t, x^-) \& Firstf(x^-, t_2, av_0 a, t_3) \& t_0 b < t_2 \&$$
$$\& \neg \exists w(w \varepsilon x_1 \& w \varepsilon x^-).$$

$\Rightarrow$ from $M \vDash I^{**}(x_1)$ & $I^{**}(x^-)$, by (5.46),

$$M \vDash Env(t, x_0^{**}) \& \forall w(w \varepsilon x_0^{**} \leftrightarrow w \varepsilon x_1 \vee w \varepsilon x^- \leftrightarrow w=y \vee w \varepsilon x^-),$$

$\Rightarrow M \vDash \forall w(w \varepsilon x_0^{**} \leftrightarrow w=y \vee w \varepsilon x^- \leftrightarrow w \varepsilon x^{+*})$,

$\Rightarrow$ by earlier proof, $M \vDash \neg \exists w(w \varepsilon x^{++} \& w \varepsilon x^{+*})$,

$\Rightarrow M \vDash \neg \exists w(w \varepsilon x^{++} \& w \varepsilon x_0^{**})$.

Now, we have that

$$M \vDash x'' = t_1 au_0 at_2 aw_2 \& Firstf(x'', t_1, au_0 a, t_2) \& x_0^{**} = t_0 bayat_2 aw_2 \& t_0 b < t_2 \&$$
$$\& \neg(y \varepsilon x'') \& Pref(aya, t_0 b) \& MinSet(x'').$$

$\Rightarrow$ from $M \vDash I^{**}(x'')$ & $I^{**}(x_0^{**})$, by (10.13), $M \vDash MinSet(x_0^{**})$.

So we have

$$M \vDash z = w_1 at_1 au_0 at_0 bayat_2 aw_2 = w_1 at_1 au_0 at_1 t'' bayat_2 aw_2 = x^{++} t'' bayat_2 aw_2 =$$
$$= t*w_0*t_1 t'' bayat_2 aw_2 = t*w_0*t_0 bayat_2 aw_2 = t*w_0*x_0^{**},$$

and

$$M \vDash Env(t_1, x^{++}) \& x^{++} = t*w_0*t_1 \& aBw_0^* \& aEw_0^* \& z = t*w_0*x_0^{**} \&$$

   & Env(t,$x_0$**) & $t_1$≺$t_0$b & Firstf($x_0$**,$t_0$b,aya,$t_2$) & ¬∃w(w ε $x^{++}$ & w ε $x_0$**) &

                                               & MinSet($x^{++}$) & MinSet($x_0$**).

⟹ from M ⊨ I**(z), by (10.6), M ⊨ MinSet(z) & Env(t,z).

Since M ⊨ ∀w(w ε $x_1$ ↔ w=y), we have that, trivially, M ⊨ $Lex^+$($x_1$).

We have that M ⊨ ¬∃w(w ε $x_1$ & w ε $x^-$) and M ⊨ ∀w(w ε $x^-$ → y≺w).

⟹ M ⊨ ∀u,v (u ε $x_1$ & v ε $x^-$ → u≺v).

We thus have

  M ⊨ Env(t,$x_0$**) & $x_0$**=$t_0$baya$x^-$ & Env(t,$x^-$) & Firstf($x^-$,$t_2$,a$v_0$a,$t_3$) & $t_0$b≺$t_2$ &

& ¬∃w(w ε $x_1$ & w ε $x^-$) & ∀u,v (u ε $x_1$ & v ε $x^-$ → u≺v) & $Lex^+$($x_1$) & $Lex^+$($x^-$).

⟹ from M ⊨ I**($x_0$**), by (9.29), M ⊨ $Lex^+$($x_0$**).

Now, we have that

    M ⊨ ∀u,v (u ε $x^{++}$ & v ε $x^{+*}$ → u≺v) & ∀w(w ε $x^{+*}$ ↔ w ε $x_0$**).

⟹ M ⊨ ∀u,v (u ε $x^{++}$ & v ε $x_0$** → u≺v).

So we have

  M ⊨ Env(t,z) & z=t*$w_0$*$x_0$** & Env($t_1$,$x^{++}$) & $x^{++}$=t*$w_0$*$t_1$ & aB$w_0$* & aE$w_0$* &

   & Env(t,$x_0$**) & $t_1$≺$t_0$b & Firstf($x_0$**,$t_0$b,aya,$t_2$) & ¬∃w(w ε $x^{++}$ & w ε $x_0$**) &

                    & ∀u,v (u ε $x^{++}$ & v ε $x_0$** → u≺v) & $Lex^+$($x^{++}$) & $Lex^+$($x_0$**),

⟹ from M ⊨ I**(z), by (9.29), M ⊨ $Lex^+$(z).

Let $x_0^{++}$=t*$w_0$*$x_1$.

⟹ from M ⊨ I**(t*) & I**($w_0$*) & I**($x_1$), since I** may be assumed to be closed w.r. to *,    M ⊨ I**($x_0^{++}$).

Then we have

$$M \vDash \text{Env}(t_1,x^{++}) \,\&\, x^{++}=t*w_0*t_1 \,\&\, aBw_0* \,\&\, aEw_0* \,\&\, \text{Env}(t_0b,x_1) \,\&$$
$$\&\, \text{Firstf}(x_1,t_0b,aya,t_0b) \,\&\, t_1<t_0b \,\&\, \neg\exists w(w\,\varepsilon\,x^{++} \,\&\, w\,\varepsilon\,x_1).$$

$\Longrightarrow$ from $M \vDash I^{**}(x^{++}) \,\&\, I^{**}(x_1)$, by (5.46),

$$M \vDash \text{Env}(t_0b, x_0^{++}) \,\&\, \forall w(w\,\varepsilon\,x_0^{++} \leftrightarrow w\,\varepsilon\,x^{++} \vee w=y).$$

We then have

$$M \vDash \text{Env}(t_0b,x_0^{++}) \,\&\, x_0^{++}=t*w_0*x_1 \,\&\, \text{Env}(t_1,x^{++}) \,\&\, x^{++}=t*w_0*t_1 \,\&\, aBw_0* \,\&$$
$$\&\, aEw_0* \,\&\, x_1=t_0bayat_0b \,\&\, \text{Max}^+T_b(t_0b,y) \,\&\, \neg(y\,\varepsilon\,x^{++}) \,\&$$
$$\&\, \forall w\,(w\,\varepsilon\,x^{++} \to w\prec y) \,\&\, \text{Special}(x^{++}).$$

$\Longrightarrow$ from $M \vDash I^{**}(x_0^{++})$, by (10.18), $M \vDash \text{Special}(x_0^{++})$.

Let $w_0 = w_0 * t_0 \text{baya}$.

$\Longrightarrow M \vDash x_0^{++} = t*w_0*x_1 = t*w_0*t_0bayat_0b = t*w_0 t_0 b$.

We claim that $M \vDash \neg\exists w(w\,\varepsilon\,x_0^{++} \,\&\, w\,\varepsilon\,x^-)$.

Suppose, for a reductio, that $M \vDash w\,\varepsilon\,x_0^{++} \,\&\, w\,\varepsilon\,x^-$.

$\Longrightarrow M \vDash (w\,\varepsilon\,x^{++} \vee w=y) \vee w\,\varepsilon\,x^-$,

$\Longrightarrow M \vDash (w\,\varepsilon\,x^{++} \vee w\,\varepsilon\,x^-) \vee (w\,\varepsilon\,x_1 \vee w\,\varepsilon\,x^-)$.

But earlier proofs ruled out both of these cases.

Hence $M \vDash \neg\exists w(w\,\varepsilon\,x_0^{++} \,\&\, w\,\varepsilon\,x^-)$, as claimed.

Note we have that

$$M \vDash z = t*w_0*t_0 bayat_2 aw_2 = t*w_0\,x^-.$$

Furthermore, we claim that

$$M \vDash \forall u,v\,(u\,\varepsilon\,x_0^{++} \,\&\, v\,\varepsilon\,x^- \to u\prec v).$$

By earlier proofs, we have that

$$M \vDash \forall w\,(w\;\varepsilon\;x^{++} \to w \prec y)\;\&\;\forall w\,(w\;\varepsilon\;x^{-} \to y \prec w).$$

$\Rightarrow$ from $M \vDash I^{**}(y)$, by (8.3), $M \vDash \forall u,v\,((u\;\varepsilon\;x^{++} \vee u=y)\;\&\;v\;\varepsilon\;x^{-} \to u \prec v)$,

$\Rightarrow M \vDash \forall u,v\,(u\;\varepsilon\;x_0^{++}\;\&\;v\;\varepsilon\;x^{-} \to u \prec v)$, as claimed.

So we have

$M \vDash \text{Env}(t,z)\;\&\;z=t*w_0\;x^{-}\;\&\;\text{Env}(t_0b,x_0^{++})\;\&\;x_0^{++}=t*w_0t_0b\;\&\;aBw_0\;\&$

$\&\;aEw_0\;\&\;\text{Env}(t,x^{-})\;\&\;\text{Firstf}(x^{-},t_2,av_0a,t_3)\;\&\;t_0b<t_2\;\&\;\neg\exists w(w\;\varepsilon\;x_0^{++}\;\&\;w\;\varepsilon\;x^{-})\;\&$

$\&\;\forall u,v\,(u\;\varepsilon\;x_0^{++}\;\&\;v\;\varepsilon\;x^{-} \to u \prec v)\;\&\;\text{Special}(x_0^{++})\;\&\;\text{Special}(x^{-}),$

$\Rightarrow$ from $M \vDash I^{**}(z)$, by (10.17), $M \vDash \text{Special}(z)$.

Hence $M \vDash V^{**}(z)$.

It remains to show that

$$M \vDash \forall w(w\;\varepsilon\;z \leftrightarrow w\;\varepsilon\;x \vee w=y).$$

Now, we have that

$M \vDash \text{Env}(t_0b,x_0^{++})\;\&\;x_0^{++}=t*w_0*x_1\;\&\;\text{Env}(t,x^{-})\;\&\;\text{Firstf}(x^{-},t_2,av_0a,t_3)\;\&$

$\&\;t_0b<t_2\;\&\;\neg\exists w(w\;\varepsilon\;x_0^{++}\;\&\;w\;\varepsilon\;x^{-}).$

$\Rightarrow$ from $M \vDash I^{**}(x_0^{++})\;\&\;I^{**}(x^{-})$, by (5.46),

$$M \vDash \forall w(w\;\varepsilon\;z \leftrightarrow w\;\varepsilon\;x_0^{++} \vee w\;\varepsilon\;x^{-}).$$

Assume that $M \vDash w\;\varepsilon\;z$.

$\Rightarrow M \vDash w\;\varepsilon\;x_0^{++} \vee w\;\varepsilon\;x^{-}$.

If $M \vDash w\;\varepsilon\;x^{-}$, then $M \vDash w\;\varepsilon\;x"$, whence $M \vDash w\;\varepsilon\;x$.

If $M \vDash w\;\varepsilon\;x_0^{++}$, then $M \vDash w\;\varepsilon\;x^{++} \vee w=y$. But $M \vDash x^{++}=x^*$, and

$M \vDash \forall u(u\;\varepsilon\;x^* \to u\;\varepsilon\;x)$. Hence $M \vDash w\;\varepsilon\;x \vee w=y$.

Therefore, $M \vDash \forall w(w \, \varepsilon \, z \to w \, \varepsilon \, x \vee w=y)$.

Conversely, assume $M \vDash w \, \varepsilon \, x \vee w=y$.

If $M \vDash w=y$, then $M \vDash w \, \varepsilon \, x_0^{++}$, so $M \vDash w \, \varepsilon \, z$.

On the other hand, if $M \vDash w \, \varepsilon \, x$, then $M \vDash w \, \varepsilon \, x' \vee w \, \varepsilon \, x''$.

But then, from $M \vDash w \, \varepsilon \, x'$ we have $M \vDash w \, \varepsilon \, x^*$, whence $M \vDash w \, \varepsilon \, x^{++}$,

and further $M \vDash w \, \varepsilon \, x_0^{++}$. So $M \vDash w \, \varepsilon \, z$.

On the other hand, from $M \vDash w \, \varepsilon \, x''$ & $w \neq u_0$, we have $M \vDash w \, \varepsilon \, x^-$, whence

$M \vDash w \, \varepsilon \, z$, whereas from $M \vDash w \, \varepsilon \, x''$ & $w=u_0$, we obtain $M \vDash w \, \varepsilon \, x_0$, whence

$M \vDash w \, \varepsilon \, x^*$, and also $M \vDash w \, \varepsilon \, z$.

Thus also $M \vDash \forall w(w \, \varepsilon \, x \vee w=y \to w \, \varepsilon \, z)$.

(1c)  $M \vDash t_0 < t_1$.

  (1ci) $M \vDash t_1 < t^{++}$.

We have that

  $M \vDash \mathrm{Env}(t,x)$ & $x = t^* w_0^* x^-$ & $aBw_0^*$ & $aEw_0^*$ & $\mathrm{Env}(t_1, x^{++})$ & $x^{++} = t^* w_0^* t_1$ &

      & $\mathrm{Env}(t, x^-)$ & $\mathrm{Lastf}(x^{++}, t_1, au_0a, t_1)$ & $\mathrm{Firstf}(x^-, t_2, av_0a, t_3)$ &

          & $\mathrm{MinMax}^+ T_b(t^{++}, v_0)$ & $t_1 < t^{++}$ & $\mathrm{Special}(x)$.

$\Rightarrow$ by (10.25), $M \vDash t^{++} = t_2$.

  (1cia) $M \vDash t_1 b < t_2$.

Let $z = w_1 a t_1 a u_0 a t_1 b a y a t_2 a w_2$.

Exactly analogous to (1bii) with $t_1 b$ in place of $t_0 b$ throughout the argument.

  (1cib) $M \vDash t_1 b = t_2$.

Let $x_2 = t_1 a y a t_1$.

$\Rightarrow$ from $M \vDash I^{**}(t_1)$ & $I^{**}(y)$, since $I^{**}$ may be assumed to be closed w.r. to *,

$$M \vDash I^{**}(x_2),$$

$\Rightarrow$ from $M \vDash MinMax^+T_b(t_0,aya)$ & $t_0<t_1$, $M \vDash Max^+T_b(t_1,aya)$,

$\Rightarrow$ $M \vDash Pref(aya,t_1)$,

$\Rightarrow$ $M \vDash Firstf(x_2,t_1,aya,t_1)$ & $Lastf(x_2,t_1,aya,t_1)$,

$\Rightarrow$ from $M \vDash I^{**}(x_2)$, by (5.22), $M \vDash Env(t_1,x_2)$ & $\forall w(w \, \varepsilon \, x_2 \leftrightarrow w=y)$,

$\Rightarrow$ by (10.5), $M \vDash MinSet(x_2)$,

$\Rightarrow$ by (10.11), $M \vDash MinSet(x^-)$.

Then we have

$M \vDash Env(t_1,x_2)$ & $x_2=t_1ayat_1$ & $x_1^+=t_1ayax^-$ & $Env(t,x^-)$ & $t_1<t_2$ &

& $Firstf(x^-,t_2,av_0a,t_3)$ & $\neg \exists w(w \, \varepsilon \, x_2 \, \& \, w \, \varepsilon \, x^-)$ & $MinSet(x_2)$ & $MinSet(x^-)$.

$\Rightarrow$ from $M \vDash I^{**}(t_1)$ & $I^{**}(y)$ & $I^{**}(x^-)$, since $I^{**}$ may be assumed to be closed w.r. to *, $M \vDash I^{**}(x_1^+)$,

$\Rightarrow$ by (10.6), $M \vDash MinSet(x_1^+)$ & $Env(t_1, x_1^+)$.

Then, as in (1bi), we have, since $M \vDash J^{**}(t)$,

$M \vDash \exists!x_1^{+*}\in I^{**} \, \exists t'\in J^{**} \, (Env(t',x_1^{+*})$ & $x_1^+ \sim x_1^{+*}$ &

& $MinSet(x_1^{+*})$ & $\forall u,v \, (u\leq_{x1+}v \leftrightarrow u\leq_{x1+*}v)$ &

& $\forall v,t_3,t_4 \, (Free^+(x_1^+,t_3,ava,t_4) \rightarrow \exists t_5 Fr(x_1^{+*},t_3b,ava,t_5))$ &

& $\forall v,t_3,t_4 \, (Bound(x_1^+,t_3,ava,t_4) \vee Free^-(x_1^+,t_3,ava,t_4) \rightarrow \exists t_5 Fr(x_1^{+*},t_3,ava,t_5)))$.

Let $z=w_1at_1au_0ax_1^{+*}$.

Then, exactly analogously to (1bi), we show that

$M \vDash I^{**}(z)$ & $MinSet(z)$ & $Env(t',z)$.

Exactly analogously to (1bi), with $x_2$ in place of $x_0{}^*$, we show that

$$M \vDash Lex^+(x_1{}^+),$$

and then that $M \vDash Lex^+(x_1{}^{+*})$.

Then, again analogously to (1bi), we derive $M \vDash Lex^+(z)$.

Finally, the same argument as in (1bi), with $x_1{}^+$ in place of $x^+$, is used to derive $M \vDash Special(z)$.

Therefore, $M \vDash V^{**}(z)$.

Exactly analogously to (1bi) we also have that $M \vDash \forall w(w\ \varepsilon\ z \rightarrow w\ \varepsilon\ x \lor w = y)$.

  (1cii) $M \vDash t_1 = t^{++}$.

$\Rightarrow M \vDash t^{++} = t_1 < t_2$.

   (1ciia) $M \vDash t_1 b < t_2$.

Same as (1cia).

   (1ciib) $M \vDash t_1 b = t_2$.

Same as (1cib).

    (*Biib2.2): $M \vDash u_0 \triangleleft_{Tb} y\ \&\ (y \approx_{Tb} v_0\ \&\ y \ll v_0)$.

$\Rightarrow M \vDash t^+ < t_0 < t^{++}$.

(2a)   $M \vDash t_1 < t_0$.

 (2ai) $M \vDash t_1 b = t_0$.

  (2aia) $M \vDash t_0 = t_2$.

We have that $M \vDash MinSet(x^-)$.

$\Rightarrow$ from $M \vDash I^{**}(x^-)$, by (10.29), $M \vDash \exists! x^{-*} \in I^* \ \exists t' \in J^{**}(Env(t', x^{-*})\ \&\ x^- \sim x^{-*}\ \&$

      $\&\ MinSet(x^{-*})\ \&\ \forall u, v\ (u \leq_{x^-} v \leftrightarrow u \leq_{x^{-*}} v)\ \&\ Env(tb, x^{-*})\ \&$

$$\& \ \forall v, t_3, t_4 \ (\text{Free}^+(x^-, t_3, ava, t_4) \to \exists t_5 \text{Fr}(x^{-*}, t_3 b, ava, t_5)) \ \&$$

$$\& \ \forall v, t_3, t_4 \ (\text{Bound}(x^-, t_3, ava, t_4) \ v \ \text{Free}^-(x^-, t_3, ava, t_4) \to \exists t_5 \text{Fr}(x^{-*}, t_3, ava, t_5))).$$

Let $z = w_1 a t_1 a u_0 a t_0 a y a x^{-*}$.

$\implies$ from $M \vDash I^{**}(w_1) \& I^{**}(t_1) \& I^{**}(u_0) \& I^{**}(t_0) \& I^{**}(y) \& I^{**}(x^{-*})$, since $I^{**}$ may be assumed to be closed w.r. to $*$, $\quad M \vDash I^{**}(z)$,

$\implies$ from $M \vDash \text{MinMax}^+ T_b(t_0, aya)$, $M \vDash \text{Pref}(aya, t_0)$,

$\implies$ from $M \vDash I^{**}(x_0^*)$, by (10.5), $M \vDash \text{MinSet}(x_0^*)$, where $x_0^* = t_0 ayat_0$.

So we have

$$M \vDash \text{Env}(t_1, x^{++}) \& x^{++} = t^* w_0^* t_1 \& aBw_0^* \& aEw_0^* \& x^{**} = t^* w_0^* x_0^* \&$$
$$\& \text{Env}(t_0, x_0^*) \& t_1 < t_0 \& \text{Firstf}(x_0^*, t_0, aya, t_0) \& \neg \exists w (w \ \varepsilon \ x^{++} \ \& \ w \ \varepsilon \ x_0^{**}) \&$$
$$\& \text{MinSet}(x^{**}) \& \text{MinSet}(x_0^*).$$

$\implies$ from $M \vDash I^{**}(t^*) \& I^{**}(w_0^*) \& I^{**}(x_0^*)$, $M \vDash I^{**}(x^{**})$,

$\implies$ by (10.6), $M \vDash \text{MinSet}(x^{**}) \& \text{Env}(t_0, x^{**})$.

From $M \vDash \text{Firstf}(x^-, t_2, av_0 a, t_3)$ we have that $M \vDash \exists t_5 \text{Firstf}(x^{-*}, t_2, av_0 a, t_5)$

as in (1bi).

By an earlier proof, $M \vDash \neg \exists w (w \ \varepsilon \ x^{**} \ \& \ w \ \varepsilon \ x^-)$.

$\implies M \vDash \neg \exists w (w \ \varepsilon \ x^{**} \ \& \ w \ \varepsilon \ x^{-*})$.

So we have

$$M \vDash \text{Env}(t_0, x^{**}) \& x^{**} = t^* w^{**} t_0 \& aBw^{**} \& aEw^{**} \& z = t^* w^{**} x^{-*} \&$$
$$\& \text{Env}(tb, x^{-*}) \& t_0 < t_2 b \& \text{Firstf}(x^{-*}, t_2 b, av_0 a, t_5) \& \neg \exists w (w \ \varepsilon \ x^{**} \ \& \ w \ \varepsilon \ x^{-*}) \&$$
$$\& \text{MinSet}(x^{**}) \& \text{MinSet}(x^{-*}).$$

$\implies$ from $M \vDash I^{**}(z)$, by (10.6), $M \vDash \text{MinSet}(z) \& \text{Env}(t', z)$.

Now, we have $M \vDash Lex^+(x^{++})$ & $Lex^+(x_0^{**})$.

$\Rightarrow$ from $M \vDash I^{**}(x^{**})$, by (9.29), $M \vDash Lex^+(x^{**})$,

$\Rightarrow$ from $M \vDash Lex^+(x^-)$, $M \vDash Lex^+(x^{-*})$,

$\Rightarrow$ since $M \vDash \forall w',w''(w' \varepsilon x^{**}$ & $w'' \varepsilon x^- \to w' \prec w'')$ & $x^- \sim x^{-*}$,

from $M \vDash I^{**}(z)$, by (9.29), $\qquad M \vDash Lex^+(z)$.

Now, we also have that

$M \vDash Env(t',z)$ & $z=t^*w^{**}x^{-*}$ & $aBw^{**}$ & $aEw^{**}$ & $Set(x^{-*})$ & $x^{**}=t^*w^{**}t_0$ &

$\qquad$ & $t_0=t_2$ & $Env(t_0,x^{**})$ & $Firstf(x^{-*},t_2b,av_0a,t_5)$ &

& $\forall u,v$ $(u \varepsilon x^{**}$ & $v \varepsilon x^{-*} \to u \prec v)$ & $\neg\exists u(u \varepsilon x^{**}$ & $u \varepsilon x^{-*})$ & $Set(x^-)$ &

$\qquad$ & $x^- \sim x^{-*}$ & $\forall u,v(u\prec_{x\text{-}}v \leftrightarrow u\prec_{x\text{-}*}v)$ &

$\qquad$ & $\forall v,t_6t_7$ $(Free^+(x^-,t_6,ava,t_7) \to \exists t_8 Fr(x^{-*},t_6b,ava,t_8))$ &

& $\forall v,t_6,t_7$ $(Bound(x^-,t_6,ava,t_7) \vee Free^-(x^-,t_6,ava,t_7) \to \exists t_8 Fr(x^{-*},t_6,ava,t_8))$ &

$\qquad\qquad\qquad\qquad\qquad\qquad\qquad\qquad$ & $Special(x^{**})$ & $Special(x^-)$.

$\Rightarrow$ from $M \vDash I^{**}(x^-)$ & $I^{**}(z)$, by (10.28), $M \vDash Special(z)$.

Therefore, $M \vDash V^{**}(z)$.

Note that $M \vDash \forall w(w \varepsilon z \leftrightarrow w \varepsilon x^{**} \vee w \varepsilon x^{-*})$.

This is proved just as in (1a), using the hypothesis $M \vDash x^- \sim x^{-*}$.

Then an exactly analogous proof to the one given in (1a), with $w \varepsilon x^{-*}$

replacing $w \varepsilon x^-$, shows that

$$M \vDash \forall w(w \varepsilon z \leftrightarrow w \varepsilon x \vee w=y),$$

as required.

(2aib) $M \vDash t_0 < t_2$.

$\Rightarrow$ $M \vDash t_1 < t_0 = t^{++} < t_2$.

This is ruled out by (10.25) as in (1ci) since $M \vDash t_2 \in I \subseteq I_0$.

(2aic) $M \vDash t_2 < t_0$.

$\Rightarrow$ $M \vDash t_0 = t^{++} \leq t_2 < t_0$, contradicting $M \vDash t_0 \in I \subseteq I_0$.

(2aii) $M \vDash t_1 b < t_0$.

$\Rightarrow$ $M \vDash t_1 < t_1 b < t_0 = t^{++}$,

$\Rightarrow$ as in (1ci), $M \vDash t^{++} = t_2$.

We then proceed exactly as in (2aia).

(2aiii) $M \vDash t_0 < t_1 b$.

$\Rightarrow$ $M \vDash t_0 \leq t_1$, by (1.13).

But this, along with (2a), contradicts $M \vDash t_0 \in I \subseteq I_0$.

(2b) $M \vDash t_0 \leq t_1$.

$\Rightarrow$ $M \vDash t^{++} = t_0 \leq t_1 < t_2$.

(2bi) $M \vDash t_1 b = t_2$.

We obtain $x^{-*}$ from $x^-$ just as in (2aia).

Let $z = w_1 a t_1 a u_0 a t_1 b a y a x^{-*}$.

$\Rightarrow$ from $M \vDash I^{**}(w_1)$ & $I^{**}(t_1)$ & $I^{**}(u_0)$ & $I^{**}(y)$ & $I^{**}(x^{-*})$, since $I^{**}$

may be assumed to be closed w.r. to *, $M \vDash I^{**}(z)$.

Let $x_3 = t_1 b a y a t_1 b$.

Then $M \vDash I^{**}(x_3)$.

$\Rightarrow$ from $M \vDash \text{MinMax}^+ T_b(t_0, aya)$ & $t_0 \leq t_1 < t_1 b$, $M \vDash \text{Pref}(aya, t_1 b)$,

$\Rightarrow$ M ⊨ Firstf($x_3$,$t_1$b,aya,$t_1$b) & Lastf($x_3$,$t_1$b,aya,$t_1$b),

$\Rightarrow$ by (5.22), M ⊨ Env($t_1$b,$x_3$) & $\forall w(w \, \varepsilon \, x_3 \leftrightarrow w=y)$,

$\Rightarrow$ by (10.5), M ⊨ MinSet($x_3$).

We then have, for $x_3^{++}=t*w_0*x_3$, that

$\quad$ M ⊨ Env($t_1$,$x^{++}$) & $x^{++}=t*w_0*t_1$ & aB$w_0$* & aE$w_0$* & $x_3^{++}=t*w_0*x_3$ &

$\quad\quad$ & Env($t_1$b,$x_3$) & $t_1<t_1$b & Firstf($x_3$,$t_1$b,aya,$t_1$b) & $\neg\exists w(w \, \varepsilon \, x^{++} \, \& \, w \, \varepsilon \, x_3)$ &

$\quad\quad\quad\quad\quad\quad\quad\quad\quad\quad\quad\quad\quad\quad\quad$ & MinSet($x^{++}$) & MinSet($x_3$).

From M ⊨ I**(t*) & I**($w_0$*) & I**($x_3$), since I** may be assumed to be closed

w.r. to *, $\quad$ M ⊨ I**($x_3^{++}$).

$\Rightarrow$ by (10.6), M ⊨ MinSet($x_3^{++}$) & Env($t_1$b,$x_3^{++}$).

We then have that

$\quad$ M ⊨ $\forall w(w \, \varepsilon \, x_3^{++} \leftrightarrow w \, \varepsilon \, x^{++} \vee w \, \varepsilon \, x_3 \leftrightarrow (w \, \varepsilon \, x^{++} \vee w=y) \leftrightarrow w \, \varepsilon \, x^{**})$.

Also, M ⊨ $x_3^{++}=t*w_0*x_3=t*w_0*t_1$bayа$t_1$b.

Let $w_3*=w_0*t_1$baya.

So we have

$\quad$ M ⊨ Env($t_1$b,$x_3^{++}$) & $x_3^{++}=t*w_3*t_1$b & aB$w_3$* & aE$w_3$* & $z=t*w_3*x^{-*}$ &

$\quad$ & Env(tb,$x^{-*}$) & $t_1$b$<t_2$b & Firstf($x^{-*}$,$t_2$b,a$v_0$a,$t_5$) & $\neg\exists w(w \, \varepsilon \, x_3^{++} \, \& \, w \, \varepsilon \, x^{-*})$ &

$\quad\quad\quad\quad\quad\quad\quad\quad\quad\quad\quad\quad\quad\quad\quad$ & MinSet($x_3^{++}$) & MinSet($x^{-*}$).

$\Rightarrow$ from M ⊨ I**(z), by (10.6), M ⊨ MinSet(z) & Env(t',z).

Now, we have M ⊨ Lex$^+$($x^{++}$), and, trivially, M ⊨ Lex$^+$($x_3$).

$\Rightarrow$ from M ⊨ I**($x_3^{++}$), by (9.29), M ⊨ Lex$^+$($x_3^{++}$).

Then, just as in (2aia), we obtain  $M \vDash Lex^+(z)$.

Now, we have

$M \vDash Env(t_1b, x_3^{++})$ & $x_3^{++}=t*w_0*x_3$ & $Env(t_1, x^{++})$ & $x^{++}=t*w_0*t_1$ & $aBw_0^*$ &

& $aEw_0^*$ & $x_3 = t_1bayat_1b$ & $Max^+T_b(t_1b, y)$ & $\neg(y \,\varepsilon\, x^{++})$ &

& $\forall u(u \,\varepsilon\, x^{++} \to u \prec v)$ & $Special(x^{++})$.

$\Rightarrow$ from  $M \vDash I^{**}(x_3^{++})$,  by (10.18),  $M \vDash Special(x_3^{++})$.

From  $M \vDash \forall w(w \,\varepsilon\, x_3^{++} \leftrightarrow w \,\varepsilon\, x^{**})$  and the previously established

$M \vDash \forall u,v\, (u \,\varepsilon\, x^{**}$ & $v \,\varepsilon\, x^{-*} \to u \prec v)$,

it follows that   $M \vDash \forall u,v\, (u \,\varepsilon\, x_3^{++}$ & $v \,\varepsilon\, x^{-*} \to u \prec v)$.

We then have

$M \vDash Env(t', z)$ & $z = t*w_3*x^{-*}$ & $aBw_3^*$ & $aEw_3^*$ & $Set(x^{-*})$ & $x_3^{++} = t*w_3*t_1b$ &

& $t_1b < t_2b$ & $Env(t_1b, x_3^{++})$ & $Firstf(x^{-*}, t_2b, av_0a, t_5)$ &

& $\forall u,v\, (u \,\varepsilon\, x_3^{++}$ & $v \,\varepsilon\, x^{-*} \to u \prec v)$ & $\neg\exists u(u \,\varepsilon\, x_3^{++}$ & $u \,\varepsilon\, x^{-*})$ & $Set(x^-)$ &

& $x^- \sim x^{-*}$ & $\forall u,v(u <_{x^-} v \leftrightarrow u <_{x^{-*}} v)$ &

& $\forall v, t_6 t_7\, (Free^+(x^-, t_6, ava, t_7) \to \exists t_8 Fr(x^{-*}, t_6b, ava, t_8))$ &

& $\forall v, t_6 t_7\, (Bound(x^-, t_6, ava, t_7)$ v $Free^-(x^-, t_6, ava, t_7) \to \exists t_8 Fr(x^{-*}, t_6, ava, t_8))$ &

& $Special(x_3^{++})$ & $Special(x^-)$.

$\Rightarrow$  from $M \vDash I^{**}(x^-)$ & $I^{**}(z)$,  by (10.28),  $M \vDash Special(z)$.

Hence   $M \vDash V^{**}(z)$.

Finally, to show that   $M \vDash \forall w(w \,\varepsilon\, z \leftrightarrow w \,\varepsilon\, x$ v $w = y)$,

note first that we have

$M \vDash Env(t_1b, x_3^{++})$ & $x_3^{++} = t*w_0*t_1b$ & $Env(t', x^{-*})$ & $Firstf(x^{-*}, t_2b, av_0a, t_5)$ &

$$\& \ t_1b<t_2b \ \& \ \neg\exists w(w \ \varepsilon \ x_3^{++} \ \& \ w \ \varepsilon \ x^{-*}).$$

$\Rightarrow$ from $M \vDash I^{**}(x_3^{++}) \ \& \ I^{**}(x^{-})$, by (5.46),

$$M \vDash \forall w(w \ \varepsilon \ z \leftrightarrow w \ \varepsilon \ x_3^{++} \ v \ w \ \varepsilon \ x^{-*}).$$

From (2aia) we have that

$$M \vDash \forall w(w \ \varepsilon \ x^{**} \ \& \ w \ \varepsilon \ x^{-} \rightarrow w \ \varepsilon \ x \ v \ w=y).$$

Since $M \vDash \forall w(w \ \varepsilon \ x_3^{++} \leftrightarrow w \ \varepsilon \ x^{**}) \ \& \ x^{-} \sim x^{-*}$, it follows that

$$M \vDash \forall w(w \ \varepsilon \ z \rightarrow w \ \varepsilon \ x \ v \ w=y).$$

Conversely, assume that $M \vDash w \ \varepsilon \ x \ v \ w=y$.

If $M \vDash w=y$, then $M \vDash w \ \varepsilon \ x_0^{*}$, whence $M \vDash w \ \varepsilon \ x^{**}$. But then $M \vDash w \ \varepsilon \ x_3^{++}$, so $M \vDash w \ \varepsilon \ z$.

If $M \vDash w \ \varepsilon \ x$, then $M \vDash w \ \varepsilon \ x' \ v \ w \ \varepsilon \ x''$. Then:

$M \vDash w \ \varepsilon \ x' \Rightarrow M \vDash w \ \varepsilon \ x^{*}, \Rightarrow M \vDash w \ \varepsilon \ x^{++}, \Rightarrow M \vDash w \ \varepsilon \ x_3^{++}, \Rightarrow M \vDash w \ \varepsilon \ z.$

On the other hand,

$M \vDash w \ \varepsilon \ x'' \ \& \ w \neq u_0 \Rightarrow M \vDash w \ \varepsilon \ x^{-}, \Rightarrow M \vDash w \ \varepsilon \ x^{-*}, \Rightarrow M \vDash w \ \varepsilon \ z,$

whereas

$M \vDash w \ \varepsilon \ x'' \ \& \ w=u_0 \Rightarrow M \vDash w \ \varepsilon \ x_0, \Rightarrow M \vDash w \ \varepsilon \ x^{*}, \Rightarrow M \vDash w \ \varepsilon \ x^{++},$

$$\Rightarrow M \vDash w \ \varepsilon \ x_3^{++}, \Rightarrow M \vDash w \ \varepsilon \ z.$$

Therefore, we also have $M \vDash w \ \varepsilon \ x \ v \ w=y \rightarrow w \ \varepsilon \ z$.

(2bii) $M \vDash t_1b<t_2$.

$\Rightarrow M \vDash t_0=t^{++}=t_1<t_1b<t_2.$

We also have

$M \vDash \text{Env}(t,x) \ \& \ x=t^*w_0^*x^{-} \ \& \ aBw_0^* \ \& \ aEw_0^* \ \& \ \text{Env}(t_1,x^{++}) \ \& \ x^{++}=t^*w_0^*t_1 \ \&$

& Env(t,x$^-$) & Lastf(x$^{++}$,t$_1$,au$_0$a,t$_1$) & Firstf(x$^-$,t$_2$,av$_0$a,t$_3$),

and M ⊨ x=t*w$_0$*x$^-$=t*w$_0$*t$_1$t"t$_2$'w$^-$t where M ⊨ t$_1$t"t$_2$'=t$_2$ & aBw$^-$ & aEw$^-$.

$\Rightarrow$ by (5.41), M ⊨ ¬∃w(w ε x$^{++}$ & w ε x$^-$),

$\Rightarrow$ from M ⊨ I**(x$^{++}$) & I**(x$^-$), by (5.46), M ⊨ ∀w(w ε x ↔ w ε x$^{++}$ v w ε x$^-$).

We now establish Claims 1-3 exactly as in the proof of (10.25), and derive,

from hypothesis M ⊨ t$_1$<t$_1$b<t$_2$, that M ⊨ Max$^+$(t$_1$b,v$_0$,x).

$\Rightarrow$ from M ⊨ Max$^+$T$_b$(t$^{++}$,v$_0$) & t$^{++}$<t$_1$b, M ⊨ Max$^+$T$_b$(t$_1$b,v$_0$).

But from M ⊨ Special(x) & Fr(x,t$_2$,av$_0$a,t$_3$) we have M ⊨ MMax$^+$(t$_2$,v$_0$,x).

$\Rightarrow$ M ⊨ t$_2$≤t$_1$b<t$_2$, which contradicts M ⊨ t$_2$∈I⊆I$_0$.

Hence (2bii) is ruled out.

  (2biii) M ⊨ t$_2$<t$_1$b.

$\Rightarrow$ by (1.13), M ⊨ t$_2$≤t$_1$,

$\Rightarrow$ M ⊨ t$_1$<t$_2$≤t$_1$, contradicting M ⊨ t$_1$∈I⊆I$_0$.

So (2biii) is also ruled out.

   (*Biib2.3): M ⊨ (u$_0$≈$_{Tb}$ y & u$_0$≪y) & y◁$_{Tb}$ v$_0$.

$\Rightarrow$ M ⊨ t$^+$=t$_0$<t$^{++}$,

$\Rightarrow$ M ⊨ t$_0$=t$^+$≤t$_1$.

(3a) M ⊨ t$_0$=t$_1$.

$\Rightarrow$ M ⊨ t$_1$=t$_0$<t$^{++}$.

Exactly as (1b).

(3b) M ⊨ t$_0$<t$_1$.

Exactly as (1c).

(*Biib2.4): $M \vDash (u_0 \approx_{Tb} y\ \&\ u_0 \ll y)\ \&\ (y \approx_{Tb} v_0\ \&\ y \ll v_0)$.

$\Rightarrow M \vDash t^+ = t_0 = t^{++}$,

$\Rightarrow M \vDash t_0 = t^+ \leq t_1\ \&\ t_0 = t^{++} \leq t_2$.

(4a)　$M \vDash t_0 = t_1$.

$\Rightarrow M \vDash t_1 = t_0 = t^{++} \leq t_2$.

　(4ai)　$M \vDash t^{++} = t_2$.

$\Rightarrow M \vDash t_1 = t_0 = t^{++} = t_2$,

$\Rightarrow M \vDash t_1 < t_2 = t_1$,　contradicting　$M \vDash t_1 \in I \subseteq I_0$.

So (4ai) is ruled out.

　(4aii)　$M \vDash t^{++} < t_2$.

Exactly as (2b).

(4b)　$M \vDash t_0 < t_1$.

$\Rightarrow M \vDash t^{++} < t_0 < t_1 < t_2$.

Exactly as in (2bi) and (2bii).

　　(*Biib1):　$M \vDash \text{Firstf}(x, t_1, au_0a, t_2)$.

$\Rightarrow$ by (9.1), $M \vDash \forall w(w \leq_x u_0 \rightarrow w = u_0)$,

$\Rightarrow M \vDash \text{Pref}(au_0a, t_1)\ \&\ \text{Tally}_b(t_2)\ \&\ ((t_1 = t_2\ \&\ x = t_1 au_0 a t_2)\ \text{v}$

$\hspace{6cm}\text{v}\ (t_1 < t_2\ \&\ (t_1 au_0 a t_2) B x))$.

(iib1i):　$M \vDash t_1 = t_2\ \&\ x = t_1 au_0 a t_2$.

$\Rightarrow M \vDash \text{Lastf}(x, t_1, au_0 a, t_2)$,

$\Rightarrow$ from $M \vDash I^{**}(x)$, by (5.22), $M \vDash \forall w(w\ \varepsilon\ x \rightarrow w = u_0)$,

$\Rightarrow$ since $M \vDash u_0 \ll y$, $M \vDash \forall w(w\ \varepsilon\ x \rightarrow w \ll y)$.

But by (*Bii), $M \vDash \exists u(u \,\varepsilon\, x \,\&\, y \prec u)$.

$\Rightarrow M \vDash u \prec y \,\&\, y \prec u$, contradicting $M \vDash I^{**}(y)$ by (8.2).

(iib1ii): $M \vDash t_1 < t_2 \,\&\, (t_1 a u_0 a t_2) Bx$.

$\Rightarrow M \vDash \exists w_2 \, t_1 a u_0 a t_2 a w_2 = x$,

$\Rightarrow$ from $M \vDash I^{**}(x)$ and the proof of (6.8), $M \vDash \exists! t^+ \in I^{**} \, MinMax^+ T_b(t^+, u_0)$,

$\Rightarrow$ from hypothesis $M \vDash Set^*(x)$, $M \vDash Special(x)$,

$\Rightarrow$ by (10.16), $M \vDash MinMax^+ T_b(t_1, u_0)$,

$\Rightarrow M \vDash t_1 = t^+$.

We now follow the general pattern of proof in (*Biib2):

(1) $M \vDash u_0 \triangleleft_{T_b} y \,\&\, y \triangleleft_{T_b} v_0$.

$\Rightarrow M \vDash t_1 = t^+ < t_0 < t^{++}$.

Let $z = t_1 a u_0 a t_0 a y a t_2 a w_2$.

Since $I^{**}$ may be assumed to be closed w.r. to *,

$$M \vDash I^{**}(z).$$

Let $x_0 = t_1 a u_0 a t_1$.

Then, for the same reason, $M \vDash I^{**}(x_0)$.

Then $M \vDash MinSet(x_0)$ is proved as in (10.5).

We prove, as in (*Biib2)(1), that $M \vDash MinSet(x^+)$, where $x^+ = t_0 a y a t_2 a w_2$, by replacing x" with x throughout the argument there.

By earlier proofs, we have that

$M \vDash Env(t_1, x_0) \,\&\, x_0 = t_1 a u_0 a t_1 \,\&\, z = t_1 a u_0 a x^+ \,\&\, Env(t, x^+) \,\&\, t_1 < t_0 \,\&$

$\&\, Firstf(x^+, t_0, aya, t_2) \,\&\, \neg\exists w(w \,\varepsilon\, x_0 \,\&\, w \,\varepsilon\, x^+) \,\&\, MinSet(x_0) \,\&\, MinSet(x^+)$.

$\Rightarrow$ from $M \vDash I^{**}(z)$, by (10.6), $M \vDash \text{MinSet}(z)\ \&\ \text{Env}(t,z)$.

Now, as in (*Biib2)(1a), we have that $M \vDash \text{Lex}^+(x_0)\ \&\ \text{Lex}^+(x^+)$.

Also, from $M \vDash I^{**}(x_0)\ \&\ I^{**}(x^+)$, by (5.46), $M \vDash \forall w(w\ \varepsilon\ z \leftrightarrow w\ \varepsilon\ x_0 \vee w\ \varepsilon\ x^+)$,

whence $M \vDash \forall w(w\ \varepsilon\ z \leftrightarrow w=u_0 \vee w\ \varepsilon\ x^+)$.

We also have $M \vDash \forall w(w\ \varepsilon\ x_0 \leftrightarrow w=u_0)\ \&\ u_0 \prec y$,

and, as in (*Biib2)(1a), that $M \vDash \forall w(w\ \varepsilon\ x^+ \leftrightarrow (w\ \varepsilon\ x \vee w \neq u_0) \vee w=y)$

and $M \vDash \forall w(w\ \varepsilon\ x \vee w \neq u_0 \rightarrow y \prec w)$.

Assume that $M \vDash u\ \varepsilon\ x_0\ \&\ v\ \varepsilon\ x^+$.

$\Rightarrow M \vDash u=u_0\ \&\ ((v\ \varepsilon\ x \vee v \neq u_0) \vee v=y)$,

$\Rightarrow M \vDash u \prec y\ \&\ (y \prec v \vee v=y)$,

$\Rightarrow$ from $M \vDash I^{**}(y)$, by (8.3), $M \vDash u \prec v$.

Therefore, $M \vDash \forall u,v(u\ \varepsilon\ x_0\ \&\ v\ \varepsilon\ x^+ \rightarrow u \prec v)$.

$\Rightarrow$ from $M \vDash I^{**}(z)$, by (9.29), $M \vDash \text{Lex}^+(z)$.

Now, we have that

$\qquad M \vDash \text{Env}(t,x)\ \&\ x=t_1 a u_0 a x^-\ \&\ \text{Env}(t_1,x_0)\ \&\ \text{Lastf}(x_0,t_1,a u_0 a,t_1)\ \&$

$\qquad\qquad \&\ \text{Firstf}(x^-,t_2,a v_0 a,t_3)\ \&\ \text{MinMax}^+ T_b(t^{++},v_0)\ \&\ t_1 \prec t^{++}\ \&\ \text{Special}(x)$.

$\Rightarrow$ by (10.25), $M \vDash t_2 = t^{++}$,

$\Rightarrow$ by (10.24), $M \vDash \text{Special}(x^-)$,

$\Rightarrow$ from $M \vDash I^{**}(x^+)$, by (10.20), $M \vDash \text{Special}(x^+)$,

$\Rightarrow$ from $M \vDash I^{**}(u_0)$, by (10.15), $M \vDash \text{Special}(x_0)$,

$\Rightarrow$ from $M \vDash I^{**}(z)$, by (10.17), $M \vDash \text{Special}(z)$.

Since $M \vDash V^{**}(x)$ by hypothesis, we have from $M \vDash \text{Env}(t,x)$ that $M \vDash J^{**}(t)$.

Hence $M \vDash V^{**}(z)$ follows from $M \vDash Env(t,z)$.

Then we have from

$M \vDash \forall w(w \ \varepsilon \ z \leftrightarrow w=u_0 \ v \ w \ \varepsilon \ x^+) \ \& \ \forall w(w \ \varepsilon \ x^+ \leftrightarrow (w \ \varepsilon \ x \ v \ w \neq u_0) \ v \ w=y)$

and $M \vDash u_0 \ \varepsilon \ x$

that $M \vDash \forall w(w \ \varepsilon \ z \leftrightarrow w \ \varepsilon \ x \ v \ w=y)$, as required.

(2) $M \vDash u_0 \triangleleft_{Tb} y \ \& \ (y \approx_{Tb} v_0 \ \& \ y \ll v_0)$.

$\Rightarrow M \vDash t_1=t^+ < t_0=t^{++}$.

(2ai) $M \vDash t_1 b = t_0$.

(2aia) $M \vDash t_0 = t_2$.

We follow the proof of (2aia) in (*Biib2) and obtain $x^{-*}$ from $x^-$ where

$M \vDash I^{**}(x^{-*}) \ \& \ MinSet(x^{-*}) \ \& \ \exists t' \in J^{**} \ Env(t', x^{-*})$.

Let $z = t_1 a u_0 a t_0 a y a x^{-*}$.

$\Rightarrow$ from $M \vDash I^{**}(t_1) \ \& \ I^{**}(u_0) \ \& \ I^{**}(t_0) \ \& \ I^{**}(y) \ \& \ I^{**}(x^{-*})$, since $I^{**}$ may be assumed to be closed w.r. to *, $M \vDash I^{**}(z)$.

For $x_0 = t_1 a u_0 a t_1$ and $x_0^* = t_0 a y a t_0$ we have

$M \vDash Pref(au_0 a, t_1) \ \& \ Pref(aya, t_0) \ \& \ t_1 < t_0 \ \& \ u_0 \neq y$.

$\Rightarrow$ from $M \vDash I^{**}(x_1^*)$, by (10.7), $M \vDash MinSet(x_1^*) \ \& \ Env(t_0, x_1^*)$

where $x_1^* = t_1 a u_0 a t_0 a y a t_0$.

As in (*Biib2)(1a), we have that $M \vDash t_2 = t^{++}$.

By (5.58), $M \vDash Env(t_0, x_1^*) \ \& \ \forall w(w \ \varepsilon \ x_1^* \leftrightarrow w=u_0 \ v \ w=y)$.

Now, we have $M \vDash \neg \exists w(w \ \varepsilon \ x^{**} \ \& \ w \ \varepsilon \ x^-) \ \& \ \forall w(w \ \varepsilon \ x^{-*} \leftrightarrow w \ \varepsilon \ x^-)$.

Assume $M \vDash w \ \varepsilon \ x_1^*$.

$\implies\ M \vDash w=u_0 \lor w=y$,

$\implies\ M \vDash w\ \varepsilon\ x^{++} \lor w=y$,

$\implies\ M \vDash w\ \varepsilon\ x^{**}$.

Therefore, $M \vDash \forall w(w\ \varepsilon\ x_1^* \to w\ \varepsilon\ x^{**})$.

It follows that $M \vDash \neg \exists w(w\ \varepsilon\ x_1^* \,\&\, w\ \varepsilon\ x^{-*})$.

So we have

$M \vDash \mathrm{Env}(t_0,x_1^*) \,\&\, x_1^*=t_1au_0at_0ayat_0 \,\&\, z=t_1au_0at_0ayax^{-*} \,\&\, \mathrm{Env}(tb,x^{-*}) \,\&$

$\quad \&\ t_0 < t_2 b \,\&\, \mathrm{Firstf}(x^{-*},t_2b,av_0a,t_5) \,\&\, \neg\exists w(w\ \varepsilon\ x_1^* \,\&\, w\ \varepsilon\ x^{-*}) \,\&\, \mathrm{MinSet}(x_1^*) \,\&$

$\&\ \mathrm{MinSet}(x^{-*})$.

$\implies$ from $M \vDash I^{**}(z)$, by (10.6), $M \vDash \mathrm{MinSet}(z) \,\&\, \mathrm{Env}(t',z)$.

Also, from $M \vDash I^{**}(x_1^*) \,\&\, I^{**}(x^{-*})$, by (5.46),

$\qquad M \vDash \forall w(w\ \varepsilon\ z \leftrightarrow w\ \varepsilon\ x_1^* \lor w\ \varepsilon\ x^{-*})$.

Now, we have

$\quad M \vDash \mathrm{Max}^+T_b(t_1,u_0) \,\&\, \mathrm{Max}^+T_b(t_0,y) \,\&\, x_1^*= t_1au_0at_0ayat_0 \,\&\, t_1<t_0 \,\&\, u_0 \prec y$,

$\implies$ from $M \vDash I^{**}(u_0) \,\&\, I^{**}(y) \,\&\, I^{**}(t_1) \,\&\, I^{**}(t_0)$, by (9.30), $M \vDash \mathrm{Lex}^+(x_1^*)$.

$\implies$ from $M \vDash \forall u,v\ (u\ \varepsilon\ x^{**} \,\&\, v\ \varepsilon\ x^- \to u \prec v)$,

$\qquad M \vDash \forall u,v\ (u\ \varepsilon\ x_1^* \,\&\, v\ \varepsilon\ x^{-*} \to u \prec v)$.

We also have $M \vDash \mathrm{Lex}^+(x^{-*})$ as in (*Biib2)(2aia).

$\implies$ from $M \vDash I^{**}(z)$, by (9.29), $M \vDash \mathrm{Lex}^+(z)$,

$\implies$ from $M \vDash I^{**}(y)$, by (10.15), $M \vDash \mathrm{Special}(x_0^*)$.

So we have

$M \vDash \mathrm{Env}(t_0,x_1^*) \,\&\, x_1^*=t_1au_0ax_0^* \,\&\, \mathrm{MinMax}^+T_b(t_1,u_0) \,\&\, \neg(u_0\ \varepsilon\ x_0^*) \,\&$

$$\text{\& Firstf}(x_0^*,t_0,aya,t_0) \text{ \& } \forall u(u \; \varepsilon \; x_0^* \rightarrow u_0 \prec v) \text{ \& Special}(x_0^*).$$

$\Rightarrow$ from $M \vDash I^{**}(x_1^*)$, by (10.20), $M \vDash \text{Special}(x_1^*)$.

Hence we have, for $w'=au_0at_0aya$,

$$M \vDash \text{Env}(t',z) \text{ \& } z=t_1w'x^{-*} \text{ \& } aBw' \text{ \& } aEw' \text{ \& Set}(x^{-*}) \text{ \& } x_1^*=t_1w't_0 \text{ \& } t_0=t_2 \text{ \&}$$
$$\text{\& Env}(t_0,x_1^*) \text{ \& Firstf}(x^{-*},t_2b,av_0a,t_5) \text{ \& } \forall u,v \; (u \; \varepsilon \; x_1^* \text{ \& } v \; \varepsilon \; x^{-*} \rightarrow u \prec v) \text{ \&}$$
$$\text{\& } \neg \exists u(u \; \varepsilon \; x_1^* \text{ \& } u \; \varepsilon \; x^{-*}) \text{ \& Set}(x) \text{ \& } x^- \sim x^{-*} \text{ \& } \forall u,v(u \prec_{x^-} v \leftrightarrow u \prec_{x^{-*}} v) \text{ \&}$$
$$\text{\& } \forall v,t_6t_7 \; (\text{Free}^+(x^-,t_6,ava,t_7) \rightarrow \exists t_8 \text{Fr}(x^{-*},t_6b,ava,t_8)) \text{ \&}$$
$$\text{\& } \forall v,t_6,t_7 \; (\text{Bound}(x^-,t_6,ava,t_7) \vee \text{Free}^-(x^-,t_6,ava,t_7) \rightarrow \exists t_8 \text{Fr}(x^{-*},t_6,ava,t_8)) \text{ \&}$$
$$\text{\& Special}(x_1^*) \text{ \& Special}(x^-).$$

$\Rightarrow$ from $M \vDash I^{**}(x^-) \text{ \& } I^{**}(z)$, by (10.28), $M \vDash \text{Special}(z)$.

It follows that $M \vDash V^{**}(z)$.

Finally, from

$$M \vDash \forall w(w \; \varepsilon \; z \leftrightarrow w \; \varepsilon \; x_1^* \vee w \; \varepsilon \; x^{-*}) \text{ \& } \forall w(w \; \varepsilon \; x_1^* \leftrightarrow w=u_0 \vee w=y) \text{ \& } x^- \sim x^{-*} \text{ \&}$$
$$\text{\& } \forall w(w \; \varepsilon \; x^- \leftrightarrow w \; \varepsilon \; x \text{ \& } w \neq u_0) \text{ \& } u_0 \; \varepsilon \; x,$$

we have that $M \vDash \forall w(w \; \varepsilon \; z \leftrightarrow w \; \varepsilon \; x \vee w=y)$.

(2aib) $M \vDash t_0 < t_2$.

(2aic) $M \vDash t_2 < t_0$.

These are both ruled out as in (*Biib2).

(2aii) $M \vDash t_1b < t_0$.

$\Rightarrow M \vDash t_1 < t_1b < t_0 = t^{++}$,

$\Rightarrow$ as in (1), $M \vDash t^{++}=t_2$.

We then proceed exactly as in (2aia).

(2aiii)  $M \vDash t_0 < t_1 b$  is ruled out as in (*Biib2).

(2b)  $M \vDash t_0 \leq t_1$.

$\Rightarrow$ from (2), $M \vDash t_0 \leq t_1 = t^+ < t_0$, contradicting $M \vDash t_0 \in I \subseteq I_0$.

So (2b) is ruled out.

(3)  $M \vDash (u_0 \approx_{Tb} y \ \& \ u_0 \ll y) \ \& \ y \triangleleft_{Tb} v_0$.

$\Rightarrow M \vDash t_1 = t^+ < t_0 = t^{++}$.

As in (1), $M \vDash t_2 = t^{++}$.

(3a)  $M \vDash t_0 b = t^{++}$.

We proceed as in (1bi) of (*Biib2) and obtain $x^{+*}$ from $x^+ = t_0 a y a x^-$, where $M \vDash I^{**}(x^{+*}) \ \& \ \text{MinSet}(x^{+*})$  and  $M \vDash \text{Env}(t', x^{+*})$ with $M \vDash J^{**}(t')$.

Let $z = t_1 a u_0 a x^{+*}$.

$\Rightarrow$ from $M \vDash I^{**}(t_1) \ \& \ I^{**}(u_0) \ \& \ I^{**}(x^{+*})$, since $I^{**}$ may be assumed to be closed w.r. to *,       $M \vDash I^{**}(z)$.

We then have:

 $M \vDash \text{Env}(t_1, x_0) \ \& \ x_0 = t_1 a u_0 a t_1 \ \& \ z = t_1 a u_0 a x^{+*} \ \& \ \text{Env}(t', x^{+*}) \ \& \ t_1 < t_0 b \ \&$

   $\& \ \text{Firstf}(x^{+*}, t_0 b, a y a, t_6) \ \& \ \neg \exists w (w \ \varepsilon \ x_0 \ \& \ w \ \varepsilon \ x^{+*}) \ \& \ \text{MinSet}(x_0) \ \&$

$\& \ \text{MinSet}(x^{+*})$.

$\Rightarrow$ from  $M \vDash I^{**}(z)$, by (10.6), $M \vDash \text{MinSet}(z) \ \& \ \text{Env}(t', z)$.

Also, from  $M \vDash I^{**}(x_0) \ \& \ I^{**}(x^{+*})$, by (5.46),

$$M \vDash \forall w(w \ \varepsilon \ z \leftrightarrow w \ \varepsilon \ x_0 \lor w \ \varepsilon \ x^{+*}),$$

whence  $M \vDash \forall w(w \ \varepsilon \ z \leftrightarrow w = u_0 \lor w \ \varepsilon \ x^{+*})$.

We have  $M \vDash Lex^+(x_0)$,  and, as in (*Biib2)(1bi),  $M \vDash Lex^+(x^{+*})$.

It was proved there that  $M \vDash \forall u,v\ (u\ \varepsilon\ x^{++}\ \&\ v\ \varepsilon\ x^{+*} \to u\prec v)$.

$\Rightarrow$  since  $M \vDash \forall w(w\ \varepsilon\ x_0 \to w\ \varepsilon\ x^{++})$,  $M \vDash \forall u,v\ (u\ \varepsilon\ x_0\ \&\ v\ \varepsilon\ x^{+*} \to u\prec v)$,

$\Rightarrow$  from  $M \vDash I^{**}(z)\ \&\ Lex^+(x_0)\ \&\ Lex^+(x^{+*})$, by (9.29),  $M \vDash Lex^+(z)$.

Finally, as in (1bi) of (*Biib2), we have  $M \vDash Special(x^+)$.

We then have that

$M \vDash Env(t',z)\ \&\ z=t_1au_0ax^{+*}\ \&\ Set(x^{+*})\ \&\ x_0=t_1au_0at_1\ \&\ t_1=t_0\ \&\ Env(t_1,x_0)\ \&$

$\&\ Firstf(x^{+*},t_0b,aya,t_6)\ \&\ \forall u,v\ (u\ \varepsilon\ x_0\ \&\ v\ \varepsilon\ x^{+*} \to u\prec v)\ \&$

$\&\ \neg\exists u(u\ \varepsilon\ x_0\ \&\ u\ \varepsilon\ x^{+*})\ \&\ Set(x^+)\ \&\ x^+ \sim x^{+*}\ \&\ \forall u,v(u<_{x^+}v \leftrightarrow u<_{x^{+*}}v)\ \&$

$\&\ \forall v,t_7t_8\ (Free^+(x^+,t_7,ava,t_8) \to \exists t_9 Fr(x^{+*},t_7b,ava,t_9))\ \&$

$\&\ \forall v,t_7,t_8\ (Bound(x^+,t_7,ava,t_8)\ \vee\ Free^-(x^+,t_7,ava,t_8) \to \exists t_9 Fr(x^{+*},t_7,ava,t_9))\ \&$

$\&\ Special(x_0)\ \&\ Special(x^+)$.

$\Rightarrow$ from  $M \vDash I^{**}(z)\ \&\ I^{**}(x^+)$,  by (10.28),  $M \vDash Special(z)$.

Thus  $M \vDash V^{**}(z)$.

Finally, to see that  $M \vDash \forall w(w\ \varepsilon\ z \leftrightarrow w\ \varepsilon\ x\ \vee\ w=y)$,  note that we have that

$$M \vDash \forall w(w\ \varepsilon\ z \leftrightarrow w=u_0\ \vee\ w\ \varepsilon\ x^{+*}),$$

along with  $M \vDash x^+ \sim x^{+*}\ \&\ \forall w(w\ \varepsilon\ x^+ \leftrightarrow ((w\ \varepsilon\ x\ \&\ w \neq u_0)\ \vee\ w=y))$

and  $M \vDash u_0\ \varepsilon\ x$.

Then the desired claim easily follows.

  (3b)  $M \vDash t_0b<t^{++}$.

Let  $z=t_1au_0at_0bayat_2aw_2$.

Then  $M \vDash I^{**}(z)$.

For $x_0^{**}=t_0bayat_2aw_2$, we have, as in (*Biib2)(1bii), that

$$M \vDash I^{**}(x_0^{**}) \,\&\, MinSet(x_0^{**}) \,\&\, Env(t,x_0^{**}).$$

$\Rightarrow$ by (*Biib2)(1bii), $M \vDash \forall w(w \,\varepsilon\, x_0^{**} \leftrightarrow w \,\varepsilon\, x^+ \leftrightarrow ((w \,\varepsilon\, x \,\&\, w \neq u_0) \,v\, w=y))$,

$\Rightarrow$ since $M \vDash u_0 \neq y$, $M \vDash \neg(u_0 \,\varepsilon\, x^+)$,

$\Rightarrow M \vDash \neg \exists w(w \,\varepsilon\, x_0 \,\&\, w \,\varepsilon\, x_0^{**})$.

So we have

$M \vDash Env(t_1,x_0) \,\&\, x_0=t_1au_0at_1 \,\&\, z=t_1au_0ax_0^{**} \,\&\, Env(t,x_0^{**}) \,\&\, t_1<t_0b \,\&$

$\&\, Firstf(x_0^{**},t_0b,aya,t_2) \,\&\, \neg\exists w(w \,\varepsilon\, x_0 \,\&\, w \,\varepsilon\, x_0^{**}) \,\&\, MinSet(x_0) \,\&$

$\&\, MinSet(x_0^{**})$.

$\Rightarrow$ from $M \vDash I^{**}(z)$, by (10.6), $M \vDash MinSet(z) \,\&\, Env(t,z)$.

From $M \vDash I^{**}(x_0) \,\&\, I^{**}(x_0^{**})$, by (5.46), $M \vDash \forall w(w \,\varepsilon\, z \leftrightarrow w \,\varepsilon\, x_0 \,v\, w \,\varepsilon\, x_0^{**})$,

whence $M \vDash \forall w(w \,\varepsilon\, z \leftrightarrow w=u_0 \,v\, w \,\varepsilon\, x_0^{**})$.

Now, as in (*Biib2)(1bi), we have

$M \vDash Lex^+(x_0^{**})$ and $M \vDash \forall u,v \,(u \,\varepsilon\, x^{++} \,\&\, v \,\varepsilon\, x_0^{**} \to u \prec v)$.

$\Rightarrow$ since $M \vDash \forall w(w \,\varepsilon\, x_0 \to w \,\varepsilon\, x^{++})$, $M \vDash \forall u,v(u \,\varepsilon\, x_0 \,\&\, v \,\varepsilon\, x_0^{**} \to u \prec v)$.

So we have

$M \vDash Env(t,z) \,\&\, z=t_1au_0ax_0^{**} \,\&\, Env(t_1,x_0) \,\&\, x_0=t_1au_0at_1 \,\&\, Env(t,x_0^{**}) \,\&$

$\&\, Firstf(x_0^{**},t_0b,aya,t_2) \,\&\, t_1<t_0b \,\&\, \neg\exists w(w \,\varepsilon\, x_0 \,\&\, w \,\varepsilon\, x_0^{**}) \,\&$

$\&\, \forall u,v(u \,\varepsilon\, x_0 \,\&\, v \,\varepsilon\, x_0^{**} \to u \prec v) \,\&\, Lex^+(x_0) \,\&\, Lex^+(x_0^{**})$,

$\Rightarrow$ from $M \vDash I^{**}(z)$, by (9.29), $M \vDash Lex^+(z)$.

Let $x_2^* = t_1au_0at_0bayat_0b$.

$\Rightarrow$ from $M \vDash I^{**}(t_1) \,\&\, I^{**}(u_0) \,\&\, I^{**}(t_0) \,\&\, I^{**}(y)$, since $I^{**}$ may be assumed to

be closed w.r. to *,  $M \vDash I^{**}(x_2^*)$.

We have

  $M \vDash \text{Pref}(au_0a,t_1)$ & $\text{Pref}(aya,t_0b)$ & $t_1 < t_0b$ & $u_0 \neq y$.

$\Rightarrow$ by (5.58),  $M \vDash \text{Env}(t_0b, x_2^*)$ & $\forall w(w \varepsilon x_2^* \leftrightarrow w=u_0 \vee w=y)$.

So we have

  $M \vDash \text{Env}(t_0b, x_2^*)$ & $x_2^* = t_1au_0at_0bayat_0b$ & $\text{Env}(t_1, x_0)$ & $x_0 = t_1au_0at_1$ &

    & $\text{Max}^+T_b(t_0, y)$ & $\neg(y \varepsilon x_0)$ & $\forall u(u \varepsilon x_0 \to u \prec y)$ & $\text{Special}(x_0)$.

$\Rightarrow$ from $M \vDash I^{**}(x_2^*)$, by (10.18),  $M \vDash \text{Special}(x_2^*)$.

By earlier proofs,

  $M \vDash \neg \exists w(w \varepsilon x_2^* \& w \varepsilon x^-)$ & $\forall u,v(u \varepsilon x_2^* \& v \varepsilon x^- \to u \prec v)$.

So we have

  $M \vDash \text{Env}(t, z)$ & $z = t_1au_0at_0bayax^-$ & $\text{Env}(t_0b, x_2^*)$ & $x_2^* = t_1au_0at_0bayat_0b$ &

    & $\text{Env}(t, x^-)$ & $\text{Firstf}(x^-, t_2, av_0a, t_3)$ & $t_0b < t_2$ & $\neg \exists w(w \varepsilon x_2^* \& w \varepsilon x^-)$ &

      & $\forall u,v(u \varepsilon x_2^* \& v \varepsilon x^- \to u \prec v)$ & $\text{Special}(x_2^*)$ & $\text{Special}(x^-)$.

$\Rightarrow$ from $M \vDash I^{**}(z)$, by (10.17),  $M \vDash \text{Special}(z)$,  as required.

Hence  $M \vDash V^{**}(z)$.

Now, we have that

 $M \vDash \forall w(w \varepsilon z \leftrightarrow w=u_0 \vee w \varepsilon x_0^{**} \leftrightarrow w=u_0 \vee w \varepsilon x^+ \leftrightarrow$

                $\leftrightarrow w=u_0 \vee ((w \varepsilon x \& w \neq u_0) \vee w=y))$.

$\Rightarrow$ since  $M \vDash u_0 \varepsilon x$,  $M \vDash \forall w(w \varepsilon z \leftrightarrow w \varepsilon x \vee w=y)$,  as required.

(4) $M \vDash (u_0 \approx_{Tb} y \,\&\, u_0 \ll y) \,\&\, (y \approx_{Tb} v_0 \,\&\, y \ll v_0)$.

$\Longrightarrow M \vDash t^+ = t_0 = t^{++}$,

$\Longrightarrow M \vDash t_1 = t^+ = t^{++} \,\&\, t_1 = t_0$,

$\Longrightarrow$ since $M \vDash t_1 < t_2$, $M \vDash t^{++} < t_2$,

$\Longrightarrow M \vDash t_0 < t_2$.

(4a) $M \vDash t_1 b = t_2$.

Let $z = t_1 a u_0 a t_1 b a y a x^{-*}$ where $x^{-*}$ is as in (2).

Then $M \vDash \text{Env}(t', x^{-*})$, where $M \vDash J^{**}(t')$.

$\Longrightarrow$ from $M \vDash I^{**}(t_1) \,\&\, I^{**}(u_0) \,\&\, I^{**}(y) \,\&\, I^{**}(x^{-*})$, since $I^{**}$ may be assumed to be closed w.r. to *, $M \vDash I^{**}(z)$.

Let $x_3^* = t_1 a u_0 a t_1 b a y a t_1 b$.

As in (2aia), we have that

$M \vDash I^{**}(x_3^*) \,\&\, \text{Env}(t_1 b, x_3^*) \,\&\, \forall w(w \,\varepsilon\, x_3^* \leftrightarrow w = u_0 \lor w = y)$,

and, further, that $M \vDash \text{MinSet}(x_3^*)$.

Since $M \vDash \text{MinSet}(x^{-*})$, we may proceed exactly analogously to (2aia) to show that

$M \vDash \text{MinSet}(z) \,\&\, \text{Env}(t', z)$.

From $M \vDash I^{**}(x_3^*) \,\&\, I^{**}(x^{-*})$, by (5.46), we have that

$M \vDash \forall w(w \,\varepsilon\, z \leftrightarrow w \,\varepsilon\, x_3^* \lor w \,\varepsilon\, x^{-*})$.

Then $M \vDash \text{Lex}^+(z)$ is proved exactly as in (2aia).

With $x_0 = t_1 a u_0 a t_1$, we then have

$M \vDash \text{Env}(t_1 b, x_3^*) \,\&\, x_3^* = t_1 a u_0 a t_1 b a y a t_1 b \,\&\, \text{Env}(t_1, x_0) \,\&\, x_0 = t_1 a u_0 a t_1 \,\&\,$

    & $Max^+T_b(t_1b,y)$ & $\neg(y \; \varepsilon \; x_0)$ & $\forall u(u \; \varepsilon \; x_0 \rightarrow u_0 \prec y)$ & $Special(x_0)$.

$\Rightarrow$ from $M \vDash I^{**}(x_3^*)$, by (10.18), $M \vDash Special(x_3^*)$.

We also have that $M \vDash Special(x^{-*})$.

Then, analogously to (2aia), again with $x_3^*$ replacing $x_1^*$ throughout the argument, we show, using (10.28), that $M \vDash Special(z)$.

Thus $M \vDash V^{**}(z)$.

Finally, that    $M \vDash \forall w(w \; \varepsilon \; z \leftrightarrow w \; \varepsilon \; x \lor w=y)$

follows analogously to (2aia).

  (4b) $M \vDash t_1b < t_2$.

We claim that $M \vDash Max^+T_b(t_1b,v_0,x)$ where $M \vDash Firstf(x^-,t_2,av_0a,t_3)$.

Note that, as in (*Biib2), $M \vDash Fr(x,t_2,av_0a,t_3)$. So $M \vDash v_0 \; \varepsilon \; x$.

Since $M \vDash Tally_b(t_1)$ & $Tally_b(t_2)$ & $t_1<t_1b<t_2$, we have $M \vDash Tally_b(t_1b)$.

Assume that $M \vDash Fr(x,t_4,aua,t_5)$ & $u<_x v_0$.

$\Rightarrow$ from $M \vDash I^{**}(x)$, by (9.19), $M \vDash u=u_0$,

$\Rightarrow$ from $M \vDash Fr(x,t_1,au_0a,t_2)$ & $Env(t,x)$, $M \vDash t_4=t_1$,

$\Rightarrow$ $M \vDash t_4=t_1<t_1b$.

So $M \vDash \forall u,t_4,t_5(Fr(x,t_4,aua,t_5)$ & $u<_x v_0 \rightarrow t_4<t_1b)$.

Along with $M \vDash Set(x)$ & $v_0 \; \varepsilon \; x$ & $Tally_b(t_1b)$, we then have

$$M \vDash Max^+(t_1b,v_0,x)$$

as claimed.

$\Rightarrow$ from hypothesis $M \vDash t_1=t^{++}$, $M \vDash Max^+T_b(t_1b,av_0a)$,

$\Rightarrow$ from hypothesis $M \vDash Special(x)$ and $M \vDash Fr(x,t_2,av_0a,t_3)$,

$$M \vDash \text{MMax}^+(t_2, v_0, x),$$

$\implies M \vDash t_2 \leq t_1 b,$

$\implies M \vDash t_2 \leq t_1 b < t_2,$ contradicting $M \vDash t_2 \in I \subseteq I_0.$

(4c) $M \vDash t_2 < t_1 b.$

$\implies$ by (1.13), $M \vDash t_2 \leq t_1,$

$\implies M \vDash t_1 \leq t_2 \leq t_1,$ contradicting $M \vDash t_1 \in I \subseteq I_0.$

(*Biib3): $M \vDash \text{Lastf}(x, t_1, au_0 a, t_2).$

$\implies M \vDash \text{Env}(t_1, x),$

$\implies M \vDash t_1 = t_2 \ \& \ (x = t_1 au_0 at_2 \ \lor \ \exists w_1(x = w_1 at_1 au_0 at_2 \ \& \ \text{Max}^+ T_b(t_1, w_1))).$

(*Biib3i): $M \vDash t_1 = t_2 \ \& \ x = t_1 au_0 at_2.$

This is ruled out as in (*Biib1i).

(*Biib3ii): $M \vDash \exists w_1(x = w_1 at_1 au_0 at_2 \ \& \ \text{Max}^+ T_b(t_1, w_1)).$

$\implies M \vDash t_1 = t_2.$

We have $M \vDash u_0 \prec y.$

$\implies M \vDash (u_0 \triangleleft_{T_b} y \ \lor \ (u_0 \approx_{T_b} y \ \& \ u_0 \ll y)).$

(1) $M \vDash u_0 \triangleleft_{T_b} y.$

$\implies M \vDash t^+ < t_0.$

(1a) $M \vDash t_1 < t_0.$

Let $z = w_1 at_1 au_0 at_0 ayat_0.$

Then $z = x^{**}.$

Then, as in (*Biib2)(2aia), we have that

$$M \vDash I^{**}(x^{**}) \ \& \ \text{Env}(t_0, x^{**}) \ \& \ \text{MinSet}(x^{**}) \ \& \ \text{Lex}^+(x^{**}).$$

As in (*Biib2)(1a), $M \vDash \text{Special}(x^{**})$.

On the other hand, we have that $M \vDash \text{MinMax}^+T_b(t_0,y)$, and, as observed earlier, from $M \vDash V^{**}(y)$ it follows that $M \vDash J^{**}(t_0)$.

Therefore $M \vDash V^{**}(z)$.

Finally, from $M \vDash I^{**}(x)\ \&\ I^{**}(y)$, the proof of (7.1)(2),

$$M \vDash \forall w(w\ \varepsilon\ z \leftrightarrow w\ \varepsilon\ x \vee w=y).$$

(1b) $M \vDash t_0 \leq t_1$.

Let $z = w_1 a t_1 a u_0 a t_1 b a y a t_1 b$.

$\Rightarrow$ from $M \vDash I^{**}(w_1)\ \&\ I^{**}(t_1)\ \&\ I^{**}(u_0)\ \&\ I^{**}(y)$, since $I^{**}$ may be assumed to be closed w.r. to *, $M \vDash I^{**}(z)$.

Recalling (*Biib2)(2bi), we have

$M \vDash x_3^{++} = t^*w_0^*x_3 = t^*w_0^*t_1 b a y a t_1 b = x^{++} b a y a t_1 b = w_1 a t_1 a u_0 a t_1 b a y a t_1 b = z$.

Then we have, as in (*Biib2)(2bi), that

$$M \vDash \text{Env}(t_1 b, x_3^{++})\ \&\ \text{MinSet}(x_3^{++})\ \&\ \text{Lex}^+(x_3^{++})\ \&\ \text{Special}(x_3^{++}).$$

$\Rightarrow$ from $M \vDash \text{Env}(t_1,x)\ \&\ V^{**}(x)$, $M \vDash J^{**}(t_1)$,

$\Rightarrow$ since $J^{**}$ is a string concept, $M \vDash J^{**}(t_1 b)$.

Hence $M \vDash V^{**}(z)$.

Again, by the proof of (7.1)(2), $M \vDash \forall w(w\ \varepsilon\ z \leftrightarrow w\ \varepsilon\ x \vee w=y)$.

(2) $M \vDash u_0 \approx_{T_b} y\ \&\ u_0 \ll y$.

$\Rightarrow M \vDash t^+ = t_0$,

$\Rightarrow M \vDash t_0 \leq t_1$.

Exactly as (1b).

This completes the proof of the STRONG SET ADJUNCTION LEMMA.

Let $\sigma^*(x,y,z)$ abbreviate the formula

$\exists! t_0(\text{MinMax}^+T_b(t_0,y)\ \&\ ((*A)[x,y,z,t_0]\ \vee\ (*B)[x,y,z,t_0])$

where

$(*A) \equiv x=aa\ \vee\ z=t_0ayat_0$,

$(*B) \equiv \exists t(\text{Env}(t,x)\ \&\ ((y\ \varepsilon\ x\ \&\ z=x)\ \vee\ (\neg(y\ \varepsilon\ x)\ \&$

$\&\ ((*Bi)[x,y,z,t_0])\ \vee\ (*Bii)[x,y,z,t_0])))).$

Here

$(*Bi) \equiv \forall u(u\ \varepsilon\ x \to u\prec y)\ \&\ ((*Bia)\ \vee\ (*Bib))$,

$(*Bii) \equiv \exists u(u\ \varepsilon\ x\ \&\ \neg(u\prec y))\ \&\ ((*Biia)\ \vee\ (*Biib)).$

Now,

$(*Bia) \equiv t<t_0\ \&\ \exists!t^*(tt^*=t_0\ \&\ z=xt^*ayat_0)$,

$(*Bib) \equiv t_0\leq t\ \&\ z=xbayatb$,

$(*Biia) \equiv \forall u(u\ \varepsilon\ x \to u\prec y)\ \&\ \exists! v_0,t_1,t_2(\text{Firstf}(x,t_1,av_0a,t_2)\ \&$

$\&\ (y\triangleleft_{Tb}v_0\ \vee\ (y\approx_{Tb}v_0\ \&\ y\ll v_0)))$,

$(*Biib) \equiv \exists u(u\ \varepsilon\ x\ \&\ u\lessapprox y)\ \&\ \exists! u_0(\text{Max}_\lessapprox(u_0,x,y)\ \&\ \exists t_1,t_2(\text{Fr}(x,t_1,au_0a,t_2)\ \&$

$\&\ ((*Biib1)\ \vee\ (*Biib2)\ \vee\ (*Biib3)))).$

Here we have

$(*Biib1) \equiv \text{Firstf}(x,t_1,au_0a,t_2)\ \&\ \exists w_2(x=t_1au_0at_2aw_2\ \&$

$\&\ ((*Biib1ii1)\ \vee\ (*Biib1ii2)\ \vee\ (*Biib1ii3)\ \vee\ (*Biib1ii4)))$,

$(*Biib2) \equiv \exists w_1,w_2(Intf(x,w_1,t_1,au_0a,t_2)\ \&\ x=w_1at_1au_0at_2aw_2\ \&$

$\&\ \exists x^-(Env(t,x^-)\ \&\ x^-=t_2aw_2\ \&\ \exists v_0,t_3((Firstf(x^-,t_2,av_0a,t_3)\ \&$

$\&\ \exists t^+,t^{++}((MinMax^+T_b(t^+,u_0)\ \&\ MinMax^+T_b(t^{++},v_0)\ \&$

$\&\ ((*Biib2.1)\ v\ (*Biib2.3)\ v\ (*Biib2.2)\ v\ (*Biib2.4)))))),$

$(*Biib3) \equiv Lastf(x,t_1,au_0a,t_2)\ \&\ \exists w_1(x=w_1at_1au_0at_2\ \&$

$\&\ ((*Biib3ii1a)\ v\ (*Biib3ii1b)\ v\ (*Biib3ii2))).$

Furthermore,

$(*Biib1ii1) \equiv (u_0 \triangleleft_{T_b} y\ \&\ y \triangleleft_{T_b} v_0)\ \&\ z=t_1au_0at_0ayat_2aw_2,$

$(*Biib1ii2) \equiv (u_0 \triangleleft_{T_b} y\ \&\ (y \approx_{T_b} v_0\ \&\ y \ll v_0))\ \&\ \exists x^-(Env(t,x^-)\ \&\ x^-=t_2aw_2\ \&$

$\&\ \exists!x^{-*}(R_{10.29}(x^-,x^{-*})\ \&\ z=t_1au_0at_0ayax^{-*})),$

$(*Biib1ii3) \equiv ((u_0 \approx_{T_b} y\ \&\ u_0 \ll y)\ \&\ y \triangleleft_{T_b} v_0)\ \&\ ((*Biib1ii3a)\ v\ (*Biib1ii3b)),$

$(*Biib1ii4) \equiv ((u_0 \approx_{T_b} y\ \&\ u_0 \ll y)\ \&\ (y \approx_{T_b} v_0\ \&\ y \ll v_0))\ \&$

$\&\ \exists!x^{-*}(R_{10.29}(x^-,x^{-*})\ \&\ z=t_1au_0at_1bayax^{-*}).$

On the other hand, we let

$(*Biib2.1)\ v\ (*Biib2.3) \equiv (u_0 \triangleleft_{T_b} y\ \&\ y \triangleleft_{T_b} v_0)\ v\ (u_0 \approx_{T_b} y\ \&\ u_0 \ll y\ \&\ y \triangleleft_{T_b} v_0))\ \&$

$\&\ ((*Biib2.1a)\ v\ (*Biib2.1b)\ v\ (*Biib2.1c)),$

$(*Biib2.2) \equiv (u_0 \triangleleft_{T_b} y\ \&\ (y \approx_{T_b} v_0\ \&\ y \ll v_0))\ \&$

$\&\ \exists!x^{-*}(R_{10.29}(x^-,x^{-*})\ \&\ z=w_1at_1au_0at_2ayax^{-*}),$

$(*Biib2.4) \equiv ((u_0 \approx_{T_b} y\ \&\ u_0 \ll y)\ \&\ (y \approx_{T_b} v_0\ \&\ y \ll v_0))\ \&$

$\&\ \exists!x^{-*}(R_{10.29}(x^-,x^{-*})\ \&\ z=w_1at_1au_0at_1bayax^{-*}).$

Now, we have

(*Biib3ii1a) ≡ $t_1 < t_0$ & $z = w_1 a t_1 a u_0 a t_0 a y a t_0$,

(*Biib3ii1b) v (*Biib3ii2) ≡ $t_0 \leq t_1$ & $z = w_1 a t_1 a u_0 a t_1 b a y a t_1 b$,

whereas

(*Biib1ii3a) ≡ $t_0 b = t^{++}$ & $\exists x^+ (Env(t, x^+)$ & $x^+ = t_0 a y a x^-$ &

& $\exists ! x^{+*}(R_{10.29}(x^+, x^{+*})$ & $z = t_1 a u_0 a x^{+*}))$.

(*Biib1ii3b) ≡ $t_0 b < t^{++}$ & $z = t_1 a u_0 a t_0 b a y a t_2 a w_2$.

Also,

(*Biib2.1a) ≡ $t_1 < t_0$ & $z = w_1 a t_1 a u_0 a t_0 a y a t_2 a w_2$,

(*Biib2.1b) ≡ $t_1 = t_0$ & ((*Biib2.1bi) v (*Biib2.1bii)),

(*Biib2.1c) ≡ $t_0 = t_1$ & $t_1 \leq t^{++}$ & ((*Biib2.1cia-ciia) v (*Biib2.1cib-ciib)),

where

(*Biib2.1bi) ≡ $t_0 b = t^{++}$ & $\exists x^+ (Env(t, x^+)$ & $x^+ = t_0 a y a t_2 a w_2$ &

& $\exists ! x^{+*}(R_{10.29}(x^+, x^{+*})$ & $z = w_1 a t_1 a u_0 a x^{+*}))$,

(*Biib2.1bii) ≡ $t_0 b < t^{++}$ & $z = w_1 a t_1 a u_0 a t_0 b a y a t_2 a w_2$.

Finally,

(*Biib2.1cia-ciia) ≡ $t_1 b < t_2$ & $z = w_1 a t_1 a u_0 a t_1 b a y a t_2 a w_2$,

(*Biib2.1cib-ciib) ≡ $t_1 b = t_2$ & $\exists x_1^+ (Env(t, x_1^+)$ & $x_1^+ = t_0 a y a x^-$ &

& $\exists ! x_1^{+*}(R_{10.29}(x^+, x_1^{+*})$ & $z = w_1 a\ t_1 a u_0 a x_1^{+*}))$.

In (*Biib1ii2), (*Biib1ii4), (*Biib2.2), (*Biib2.4), (*Biib1ii3a), (*Biib2.1bi), (*Biib2.1cib-ciib), the formula $R_{10.29}(x,y)$ abbreviates

$$\exists t' \text{Env}(t',y) \& x \sim y \& \forall v,w \ (v <_x w \leftrightarrow v <_y w) \&$$

$$\& \forall w,t_1,t_2 \ (\text{Free}^+(x,t_1,awa,t_2) \rightarrow \exists t_3 \text{Fr}(y,t_1 b,awa,t_3)) \&$$

$$\& \forall w,t_1,t_2 \ (\text{Bound}(x,t_1,awa,t_2) \lor \text{Free}^-(x,t_1,awa,t_2) \rightarrow \exists t_3 \text{Fr}(y,t_1,awa,t_3))$$

expressing the relationship between z and z' in (10.29).

Remark:  The proof of the STRONG SET ADJUNCTION LEMMA shows that

QT$^+$ ⊦ ∀x,y (V**(x) & V**(y) → ∃!z (V**(z) & σ*(x,y,z)) &

    & ∀z(σ*(x,y,z) → ∀w (w ε z ↔ (w ε x v w=y)) & (σ*(x,y,x) ↔ y ε x)).

# 13. Interpreting AST in QT⁺

The theory AST (Adjunctive Set Theory) is formulated in the first-order language with identity whose sole non-logical symbol is the binary relation symbol $\in$. Its non-logical axioms are

NULL.  $\exists x \forall y \neg (y \in x)$

ADJ.    $\forall x,y \exists z \forall w (w \in z \leftrightarrow w \in x \lor w=y)$.

We introduce the following abbreviations for certain formulae of QT⁺:

$\quad \text{Set}^+(x) \equiv I_{7.1}(x) \,\&\, \text{Set}(x)$

$\quad [x \in y]^+ \equiv \text{Set}^+(x) \,\&\, x \,\varepsilon\, y$

$\quad [x=y]^+ \equiv x=y$.

Thus each atomic formula A of AST is associated with a unique formula $[A]^+$ of QT⁺. We extend this map in a natural way to arbitrary formulae of AST by setting

$\quad [\neg A]^+ \equiv \neg [A]^+ \quad \text{and} \quad [A \,\&\, B]^+ \equiv [A]^+ \,\&\, [B]^+,$

and further

$\quad [\exists x\, A]^+ \equiv \exists x(\text{Set}^+(x) \,\&\, [A]^+),$

with the remaining logical operations defined in terms of $\neg$, & and $\exists$ in the usual way.

Claim 13.1.  $\text{QT}^+ \vdash [\exists x \forall y \neg (y \in x)]^+$.

We show that

$\quad \text{QT}^+ \vdash \exists x(\text{Set}^+(x) \,\&\, \forall y(\text{Set}^+(y) \to \neg(\text{Set}^+(y) \,\&\, y \,\varepsilon\, x)))$.

Consider $x = aa$.

We have that $QT^+ \vdash Set^+(aa)$.

By (5.18), $QT^+ \vdash \forall y \neg(y \in aa)$.

$\Rightarrow QT^+ \vdash \forall y(Set^+(y) \to \neg(Set^+(y) \& y \, \varepsilon \, aa))$,

$\Rightarrow QT^+ \vdash \exists x(Set^+(x) \& \forall y(Set^+(y) \to \neg(Set^+(y) \& y \, \varepsilon \, aa)))$.

That is, $QT^+ \vdash [\exists x \forall y \neg(y \in x)]^+$.

Claim 13.2. $QT^+ \vdash [\forall x,y \exists z \forall w(w \in z \leftrightarrow w \in x \lor w = y)]^+$.

We show that

$\quad QT^+ \vdash \forall x,y(Set^+(x) \& Set^+(y) \to \exists z(Set^+(z) \& \forall w(Set^+(w) \to$

$\qquad\qquad \to (Set^+(w) \& w \, \varepsilon \, z \leftrightarrow (Set^+(w) \& w \, \varepsilon \, x) \lor w = y)))))$.

Assume $M \vDash Set^+(x) \& Set^+(y)$.

By SET ADJUNCTION LEMMA (7.1),

(*) $\quad M \vDash \exists z(Set^+(z) \& \forall w(w \, \varepsilon \, z \leftrightarrow w \, \varepsilon \, x \lor w = y))$.

Assume $M \vDash Set^+(w)$.

Suppose that $M \vDash Set^+(w) \& w \, \varepsilon \, z$.

$\Rightarrow$ from (*), $M \vDash w \, \varepsilon \, x \lor w = y$,

$\Rightarrow M \vDash (Set^+(w) \& w \, \varepsilon \, x) \lor w = y$.

So $M \vDash Set^+(w) \& w \, \varepsilon \, z \to (Set^+(w) \& w \, \varepsilon \, x) \lor w = y$.

Conversely, suppose that $M \vDash (Set^+(w) \& w \, \varepsilon \, x) \lor w = y$.

If $M \vDash Set^+(w) \& w \, \varepsilon \, x$, then from (*), $M \vDash w \, \varepsilon \, z$.

But then $M \vDash Set^+(w) \& w \, \varepsilon \, z$.

If  M ⊨ w=y,  then from (*),  M ⊨ w ε z.   Hence  M ⊨ Set⁺(w) & w ε z.

Therefore,  M ⊨ (Set⁺(w) & w ε x) v w=y → Set⁺(w) & w ε z.

This completes the proof of

QT⁺ ⊢ ∀x,y(Set⁺(x) & Set⁺(y) → ∃z(Set⁺(z) & ∀w(Set⁺(w) →

$\quad\quad\quad\quad\quad$ → (Set⁺(w) & w ε z ↔ (Set⁺(w) & w ε x) v w=y))))).

But then  QT⁺ ⊢ [∀x,y∃z∀w(w ∈ z ↔ w ∈ x v w=y)]⁺.

## 14. Interpreting $PS_0$ in $QT^+$

The theory $PS_0$ is formulated in the first-order language $\{\,0,\,;\,\}$, where "0" is an individual constant and ";" a binary function symbol.

Its axioms are

(1) $0;x \neq 0$

(2) $x;y,y = x;y$

(3) $x;y,z = x;z,y$

(4) $x;y,z = x;y \leftrightarrow x;z = x \lor z = y$

where we write "$x;y,z$" for "$(x;y);z$".

Note that

   (2) & (3) $\vdash x;z = x \lor z = y \rightarrow x;y,z = x;y$.

For, assume $z = y$.

$\Rightarrow x;y,z = x;y,y = x;y$ by (2).

Hence (2) $\vdash z = y \rightarrow x;y,z = x;y$.

Assume $x;z = x$.

$\Rightarrow x;z,y = (x;z);y = x;y$.

But $x;z,y = x;y,z$ by (3).

$\Rightarrow x;y,z = x;y$.

Hence (3) $\vdash x;z = x \rightarrow x;y,z = x;y$.

Therefore

   (2) & (3) $\vdash x;z = x \lor z = y \rightarrow x;y,z = x;y$, as claimed.

The theory $PS_0$ was introduced by Kirby as the quantifier-free theory of finite sets. Letting

$$x \in y \equiv x;y = x,$$

we may informally express the meaning of (1)-(4) in more familiar notation as:

$$\neg(x \in \emptyset)$$

$$y \in x \cup \{y\}$$

$$(x \cup \{y\}) \cup \{z\} = (x \cup \{z\}) \cup \{y\}$$

$$z \in x \cup \{y\} \leftrightarrow z \in x \vee z = y.$$

("PS" is meant to stand for Peano Set Theory.)

$PS_0$ may be thought of as a minimal theory of the <u>adjunction operation</u> $x \cup \{y\}$. We shall reformulate $PS_0$ as a first-order theory $PS_0'$ in the language $\{0, S\}$ where "0" is an individual constant and "S" a ternary relation symbol:

(1') $\neg S(0,x,0)$

(2') $S(x,y,z_1) \,\&\, S(z_1,y,z_2) \rightarrow z_1 = z_2$

(3') $S(x,y,z_1) \,\&\, S(z_1,z,z_2) \,\&\, S(x,z,z_3) \,\&\, S(z_3,y,z_4) \rightarrow z_2 = z_4$

(4') $S(x,y,z) \,\&\, S(z,w,z) \rightarrow S(x,w,x) \vee w = y$

(5') $\exists z\, S(x,y,z)$

(6') $S(x,y,z_1) \,\&\, S(x,y,z_2) \rightarrow z_1 = z_2.$

Let $[0]^* \equiv aa$, and let $[S(x,y,z)]^* \equiv \sigma^*(x,y,z)$,

where $\sigma^*(x,y,z)$ is as in the Remark following the STRONG SET ADJUNCTION LEMMA (12.2).

Let $[x = y]^* \equiv x=y$.

(14.1) <u>Claim</u>: $QT^+ \vdash V^{**}(x) \rightarrow \neg[S(0,x,0)]^*$.

Assume $M \vDash V^{**}(x)$.

Suppose, for a reductio, that $M \vDash \sigma^*(aa,x,aa)$.

$\Rightarrow$ by (12.2), $M \vDash \forall w(w \, \varepsilon \, aa \leftrightarrow w \, \varepsilon \, aa \vee w=x)$,

$\Rightarrow M \vDash x \, \varepsilon \, aa \leftrightarrow x \, \varepsilon \, aa \vee x=x$.

But, by (5.18), $M \vDash \neg(x \, \varepsilon \, aa)$, whereas $M \vDash x=x \vee x \, \varepsilon \, aa$, a contradiction.

Therefore, $M \vDash \neg\sigma^*(aa,x,aa)$, that is, $M \vDash \neg[S(0,x,0)]^*$.

Hence $M \vDash V^{**}(x) \rightarrow \neg\sigma^*(aa,x,aa)$.

(14.2) <u>Claim</u>: $QT^+ \vdash V^{**}(x) \, \& \, V^{**}(y) \, \& \, V^{**}(z_1) \, \& \, V^{**}(z_2) \rightarrow$

$$\rightarrow [S(x,y,z_1) \, \& \, S(z_1,y,z_2) \rightarrow z_1=z_2]^*.$$

Assume $M \vDash V^{**}(x) \, \& \, V^{**}(y) \, \& \, V^{**}(z_1) \, \& \, V^{**}(z_2)$ along with

$M \vDash \sigma^*(x,y,z_1) \, \& \, \sigma^*(z_1,y,z_2)$.

$\Rightarrow M \vDash y \, \varepsilon \, x \vee \neg(y \, \varepsilon \, x)$.

(a) $M \vDash y \, \varepsilon \, x$.

$\Rightarrow$ from $M \vDash \sigma^*(x,y,z_1)$, $M \vDash x=z_1$,

$\Rightarrow M \vDash y \, \varepsilon \, z_1$,

$\Rightarrow$ from $M \vDash \sigma^*(z_1,y,z_2)$, $M \vDash z_1=z_2$.

(b) $M \vDash \neg(y \ \varepsilon \ x)$.

$\Rightarrow$ from (12.2), $M \vDash y \ \varepsilon \ z_1$,

$\Rightarrow$ from $M \vDash \sigma^*(z_1,y,z_2)$, $M \vDash z_1=z_2$.

Hence we have

$M \vDash V^{**}(x) \ \& \ V^{**}(y) \ \& \ V^{**}(z_1) \ \& \ V^{**}(z_2) \ \to \ (\sigma^*(x,y,z_1) \ \& \ \sigma^*(z_1,y,z_2) \to z_1=z_2)$,

that is,

$QT^+ \vdash V^{**}(x) \ \& \ V^{**}(y) \ \& \ V^{**}(z_1) \ \& \ V^{**}(z_2) \ \to \ [S(x,y,z_1) \ \& \ S(z_1,y,z_2) \to z_1=z_2]^*$.

(14.3) <u>Claim</u>:

$QT^+ \vdash V^{**}(x) \ \& \ V^{**}(y) \ \& \ V^{**}(z) \ \& \ V^{**}(z_1) \ \& \ V^{**}(z_2) \ \& \ V^{**}(z_3) \ \& \ V^{**}(z_4) \to$

$\to [S(x,y,z_1) \ \& \ S(z_1,z,z_2) \ \& \ S(x,z,z_3) \ \& \ S(z_3,y,z_4) \to z_2=z_4]^*$.

Assume that $M \vDash V^{**}(x) \ \& \ V^{**}(y) \ \& \ V^{**}(z_1) \ \& \ V^{**}(z_2) \ \& \ V^{**}(z_3) \ \& \ V^{**}(z_4)$

and $\qquad M \vDash \sigma^*(x,y,z_1) \ \& \ \sigma^*(z_1,z,z_2) \ \& \ \sigma^*(x,z,z_3) \ \& \ \sigma^*(z_3,y,z_4)$.

We then have that

$M \vDash \forall w(w \ \varepsilon \ z_2 \leftrightarrow w \ \varepsilon \ z_1 \lor w=z \leftrightarrow (w \ \varepsilon \ x \lor w=y) \lor w=z \leftrightarrow$

$\leftrightarrow (w \ \varepsilon \ x \lor w=z) \lor w=y \leftrightarrow w \ \varepsilon \ z_3 \lor w=y \leftrightarrow w \ \varepsilon \ z_4)$,

$\Rightarrow M \vDash z_2 \sim z_4$.

But from $M \vDash V^{**}(z_2) \ \& \ V^{**}(z_4)$ we have $M \vDash Set^*(z_2) \ \& \ Set^*(z_4)$.

Hence, by THE UNIQUENESS LEMMA (11.4), $M \vDash z_2=z_4$,

which suffices to prove the claim.

(14.4) <u>Claim</u>: QT⁺ ⊢ V**(x) & V**(y) & V**(z) & V**(w) →

$$\to [S(x,y,z) \mathbin{\&} S(z,w,z) \to S(x,w,x) \lor w=y]^*.$$

Assume   M ⊨ V**(x) & V**(y) & V**(z) & V**(w)

along with        M ⊨ σ*(x,y,z) & σ*(z,w,z).

From STRONG SET ADJUNCTION LEMMA we have

   M ⊨ ∀v(v ε z ↔ w ε x ∨ v=y)  and  M ⊨ ∀v(v ε z ↔ v ε z ∨ v=w).

⟹ M ⊨ w ε z.

Assume that   M ⊨ w≠y.

⟹ M ⊨ w ε x,

⟹ from STRONG SET ADJUNCTION LEMMA,   M ⊨ σ*(x,w,x),  as required.

(14.5) <u>Claim</u>: QT⁺ ⊢ V**(x) & V**(y) → ∃z(V**(z) & [S(x,y,z)]*).

(14.6) <u>Claim</u>: QT⁺ ⊢ V**(x) & V**(y) & V**($z_1$) & V**($z_2$) →

$$\to [S(x,y,z_1) \mathbin{\&} S(x,y,z_2) \to z_1=z_2]^*.$$

This follows immediately from THE STRONG SET ADJUNCTION LEMMA.

## 15. Interpreting $PS_0 + EXT$ in $QT^+$

The axioms of $PS_0 + EXT$ are those of $PS_0$ together with

   EXT.    $\forall x,y(\forall z(S(x,z,x) \leftrightarrow S(y,z,y)) \rightarrow x=y)$.

We let $V^{**}(x)$ define the domain of the interpretation $**$.

Let $[0]^{**} \equiv aa$.

Let $[S(x,y,z)]^{**} \equiv \sigma^*(x,y,z)$.

Let $[x=y]^{**} \equiv x \stackrel{*}{=} y$

where       $x \stackrel{*}{=} y \equiv \forall z(V^{**}(z) \rightarrow (z \, \varepsilon \, x \leftrightarrow z \, \varepsilon \, y))$.

(15.1) <u>Claim</u>: $QT^+ \vdash V^{**}(x) \rightarrow \neg[S(0,x,0)]^{**}$.

Same as (14.1).

(15.2) <u>Claim</u>: $QT^+ \vdash V^{**}(x) \, \& \, V^{**}(y) \, \& \, V^{**}(z_1) \, \& \, V^{**}(z_2) \rightarrow$

$$\rightarrow [S(x,y,z_1) \, \& \, S(z_1,y,z_2) \rightarrow z_1=z_2]^{**}.$$

Assume $M \vDash V^{**}(x) \, \& \, V^{**}(y) \, \& \, V^{**}(z_1) \, \& \, V^{**}(z_2)$ along with

   $M \vDash \sigma^*(x,y,z_1) \, \& \, \sigma^*(z_1,y,z_2)$.

Exactly the same argument applies as in (14.2) to obtain $M \vDash z_1 = z_2$. The desired conclusion follows immediately.

(15.3) <u>Claim</u>:

$QT^+ \vdash V^{**}(x) \& V^{**}(y) \& V^{**}(z) \& V^{**}(z_1) \& V^{**}(z_2) \& V^{**}(z_3) \& V^{**}(z_4) \rightarrow$

$\rightarrow [S(x,y,z_1) \& S(z_1,z,z_2) \& S(x,z,z_3) \& S(z_3,y,z_4) \rightarrow z_2=z_4]^{**}$.

Assume $M \vDash V^{**}(x) \& V^{**}(y) \& V^{**}(z_1) \& V^{**}(z_2) \& V^{**}(z_3) \& V^{**}(z_4)$

and $M \vDash \sigma^*(x,y,z_1) \& \sigma^*(z_1,z,z_2) \& \sigma^*(x,z,z_3) \& \sigma^*(z_3,y,z_4)$.

We argue just as in (14.3) and obtain

$$M \vDash \forall w(w \, \varepsilon \, z_2 \leftrightarrow w \, \varepsilon \, z_4),$$

whence $M \vDash z_2 \stackrel{*}{=} z_4$.

(15.4) <u>Claim</u>: $QT^+ \vdash V^{**}(x) \& V^{**}(y) \& V^{**}(z) \& V^{**}(w) \rightarrow$

$\rightarrow [S(x,y,z) \& S(z,w,z) \rightarrow S(x,w,x) \lor w=y]^{**}$.

Assume $M \vDash V^{**}(x) \& V^{**}(y) \& V^{**}(z) \& V^{**}(w)$

along with $M \vDash \sigma^*(x,y,z) \& \sigma^*(z,w,z)$.

From STRONG SET ADJUNCTION LEMMA we have

$M \vDash \forall v(v \, \varepsilon \, z \leftrightarrow w \, \varepsilon \, x \lor v=y)$ and $M \vDash \forall v(v \, \varepsilon \, z \leftrightarrow v \, \varepsilon \, z \lor v=w)$.

$\Longrightarrow M \vDash w \, \varepsilon \, z$,

$\Longrightarrow M \vDash w \, \varepsilon \, x \lor w=y$,

$\Longrightarrow$ from STRONG SET ADJUNCTION LEMMA,

$$M \vDash \sigma^*(x,w,x) \lor \forall v(V^{**}(v) \rightarrow (v \, \varepsilon \, w \leftrightarrow v \, \varepsilon \, y)),$$

as required.

(15.5) <u>Claim</u>: $QT^+ \vdash V^{**}(x) \& V^{**}(y) \to \exists z(V^{**}(z) \& [S(x,y,z)]^{**})$.

Immediate from the STRONG SET ADJUNCTION LEMMA.

(15.6) <u>Claim</u>: $QT^+ \vdash V^{**}(x) \& V^{**}(y) \& V^{**}(z_1) \& V^{**}(z_2) \to$

$$\to [S(x,y,z_1) \& S(x,y,z_2) \to z_1 = z_2]^{**}.$$

Assume $\quad M \vDash V^{**}(x) \& V^{**}(y) \& V^{**}(z_1) \& V^{**}(z_2)$

along with $\quad M \vDash \sigma^*(x,y,z_1) \& \sigma^*(x,y,z_2)$.

From the STRONG SET ADJUNCTION LEMMA we have that

$\quad M \vDash \forall w(w \, \varepsilon \, z_1 \leftrightarrow w \, \varepsilon \, x \lor w=y)$ and $M \vDash \forall w(w \, \varepsilon \, z_2 \leftrightarrow w \, \varepsilon \, x \lor w=y)$,

$\Longrightarrow M \vDash \forall w(w \, \varepsilon \, z_1 \leftrightarrow w \, \varepsilon \, z_2)$,

$\Longrightarrow M \vDash \forall w(V^{**}(w) \to (w \, \varepsilon \, z_1 \leftrightarrow w \, \varepsilon \, z_2))$, as required.

(EXT**) <u>Claim</u>: $QT^+ \vdash V^{**}(x) \& V^{**}(y) \to [\forall z(S(x,z,x) \leftrightarrow S(y,z,y)) \to x=y]^{**}$.

Assume $\quad M \vDash V^{**}(x) \& V^{**}(y)$

along with $\quad M \vDash \forall z(V^{**}(z) \to (\sigma^*(x,z,x) \leftrightarrow \sigma^*(y,z,y)))$.

From the Remark following the proof of the STRONG SET ADJUNCTION

LEMMA we have that

$$M \vDash \forall w(V^{**}(w) \to (\sigma^*(x,w,x) \leftrightarrow w \, \varepsilon \, x))$$

and also $\quad M \vDash \forall w(V^{**}(w) \to (\sigma^*(y,w,y) \leftrightarrow w \, \varepsilon \, y))$.

Hence

$$M \vDash \forall z(V^{**}(z) \to (z \, \varepsilon \, x \leftrightarrow z \, \varepsilon \, y)), \quad \text{as required.}$$

## 16. Interpreting AST$_1$+EXT in QT$^+$

Let AST$_1$ be the first-order theory formulated in the same language as AST consisting of the following non-logical axioms:

NULL.  $\exists x \forall y \neg (y \in x)$

ADJ$_1$.  $\forall x,y \exists z (y \in z \ \& \ \forall w(w \in x \rightarrow w \in z) \ \& \ \forall w(w \in z \rightarrow w \in x \lor w=y))$.

We consider the theory AST$_1$+EXT with the additional axiom

EXT.  $\forall x,y(\forall z(z \in x \leftrightarrow z \in y) \rightarrow x=y)$.

We let the domain of the interpretation $^{+*}$ be defined by $V^{**}(x)$.

Let  $[x=y]^{+*} \equiv x \overset{*}{=} y$.

Let  $[x \in y]^{+*} \equiv x \ \varepsilon \ y$.

The verification of [NULL]$^{+*}$ in QT$^+$ is the same as for Claim 13.1.

Now [EXT]$^{+*}$ is the formula

$\forall x(V^{**}(x) \rightarrow \forall y(V^{**}(y) \rightarrow (\forall z(V^{**}(z) \rightarrow (z \ \varepsilon \ x \leftrightarrow z \ \varepsilon \ y)) \rightarrow$
$\rightarrow \forall z(V^{**}(z) \rightarrow (z \ \varepsilon \ x \leftrightarrow z \ \varepsilon \ y)))))$,

and we immediately have that  QT ⊢ [EXT]$^{+*}$.

Claim.  QT$^+$ ⊢ [ADJ$_1$]$^{+*}$.

We have that [ADJ$_1$]$^{+*}$ is the formula

$$\forall x(V^{**}(x) \to \forall y(V^{**}(y) \to \exists z(V^{**}(z) \ \& \ y \ \varepsilon \ z \ \&$$

$$\& \ \forall w(V^{**}(w) \to (w \ \varepsilon \ x \to w \ \varepsilon \ z)) \ \&$$

$$\& \ \forall w(V^{**}(w) \to (w \ \varepsilon \ z \to w \ \varepsilon \ x \ \lor \ \forall v(V^{**}(v) \to (v \ \varepsilon \ w \leftrightarrow v \ \varepsilon \ y)))))))).$$

Assume $M \vDash V^{**}(x) \ \& \ V^{**}(y)$.

From the STRONG SET ADJUNCTION LEMMA we have that

$$M \vDash \exists z(V^{**}(z) \ \& \ \sigma^*(x,y,z) \ \& \ \forall w(w \ \varepsilon \ z \leftrightarrow w \ \varepsilon \ x \ \lor \ w=y)).$$

$\Rightarrow M \vDash y \ \varepsilon \ z \ \& \ \forall w(V^{**}(w) \to (w \ \varepsilon \ x \to w \ \varepsilon \ z))$.

Assume now that $M \vDash V^{**}(w) \ \& \ w \ \varepsilon \ z$.

$\Rightarrow M \vDash w \ \varepsilon \ x \ \lor \ w=y$,

$\Rightarrow M \vDash w \ \varepsilon \ x \ \lor \ \forall v(V^{**}(v) \to (v \ \varepsilon \ w \leftrightarrow v \ \varepsilon \ y))$.

Hence we also have

$$M \vDash \forall w(V^{**}(w) \to (w \ \varepsilon \ z \to w \ \varepsilon \ x \ \lor \ \forall v(V^{**}(v) \to (v \ \varepsilon \ w \leftrightarrow v \ \varepsilon \ y))))$$

as required.

Therefore $QT^+ \vdash [ADJ_1]^{+*}$.

Hence we have:

<u>The Interpretation Theorem.</u>     $AST+EXT \leq_I QT^+$.

# Index of definitions